# NEUTROSOPHIC BILINEAR ALGEBRAS AND THEIR GENERALIZATIONS

**W. B. Vasantha Kandasamy**
**Florentin Smarandache**

2010

# NEUTROSOPHIC BILINEAR ALGEBRAS AND THEIR GENERALIZATIONS

# CONTENTS

# DEDICATION

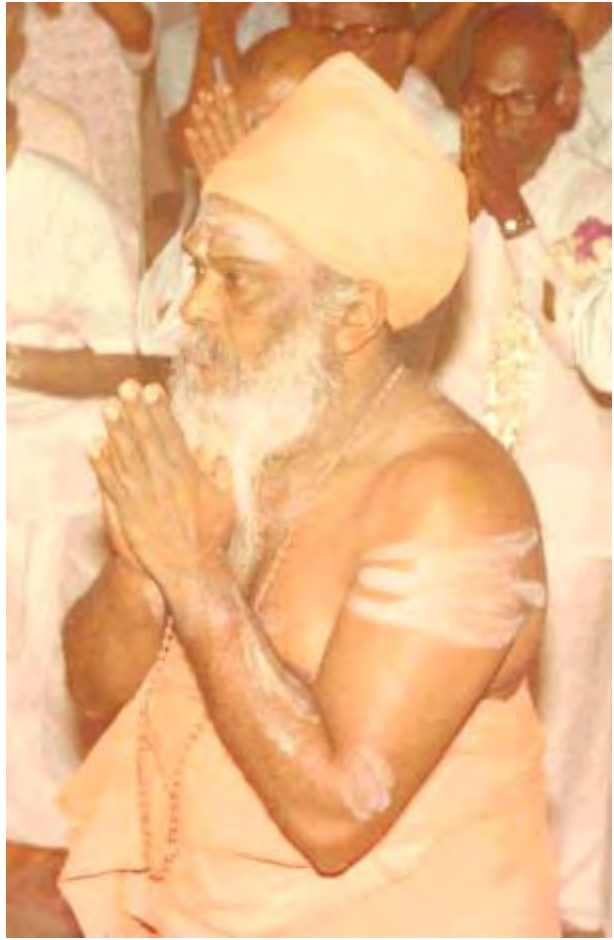

This book is dedicated to the memory of Thavathiru Kundrakudi Adigalar (11 July 1925 - 15 April 1995), a spiritual leader who worked tirelessly for social development and communal harmony. His powerful writings and speeches provided a massive impetus to the promotion of Tamil literature and culture. Moving beyond the realm of religion, he was also actively involved in bringing change at the grassroots level. During his lifetime and after, he was celebrated for his successful efforts to develop poverty-stricken villages around Kundrakudi through planning, support and continued intervention. His initiatives for rural development earned praise even from Indira Gandhi, the then Prime Minister of India. Although he headed a famed religious center, Kundrakudi Adigalar maintained a scientific outlook towards the world and paid special attention to the educational uplift of poor people across caste, community or religion. The services rendered by him to society are endless, and the dedication of this book is a small gesture to pay homage to that great man.

# PREFACE

This book introduces the concept of neutrosophic bilinear algebras and their generalizations to n-linear algebras, n>2.

This book has five chapters. The reader should be well-versed with the notions of linear algebras as well as the concepts of bilinear algebras and n- linear algebras. Further the reader is expected to know about neutrosophic algebraic structures as we have not given any detailed literature about it.

The first chapter is introductory in nature and gives a few essential definitions and references for the reader to make use of the literature in case the reader is not thorough with the basics. The second chapter deals with different types of neutrosophic bilinear algebras and bivector spaces and proves several results analogous to linear bialgebra.

In chapter three the authors introduce the notion of n-linear algebras and prove several theorems related to them. Many of the classical theorems for neutrosophic algebras are proved with appropriate modifications. Chapter four indicates the probable applications of these algebraic structures. The final chapter suggests about 80 innovative problems for the reader to solve.

The interesting feature of this book is that it has over 225 illustrative examples, this is mainly provided to make the reader understand these new concepts. This book contains over 60 theorems and has introduced over 100 new concepts.

The authors deeply acknowledge Dr. Kandasamy for the proof reading and Meena and Kama for the formatting and designing of the book.

W.B.VASANTHA KANDASAMY
FLORENTIN SMARANDACHE

**Chapter One**

# INTRODUCTION TO BASIC CONCEPTS

This chapter has two sections. In section one basic notions about bilinear algebras and n-linear algebras are recalled. In section two an introduction to indeterminacy and algebraic neutrosophic structures essential for this book are given.

## 1.1 Introduction to Bilinear Algebras and their Generalizations

In this section we just recall some necessary definitions about bilinear algebras.

**DEFINITION 1.1.1:** *Let $(G, +, \bullet)$ be a bigroup where $G = G_1 \cup G_2$; bigroup G is said to be commutative if both $(G_1, +)$ and $(G_2, \bullet)$ are commutative.*

**DEFINITION 1.1.2:** *Let $V = V_1 \cup V_2$ where $V_1$ and $V_2$ are two proper subsets of V and $V_1$ and $V_2$ are vector spaces over the same field F that is V is a bigroup, then we say V is a bivector space over the field F.*
*If one of $V_1$ or $V_2$ is of infinite dimension then so is V. If $V_1$ and $V_2$ are of finite dimension so is V; to be more precise if $V_1$ is* *of dimension n and $V_2$ is of dimension m then we define dimension of the bivector space $V = V_1 \cup V_2$ to be of dimension m + n. Thus there exists only m + n elements which are linearly independent and has the capacity to generate $V = V_1 \cup V_2$.*

The important fact is that same dimensional bivector spaces are in general not isomorphic.

***Example 1.1.1:*** Let $V = V_1 \cup V_2$ where $V_1$ and $V_2$ are vector spaces of dimension 4 and 5 respectively defined over rationals where $V_1 = \{(a_{ij}) \mid a_{ij} \in Q\}$, collection of all $2 \times 2$ matrices with entries from Q. $V_2 = \{$Polynomials of degree less than or equal to 4 with coefficients from $Q\}$. Clearly V is a finite dimensional bivector space over Q of dimension 9. In order to avoid confusion we can follow the following convention whenever essential. If $v \in V = V_1 \cup V_2$ then $v \in V_1$ or $v \in V_2$ if $v \in V_1$ then v has a representation of the form $(x_1, x_2, x_3, x_4, 0, 0, 0, 0, 0)$ where $(x_1, x_2, x_3, x_4) \in V_1$ if $v \in V_2$ then $v = (0, 0, 0, 0, y_1, y_2, y_3, y_4, y_5)$ where $(y_1, y_2, y_3, y_4, y_5) \in V_2$.

**DEFINITION 1.1.3:** *Let $V = V_1 \cup V_2$ be a bigroup. If $V_1$ and $V_2$ are linear algebras over the same field F then we say V is a linear bialgebra over the field F.*

*If both $V_1$ and $V_2$ are of infinite dimensional linear algebras over F then we say V is an infinite dimensional linear bialgebra over F. Even if one of $V_1$ or $V_2$ is infinite dimension then we say V is an infinite dimensional linear bialgebra. If both $V_1$ and $V_2$ are finite dimensional linear algebra over F then we say $V = V_1 \cup V_2$ is a finite bidimensional linear bialgebra.*

***Examples 1.1.2:*** Let $V = V_1 \cup V_2$ where $V_1 = \{$set of all $n \times n$ matrices with entries from $Q\}$ and $V_2$ be the polynomial ring $Q[x]$. $V = V_1 \cup V_2$ is a linear bialgebra over Q and the linear bialgebra is an infinite dimensional linear bialgebra.

***Example 1.1.3:*** Let $V = V_1 \cup V_2$ where $V_1 = Q \times Q \times Q$ abelian group under '+', $V_2 = \{$set of all $3 \times 3$ matrices with entries from $Q\}$ then $V = V_1 \cup V_2$ is a bigroup. Clearly V is a linear bialgebra over Q. Further dimension of V is 12; V is a 12 dimensional linear bialgebra over Q.

The standard basis is {(0 1 0), (1 0 0), (0 0 1)} ∪

$$\left\{ \begin{pmatrix} 1 & 0 & 0 \\ 0 & 0 & 0 \\ 0 & 0 & 0 \end{pmatrix}, \begin{pmatrix} 0 & 1 & 0 \\ 0 & 0 & 0 \\ 0 & 0 & 0 \end{pmatrix}, \begin{pmatrix} 0 & 0 & 1 \\ 0 & 0 & 0 \\ 0 & 0 & 0 \end{pmatrix}, \begin{pmatrix} 0 & 0 & 0 \\ 1 & 0 & 0 \\ 0 & 0 & 0 \end{pmatrix}, \begin{pmatrix} 0 & 0 & 0 \\ 0 & 1 & 0 \\ 0 & 0 & 0 \end{pmatrix}, \right.$$

$$\left. \begin{pmatrix} 0 & 0 & 0 \\ 0 & 0 & 1 \\ 0 & 0 & 0 \end{pmatrix}, \begin{pmatrix} 0 & 0 & 0 \\ 0 & 0 & 0 \\ 1 & 0 & 0 \end{pmatrix}, \begin{pmatrix} 0 & 0 & 0 \\ 0 & 0 & 0 \\ 0 & 1 & 0 \end{pmatrix}, \begin{pmatrix} 0 & 0 & 0 \\ 0 & 0 & 0 \\ 0 & 0 & 1 \end{pmatrix} \right\}$$

**DEFINITION 1.1.4:** *Let $V = V_1 \cup V_2$ be a bigroup. Suppose V is a linear bialgebra over F. A non empty proper subset W of V is said to be a linear subbialgebra of V over F if*
  i.  *$W = W_1 \cup W_2$ is a subbigroup of $V = V_1 \cup V_2$.*
  ii. *$W_1$ is a linear subalgebra over F.*
  iii. *$W_2$ is a linear subalgebra over F.*

For more refer [48, 51-2]. For n-linear algebra of type I and II, refer[54-5].

## 1.2 Introduction to Neutrosophic Algebraic Structures

In this section we just recall some basic neutrosophic algebraic structures essential to make this book a self contained one. For more refer [36-43, 53].

In this section we assume fields to be of any desired characteristic and vector spaces are taken over any field. We denote the indeterminacy by 'I', as i will make a confusion, as it denotes the imaginary value, viz. $i^2 = -1$ that is $\sqrt{-1} = i$. The indeterminacy I is such that $I \cdot I = I^2 = I$.

Here we recall the notion of neutrosophic groups. Neutrosophic groups in general do not have group structure.

**DEFINITION 1.2.1:** *Let (G, \*) be any group, the neutrosophic group is generated by I and G under \* denoted by N(G) = {⟨G ∪ I⟩, \*}.*

**Example 1.2.1:** *Let $Z_7$ = {0, 1, 2, …, 6} be a group under addition modulo 7. N(G) = {⟨$Z_7$ ∪ I⟩, '+' modulo 7} is a neutrosophic group which is in fact a group. For N(G) = {a + bI / a, b ∈ $Z_7$} is a group under '+' modulo 7. Thus this neutrosophic group is also a group.*

**Example 1.2.2:** *Consider the set G = $Z_5$ \ {0}, G is a group under multiplication modulo 5. N(G) = {⟨G ∪ I⟩, under the binary operation, multiplication modulo 5}. N(G) is called the neutrosophic group generated by G ∪ I. Clearly N(G) is not a group, for $I^2$ = I and I is not the identity but only an indeterminate, but N(G) is defined as the neutrosophic group.*

Thus based on this we have the following theorem:

**THEOREM 1.2.1:** *Let (G, \*) be a group, N(G) = {⟨G ∪ I⟩, \*} be the neutrosophic group.*
  1. *N(G) in general is not a group.*
  2. *N(G) always contains a group.*

*Proof:* To prove N(G) in general is not a group it is sufficient we give an example; consider ⟨$Z_5$ \ {0} ∪ I⟩ = G = {1, 2, 4, 3, I, 2 I, 4 I, 3 I}; G is not a group under multiplication modulo 5. In fact {1, 2, 3, 4} is a group under multiplication modulo 5. N(G) the neutrosophic group will always contain a group because we generate the neutrosophic group N(G) using the group G and I. So G $\subsetneq$ N(G); hence N(G) will always contain a group.

Now we proceed onto define the notion of neutrosophic subgroup of a neutrosophic group.

**DEFINITION 1.2.2:** *Let N(G) = ⟨G ∪ I⟩ be a neutrosophic group generated by G and I. A proper subset P(G) is said to be a neutrosophic subgroup if P(G) is a neutrosophic group i.e. P(G) must contain a (sub) group of G.*

***Example 1.2.3***: Let $N(Z_2) = \langle Z_2 \cup I \rangle$ be a neutrosophic group under addition. $N(Z_2) = \{0, 1, I, 1 + I\}$. Now we see $\{0, I\}$ is a group under + in fact a neutrosophic group $\{0, 1 + I\}$ is a group under '+' but we call $\{0, I\}$ or $\{0, 1 + I\}$ only as pseudo neutrosophic groups for they do not have a proper subset which is a group. So $\{0, I\}$ and $\{0, 1 + I\}$ will be only called as pseudo neutrosophic groups (subgroups).

We can thus define a pseudo neutrosophic group as a neutrosophic group, which does not contain a proper subset which is a group. Pseudo neutrosophic subgroups can be found as a substructure of neutrosophic groups. Thus a pseudo neutrosophic group though has a group structure is not a neutrosophic group and a neutrosophic group cannot be a pseudo neutrosophic group. Both the concepts are different.

Now we see a neutrosophic group can have substructures which are pseudo neutrosophic groups which is evident from the following example.

***Example 1.2.4:*** Let $N(Z_4) = \langle Z_4 \cup I \rangle$ be a neutrosophic group under addition modulo 4. $\langle Z_4 \cup I \rangle = \{0, 1, 2, 3, I, 1 + I, 2I, 3I, 1 + 2I, 1 + 3I, 2 + I, 2 + 2I, 2 + 3I, 3 + I, 3 + 2I, 3 + 3I\}$. $o(\langle Z_4 \cup I \rangle) = 4^2$.

Thus neutrosophic group has both neutrosophic subgroups and pseudo neutrosophic subgroups. For $T = \{0, 2, 2 + 2I, 2I\}$ is a neutrosophic subgroup as $\{0\ 2\}$ is a subgroup of $Z_4$ under addition modulo 4. $P = \{0, 2I\}$ is a pseudo neutrosophic group under '+' modulo 4.

**DEFINITION 1.2.3:** *Let K be the field of reals. We call the field generated by $K \cup I$ to be the neutrosophic field for it involves the indeterminacy factor in it. We define $I^2 = I$, $I + I = 2I$ i.e., $I + \ldots + I = nI$, and if $k \in K$ then $k.I = kI$, $0I = 0$. We denote the neutrosophic field by K(I) which is generated by $K \cup I$ that is $K(I) = \langle K \cup I \rangle$. $\langle K \cup I \rangle$ denotes the field generated by K and I.*

***Example 1.2.5:*** Let R be the field of reals. The neutrosophic field of reals is generated by R and I denoted by $\langle R \cup I \rangle$ i.e. R(I) clearly $R \subset \langle R \cup I \rangle$.

*Example 1.2.6:* Let Q be the field of rationals. The neutrosophic field of rationals is generated by Q and I denoted by Q(I).

**DEFINITION 1.2.4:** *Let K(I) be a neutrosophic field we say K(I) is a prime neutrosophic field if K(I) has no proper subfield, which is a neutrosophic field.*

*Example 1.2.7:* Q(I) is a prime neutrosophic field where as R(I) is not a prime neutrosophic field for Q(I) $\subset$ R(I).

**DEFINITION 1.2.5:** *Let K(I) be a neutrosophic field, P $\subset$ K(I) is a neutrosophic subfield of P if P itself is a neutrosophic field. K(I) will also be called as the extension neutrosophic field of the neutrosophic field P.*

We can also define neutrosophic fields of prime characteristic p (p is a prime).

**DEFINITION 1.2.6:** *Let $Z_p = \{0, 1, 2, ..., p-1\}$ be the prime field of characteristic p. $\langle Z_p \cup I \rangle$ is defined to be the neutrosophic field of characteristic p. Infact $\langle Z_p \cup I \rangle$ is generated by $Z_p$ and I and $\langle Z_p \cup I \rangle$ is a prime neutrosophic field of characteristic p.*

*Example 1.2.8:* $Z_7 = \{0, 1, 2, 3, ..., 6\}$ be the prime field of characteristic 7. $\langle Z_7 \cup I \rangle = \{0, 1, 2, ..., 6, I, 2I, ..., 6I, 1 + I, 1 + 2I, ..., 6 + 6I \}$ is the prime field of characteristic 7.

**DEFINITION 1.2.7:** *Let G(I) by an additive abelian neutrosophic group and K any field. If G(I) is a vector space over K then we call G(I) a neutrosophic vector space over K.*

Elements of these neutrosophic fields will also be known as neutrosophic numbers. For more about neutrosophy please refer [36-43]. We see $Z_nI = \{aI \mid a \in Z_n\}$ is a neutrosophic field called pure neutrosophic field. Likewise QI, RI and $Z_pI$ are neutrosophic fields where p is a prime. Thus $Z_5I = \{0, I, 2I, 3I, 4I\}$ is a pure neutrosophic field. For more about neutrosophic vector spaces please refer [53].

**Chapter Two**

# NEUTROSOPHIC LINEAR BIALGEBRA

In this chapter we introduce the notion of neutrosophic linear bialgebras and describe a few properties about them. Strong neutrosophic linear bialgebra are also introduced. This chapter has four sections. In section one, we introduce the new notion of neutrosophic bivector space. Strong neutrosophic bivector spaces are introduced in section two. Section three introduces the notion of neutrosophic bivector space of type III. Section four studies the biinner product in strong neutrosophic bivector space.

## 2.1 Neutrosophic Bivector Spaces

In this section we introduce the notion of neutrosophic bivector spaces and study their properties.

**DEFINITION 2.1.1:** *Let $V = V_1 \cup V_2$ where each $V_i$ is a neutrosophic vector space over the same field F and $V_i \neq V_j$, $V_i \not\subseteq V_j$ and $V_j \not\subseteq V_i$; $1 \leq i, j \leq 2$, then we define V to be a neutrosophic bivector space over the real field F.*

*Note:* We assume here F is just a real field that is F is Q or $Z_n$ or R or C. (n a prime $n < \infty$).

We will illustrate this by some simple examples.

*Example 2.1.1:* Let $V_1 = \langle Q \cup I \rangle = N(Q) = \{a + bI \mid a, b \in Q\}$ be a neutrosophic vector space over Q. Take

$$V_2 = \left\{ \begin{pmatrix} a & b \\ c & d \end{pmatrix} \middle| a, b, c, d \in N(Q) \right\},$$

a neutrosophic vector space over Q. $V = V_1 \cup V_2$ is a neutrosophic bivector space over Q.

*Example 2.1.2:* Let $V = V_1 \cup V_2 = \{N(Q)[x]\} \cup \{(a, b, c) \mid a, b, c \in N(Q)\}$. V is a neutrosophic bivector space over Q.

Now we will define a quasi neutrosophic bivector space.

**DEFINITION 2.1.2:** *Let $V = V_1 \cup V_2$ be such that $V_1$ is a vector space over the real field F and $V_2$ is a neutrosophic vector space over F. We define $V = V_1 \cup V_2$ to be a quasi- neutrosophic bivector space over F.*

We will give some examples of quasi neutrosophic bivector spaces.

*Example 2.1.3:* Let $V = V_1 \cup V_2$ where $V_1 = \{Z_7[x] \mid$ all polynomials in the variable x with coefficients from $Z_7\}$ is a vector space over $Z_7$ and $V_2 = \{(Z_7 I \times Z_7 I \times Z_7 I \times Z_7 I) = \{(a, b, c, d) \mid a, b, c, d \in Z_7 I\}\}$ is a neutrosophic vector space over $Z_7$. Then $V = V_1 \cup V_2$ is a quasi neutrosophic bivector space over $Z_7$.

***Example 2.1.4:*** Let $V = V_1 \cup V_2$ where $V_1 = \{Q \times Q \times Q \times Q \times R\} = \{(a, b, c, d, e) \mid a, b, c, d \in Q$ and $e \in R\}$ is a vector space over Q and $V_2 = \{QI \times Q \times QI \times Q \times QI\} = \{(a, b, c, d, e) \mid a, c, e \in QI$ and $b, d \in Q\}$ is a neutrosophic vector space over Q. Thus $V = V_1 \cup V_2$ is a quasi neutrosophic bivector space over Q.

Now we define substructures of these structures.

**DEFINITION 2.1.3:** *Let $V = V_1 \cup V_2$ be a neutrosophic bivector space over the field F. Let $W = W_1 \cup W_2 \subseteq V_1 \cup V_2$ be such that W is a neutrosophic bivector space over F, then we define W to be a neutrosophic bivector subspace of V over F.*

We will illustrate this by some examples.

***Example 2.1.5:*** Let $V = V_1 \cup V_2$

$$= \{Z_3I \times Z_3I \times Z_3I \times Z_3I\} \cup \left\{ \begin{pmatrix} a & b \\ c & d \end{pmatrix} \middle| a,b,c,d \in N(Z_3) \right\}$$

be a neutrosophic bivector space over the field $Z_3$. Let $W = W_1 \cup W_2$

$$= \{(a, b, 0, 0) \mid a, b \in Z_3I\} \cup \left\{ \begin{pmatrix} a & b \\ c & d \end{pmatrix} \middle| a,b,c,d \in Z_3I \right\}$$

$\subseteq V_1 \cup V_2$; W is a neutrosophic bivector subspace of V over the field $Z_3$.

***Example 2.1.6:*** Let $V = V_1 \cup V_2$

$$= N(Q)[x] \cup \left\{ \begin{pmatrix} a & b & e & g \\ c & d & f & h \end{pmatrix} \middle| a,b,c,d,e,f,g,h \in QI \right\}$$

be a neutrosophic bivector space over the field Q. Let $W = W_1 \cup W_2$

$$= QI[x] \cup \left\{ \begin{pmatrix} a & b & e & g \\ 0 & 0 & 0 & 0 \end{pmatrix} \middle| a,b,e,g \in QI \right\}$$

$\subseteq V_1 \cup V_2$, W is a neutrosophic bivector subspace of V over the field Q.

**DEFINITION 2.1.4:** *Let $V = V_1 \cup V_2$ be a neutrosophic bivector space over the field F. Let $W = W_1 \cup W_2 \subseteq V_1 \cup V_2$ be such that W is only a quasi neutrosophic bivector space over F; that is one of $W_1$ or $W_2$ is only a neutrosophic vector space over F and other is just a vector space over the field F; then we call W to be a pseudo quasi neutrosophic bivector subspace of V over the field F.*

*Example 2.1.7:* Let $V = V_1 \cup V_2$ where $V_1 = (Z_5 I \times Z_5 I \times Z_5 I)$ a neutrosophic vector space over $Z_5$ and $V_2 = N(Z_5)[x]$ a neutrosophic bivector space over $Z_5$. $V = V_1 \cup V_2$ is a neutrosophic bivector space over the field $Z_5$.

Take $W = W_1 \cup W_2 = \{Z_5 I \times \{0\} \times \{0\}\} \cup \{Z_5 I [x]\} \subseteq V_1 \cup V_2$; W is a pseudo quasi neutrosophic bivector subspace of V over $Z_5$.

*Example 2.1.8:* Let $V = V_1 \cup V_2$

$$= \left\{ \begin{pmatrix} a & b & c \\ d & e & f \\ g & h & i \end{pmatrix} \middle| a,b,c,d,e,f,g,h,i \in RI \right\} \cup$$

$$\left\{ \begin{pmatrix} a_1 & a_2 & a_3 & a_4 & a_5 \\ b_1 & b_2 & b_3 & b_4 & b_5 \end{pmatrix} \middle| a_i, b_i \in N(Q); 1 \leq i \leq 5 \right\}$$

be a neutrosophic bivector space over the field Q.

Take $W = W_1 \cup W_2$

$$= \left\{ \begin{pmatrix} a & b & c \\ 0 & d & e \\ 0 & 0 & f \end{pmatrix} \middle| a,b,c,d,e,f \in QI \right\} \cup$$

$$\left\{ \begin{pmatrix} a_1 & a_2 & a_3 & a_4 & a_5 \\ 0 & 0 & 0 & 0 & 0 \end{pmatrix} \middle| a_i \in Q; 1 \le i \le 5 \right\}$$

$\subseteq V_1 \cup V_2$; W is a pseudo quasi neutrosophic bivector subspace of V over the field Q.

**DEFINITION 2.1.5:** *Let $V = V_1 \cup V_2$ be a neutrosophic bivector space over the field F. Suppose $W = W_1 \cup W_2 \subseteq V_1 \cup V_2$ is such that W is just a bivector space over the field F then we define W to be a pseudo bivector subspace of V over the field F.*

We will give some examples of this notion.

***Example 2.1.9:*** Let $V = V_1 \cup V_2 =$

$$\left\{ \begin{pmatrix} a & b \\ c & d \end{pmatrix} \middle| a,b,c,d \in N(Z_{11}) \right\}$$

$\cup$ ($\{N(Z_{11}) \times N(Z_{11}) \times N(Z_{11}) \times N(Z_{11})\}$) be a neutrosophic bivector space over the field $Z_{11}$. Let $W = W_1 \cup W_2 =$

$$\left\{ \begin{pmatrix} a & b \\ c & d \end{pmatrix} \middle| a,b,c,d \in Z_{11} \right\}$$

$\cup \{(a, b, c, d) \mid a, b, c, d \in Z_{11}\} \subseteq V_1 \cup V_2$ be a bivector space over $Z_{11}$. Thus W is a pseudo bivector subspace of V over the field $Z_{11}$.

***Example 2.1.10:*** Let $V = V_1 \cup V_2 =$

$$\left\{ \begin{pmatrix} a_1 & a_2 & a_3 \\ a_4 & a_5 & a_6 \end{pmatrix} \middle| a_i \in N(Q); 1 \le i \le 6 \right\}$$

$\cup \{N(Q) \times N(Q) \times N(Q) \times N(Q)\}$ be a neutrosophic bivector spaces over the field Q.
Take $W = W_1 \cup W_2 =$

$$\left\{ \begin{pmatrix} a_1 & a_2 & a_3 \\ a_4 & a_5 & a_6 \end{pmatrix} \middle| a_i \in Q; 1 \le i \le 6 \right\}$$

$\cup \ (Q \times Q \times \{0\} \times \{0\})\} \subseteq V_1 \cup V_2$; W is a pseudo bivector subspace of V over the field Q.

**DEFINITION 2.1.6:** *Let $V = V_1 \cup V_2$ be a neutrosophic bivector space over the field F. Let $W = W_1 \cup W_2 \subseteq V_1 \cup V_2$ be such that W is a neutrosophic bivector space over the subfield $K \subseteq F$. Then we call W to be a neutrosophic special bivector subspace of V over the subfield K of F.*

We will give some examples.

***Example 2.1.11:*** Let $V = V_1 \cup V_2 = \{RI[x]\} \cup$

$$\left\{ \begin{pmatrix} a & b \\ c & d \end{pmatrix} \middle| a,b,c,d \in RI \right\}$$

be a neutrosophic bivector space over the field R. Take $W = W_1 \cup W_2 = \{QI[x]\} \cup$

$$\left\{ \begin{pmatrix} a & b \\ 0 & d \end{pmatrix} \middle| a,b,d \in RI \right\}$$

$\subseteq V_1 \cup V_2$ ; W is a neutrosophic special bivector subspace of V over the subfield Q of R.

**Example 2.1.12:** Let $V = V_1 \cup V_2 = (RI \times RI \times RI) \cup$

$$\left\{ \begin{pmatrix} a_1 & a_2 & a_3 & a_4 \\ a_5 & a_6 & a_7 & a_8 \end{pmatrix} \middle| a_i \in RI; 1 \leq i \leq 8 \right\}$$

be a neutrosophic bivector space over the field $Q(\sqrt{2})$. Take $W = W_1 \cup W_2 = \{(QI \times QI \times QI)\} \cup$

$$\left\{ \begin{pmatrix} a_1 & a_2 & a_3 & a_4 \\ 0 & 0 & 0 & 0 \end{pmatrix} \middle| a_i \in RI; 1 \leq i \leq 4 \right\}$$

$\subseteq V_1 \cup V_2$, W is a neutrosophic special bivector subspace of V over the subfield $Q \subseteq Q(\sqrt{2})$.

**DEFINITION 2.1.7:** *Let $V = V_1 \cup V_2$ be a neutrosophic bivector space over the field F. If V has no neutrosophic special bivector subspace then we call V to be a neutrosophic special simple bivector space over F.*

**Example 2.1.13:** Let $V = V_1 \cup V_2 = \{QI \times QI\} \cup \{QI[x]\}$ be a neutrosophic bivector space over the field $F = Q$. V is a neutrosophic special simple bivector space over Q as Q has no proper subfield.

**Example 2.1.14:** Let $V = V_1 \cup V_2 = Z_{23}I \times Z_{23}I \times Z_{23}I \times Z_{23}I\} \cup$

$$\left\{ \begin{pmatrix} a & b \\ c & d \end{pmatrix} \middle| a,b,c,d \in Z_{23}I \right\}$$

be a neutrosophic bivector space over the field $Z_{23}$. V is a neutrosophic special simple bivector space over $Z_{23}$ as $Z_{23}$ is a prime field.

In view of these examples we have the following interesting theorem.

**THEOREM 2.1.1:** *Let $V = V_1 \cup V_2$ be a neutrosophic bivector space over the field F. If F is a prime field of characteristic zero or a prime p then V is a neutrosophic special simple bivector space over F.*

*Proof:* Given F is a prime field of characteristic zero or a prime p, so F has no subfields. Thus for no $W = W_1 \cup W_2 \subseteq V_1 \cup V_2$ can be neutrosophic special bivector subspace of $V = V_1 \cup V_2$ as F has no subfield. Hence the claim.

Now we proceed onto define the notion of neutrosophic bilinear algebra.

**DEFINITION 2.1.8:** *Let $V = V_1 \cup V_2$ where both $V_1$ and $V_2$ are neutrosophic linear algebras over the field F, then we define V to be a neutrosophic bilinear algebra over F.*

We will illustrate this by some simple examples.

*Example 2.1.15:* Let $V = V_1 \cup V_2 =$

$$\left\{ \begin{pmatrix} a & b \\ c & d \end{pmatrix} \middle| a, b, c, d \in QI \right\}$$

$\cup$ (QI × QI × QI × QI × QI)} be a neutrosophic bilinear algebra over Q.

*Example 2.1.16:* Let $V = V_1 \cup V_2 =$

$$\left\{ \begin{pmatrix} a & 0 \\ b & c \end{pmatrix} \middle| a, b, c \in Z_{29}I \right\} \cup \{Z_{29}I [x]\},$$

V is a neutrosophic bilinear algebra over the field $Z_{29}$.

Now as in case of neutrosophic bivector spaces we can define the following substructures.

**DEFINITION 2.1.9:** *Let $V = V_1 \cup V_2$ be a neutrosophic bilinear algebra over the field F. If $W = W_1 \cup W_2 \subseteq V_1 \cup V_2$ is a neutrosophic bilinear algebra over the field F then we call W to be a neutrosophic bilinear subalgebra of V over the field F.*

We give an example.

*Example 2.1.17:* Let $V = V_1 \cup V_2 =$

$$\left\{ \begin{pmatrix} a & b & c \\ 0 & d & e \\ 0 & 0 & f \end{pmatrix} \middle| a,b,c,d,e,f \in QI \right\}$$

$\cup$ {QI $\times$ QI $\times$ QI $\times$ QI} be a neutrosophic bilinear algebra over the field Q.
   Choose

$$W = \left\{ \begin{pmatrix} a & 0 & 0 \\ 0 & b & 0 \\ 0 & 0 & c \end{pmatrix} \middle| a,b,c \in QI \right\}$$

$\cup$ {QI $\times$ {0} $\times$ {0} $\times$ QI} = $W_1 \cup W_2 \subseteq V_1 \cup V_2$; W is a neutrosophic bilinear subalgebra of V over Q.

**DEFINITION 2.1.10:** *Let $V = V_1 \cup V_2$ be a neutrosophic bilinear algebra over the field F. Let $W = W_1 \cup W_2 \subseteq V_1 \cup V_2$ be a neutrosophic bilinear algebra over a subfield $K \subseteq F$; then we define W to be a neutrosophic special bilinear subalgebra of V over the subfield K of F. If V has no special neutrosophic bilinear subalgebra's then we call V to be a special neutrosophic simple bilinear algebra or neutrosophic special simple bilinear algebra.*

We will illustrate this by some simple examples.

***Example 2.1.18:*** Let $V = V_1 \cup V_2 =$

$$\left\{ \begin{pmatrix} a & b \\ c & d \end{pmatrix} \middle| a,b,c,d \in RI \right\}$$

$\cup$ {RI[x]} be a neutrosophic bilinear algebra over the field R. Take $W = W_1 \cup W_2 =$

$$\left\{ \begin{pmatrix} a & b \\ c & d \end{pmatrix} \middle| a,b,c,d \in QI \right\}$$

$\cup$ {QI[x]} $\subseteq V_1 \cup V_2$ ; W is a neutrosophic special bilinear algebra over the subfield Q of the field R.

***Example 2.1.19:*** Let $V = V_1 \cup V_2$

$$= \left\{ \begin{pmatrix} a_1 & a_2 & a_3 \\ a_4 & a_5 & a_6 \\ a_7 & a_8 & a_9 \end{pmatrix} \middle| a_i \in QI; 1 \le i \le 9 \right\}$$

$\cup$ {RI [x]} be a neutrosophic bilinear algebra over Q. Clearly V is a neutrosophic special simple bilinear algebra.

***Example 2.1.20:*** Let $V = V_1 \cup V_2 = \{Z_7I\ [x]\} \cup \{Z_7I \times Z_7I \times Z_7I\}$ be a neutrosophic bilinear algebra over the field $Z_7$. V is a neutrosophic simple bilinear algebra.

In view of these examples we have the following theorem, the proof of which is left as an exercise for the reader.

**THEOREM 2.1.2:** *Let $V = V_1 \cup V_2$ be a neutrosophic bilinear algebra over the field F, where F is a prime field (i.e., F has no subfields other than itself). V is a neutrosophic special simple bilinear algebra.*

**DEFINITION 2.1.11:** *Let $V = V_1 \cup V_2$ be a neutrosophic bilinear algebra over a field F. Suppose $W = W_1 \cup W_2 \subseteq V_1 \cup V_2$ and if W is only a bilinear algebra over the field F, then we call W to be a pseudo bilinear subalgebra of V over the field F.*

*Example 2.1.21:* Let $V = V_1 \cup V_2 =$

$$\left\{ \begin{pmatrix} a & b \\ c & d \end{pmatrix} \middle| a, b, c, d \in N(Q) \right\}$$

$\cup \{Q \times QI \times QI \times Q\}$ be a neutrosophic bilinear algebra over Q, where $W = W_1 \cup W_2 =$

$$\left\{ \begin{pmatrix} a & b \\ c & d \end{pmatrix} \middle| a, b, c, d \in Q \right\}$$

$\cup \{Q \times \{0\} \times \{0\} \times \{0\}\} \subseteq V_1 \cup V_2$; W is a pseudo bilinear subalgebra of V over the field F.

*Example 2.1.22:* Let $V = V_1 \cup V_2 =$

$$\left\{ \begin{pmatrix} a & b & c \\ d & e & f \\ g & h & k \end{pmatrix} \middle| a, b, c, d, e, f, g, h, k \in N(R) \right\}$$

$\cup \{(a, b, c, d) \mid a, b, c, d \in N(Q)\}$ be a neutrosophic bilinear algebra over the field Q. Take $W = W_1 \cup W_2 =$

$$\left\{ \begin{pmatrix} a & b & c \\ 0 & d & e \\ 0 & 0 & f \end{pmatrix} \middle| a, b, c, d, e, f \in R \right\}$$

∪ {(a, b, 0, d) | a, b, d ∈ Q} ⊆ V$_1$ ∪ V$_2$; W is a pseudo bilinear subalgebra of V over Q.

**DEFINITION 2.1.12:** *Let V = V$_1$ ∪ V$_2$ be a neutrosophic bilinear algebra over the field F. Let W = W$_1$ ∪ W$_2$ ⊆ V$_1$ ∪ V$_2$ be a proper bisubset of V which is just a neutrosophic bivector space over the field F. We define W to be a pseudo neutrosophic bivector subspace of V over F.*

We will illustrate this by some simple examples.

*Example 2.1.23:* Let V = V$_1$ ∪ V$_2$ =

$$\left\{ \sum_{i=0}^{n} a_i x^i \,\middle|\, a_i \in QI, i = 0,1,2,\ldots, n \leq \infty \right\} \cup$$

$$\left\{ \begin{pmatrix} a & b \\ c & d \end{pmatrix} \,\middle|\, a,b,c,d \in QI \right\}$$

be a neutrosophic bilinear algebra over the field Q. Let W = W$_1$ ∪ W$_2$ =

$$\left\{ \sum_{i=0}^{5} a_i x^i \,\middle|\, a_i \in QI, i = 0,1,2,\ldots,5 \right\} \cup \left\{ \begin{pmatrix} 0 & a \\ b & 0 \end{pmatrix} \,\middle|\, a,b \in QI \right\}$$

⊆ V$_1$ ∪ V$_2$, W is a pseudo neutrosophic bivector subspace of V over Q.

*Example 2.1.24:* Let V = V$_1$ ∪ V$_2$ =

$$\left\{ \begin{pmatrix} a & b & c \\ d & e & f \\ g & h & i \end{pmatrix} \,\middle|\, a,b,c,d,e,f,g,h,i \in Z_{11}I \right\}$$

$$\cup \left\{ \sum_{i=0}^{n} a_i x^i \,\middle|\, a_i \in Z_{11}I, i = 0,1,2,\ldots,n; n \leq \infty \right\}$$

be a neutrosophic bilinear algebra over the field $Z_{11}$. Let $W = W_1 \cup W_2 =$

$$\left\{ \begin{pmatrix} a & b & c \\ d & 0 & 0 \\ 0 & 0 & 0 \end{pmatrix} \,\middle|\, a,b,c,d \in Z_{11}I \right\} \cup \left\{ \sum_{i=1}^{6} a_0 x^i \,\middle|\, a_i \in Z_{11}I, 1 \leq i \leq 6 \right\}$$

$\subseteq V_1 \cup V_2$, $W$ is only a pseudo neutrosophic bivector subspace of $V$ as we see

$$\begin{pmatrix} a & b & c \\ d & 0 & 0 \\ 0 & 0 & 0 \end{pmatrix} \begin{pmatrix} a & b & c \\ d & 0 & 0 \\ 0 & 0 & 0 \end{pmatrix} = \begin{pmatrix} a^2 + bd & ab & ac \\ ad & bd & cd \\ 0 & 0 & 0 \end{pmatrix} \notin W_1.$$

Similarly if we take $a = 2x^6 + 3x + 1$ and $b = (2Ix^6 + 3Ix + I)(3x^4 + 2x^2 + I) = 6Ix^{10} + 4Ix^8 + 2Ix^6 + 9Ix^5 + 6Ix^3 + 3Ix + 3Ix^4 + 2Ix^2 + I \notin W_2$ but $a, b \in W_2$. Thus $W$ is only a neutrosophic bivector subspace of $V$ and not a neutrosophic bilinear subalgebra of $V$ over $Z_{11}$.

We have the following interesting theorem, the proof of which is left as an exercise for the reader.

**THEOREM 2.1.3:** *Let $V = V_1 \cup V_2$ be a neutrosophic bilinear algebra over the field $F$. $V$ is clearly a neutrosophic bivector space over the field $F$. If $V$ is a neutrosophic bivector space over the field $F$ then in general $V$ is not a neutrosophic bilinear algebra over the field $F$.*

**DEFINITION 2.1.13:** *Let $V = V_1 \cup V_2$ where $V_1$ is only a neutrosophic linear algebra over the field $F$ and $V_2$ is just a linear algebra over $F$ then we define $V = V_1 \cup V_2$ to be a quasi neutrosophic bilinear algebra over the field $F$.*

We illustrate this by some examples.

*Example 2.1.25:* Let $V = V_1 \cup V_2 =$

$$\left\{ \begin{pmatrix} a & b \\ c & d \end{pmatrix} \middle| a,b,c,d \in R \right\} \cup \{(a\ b\ c\ d\ e\ f) \mid a, b, c, d, e, f \in QI\}$$

be a quasi neutrosophic bilinear algebra over the field Q.

*Example 2.1.26:* Let $V = V_1 \cup V_2 =$

$$\left\{ \begin{pmatrix} a & b \\ c & d \end{pmatrix} \middle| a,b,c,d \in QI \right\} \cup \left\{ \sum_{i=0}^{n} a_i x^i \middle| a_i \in Q \right\};$$

V is a quasi neutrosophic bilinear algebra over the field Q.

**DEFINITION 2.1.14:** *Let $V = V_1 \cup V_2$ where $V_1$ is a neutrosophic vector space over the field F and $V_2$ is a neutrosophic linear algebra over the field F. $V = V_1 \cup V_2$ is defined to be a pseudo neutrosophic quasi bilinear algebra over F.*

We will illustrate this by some simple examples.

*Example 2.1.27:* Let $V = V_1 \cup V_2 =$

$$\left\{ \begin{pmatrix} a & b & e \\ c & d & f \end{pmatrix} \middle| a,b,c,d,e,f \in QI \right\}$$

$$\cup \left\{ \begin{pmatrix} a & b \\ c & d \end{pmatrix} \middle| a,b,c,d \in N(Q) \right\}.$$

V is a pseudo neutrosophic quasi bilinear algebra over the field Q.

*Example 2.1.28:* Let V = $V_1 \cup V_2$ =

$$\left\{ \sum_{i=0}^{8} a_i x^i \,\middle|\, a_i \in QI; 0 \leq i \leq 8 \right\} \cup$$

$$\left\{ \begin{pmatrix} a & b & c \\ d & e & f \\ g & h & i \end{pmatrix} \,\middle|\, a,b,c,d,e,f,g,h,i \in RI \right\};$$

V is a pseudo neutrosophic quasi bilinear algebra over the field Q.

We can have for any neutrosophic bilinear algebra V a substructure which is a pseudo neutrosophic quasi bilinear subalgebra of V.

**DEFINITION 2.1.15:** *Let V = $V_1 \cup V_2$ be a neutrosophic bilinear algebra over a field F. Suppose W = $W_1 \cup W_2 \subseteq V_1 \cup V_2$ such that $W_1$ is a neutrosophic vector space over F and $W_2$ is a neutrosophic linear algebra over F then we define W to be a pseudo neutrosophic quasi bilinear subalgebra of V over the field F.*

We will illustrate this by some simple examples.

*Example 2.1.29:* Let V = $V_1 \cup V_2$ =

$$\left\{ \begin{pmatrix} a & b \\ c & d \end{pmatrix} \,\middle|\, a,b,c,d \in QI \right\}$$

$\cup$ {(a, b, c, d, e, f) | a, b, c, d, e, f $\in$ QI} be a neutrosophic bilinear algebra over the field Q. Take W = $W_1 \cup W_2$ =

$$\left\{ \begin{pmatrix} 0 & a \\ b & 0 \end{pmatrix} \,\middle|\, a,b \in QI \right\} \cup \{(a\ 0\ c\ 0\ e\ 0) \mid a, c, e \in QI\}$$

$\subseteq V_1 \cup V_2$ W is a pseudo neutrosophic quasi bilinear subalgebra of V over the field Q.

**Example 2.1.30:** Let $V = V_1 \cup V_2 =$

$$\left\{ \sum_{i=0}^{n} a_i x^i \,\middle|\, a_i \in Z_{17}I; 0 \leq i \leq n \leq \infty \right\} \cup$$

$$\left\{ \begin{pmatrix} a & b & c \\ d & e & f \\ g & h & i \end{pmatrix} \,\middle|\, a,b,c,d,e,f,g,h,i \in Z_{17}I \right\}$$

be a neutrosophic bilinear algebra over the field $Z_{17}$. Choose $W = W_1 \cup W_2 =$

$$\left\{ \sum_{i=0}^{n} a_i x^{2i} \,\middle|\, a_i \in Z_{17}I; 0 \leq i \leq n \leq \infty \right\} \cup$$

$$\left\{ \begin{pmatrix} a & 0 & 0 \\ b & 0 & 0 \\ c & e & 0 \end{pmatrix} \,\middle|\, a,b,c,e \in Z_{17}I \right\}$$

$\subseteq V_1 \cup V_2$; W is only a pseudo neutrosophic quasi bilinear subalgebra of V over the field $Z_{17}$.

We see

$$\begin{pmatrix} a & 0 & 0 \\ b & 0 & 0 \\ c & e & 0 \end{pmatrix} \in W_2$$

But

$$\begin{pmatrix} a & 0 & 0 \\ b & 0 & 0 \\ c & e & 0 \end{pmatrix} \begin{pmatrix} a & 0 & 0 \\ b & 0 & 0 \\ c & e & 0 \end{pmatrix} \notin W_2.$$

Now we proceed onto define linear bitransformation of a neutrosophic bivector space and neutrosophic bilinear algebra over the field F.

**DEFINITION 2.1.16:** *Let $V = V_1 \cup V_2$ be a neutrosophic bivector space over a field F and $W = W_1 \cup W_2$ be another neutrosophic bivector space over the same field F.*

*Define $T = T_1 \cup T_2 : V = V_1 \cup V_2 \to W = W_1 \cup W_2$ as follows $T_i: V_i \to W_i$, $i = 1, 2$ is just a neutrosophic linear transformation from $V_i$ to $W_i$. This $T = T_1 \cup T_2$ is a neutrosophic linear bitransformation of V into W. If W = V then we call the neutrosophic linear bitransformation as neutrosophic linear bioperator. We denote it by $BN_F(V,W) = \{$set of all neutrosophic linear bitransformations of $V = V_1 \cup V_2$ to $W = W_1 \cup W_2\}$; $BN_F(V,W) = BN_F(V_1, W_1) \cup BN_F(V_2, W_2)$. $BN_F(V,V) = \{$set of all neutrosophic linear bioperators of V to V and $BN_F(V_1, V_1) \cup BN_F(V_2, V_2) = BN_F(V,V)$.*

Interested reader can study the algebraic structures of $BN_F(V,W)$ and $BN_F(V,V)$. However we give an example of each.

***Example 2.1.31:*** Let $V = V_1 \cup V_2 =$

$$\left\{ \begin{pmatrix} a & b \\ c & 0 \end{pmatrix} \middle| a,b,c \in QI \right\} \cup \left\{ \begin{pmatrix} a & b & c \\ d & e & f \end{pmatrix} \middle| a,b,c,d,e,f \in QI \right\}$$

and $W = W_1 \cup W_2 = \{QI \times QI \times QI\} \cup$

$$\left\{ \sum_{i=0}^{5} a_i x^i \middle| a_i \in QI; 1 \leq i \leq 5 \right\}$$

be neutrosophic bivector spaces over the field Q. Define $T = T_1 \cup T_2 : V = V_1 \cup V_2 \to W = W_1 \cup W_2$ where $T_1 : V_1 \to W_1$ and $T_2 : V_2 \to W_2$ as follows:

$$T_1 \begin{pmatrix} a & b \\ c & 0 \end{pmatrix} = (a, b, c)$$

and

$$T_2 \begin{pmatrix} a & b & c \\ d & e & f \end{pmatrix} = \{a + bx + cx^2 + dx^3 + ex^4 + fx^5 \mid a, b, c, d, e, f \in QI\}.$$

Clearly $T = T_1 \cup T_2$ is a neutrosophic linear bitransformation of V to W.

*Example 2.1.32:* Let $V = V_1 \cup V_2 =$

$$\left\{ \begin{pmatrix} a & b \\ c & d \end{pmatrix} \middle| a, b, c, d \in Z_7 I \right\}$$

$\cup \{Z_7 I \times Z_7 I \times Z_7 I \times Z_7 I\}$ be a neutrosophic bivector space over the field $Z_7$. $T = T_1 \cup T_2 : V = V_1 \cup V_2 \to V = V_1 \cup V_2$, where

$$T_1 : V_1 \to V_1$$

and

$$T_2 : V_2 \to V_2$$

is as follows.

$$T_1 \begin{pmatrix} a & b \\ c & d \end{pmatrix} = \begin{pmatrix} c & d \\ a & b \end{pmatrix}$$

and

$$T_2(a, b, c, d) = (a, b + c, d, a + d).$$

It is easily verified. T is a neutrosophic linear bioperator on V.

## 2.2 Strong Neutrosophic Bivector Spaces

In this section we for the first time introduce the notion of strong neutrosophic bivector spaces and study them.

**DEFINITION 2.2.1:** *Let $V = V_1 \cup V_2$ where $V_1$ and $V_2$ are neutrosophic additive abelian groups. Suppose $V = V_1 \cup V_2$ is a neutrosophic bivector space over a neutrosophic field F then we call V to be a strong neutrosophic bivector space.*

We will give some examples.

**Example 2.2.1:** Let $V = V_1 \cup V_2 = \{Z_5I\,[x]\} \cup \{Z_5I \times Z_5I \times Z_5I\}$ be a strong neutrosophic bivector space over the neutrosophic field $Z_5I$.

**Example 2.2.2:** Let $V = V_1 \cup V_2 =$

$$\left\{ \begin{pmatrix} a & b \\ c & d \end{pmatrix} \middle| a,b,c,d \in RI \right\} \cup$$

$$\left\{ \begin{pmatrix} a & b & c & d & e & f \\ g & h & i & j & k & l \end{pmatrix} \middle| a,b,c,d,e,f,g,h,i,j,k,l \in QI \right\}$$

be a strong neutrosophic bivector space over the neutrosophic field QI.

**DEFINITION 2.2.2:** *Let $V = V_1 \cup V_2$ be a strong neutrosophic bivector space over the neutrosophic field K. If $W = W_1 \cup W_2 \subseteq V_1 \cup V_2$; and if W is a strong neutrosophic bivector space over the neutrosophic field K, then we call W to be a strong neutrosophic bivector subspace of V over the neutrosophic field F.*

We will illustrate this by some simple examples.

**Example 2.2.3:** Let $V = V_1 \cup V_2 =$

$$\left\{ \begin{pmatrix} a_1 & a_2 & a_3 \\ a_4 & a_5 & a_6 \end{pmatrix} \middle| a_i \in QI; 1 \le i \le 6 \right\}$$

$$\cup \left\{ \begin{pmatrix} a_1 & a_2 \\ a_3 & a_4 \\ a_5 & a_6 \\ a_7 & a_8 \end{pmatrix} \middle| a_i \in QI; 1 \le i \le 8 \right\}$$

be a strong neutrosophic bivector space over the neutrosophic field QI. Take $W = W_1 \cup W_2$

$$= \left\{ \begin{pmatrix} a_1 & 0 & a_3 \\ 0 & a_4 & 0 \end{pmatrix} \middle| a_1, a_3, a_4 \in QI \right\} \cup$$

$$\left\{ \begin{pmatrix} a_1 & a_2 \\ 0 & 0 \\ a_5 & a_6 \\ 0 & 0 \end{pmatrix} \middle| a_1, a_2, a_5, a_6 \in QI \right\}$$

$\subseteq V_1 \cup V_2$; W is a strong neutrosophic bivector subspace of V over the field QI.

***Example 2.2.4:*** Let $V = V_1 \cup V_2 = \{QI \times QI \times QI\} \cup \{QI[x]\}$ be a strong neutrosophic bivector space over the neutrosophic field QI. Take $W = W_1 \cup W_2 = \{(QI \times \{0\} \times QI\} \cup$

$$\left\{ \sum_{i=0}^{8} a_i x^i \middle| a_i \in QI; 0 \le i \le 8 \right\}$$

$\subseteq V_1 \cup V_2$, W is a strong neutrosophic bivector subspace of V over QI.
Let us define the notion of strong neutrosophic bilinear algebra.

**DEFINITION 2.2.3:** *Let $V = V_1 \cup V_2$ be a neutrosophic bilinear algebra over the neutrosophic field K, we define V to be strong neutrosophic bilinear algebra over K.*

We will illustrate this by examples.

*Example 2.2.5:* Let $V = V_1 \cup V_2 = \{QI[x]\} \cup$

$$\left\{ \begin{pmatrix} a & b \\ 0 & d \end{pmatrix} \middle| a, b, d \in QI \right\}$$

be the strong neutrosophic bilinear algebra over the neutrosophic field QI.

*Example 2.2.6:* Let $V = V_1 \cup V_2 = \{QI \times QI \times QI \times QI \times QI\} \cup$

$$\left\{ \begin{pmatrix} a & b \\ c & d \end{pmatrix} \middle| a, b, c, d \in N(Q) \right\}$$

be a strong neutrosophic bilinear algebra over the neutrosophic field QI.

*Example 2.2.7:* Let $V = V_1 \cup V_2 = \{Z_{11}I \times Z_{11}I \times Z_{11}I\} \cup$

$$\left\{ \begin{pmatrix} a & b & c \\ d & e & f \\ g & h & i \end{pmatrix} \middle| a, b, c, d, e, f, g, h, i \in N(Z_{11}) \right\}$$

be a strong neutrosophic bilinear algebra over the neutrosophic field $Z_{11}I$.

**DEFINITION 2.2.4:** *Let $V = V_1 \cup V_2$ be a strong neutrosophic bilinear algebra over the neutrosophic field K. If $W = W_1 \cup W_2 \subseteq V_1 \cup V_2$ be a strong neutrosophic bilinear algebra over K then we define W to be a strong neutrosophic bilinear subalgebra of V over the neutrosophic field K.*

*Example 2.2.8:* Let $V = V_1 \cup V_2 =$

$$\left\{ \begin{pmatrix} a & b \\ c & d \end{pmatrix} \middle| a, b, c, d \in QI \right\}$$

∪ {(a, b, c, d, e, f) | a, b, c, d, e, f ∈ QI} be a strong neutrosophic bilinear algebra over the neutrosophic field K = QI. Take W = $W_1 \cup W_2$ =

$$\left\{ \begin{pmatrix} a & 0 \\ 0 & d \end{pmatrix} \middle| a, d \in QI \right\}$$

∪ {(a, 0, 0, d, 0, f) | a, d, f ∈ QI} ⊆ $V_1 \cup V_2$; W is a strong neutrosophic bilinear subalgebra of V over K = QI.

*Example 2.2.9:* Let V = $V_1 \cup V_2$ = {$Z_7I$ [x]} ∪

$$\left\{ \begin{pmatrix} a & b \\ c & d \end{pmatrix} \middle| a, b, c, d \in Z_7I \right\}$$

be a strong neutrosophic bilinear algebra over the neutrosophic field $Z_7I$.
Take W = $W_1 \cup W_2$ =

$$\left\{ \sum_{i=0}^{n} a_i x^{2i} \middle| a_i \in Z_7I; i = 0, 1, 2, ..., \infty \right\} \cup \left\{ \begin{pmatrix} a & a \\ a & a \end{pmatrix} \middle| a \in Z_7I \right\}$$

⊆ $V_1 \cup V_2$; W is a strong neutrosophic bilinear subalgebra of V over the neutrosophic field $Z_7I$.

**DEFINITION 2.2.5:** *Let V = $V_1 \cup V_2$ be a strong neutrosophic bilinear algebra over the neutrosophic field K. Let W = $W_1 \cup W_2 \subseteq V_1 \cup V_2$ where $W_1$ is just a neutrosophic vector space over the neutrosophic field K and $W_2$ is a neutrosophic linear subalgebra over the neutrosophic field K. W = $W_1 \cup W_2$ is defined to be a pseudo strong neutrosophic linear subalgebra of V over K.*

We will illustrate this situation by some examples.

***Example 2.2.10:*** Let $V = V_1 \cup V_2 =$

$$\left\{ \begin{pmatrix} a_1 & a_2 & a_3 \\ a_4 & a_5 & a_6 \\ a_7 & a_8 & a_9 \end{pmatrix} \middle| a_i \in Z_7I; 1 \leq i \leq 9 \right\}$$

$\cup$ $\{Z_7I\ [x]\}$ be a neutrosophic bilinear algebra over the neutrosophic field $Z_7I$.
Take $W = W_1 \cup W_2 =$

$$\left\{ \begin{pmatrix} a & 0 & 0 \\ 0 & 0 & b \\ 0 & c & 0 \end{pmatrix} \middle| a, b, c \in Z_7I \right\} \cup$$

$$\left\{ \sum_{i=0}^{n} a_i x^{2i} \middle| 0 \leq i \leq n; i = 0, 1, 2, ..., \infty \right\}$$

$\subseteq V_1 \cup V_2$; W is a pseudo strong neutrosophic bilinear subalgebra of V over the field $Z_7I$.

***Example 2.2.11:*** Let $V = V_1 \cup V_2 =$

$$\left\{ \begin{pmatrix} a & b \\ c & d \end{pmatrix} \middle| a, b, c, d \in Z_{23}I \right\}$$

$\cup$ $\{Z_{23}I \times Z_{23}I \times Z_{23}I \times Z_{23}I\}$ be a strong neutrosophic bilinear algebra over the neutrosophic field $Z_{23}I$. Take $W = W_1 \cup W_2 =$

$$\left\{ \begin{pmatrix} 0 & d \\ a & 0 \end{pmatrix} \middle| a, d \in Z_{23}I \right\}$$

$\cup$ $\{(0\ 0\ 0\ b)\ |\ a, b \in Z_{23}I\} \subseteq V_1 \cup V_2$, W is a pseudo strong neutrosophic bilinear subalgebra of V over $Z_{23}I$.

***Example 2.2.12:*** Let $V = V_1 \cup V_2 = \{QI[x]\} \cup$

$$\left\{ \begin{pmatrix} a & b & c \\ d & e & f \\ g & h & i \end{pmatrix} \middle| a,b,c,d,e,f,g,h,i \in QI \right\}$$

be a strong neutrosophic bilinear algebra over the neutrosophic field QI. Take $W = W_1 \cup W_2 =$

$$\left\{ \sum_{i=0}^{n} a_i x \middle| a_i \in QI \right\} \cup \left\{ \begin{pmatrix} a & 0 & 0 \\ 0 & b & 0 \\ 0 & 0 & d \end{pmatrix} \middle| a,b,d \in QI \right\}$$

$\subseteq V_1 \cup V_2$; W is a pseudo strong neutrosophic bilinear subalgebra of V over the field QI.

**DEFINITION 2.2.6:** *Let $V = V_1 \cup V_2$ be a strong neutrosophic bilinear algebra over the neutrosophic field K. Let $W = W_1 \cup W_2 \subseteq V_1 \cup V_2$ be a strong neutrosophic bivector space over the neutrosophic field K. W is defined as the strong pseudo neutrosophic bivector subspace of V over the field K.*

We illustrate this by some examples.

***Example 2.2.13:*** Let $V = V_1 \cup V_2 =$

$$\left\{ \begin{pmatrix} a & b & c \\ d & e & f \\ g & h & i \end{pmatrix} \middle| a,b,c,d,e,f,g,h,i \in Z_{13}I \right\} \cup$$

$$\left\{ \sum_{i=0}^{\infty} a_i x^i \middle| a_i \in Z_{13}I; i = 1,2,...,\infty \right\}$$

be a strong neutrosophic bilinear algebra over the neutrosophic field $Z_{23}I$. Take $W = W_1 \cup W_2 =$

$$\left\{ \begin{pmatrix} a & b & c \\ d & e & f \\ 0 & 0 & 0 \end{pmatrix} \middle| a,b,c,d,e,f \in Z_{23}I \right\} \cup$$

$$\left\{ \sum_{i=0}^{9} a_i x^i \middle| a_i \in Z_{23}I; 0 \leq i \leq 9 \right\}$$

$\subseteq V_1 \cup V_2$, W is a strong pseudo neutrosophic bivector subspace of V over the field $Z_{23}I$.

*Example 2.2.14:* Let $V = V_1 \cup V_2 =$

$$\left\{ \begin{pmatrix} a_1 & a_2 & a_3 & a_4 \\ a_5 & a_6 & a_7 & a_8 \\ a_9 & a_{10} & a_{11} & a_{12} \\ a_{13} & a_{14} & a_{15} & a_{16} \end{pmatrix} \middle| a_i \in Z_{23}I; 1 \leq i \leq 16 \right\}$$

$$\cup \left\{ \sum_{i=0}^{\infty} a_i x^i \middle| a_i \in Z_{11}I; i = 1,2,...,\infty \right\}$$

be a strong neutrosophic bilinear algebra over the neutrosophic field $Z_{11}I$. Take $W = W_1 \cup W_2 =$

$$\left\{ \begin{pmatrix} a_1 & a_2 & a_3 & a_4 \\ a_5 & a_6 & a_7 & a_8 \\ 0 & 0 & 0 & 0 \\ 0 & 0 & 0 & 0 \end{pmatrix} \middle| a_i \in Z_{11}I; 1 \leq i \leq 8 \right\} \cup$$

$$\left\{\sum_{i=0}^{6} a_i x^i \,\middle|\, a_i \in Z_{11}I;\ i = 1, 2, \ldots, 6\right\}$$

$\subseteq V_1 \cup V_2$, W is pseudo strong neutrosophic bivector subspace of V over $Z_{11}I$.

**DEFINITION 2.2.7:** *Let $V = V_1 \cup V_2$ be a strong neutrosophic bilinear algebra over the neutrosophic field K. Take $W = W_1 \cup W_2 \subseteq V_1 \cup V_2$ and $F \subseteq K$ (F a field and is not a neutrosophic subfield of K). If W is a neutrosophic bilinear algebra over the field F then we define W to be a pseudo strong neutrosophic bilinear subalgebra of V over the subfield F of the neutrosophic field K.*

We will illustrate this by some examples.

*Example 2.2.15:* Let $V = V_1 \cup V_2 =$

$$\left\{\begin{pmatrix} a & b \\ c & d \end{pmatrix} \,\middle|\, a, b, c, d \in N(Q)\right\} \cup \left\{\sum_{i=0}^{n} a_i x^i \,\middle|\, a_i \in N(Q); 0 \leq i \leq \infty\right\}$$

be a strong neutrosophic bilinear subalgebra of V over the neutrosophic field QI. Take $W = W_1 \cup W_2 =$

$$\left\{\begin{pmatrix} a & a \\ a & a \end{pmatrix} \,\middle|\, a \in QI\right\} \cup \left\{\sum_{i=0}^{\infty} a_i x^i \,\middle|\, a_i \in QI; i = 0, 1, \ldots, \infty\right\}$$

$\subseteq V_1 \cup V_2$, W is a pseudo strong neutrosophic bilinear subalgebra of V over the subfield Q of N(Q).

**DEFINITION 2.2.8:** *Let $V = V_1 \cup V_2$ be a strong neutrosophic bilinear algebra over the neutrosophic field K. Let $W = W_1 \cup W_2 \subseteq V$ be a bivector space over the real field $F \subseteq K$. We call $W = W_1 \cup W_2 \subseteq V_1 \cup V_2$ as a pseudo bivector subspace of V over the real subfield F of K.*

We will illustrate this by some simple examples.

***Example 2.2.16:*** Let $V = V_1 \cup V_2 =$

$$\left\{ \begin{pmatrix} a & b \\ c & d \end{pmatrix} \middle| a, b, c, d \in N(Q) \right\} \cup \left\{ \sum_{i=0}^{\infty} a_i x^i \middle| a_i \in N(Q); i = 0, \ldots, \infty \right\}$$

be a neutrosophic bilinear algebra over the neutrosophic field $N(Q)$. Take $W = W_1 \cup W_2 =$

$$\left\{ \begin{pmatrix} a & b \\ 0 & 0 \end{pmatrix} \middle| a, b \in Q \right\} \cup \left\{ \sum_{i=0}^{8} a_i x^i \middle| a_i \in Q; i = 0, 1, 2, \ldots, 8 \right\}$$

$\subseteq V_1 \cup V_2$, $W$ is a pseudo bivector space over the field $Q$.

***Example 2.2.17:*** Let $V = V_1 \cup V_2 =$

$$\left\{ \begin{pmatrix} a & 0 & 0 \\ 0 & b & 0 \\ 0 & 0 & c \end{pmatrix} \middle| a, b, c \in Z_{29}I \right\}$$

$\cup \ \{(N(Z_{29}) \times N(Z_{29}) \times N(Z_{29}) \times N(Z_{29}))\}$ be a strong neutrosophic bilinear algebra over the neutrosophic field $N(Z_{29})$. Take $W = W_1 \cup W_2 =$

$$\left\{ \begin{pmatrix} a & 0 & 0 \\ 0 & b & 0 \\ 0 & 0 & c \end{pmatrix} \middle| a, b, c \in Z_{29} \right\} \cup \{(Z_{29} \times Z_{29} \times \{0\} \times \{0\})\}$$

$\subseteq V_1 \cup V_2$; $W$ is a pseudo bivector subspace of $V$ over the field $Z_{29}$.

**DEFINITION 2.2.9:** *Let $V = V_1 \cup V_2$ be a strong neutrosophic bilinear algebra over the neutrosophic field $K$. Let $W = W_1 \cup$* $W_2 \subseteq V_1 \cup V_2$ *be a bilinear algebra over a real subfield F of K. We define W to be a pseudo bilinear subalgebra of V over the field F.*

We will illustrate this by some simple examples.

***Example 2.2.18:*** Let $V = V_1 \cup V_2 =$

$$\left\{ \begin{pmatrix} a & b \\ c & d \end{pmatrix} \middle| a,b,c,d \in N(Q) \right\} \cup$$

$$\left\{ \sum_{i=0}^{\infty} a_i x^i \middle| a_i \in N(Q); i = 0, 1, \ldots, \infty \right\}$$

be a strong neutrosophic bilinear algebra over the neutrosophic field $N(Q)$. Take $W = W_1 \cup W_2 =$

$$\left\{ \begin{pmatrix} a & b \\ c & d \end{pmatrix} \middle| a,b,c,d \in Q \right\} \cup \left\{ \sum_{i=0}^{\infty} a_i x^i \middle| a_i \in Q; i = 0, 1, \ldots, \infty \right\}$$

$\subseteq V_1 \cup V_2$; W is a pseudo bilinear subalgebra of V over Q.

***Example 2.2.19:*** Let $V = V_1 \cup V_2 =$

$$\left\{ \begin{pmatrix} a & b & c & d \\ 0 & d & e & f \\ 0 & 0 & g & h \\ 0 & 0 & 0 & i \end{pmatrix} \middle| a,b,c,d,e,f,g,h,i \in N(Z_{17}) \right\}$$

$\cup \{(N(Z_{17}) \times N(Z_{17}) \times N(Z_{17}) \times N(Z_{17}) \times N(Z_{17}))\}$ be a strong neutrosophic bilinear algebra over the neutrosophic field $N(Z_{17})$. Take $W = W_1 \cup W_2 =$

$$\left\{ \begin{pmatrix} a & 0 & 0 & 0 \\ 0 & b & 0 & 0 \\ 0 & 0 & c & 0 \\ 0 & 0 & 0 & d \end{pmatrix} \middle| \, a,b,c,d \in Z_{17} \right\}$$

$\cup \, \{Z_{17} \times \{0\} \times Z_{17} \times \{0\} \times Z_{17}\} \subseteq V_1 \cup V_2$; W is a pseudo bilinear subalgebra of V over the field $Z_{17} \subseteq N(Z_{17})$.

**DEFINITION 2.2.10:** *Let $V = V_1 \cup V_2$ be a strong neutrosophic bivector space over the neutrosophic field K. Let $W = W_1 \cup W_2$ be a strong neutrosophic bivector space over the same neutrosophic field K. Let $T : V \to W$ i.e., $T = T_1 \cup T_2 : V_1 \cup V_2 \to W_1 \cup W_2$ be a bimap such that $T_i : V_i \to W_i$ is a strong neutrosophic linear transformation from $V_i$ to $W_i$; $i = 1,2$. We define $T = T_1 \cup T_2$ to be a strong neutrosophic linear bitransformation from V to W. If $W = V$ then we call T to be a strong neutrosophic linear bioperator on V.*

*$SNHom_K (V, W)$ denotes the set of all strong neutrosophic linear bitransformations from V to W.*

*$SNHom_K (V, V)$ denotes the set of all strong neutrosophic linear bioperator from V to V.*

Interested reader is requested to give examples.
Also the study of substructure preserving strong neutrosophic linear bitransformations (bioperators) is an interesting field of research.

Now we proceed onto define bilinearly independent bivectors and other related properties.

**DEFINITION 2.2.11:** *Let $V = V_1 \cup V_2$ be a strong neutrosophic bivector space over the neutrosophic field K. A proper bisubset $S = S_1 \cup S_2 \subseteq V_1 \cup V_2$ is said to be a bibasis of V if S is a bilinearly independent biset and each $S_i \subseteq V_i$ generates $V_i$; that is $S_i$ is a basis of $V_i$ true for $i = 1, 2$.*

**DEFINITION 2.2.12:** *Let $V = V_1 \cup V_2$ be a strong neutrosophic bivector space over the neutrosophic field K. Let $X = X_1 \cup X_2 \subseteq V_1 \cup V_2$ be a biset of V, we say X is a linearly biindependent bisubset of V over K if each of the subsets $X_i$ contained in $V_i$ is a linearly independent subset of $V_i$ over the K; i = 1, 2.*

The reader is expected to prove the following:

**THEOREM 2.2.1:** *Let $V = V_1 \cup V_2$ be a strong neutrosophic bivector space over the neutrosophic field K. Let $B = B_1 \cup B_2$ be a bibasis of V over K then B is a linearly biindependent subset of V over K. If $X = X_1 \cup X_2$ be a bisubset of V which is bilinearly independent bisubset of V then X in general need not be a bibasis of V over K.*

We will explain this by some examples.

***Example 2.2.20:*** Let $V = V_1 \cup V_2 = \{(QI \times QI \times QI)\} \cup$

$$\left\{ \sum_{i=0}^{\infty} a_i x^i \,\middle|\, 0 \leq i \leq n \leq \infty; a_i \in QI \right\}$$

be a strong neutrosophic bivector space over the neutrosophic field QI. Let $B = B_1 \cup B_2 = \{(I, 0, 0), (0, I, 0), (0, 0, I)\} \cup \{I, Ix, Ix^2, \ldots, Ix^n, \ldots, Ix^\infty\} \subseteq V_1 \cup V_2$ be a bibasis of V over the neutrosophic field QI. Take $X = X_1 \cup X_2 = \{I, 0, 2I), (0, 3I, I)\} \cup \{I, Ix, Ix^2, Ix^3, Ix^7\} \subseteq V_1 \cup V_2$; X is a linearly independent bisubset of V but is not a bibasis of V over QI.

***Example 2.2.21:*** Let $V = V_1 \cup V_2 =$

$$\left\{ \begin{pmatrix} a & b \\ c & d \end{pmatrix} \,\middle|\, a,b,c,d \in Z_{13}I \right\} \cup \left\{ \begin{pmatrix} a_1 & a_2 & a_3 \\ a_4 & a_5 & a_6 \end{pmatrix} \,\middle|\, a_i \in Z_{13}I; 1 \leq i \leq 6 \right\}$$

be a strong neutrosophic bivector space over the neutrosophic field $Z_{13}I$. Let $B = B_1 \cup B_2 =$

$$\left\{ \begin{pmatrix} I & 0 \\ 0 & 0 \end{pmatrix}, \begin{pmatrix} 0 & I \\ 0 & 0 \end{pmatrix}, \begin{pmatrix} 0 & 0 \\ I & 0 \end{pmatrix}, \begin{pmatrix} 0 & 0 \\ 0 & I \end{pmatrix} \right\} \cup$$

$$\left\{ \begin{pmatrix} I & 0 & 0 \\ 0 & 0 & 0 \end{pmatrix}, \begin{pmatrix} 0 & I & 0 \\ 0 & 0 & 0 \end{pmatrix}, \begin{pmatrix} 0 & 0 & I \\ 0 & 0 & 0 \end{pmatrix}, \right.$$

$$\left. \begin{pmatrix} 0 & 0 & 0 \\ I & 0 & 0 \end{pmatrix}, \begin{pmatrix} 0 & 0 & 0 \\ 0 & I & 0 \end{pmatrix}, \begin{pmatrix} 0 & 0 & 0 \\ 0 & 0 & I \end{pmatrix} \right\}$$

$\subseteq V_1 \cup V_2$, B is a bibasis of V over $Z_{13}I$. Take X =

$$\left\{ \begin{pmatrix} I & I \\ 0 & 0 \end{pmatrix}, \begin{pmatrix} 0 & 0 \\ I & I \end{pmatrix} \right\} \cup$$

$$\left\{ \begin{pmatrix} 3I & 0 & I \\ 0 & 0 & 0 \end{pmatrix}, \begin{pmatrix} 0 & I & 4I \\ I & I & 0 \end{pmatrix}, \begin{pmatrix} 0 & I & 0 \\ 2I & 0 & 4I \end{pmatrix} \right\}$$

$= X_1 \cup X_2 \subseteq V_1 \cup V_2$, X is only a linearly independent biset of V but is not a bibasis of V over $Z_{13}I$.

**DEFINITION 2.2.13:** *Let $V = V_1 \cup V_2$ be a strong neutrosophic bivector space over the neutrosophic field K. Let $X = X_1 \cup X_2 \subseteq V_1 \cup V_2$, if X is not a bilinearly independent bisubset of V then we say X is a bilinearly dependent bisubset of V.*

*Example 2.2.22:* Let $V = V_1 \cup V_2 = \{QI \times QI \times QI \times QI\} \cup$

$$\left\{ \sum_{i=0}^{5} a_i x^i \,\middle|\, a_i \in QI; 0 \le i \le 5 \right\}$$

be a strong neutrosophic bivector space over the neutrosophic field QI. Let $X = X_1 \cup X_2 = \{(I, I, 0, 0), (0, I, I, 0), (0, 0, I, I), (I, I, I, I), (3I, 2I, I, 0)\} \cup \{I, Ix^2, 1 + 3Ix^3, 5Ix^3 + 3Ix^2, Ix^5 + 3Ix +$ $5Ix^2 + 3Ix^4\} \subseteq V_1 \cup V_2$. It is easily verified X is a linearly dependent bisubset of V over QI.

***Example 2.2.23:*** Let $V = V_1 \cup V_2 =$

$$\left\{ \begin{pmatrix} a & b \\ c & d \end{pmatrix} \middle| a,b,c,d \in Z_2I \right\} \cup \left\{ \sum_{i=0}^{\infty} a_i x^i \middle| 0 \leq i \leq \infty; a_i \in Z_2I \right\}$$

be a strong neutrosophic bilinear algebra over the neutrosophic field $Z_2I$. $B = B_1 \cup B_2 =$

$$\left\{ \begin{pmatrix} I & 0 \\ 0 & 0 \end{pmatrix}, \begin{pmatrix} 0 & I \\ 0 & 0 \end{pmatrix}, \begin{pmatrix} 0 & 0 \\ I & 0 \end{pmatrix}, \begin{pmatrix} 0 & 0 \\ 0 & I \end{pmatrix} \right\} \cup \{I, Ix, Ix^2, \ldots, Ix^n, \ldots\}$$

is a bibasis of B.

$$\left\{ \begin{pmatrix} I & I \\ 0 & I \end{pmatrix}, \begin{pmatrix} I & 0 \\ 0 & 0 \end{pmatrix}, \begin{pmatrix} I & I \\ I & I \end{pmatrix}, \begin{pmatrix} I & I \\ I & I \end{pmatrix}, \begin{pmatrix} 0 & I \\ 0 & 0 \end{pmatrix} \right\}$$

$\cup \{I + Ix^2 + Ix^3 + Ix^2, Ix^2, I, Ix, I + Ix^2\}$ is a linearly dependent bisubset of V over $Z_2I$. The number of bielements in the bibasis $B = B_1 \cup B_2$ is the bidimension of $V = V_1 \cup V_2$, denoted by $|B| = (|B_1|, |B_2|)$.

If $|B| = (|B_1|, |B_2|) = (n, m)$ and if $n < \infty$ and $m < \infty$ then we say V is a finite bidimensional strong neutrosophic bilinear algebra (bivector space) over the neutrosophic field K. Even if one of m or n is $\infty$ or both m and n is infinite then we say the bidimension of V is infinite.

***Example 2.2.24:*** Let $V = V_1 \cup V_2 =$

$$\left\{ \begin{pmatrix} a & b \\ c & d \end{pmatrix} \middle| a,b,c,d \in QI \right\} \cup \{QI \times QI \times QI\}$$

be a strong neutrosophic bivector space over the neutrosophic field QI. $B = B_1 \cup B_2 =$

$$\left\{ \begin{pmatrix} I & 0 \\ 0 & 0 \end{pmatrix}, \begin{pmatrix} 0 & I \\ 0 & 0 \end{pmatrix}, \begin{pmatrix} 0 & 0 \\ I & 0 \end{pmatrix}, \begin{pmatrix} 0 & 0 \\ 0 & I \end{pmatrix} \right\}$$

$\cup$ {(I 0 0), (0, I, 0) (0, 0, I)} $\subseteq V_1 \cup V_2$; B is a bibasis of V over QI and the bidimension of V is finite (4, 3).

*Example 2.2.25:* Let $V = V_1 \cup V_2 =$

$$\left\{ \sum_{i=0}^{\infty} a_i x^i \,\middle|\, a_i \in Z_2 I \right\} \cup \{(Z_2 I \times Z_2 I)\}$$

be a strong neutrosophic bivector space over $Z_2 I$. $B = B_1 \cup B_2 = \{I, Ix, Ix^2, \ldots, Ix^n, \ldots \infty\} \cup \{(I, 0), (0, I)\} \subseteq V_1 \cup V_2$ is a bibasis of V over $Z_2 I$. The bidimension of V is $(\infty, 2)$.

*Example 2.2.26:* Let $V = V_1 \cup V_2 =$

$$\left\{ \sum_{i=0}^{\infty} a_i x^i \,\middle|\, a_i \in QI; i = 1, 2, \ldots, \infty \right\} \cup$$

$$\left\{ \begin{pmatrix} a & b \\ c & d \end{pmatrix} \,\middle|\, a, b, c, d \in RI \right\}$$

be a strong neutrosophic bilinear algebra over the neutrosophic field QI. $B = B_1 \cup B_2 = \{I, Ix, Ix^2, \ldots, Ix^n, \ldots\} \cup$ {an infinite basis for $V_2$} is a bibasis of V over QI. Thus the bidimension of V is infinite and $|B| = (\infty, \infty)$.

*Example 2.2.27:* Let $V = V_1 \cup V_2 =$

$$\left\{ \begin{pmatrix} a & b \\ c & d \end{pmatrix} \,\middle|\, a, b, c, d \in QI \right\} \cup \{RI \times RI\}$$

be a strong neutrosophic bivector space over the neutrosophic field QI. Take $B = B_1 \cup B_2 =$

$$\left\{ \begin{pmatrix} I & 0 \\ 0 & 0 \end{pmatrix}, \begin{pmatrix} 0 & I \\ 0 & 0 \end{pmatrix}, \begin{pmatrix} 0 & 0 \\ I & 0 \end{pmatrix}, \begin{pmatrix} 0 & 0 \\ 0 & I \end{pmatrix} \right\}$$

$\cup$ {An infinite set}, B is a bibasis of V over QI. The bidimension of V is $(4, \infty)$; thus the bidimension of V is infinite.
It is interesting to note that if V and W are strong neutrosophic bivector spaces over the neutrosophic field K. Suppose bidimension of V is $(n_1, n_2)$ then we say the bidimension of V and W are the same if and only if W is just of bidimension $(n_1, n_2)$ or $(n_2, n_1)$.

**DEFINITION 2.2.14:** *Let $V = V_1 \cup V_2$ and $W = W_1 \cup W_2$ be two strong neutrosophic bivector spaces over the neutrosophic field K. Let $T = T_1 \cup T_2$ be a bilinear transformation (linear bitransformation) from V to W defined by $T_i : V_i \to W_j$, $i = 1, 2$, $j = 1, 2$, such that $T_1 : V_1 \to W_1$ and $T_2 : V_2 \to W_2$ or $T_1 : V_1 \to W_2$ and $T_2 : V_2 \to W_1$. The bikernel of T denoted by $kerT = kerT_1 \cup kerT_2$ where $ker\, T_i = \{v^j \in V_i \,|\, T(v^j) = 0; i = 1, 2\}$. Thus biker $T = \{(v^1, v^2) \in V_1 \cup V_2 / T(v^1, v^2) = T_1(v^1) \cup T(v^2) = 0 \cup 0\}$.*

It is easily verified ker T is a proper neutrosophic bisubgroup of V. Further ker T is a strong neutrosophic bisubspace of V.

*Example 2.2.28:* Let $V = V_1 \cup V_2 =$

$$\left\{ \begin{pmatrix} a_1 & a_2 & a_3 \\ a_4 & a_5 & a_6 \end{pmatrix} \middle| a_i \in QI, 1 \le i \le 6 \right\} \cup$$

$$\left\{ \begin{pmatrix} a_1 & a_2 \\ a_3 & a_4 \\ a_5 & a_6 \\ a_7 & a_8 \end{pmatrix} \middle| a_i \in QI, 1 \le i \le 8 \right\}$$

be a strong neutrosophic bivector space over QI. $W = W_1 \cup W_2$ =

$$\left\{\sum_{i=0}^{7} a_i x^i \,\middle|\, a_i \in QI; 0 \leq i \leq 7\right\} \cup \left\{\begin{pmatrix} a_1 & a_2 & a_3 \\ 0 & a_4 & a_5 \\ 0 & 0 & a_6 \end{pmatrix} \,\middle|\, a_i \in QI, 0 \leq i \leq 6\right\}$$

be a strong neutrosophic bivector space over QI. Define a bimap $T = T_1 \cup T_2 : V_1 \cup V_2 \to W_1 \cup W_2$ by $T_1 : V_1 \to W_2$ and $T_2 : V_2 \to W_1$ such that

$$T_1 \begin{pmatrix} a_1 & a_2 & a_3 \\ a_4 & a_5 & a_6 \end{pmatrix} = \begin{pmatrix} a_1 & a_2 & a_3 \\ 0 & a_4 & a_5 \\ 0 & 0 & a_6 \end{pmatrix}$$

and

$$T_2 \begin{pmatrix} a_1 & a_2 \\ a_3 & a_4 \\ a_5 & a_6 \\ a_7 & a_8 \end{pmatrix} = \sum_{i=0}^{7} a_i x^i$$

where $a_1 \to a_0$, $a_2 \to a_1$, $a_3 \to a_2$, $a_4 \to a_3$, $a_5 \to a_4$, $a_6 \to a_5$, $a_7 \to a_6$ and $a_8 \to a_7$.

$T = T_1 \cup T_2$ is a linear bimap.

$$\text{biker } T = \left\{\begin{pmatrix} 0 & 0 & 0 \\ 0 & 0 & 0 \\ 0 & 0 & 0 \end{pmatrix}\right\} \cup \left\{\begin{pmatrix} 0 & 0 \\ 0 & 0 \\ 0 & 0 \\ 0 & 0 \end{pmatrix}\right\}.$$

Interested reader can construct more examples in which biker T is a proper non zero neutrosophic bisubspace of V. We will prove results when we define strong neutrosophic n-vector spaces $n > 2$, for $n = 2$ gives the strong neutrosophic bivector space. Further neutrosophic bivector spaces (bilinear algebras) and the strong neutrosophic bivector spaces (bilinear algebras) which we have defined in sections 2.1 and 2.2 are type I neutrosophic bivector spaces and strong neutrosophic bivector spaces respectively. In the following section we define type II neutrosophic bivector spaces (bilinear algebras).

## 2.3 Neutrosophic Bivector Spaces of Type II

In this section we proceed onto define neutrosophic bivector spaces of type II and neutrosophic linear bialgebras (or bilinear algebras) of type II. We discuss several interesting properties about them. We also give the difference between type I and type II neutrosophic bivector spaces.

**DEFINITION 2.3.1:** *Let $V = V_1 \cup V_2$ where $V_1$ is a neutrosophic vector space over the real field $F_1$ and $V_2$ is a neutrosophic vector space over the real field $F_2$ such that $F_1 \neq F_2$, $F_1 \not\subseteq F_2$, $F_2 \not\subseteq F_1$ and $V_1 \neq V_2$, $V_1 \not\subseteq V_2$ and $V_2 \not\subseteq V_1$.*
*We call V to be a neutrosophic bivector space over the bifield $F = F_1 \cup F_2$ of type II.*

We will illustrate this by some simple examples.

*Example 2.3.1:* Let $V = V_1 \cup V_2$ where

$$V_1 = \left\{ \begin{pmatrix} a & b \\ c & d \end{pmatrix} \middle| a,b,c,d \in Z_7 I \right\}$$

be a neutrosophic vector space over the field $Z_7$ and

$$V_2 = \left\{ \begin{pmatrix} a_1 & a_2 \\ a_3 & a_4 \\ a_5 & a_6 \end{pmatrix} \middle| a_i \in N(Q), 1 \leq i \leq 6 \right\}$$

is a neutrosophic vector space over the field Q. $V = V_1 \cup V_2$ is a neutrosophic bivector space over the bifield $F = Z_7 \cup Q$ of type II.

***Example 2.3.2:*** Let $V = V_1 \cup V_2$ where $V_1 = \{QI[x]\}$ a neutrosophic vector space over the field Q and $V_2 =$

$$\left\{ \begin{pmatrix} a & b & 0 \\ 0 & d & e \\ 0 & 0 & f \end{pmatrix} \middle| a,b,c,d,e,f \in Z_{11} \right\}$$

be a neutrosophic vector space over the field $Z_{11}$. $V = V_1 \cup V_2$ is a neutrosophic bivector space over the bifield $F = Q \cup Z_{11}$ of type II.

***Example 2.3.3:*** Let $V = V_1 \cup V_2$ where $V_1 = \{Z_{13}I \times Z_{13}I \times Z_{13}I \times Z_{13}I\}$ is a neutrosophic vector space over the field $Z_{13}$ and

$$V_2 = \left\{ \begin{pmatrix} a_1 & a_2 & a_3 & a_4 \\ a_5 & a_6 & a_7 & a_8 \end{pmatrix} \middle| a_i \in Z_{23}, 1 \le i \le 8 \right\}$$

be a neutrosophic vector space over the field $Z_{23}$. $V = V_1 \cup V_2$ is a neutrosophic bivector space over the bifield $F = Z_{13} \cup Z_{23}$ of type II.

**DEFINITION 2.3.2:** *Let $V = V_1 \cup V_2$ be a neutrosophic bivector space over the bifield $F = F_1 \cup F_2$ of type II. Let $W = W_1 \cup W_2 \subseteq V_1 \cup V_2$, if W is a neutrosophic bivector space over the bifield $F = F_1 \cup F_2$ of type II, then we call W to be a neutrosophic bivector subspace of V over the bifield $F = F_1 \cup F_2$ of type II.*

We will illustrate this by examples.

***Example 2.3.4:*** Let $V = V_1 \cup V_2 =$

$$\left\{ \begin{pmatrix} a & b \\ c & d \end{pmatrix} \middle| a,b,c,d \in Z_7 I \right\} \cup$$

$$\left\{ \begin{pmatrix} a_1 & a_2 & a_3 & a_4 \\ a_5 & a_6 & a_7 & a_8 \end{pmatrix} \middle| a_i \in Z_{11}I, 1 \le i \le 8 \right\}$$

be a neutrosophic bivector space of V over the bifield $F = Z_7 \cup Z_{11}$ of type II. Take $W = W_1 \cup W_2 =$

$$\left\{ \begin{pmatrix} a & b \\ 0 & 0 \end{pmatrix} \middle| a, b \in Z_7 I \right\} \cup$$

$$\left\{ \begin{pmatrix} a_1 & a_2 & a_3 & a_4 \\ a_5 & a_6 & 0 & 0 \end{pmatrix} \middle| a_i \in Z_{11}I, i = 1, 2, 4, 5, 6, 3 \right\}$$

$\subseteq V_1 \cup V_2$; W is a neutrosophic bivector subspace of V over the bifield $Z_7 \cup Z_{11}$ of type II.

***Example 2.3.5:*** Let $V = V_1 \cup V_2 =$

$$\left\{ \begin{pmatrix} a & b \\ c & d \end{pmatrix} \middle| a, b, c, d \in QI \right\}$$

$\cup \{Z_{13}I \times Z_{13}I \times Z_{13}I \times Z_{13}I \times Z_{13}I\}$ be a neutrosophic bivector space over the bifield $F = Q \cup Z_{13}$ of type II. Take $W = W_1 \cup W_2 =$

$$\left\{ \begin{pmatrix} a & a \\ a & a \end{pmatrix} \middle| a \in QI \right\} \cup \{(a\ a\ a\ a\ a) \mid a \in Z_{13}I\}$$

$\subseteq V_1 \cup V_2$; W is a neutrosophic bivector subspace of V over the bifield $F = Q \cup Z_{13}$.

Now we define a substructures on these neutrosophic bivector spaces over the bifield. It is pertinent to mention here that the term type II will be suppressed as one can easily understand by the very definition it is distinct from type I.

**DEFINITION 2.3.3:** *Let $V = V_1 \cup V_2$ be a neutrosophic bivector space over the real bifield $F = F_1 \cup F_2$. Let $W = W_1 \cup W_2 \subseteq V_1 \cup V_2$ and $K = K_1 \cup K_2 \subseteq F_1 \cup F_2 = F$. If $W$ is a neutrosophic bivector space over the bifield $K = K_1 \cup K_2$ then we call $W$ to be a special subneutrosophic bivector subspace of $V$ over the bisubfield $K$ of $F$.*

We will give an example of this definition.

***Example 2.3.6:*** Let $V = V_1 \cup V_2 = \{Q(\sqrt{2}, \sqrt{3})I \times Q(\sqrt{2}, \sqrt{3})I\} \cup$

$$\left\{ \begin{pmatrix} a & b \\ c & d \end{pmatrix} \middle| a,b,c,d \in Q(\sqrt{5}, \sqrt{7}) \right\}$$

be a neutrosophic bivector space over the bifield $F = Q(\sqrt{2}, \sqrt{3}) \cup Q(\sqrt{5}, \sqrt{7}) = F$. Take $W = \{Q(\sqrt{2})I \times Q(\sqrt{2})I\} \cup$

$$\left\{ \begin{pmatrix} a & a \\ a & a \end{pmatrix} \middle| a \in Q(\sqrt{5}) \right\}$$

$= W_1 \cup W_2 \subseteq V_1 \cup V_2$, $W$ is a special subneutrosophic bivector subspace of $V$ over the subfield $Q(\sqrt{2}) \cup Q(\sqrt{5}) = K_1 \cup K_2 \subseteq Q(\sqrt{2}, \sqrt{3}) \cup Q(\sqrt{5}, \sqrt{7}) = F$.

Now we define the neutrosophic bivector space $V$ to be bisimple if $V$ has no proper special subneutrosophic bivector subspace over a bisubfield.

We will illustrate this by some examples.

***Example 2.3.7:*** Let $V = V_1 \cup V_2 =$

$$\{RI\} \cup \left\{ \begin{pmatrix} a & b \\ c & d \end{pmatrix} \middle| a,b,c,d \in Z_7 I \right\}$$

be a neutrosophic bivector space over the real bifield $F = Q \cup Z_7$. We see the real bifield is bisimple; i.e., it has no subbifields or bisubfields. So V is a bisimple neutrosophic bivector space over F.

***Example 2.3.8:*** Let $V = V_1 \cup V_2 = \{Z_2I \times Z_2I \times Z_2I \times Z_2I\} \cup$

$$\left\{ \begin{pmatrix} a & b \\ c & d \end{pmatrix} \middle| a, b, c, d \in Z_3I \right\}$$

be a neutrosophic bivector space over the real bifield $F = Z_2 \cup Z_3$. V is a bisimple neutrosophic bivector space over F. We see both $Z_2$ and $Z_3$ are prime fields of characteristic two and three respectively.

In view of this we have the following theorem.

**THEOREM 2.3.1:** *Let $V = V_1 \cup V_2$ be a neutrosophic bivector space over a bifield $F = F_1 \cup F_2$. If both $F_1$ and $F_2$ are prime fields then V is a bisimple neutrosophic bivector space over the real bifield $F = F_1 \cup F_2$.*

The proof of the above theorem is left as an exercise to the reader. A natural question arise; if one of the fields $F_1$ and $F_2$ alone is a prime field can we have some special type of substructures. In view of this we have the following definition.

**DEFINITION 2.3.4**: *Let $V = V_1 \cup V_2$ be a neutrosophic bivector space over the real bifield $F = F_1 \cup F_2$ where one of $F_1$ or $F_2$ is a prime field. Let $W = W_1 \cup W_2$ be such that $W_1$ is a neutrosophic vector subspace of $V_1$ over $K_1 \subseteq F_1$ ($F_2$ is a prime field ) and $W_2$ is a neutrosophic vector subspace of $V_2$ over $F_2$; then we call $W = W_1 \cup W_2$ to be a quasi special neutrosophic bivector subspace of V over the quasi bisubfield $K_1 \cup F_2$.*

(If $F = F_1 \cup F_2$ is a bifield, $K_1 \subseteq F_1$ is a proper subfield of $F_1$ then $K_1 \cup F_2$ is called the quasi bisubfield of the bifield $F = F_1 \cup F_2$). We will illustrate this by some examples.

***Example 2.3.9:*** Let $V = V_1 \cup V_2 =$

$$\left\{ \begin{pmatrix} a & b \\ c & d \end{pmatrix} \middle| a,b,c,d \in Z_{17}I \right\}$$

$\cup \{RI \times RI \times RI\}$ be a neutrosophic bivector space over the bifield $F = Z_{17} \cup Q(\sqrt{2},\sqrt{3},\sqrt{5},\sqrt{7},\sqrt{11})$. Take $W = W_1 \cup W_2$

$$= \left\{ \begin{pmatrix} a & b \\ c & d \end{pmatrix} \middle| a,b,c,d \in Z_{17}I \right\} \cup$$

$\{RI \times \{0\} \times RI\} \subseteq V_1 \cup V_2$, W is a quasi special neutrosophic bivector subspace of V over the quasi bisubfield $Z_{17} \cup Q(\sqrt{2})$ of the bifield F.

***Example 2.3.10:*** Let $V = V_1 \cup V_2 = \{Z_7I[x]\} \cup \{RI \times RI \times RI\}$ be a neutrosophic bivector space over the real bifield $Z_7 \cup R$. Let $W = W_1 \cup W_2 =$

$$\left\{ \sum_{i=0}^{9} a_i x^i \middle| 0 \leq i \leq 9; a_i \in Z_7I \right\} \cup \{QI \times QI \times QI\}$$

$\subseteq V_1 \cup V_2$; W is a quasi special neutrosophic bivector subspace of V over the real quasi bifield $Z_7 \cup Q \subseteq Z_7 \cup R$.

Now we proceed on to define the notion of bibasis of the neutrosophic bivector space of type II.

**DEFINITION 2.3.5:** *Let $V = V_1 \cup V_2$ be a neutrosophic bivector space of type II over the bifield $F = F_1 \cup F_2$. Let $B = B_1 \cup B_2 \subseteq V_1 \cup V_2$ be a bisubset of V such that $B_i$ is a linearly independent bisubset of $V_i$; and generates $V_i$ for i = 1, 2, then we call B to be bibasis of V over the bifield $F_1 \cup F_2 = F$.*

We will illustrate this by some simple examples.

***Example 2.3.11:*** Let $V = V_1 \cup V_2 =$

$$\left\{ \begin{pmatrix} a_1 & a_2 & a_3 \\ a_4 & a_5 & a_6 \end{pmatrix} \middle| a_i \in Z_7I;\ 0 \leq i \leq 6 \right\} \cup$$

$$\left\{ \begin{pmatrix} a_1 & a_2 \\ a_3 & a_4 \\ a_5 & a_6 \\ a_7 & a_8 \end{pmatrix} \middle| a_i \in Z_5I;\ 1 \leq i \leq 8 \right\}$$

be a neutrosophic bivector space of type II over the bifield $F = Z_7 \cup Z_5$.

Take $B = B_1 \cup B_2$

$$= \left\{ \begin{pmatrix} I & 0 & 0 \\ 0 & 0 & 0 \end{pmatrix}, \begin{pmatrix} 0 & I & 0 \\ 0 & 0 & 0 \end{pmatrix}, \begin{pmatrix} 0 & 0 & I \\ 0 & 0 & 0 \end{pmatrix}, \right.$$

$$\left. \begin{pmatrix} 0 & 0 & 0 \\ I & 0 & 0 \end{pmatrix}, \begin{pmatrix} 0 & 0 & 0 \\ 0 & I & 0 \end{pmatrix}, \begin{pmatrix} 0 & 0 & 0 \\ 0 & 0 & I \end{pmatrix} \right\}$$

$$\cup \left\{ \begin{pmatrix} I & 0 \\ 0 & 0 \\ 0 & 0 \\ 0 & 0 \end{pmatrix}, \begin{pmatrix} 0 & I \\ 0 & 0 \\ 0 & 0 \\ 0 & 0 \end{pmatrix}, \begin{pmatrix} 0 & 0 \\ I & 0 \\ 0 & 0 \\ 0 & 0 \end{pmatrix}, \begin{pmatrix} 0 & 0 \\ 0 & I \\ 0 & 0 \\ 0 & 0 \end{pmatrix}, \right.$$

$$\left. \begin{pmatrix} 0 & 0 \\ 0 & 0 \\ I & 0 \\ 0 & 0 \end{pmatrix}, \begin{pmatrix} 0 & 0 \\ 0 & 0 \\ 0 & I \\ 0 & 0 \end{pmatrix}, \begin{pmatrix} 0 & 0 \\ 0 & 0 \\ 0 & 0 \\ I & 0 \end{pmatrix}, \begin{pmatrix} 0 & 0 \\ 0 & 0 \\ 0 & 0 \\ 0 & I \end{pmatrix} \right\}$$

$\subseteq V_1 \cup V_2$, B is a bibasis of V over the bifield.

***Example 2.3.12:*** Let $V = V_1 \cup V_2 =$

$$\left\{ \begin{pmatrix} a & a \\ a & a \end{pmatrix} \middle| a \in QI \right\} \cup \left\{ \begin{pmatrix} a & a \\ a & a \\ a & a \\ a & a \end{pmatrix} \middle| a \in Z_{37}I \right\}$$

be a neutrosophic bivector space over the bifield $F = F_1 \cup F_2 = Q \cup Z_{37}$. Take $B = B_1 \cup B_2 =$

$$\left\{ \begin{pmatrix} I & I \\ I & I \end{pmatrix} \right\} \cup \left\{ \begin{pmatrix} I & I \\ I & I \\ I & I \\ I & I \end{pmatrix} \right\}$$

$\subseteq V_1 \cup V_2$, $B$ is a bibasis of $V$ over the bifield $F = Q \cup Z_{37}$.

**DEFINITION 2.3.6:** *Let $V = V_1 \cup V_2$ be a neutrosophic bivector space over the bifield $F = F_1 \cup F_2$. Let $P = P_1 \cup P_2 \subseteq V_1 \cup V_2$ be a proper bisubset of $V$ such that each $P_i$ is a linearly independent subset of $V_i$ over $F_i$; $i = 1, 2$; then we define $P = P_1 \cup P_2$ to be a bilinearly independent bisubset of $V$ over the bifield $F = F_1 \cup F_2$ or $P$ is defined to be the linearly biindependent bisubset of $V$ over the bifield $F = F_1 \cup F_2$.*

It is interesting and important to note that every bibasis is a linearly biindependent bisubset, but a linearly biindependent bisubset need not in general to be a bibasis of $V$ over the bifield $F$.
We will illustrate this situation by an example.

***Example 2.3.13:*** Let $V = V_1 \cup V_2 =$

$$\left\{ \begin{pmatrix} a & b \\ c & d \end{pmatrix} \middle| a,b,c,d \in RI \right\} \cup$$

$$\left\{ \begin{pmatrix} a_1 & a_2 & a_3 \\ b_1 & b_2 & b_3 \end{pmatrix} \middle| a_i, b_i \in Z_7 I; 1 \leq i \leq 3 \right\}$$

be a neutrosophic bivector space of type II over the bifield $F = Q \cup Z_7$. Take $B = B_1 \cup B_2 =$

$$\left\{ \begin{pmatrix} I & 0 \\ 0 & 0 \end{pmatrix}, \begin{pmatrix} 0 & I \\ 0 & 0 \end{pmatrix}, \begin{pmatrix} 0 & 0 \\ I & 0 \end{pmatrix}, \begin{pmatrix} 0 & 0 \\ 0 & I \end{pmatrix} \right\} \cup$$

$$\left\{ \begin{pmatrix} I & I & 0 \\ 0 & 0 & 0 \end{pmatrix}, \begin{pmatrix} 0 & 0 & I \\ I & 0 & 0 \end{pmatrix}, \begin{pmatrix} I & I & I \\ 0 & I & I \end{pmatrix} \right\} \subseteq V_1 \cup V_2.$$

Clearly B is a linearly biindependent bisubset of V but is not a bibasis of V. Thus in general every linearly biindependent bisubset of V need not be a bibasis of V.

**DEFINITION 2.3.7:** *Let $V = V_1 \cup V_2$ be a neutrosophic bivector space over the bifield $F = F_1 \cup F_2$, and $W = W_1 \cup W_2$ be another neutrosophic bivector space over the same bifield $F = F_1 \cup F_2$ that is $V_i$ and $W_i$ are vector spaces over the field $F_i$, $i = 1, 2$. Let $T = T_1 \cup T_2$ be a bimap from V to W; where $T_i : V_i \to W_i$ is a linear transformation from $V_i$ to $W_i$, $i = 1, 2$. We define $T = T_1 \cup T_2 : V = V_1 \cup V_2 \to W = W_1 \cup W_2$ to be a neutrosophic linear bitransformation of V to W of type II.*

We will illustrate this by a simple example.

***Example 2.3.14:*** Let $V = V_1 \cup V_2 =$

$$\left\{ \begin{pmatrix} a & b \\ c & d \end{pmatrix} \middle| a,b,c,d \in Z_7 I \right\} \cup$$

$$\left\{ \begin{pmatrix} a_1 & a_2 & a_3 & a_4 \\ a_5 & a_6 & a_7 & a_8 \end{pmatrix} \middle| a_i \in Z_{13}I; 1 \leq i \leq 8 \right\}$$

be a neutrosophic bivector space of type II over the bifield $F = F_1 \cup F_2 = Z_7 \cup Z_{13}$. Let $W = W_1 \cup W_2 = \{Z_7I \times Z_7I \times Z_7I \times Z_7I\}$

$$\cup \left\{ \begin{pmatrix} a_1 & a_2 \\ a_3 & a_4 \\ a_5 & a_6 \\ a_7 & a_8 \\ a_9 & a_{10} \end{pmatrix} \middle| a_i \in Z_{13}I; 1 \leq i \leq 10 \right\}$$

be a neutrosophic bivector space of type II over the bifield $F = F_1 \cup F_2 = Z_7 \cup Z_{13}$.

Define $T = T_1 \cup T_2: V = V_1 \cup V_2 \to W = W_1 \cup W_2$ where $T_1: V_1 \to W_1$ and $T_2: V_2 \to W_2$ is defined by

$$T_1 \begin{pmatrix} a & b \\ c & d \end{pmatrix} = (a, b, c, d)$$

and

$$T_2 \begin{pmatrix} a_1 & a_2 & a_3 & a_4 \\ a_5 & a_6 & a_7 & a_8 \end{pmatrix} = \begin{pmatrix} a_1 & a_2 \\ a_3 & a_4 \\ a_5 & a_6 \\ a_7 & a_8 \\ 0 & 0 \end{pmatrix}.$$

It is easily verified T is a linear bitransformation of V to W.

*Note:* If we take in the definition 2.3.7; $W = V$ then we call T to be a linear bioperator on V of type II. We will denote by $\text{N Hom}_{F_1 \cup F_2}(V, W)$, the collection of all neutrosophic linear bitransformations of V to W. $\text{N Hom}_{F_1 \cup F_2}(V, V)$ denotes the collection of all neutrosophic linear bioperators of V to V.

***Example 2.3.15:*** Let $V = V_1 \cup V_2 =$

$$\left\{ \begin{pmatrix} a & b \\ c & d \end{pmatrix} \middle| a,b,c,d \in QI \right\} \cup$$

$$\left\{ \begin{pmatrix} a_1 & a_2 & a_3 & a_4 & a_5 \\ a_6 & a_7 & a_8 & a_9 & a_{10} \end{pmatrix} \middle| a_i \in Z_{19}I; 1 \le i \le 10 \right\};$$

be a neutrosophic bivector space over the bifield $F = Q \cup Z_{19}$.

Define $T = T_1 \cup T_2 : V = V_1 \cup V_2 \to V = V_1 \cup V_2$ where $T_1: V_1 \to V_1$ and $T_2: V_2 \to V_2$ such that,

$$T_1 \begin{pmatrix} a & b \\ c & d \end{pmatrix} = \begin{pmatrix} a & a \\ a & a \end{pmatrix}$$

and

$$T_2 \begin{pmatrix} a_1 & a_2 & a_3 & a_4 & a_5 \\ a_6 & a_7 & a_8 & a_9 & a_{10} \end{pmatrix} = \begin{pmatrix} a & a & a & a & a \\ b & b & b & b & b \end{pmatrix}$$

where $a, b \in Z_{19}I$.

It is easily verified T is a neutrosophic bilinear operator on V of type II.

**DEFINITION 2.3.8:** *Let $V = V_1 \cup V_2$ and $W = W_1 \cup W_2$ be two neutrosopic bivector spaces over the bifield $F = F_1 \cup F_2$.*

*Let $T = T_1 \cup T_2: V_1 \cup V_2 = V \to W_1 \cup W_2 = W$ be a linear bitransformation of V to W. The bikernel of T denoted by ker T $= ker\ T_1 \cup ker\ T_2$ where $ker\ T_i = \{v^i \in V_i / T_i(v^i) = \overline{0}\}$; $i = 1, 2$. Thus ker $T = \{(v^1, v^2) \in V_1 \cup V_2 / T(v^1, v^2)\} = \{T_1(v^1) \cup T_2(v^2) = 0 \cup 0\}$.*

It is easily verified that ker T is a proper neutrosophic bisubgroup of V. Futher ker T is a neutrosophic bisubspace of V.

The reader is expected to give some examples.

**THEOREM 2.3.2:** *Let $V = V_1 \cup V_2$ and $W = W_1 \cup W_2$ be two neutrosophic bivector spaces over the bifield $F = F_1 \cup F_2$ of type II and suppose V is finite bidimensional. Let $T = T_1 \cup T_2$ be a neutrosophic bilinear transformation (linear bitransformation) of V into W. ($T_i: V_i \to W_i$; $i = 1, 2$).*

*Then*

*birank $T$ + binullity $T$ = $(n_1, n_2)$ dim V*
*= bidimension V;*

*that is (birank $T$ =) rank $T_1 \cup$ rank $T_2$ + (binullity $T$ =) nullity $T_1 \cup$ nullity $T_2$ = (bidimension $V$ = ) dim$V_1 \cup$ dim $V_2 = (n_1, n_2)$. (Here dim $V_i = n_i$; $i = 1, 2$).*

The proof is left as an exercise to the reader. Further the following theorem is also left as an exercise to the reader.

**THEOREM 2.3.3:** *Let $V = V_1 \cup V_2$ and $W = W_1 \cup W_2$ be two neutrosophic bivector spaces over the bifield $F = F_1 \cup F_2$. Let $T = T_1 \cup T_2$ and $S = S_1 \cup S_2$ be neutrosophic bilinear transformations from V into W. The bifunction*

$(T + S) = (T_1 \cup T_2 + S_1 \cup S_2)$
$= (T_1 + S_1) \cup (T_2 + S_2)$

*is defined by*

$(T + S)(\alpha) = (T_1 + S_1) \cup (T_2 + S_2)(\alpha_1 \cup \alpha_2)$
$= (T_1 + S_1)\alpha_1 \cup (T_2 + S_2)\alpha_2$
$= (T_1\alpha_1 + S_1\alpha_1) \cup (T_2\alpha_2 + S_2\alpha_2)$

*is a neutrosophic linear bitransformation from $V = V_1 \cup V_2$ to $W_1 \cup W_2$. ($\alpha = \alpha_1 \cup \alpha_2 \in V_1 \cup V_2$). If $C = C_1 \cup C_2$ is a biscalar from the bifield then the bifunction*

$(CT)\alpha = (C_1 \cup C_2)(T_1 \cup T_2)(\alpha_1 \cup \alpha_2)$
$= C_1T_1\alpha_1 \cup C_2T_2\alpha_2$

*is a bilinear transformation (linear bitransformation ) from V into W. Thus the set of all linear bitransfomations defined by biaddition and scalar bimultiplication is a neutrosophic*

*bivector space (vector bispace) over the same bifield $F = F_1 \cup F_2$.*

*Let $NL(V, W) = NL^1(V_1, W_1) \cup NL^2(V_2, W_2)$ be a neutrosophic bivector space over the bifield $F = F_1 \cup F_2$.*

*Further if $V = V_1 \cup V_2$ is a neutrosophic bivector space over the bifield $F = F_1 \cup F_2$ of finite bidimension $(n_1, n_2)$ and $W = W_1 \cup W_2$ is a neutrosophic bivector space of finite dimension $(m_1, m_2)$ over the same bifield $F = F_1 \cup F_2$. Then $NL(V, W)$ is of finite bidimension and has $(m_1 n_1, m_2 n_2)$ bidimension over the same bifield $F = F_1 \cup F_2$.*

Further we have another interesting property about these neutrosophic bivector spaces.

Let $V = V_1 \cup V_2$, $W = W_1 \cup W_2$ and $Z = Z_1 \cup Z_2$ be three neutrosophic bivector spaces over the same bifield $F = F_1 \cup F_2$. Let T be a neutrosophic bilinear transformation from V into W and S be a neutrosophic linear bitransformation from W into Z. Then the bicomposed bifunction S o T = ST defined by $ST(\alpha) = S(T(\alpha))$ is a neutrosophic bilinear transformation from V into Z. The reader is expected to prove the above claim.

Now we proceed on to define the notion of neutrosophic bilinear algebra or neutrosophic linear bialgebra of type II over the bifield $F = F_1 \cup F_2$.

**DEFINITION 2.3.9:** *Let $V = V_1 \cup V_2$ be a neutrosophic bivector space of type II over the bifield $F = F_1 \cup F_2$. If each $V_i$ is a neutrosophic linear algebra over $F_i$, $i = 1, 2$, then we call V to be a neutrosophic bilinear algebra over the bifield $F = F_1 \cup F_2$ of type II.*

We will illustrate this by some simple examples.

***Example 2.3.16:*** Let

$$V = \left\{ \begin{pmatrix} a & b \\ c & d \end{pmatrix} \middle| a, b, c, d \in QI \right\}$$

∪ {(a, b, c, d, e) | a, b, c, d, e ∈ $Z_7I$} be a neutrosophic bivector space over the bifield F = Q ∪ $Z_7$. V is clearly a neutrosophic bilinear algebra over F.

***Example 2.3.17:*** Let V = $V_1 \cup V_2$

$$= \left\{ \begin{pmatrix} a_1 & a_2 & a_3 \\ a_4 & a_5 & a_6 \\ a_7 & a_8 & a_9 \end{pmatrix} \middle| a_i \in Z_{11}I; 1 \le i \le 9 \right\}$$

{$Z_{13}I[x]$; all polynomials in the variable x with coefficients from $Z_{13}I$}; V is a neutrosophic bilinear algebra over the bifield F = $Z_{11} \cup Z_{13}$.

***Example 2.3.18:*** Let V = $V_1 \cup V_2$ =

$$\left\{ \begin{pmatrix} a_1 & a_2 \\ a_3 & a_4 \\ a_5 & a_6 \end{pmatrix} \middle| a_i \in Z_{17}I; 1 \le i \le 6 \right\} \cup \left\{ \begin{pmatrix} a & b \\ c & d \end{pmatrix} \middle| a, b, c, d \in QI \right\};$$

V is only a neutrosophic bivector space over the bifield $Z_{17}$ ∪ Q. Clearly V is not a neutrosophic bilinear algebra over the bifield of type II as $V_1$ is not a neutrosophic linear algebra over the field $Z_{17}$.

Thus we have the following interesting result, the proof of which is left as an exercise for the reader.

**THEOREM 2.3.4:** *Let V = $V_1 \cup V_2$ be a neutrosophic bilinear algebra over a bifield F = $F_1 \cup F_2$ of type II. Clearly V is a neutrosophic bivector space over the bifield F. However a neutrosophic bivector space of type II need not in general be a neutrosophic bilinear algebra of type II.*

Now we proceed on to define the new notion of neutrosophic linear bisubalgebra or neutrosophic bilinear sub algebra of type II

**DEFINITION 2.3.10:** *Let $V = V_1 \cup V_2$ be a neutrosophic bilinear algebra over a bifield $F = F_1 \cup F_2$ of type II. Take $W = W_1 \cup W_2 \subseteq V_1 \cup V_2$; W is a neutrosophic sub bilinear algebra or neutrosophic bilinear subalgebra of V if W is itself a neutrosophic linear bialgebra of type II over the bifield $F = F_1 \cup F_2$.*

We will illustrate this situation by some examples.

*Example 2.3.19:* Let $V = V_1 \cup V_2 =$

$$\left\{ \begin{pmatrix} a & b \\ c & d \end{pmatrix} \middle| a,b,c,d \in QI \right\}$$

$\cup \{(a_1\ a_2\ a_3\ a_4\ a_5\ a_6) \mid a_i \in Z_2I;\ 1 \leq i \leq 6\}$ be a neutrosophic bilinear algebra of type II over the bifield $F = Q \cup Z_2$.
   Take $W = W_1 \cup W_2 =$

$$\left\{ \begin{pmatrix} a & a \\ a & a \end{pmatrix} \middle| a \in QI \right\} \cup \{(a\ a\ a\ a\ a\ a) \mid a \in Z_2I\}$$

$\subseteq V_1 \cup V_2$, W is a neutrosophic bilinear subalgebra of V over the bifield $F = Q \cup Z_2$ of type II.

*Example 2.3.20:* Let $V = V_1 \cup V_2 =$

$$\left\{ \begin{pmatrix} a_1 & a_2 & a_3 \\ a_4 & a_5 & a_6 \\ a_7 & a_8 & a_9 \end{pmatrix} \middle| a_i \in Z_3I;\ 1 \leq i \leq 9 \right\}$$

∪ {QI[x]; all polynomials in the variable x with coefficients from QI} be a neutrosophic bilinear algebra of type II over the bifield $F = Z_3 \cup Q$.

Take $W = W_1 \cup W_2 =$

$$\left\{ \begin{pmatrix} a_1 & a_2 & a_3 \\ 0 & a_4 & a_5 \\ 0 & 0 & a_6 \end{pmatrix} \middle| a_i \in Z_3 I; 1 \le i \le 6 \right\}$$

∪ {All polynomials of even degree with coefficients from the field QI} $\subseteq V_1 \cup V_2$; W is a neutrosophic bilinear subalgebra of V of type II over the bifield $Z_3 \cup Q$.

**DEFINITION 2.3.11:** *Let $V = V_1 \cup V_2$ be a neutrosophic bilinear algebra over a bifield $F = F_1 \cup F_2$ of type II. Let $W = W_1 \cup W_2 \subseteq V_1 \cup V_2$, suppose W is only a neutrosophic bivector space of type II over the bifield $F = F_1 \cup F_2$ and is not a neutrosophic bilinear subalgebra of V of type II over the bifield F then we say W is a pseudo neutrosophic bivector subspace of V over the bifield $F = F_1 \cup F_2$ of type II.*

We will illustrate this by some examples.

*Example 2.3.21:* Let

$$V = \left\{ \begin{pmatrix} a & b \\ c & d \end{pmatrix} \middle| a,b,c,d \in QI \right\}$$

∪ {$Z_7 I[x]$; all polynomials in the variable x with coefficients from $Z_7 I$} be a neutrosophic bilinear algebra over the bifield $F = F_1 \cup F_2 = Q \cup Z_7$. Take $W = W_1 \cup W_2 =$

$$\left\{ \begin{pmatrix} 0 & a \\ b & 0 \end{pmatrix} \middle| a,b \in QI \right\} \cup \left\{ \sum_{i=0}^{20} a_i x^i \middle| a_i \in Z_7 I; 1 \le i \le 20 \right\}$$

$\subseteq V_1 \cup V_2$. Clearly W is only a bivector space over the bifield $F = Q \cup Z_7$. For product of two elements is not defined in both $W_1$ and $W_2$. Thus W is a pseudo neutrosophic bivector subspace of V over the bifield $F = Q \cup Z_7$.

*Example 2.3.22:* Let $V = V_1 \cup V_2 =$

$$\left\{ \begin{pmatrix} a & b & c \\ d & e & f \\ g & h & i \end{pmatrix} \middle| a,b,c,d,e,f,g,h,i \in Z_{11}I \right\}$$

$\cup \{(a, b, c) \mid a, b, c \in Z_{29}I\}$ be a neutrosophic bilinear algebra over the bifield $F = Z_{11} \cup Z_{29}$. Take $W = W_1 \cup W_2 =$

$$\left\{ \begin{pmatrix} a & 0 & b \\ 0 & c & 0 \\ d & 0 & 0 \end{pmatrix} \middle| a,b,c,d \in Z_{11}I \right\} \cup \{(a\ a\ a) \mid a \in Z_{29}I\}$$

$\subseteq V_1 \cup V_2$; W is a pseudo neutrosophic bivector subspace of V as $W_1$ is only a neutrosophic vector space over $Z_{11}$ which is not a neutrosophic linear algebra over $Z_{11}$, but $W_2$ is neutrosophic linear algebra over $Z_{29}$. Thus W is only a pseudo neutrosophic bivector subspace of V.

Let $V = V_1 \cup V_2$ be a bivector space over the bifield $F = F_1 \cup F_2$. A linear bitransformation $f = f_1 \cup f_2$ from $V = V_1 \cup V_2$ into the bifield $F = F_1 \cup F_2$ of biscalars is called as a linear bifunctional or bilinear functional.

However when the bivector space which are neutrosophic bivector spaces are defined over a real bifield $F = F_1 \cup F_2$ we see the notion of linear bifunctional is not possible. Hence to have the concept of linear bifunctional we need the bivector spaces to be defined over neutrosophic bifields.

However we can define neutrosophic hyper bispace of a neutrosophic bivector space.

**DEFINITION 2.3.12:** *Let $V = V_1 \cup V_2$ be a finite dimensional neutrosophic bivector space of type II over the bifield $F = F_1 \cup F_2$ of dimension $(n_1, n_2)$. Let $W = W_1 \cup W_2 \subseteq V_1 \cup V_2$ be a neutrosophic bivector subspace of V of dimension $((n_1 - 1), (n_2 - 1))$ over the bifield $F = F_1 \cup F_2$. Then we call W to be a neutrosophic hyper bispace of V.*

We will illustrate this situation by some examples.

*Example 2.3.23:* Let $V = V_1 \cup V_2 =$

$$\left\{ \begin{pmatrix} a & b \\ c & d \end{pmatrix} \middle| a, b, c, d \in QI \right\}$$

$\cup \{(a\ b\ c) \mid a, b, c \in Z_{17}I\}$ be a neutrosophic bivector space of finite bidimension over the bifield $F = Q \cup Z_{17}$. Take $W = W_1 \cup W_2 =$

$$\left\{ \begin{pmatrix} a & b \\ 0 & d \end{pmatrix} \middle| a, b, d \in QI \right\}$$

$\cup \{(a, b, 0) \mid a, b \in Z_{17}I\} \subseteq V_1 \cup V_2$; W is a neutrosophic hyper bisubspace of V over the bifield $F = Q \cup Z_{17}$.

*Example 2.3.24:* Let $V = V_1 \cup V_2 =$

$$\left\{ \sum_{i=0}^{12} a_i x^i \middle| a_i \in Z_2 I \right\} \cup \left\{ \begin{pmatrix} a_1 & a_2 \\ a_3 & a_4 \\ a_5 & a_6 \end{pmatrix} \middle| a_i \in Z_3 I \right\}$$

be a neutrosophic bivector space over the bifield $F = Z_2 \cup Z_3$. Take $W =$

$$\left\{ \sum_{i=0}^{11} a_i x^i \middle| a_i \in Z_2 I \right\} \cup$$

$$\left\{ \begin{pmatrix} a_1 & a_2 \\ a_3 & 0 \\ a_5 & a_6 \end{pmatrix} \middle| a_i \in Z_3 I; i = 1, 2, 3, 5, 6 \right\}$$

$\subseteq V_1 \cup V_2$; W is neutrosophic hyper bispace of V over the bifield $Z_2 \cup Z_3$. Clearly the bidimension of V is (13, 6) where the bidimension of W is (12, 5).

The notion of biannihilator of a biset S of a neutrosophic bivector space over a real bifield cannot be defined as the notion of linear functional is undefined for these bispaces.
    We can define neutrosophic bipolynomial ring over the bifield F. Let $F[x] = F_1[x] \cup F_2[x]$ be such that $F_i[x]$ is a polynomial ring over $F_i$ then we cannot call $F[x] = F_1[x] \cup F_2[x]$ to be a neutrosophic bipolynomial biring over $F_1 \cup F_2$ as $F_1$ and $F_2$ are not neutrosophic fields they are only real fields.

Now we can define yet another new substructure.

**DEFINITION 2.3.13:** *Let $V = V_1 \cup V_2$ be a neutrosophic bivector space over the bifield $F = F_1 \cup F_2$. Let $W = W_1 \cup W_2 \subseteq V_1 \cup V_2$ be such that W is just a bivector space over the real bifield $F = F_1 \cup F_2$; i.e., W is not a neutrosophic bivector space over the bifield F; then we call W to be a pseudo bivector subspace of V over the bifield F.*

We will illustrate this by the following examples.

***Example 2.3.25:*** Let $V = V_1 \cup V_2 =$

$$\left\{ \begin{pmatrix} a & b \\ c & d \end{pmatrix} \middle| a, b, c, d \in N(Q) \right\}$$

$\cup \{(a, b, c, d, e, f) \mid a, b, c, d, e, f \in N(Z_2)\}$ be a neutrosophic bivector space over the bifield $F = Q \cup Z_2$. Take $W = W_1 \cup W_2$

$$= \left\{ \begin{pmatrix} a & b \\ c & d \end{pmatrix} \middle| a, b, c, d \in Q \right\}$$

$\cup \{(a, b, c, d, e, f) \mid a, b, c, d, e, f \in Z_2\} \subseteq V_1 \cup V_2$. Clearly W is only a bivector space over the bifield F. Thus W is a pseudo bivector subspace of V over the bifield $F = Q \cup Z_2$.

*Example 2.3.26:* Let $V = V_1 \cup V_2 = \{Z_{17} I[x]$; all polynomials in the variable x with coefficients from $Z_{17}I\} \cup$

$$\left\{ \begin{pmatrix} a_1 & a_2 & a_3 \\ a_4 & a_5 & a_6 \end{pmatrix} \middle| a_i \in N(Z_{13}) \right\}$$

be a neutrosophic bivector space over the bifield $F = Z_{17} \cup Z_{13}$. We see there does not exist a $W = W_1 \cup W_2 \subseteq V_1 \cup V_2$ such that W is a bivector space over the bifield $F = Z_{17} \cup Z_{13}$.

Thus we see from this example that all neutrosophic bivector spaces need not in general contain pseudo bivector subspaces.

In view of this we have the following result which proves the existence of a class of neutrosophic bivector spaces which do not contain pseudo bivector subspaces.

**THEOREM 2.3.5:** *Let $V = V_1 \cup V_2$ be a neutrosophic bivector space over the real bifield $F = F_1 \cup F_2$. Even if one of $V_1$ (or $V_2$) has its entries from the neutrosophic field $F_1I$ (or $F_2I$) then we see V has no pseudo bivector subspaces.*

*Proof:* We see in $V = V_1 \cup V_2$ the entries are in one of $V_1$ or $V_2$ or in both $V_1$ and $V_2$, the entries are taken from $F_1I$ ($F_2I$) or from $F_1I$ and $F_2I$. Since $F_i \not\subseteq F_iI$; $i = 1, 2$ we see $V_i$ can never be a vector space over $F_i$ but only a neutrosophic vector space over $F_i$, $i = 1, 2$. Hence the claim.

We say a neutrosophic bivector space $V = V_1 \cup V_2$ is a pseudo simple neutrosophic bivector space if V has no proper pseudo bivector subspace.

*Example 2.3.27:* Let $V = V_1 \cup V_2 =$

$$\left\{ \begin{pmatrix} a & b \\ c & d \end{pmatrix} \middle| a,b,c,d \in Z_7 I \right\} \cup \left\{ \begin{pmatrix} a_1 \\ a_2 \\ a_3 \\ a_4 \end{pmatrix} \middle| a_i \in QI; 1 \leq i \leq 4 \right\}$$

be a neutrosophic bivector space over the bifield $F = Z_7 \cup Q$. V is a pseudo simple neutrosophic bivector space.

*Example 2.3.28:* Let $V = V_1 \cup V_2 = \{Z_{11}I[x];$ all polynomials in the variable x with coefficients from the field $Z_{11}I\} \cup \{(a_1, a_2, a_3, a_4, a_5, a_6, a_7) \mid a_i \in N(Q); 1 \leq i \leq 7\}$ be a neutrosophic bivector space over the bifield $F = Z_{11} \cup Q$. V is a pseudo simple neutrosophic bivector space over the bifield F.

Now we proceed onto define the notion of quasi pseudo bivector subspace of a neutrosophic bivector space.

**DEFINITION 2.3.14:** *Let $V = V_1 \cup V_2$ be a neutrosophic bivector space over the bifield $F = F_1 \cup F_2$. Let $W = W_1 \cup W_2 \subseteq V_1 \cup V_2$ where only one of $W_1$ or $W_2$ is a vector space over $F_1$ or $F_2$ and the other is a neutrosophic vector space over $F_1$ or $F_2$ then we call W to be a quasi pseudo bivector subspace of V over the bifield $F = F_1 \cup F_2$.*

We will illustrate this situation by some examples.

*Example 2.3.29:* Let $V = V_1 \cup V_2 = \{(a, b, c, d) \mid a, b, c, d \in Z_{13}I\} \cup$

$$\left\{ \begin{pmatrix} a_1 & a_2 & a_3 \\ a_4 & a_5 & a_6 \\ a_7 & a_8 & a_9 \end{pmatrix} \middle| a_i \in N(Z_5); 1 \le i \le 9 \right\}$$

be a neutrosophic bivector space over the bifield $F = Z_{13} \cup Z_5$. Take $W = W_1 \cup W_2 = \{(a, a, a, a) \mid a \in Z_{13}I\} \cup$

$$\left\{ \begin{pmatrix} a_1 & a_2 & a_3 \\ a_4 & a_5 & a_6 \\ a_7 & a_8 & a_9 \end{pmatrix} \middle| a_i \in Z_5; 1 \le i \le 9 \right\}$$

$\subseteq V_1 \cup V_2$; $W_1$ is a neutrosophic vector subspace of $V_1$ over $Z_{13}$ and $W_2$ is just vector space of $V_2$ over $Z_5$. We see $W_2$ is not a neutrosophic vector subspace of $V_2$ over $Z_5$. Thus $W = W_1 \cup W_2$ is a quasi pseudo bivector subspace of V over the bifield $F = Z_{13} \cup Z_5$.

*Example 2.3.30:* Let $V = V_1 \cup V_2 = \{N(Z_{19})[x]$; all polynomials in the variable x with coefficients from the field $N(Z_{19})\} \cup$

$$\left\{ \begin{pmatrix} a & b & c \\ d & e & f \end{pmatrix} \middle| a, b, c, d, e, f \in Z_{43}I \right\}$$

be a neutrosophic bivector space over the bifield $F = Z_{19} \cup Z_{43}$. Take $W = W_1 \cup W_2 = \{Z_{19}[x]$; the set of all polynomials in the variable x with coefficients from $Z_{19}\} \cup$

$$\left\{ \begin{pmatrix} a & a & a \\ a & a & a \end{pmatrix} \middle| a \in Z_{43}I \right\}$$

$\subseteq V_1 \cup V_2$; W is a quasi pseudo bivector subspace of V over the bifield $F = Z_{19} \cup Z_{43}$.

Now it may so happen we can have for some neutrosophic bivector subspace both pseudo bivector subspace as well as quasi pseudo bivector subspaces.

We will illustrate this situation by an example.

***Example 2.3.31:*** Let $V = V_1 \cup V_2 =$

$$\left\{ \begin{pmatrix} a & b \\ c & d \end{pmatrix} \middle| a,b,c,d \in N(Q) \right\}$$

$\cup$ {$N(Z_{47})[x]$; all polynomials in the variable x with coefficients from the neutrosophic field $N(Z_{47})$} be a neutrosophic bivector space over the bifield $F = Q \cup Z_{47}$. Take $W = W_1 \cup W_2 =$

$$\left\{ \begin{pmatrix} a & a \\ a & a \end{pmatrix} \middle| a \in Q \right\}$$

$\cup$ [$Z_{47}[x]$; all polynomials in the variable x with coefficients from the field $Z_{47}$} $\subseteq V_1 \cup V_2$; clearly W is a pseudo bivector subspace of V over the bifield $F = Q \cup Z_{47}$.

Let $S = S_1 \cup S_2 =$

$$\left\{ \begin{pmatrix} a & a \\ a & a \end{pmatrix} \middle| a \in QI \right\}$$

$\cup$ {$Z_{47}I[x]$; all polynomials in the variable x with coefficients from $Z_{47}I$} $\subseteq V_1 \cup V_2$. S is a quasi pseudo bivector subspace of V. Thus V can have both types of bivector subspaces.

Finally we define subneutrosophic bivector subspace.

**DEFINITION 2.3.15:** *Let $V = V_1 \cup V_2$ be a neutrosophic bivector space over the bifield $F = F_1 \cup F_2$. Let $W = W_1 \cup W_2 \subseteq V_1 \cup V_2$ be a neutrosophic bivector space over the bisubfield $K = K_1 \cup K_2 \subseteq F_1 \cup F_2$; $K_i \subseteq F_i$; $K_i$ is a proper subfield of $F_i$; i = 1, 2. We* *then call W to be a subneutrosophic bivector subspace of V over the subbifield K of the bifield F. If $V = V_1 \cup V_2$ has no subneutrosophic bivector subspace then we call V to be a sub bisimple neutrosophic bivector space.*

We will illustrate this situation by some examples.

***Example 2.3.32:*** Let $V = V_1 \cup V_2 =$

$$\left\{ \begin{pmatrix} a & b \\ c & d \end{pmatrix} \middle| a,b,c,d \in RI \right\} \cup \{(a, b, c, d, e) \mid a, b, c, d, e \in RI\}$$

be a neutrosophic bivector space over the bifield,

$$F = Q(\sqrt{2},\sqrt{3},\sqrt{7},\sqrt{11},\sqrt{17}) \cup Q(\sqrt{19},\sqrt{23},\sqrt{43},\sqrt{41},\sqrt{7}).$$

Take $W = W_1 \cup W_2 =$

$$\left\{ \begin{pmatrix} a & a \\ a & a \end{pmatrix} \middle| a \in RI \right\} \cup \{(a, a, a, a, a) \mid a \in RI\} \subseteq V_1 \cup V_2,$$

$W = W_1 \cup W_2$ is a neutrosophic bivector space over the subbifield
$$\begin{aligned} K &= Q(\sqrt{2},\sqrt{11},\sqrt{17}) \cup Q(\sqrt{19},\sqrt{41}) \\ &= K_1 \cup K_2 \subseteq F_1 \cup F_2. \end{aligned}$$

Thus W is a subneutrosophic bivector subspace of V over the subbifield $K = K_1 \cup K_2$.

***Example 2.3.33:*** Let $V = V_1 \cup V_2 =$

$$\left\{ \begin{pmatrix} a_1 & a_2 & a_3 & a_4 & a_5 \\ a_6 & a_7 & a_8 & a_9 & a_{10} \end{pmatrix} \middle| a_i \in RI \right\} \cup$$

$$\left\{ \begin{pmatrix} a_1 & a_2 \\ a_3 & a_4 \\ a_5 & a_6 \end{pmatrix} \middle| a_i \in \left[ \frac{Z_2[x]}{\langle x^2 + x + 1 \rangle} \right] I \right\}$$

be a neutrosophic bivector space over the bifield $R \cup \frac{Z_2[x]}{\langle x^2 + x + 1 \rangle}$. Take $W = W_1 \cup W_2 =$

$$\left\{ \begin{pmatrix} a & a & a & a & a \\ a & a & a & a & a \end{pmatrix} \middle| a_i \in QI \right\} \cup \left\{ \begin{pmatrix} a & a \\ a & a \\ a & a \end{pmatrix} \middle| a_i \in Z_2 I \right\}$$

$\subseteq V_1 \cup V_2$; W is a subneutrosophic bivector subspace of V over the subbifield $K = K_1 \cup K_2 = Q \cup Z_2 \subseteq R \cup \frac{Z_2[x]}{\langle x^2 + x + 1 \rangle}$.

Now we proceed onto define the notion of strong neutrosophic bivector space and discuss a few important properties about them.

**DEFINITION 2.3.16:** *Let $V = V_1 \cup V_2$ be a neutrosophic bivector space over the neutrosophic bifield $F = F_1 \cup F_2$, then we call V to be a strong neutrosophic bivector space of type II.*

We will illustrate this by some examples.

***Examples 2.3.34:*** Let $V = V_1 \cup V_2 =$

$$\left\{ \begin{pmatrix} a & b \\ c & d \end{pmatrix} \middle| a, b, c, d \in QI \right\}$$

$\cup$ {(a, b, c, d, e, f, g) | a, b, c, d, e, f, g $\in Z_7$ I} be a strong neutrosophic bivector space over the neutrosophic bifield $F = F_1 \cup F_2 = QI \cup Z_7 I$.

***Example 2.3.35:*** Let $V = V_1 \cup V_2 =$

$$\left\{ \begin{pmatrix} a_1 & a_2 & a_3 & a_4 \\ a_5 & a_6 & a_7 & a_8 \end{pmatrix} \middle| a_i \in N(Q); 1 \le i \le 8 \right\} \cup$$

$$\left\{ \begin{pmatrix} a_1 & a_2 & a_3 \\ a_4 & a_5 & a_6 \\ a_7 & a_8 & a_9 \end{pmatrix} \middle| a_i \in N(Z_{11}); 1 \le i \le 9 \right\}$$

be a strong neutrosophic bivector space over the neutrosophic bifield $F = QI \cup Z_{11}I$.

We see strong neutrosophic bivector spaces are defined over neutrosophic bifields but neutrosophic bivector spaces are defined over real bifields. We see only incase of strong neutrosophic bispaces we can define neutrosophic bifunctionals but incase of neutrosophic bivector spaces we cannot define neutrosophic bifunctionals.

Now we will proceed onto define substructures in strong neutrosophic bivector spaces.

**DEFINITION 2.3.17:** *Let $V = V_1 \cup V_2$ be strong a neutrosophic bivector space over the neutrosophic bifield $F = F_1 \cup F_2$. Let $W = W_1 \cup W_2 \subseteq V_1 \cup V_2$, if W is a strong neutrosophic bivector space over the neutrosophic bifield $F = F_1 \cup F_2$, then we call W to be a strong neutrosophic bivector subspace of V over the neutrosophic bifield $F = F_1 \cup F_2$.*

We will illustrate this by the following examples.

***Example 2.3.36:*** Let $V = V_1 \cup V_2 =$

$$\left\{ \begin{pmatrix} a & b & c \\ d & e & f \end{pmatrix} \middle| a,b,c,d,e,f \in N(Z_7) \right\} \cup$$

$$\left\{ \begin{pmatrix} a_1 \\ a_2 \\ a_3 \\ a_4 \\ a_5 \end{pmatrix} \middle| a_i \in N(Z_{11}); 1 \le i \le 5 \right\}$$

be a strong neutrosophic bivector space over the neutrosophic bifield $Z_7 I \cup Z_{11} I$. Take $W = W_1 \cup W_2 =$

$$\left\{ \begin{pmatrix} a & a & a \\ a & a & a \end{pmatrix} \middle| a \in Z_7 I \right\} \cup \left\{ \begin{pmatrix} a \\ a \\ a \\ a \\ a \\ a \end{pmatrix} \middle| a \in Z_{11} I \right\}$$

$\subseteq V_1 \cup V_2$; W is a strong neutrosophic bivector subspace of V over the neutrosophic bifield $Z_7 I \cup Z_{11} I$.

*Example 2.3.37:* Let $V = V_1 \cup V_2 =$

$$\left\{ \begin{pmatrix} a_1 & a_2 & a_3 & a_4 \\ 0 & a_5 & a_6 & a_7 \\ 0 & 0 & a_8 & a_9 \\ 0 & 0 & 0 & a_{10} \end{pmatrix} \middle| a_i \in Z_{23} I; 1 \le i \le 10 \right\} \cup$$

$$\left\{ \begin{pmatrix} a_1 & 0 & a_2 & 0 & a_3 & 0 & a_4 \\ 0 & a_5 & 0 & a_6 & 0 & a_7 & 0 \end{pmatrix} \middle| a_i \in Z_{17} I; 1 \le i \le 7 \right\}$$

be a strong neutrosophic bivector space over the neutrosophic bifield $F = F_1 \cup F_2 = Z_{23} I \cup Z_{17} I$.

Let $W = W_1 \cup W_2 =$

$$\left\{ \begin{pmatrix} a & a & a & a \\ 0 & a & a & a \\ 0 & 0 & a & a \\ 0 & 0 & 0 & a \end{pmatrix} \middle| a \in Z_{23}I \right\} \cup$$

$$\left\{ \begin{pmatrix} a & 0 & a & 0 & a & 0 & a \\ 0 & 0 & 0 & 0 & 0 & 0 & 0 \end{pmatrix} \middle| a \in Z_{17}I \right\}$$

$\subseteq V_1 \cup V_2$; W is a strong neutrosophic bivector subspace of V over the neutrosophic bifield $F = F_1 \cup F_2$.

**DEFINITION 2.3.18:** *Let $V = V_1 \cup V_2$ be a strong neutrosophic bivector space over the neutrosophic bifield $F = F_1 \cup F_2$. Let $W = W_1 \cup W_2 \subseteq V_1 \cup V_2$; is defined to be a pseudo strong neutrosophic bivector subspace of V if W is a neutrosophic bivector space over the real bifield $K = K_1 \cup K_2 \subseteq F_1 \cup F_2$.*

We will illustrate this by some examples.

*Example 2.3.38:* Let $V = V_1 \cup V_2 =$

$$\left\{ \begin{pmatrix} a & b \\ c & d \end{pmatrix} \middle| a,b,c,d \in N(Z_{11}) \right\} \cup$$

$\{(a_1, a_2, a_3, a_4, a_5, a_6, a_7) \mid a_i \in N(Z_{11}); 1 \leq i \leq 7\}$ be a strong neutrosophic bivector space over the neutrosophic bifield $F = N(Z_{11}) \cup N(Z_{17}) = F_1 \cup F_2$. Choose $W = W_1 \cup W_2 =$

$$\left\{ \begin{pmatrix} a & b \\ c & c \end{pmatrix} \middle| a,b,c \in Z_{11}I \right\} \cup$$

$\{(a_1\ 0\ a_3\ 0\ a_5\ 0\ a_7) \mid a_1, a_3, a_5, a_7 \in N(Z_{11})\} \subseteq V_1 \cup V_2$. W is a pseudo strong neutrosophic bivector subspace of V over the real bifield $Z_{11} \cup Z_{17} \subseteq F_1 \cup F_2$.

***Example 2.3.39:*** Let $V = V_1 \cup V_2 = \{N(Z_{19})[x]$; all polynomials in the variable x with coefficients from $N(Z_{19})\} \cup \{(x_1, x_2, x_3, x_4, x_5) \mid x_i \in N(Z_{23}); 1 \leq i \leq 5\}$ be strong neutrosophic space over the neutrosophic bifield $F = N(Z_{19}) \cup Z_{23}I$.

Take $W = W_1 \cup W_2 = \{Z_{19}I[x]$; all polynomials in the variable x with coefficients from $Z_{19}I\} \cup \{(a\ a\ a\ a\ b) \mid a, b \in Z_{23}I\} \subseteq V_1 \cup V_2$; W is a pseudo strong neutrosophic bivector subspace of V over the real bifield $K = K_1 \cup K_2 = Z_{19} \cup Z_{23} \subseteq F_1 \cup F_2 = N(Z_{19}) \cup Z_{23}I$.

Recall a bifield $F = F_1 \cup F_2$ is said to be a quasi neutrosophic bifield if one of $F_1$ or $F_2$ is a neutrosophic field and the other is just a real field. $F = QI \cup Z_{17}$ is a quasi neutrosophic bifield. $F = Q \cup Z_{11}I$ is a quasi neutrosophic bifield. $F = N(Z_2) \cup Z_3$ is a quasi neutrosophic bifield.

**DEFINITION 2.3.19:** *Let $V = V_1 \cup V_2$ be a neutrosophic bivector space over the bifield $F = F_1 \cup F_2$. If $F = F_1 \cup F_2$ is only a quasi neutrosophic bifield then we call V to be a quasi strong neutrosophic bivector space over the quasi neutrosophic bifield.*

We will illustrate this situation by some simple examples.

***Example 2.3.40:*** Let $V = V_1 \cup V_2 =$

$$\left\{ \begin{pmatrix} a & b \\ c & d \end{pmatrix} \middle| a,b,c,d \in QI \right\} \cup$$

$$\left\{ \begin{pmatrix} a & 0 & 0 & 0 \\ b & c & 0 & 0 \\ 0 & d & 0 & 0 \\ 0 & 0 & e & f \end{pmatrix} \middle| a,b,c,d,e,f \in Z_{17}I \right\}$$

be a quasi strong neutrosophic bivector space over the quasi neutrosophic bifield $F = Q \cup Z_{17}I$.

*Example 2.3.41:* Let $V = V_1 \cup V_2 =$

$$\left\{ \begin{pmatrix} a_1 & a_2 & a_3 & a_4 \\ a_5 & a_6 & a_7 & a_8 \\ a_9 & a_{10} & a_{11} & a_{12} \end{pmatrix} \middle| a_i \in N(Z_{23}); 1 \le i \le 12 \right\} \cup$$

$\{Z_{11}I[x]$; all polynomials in the variable x with coefficients from the field $Z_{11}I\}$ be a quasi strong neutrosophic bivector space over the quasi neutrosophic bifield $F = Z_{23} \cup Z_{11}I$.

**DEFINITION 2.3.20:** *Let $V = V_1 \cup V_2$ be a neutrosophic bivector space over the bifield $F = F_1 \cup F_2$. Let $W = W_1 \cup W_2 \subseteq V_1 \cup V_2$ be a strong neutrosophic bivector space over the neutrosophic subbifield $K = K_1 \cup K_2 \subseteq F_1 \cup F_2$; then we call W to be a strong neutrosophic bivector subspace of V over the neutrosophic bisubfield $K = K_1 \cup K_2 \subseteq F_1 \cup F_2$.*

We will illustrate this situation by some examples.

*Example 2.3.42:* Let $V = V_1 \cup V_2 =$

$$\left\{ \begin{pmatrix} a_1 & a_2 & a_3 & a_4 & a_5 \\ 0 & a_6 & 0 & a_7 & 0 \\ a_8 & 0 & 0 & 0 & a_9 \\ 0 & a_{10} & a_{11} & a_{12} & 0 \\ a_{13} & 0 & a_{14} & 0 & a_{15} \end{pmatrix} \middle| a_i \in N(Q); 1 \le i \le 15 \right\} \cup$$

$$\left\{ \begin{pmatrix} a_1 & a_2 \\ 0 & a_3 \\ a_4 & 0 \\ 0 & a_5 \\ a_6 & 0 \\ 0 & a_7 \\ a_8 & a_9 \end{pmatrix} \middle| a_i \in N(Z_{11}); 1 \le i \le 9 \right\}$$

be a strong neutrosophic bivector space over the neutrosophic bifield $F = F_1 \cup F_2 = N(Q) \cup N(Z_{11})$. Take $W = W_1 \cup W_2 =$

$$\left\{ \begin{pmatrix} a & a & a & a & a \\ 0 & a & 0 & a & 0 \\ a & 0 & 0 & 0 & a \\ 0 & a & a & a & 0 \\ a & 0 & a & 0 & a \end{pmatrix} \middle| a \in N(Q) \right\} \cup \left\{ \begin{pmatrix} a & a \\ 0 & a \\ a & 0 \\ 0 & a \\ a & 0 \\ 0 & a \\ a & a \end{pmatrix} \middle| a \in Z_{11}I \right\}$$

$\subseteq V_1 \cup V_2$; W is a strong subneutrosophic bivector subspace of V over the neutrosophic bisubfield $K = K_1 \cup K_2 = QI \cup Z_{11}I \subseteq N(Q) \cup N(Z_{11})$.

*Example 2.3.43:* Let $V = V_1 \cup V_2 =$

$$\left\{ \begin{pmatrix} a_1 & a_2 & a_3 \\ a_4 & a_5 & a_6 \\ a_7 & a_8 & a_9 \\ a_{10} & a_{11} & a_{12} \\ a_{13} & a_{14} & a_{15} \end{pmatrix} \middle| a_i \in N(Z_{47}); 1 \leq i \leq 15 \right\} \cup$$

$$\left\{ \begin{pmatrix} a_1 & a_2 & a_3 & a_4 & a_5 & a_6 & a_7 \\ a_8 & a_9 & a_{10} & a_{11} & a_{12} & a_{13} & a_{14} \\ a_{15} & a_{16} & a_{17} & a_{18} & a_{19} & a_{20} & a_{21} \end{pmatrix} \middle| a_i \in N(Z_3); 1 \leq i \leq 21 \right\}$$

be a strong neutrosophic bivector space over the neutrosophic bifield $F = F_1 \cup F_2 = N(Z_{47}) \cup N(Z_3)$.
   Take $W = W_1 \cup W_2 =$

$$\left\{\begin{pmatrix} a & a & a \\ a & a & a \\ a & a & a \\ a & a & a \\ a & a & a \end{pmatrix} \middle| a \in Z_{47}I\right\} \cup \left\{\begin{pmatrix} a & a & a & a & a & a & a \\ 0 & 0 & 0 & 0 & 0 & 0 & 0 \\ a & 0 & a & 0 & a & 0 & 0 \end{pmatrix} \middle| a \in Z_3I\right\}$$

$\subseteq V_1 \cup V_2$; W is a strong subneutrosophic bivector subspace of V over the neutrosophic bisubfield $K = K_1 \cup K_2 = Z_{47}I \cup Z_3I \subseteq N(Z_{47}) \cup N(Z_3) = F_1 \cup F_2$.

Now we state a result which will prove the existence of strong subneutrosophic bivector subspaces.

**THEOREM 2.3.6:** *Let $V = V_1 \cup V_2$ be a strong neutrosophic bivector space over a neutrosophic field $F = F_1 \cup F_2$ where both $F_1$ and $F_2$ are of the form $F_i = N(K_i)$ where $K_i$ is a real field; $i = 1, 2$ then V has a strong subneutrosophic bivector subspace provided V has neutrosophic bivector subspaces.*

The proof of this theorem is left as an exercise for the reader.

**DEFINITION 2.3.21:** *Let $V = V_1 \cup V_2$ be a strong neutrosophic bivector space over the neutrosophic bifield $F = F_1 \cup F_2$. If V has no strong sub neutrosophic bivector subspaces then we call V to be a bisimple strong neutrosophic bivector space.*

We will illustrate this by some simple examples.

*Example 2.3.44:* Let $V = V_1 \cup V_2 =$

$$\left\{\begin{pmatrix} x_1 & x_2 & x_3 & x_4 & x_5 & x_6 \\ x_7 & x_8 & x_9 & x_{10} & x_{11} & x_{12} \end{pmatrix} \middle| x_i \in Z_{13}I; 1 \leq i \leq 12\right\} \cup$$

$$\left\{ \begin{pmatrix} a & b & a \\ b & a & b \\ a & a & a \\ b & b & b \\ a & a & b \\ b & b & a \\ a & b & b \\ b & a & a \end{pmatrix} \middle| a, b \in Z_5 I \right\}$$

be a strong space over a neutrosophic bivector space over the neutrosophic bifield $F = F_1 \cup F_2 = Z_{13}I \cup Z_5I$. We see there exists no strong neutrosophic bivector subspace for V. This is true as $F = F_1 \cup F_2 = Z_{13}I \cup Z_5I$ has no neutrosophic subbifield. Hence the claim that V is a bisimple strong subneutronsophic bivector space.

*Example 2.3.45:* Let $V = V_1 \cup V_2 = \{N(Z_{19})[x]$; all polynomial in the variable x with coefficients from $N(Z_{19})\} \cup$

$$\left\{ \begin{pmatrix} a_1 \\ a_2 \\ a_3 \\ a_4 \\ a_5 \\ a_6 \end{pmatrix} \middle| a_i \in N(Z_{23}); 1 \le i \le 6 \right\}$$

be a strong neutrosophic bivector space over the neutrosophic bifield $Z_{19}I \cup Z_{23}I$. Take any $W = W_1 \cup W_2 \subseteq V_1 \cup V_2$; we see as F has no subbifield which is neutrosophic, V has no strong subneutrosophic bivector spaces; so V is a bisimple strong subneutrosophic bivector space.

Now we give a theorem which guarantees the existence of bisimple strong subneutrosophic bivector spaces.

**THEOREM 2.3.7:** *Let $V = V_1 \cup V_2$ be a strong neutrosophic bivector space over the neutrosophic bifield $F = F_1 \cup F_2$, where $F_1$ and $F_2$ are of the form KI where K is the prime real field of characteristic zero or a prime p. Then V is a bisimple strong subneutrosophic bivector space over the neutrosophic bifield $F = F_1 \cup F_2$.*

*Proof:* Given $V = V_1 \cup V_2$ is a strong neutrosophic bivector space over the neutrosophic bifield $F = F_1 \cup F_2$ and $F_i = K_i I$ where $K_i$ is a prime field, $i = 1, 2$. So $F_i$ has no proper neutrosophic subfield for $i = 1, 2$. Hence V cannot have a strong subneutrosophic bivector space over any subfield of the bifield F. Hence $V = V_1 \cup V_2$ is a bisimple strong subneutrosophic bivector space over F.

Thus we have proved the existence of bisimple strong subneutrosophic bivector spaces.

Now we proceed on to define the concept of linearly independent bisubset and the basis for the strong neutrosophic bivector spaces.

**DEFINITION 2.3.22:** *Let $V = V_1 \cup V_2$ be a strong neutrosophic bivector space defined over the neutrosophic bifield $F = F_1 \cup F_2$. A bisubset $S = S_1 \cup S_2 \subseteq V_1 \cup V_2$ is said to be a linearly biindependent or bilinearly independent over F if each $S_i$ is a linearly independent subset of $V_i$ over $F_i$; $i = 1, 2$. If $S = S_1 \cup S_2$ be a linearly biindependent bisubset of V and if each $S_i$ generates $V_i$ over $F_i$ for $i = 1, 2$ then we say S is a bibasis of $V = V_1 \cup V_2$ over $F = F_1 \cup F_2$.*

We will illustrate this situation by some examples.

***Example 2.3.46:*** Let $V = V_1 \cup V_2 =$

$$\left\{ \begin{pmatrix} a & a \\ a & a \\ a & a \end{pmatrix} \middle| a \in N(Z_{11}) \right\} \cup \{(a\ a\ a) \mid a \in N(Z_{17})\}$$

be a strong neutrosophic bivector space over the neutrosophic bifield $F = N(Z_{11}) \cup N(Z_{17})$. Take $S = S_1 \cup S_2 =$

$$\left\{ \begin{pmatrix} 1 & 1 \\ 1 & 1 \\ 1 & 1 \end{pmatrix} \right\} \cup \{(1\ 1\ 1)\} \subseteq V_1 \cup V_2.$$

S is a bibasis of V we say the bidimension of V is the (number of elements in $S_1$) $\cup$ (number of elements in $S_2$) where $S = S_1 \cup S_2$ is a bibasis of V over the neutrosophic bifield $F = F_1 \cup F_2$ and it is denoted by $(|S_1|, |S_2|)$ or $|S_1| \cup |S_2|$. We see the bidimension of $V = V_1 \cup V_2$ in example 2.3.46 is (1, 1).

***Example 2.3.47:*** Let $V = V_1 \cup V_2 = \{Z_{17}I[x]$; all polynomials in the variable x with coefficients from the neutrosophic field $Z_{17}I\} \cup \{(N(Q) \times N(Q) \times N(Q))\}$ be a strong neutrosophic bivector space over the neutrosophic bifield $F = Z_{17}I \cup N(Q)$.

Take $S = S_1 \cup S_2 = \{I, Ix, Ix^2, \ldots, Ix^n, \ldots\} \cup \{(100), (010), (001)\} \subseteq V_1 \cup V_2$; S is a bibasis of V over the bifield $F = Z_{17}I \cup N(Q)$ and bidimension of V over F is $(\infty, 3)$.

We say the bidimension is bifinite if both $|S_1|$ and $|S_2|$ are finite; even if one of $|S_1|$ or $|S_2|$ is not finite we say the bidimension of V is biinfinite over F. We see the bidimension of V given in example 2.3.47 is biinfinite.

Next we will prove that in general every linearly biindependent bisubset of a strong neutrosophic bisubset of a strong neutrosophic bivector space need not form a bibasis of $V = V_1 \cup V_2$ over $F = F_1 \cup F_2$.

We will illustrate this by some examples.

***Examples 2.3.48:*** Let $V = V_1 \cup V_2 =$

$$\left\{ \begin{pmatrix} a_1 & a_2 \\ a_3 & a_4 \\ a_5 & a_6 \\ a_7 & a_8 \end{pmatrix} \middle| a_i \in Z_{19}I; 1 \le i \le 8 \right\} \cup$$

$\{(a_1, a_2, a_3, a_4, a_5) \mid a_i \in N(Z_{11}); 1 \le i \le 5\}$ be a strong neutrosophic bivector space over the bifield $F = F_1 \cup F_2 = Z_{19}I \cup N(Z_{11})$. Take $S = S_1 \cup S_2 =$

$$\left\{ \begin{pmatrix} I & I \\ 0 & 0 \\ 0 & 0 \\ 0 & 0 \end{pmatrix}, \begin{pmatrix} 0 & 0 \\ I & 0 \\ 0 & I \\ 0 & 0 \end{pmatrix}, \begin{pmatrix} 0 & 0 \\ 0 & I \\ 0 & 0 \\ I & 0 \end{pmatrix}, \begin{pmatrix} 0 & 0 \\ 0 & 0 \\ I & 0 \\ 0 & I \end{pmatrix} \right\} \cup$$

$\{(I, 0, 0, 0, 0), (0, I, I, 0, 0), (0, 0, I, 0, I)\} \subseteq V_1 \cup V_2$; $S$ is a linearly biindependent bisubset of $V$ over the bifield $F = Z_{19}I \cup N(Z_{11})$. Clearly $S$ is not a bibasis of $V = V_1 \cup V_2$ over $F = Z_{19}I \cup N(Z_{11})$.

*Example 2.3.49:* Let $V = V_1 \cup V_2 =$

$$\left\{ \begin{pmatrix} a_1 & a_2 & a_3 & a_4 \\ a_5 & a_6 & a_7 & a_8 \end{pmatrix} \middle| a_i \in QI; 1 \le i \le 8 \right\} \cup$$

$$\left\{ \begin{pmatrix} a & b \\ b & a \\ a & a \\ b & b \\ c & c \\ c & a \end{pmatrix} \middle| a, b, c \in Z_7 I \right\}$$

be a strong neutrosophic bivector space over the bifield $F = F_1 \cup F_2 = QI \cup Z_7I$. Take $S = S_1 \cup S_2 =$

$$\left\{ \begin{pmatrix} I & I & 0 & 0 \\ 0 & 0 & 0 & 0 \end{pmatrix}, \begin{pmatrix} 0 & 0 & 3I & I \\ 0 & 0 & 0 & 0 \end{pmatrix}, \begin{pmatrix} 0 & 0 & 0 & 0 \\ -3I & I & 0 & 0 \end{pmatrix}, \right.$$

$$\left. \begin{pmatrix} 0 & 0 & 0 & 0 \\ 0 & 0 & 7I & 0 \end{pmatrix}, \begin{pmatrix} 0 & 0 & 0 & 0 \\ 0 & 0 & 0 & 2I \end{pmatrix} \right\} \cup$$

$$\left\{ \begin{bmatrix} I & I \\ 0 & 0 \\ 0 & 0 \\ 0 & 0 \\ 0 & 0 \\ I & I \end{bmatrix}, \begin{bmatrix} 0 & 0 \\ 0 & I \\ 0 & I \\ 0 & 0 \\ I & I \\ 0 & 0 \end{bmatrix}, \begin{bmatrix} 0 & 0 \\ 0 & 0 \\ I & 0 \\ I & 0 \\ 0 & 0 \\ I & 0 \end{bmatrix} \right\}$$

$\subseteq V_1 \cup V_2$; S is a linearly biindependent bisubset of $V = V_1 \cup V_2$ over $F = F_1 \cup F_2 = QI \cup Z_7I$.

Now we will proceed on to define the notion of strong neutrosophic bilinear algebra or strong neutrosophic linear bialgebra.

**DEFINITION 2.3.23:** *Let $V = V_1 \cup V_2$ be a strong neutrosophic bivector space over the neutrosophic bifield $F = F_1 \cup F_2$. If each $V_i$ is a neutrosophic linear algebra over the field $F_i$, $i = 1, 2$ then we call $V = V_1 \cup V_2$ to be a strong neutrosophic bilinear algebra over the neutrosophic bifield $F = F_1 \cup F_2$.*

We will illustrate this by some simple examples.

*Example 2.3.50:* Let $V = V_1 \cup V_2 =$

$$\left\{ \begin{pmatrix} a & b \\ c & d \end{pmatrix} \middle| a, b, c, d \in QI \right\} \cup$$

{(a, b, c, d, e, f, g, h, i) | a, b, c, d, e, f, g, h, i $\in Z_{11}I$} be a strong neutrosophic bilinear algebra over the bifield $F = QI \cup Z_{11}I$.

*Example 2.3.51:* Let $V = V_1 \cup V_2 = \{Z_{13}I[x];$ all polynomials in the variable x with coefficients from $Z_{13}I\} \cup$

$$\left\{ \begin{pmatrix} a & a \\ a & a \end{pmatrix} \middle| a \in Z_{23}I \right\}$$

be a strong neutrosophic bilinear algebra over the neutrosophic bifield $F = F_1 \cup F_2 = Z_{13}I \cup Z_{23}I$.

We see in general all strong neutrosophic bivector spaces are not strong neutrosophic bilinear algebras. But all strong neutrosophic bilinear algebras are strong neutrosophic bivector spaces.

We will illustrate the former one by an example as the latter claim simply follows from the very definition of strong neutrosophic bilinear algebra.

*Example 2.3.52:* Let $V = V_1 \cup V_2 =$

$$\left\{ \begin{pmatrix} a \\ b \\ c \\ d \\ e \end{pmatrix} \middle| a,b,c,d,e \in Z_{13}I \right\} \cup$$

$$\left\{ \begin{pmatrix} a_1 & a_2 & a_3 & a_4 & a_5 \\ a_6 & a_7 & a_8 & a_9 & a_{10} \end{pmatrix} \middle| a_i \in Z_7I; 1 \le i \le 10 \right\}$$

be a strong neutrosophic bivector space over the neutrosophic bifield $F = Z_{13}I \cup Z_7I$. We see $V = V_1 \cup V_2$ is not a strong neutrosophic bilinear algebra over the bifield $F = F_1 \cup F_2 = Z_{13}I \cup Z_7I$ as we see multiplication of elements within $V_i$ are not defined for $i = 1, 2$.

Now we define yet a new concept called quasi strong neutrosophic bilinear algebra.

**DEFINITION 2.3.24:** *Let $V = V_1 \cup V_2$ where $V_1$ is a strong neutrosophic vector space over the neutrosophic field $F_1$ ($V_1$ is only a vector space and $V_2$ is a strong neutrosophic linear algebra) over the neutrosophic field $F_2$ then we call $V = V_1 \cup V_2$ to be a quasi strong neutrosophic bilinear algebra over the neutrosophic bifield $F = F_1 \cup F_2$.*

We will illustrate this by the following examples.

*Example 2.3.53:* Let $V = V_1 \cup V_2 =$

$$\left\{ \begin{pmatrix} a_1 & a_2 & a_3 & a_4 \\ a_5 & a_6 & a_7 & a_8 \\ a_9 & a_{10} & a_{11} & a_{12} \end{pmatrix} \middle| a_i \in Z_{29}I; 1 \le i \le 12 \right\} \cup$$

{QI[x]; all polynomials in the variable x with coefficients from QI} be a quasi strong neutrosophic bilinear algebra over the neutrosophic bifield $F = Z_{29}I \cup QI$.

*Example 2.3.54:* Let $V = V_1 \cup V_2 =$

$$\left\{ \begin{pmatrix} a_1 & a_2 & a_3 & a_4 & a_5 & a_6 & a_7 \\ a_8 & a_9 & a_{10} & a_{11} & a_{12} & a_{13} & a_{14} \end{pmatrix} \middle| a_i \in Z_5I; 1 \le i \le 14 \right\} \cup$$

$$\left\{ \begin{pmatrix} a & b \\ c & d \end{pmatrix} \middle| a,b,c,d \in Z_2I \right\}$$

be a quasi strong neutrosophic bilinear algebra over the neutrosophic bifield $F = Z_5I \cup Z_2I$.

**DEFINITION 2.3.25:** *Let $V = V_1 \cup V_2$ be a strong neutrosophic bivector space over the neutrosophic bifield $F = F_1 \cup F_2$. Suppose $W = W_1 \cup W_2 \subseteq V_1 \cup V_2$ is such that W is a strong neutrosophic bilinear algebra over the neutrosophic field $F = F_1 \cup F_2$ then we call W to be a pseudo strong neutrosophic bilinear subalgebra of V over $F = F_1 \cup F_2$.*

**Example 2.3.55:** Let $V = V_1 \cup V_2 =$

$$\left\{ \begin{pmatrix} 0 & x \\ y & 0 \end{pmatrix} \middle| x, y \in Z_7 I \right\} \cup$$

$$\left\{ \begin{pmatrix} a_1 & 0 & a_2 \\ 0 & a_3 & 0 \\ a_4 & a_5 & 0 \end{pmatrix} \middle| a_1, a_2, a_3, a_4, a_5 \in Z_5 I \right\}$$

be a strong neutrosophic bivector space over the bifield $F = F_1 \cup F_2 = Z_7 I \cup Z_5 I$. Take $W = W_1 \cup W_2 =$

$$\left\{ \begin{pmatrix} 0 & x \\ 0 & 0 \end{pmatrix} \middle| x \in Z_7 I \right\} \cup \left\{ \begin{pmatrix} a_1 & 0 & a_2 \\ 0 & 0 & 0 \\ 0 & 0 & 0 \end{pmatrix} \middle| a_1, a_2 \in Z_5 I \right\}$$

to be a strong neutrosophic bilinear algebra over the neutrosophic bifield $F = Z_7 I \cup Z_5 I$, both $W_1$ and $W_2$ is closed under matrix multiplication. Thus $W = W_1 \cup W_2 \subseteq V_1 \cup V_2$ is a pseudo strong neutrosophic bilinear subalgebra of $V$ over the neutrosophic bifield $F$.

**Example 2.3.56:** Let $V = V_1 \cup V_2 =$

$$\left\{ \begin{pmatrix} a & b \\ c & d \end{pmatrix} \middle| a, b, c, d \in Z_{19} I \right\} \cup$$

$$\left\{ \begin{pmatrix} a & a & a \\ 0 & 0 & b \\ c & 0 & d \end{pmatrix} \middle| a, b, c, d \in Z_{41} I \right\}$$

be a strong neutrosophic bivector space over the neutrosophic bifield $F = F_1 \cup F_2 = Z_{19} I \cup Z_{41} I$. Take $W = W_1 \cup W_2 =$

$$\left\{ \begin{pmatrix} 0 & 0 \\ d & 0 \end{pmatrix} \middle| d \in Z_{19}I \right\} \cup \left\{ \begin{pmatrix} a & 0 & 0 \\ 0 & 0 & 0 \\ 0 & 0 & d \end{pmatrix} \middle| a,d \in Z_{41}I \right\}$$

$\subseteq V_1 \cup V_2$. W is a pseudo strong neutrosophic bilinear sub algebra of V over the bifield $F = Z_{19}I \cup Z_{41}I$.

**DEFINITION 2.3.26:** *Let $V = V_1 \cup V_2$ be a strong neutrosophic bilinear algebra over the neutrosophic bifield $F = F_1 \cup F_2$. Let $W = W_1 \cup W_2 \subseteq V_1 \cup V_2$, if W is a strong neutrosophic bivector space over F then we call W to be pseudo strong neutrosophic bivector subspace of V over F provided W is not a strong neutrosophic bilinear subalgebra of V over F.*

We will illustrate this situation by some Examples.

*Example 2.3.57:* Let $V = V_1 \cup V_2 =$

$$\left\{ \begin{pmatrix} a & b & c \\ d & e & f \\ g & h & i \end{pmatrix} \middle| a,b,c,d,e,f,g,h,i \in N(Q) \right\} \cup$$

$$\left\{ \begin{pmatrix} a & b \\ c & d \end{pmatrix} \middle| a,b,c,d \in Z_7I \right\}$$

be a strong neutrosophic bilinear algebra over the neutrosophic bifield $F = F_1 \cup F_2 = N(Q) \cup Z_7I$.
Take $W = W_1 \cup W_2 =$

$$\left\{ \begin{pmatrix} 0 & 0 & a \\ 0 & b & 0 \\ c & 0 & 0 \end{pmatrix} \middle| a,b,c \in QI \right\} \cup \left\{ \begin{pmatrix} 0 & a \\ b & 0 \end{pmatrix} \middle| a,b \in Z_7I \right\}$$

⊆ $V_1 \cup V_2$; W is a pseudo strong neutrosophic bivector subspace of V over F.

*Example 2.3.58:* Let $V = V_1 \cup V_2 =$

$$\left\{ \sum_{i=0}^{\infty} a_i x_i \mid a_i \in Z_2 I ; x \text{ is a variable or indeterminate} \right\} \cup$$

$$\left\{ \begin{pmatrix} a & b \\ c & d \end{pmatrix} \middle| a, b, c, d \in Z_3 I \right\}$$

be a strong neutrosophic bilinear algebra over the neutrosophic bifield $F = Z_2 I \cup Z_3 I$. Take $W = W_1 \cup W_2 =$

$$\left\{ \sum_{i=0}^{9} a_i x^i \middle| a_i \in Z_2 I; 0 \le i \le 9 \right\} \cup \left\{ \begin{pmatrix} 0 & b \\ c & 0 \end{pmatrix} \middle| b, c \in Z_3 I \right\}$$

⊆ $V_1 \cup V_2$; W is a pseudo strong neutrosophic bivector subspace of V over the bifield F.

**DEFINITION 2.3.27:** *Let $V = V_1 \cup V_2$ be a strong neutrosophic bilinear algebra over the neutrosophic bifield $F = F_1 \cup F_2$. Let $W = W_1 \cup W_2 \subseteq V_1 \cup V_2$, where one of $W_1$ or $W_2$ is alone a strong neutrosophic linear subalgebra and the other is just a strong neutrosophic vector subspace; then we call W to be a quasi strong neutrosophic bilinear subalgebra of V over the neutrosophic bifield $F = F_1 \cup F_2$.*

We will illustrate this situation by some examples.

*Example 2.3.59:* Let $V = V_1 \cup V_2 =$

$$\left\{ \begin{pmatrix} a & b \\ c & d \end{pmatrix} \middle| a, b, c, d \in Z_7 I \right\} \cup \left\{ \sum_{i=0}^{\infty} a_i x^i \mid a_i \in Z_2 I; 0 \le i \le \infty \right\}$$

be a strong neutrosophic bilinear algebra over the neutrosophic bifield $F = F_1 \cup F_2 = Z_7I \cup Z_2I$. Take $W = W_1 \cup W_2 =$

$$\left\{ \begin{pmatrix} a & a \\ a & a \end{pmatrix} \middle| a \in Z_7I \right\} \cup \left\{ \sum_{i=0}^{8} a_i x^i \middle| a_i \in Z_2I; 0 \le i \le 8 \right\}$$

$\subseteq V_1 \cup V_2$; W is a quasi strong neutrosophic bilinear subalgebra of V over $F = Z_7I \cup Z_2I$.

***Example 2.3.60:*** Let $V = V_1 \cup V_2 =$

$$\left\{ \begin{pmatrix} a & b & c \\ d & e & f \\ g & h & i \end{pmatrix} \middle| a,b,c,d,e,f,g,h,i \in Z_7I \right\} \cup$$

{N(Q)[x]; all polynomials in the variable x with coefficients from N(Q)} be a strong neutrosophic bilinear algebra over the bifield $F = Z_7I \cup QI$. Take $W = W_1 \cup W_2 =$

$$\left\{ \begin{pmatrix} 0 & 0 & a \\ 0 & b & 0 \\ c & 0 & 0 \end{pmatrix} \middle| a,b,c \in Z_7I \right\} \cup$$

{QI[x]; all polynomials in the variable x with coefficients from QI} $\subseteq V_1 \cup V_2$; W is a quasi strong neutrosophic bilinear subalgebra of V over F.

Now we proceed onto define the notion of strong neutrosophic bilinear subalgebra.

**DEFINITION 2.3.28:** *Let $V = V_1 \cup V_2$ be a strong neutrosophic bilinear algebra over the neutrosophic bifield $F = F_1 \cup F_2$. Let $W = W_1 \cup W_2 \subset V_1 \cup V_2$ be a proper bisubset of V; if W is a strong neutrosophic bilinear algebra over the bifield $F = F_1 \cup$* $F_2$; then we call $W = W_1 \cup W_2$ to be a *strong neutrosophic bilinear subalgebra of V over the bifield* $F = F_1 \cup F_2$.

We will illustrate this situation by some examples.

***Example 2.3.61:*** Let $V = V_1 \cup V_2 =$

$$\left\{ \begin{pmatrix} a & b \\ c & d \end{pmatrix} \middle| a,b,c,d \in N(Z_{41}) \right\} \cup$$

$\{N(Z_{11})[x]$; all polynomials in the variable x with coefficients from $N(Z_{11})\}$ be a strong neutrosophic bilinear algebra over the neutrosophic bifield $F = Z_{41}I \cup Z_{11}I$. Take $W = W_1 \cup W_2 =$

$$\left\{ \begin{pmatrix} a & b \\ c & d \end{pmatrix} \middle| a,b,c,d \in Z_{41}I \right\} \cup$$

$\{Z_{11}I[x]$; all polynomials in the variable x with coefficients from the neutrosophic field $Z_{11}I\} \subseteq V_1 \cup V_2$; W is a strong neutrosophic bilinear subalgebra of V over the bifield $F = Z_{41}I \cup Z_{11}I$.

***Example 2.3.62:*** Let $V = V_1 \cup V_2 =$

$$\left\{ \begin{pmatrix} a & b & c \\ d & e & f \\ g & h & i \end{pmatrix} \middle| a,b,c,d,e,f,g,h,i \in Z_{29}I \right\} \cup$$

$\{(a\ b\ c\ d\ e\ f) \mid a, b, c, d, e, f \in Z_{53}I\}$ be a strong neutrosophic bilinear algebra over the neutrosophic bifield $F = Z_{29}I \cup Z_{53}I$. Let $W = W_1 \cup W_2 =$

$$\left\{ \begin{pmatrix} a & b & c \\ 0 & d & e \\ 0 & 0 & f \end{pmatrix} \middle| a,b,c,d,e,f \in Z_{29}I \right\} \cup$$

$\{(a\ 0\ b\ 0\ d\ 0) \mid a, b, d \in Z_{53}I\} \subseteq V_1 \cup V_2$; W is a strong neutrosophic bilinear subalgebra of V over the neutrosophic bifield $F = Z_{29}I \cup Z_{53}I$.

**DEFINITION 2.3.29:** *Let $V = V_1 \cup V_2$ be a strong neutrosophic bilinear algebra over the neutrosophic bifield $F = F_1 \cup F_2$. Let $W = W_1 \cup W_2 \subseteq V_1 \cup V_2$; be such that W is a strong neutrosophic bilinear algebra over the proper neutrosophic bisubfield $K = K_1 \cup K_2 \subseteq F_1 \cup F_2$; $K_i$ is a proper neutrosophic subfield of $F_i$, $i = 1, 2$. We call $W = W_1 \cup W_2$ to be a strong subneutrosophic bilinear subalgebra of V over the neutrosophic subbifield $K = K_1 \cup K_2 \subseteq F_1 \cup F_2$.*

We will illustrate this by some examples.

***Example 2.3.63:*** Let $V = V_1 \cup V_2 =$

$$\left\{ \begin{pmatrix} a & b & c \\ d & e & f \\ g & h & i \end{pmatrix} \middle| a,b,c,d,e,f,g,h,i \in N(Q) \right\} \cup$$

$\{N(Z_{47})[x]$; all polynomials in the variable the x with coefficients from the neutrosophic field $N(Z_{47})\}$ be a strong neutrosophic bilinear algebra over the neutrosophic bifield $F = N(Q) \cup N(Z_{47})$. Take $W = W_1 \cup W_2 =$

$$\left\{ \begin{pmatrix} a & a & a \\ a & a & a \\ a & a & a \end{pmatrix} \middle| a \in QI \right\} \cup$$

$\{Z_{47}I[x]$; all polynomials in the variable x with coefficients from the neutrosophic field $Z_{47}I\} \subseteq V_1 \cup V_2$; W is a strong subneutrosophic bilinear subalgebra of V over the neutrosophic subbifield $K = K_1 \cup K_2 = QI \cup Z_{47}I \subseteq N(Q) \cup N(Z_{47})$.

*Example 2.3.64:* Let $V = V_1 \cup V_2 =$

$$\left\{ \begin{pmatrix} a & b & c & d \\ 0 & e & f & g \\ 0 & 0 & h & i \\ 0 & 0 & 0 & j \end{pmatrix} \middle| a,b,c,d,e,f,g,h,i,j \in N(Z_{11}) \right\} \cup$$

$$\left\{ \begin{pmatrix} a & 0 & 0 \\ b & c & 0 \\ d & e & f \end{pmatrix} \middle| a,b,c,d,e,f \in N(Z_{17}) \right\}$$

be a strong neutrosophic bilinear algebra over the neutrosophic bifield $F = F_1 \cup F_2 = N(Z_{11}) \cup N(Z_17)$.

Take $W = W_1 \cup W_2 =$

$$\left\{ \begin{pmatrix} a & a & a & a \\ 0 & a & a & a \\ 0 & 0 & a & a \\ 0 & 0 & 0 & a \end{pmatrix} \middle| a \in Z_{11}I \right\} \cup$$

$$\left\{ \begin{pmatrix} a & 0 & 0 \\ b & c & 0 \\ d & e & f \end{pmatrix} \middle| a,b,c,d,e,f \in Z_{17}I \right\}$$

$\subseteq V_1 \cup V_2$; W is a strong subneutrosophic bilinear subalgebra of V over the neutrosophic subbifield $K = K_1 \cup K_2 = Z_{11}I \cup Z_{17}I \subseteq N(Z_{11}) \cup N(Z_{17})$.

A neutrosophic bifield $F = F_1 \cup F_2$ is said to be neutrosophic biprime if both $F_1$ and $F_2$ have no proper neutrosophic subbifields contained in them. $F = Z_{11}I \cup Z_2I$ is neutrosophic biprime. $F = QI \cup Z_3I$ is neutrosophic biprime.

We see if $F_1$ is a neutrosophic prime field then it is of the form $QI$ or $Z_pI$; p a prime.

Now we will define bisimple strong subneutrosophic linear bialgebra.

**DEFINITION 2.3.30:** *Let $V = V_1 \cup V_2$ be a strong neutrosophic bilinear algebra over the neutrosophic bifield $F = F_1 \cup F_2$. If V has no proper strong subneutrosophic bilinear subalgebra then we define V to be a bisimple strong subneutrosophic linear bialgebra.*

We will illustrate this by some simple examples.

*Example 2.3.65:* Let $V = V_1 \cup V_2 =$

$$\left\{ \begin{pmatrix} a & b \\ c & d \end{pmatrix} \middle| a,b,c,d \in N(Q) \right\} \cup$$

$\{(x_1, x_2, x_3, x_4, x_5, x_6) \mid x_i \in Z_7I, 1 \leq i \leq 6\}$ be a strong neutrosophic linear bialgebra over the neutrosophic bifield $F = F_1 \cup F_2 = QI \cup Z_7I$. Since F has no neutrosophic subbifield V is a bisimple strong subneutrosophic linear bialgebra over F.

*Example 2.3.66:* Let $V = V_1 \cup V_2 =$

$$\left\{ \begin{pmatrix} a & b & c \\ d & e & f \\ g & h & i \end{pmatrix} \middle| a,b,c,d,e,f,g,h,i \in Z_{17}I \right\} \cup$$

$\{N(Z_{11})[x]$; all polynomials in the variable x with coefficients from $N(Z_{11})\}$ be a strong neutrosophic linear bialgebra over the

neutrosophic bifield $F = F_1 \cup F_2 = Z_{17}I \cup Z_{11}I$. Clearly V is a bisimple strong subneutrosophic bilinear algebra over F.

In view of this we have the following theorem.

**THEOREM 2.3.8:** *Let $V = V_1 \cup V_2$ be a strong neutrosophic bilinear algebra over the bifield $F = F_1 \cup F_2$ if each $F_i$ is of the form $K_iI$ where $K_i$ is a prime field $i = 1, 2$ then V is a bisimple strong subneutrosophic bilinear algebra over F.*

*Proof:* Follows from the fact that $F = F_1 \cup F_2$, the neutrosophic bifield has no proper neutrosophic subbifield.
   Now as in case of strong neutrosophic bivector spaces we can define the bibasis of a strong neutrosophic bilinear algebra and linearly biindependent bisubset. This task is left as an exercise for the interested reader.

We define linear bitransformation of a strong neutrosophic bilinear algebra into a strong neutrosophic bilinear algebra which we choose to call as strong neutrosophic linear bitransformation or when the context of reference is clear we just call it as strong bilinear transformation or in short just bilinear transformation or linear bitransformation.

**DEFINITION 2.3.31:** *Let $V = V_1 \cup V_2$ and $W = W_1 \cup W_2$ be two strong neutrosophic bilinear algebras over the same neutrosophic bifield $F = F_1 \cup F_2$. A bimap $T = T_1 \cup T_2 : V = V_1 \cup V_2 \to W = W_1 \cup W_2$ is defined to be a strong neutrosophic bilinear transformation or strong bilinear transformation or just bilinear transformation if each $T_i : V_i \to W_i$ is a linear transformation of $V_i$ to $W_i$ over $F_i$ for $i = 1, 2$.*

We will first illustrate this by some examples.

***Example 2.3.67:*** Let $V = V_1 \cup V_2 =$

$$\left\{ \begin{pmatrix} a & b \\ c & d \end{pmatrix} \middle| a,b,c,d \in Z_{17}I \right\} \cup$$

$$\left\{ \begin{pmatrix} a & b & c \\ d & e & f \\ g & h & i \end{pmatrix} \middle| a,b,c,d,e,f,g,h,i \in Z_{11}I \right\}$$

be a strong neutrosophic bilinear algebra over the neutrosophic bifield $F = F_1 \cup F_2 = Z_{17}I \cup Z_{11}I$. Take $W = W_1 \cup W_2 = \{(a, b, c, d) \mid a, b, c, d \in Z_{17}I\} \cup \{(a_1, a_2, a_3, a_4, a_5, a_6, a_7, a_8, a_9) \mid a_i \in Z_{11}I; 1 \le i \le 9\}$ to be a strong neutrosophic bilinear algebra over the same neutrosophic bifield $F = Z_{17}I \cup Z_{11}I$. The bimap $T = T_1 \cup T_2: V = V_1 \cup V_2 \to W_1 \cup W_2 = W$ where $T_1: V_1 \to W_1$ and $T_2: V_2 \to W_2$ defined by

$$T_1 \begin{pmatrix} a & b \\ c & d \end{pmatrix} = (a, b, c, d)$$

and

$$T_2 \begin{pmatrix} a & b & c \\ d & e & f \\ g & h & i \end{pmatrix} = (a, b, c, d, e, f, g, h, i)$$

is a strong neutrosophic bilinear transformation of V to W.

***Example 2.3.68:*** Let $V = V_1 \cup V_2 = \{N(Q)[x]$, all polynomials in the variable x with coefficients from $N(Q)\} \cup$

$$\left\{ \begin{pmatrix} a & b & c \\ d & e & f \\ g & h & i \end{pmatrix} \middle| a,b,c,d,e,f,g,h,i \in Z_2I \right\}$$

be a strong neutrosophic bilinear algebra over the neutrosophic bifield $F = QI \cup Z_2I$. Let $W = W_1 \cup W_2 =$

$$\left\{ \sum_{i=0}^{\infty} a_i x^{2i} \mid a_i \in N(Q); 0 \le i \le \infty \right\} \cup$$

$\{(a_1, a_2, a_3, a_4, a_5, a_6) \mid a_i \in Z_2I; 1 \leq i \leq 6\}$ be a strong neutrosophic bilinear space over the same neutrosophic bifield $F = QI \cup Z_2I$. Define the bimap

$$T = T_1 \cup T_2 : V = V_1 \cup V_2 \to W = W_1 \cup W_2$$

where $T_1: V_1 \to W_1$ and $T_2: V_2 \to W_2$ are defined by

$$T_1\left(\sum_{i=0}^{\infty} a_i x^i\right) = \left(\sum_{i=0}^{\infty} a_i x^{2i}\right) \text{ that is } x \to x^2$$

and

$$T_2 \begin{pmatrix} a & b & c \\ g & d & e \\ h & i & f \end{pmatrix} \to (a, b, c, d, e, f)$$

$T = T_1 \cup T_2$ is a strong neutrosophic bilinear transformation of V into W.

If in the definition of a bilinear transformation we put $W = V$ i.e., $W = W_1 \cup W_2 = V_1 \cup V_2$ i.e., $V_i = W_i$; $i = 1, 2$. That is the range bispace W is the same as the domain bispace then we call the strong neutrosophic bilinear transformation as the strong neutrosophic bilinear operator or strong neutrosophic linear bioperator on V.

We will illustrate this by some examples.

***Example 2.3.69:*** Let $V = V_1 \cup V_2 =$

$$\left\{\begin{pmatrix} a & b \\ c & d \end{pmatrix} \middle| a,b,c,d \in Z_5I\right\} \cup \left\{\sum_{i=0}^{\infty} a_i x^i \mid a_i \in N(Q); 0 \leq i \leq \infty\right\}$$

be a strong neutrosophic bilinear algebra over the neutrosophic bifield $F = Z_5I \cup QI$. Define $T = T_1 \cup T_2$ a bimap from $V = V_1 \cup V_2$ into $V = V_1 \cup V_2$ where $T_1: V_1 \to V_1$ and $T_2: V_2 \to V_2$ is given by

$$T_1\begin{pmatrix} a & b \\ c & d \end{pmatrix} = \begin{pmatrix} d & c \\ b & a \end{pmatrix}$$

and

$$T_2\left(\sum_{i=0}^{\infty} a_i x^i\right) = \left(\sum_{i=0}^{\infty} a_i x^{2i}\right)$$

i.e., $x \to x^2$. T is a strong neutrosophic linear bioperator on V.

*Example 2.3.70:* Let $V = V_1 \cup V_2 =$

$$\left\{ \begin{pmatrix} a_1 & a_2 & a_3 \\ 0 & a_4 & a_5 \\ 0 & 0 & a_6 \end{pmatrix} \middle| a_i \in N(Q); 1 \le i \le 6 \right\} \cup$$

$\{(a_1, a_2, a_3, a_4, a_5) \mid a_i \in Z_{11}I; 1 \le i \le 5\}$ be a strong neutrosophic bilinear algebra over the bifield $F = QI \cup Z_{11}I$. Define $T = T_1 \cup T_2 : V = V_1 \cup V_2 \to V = V_1 \cup V_2$ and $T_1 : V_1 \to V_1$ and $T_2 : V_2 \to V_2$ given by

$$T_1\begin{pmatrix} a_1 & a_2 & a_3 \\ 0 & a_4 & a_5 \\ 0 & 0 & a_6 \end{pmatrix} = \begin{pmatrix} a_6 & a_5 & a_3 \\ 0 & a_4 & a_2 \\ 0 & 0 & a_1 \end{pmatrix}$$

and

$$T_2(a_1, a_2, a_3, a_4, a_5) = (a_5, a_3, a_4, a_2, a_1)$$

$T = T_1 \cup T_2$ is a strong neutrosophic bilinear operator on V.

It is interesting to study the collection of all strong neutrosophic linear transformation of strong neutrosophic bilinear algebra $V = V_1 \cup V_2$ into a strong neutrosophic bilinear algebra $W = W_1 \cup W_2$ defined over the bifield $F = F_1 \cup F_2$.

We will denote this collection by

$$SNH_{F=F_1 \cup F_2}(V, W) = SNH_{F_1}(V_1, W_1) \cup SNH_{F_2}(V_2, W_2)$$

= {Collection of all bilinear transformation of $V_1 \cup V_2$ into $W_1 \cup W_2$}
= {Collection of all linear transformation of $V_1$ into $W_1$} $\cup$ {Collection of all linear transformation of $V_2$ into $W_2$}.

Interested reader can study and analyse the algebraic structure of $\text{SNH}_{F_1 \cup F_2}(V, W)$. On similar lines the set of all strong neutrosophic bilinear operators (linear bioperators) of a strong neutrosophic linear bialgebra over the neutrosophic bifield $F = F_1 \cup F_2$ is denoted by

$$\text{SNH}_{F_1 \cup F_2}(V, V) = \text{SNH}_{F_1}(V_1, V_1) \cup \text{SNH}_{F_2}(V_2, V_2)$$
$$= \text{SNH}_{F_1 \cup F_2}(V_1 \cup V_2, V_1 \cup V_2)$$

= {Collection of all strong neutrosophic linear bioperators of $V = V_1 \cup V_2$ into $V = V_1 \cup V_2$}.
= {Collection of all strong neutrosophic linear operators of $V_1$ into $V_1$} $\cup$ {Collection of all strong neutrosophic linear operators on $V_2$ into $V_2$}.

Interested reader is requested to study the algebraic structure of $\text{SNH}_{F_1 \cup F_2}(V, V)$. We will prove the following interesting property about strong neutrosophic linear bitransformation.

**THEOREM 2.3.9:** *Let $V = V_1 \cup V_2$ be a $(n_1, n_2)$ bidimensional finite strong neutrosophic bivector space over the neutrosophic bifield $F = F_1 \cup F_2$. Let $\{\alpha_1^1 \ldots \alpha_{n_1}^1\} \cup \{\alpha_1^2 \ldots \alpha_{n_2}^2\}$ be a bibasis of V over $F = F_1 \cup F_2$. Let $W = W_1 \cup W_2$ be a strong neutrosophic bivector space over the same neutrosophic bifield $F = F_1 \cup F_2$.*

*Let $\{\beta_1^1 \ldots \beta_{n_1}^1\} \cup \{\beta_1^2 \ldots \beta_{n_2}^2\}$ be any bivector in W. Then there is precisely a bilinear transformation $T = T_1 \cup T_2$ from $V = V_1 \cup V_2$ into $W = W_1 \cup W_2$ such that $T_i(\alpha_j^i) = (\beta_j^i)$ for $j = 1, 2, \ldots, n_i$ and $i = 1, 2$.*

**Proof:** Given $V = V_1 \cup V_2$ and $W = W_1 \cup W_2$ are two strong neutrosophic bivector spaces defined over the neutrosophic bifield $F = F_1 \cup F_2$. Let $T = T_1 \cup T_2 : V = V_1 \cup V_2 \to W = W_1 \cup W_2$.

Let $\{\alpha_1^1, \alpha_2^1, \ldots, \alpha_{n_i}^1\} \cup \{\alpha_1^2, \alpha_2^2, \ldots, \alpha_{n_2}^2\}$ be a bibasis of V. Given $\{\beta_1^1, \beta_2^1, \ldots, \beta_{n_i}^1\} \cup \{\beta_1^2, \beta_2^2, \ldots, \beta_{n_2}^2\}$ is a bivector in $W = W_1 \cup W_2$. To prove there is a bilinear transformation $T = T_1 \cup T_2$ with $T_i(\alpha_j^i) = (\beta_j^i)$ for each $j = 1, 2, \ldots, n_i$ and $i = 1, 2$. For every $\alpha = \alpha^1 \cup \alpha^2$ in $V = V_1 \cup V_2$ we have for every $\alpha^i \in V_i$ ($i = 1, 2$) a unique $x_1^i, x_2^i, \ldots, x_{n_i}^i$ such that

$$\alpha^i = x_1^i \alpha_1^i + x_2^i \alpha_2^i + \ldots + x_{n_i}^i \alpha_{n_i}^i.$$

This is true for every $i$; $i = 1, 2$. For this vector $\alpha^i$ define

$$T_i(\alpha^i) = x_1^i \beta_1^i + x_2^i \beta_2^i + \ldots + x_{n_i}^i \beta_{n_i}^i$$

true for $i = 1, 2$. Thus $T_i$ is well defined for associating with each vector $\alpha^i$ in $V_i$ a vector $T_i \alpha^i$ in $W_i$ ($i = 1, 2$). This rule for $T = T_1 \cup T_2$ is a well defined rule for each $T_i : V_i \to W_i$; $i = 1, 2$.

From the definition it is clear that $T_i \alpha_j^i = \beta_j^i$ for each $j$. To see T is bilinear. Let

$$\beta^i = y_1^i \alpha_1^i + y_2^i \alpha_2^i + \ldots + y_{n_i}^i \alpha_{n_i}^i$$

be in V and let $C^i$ be any scalar from $F_i$. Now

$$C^i \alpha^i + \beta^i = \left(C^i x_1^i + y_1^i\right)\beta_1^i + \ldots + \left(C^i x_{n_i}^i + y_{n_i}^i\right)\beta_{n_i}^i \,;$$

$i = 1, 2$.

On the other hand

$$T_i(C^i \alpha^i + \beta^i) = C^i \sum_{j=1}^{n_1} x_j^i \beta_j^i + \sum_{j=1}^{n_2} y_j^i \beta_j^i$$

true for each $i = 1, 2$; i.e., true for every linear transformation $T_i$ in T.

$$T_i(C^i \alpha^i + \beta^i) = C_i T_i(\alpha^i) + T_i(\beta^i)$$

true for every i.

Thus

$$T(C\alpha + \beta) = T_1(C_1 \alpha_1 + \beta_1) + T_2(C_2 \alpha_2 + \beta_2).$$

If $S = S_1 \cup S_2$ is a bilinear transformation from $V = V_1 \cup V_2$ into $W = W_1 \cup W_2$ with $S_i \alpha_j^i = \beta_j^i$; $j = 1, 2, \ldots, n_i$, $i = 1, 2$ then for any bivector $\alpha = \alpha^1 \cup \alpha^2$ we have for every $\alpha^i$ in $\alpha$ (i = 1, 2);

$$\alpha^i = \sum_{j=1}^{n_i} x_j^i \alpha_j^i$$

We have

$$\begin{aligned}
S_i \alpha^i &= S_i \sum_{j=1}^{n_i} y_j^i \alpha_j^i \\
&= \sum_{j=1}^{n_i} x_j^i S_i\left(\alpha_j^i\right) \\
&= \sum_{j=1}^{n_i} x_j^i \beta_j^i
\end{aligned}$$

so that S is exactly the rule T which we have defined. The prove $T\alpha = \beta$; i.e., if $\alpha = \alpha^1 \cup \alpha^2$ and $\beta = \beta^1 \cup \beta^2$ then $T_i \alpha_j^i = \beta_j^i$; $1 \leq j \leq n_i$; i = 1, 2.

The reader is requested to make the bimatrix analogue of the linear bitransformation from a strong neutrosophic bivector space V into a strong neutrosophic bivector space defined over the same neutrosophic bifield $F = F_1 \cup F_2$.

Now we proceed onto define the notion of binull space or null bispace and birank of a bilinear transformation T.

**DEFINITION 2.3.32:** *Let $V = V_1 \cup V_2$ and $W = W_1 \cup W_2$ be two strong neutrosophic bivector spaces defined over the neutrosophic bifield $F = F_1 \cup F_2$ of bidimensions $(n_1, n_2)$ and $(m_1, m_2)$ respectively. Let $T = T_1 \cup T_2 : V = V_1 \cup V_2 \to W = W_1 \cup W_2$ be a bilinear transformation. The binull space or null bispace of $T = T_1 \cup T_2$ is the set of all bivectors $a = \alpha_1 \cup \alpha_2$ in V such that $T_i \alpha^i = 0$; i = 1, 2.*

*If V is finite dimensional the birank of T is the dimension of the birange of $T = T_1 \cup T_2$ and binullity of T is the dimension of the null bispace of T.*

We have the following interesting relation between the birank of T and binullity of T.

**THEOREM 2.3.10:** *Let $V = V_1 \cup V_2$ and $W = W_1 \cup W_2$ be strong neutrosophic bivector space defined over the neutrosophic bifield $F = F_1 \cup F_2$ and suppose V is finite say $(n_1, n_2)$ dimensional T is a linear bitransformation from V into W. Then birank T + binullity T = bidimension V = $(n_1, n_2)$. Thus (rank $T_1$ $\cup$ rank $T_2$) + (nullity $T_1$ $\cup$ nullity $T_2$) = $(n_1, n_2)$.*

The proof is left as an exercise to the reader.

Now as in case of usual neutrosophic bivector spaces we have in case of strong neutrosophic bivector spaces the following result to be true.

Suppose $V = V_1 \cup V_2$ and $W = W_1 \cup W_2$ be any two strong neutrosophic bivector spaces over the bifield $F = F_1 \cup F_2$.

Let T and S be strong neutrosophic linear bitransformations from V into W. The bifunction

$(T + S) = (T_1 \cup T_2 + S_1 \cup S_2) = (T_1 + S_1) \cup (T_2 + S)$

is defined by

$(T + S)\alpha = T\alpha + S\alpha;$

i.e.,

$(T_1 \cup T_2 + S_1 \cup S_2)(\alpha_1 \cup \alpha_2) = (T_1 + S_1)(\alpha_1) \cup (T_2 + S_2)(\alpha_2)$
$= (T_1\alpha_1 + S_1\alpha_1) \cup T_2\alpha_2 + S_2\alpha_2.$

For any $C \in F_1 \cup F_2 = F$ the bifunction CT is defined by $(CT)\alpha = C(T\alpha)$ is a linear bitransformation from V into W. Further it can be proved that the set of all linear bitransformations from $V = V_1 \cup V_2$ into $W = W_1 \cup W_2$ together with addition and scalar multiplication defined above is a strong neutrosophic bivector space over the same neutrosophic bifield $F = F_1 \cup F_2$. Further it can be proved that if $V = V_1 \cup V_2$ be a finite bidimension $(n_1, n_2)$ strong neutrosophic bivector space over the bifield $F = F_1 \cup F_2$ and $W = W_1 \cup W_2$ be a finite $(m_1, m_2)$ bidimension strong neutrosophic bivector space over the same neutrosophic bifield $F = F_1 \cup F_2$, then the bispace $SNH_{F_1 \cup F_2}(V, W) = SNL_2(V, W)$ is a finite bidimensional bispace of bidimension $(m_1 n_1, m_2 n_2)$ over the same neutrosophic bifield $F = F_1 \cup F_2$. These results hold good when the strong neutrosophic bivector spaces are strong neutrosophic bilinear algebras.

We now proceed onto define biinvertible bilinear transformation.

**DEFINITION 2.3.33:** *Let $V = V_1 \cup V_2$ and $W = W_1 \cup W_2$ be two strong neutrosophic bivector spaces defined over the same neutrosophic bifield $F = F_1 \cup F_2$ of type II. A bilinear transformation $T = T_1 \cup T_2$ from V into W is biinvertible if and only if*
  i. *$T = T_1 \cup T_2$ is one to one that is each $T_i$ is one to one from $V_i$ into $W_i$ such that $T_i \alpha_i = T_i \beta_i$ implies $\alpha_i = \beta_i$ true for each i, i = 1, 2, ..., n.*
  ii. *T is onto, that is birange of T is all of $W = W_1 \cup W_2$ i.e., each $T_i : V_i \to W_i$ is onto and range $T_i$ is all of $W_i$ true for every i; i = 1, 2.*

We will first illustrate this situation by some examples.

*Example 2.3.71:* Let $V = V_1 \cup V_2 = $

$$\left\{ \begin{pmatrix} a & b \\ c & d \end{pmatrix} \middle| a,b,c,d \in Z_7 I \right\} \cup$$

$\{(a_1, a_2, a_3, a_4, a_5) \mid a_i \in Z_{11}I; 1 \leq i \leq 5\}$ be a strong neutrosophic bilinear algebra defined over the neutrosophic bifield $F = Z_7 I \cup Z_{11} I$. Define $T = T_1 \cup T_2 : V = V_1 \cup V_2 \to V = V_1 \cup V_2$ where
$$T_1 : V_1 \to V_1 \text{ and } T_2 : V_2 \to V_2.$$
such that
$$T_1 \begin{pmatrix} a & b \\ c & d \end{pmatrix} = \begin{pmatrix} b & a \\ d & c \end{pmatrix}$$
and
$$T_2 (a_1, a_2, a_3, a_4, a_5) = (a_5, a_4, a_3, a_2, a_1).$$
Clearly $T = T_1 \cup T_2$ is a strong neutrosophic linear bioperator of V into V.

Take $S = S_1 \cup S_2$ and define from $V = V_1 \cup V_2$ to $V = V_1 \cup V_2$ by

$$S_1 \begin{pmatrix} a & b \\ c & d \end{pmatrix} = \begin{pmatrix} b & a \\ d & c \end{pmatrix}$$

and

$$S_2 (a_1, a_2, a_3, a_4, a_5) = (a_5, a_4, a_3, a_2, a_1).$$

$S = S_1 \cup S_2 = T = T_1 \cup T_2$ such that $S$ is a strong neutrosophic linear bioperator of V into V.

We see $T_1 \cdot T_1 = T_1 \cdot S_1 = S_1 \cdot T_1$ is identity linear bioperator on V. We have $S_1 = T_1$.

For consider

$$T_1 \cdot S_1 \begin{pmatrix} a & b \\ c & d \end{pmatrix} = S_1 \left[ T_1 \begin{pmatrix} a & b \\ c & d \end{pmatrix} \right]$$

$$= S_1 \begin{pmatrix} b & a \\ d & c \end{pmatrix} = \begin{pmatrix} a & b \\ d & c \end{pmatrix}.$$

Thus

$$S_1 T_1 \begin{pmatrix} a & b \\ d & c \end{pmatrix} = \begin{pmatrix} a & b \\ d & c \end{pmatrix},$$

hence

$$S_1 \cdot T_1 = T_1 \cdot S_1 = T_1 \cdot T_1 = S_1 \cdot S_1$$

(as $S_1 = T_1$) is such that $T_1 = T_1^{-1}$.

Now consider $T_2 : V_2 \to V_2$ we see $T_2 = S_2$

$$T_2 \cdot S_2 [(a_1, a_2, a_3, a_4, a_5)] = T_2 (a_5, a_4, a_3, a_2, a_1)$$
$$= (a_1, a_2, a_3, a_4, a_5)$$
$$= \text{identity bioperator on } V_2.$$

Thus $T_2 = T_2^{-1}$. We see $T = T_1 \cup T_2$ has the inverse bioperator $T_{-1} = T_1^{-1} \cup T_2^{-1}$.

Now we can also give an example of a linear bitransformation of strong neutrosophic bivector spaces (or bilinear algebras).

**Example 2.3.72:** Let $V = V_1 \cup V_2 =$

$$\left\{ \begin{pmatrix} a & b \\ c & d \end{pmatrix} \middle| a,b,c,d \in Z_7I \right\} \cup$$

$$\left\{ \begin{pmatrix} a_1 & a_2 & a_3 \\ a_4 & a_5 & a_6 \end{pmatrix} \middle| a_i \in Z_{11}I; 1 \leq i \leq 6 \right\}$$

be a strong neutrosophic bivector space over the neutrosophic bifield $F = Z_7I \cup Z_{11}I$.

Take $W = W_1 \cup W_2 =$

$$\left\{ \begin{pmatrix} a & 0 & b \\ 0 & c & d \end{pmatrix} \middle| a,b,c,d \in Z_7I \right\} \cup$$

$$\left\{ \begin{pmatrix} a_1 & a_2 \\ a_3 & a_4 \\ a_5 & a_6 \end{pmatrix} \middle| a_i \in Z_{11}I; 1 \leq i \leq 6 \right\}$$

to be a strong neutrosophic bivector space over the same neutrosophic bifield $F = Z_7I \cup Z_{11}I$. Define $T = T_1 \cup T_2 : V = V_1 \cup V_2 \to W = W_1 \cup W_2$ where $T_1 : V_1 \to W_1$ and $T_2 : V_2 \to W_2$ such that

$$T_1 \begin{pmatrix} a & b \\ c & d \end{pmatrix} = \begin{pmatrix} a & 0 & b \\ 0 & c & d \end{pmatrix}$$

and

$$T_2 \begin{pmatrix} a_1 & a_2 & a_3 \\ a_4 & a_5 & a_6 \end{pmatrix} = \begin{pmatrix} a_1 & a_2 \\ a_3 & a_4 \\ a_5 & a_6 \end{pmatrix}.$$

$T = T_1 \cup T_2$ is a neutrosophic linear bitransformation of $V = V_1 \cup V_2$ into $W = W_1 \cup W_2$.
Define a bimap $S = S_1 \cup S_2$:

$$W = W_1 \cup W_2 \to V = V_1 \cup V_2$$

where $S_1 : W_1 \to V_1$ and $S_2 : W_2 \to V_2$ such that

$$S_1 \begin{pmatrix} a & 0 & b \\ 0 & c & d \end{pmatrix} = \begin{pmatrix} a & b \\ c & d \end{pmatrix}$$

and

$$S_2 \begin{pmatrix} a_1 & a_2 \\ a_3 & a_4 \\ a_5 & a_6 \end{pmatrix} = \begin{pmatrix} a_1 & a_2 & a_3 \\ a_4 & a_5 & a_6 \end{pmatrix}.$$

$S = S_1 \cup S_2$ is clearly a linear transformation from $W = W_1 \cup W_2$ to $V = V_1 \cup V_2$.

Now we find $T \cdot S$ and $S \cdot T$.

$$T \cdot S = T_1 \cdot S_1 \cup T_2 \cdot S_2$$

Now

$$T_1 \cdot S_1 \begin{pmatrix} a & b \\ c & d \end{pmatrix} = S_1 \cdot T_1 \begin{pmatrix} a & b \\ c & d \end{pmatrix}$$

$$= S_1 \begin{pmatrix} a & 0 & b \\ 0 & c & d \end{pmatrix} = \begin{pmatrix} a & b \\ c & d \end{pmatrix}.$$

That is $T_1 \cdot S_1$ is the identity transformation of $V_1$.

Now consider

$$S_1 \cdot T_1 \begin{pmatrix} a & 0 & b \\ 0 & c & d \end{pmatrix} = T_1 \left( S_1 \begin{pmatrix} a & 0 & b \\ 0 & c & d \end{pmatrix} \right)$$

$$= T_1 \begin{pmatrix} a & b \\ c & d \end{pmatrix} = \begin{pmatrix} a & 0 & b \\ 0 & c & d \end{pmatrix}.$$

Thus $S_1 \cdot T_1$ is the identity transformation of $W_1$.

Now consider

$$T_2 \cdot S_2 \begin{pmatrix} a_1 & a_2 & a_3 \\ a_4 & a_5 & a_6 \end{pmatrix} = S_2 \left( T_2 \begin{pmatrix} a_1 & a_2 & a_3 \\ a_4 & a_5 & a_6 \end{pmatrix} \right)$$

$$= S_2 \begin{pmatrix} a_1 & a_2 \\ a_3 & a_4 \\ a_5 & a_6 \end{pmatrix} = \begin{pmatrix} a_1 & a_2 & a_3 \\ a_4 & a_5 & a_6 \end{pmatrix}.$$

Thus $T_2 \cdot S_2$ is the identity linear transformation on $V_2$.
Consider

$$S_2 \cdot T_2 \begin{pmatrix} a_1 & a_2 \\ a_3 & a_4 \\ a_5 & a_6 \end{pmatrix} = T_2 S_2 \begin{pmatrix} a_1 & a_2 \\ a_3 & a_4 \\ a_5 & a_6 \end{pmatrix}$$

$$= T_2 \begin{pmatrix} a_1 & a_2 & a_3 \\ a_4 & a_5 & a_6 \end{pmatrix} = \begin{pmatrix} a_1 & a_2 \\ a_3 & a_4 \\ a_5 & a_6 \end{pmatrix}.$$

Thus $S_2 \cdot T_2$ is the identity linear transformation on $W_2$. Thus $T \cdot S = T_1 \cdot S_1 \cup T_2 \cdot S_2$ is the identity bilinear transformation on $V_1 \cup V_2$ and $S \cdot T = S_1 \cdot T_1 \cup S_2 \cdot T_2$ is the identity linear bitransformation on $W = W_1 \cup W_2$.

In view of this example the reader is requested to prove the following result.

Let $V = V_1 \cup V_2$ and $W = W_1 \cup W_2$ be two strong neutrosophic bivector spaces defined over the neutrosophic bifield $F = F_1 \cup F_2$ of type II. Let $T = T_1 \cup T_2$ be a strong linear bitransformation from $V = V_1 \cup V_2$ into $W = W_1 \cup W_2$. If T is biinvertible then the biinverse bifunction $T^{-1} = T_1^{-1} \cup T_2^{-1}$ is a bilinear transformation from W into V.

Suppose $T = T_1 \cup T_2$ is a linear bitransformation from the strong neutrosophic bivector spaces $V = V_1 \cup V_2$ into $W = W_1 \cup W_2$ then T is binon singular if and only if T carries each bilinearly independent bisubset of $V = V_1 \cup V_2$ into a bilinearly independent bisubset of $W = W_1 \cup W_2$.

The following nice result is left as an exercise for the reader.

**THEOREM 2.3.11:** *Let $V = V_1 \cup V_2$ and $W = W_1 \cup W_2$ be a strong neutrosophic bivector spaces defined over the same neutrosophic bifield $F = F_1 \cup F_2$ of type II. If $T = T_1 \cup T_2$ is a bilinear transformation of V into W then the following are equivalent.*
  i. *$T = T_1 \cup T_2$ is biinvertible*
  ii. *$T = T_1 \cup T_2$ is binon singular*
  iii. *$T = T_1 \cup T_2$ is onto that is the birange of $T = T_1 \cup T_2$ is $W = W_1 \cup W_2$.*

We prove an important result.

**THEOREM 2.3.12:** *Every $(n_1, n_2)$ bidimensional strong neutrosophic bivector space $V = V_1 \cup V_2$ defined over the neutrosophic bifield $F = F_1 \cup F_2$ is biisomorphic to $F_1^{n_1} \cup F_2^{n_2}$.*

*Proof:* Let $V = V_1 \cup V_2$ be a $(n_1, n_2)$ bidimensional strong neutrosophic bivector space over the neutrosophic bifield $F = F_1 \cup F_2$ of type II. Let $B = \{\alpha_1^1, \alpha_2^1, \ldots, \alpha_{n_1}^1\} \cup \{\alpha_1^2, \alpha_2^2, \ldots, \alpha_{n_2}^2\}$ be a bibasis of V. We define a bifunction $T = T_1 \cup T_2$ from $V = V_1 \cup V_2$ into $F_1^{n_1} \cup F_2^{n_2}$ is follows.

If $\alpha = \alpha_1 \cup \alpha_2$ is in $V = V_1 \cup V_2$, let $T\alpha = T_1(\alpha_1) \cup T_2(\alpha_2)$ be the $(n_1, n_2)$ pair, $\left(x_1^1, x_2^1, \ldots, x_{n_1}^1\right) \cup \left(x_1^2, x_2^2, \ldots, x_{n_2}^2\right)$ of the bicoordinate of $\alpha = \alpha_1 \cup \alpha_2$ relative to the biordered bibasis B; i.e., the $(n_1, n_2)$ pair such that
$$\alpha = x_1^1\alpha_1^1 + x_2^1\alpha_2^1 + \ldots + x_{n_1}^1\alpha_{n_1}^1 \cup x_1^2\alpha_1^2 + x_2^2\alpha_2^2 + \ldots + x_{n_2}^2\alpha_{n_2}^2.$$

Clearly T is a linear bitransformation; T is a one to one map of $V = V_1 \cup V_2$ onto $F_1^{n_1} \cup F_2^{n_2}$ or each $T_i$ is linear and one to one and maps $V_i$ to $F_i^{n_i}$; $i = 1, 2$ for every i. Thus as in case of vector space transformation by matrices give a representation of bitransformation by bimatrices where the bimatrices are neutrosophic bimatrices.

Let $V = V_1 \cup V_2$ be a bivector space of $(n_1, n_2)$ bidimension over $F = F_1 \cup F_2$. Let $W = W_1 \cup W$ be a bivector space over the same bifield $F = F_1 \cup F_2$ of $(m_1, m_2)$ dimension. Let
$$B = \{\alpha_1^1, \alpha_2^1, \ldots, \alpha_{n_1}^1\} \cup \{\alpha_1^2, \alpha_2^2, \ldots, \alpha_{n_2}^2\}$$
be a bibasis of $V = V_1 \cup V_2$ and
$$C = \{\beta_1^1, \beta_2^1, \ldots, \beta_{m_1}^1\} \cup \{\beta_1^2, \beta_2^2, \ldots, \beta_{m_2}^2\}$$
be a bibasis for W. If $T = T_1 \cup T_2$ is any bilinear transformation of type II from $V = V_1 \cup V_2$ into $W = W_1 \cup W_2$ then $T_i$ is determined by its action on $\alpha^i$, $i = 1, 2$. Each of the $(n_1, n_2)$ pair vector; $T_i \alpha_j^i$; $j = 1, 2, \ldots, n$. $i = 1, 2$ is uniquely expressible as a linear combination
$$T_i \alpha_j^i = \sum_{k_i-1}^{m_i} A_{k_i,j}^i \beta_{k_i}^i .$$

This is true for every i, $1 \leq j \leq n_i$; $i = 1, 2$; of $\beta_{k_i}^i$; the scalars being the coordinates of $A_{ij}^i, \ldots, A_{m_i,j}^i$. $T_i \alpha_j^i$ in the basis $\beta_1^i, \ldots, \beta_{m_i}^i$ of C. True for each i; $i = 1, 2$.

Accordingly the bitransformation $T = T_1 \cup T_2$ is determined by the $(m_1n_1, m_2n_2)$ scalars $A_{k_i,j}^i$. The $m_i \times n_i$ neutrosophic matrix $A^i$ defined by $A_{k_i,j}^i$ is called the component neutrosophic matrix $T_i$ of T relative to the component basis $\{\alpha_1^i, \alpha_2^i, \ldots, \alpha_{n_i}^i\}$ and $\{\beta_1^i, \beta_2^i, \ldots, \beta_{m_i}^i\}$ of B and C respectively. Since this is true for every i; $i = 1, 2$; We have $A = A_{k_1j}^1 \cup A_{k_2j}^2 = A^1 \cup A^2$ the neutrosophic bimatrix associated with $T = T_1 \cup T_2$. Each $A^i$ determines the linear transformation $T_i$; for $i = 1, 2$. If $\alpha^i = x_1^i \alpha_1^i + x_2^i \alpha_2^i + \ldots + x_{n_i}^i . \alpha_{n_i}^i$ is a neutrosophic vector in $V_i$ then
$$T_i \alpha^i = \left( T_i \sum_{j=1}^{n_i} x_j^i \alpha_j^i \right) = \left( \sum_{j=1}^{n_i} x_j^i T_i \alpha_j^i \right)$$

$$= \left( \sum_{j=1}^{n_i} x_j^i \sum_{k=1}^{m_i} A_{k,j}^i \beta_{k_i}^i \right) = \left( \sum_{k=1}^{m_i} \sum_{j=1}^{n_i} \left( A_{k,j}^i x_j^i \right) \beta_{k_i}^i \right); i = 1, 2.$$

If $X^i$ is the coordinate neutrosophic matrix of $\alpha^i$ in the component bibasis of B then the above computation shows that $A^i X^i$ is the coordinate neutrosophic matrix of the vector $T^i \alpha^i$; that is the component of the bibasis C because the scalar

$$\left( \sum_{j=1}^{n_i} A_{k,j}^i x_{k_i}^i \right)$$

is the $k^{th}$ row of the column neutrosophic matrix $A^i X^i$. This is true for every i; i = 1, 2. Let us also observe that if $A^i$ is any $m_i \times n_i$ neutrosophic matrix over the neutrosophic field $F_i$, then

$$T_i \left( \sum_{j=1}^{n_i} x_j^i \alpha_j^i \right) = \left( \sum_{k=1}^{m_i} \left( \sum_{j=1}^{n_i} A_{k,j}^i x_{k_j}^i \right) \beta_{k_i}^i \right)$$

defines a linear transformation $T_i$ from $V_i$ into $W_i$, the neutrosophic matrix of which is $A^i$ relative to $\{\alpha_1^i, \alpha_2^i, \ldots, \alpha_{n_i}^i\}$ and $\{\beta_1^i, \beta_2^i, \ldots, \beta_{m_i}^i\}$; this is true for every i and i = 1, 2. Hence $T = T_1 \cup T_2$ is a linear bitransformation from $V = V_1 \cup V_2$ into $W = W_1 \cup W_2$, the neutrosophic bimatrix which is $A = A_1 \cup A_2$ relative to the bibasis $B = \{\alpha_1^1, \alpha_2^1, \ldots, \alpha_{n_1}^1\} \cup \{\alpha_1^2, \alpha_2^2, \ldots, \alpha_{n_2}^2\}$ and $C = \{\beta_1^1, \beta_2^1, \ldots, \beta_{m_1}^1\} \cup \{\beta_1^2, \beta_2^2, \ldots, \beta_{m_2}^2\}$.

Now as in case of bivector spaces of type II we can in case of strong neutrosophic bivector spaces prove that we can construct a biisomorphism between the strong neutrosophic bispace $NL_2(V, W)$ and the neutrosophic bispace of all neutrosophic bimatrices of biorder $(m_1 \times n_1, m_2 \times n_2)$ over the same neutrosophic bifield $F = F_1 \cup F_2$ over which $V = V_1 \cup V_2$ and $W = W_1 \cup W_2$ are defined as strong neutrosophic bivector spaces of bidimensions $(n_1, n_2)$ and $(m_1, m_2)$ respectively.

Now we will proceed onto define the notion of bilinear functionals or which we may call as linear bifunctionals. We know linear bifunctionals could not be defined for neutrosophic bivectors spaces of type I.

Let $V = V_1 \cup V_2$ be a strong neutrosophic bivector space over the neutrosophic bifield $F = F_1 \cup F_2$ of type II. A bilinear transformation or linear bitransformation $f = f_1 \cup f_2$ from V into the bifield $F = F_1 \cup F_2$ is defined as the linear bifunctional or bilinear functional on V, i.e., $f = f_1 \cup f_2$ is a bifunction from $V = V_1 \cup V_2$ into $F = F_1 \cup F_2$ such that
$$f(c\alpha + \beta) = f_1(c_1\alpha_1 + \beta_1) + f_2(c_2\alpha_2 + \beta_2)$$
$$= \{c_1 f_1(\alpha_1) \cup c_2 f_2(\alpha_2)\} + \{f_1(\beta_1) \cup f_2(\beta_2)\}$$
where $c = c_1 \cup c_2$ and $\alpha = \alpha_1 \cup \alpha_2$ and $\beta = \beta_1 \cup \beta_2$, $\beta_i, \alpha_i \in V_i$, $i = 1, 2$. That is $f = f_1 \cup f_2$ where each $f_i$ is a linear functional on $V_i$; $i = 1, 2$.

The following observations are both interesting and important. Let $F = F_1 \cup F_2$ be a neutrosophic bifield and let $F_1^{n_1} \cup F_2^{n_2}$ be a strong neutrosophic bivector space of type II over the bifield $F_1 \cup F_2$. A bilinear functional $f = f_1 \cup f_2$ from $F_1^{n_1} \cup F_2^{n_2}$ to $F_1 \cup F_2$ given by
$$f_1\left(x_1^1, x_2^1, ..., x_{n_1}^1\right) \cup f_2\left(x_2^1, x_2^2, ..., x_{n_2}^2\right)$$
$$= x_1^1 \alpha_1^1 + x_2^1 \alpha_2^1 + ... + x_{n_1}^1 . \alpha_{n_1}^1 \cup x_1^2 \alpha_1^2 + x_2^2 \alpha_2^2 + ... + x_{n_2}^2 . \alpha_{n_2}^2$$
where $\alpha_j^i \in F_i$; $1 \leq j \leq n_i$ and $i = 1, 2$; is a bilinear functional of $F_1^{n_1} \cup F_2^{n_2}$.

It is the bilinear functional which is represented by the neutrosophic bimatrix
$$\left[\alpha_1^1, \alpha_2^1, ..., \alpha_{n_1}^1\right] \cup \left[\alpha_1^2, \alpha_2^2, ..., \alpha_{n_2}^2\right]$$
relative to the standard bibasis for $F_1^{n_1} \cup F_2^{n_2}$ on the bibasis $\{1\} \cup \{1\}$ or $\{I\} \cup \{I\}$ for $F = F_1 \cup F_2$ depending on $F_i = N(K_i)$ or $F_i = K_i I$ respectively; $K_i$ – real field; $i = 1, 2$.

$\alpha_j^i = f_i(E_j^i)$ ; $j = 1, 2, \ldots, n_i$ for every $i = 1, 2$. Every bilinear functional on $F_1^{n_1} \cup F_2^{n_2}$ is of this form for some biscalar $\{\alpha_1^1, \alpha_2^1, \ldots, \alpha_{n_1}^1\} \cup \{\alpha_1^2, \alpha_2^2, \ldots, \alpha_{n_2}^2\}$. This is immediate from the definition of bilinear functional of type II because we define $\alpha_j^i = f_i(E_j^i)$. Hence

$$f_1(x_1^1, x_2^1, \ldots, x_{n_1}^1) \cup f_2(x_2^1, x_2^2, \ldots, x_{n_2}^2)$$
$$= f_1\left(\sum_{j=1}^{n_1} x_j^1 E_j^1\right) \cup f_2\left(\sum_{j=1}^{n_2} x_j^2 E_j^2\right)$$
$$= \sum_{j=1}^{n_1} x_j^1 f_1(E_j^1) \cup \sum_{j=1}^{n_2} x_j^2 f_2(E_j^2)$$
$$= \sum_{j=1}^{n_1} x_j^1 \alpha_j^1 \cup \sum_{j=1}^{n_2} x_j^2 \alpha_j^2.$$

Now we proceed onto define the new notion of strong neutrosophic bidual space or equivalently strong dual bispace of the strong neutrosophic bivector space $V = V_1 \cup V_2$ defined over the neutrosophic bifield $F = F_1 \cup F_2$ of type II.

Now as in case of $SNL^2(V, W) = SNL(V_1, W_1) \cup SNL(V_2, W_2)$ we in case of bilinear functional have $SNL^2(V, F) = SNL(V_1, F_1) \cup SNL(V_2, F_2)$. We define $V^* = SNL^2(V, F) = V_1^* \cup V_2^* = SNL(V_1, F_1) \cup SNL(V_2, F_2)$

That is each $V_i^*$ is the strong neutrosophic dual space of $V_i$, $V_i$ defined over $F_i$, $i = 1, 2$. We know if the strong neutrosophic vector space $V_i$, $\dim V_i = \dim V_i^*$ for every $i$, $1 \leq i \leq 2$.
Thus
$$\dim V = \dim V_1 \cup \dim V_2$$
$$= \dim V^*$$
$$= \dim V_1^* \cup \dim V_2^*.$$

If $B = \{\alpha_1^1, \alpha_2^1, \ldots, \alpha_{n_1}^1\} \cup \{\alpha_1^2, \alpha_2^2, \ldots, \alpha_{n_2}^2\}$ is a bibasis for $V = V_1 \cup V_2$, then we know for a bilinear function of type II, $f = f_1$

∪ $f_2$ we have $f_k$ on $V_k$ is such that $f_i^k\left(\alpha_j^k\right) = \delta_{ij}^k$ true for k = 1, 2. In this way we obtain from the biset B = $B_1 \cup B_2$ a pair of $n_i$ sets of distinct bifunctionals (i = 1, 2); $\{f_1^1, f_2^1, ..., f_{n_1}^1\}$ ∪ $\{f_1^2, f_2^2, ..., f_{n_2}^2\}$ on V = $V_1 \cup V_2$. These bifunctionals are also bilinearly independent over the bifield F = $F_1 \cup F_2$, i.e., $\{f_1^i, f_2^i, ..., f_{n_i}^i\}$ is linearly independent on $V_i$ over the neutrosophic field $F_i$, for every i, $1 \leq i \leq 2$.
Thus

$$f^i = \left(\sum_{j=1}^{n_i} c_j^1 f_j^1\right), i = 1, 2.$$

That is

$$f = \sum_{j=1}^{n_1} c_j^1 f_j^1 \cup \sum_{j=1}^{n_2} c_j^2 f_j^2.$$

$$f^i\left(\alpha_j^i\right) = \sum_{k=1}^{n_i} c_k^i f_k^i\left(\alpha_j^k\right)$$

$$= \sum_{k=1}^{n_i} c_k^i \delta_{ki}$$

$$= c_j^i.$$

This is true for i = 1, 2 and $1 \leq j \leq n_i$.

In particular if each $f_i$ is a zero functional $f^i \alpha_j^1 = 0$ for each j and hence the scalar $c_j^i$ are all zero. Thus $\{f_1^i, f_2^i, ..., f_{n_i}^i\}$ are $n_i$ linearly independent linear functionals of $V_i$ defined on $F_i$, true for each i; $1 \leq i \leq 2$. Since $V_i^*$ is of dimension $n_i$; it must be that $\{f_1^i, f_2^i, ..., f_{n_i}^i\}$ is a basis of $V_i^*$ which is the dual basis of B. Thus $B^* = B_1^* \cup B_2^* = \{f_1^1, f_2^1, ..., f_{n_1}^1\} \cup \{f_1^2, f_2^2, ..., f_{n_2}^2\}$ is the bidual basis or dual bibasis of B = $\{\alpha_1^1, \alpha_2^1, ..., \alpha_{n_1}^1\} \cup \{\alpha_1^2, \alpha_2^2, ..., \alpha_{n_2}^2\}$. $B^*$ forms the bibasis of $V^* = V_1^* \cup V_2^*$.

Interested reader is left with the task of proving the following theorem.

**THEOREM 2.3.13:** *Let $V = V_1 \cup V_2$ be a finite $(n_1, n_2)$ bidimension strong neutrosophic bivector space defined over the neutrosophic bifield $F = F_1 \cup F_2$. Let $B = \{\alpha_1^1, \alpha_2^1, \ldots, \alpha_{n_1}^1\} \cup \{\alpha_1^2, \alpha_2^2, \ldots, \alpha_{n_2}^2\} = B_1 \cup B_2$ be a bibasis for $V = V_1 \cup V_2$. There is a unique bidual basis (dual bibasis) $B^* = B_1^* \cup B_2^* = \{f_1^1, f_2^1, \ldots, f_{n_1}^1\} \cup \{f_1^2, f_2^2, \ldots, f_{n_2}^2\}$ for $V^* = V_1^* \cup V_2^*$ such that $f_i^k(\alpha_j) = \delta_{ij}^k$. For each bilinear functional $f = f^1 \cup f^2$ we have*

$$f = \left(\sum_{k=1}^{n_i} f^i(\alpha_k^i) f_k^i\right).$$

*That is*

$$f = \left(\sum_{k=1}^{n_1} f^1(\alpha_k^1) f_k^1\right) \cup \left(\sum_{k=1}^{n_2} f^2(\alpha_k^2) f_k^2\right)$$

*and for each bivector $\alpha = \alpha^1 \cup \alpha^2$ in $V = V_1 \cup V_2$ we have*

$$\alpha = \left(\sum_{k=1}^{n_1} f_k^1(\alpha^1) \alpha_k^1\right) \cup \left(\sum_{k=1}^{n_2} f_k^2(\alpha^2) \alpha_k^2\right).$$

Now we proceed onto defined yet a new feature of the strong neutrosophic bivector space of type II.

Let $V = V_1 \cup V_2$ be a strong neutrosophic bivector space defined over the neutrosophic bifield $F = F_1 \cup F_2$ of type II. Let $S = S_1 \cup S_2$ be a bisubset of $V = V_1 \cup V_2$ (that is $S_i \subseteq V_i$; $i = 1, 2$); the biannihilator of $S$ is $S^\circ = S_1^\circ \cup S_2^\circ$ of bilinear functionals on $V = V_1 \cup V_2$ such that $f(\alpha) = 0 \cup 0$ i.e., if $f = f^1 \cup f^2$ for every $\alpha = \alpha_1 \cup \alpha_2 \in S = S_1 \cup S_2$ ($\alpha_1 \in S_1$, $\alpha_2 \in S_2$); $f^i(\alpha_i) = 0$ for every $\alpha_i \in S_i$; $i = 1, 2$.

It is interesting to note that $S^\circ = S_1^\circ \cup S_2^\circ$ is a strong neutrosophic bisubspace of $V^* = V_1^* \cup V_2^*$; whether $S = S_1 \cup S_2$ is a bisubspace of $V = V_1 \cup V_2$ or only just a bisubset of $V = V_1$

∪ $V_2$. If $S = (0 \cup 0)$ then $S° = V^* = V_1^* \cup V_2^*$. If $S = V$, i.e., $V_1$ ∪ $V_2 = S_1 \cup S_2$ then $S°$ is just the zero bisubspace of $V^* = V_1^* \cup V_2^*$.

We leave the following theorem for the reader to prove.

**THEOREM 2.3.14:** *Let $V = V_1 \cup V_2$ be a strong neutrosophic bivector space of $(n_1, n_2)$ bidimension over the neutrosophic bifield $F = F_1 \cup F_2$ of type II. Let $W = W_1 \cup W_2$ be a strong neutrosophic bisubspace of $V = V_1 \cup V_2$. Then $\dim W + \dim W° = \dim V$ that is $(\dim W_1 \cup \dim W_2) + \dim W_1° \cup \dim W_2° = \dim V_1 \cup \dim V_2 = (n_1, n_2)$. (That is if $\dim W = (k_1, k_2) = \dim W_1 \cup \dim W_2$ that is $(k_1, k_2) + (n_1 - k_1, n_2 - k_2) = (n_1, n_2))$.*

Now only in case of strong neutrosophic bivector spaces of finite $(n_1, n_2)$ bidimension over the neutrosophic bifield $F = F_1 \cup F_2$ we are in a position to define the strong neutrosophic bihyper subspaces of $V = V_1 \cup V_2$.

Suppose $V = V_1 \cup V_2$ be a strong neutrosophic bivector space over the neutrosophic bifield $F = F_1 \cup F_2$ of $(n_1, n_2)$ bidimension. $f = f^1 \cup f^2$ be a bilinear functional on V. The binull space of f or the null bisubspace of f denoted by $N_f = N_{f^1}^1 \cup N_{f^2}^2$.

The bidimension of $N_f = \dim N_{f^1}^1 \cup N_{f^2}^2$; but $\dim N_{f^i}^i = \dim V_i^{-1} = n_i - 1$ true for $i = 1, 2$. Thus bidimensin of $N_f = \dim (V_1 - 1) \cup \dim (V_2 - 1) = \dim N_{f^1}^1 \cup \dim N_{f^2}^2$.

We know in a vector space of dimension n a subspace of dimension $n - 1$ is called a hypersubspace likewise in a strong neutrosophic bivector space of bidimenion $(n_1, n_2)$ over the neutrosophic bifield $F = F_1 \cup F_2$ the bisubspace of bidimension $(n_1 - 1, n_2 - 1)$, we call that bisubspace to be a strong neutrosophic bihypersubspace of V.

Thus $N_f = N_{f^1}^1 \cup N_{f^2}^2$ is a strong neutrosophic bihyper subspace of the strong neutrosophic bivector space of $V = V_1 \cup V_2$.

One can prove as in case of bivector space in case of strong neutrosophic bivector spaces of type II if $W_1 = W_1^1 \cup W_1^2$ and $W_2 = W_2^1 \cup W_2^2$ be strong neutrosophic bivector subspaces of a strong neutrosophic bivector space $V = V_1 \cup V_2$ over the neutrosophic bifield $F = F_1 \cup F_2$ of bidimension $(n_1, n_2)$ then $W_1 = W_2$ if and only if $W_1^\circ = W_2^\circ$ that is if and only if $(W_1^i)^\circ = (W_2^i)^\circ$ for $i = 1, 2$.

Further as in case of bivector spaces of type II we can define the concept of dual space of a dual space $V^{**} = V_1^{**} \cup V_2^{**} = V_1 \cup V_2$.

Let $V = V_1 \cup V_2$ be a strong neutrosophic bivector space over the neutrosophic bifield $F = F_1 \cup F_2$ (i.e., each $V_i$ is a strong neutrosophic vector space over $F_i$, $i = 1, 2$) Let $V^* = V_1^* \cup V_2^*$ be the bivector space which is the bidual of V over the same bifield $F = F_1 \cup F_2$. The bidual of the bidual space $V^*$, i.e., $V^{**}$ in terms of the bibasis and bidual basis is given in the following:

Let $\alpha = \alpha^1 \cup \alpha^2$ be a strong neutrosophic bivector space V $= V_1 \cup V_2$ then $\alpha$ induces bilinear function $L_\alpha = L^1_{\alpha^1} \cup L^2_{\alpha^2}$ defined by

$$\begin{aligned} L_\alpha(f) &= L^1_{\alpha^1}(f^1) \cup L^2_{\alpha^2}(f^2) \\ &= f(\alpha) \\ &= f^1(\alpha^1) \cup f^2(\alpha^2) \end{aligned}$$

$f \in V_1^* \cup V_2^* = V^*$; $f^i \in V_i^*$; $i = 1, 2$. The fact each $L^i_{\alpha^i}$ is linear is just a reformation of the definition of the linear operators on $V_i^*$ for each $i = 1, 2$. The fact that each $L^i_{\alpha^i}$ is linear is just a reformation of the definition of linear operators in $V_i^*$; $i = 1, 2$.

$$\begin{aligned} L_\alpha(cf + g) &= L^1_{\alpha^1}(c_1 f^1 + g_1) + L^2_{\alpha^2}(c_2 f^2 + g_2) \\ &= (c_1 f^1 + g_1)(\alpha^1) \cup (c_2 f^2 + g_2)(\alpha^2) \\ &= c_1 f^1(\alpha^1) + g_1(\alpha^1) \cup c_2 f^2(\alpha^2) + g_2(\alpha^2) \\ &= c_1 L_\alpha(f) + L_\alpha(g) \end{aligned}$$

where $f = f^1 \cup f^2$ and $g = g_1 \cup g_2$. If $V = V_1 \cup V_2$ is a strong neutrosophic finite $(n_1, n_2)$ bidimensional and $\alpha \neq 0 = \alpha^1 \cup \alpha^2$ then $L_\alpha = L^2_{\alpha^1} \cup L^2_{\alpha^2} \neq 0 \cup 0$, in other words there exists a bilinear function $f = f^1 \cup f^2$ such that $f(\alpha) \neq 0$; i.e., $f(\alpha) = f^1(\alpha^1) \cup f^2(\alpha^2)$ for each $f^i(\alpha^i) \neq 0$; $i = 1, 2$. Further the bimapping $\alpha = \alpha^1 \cup \alpha^2 \to L_\alpha = L^2_{\alpha^1} \cup L^2_{\alpha^2}$ is a biisomorphism of $V = V_1 \cup V_2$ on to $V^{**} = V_1^{**} \cup V_2^{**}$.

Several properties in this direction can be analysed by any interested reader.

It can be easily proved as in case of bivector spaces of type II;

"If $S = S_1 \cup S_2$ is any biset of a $(n_1, n_2)$ bifinite dimensional strong neutrosophic bivector space $V = V_1 \cup V_2$ then $(S°)° = ((S_1°) \cup (S_2°))° = (S_1°)° \cup (S_2°)°$ is the bisubspace spanned by $S = S_1 \cup S_2$.

Thus if $V = V_1 \cup V_2$ is a strong neutrosophic bivector space of type II defined over the neutrosophic bifield $F = F_1 \cup F_2$. We define the bihypersubspace or hyperbispace of $V = V_1 \cup V_2$. Assume $V = (V_1 \cup V_2)$ is a $(n_1, n_2)$ dimension over $F = F_1 \cup F_2$. If $N = N_1 \cup N_2$ is a bihyperspace of V that is $N = N_1 \cup N_2$ is of $(n_1 - 1, n_2 - 1)$ bidimensional over $F = F_1 \cup F_2$ then we can define $N = N_1 \cup N_2$ to be a hyper space of V if

(1) $N = N_1 \cup N_2$ is a proper strong neutrosophic bivector subspace of V.
(2) If W is a strong neutrosophic bisubspace of V which contains N then either $W = N$ or $W = V$.

Condition (1) and (2) together say that $N = N_1 \cup N_2$ is a proper strong neutrosophic bisubspace and there is no larger proper strong neutrosophic bisubspace; in short $N = N_1 \cup N_2$ is a maximal proper strong neutrosophic bisubspace of V. Thus if $V = V_1 \cup V_2$ is a strong neutrosophic bivctor space over the neutrosophic bifield $F = F_1 \cup F_2$, a bihyper space in $V = V_1 \cup V_2$ is a maximal proper strong neutrosophic bisubspace of $V = V_1 \cup V_2$ over the neutrosophic bifield $F = F_1 \cup F_2$.

The following property about bihyperspace of V can be easily proved.

If $f = f^1 \cup f^2$ is a nonzero bilinear functional on the strong neutrosophic bivector space $V = V_1 \cup V_2$ over the neutrosophic bifield $F = F_1 \cup F_2$ of $(n_1, n_2)$ finite bidimension over F then the bihyperspace of V is the binull space of a non zero bilinear functional on V. It need not be unique.

We just give another interesting property about strong neutrosophic bivector spaces over a bifield of type II.

Let $V = V_1 \cup V_2$ and $W = W_1 \cup W_2$ be two strong neutrosophic bivector spaces over the same neutrosophic bifield $F = F_1 \cup F_2$. For each bilinear transformation $T = T_1 \cup T_2$ from $V = V_1 \cup V_2$ into $W = W_1 \cup W_2$ there is a unique bilinear transformation $T^t = T_1^t \cup T_2^t$ from $W^* = W_1^* \cup W_2^*$ into $V^* = V_1^* \cup V_2^*$ such that

$$\left(T_g^t\right)\alpha = \left(T_{1g_1}^t \cup T_{2g_2}^t\right)(\alpha^1 \cup \alpha^2)$$
$$= g_1(T_1\alpha^1) \cup g_2(T_2\alpha^2)$$
$$= g(T\alpha)$$

for every $g = g_1 \cup g_2 \in W^* = W_1^* \cup W_2^*$ and $\alpha = \alpha^1 \cup \alpha^2$ in $V = V_1 \cup V_2$.

We call $T^t = T_1^t \cup T_2^t$ as a bitranspose of $T = T_1 \cup T_2$. This bitransformaiton $T^t$ is also called as the biadjoint of T.

We now prove an important property about $T^t$.

**THEOREM 2.3.15:** *Let $V = V^1 \cup V^2$ and $W = W^1 \cup W^2$ be any two strong neutrosophic bivector spaces over the neutrosophic bifield $F = F^1 \cup F^2$ and let $T = T^1 \cup T^2$ be a strong neutrosophic bilinear transformation from $V = V^1 \cup V^2$ into $W = W^1 \cup W^2$. The binull space of $T^t = T_1^t \cup T_2^t$ is the biannihilator of the birange of $T = T^1 \cup T^2$. If V and W are finite bidimensional then*
i. *birank $(T^t)$ = birank T.*
ii. *The birange of $T^t = T_1^t \cup T_2^t$ is the annihilator of the binull space of $T = T_1 \cup T_2$.*

*Proof:* Let $g = g^1 \cup g^2$ be in $W^* = W_1^* \cup W_2^*$ the dual space of the strong neutrosophic bivector space $W = W_1 \cup W_2$. By definition we have $(T_g^t)\alpha = g(T\alpha)$ where for each $\alpha = \alpha^1 \cup \alpha^2 \in V_1 \cup V_2$. $T = T_1 \cup T_2 : V_1 \cup V_2 \to W_1 \cup W_2$. The statement that $g = g^1 \cup g^2$ is in the binull space of $T^t = T_1^t \cup T_2^t$ means that $g(T\alpha) = 0$; i.e., $g^1 T_1 \alpha^1 \cup g^2 T_2 \alpha^2 = 0 \cup 0$ for every $\alpha = \alpha^1 \cup \alpha^2 \in V = V_1 \cup V_2$.

Thus the binull space $T^t = T_1^t \cup T_2^t$ is precisely the biannihilator of the birange of $T = T_1 \cup T_2$. Suppose $V = V_1 \cup V_2$ and $W = W_1 \cup W_2$ are finite bidimensional, we say bidimension $V = (n_1, n_2)$ and bidimension $W = (m_1, m_2)$.

*Proof of (i):* Let $r = (r_1, r_2)$ be the birank of $T = T_1 \cup T_2$, i.e., the bidimension of the birange of T is $(r_1, r_2)$.

By earlier results the biannihilator of the birange of $T = T_1 \cup T_2$ has bidimension $(m_1 - r_1, m_2 - r_2)$. By the first statement of the theorem the binullity of $T^t = T_1^t \cup T_2^t$ must be $(m_1 - r_1, m_2 - r_2)$. Since $T^t = T_1^t \cup T_2^t$ is bilinear transformation on an $(m_1, m_2)$ bidimensional bispace the birank of $T^t = T_1^t \cup T_2^t$ is $(m_1 - (m_1 - r_1), m_2 - (m_2 - r_2))$ and so T and $T^t$ have the same birank.

*Proof for (ii):* Let $N = N_1 \cup N_2$ be the binull space of $T = T_1 \cup T_2$. Every bifunction in the birange of $T^t = T_1^t \cup T_2^t$ is in the biannihilator of $N = N_1 \cup N_2$, for suppose $f = T^t g$; i.e., $f^1 \cup f^2 = T_1^t g^1 \cup T_2^t g^2$ for some $g = g^1 \cup g^2$ in $W^* = W_1^* \cup W_2^*$ then if $\alpha = \alpha^1 \cup \alpha^2$ is in $N = N_1 \cup N_2$ ;

$$\begin{aligned} f(\alpha) &= f^1(\alpha^1) \cup f^2(\alpha^2) \\ &= (T_g^t)\alpha = (T_1^t g^1)\alpha^1 \cup (T_2^t g^2)\alpha^2 \\ &= g(T\alpha) \\ &= g^1(T_1 \alpha^1) \cup g_2(T_2 \alpha^2) \\ &= g^1(0) \cup g^2(0) \\ &= 0 \cup 0. \end{aligned}$$

Now the birange of $T^t = T_1^t \cup T_2^t$ is a bisubspace of the space $N^o = N_1^o \cup N_2^o$ and

$$\begin{aligned} \dim N^o &= (n_1 - \dim N_1) \cup (n_2 - \dim N_2) \\ &= \text{birank } T \\ &= \text{birank } T^t \end{aligned}$$

so that birange of $T^t$ must exactly be $N^o$.

**THEOREM 2.3.16:** *Let $V = V_1 \cup V_2$ and $W = W_1 \cup W_2$ be two $(n_1, n_2)$ and $(m_1, m_2)$ dimensional bivector spaces over the neutrosophic bifield $F = F_1 \cup F_2$. Let B be a bibasis of V and B\* the bidual basis of V\*. Let C be a bibasis of W with dual bibasis C\*. Let $T = T_1 \cup T_2$ be a bilinear transformation from V into W; let A be the neutrosophic bimatrix of $T = T_1 \cup T_2$ relative to B and C and let B be a neutrosophic bimatrix of $T^t$ relative to B\*, C\*. Then $B_{ij}^k = A_{ij}^k$ for $k = 1, 2$. That is*

$$A_{ij}^1 \cup A_{ij}^2 = B_{ij}^1 \cup B_{ij}^2.$$

*Proof:* Given $V = V_1 \cup V_2$ and $W = W_1 \cup W_2$ are strong neutrosophic bivector spaces over the neutrosophic bifield $F = F_1 \cup F_2$. Given $V = V_1 \cup V_2$ is $(n_1, n_2)$ bidimension and $W = W_1 \cup W_2$ is of $(m_1, m_2)$ bidimension over the bifield $F = F_1 \cup F_2$. Let

$$B = \{\alpha_1^1, \alpha_2^1, \ldots, \alpha_{n_1}^1\} \cup \{\alpha_1^2, \alpha_2^2, \ldots, \alpha_{n_2}^2\}$$

be a bibasis of $V = V_1 \cup V_2$ and the dual bibasis of B,

$$B^* = B_1^* \cup B_2^* = \{f_1^1, f_2^1, \ldots, f_{n_1}^1\} \cup \{f_1^2, f_2^2, \ldots, f_{n_2}^2\}.$$

Let

$$C = C_1 \cup C_2 = \{\beta_1^1, \beta_2^1, \ldots, \beta_{m_1}^1\} \cup \{\beta_1^2, \beta_2^2, \ldots, \beta_{m_2}^2\}$$

be a bibasis of $W = W_1 \cup W_2$.
The dual bibasis of C,

$$C^* = C_1^* \cup C_2^* = \{g_1^1, g_2^1, \ldots, g_{m_1}^1\} \cup \{g_1^2, g_2^2, \ldots, g_{m_2}^2\}.$$

Now by definition for $\alpha = \alpha^1 \cup \alpha^2$;

$$T_k \alpha_j^k = \sum_{i=1}^{m_k} A_{ij}^k \beta_i^k \; ; j = 1, 2, \ldots, n_k; k = 1, 2.$$

$$T_k^t g_j^k = \sum_{i=1}^{n_k} B_{ij}^k f_i^k \; ; j = 1, 2, \ldots, m_k \text{ and } k = 1, 2.$$

Further

$$\begin{aligned}
(T_k^t g_j^k)(\alpha_i^k) &= g_j^k (T_k^t \alpha_i^k) \\
&= g_j^k \left( \sum_{p=1}^{m_k} A_{pi}^k \beta_p^k \right) \\
&= \sum_{p=1}^{m_k} A_{pi}^k g_j^k (\beta_p^k) \\
&= \sum_{p=1}^{m_k} A_{pi} \delta_{jp} \\
&= A_{ji}^k.
\end{aligned}$$

For any bilinear functional $f = f^1 \cup f^2$ on V,

$$f^k = \sum_{i=1}^{m_k} f^k(\alpha_i^k) f_i^k \; ; k = 1, 2.$$

If we apply this formula to the functional $f^k = T_k^t g_j^k$ and use the fact $(T_k^t g_j^k) \alpha_i^k = A_{ji}^k$, we have

$$(T_k^t g_j^k) = \sum_{i=1}^{n_k} A_{ji}^k f_i^k$$

from which it follows $B_{ij}^k = A_{ij}^k$; true for $k = 1, 2$. That is

$$B_{ij}^1 \cup B_{ij}^2 = A_{ij}^1 \cup A_{ij}^2.$$

If $A = A^1 \cup A^2$ is a $(m_1 \times n_1, m_2 \times n_2)$ neutrosophic bimatrix over the neutrosophic bifield $F = F_1 \cup F_2$ then the bitranspose of A is the $(n_1 \times m_1, n_2 \times m_2)$ neutrosophic bimatrix $A^t$ defined by $\left(A_{ij}^1\right)^t \cup \left(A_{ij}^2\right)^t = A_{ij}^1 \cup A_{ij}^2.$

 We leave it as an exercise for the reader to prove the birow rank of A is equal to the bicolumn rank of A, that is for each neutrosophic matrix $A^i$ we have the column rank of $A^i$ to be equal to the row rank of $A^i$; $i = 1, 2$.

We see all these results holds good for strong neutrosophic bilinear algebras defined over the neutrosophic bifield $F_1 \cup F_2 = F$ with appropriate modification if necessary.

Now we proceed onto define the notion of neutrosophic bipolynomial over a neutrosophic bifield $F = F_1 \cup F_2$.

**DEFINITION 2.3.34:** *Let $F[x] = F_1[x] \cup F_2[x]$ be such that each $F_i[x]$ is a polynomial over $F_i$, $F_i$ a neutrosophic field; $i = 1, 2$ and $F_1 \neq F_2$ i.e., $F = F_1 \cup F_2$ is a neutrosophic bifield. We call $F[x]$ the neutrosophic bipolynomial over the neutrosophic bifield $F = F_1 \cup F_2$. Any element $p(x) \in F[x]$ will be the form $p(x) = p_1(x) \cup p_2(x)$ where $p_i(x)$ is a neutrosophic polynomial in $F_i[x]$; i.e., $p_i(x)$ is a neutrosophic polynomial in the variable x with coefficients from the neutrosophic field $F_i$, $i = 1, 2$. The bidegree of $p(x)$ is a pair given by $(n_1, n_2)$ where $n_i$ is the degree of the polynomial $p_i(x)$; $i = 1, 2$.*

We will illustrate this situation by some simple examples.

*Example 2.3.73:* Let $F = F_1 \cup F_2 = Z_3 I \cup QI$ be a neutrosophic bifield. $F[x] = F_1[x] \cup F_2[x] = Z_3I[x] \cup QI[x] = $ {all polynomials in the variable x with coefficients from the neutrosophic field $Z_3I$} $\cup$ {all polynomials in the variable x with coefficients from the neutrosophic field $QI$} is a bipolynomial strong neutrosophic bivector space over the neutrosophic bifield $F = F_1 \cup F_2 = Z_3I \cup QI$.
Let $p(x) = 2I + Ix + 2Ix^3 + Ix^7 \cup 3I + 7Ix + 270I\ x^7 - 5762I\ x^9 + 3006I\ x^{29}$; $p(x) \in F[x]$.

*Example 2.3.74:* Let $F = F_1 \cup F_2 = N(Q) \cup N(Z_{11})$ be a neutrosophic bifield. $F[x] = F_1[x] \cup F_2[x] = N(Q)[x] \cup N(Z_{11})[x] = $ {all polynomials in the variable x with coefficients from the neutrosophic field $N(Q)$} $\cup$ {all polynomials in the variable x with coefficients from the neutrosophic field $N(Z_{11})$} is a bipolynomial strong neutrosophic bivector space over the neutrosophic bifield $F = N(Q) \cup N(Z_{11})$.

Take $p(x) = p_1(x) \cup p_2(x) = 3 + 17x^2 - 245 x^5 + 346 x^7 - 93/2 x^8 + 47I x^9 - 5009 x^{11} \cup 3I + 8 + 4Ix + 5x^2 + 10I x^7 + 2x^{11} + 9x^{20} \in F[x]$ is a bipolynomials with coefficients from the bifield $N(Q) \cup N(Z_{11})$

**DEFINITION 2.3.35:** *Let $F[x] = F_1[x] \cup F_2[x]$ be a bipolynomial over the neutrosophic bifield $F = F_1 \cup F_2$. $F[x]$ is a strong neutrosophic bilinear algebra over the bifield F. Infact $F[x]$ is a bicommutative neutrosophic linear bialgebra over the bifield F. $F[x]$ the strong neutrosophic bilinear algebra may or may not have the biidentity $I_2 = 1 \cup I$ or $I \cup I$ or $I \cup 1$ or $1 \cup 1$ depending on the neutrosophic bifield $F = F_1 \cup F_2$. We call a neutrosophic bipolynomial $p(x) = p_1(x) \cup p_2(x)$ to be bimonic polynomial if each $p_i(x)$ is a monic neutrosophic polynomial in x for i = 1, 2. We will call a neutrosophic bipolynomial to be a neutrosophic monic bipolynomial if for each $p_i(x) \in F_i[x]$ the coefficient associated with the highest degree is I for i = 1, 2.*

We will first illustrate these situations before we proceed on to prove further results.

*Example 2.3.75:* Let $F = F_1 \cup F_2 = Z_7I \cup Z_{13}I$ be a neutrosophic bifield. $F[x] = F_1[x] \cup F_2[x] = Z_7I[x] \cup Z_{13}[x] = $ {all polynomials in the variable x with coefficients from the neutrosophic field $Z_7I$} $\cup$ {all polynomials in the variable x with coefficients from the neutrosophic field $Z_{13}I$} is a bipolynomial strong neutrosophic bilinear algebra over the bifield $F = F_1 \cup F_2 = Z_7I \cup Z_{13}I$.
   It is easily verified that $F[x]$ has no monic bipolynomial however $F[x]$ has neutrosophic monic bipolynomials. For take
   $p(x) = Ix^{29} + 2I x^8 + 4Ix^3 + 5I \cup Ix^{47} + 12Ix^{25} + 10I x^{12} + 7I x^4 + 5Ix + 3I \in F_1[x] \cup F_2[x]$. Clearly $p(x)$ is a neutrosophic monic bipolynomial in $F(x)$.

*Example 2.3.76:* Let $F = F_1 \cup F_2 = N(Z_{23}) \cup N(Z_{47})$ be a neutrosophic bifield. $F[x] = F_1[x] \cup F_2[x] = N(Z_{23})[x] \cup N(Z_{47})[x] = $ {all polynomials in the variable x with coefficients from the neutrosophic field $N(Z_{23})$} $\cup$ {all polynomials in the variable x with coefficients from the neutrosophic field $N(Z_{47})$}. $F[x]$ is a bipolynomial strong neutrosophic bilinear algebra over the neutrosophic bifield $F = N(Z_{23}) \cup N(Z_{47})$. Take $p(x) = p_1(x) \cup p_2(x) \in N(Z_{23})[x] \cup N(Z_{47})[x]$; $p(x) = \{x^{48} + 3I\, x^{20} + 15x^{12} + 4Ix^{26} + 13I\, x^7 + 5x^3 + 20I + 4\} \cup \{x^{104} + 46I\, x^{100} + 45Ix^{79} + 27x^{68} + 40x^{27} + 37x^5 + Ix^2 + 7I + 4\}$ is a monic bipolynomial in $F[x]$.

It is interesting to note that in case of these bipolynomial strong neutrosophic bilinear algebra the conditions under which the bilinear algebra will have monic bipolynomials and when it will never have monic bipolynomials.

**THEOREM 2.3.17:** *Let $F = F_1 \cup F_2$ be a neutrosophic bifield $F[x] = F_1[x] \cup F_2[x]$ be a bipolynomial strong neutrosophic bilinear algebra over the neutrosophic bifield $F = F_1 \cup F_2$.*
*If both $F_1$ and $F_2$ are of the form $K_1 I$ and $K_2 I$ where $K_1$ and $K_2$ are real fields then the bipolynomials strong neutrosophic bilinear algebra $F[x]$ will not contain monic bipolynomial.*

*Proof:* Given $F = F_1 \cup F_2$ is a neutrosophic bifield where $F_1 = K_1 I$ and $F_2 = K_2 I$ with $K_1$ and $K_2$ real fields.
$F[x] = F_1[x] \cup F_2[x] = K_1 I[x] \cup K_2 I[x]$ is the bipolynomial strong neutrosophic bilinear algebra over the neutrosophic bifield $F = F_1 \cup F_2$. We see clearly $F = F_1 \cup F_2 = K_1 I \cup K_2 I$ does not contain any real element; every bipair in $K_1 I \cup K_2 I$ is neutrosophic. Hence $1 \notin K_i I$ for $i = 1, 2$. Thus no bipolynomials in the variable x has real coefficients i.e., no polynomial $p(x)$ in the variable x in $K_i I[x]$ has real coefficients; for $i = 1, 2$.
Thus $F[x] = K_1 I[x] \cup K_2 I[x]$ has no monic bipolynomial.
Hence the claim.

**THEOREM 2.3.18:** *Let $F = F_1 \cup F_2$ be a neutrosophic bifield. $F[x] = F_1[x] \cup F_2[x]$ be a bipolynomials strong neutrosophic bilinear algebra over the neutrosophic bifield $F = F_1 \cup F_2$. If each $F_i$ is of the form $N(K_i)$ where $K_i$ is a real field for $i = 1, 2$ then the bipolynomial strong neutrosophic linear bialgebra has monic bipolynomials.*

*Proof:* Given $F = F_1 \cup F_2 = N(K_1) \cup N(K_2)$ (where $K_1$ and $K_2$ are real fields) is a neutrosophic bifield. $F[x] = F_1[x] \cup F_2[x] = N(K_1)[x] \cup N(K_2)[x]$ is a bipolynomial strong neutrosophic bilinear algebra over the neutrosophic bifield $F = F_1 \cup F_2$.

Now clearly $K_i \subseteq N(K_i)$ for $i = 1, 2$; that is the neutrosophic field $N(K_i)$ contains the real field $K_i$ as a proper subfield, true for $i = 1, 2$. Thus $1 \in N(K_i)$ for $i = 1, 2$. Now we can take $p(x) = p_1(x) \cup p_2(x)$ where both $p_1(x)$ is a monic polynomial of the form say $x^{219} + 7I x^{200} + 14I x^{14} + 27x^{10} + 205$ in $F_1[x]$ and $p_2(x)$ to be a monic polynomial of the form $x^3 + 7x + 21I$ in $F_2[x]$ we see $p(x)$ is a monic bipolynomial in $F[x]$.

Hence the claim

Thus we see when both the neutrosophic fields $F_i$; $i = 1, 2$, are not pure neutrosophic fields then certainly the bipolynomial strong neutrosophic bilinear algebra has monic bipolynomials.

Further even if one of the neutrosophic field $F_i$ is a pure neutrosophic field that is $F_i = K_i I$ where $K_i$ is a real field $i = 1, 2$; then $F[x] = F_1[x] \cup F_2[x]$ has no bipolynomial which is a monic bipolynomial.

The reader is expected to prove the following results.

**THEOREM 2.3.19:** *Let $F[x] = F_1[x] \cup F_2[x]$ be a strong neutrosophic bilinear algebra of bipolynomials over the neutrosophic bifield $F = F_1 \cup F_2$ then*
  i. *If $f(x) = f_1(x) \cup f_2(x)$ and $g(x) = g_1(x) \cup g_2(x)$ are two non zero bipolynomials in $F[x]$, the bipolynomial $f(x) g(x) = f_1(x) g_1(x) \cup f_2(x) g_2(x)$ is a non zero bipolynomial in $F[x]$.*
  ii. *The bidegree of $(f(x) g(x)) =$ bidegree of $f(x)$ + bidegree of $g(x)$ where bidegree of $f = (n_1, n_2)$ and bidegree of $g = (m_1, m_2)$.*
  iii. *$f(x) g(x)$ is monic bipolynomial if both $f(x)$ and $g(x)$ are monic polynomials and $F = F_1 \cup F_2$ is neutrosophic bifield of the form $F = N(K_1) \cup N(K_2)$ where $K_1$ and $K_2$ are real fields.*
  iv. *$f(x) g(x)$ is a monic neutrosophic bipolynomial if both $f(x)$ and $g(x)$ are monic neutrosophic bipolynomials. ($F = K_1 I \cup K_2 I$).*

v.  If $f + g$  $= f_1 \cup f_2 + g_1 \cup g_2$
    $= (f_1 + g_1) \cup (f_2 + g_2)$
    $\neq 0 \cup 0$
    $= \max(\text{bideg } f, \text{bideg } g)$.

vi. If $f, g, h$ are bipolynomials over the neutrosophic bifield $F = F_1 \cup F_2$. $f(x) = f_1(x) \cup f_2(x)$, $g(x) = g_1(x) \cup g_2(x)$ and $h(x) = h_1(x) \cup h_2(x)$; $g(x) \neq 0 \cup 0$, $f(x) \neq 0 \cup 0$ and $h(x) \neq 0 \cup 0$ and if $fg = fh$ then $g = h$.

As in case of polynomials we can derive most of the results in case of neutrosophic bipolynomial.

Let $A = A_1 \cup A_2$ be a strong neutrosophic bilinear algebra with biidentity $1 = 1 \cup 1$ over the neutrosophic bifield $F = F_1 \cup F_2$ where we make the convention for any real $\alpha = \alpha_1 \cup \alpha_2$ ($\alpha_1$ and $\alpha_2$ are both real).

$$\alpha^0 = \alpha_1^0 \cup \alpha_2^0 = 1 \cup 1 = 1_2.$$

We cannot derive the properties enjoyed by usual bipolynomials.

Thus to get the analogue of the Lagrange biinterpolation formula in case of bipolynomial strong neutrosophic spaces we have to make more assumptions. We can derive several results analogous to bipolynomial bilinear algebra.

Suppose $f = f_1 \cup f_2$ and $d = d_1 \cup d_2$ be any two non zero neutrosophic bipolynomials over the neutrosophic bifield $F = F_1 \cup F_2$ such that bideg $d \leq$ bideg $f$ (i.e., bideg $d = (n_1, n_2)$ and bideg $f = (m_1, m_2)$ and $n_i \leq m_i$ for $i = 1, 2$ (then we say bideg $d \leq$ bideg $f$) then there exists a bipolynomial $g = g_1 \cup g_2$ in $F[x] = F_1[x] \cup F_2[x]$ such that either $f - dg = 0$ that is $(f_1 - d_1 g_1) \cup (f_2 - d_2 g_2) = 0 \cup 0$ or bideg $(f - dg) <$ bideg $f$.

We have also the following interesting result in case of neutrosophic bipolynomials.

**THEOREM 2.3.20:** *Let $f = f_1 \cup f_2$ and $d = d_1 \cup d_2$ be bipolynomials over the neutrosophic bifield $F = F_1 \cup F_2$ and $d = d_1 \cup d_2$ is different from $0 \cup 0$ then there exists bipolynomials $q = q_1 \cup q_2$ and $r = r_1 \cup r_2$ in $F[x] = F_1[x] \cup F_2[x]$ such that $f = dq + r$; i.e., $f = f_1 \cup f_2 = (d_1 q_1 + r_1) \cup (d_2 q_2 + r_2)$.*

*(2) Either $r = r_1 \cup r_2 = (0 \cup 0)$ or bideg $r <$ bideg d. The bipolynomials $q = q_1 \cup q_2$ and $r = r_1 \cup r_2$ satisfying the conditions (1) and (2) are unique.*

The proof is direct hence left as an exercise for the reader.

**DEFINITION 2.3.36:** *Let $d = d_1 \cup d_2$ be a non zero bipolynomial over the neutrosophic bifield $F = F_1 \cup F_2$. If $f = f_1 \cup f_2$ is in $F[x] = F_1[x] \cup F_2[x]$, the proceeding theorem show there exists atmost one bipolynomial $q = q_1 \cup q_2$ in $F[x]$ such that $f = dq$ i.e., $f = f_1 \cup f_2 = d_1q_1 \cup d_2q_2$. If such a $q = q_1 \cup q_2$ exists we say that $d = d_1 \cup d_2$ bidivides $f = f_1 \cup f_2$ and f is bidivisible by $d = d_1 \cup d_2$ and $f = f_1 \cup f_2$ is a bimultiple of $d = d_1 \cup d_2$ and we call $q = q_1 \cup q_2$ to be the biquotient of f and $d = d_1 \cup d_2$ and write $q = f / d$ that is $q = q_1 \cup q_2 = f_1 / d_1 \cup f_2 / d_2$.*

The following result is direct. If $f = f_1 \cup f_2$ is a bipolynomial over the neutrosophic bifield $F = F_1 \cup F_2$ and $c = c_1 \cup c_2$ be an element of F. f is bidivisible by $x - c = (x - c_1) \cup (x - c_2)$ if and only if $f(c) = f_1(c_1) \cup f_2(c_2) = 0 \cup 0$.

We can prove the fundamental theorem of algebra namely that every polynomial of degree n has atmost n roots can be proved in case of neutrosophic bipolynomials. A bipolynomial $f = f_1 \cup f_2$ of degree $(n_1, n_2)$ over a neutrosophic bifield $F = F_1 \cup F_2$ has atmost $(n_1, n_2)$ biroots in $F = F_1 \cup F_2$. Now we will prove Taylors formula for bipolynomials over the neutrosophic bifield $F = F_1 \cup F_2$.

**THEOREM 2.3.21:** *Let $F = F_1 \cup F_2$ be a neutrosophic bifield of bicharacteristic (0, 0). $c = c_1 \cup c_2$ be an element in $F = F_1 \cup F_2$ and $(n_1, n_2)$ be a pair of positive integers. If $f = f_1 \cup f_2$ is a neutrosophic bipolynomial over the bifield $F = F_1 \cup F_2$ with bideg $f \leq (n_1, n_2)$ then*

$$f = \sum_{k_1=0}^{n_1} \frac{D^{k_1} f_1 c_1 (x - c_1)^{k_1}}{\lfloor k_1} \cup \sum_{k_2=0}^{n_2} \frac{D^{k_2} f_2 c_2 (x - c_2)^{k_2}}{\lfloor k_2}.$$

*Proof:* We know Taylors theorem is a consequence of the binomial theorem and the linearity of the operators $D^1$, $D^2$, ..., $D^n$. We know the binomial theorem

$$(a+b)^m = \sum_{k=0}^{m} \begin{bmatrix} m \\ k \end{bmatrix} a^{m-k} b^k$$

where $\begin{pmatrix} m \\ k \end{pmatrix} = \frac{m!}{k!(m-k)!} = \frac{m(m-1)...(m-k+1)}{1.2...k = \lfloor k}$ is the familiar binomial coefficient giving the number of combinations of m objects taken k at a time.

Now we apply the binomial theorem to the pair of neutrosophic polynomials

$$\begin{aligned} x^{m_1} \cup x^{m_2} &= (c_1 + (x-c_1))^{m_1} \cup (c_2 + (x-c_2))^{m_2} \\ &= \sum_{0}^{m_1} \begin{pmatrix} m_1 \\ c_1 \end{pmatrix} c_1^{m_1-k_1}(x-c_1)^{k_1} \cup \sum_{0}^{m_2} \begin{pmatrix} m_2 \\ c_2 \end{pmatrix} c_2^{m_2-k_2}(x-c_2)^{k_2} \\ &= \left\{ c_1^{m_1} + m_1 c_1^{m_1-1}(x-c_1) + ... + (x-c_1)^{m_1} \right\} \cup \\ &\quad \left\{ c_2^{m_2} + m_2 c_2^{m_2-1}(x-c_2) + ... + (x-c_2)^{m_2} \right\} \end{aligned}$$

and this is the statement of Taylor's biformula for the case
$$f = x^{m_1} \cup x^{m_2}.$$
If
$$f = \sum_{m_1=0}^{n_1} a_{m_1}^1 x^{m_1} \cup \sum_{m_2=0}^{n_2} a_{m_2}^2 x^{m_2},$$

$$D_f^k(c) = \sum_{m_1=0}^{n_1} a_{m_1}^1 D^{k_1} x^{m_1}(c_1) \cup \sum_{m_2=0}^{n_2} a_{m_2}^2 D^{k_2} x^{m_2}(c_2)$$

and

$$\sum_{k_1=0}^{m_1} \frac{D^{k_1} f_1(c_1)(x-c_1)^{k_1}}{\lfloor k_1} \cup \sum_{k_2=0}^{m_2} \frac{D^{k_2} f_2(c_2)(x-c_2)^{k_2}}{\lfloor k_2}$$

$$= \sum_{k_1} \sum_{m_1} a^1_{m_1} \frac{D^{k_1} x^{m_1}(c_1)(x-c_1)^{k_1}}{\lfloor k_1} \cup$$

$$\sum_{k_2} \sum_{m_2} a^2_{m_2} \frac{D^{k_2} x^{m_2}(c_2)(x-c_2)^{k_2}}{\lfloor k_2}$$

$$= \sum_{m_1} a^1_{m_1} \sum_{k_1} \frac{D^{k_1} x^{m_1}(c_1)(x-c_1)^{k_1}}{\lfloor k_1} \cup$$

$$\sum_{m_2} a^2_{m_2} \sum_{k_2} \frac{D^{k_2} x^{m_2}(c_2)(x-c_2)^{k_2}}{\lfloor k_2}.$$

If $c = c_1 \cup c_2$ is a biroot of the neutrosophic bipolynomial $f = f_1 \cup f_2$ with bimultiplicity $c = c_1 \cup c_2$ as a biroot of $f = f_1 \cup f_2$ is the largest bipositive integer $(r_1, r_2)$ such that $(x - c_1)^{r_1} \cup (x - c_2)^{r_2}$ bidivides $f = f_1 \cup f_2$.

Now we have still an interesting result on these neutrosophic bipolynomials and their bimultiplicity.

**THEOREM 2.3.22:** *If $F = F_1 \cup F_2$ is a neutrosophic bifield of $(0, 0)$ bicharacteristic (i.e., each $F_i$ is of characteristic zero for $i = 1, 2$) and $f = f_1 \cup f_2$ be a neutrosophic bipolynomial over the bifield $F = F_1 \cup F_2$ with bideg $f \leq (n_1, n_2)$. Then the biscalar $c = c_1 \cup c_2$ is a biroot of $f = f_1 \cup f_2$ of multiplicity $(r_1, r_2)$ if and only if $(D^{k_1} f_1)(c_1) \cup (D^{k_2} f_2)(c_2) = 0 \cup 0;\ 0 \leq k_i \leq r_i - 1;\ i = 1, 2$. $D^{r_i} f_i(c_i) \neq 0$ for every $i = 1, 2$.*

*Proof:* Suppose that $(r_1, r_2)$ is the bimultiplicity of $c = c_1 \cup c_2$ as a biroot of $f = f_1 \cup f_2$.

Then there exists a neutrosophic bipolynomial $g = g_1 \cup g_2$ such that $f = (x - c_1)^{r_1} g_1 \cup (x - c_2)^{r_2} g_2$ and $g(c) = g_1(c_1) \cup g_2(c_2) \neq 0 \cup 0$.

For otherwise $f = f_1 \cup f_2$ would be bidivisible by $(x - c_1)^{r_1 + 1} \cup (x - c_2)^{r_2 + 1}$. By Taylors biformula applied to $g = g_1 \cup g_2$,

$$f = (x-c_1)^{r_1} \sum_{m_1=0}^{n_1-r_1} \frac{(D^{m_1}g_1)(c_1)(x-c_1)^{m_1}}{\lfloor m_1}\ \cup$$

$$(x-c_1)^{r_2} \sum_{m_2=0}^{n_2-r_2} \frac{(D^{m_2}g_2)(c_2)(x-c_2)^{m_2}}{\lfloor m_2}$$

$$= \sum_{m_1=0}^{n_1-r_1} \frac{D^{m_1}g_1(x-c_1)^{r_1+m_1}}{\lfloor m_1}\ \cup\ \sum_{m_2=0}^{n_2-r_2} \frac{(D^{m_2}g_2)(x-c_2)^{r_2+m_2}}{\lfloor m_2}.$$

Since there is only one way to write $f = f_1 \cup f_2$ (i.e., only one way to write each component $f_i$ of $f$; $i = 1, 2$) as a bilinear combination of bipowers of $(x-c_1)^{k_1} \cup (x-c_2)^{k_2}$; $0 \le k_i \le n_i$, $i = 1, 2$ it follows that

$$\frac{(D^{k_i}f_i)(c_i)}{\lfloor k_i} = \begin{cases} 0 & \text{if } 0 \le k_i \le r_i - 1 \\ \dfrac{D^{k_i-r_i}g_i(c_i)}{(k_i-r_i)!} & \text{if } r_i \le k_i \le n_i. \end{cases}$$

This is true for every $i$, $i = 1, 2$. Therefore $D^{k_i} f_i(c_i) = 0$ for $0 \le k_i \le r_i - 1$; $i = 1, 2$ and $D^{r_i} f_i(c_i) \ne g_i(c_i) \ne 0$ for every $i$, $i = 1, 2$.

Conversely if these conditions are satisfied, it follows at once from Taylor's biformula that there is a neutrosophic bipolynomial $g = g_1 \cup g_2$ such that $f = f_1 \cup f_2 = (x-c_1)^{r_1} g_1 \cup (x-c_2)^{r_2} g_2$ and $g_1(c_1) \cup g_2(c_2) = g(c) \ne 0 \cup 0$.

Now suppose that $(r_1, r_2)$ is not the largest positive biinteger pair such that $(x-c_1)^{r_1} \cup (x-c_2)^{r_2}$ bidivides $f_1 \cup f_2$; i.e., each $(x-c_i)^{r_i}$ divides $f_i$; $i = 1, 2$. then there is a bipolynomial $h = h_1 \cup h_2$ such that $f = (x-c_1)^{r_1+1} h_1 \cup (x-c_2)^{r_2+1} h_2$. But this implies $g = g_1 \cup g_2 = (x-c_1)h_1 \cup (x-c_2)h_2$; hence $g(c) = g_1(c_1) \cup g_2(c_2) = 0 \cup 0$ a contradiction; hence the claim.

Now we proceed onto define principal biideal generated by the neutrosophic bipolynomial $d = d_1 \cup d_2$.

**DEFINITION 2.3.37:** *Let $F = F_1 \cup F_2$ be a neutrosophic bifield. A biideal in $F[x] = F_1[x] \cup F_2[x]$ is a strong neutrosophic bisubspace $m = m_1 \cup m_2$ of $F[x] = F_1[x] \cup F_2[x]$ such that when $f = f_1 \cup f_2$ and $g = g_1 \cup g_2$ then $fg = f_1 g_1 \cup f_2 g_2$ belongs to $m = m_1 \cup m_2$; i.e., each $f_i g_i \in m_i$ whenever f is in $F[x]$ and $g \in m$ (i = 1, 2).*

If in particular the biideal $m = d F[x]$ for some bipolynomial $d = d_1 \cup d_2$ in $F[x] = F_1[x] \cup F_2[x]$, i.e., the biset of all bimultiple $df = d_1 f_1 \cup d_2 f_2$ of $d = d_1 \cup d_2$ by arbitrary $f = f_1 \cup f_2$ in $F[x] = F_1[x] \cup F_2[x]$ is a biideal; for m is non empty; m infact contains d. If $f, g \in F[x] = F_1[x] \cup F_2[x]$ and $c = c_1 \cup c_2$ is a biscalar then $c(df) - dg = (c_1 d_1 f_1 - d_1 g_1) \cup (c_2 d_2 f_2 - d_2 g_2) = d_1(c_1 f_1 - g_1) \cup d_2(c_2 f_2 - g_2)$ belongs to $m = m_1 \cup m_2$, that is $d_i(c_i f_i - g_i) \in m_i$; i = 1, 2; so that m is strong neutrosophic bivector subspace. Finally m contains

$$(df)g = d(fg)$$
$$= (d_1 f_1)g_1 \cup (d_2 f_2)g_2$$
$$= d_1(f_1 g_1) \cup d_2(f_2 g_2)$$

as well $m = m_1 \cup m_2$ is called the principal biideal generated by $d = d_1 \cup d_2$.

We will prove the following biprincipal ideal or principal biideal of $F[x] = F_1[x] \cup F_2[x]$.

**THEOREM 2.3.23:** *Let $F = F_1 \cup F_2$ be a bifield which is a neutrosophic bifield and $m = m_1 \cup m_2$ a non zero biideal in $F[x] = F_1[x] \cup F_2[x]$. Then there is a unique monic bipolynomial $d = d_1 \cup d_2$ in $F[x]$ where each $d_i$ is a monic polynomial in $F_i[x]$; i = 1, 2 such that m is the principal biideal generated by d.*

*Proof:* Given $F = F_1 \cup F_2$ is a neutrosophic bifield and $F[x] = F_1[x] \cup F_2[x]$ be the bipolynomial strong neutrosophic bivector space over the neutrosophic bifield $F = F_1 \cup F_2$. Let $m = m_1 \cup m_2$ be a non zero biideal of $F[x] = F_1[x] \cup F_2[x]$, we call a bipolynomial p(x) to be bimonic, i.e. if in $p(x) = p_1(x) \cup p_2(x)$ every $p_i(x)$ is a monic polynomial for i = 1, 2. Similarly we call a bipolynomial to be biminimal if in $p(x) = p_1(x) \cup p_2(x)$ each polynomial $p_i(x)$ is of minimal degree. Now $m = m_1 \cup m_2$ contains a non zero bipolynomial $p(x) = p_1(x) \cup p_2(x)$ where each $p_i(x) \neq 0$ for $i = 1, 2$. Among all the non zero bipolynomials in m there is a bipolynomial $d = d_1 \cup d_2$ of minimal bidegree. Without loss in generality we may assume that minimal bipolynomial is monic, i.e., d is monic. Suppose $f = f_1 \cup f_2$ is any bipolynomial in m then we know $f = dq + r = f_1 \cup f_2 = d_1q_1 + r_1 \cup d_2q_2 + r_2$ where $r = (r_1, r_2) = (0, 0)$ or bidegree $r <$ bidegree d; i.e., $f = f_1 \cup f_2 = (d_1q_1 + r_1) \cup (d_2q_2 + r_2)$. Since d is in m, $dq = d_1q_1 \cup d_2q_2 \in m$ and $f \in m$ so $f - dg = r = r_1 \cup r_2 \in m$. But since d is a bipolynomial in m of minmal bidegree we cannot have bidegree $r <$ bidegree d so $r = 0 \cup 0$.

Thus $m = dF[x] = d_1F_1[x] \cup d_2F_2[x]$. If g is any other bimonic polynomial such that $gF[x] = m = g_1F_1[x] \cup g_2F_2[x]$ then there exists non zero bipolynomial $p = p_1 \cup p_2$ and $q = q_1 \cup q_2$ such that $d = gp$ and $g = dq$. i.e., $d = d_1 \cup d_2 = g_1p_1 \cup g_2p_2$ and $g_1 \cup g_2 = d_1q_1 \cup d_2q_2$. Thus

$$d = dpq$$
$$= d_1p_1q_1 \cup d_2p_2q_2$$
$$= (d_1 \cup d_2)pq$$

and bidegree d = bidegree d + bidegree d + bidegree p + bideg q. Hence bidegree p = bidegree q = (0, 0) and as d and g are bimonic $p = q = 1$. Thus $d = g$.

Hence the claim.

If in the biideal m we have $f = pq + r$ where p, f ∈ m; i.e., $p = p_1 \cup p_2 \in m$ and $f = f_1 \cup f_2 \in m$; $f = f_1 \cup f_2 = (p_1q_1 + r_1) \cup (p_2q_2 + r_2)$ where the biremainder $r = r_1 \cup r_2 \in m$ and is different from $0 \cup 0$ and has smaller bidegree than p.

The interested reader is requested to prove the following results.

**COROLLARY 2.3.1:** *If $p^1$, $p^2$ are bipolynomials over a neutrosophic bifield $F = F_1 \cup F_2$ not all of which are zero; i.e., $0 \cup 0$, then there is a unique bimonic polynomial $d = d_1 \cup d_2$ in $F[x] = F_1[x] \cup F_2[x]$ such that*

i. $d = d_1 \cup d_2$ is the biideal generated by $p^1$, $p^2$ where $p^1 = p_1^1 \cup p_2^1$ and $p^2 = p_1^2 \cup p_2^2$.

ii. $d = d_1 \cup d_2$ bidivides each of the bipolynomials $p^i = p_1^i \cup p_2^i$ that is $d_j / p_j^i$, $j = 1, 2$ and $i = 1, 2$.

iii. $d$ is bidivisible by every bipolynomial which bidivides each of the bipolynimial $p^1$ and $p^2$.

*Any bipolynomial satisfying (i) and (ii) necessarily satisfies (iii).*

Next we proceed on to define greatest common bidivisor or bigreatest common divisor.

**DEFINITION 2.3.38:** *If $p^1$, $p^2$ (where $p^1 = p_1^1 \cup p_2^1$ and $p_2 = p_1^2 \cup p_2^2$) are neutrosophic bipoynomials over the neutrosophic bifield $F = F_1 \cup F_2$ such that both the bipolynomials are not $0 \cup 0$. Then the monic generator $d = d_1 \cup d_2$ of the biideal $\{p_1^1 F_1[x]\} + p_1^2 F_1[x]\} \cup \{p_2^1 F_2[x]\} + p_2^2 F_2[x]\}$ is called the greatest common bidivisor or bigreatest common divisor of $p^1$ and $p^2$. This terminology is justified by the proceeding statement. We say the neutrosophic bipolynomials $p^1 = p_1^1 \cup p_2^1$ and $p_2 = p_1^2 \cup p_2^2$ are birelatively prime if their bigreatest common divisor is $(1, 1)$ or $(I, I)$ or equivalently if the biideal they generate is all of $F[x] = F_1[x] \cup F_2[x]$.*

This result and definition 2.3.38 can be extended to any arbitrary number of bipolynomials $p^1, p^2, \ldots, p^n$, $n > 2$. We will now proceed onto define bifactorization, biprime, biirreducible of neutrosophic bipolynomials over the neutrosophic bifield $F = F_1 \cup F_2$.

**DEFINITION 2.3.39:** *Let $F = F_1 \cup F_2$ be a neutrosophic bifield. A bipolynomial $f = f_1 \cup f_2$ in $F[x] = F_1[x] \cup F_2[x]$ is said to be bireducible over the bifield $F = F_1 \cup F_2$ if there exists bipolynomials $g, h \in F[x]$, $g = g_1 \cup g_2$ and $h = h_1 \cup h_2$ in $F[x] = F_1[x] \cup F_2[x]$ of bidegree $\geq (1, 1)$ or $(I, I)$ such that $f = gh = g_1 h_1 \cup g_2 h_2 = f_1 \cup f_2$ and if such g and h does not exist $f = f_1 \cup$* $f_2$ is said to be biirreducible over the bifield $F = F_1 \cup F_2$. A non biscalar, bi-irreducible neutrosophic bipolynomial over the neutrosophic bifield $F = F_1 \cup F_2$ is called the biprime polynomial over the bifield $F = F_1 \cup F_2$ and some times we say it is biprime in $F[x] = F_1[x] \cup F_2[x]$.

The following results can be proved by any interested reader.

**THEOREM 2.3.24:** *Let $p = p^1 \cup p^2$, $f = f^1 \cup f^2$ and $g = g^1 \cup g^2$ be neutrosophic bipolynomials over the neutrosophic bifield $F = F_1 \cup F_2$. Suppose that p is a biprime bipolynomial and that p bidivides the product $fg = f_1 g_1 \cup f_2 g_2$ then either p bidivides f or p bidivides g.*

**THEOREM 2.3.25:** *If $p = p^1 \cup p^2$ is a biprime bipolynomial that bidivides a biproduct $f_1$ and $f_2$ that is $f_1 f_2$ then p bidivides one of the bipolynomial $f_1$ or $f_2$.*

**THEOREM 2.3.26:** *If $F = F_1 \cup F_2$ be a neutrosophic bifield a non zero biscalar monic neutrosophic bipolynomial in $F[x] = F_1[x] \cup F_2[x]$ can be bifactored as a biproduct of bimonic primes in $F[x] = F_1[x] \cup F_2[x]$ in one and only one way except for the order.*

**THEOREM 2.3.27:** *Let $f = f_1 \cup f_2$ be a non scalar neutrosophic monic bipolynomial over the neutrosophic bifield $F = F_1 \cup F_2$ and let $f = p_1^{n_1^1} \ldots p_{k_1}^{n_{k_1}^1} \cup p_1^{n_1^2} \ldots p_{k_2}^{n_{k_2}^2}$ be the prime bifactorization of f. For each $j_t$; $1 \leq j_t \leq k_t$; $t = 1, 2$, let $f_j^t = f^t / p_j^{n_j^t} = \prod_{i \neq j} p_i^{n_{k_i}^t}$, then $f_1^t, \ldots, f_{k_t}^t$ are relatively prime for $t = 1, 2$.*

**THEOREM 2.3.28:** *If $f = f_1 \cup f_2$ is a neutrosophic bipolynomial over the bifield $F = F_1 \cup F_2$ with derivative $f' = f_1' \cup f_2'$. Then f is a biproduct of distinct irreducible bipolynomials over the bifield $F = F_1 \cup F_2$ if and only if f and f ′ are relatively biprime, that is each $f_i$ and $f_i'$ are relatively prime for $i = 1, 2$.*

Now we proceed onto define bicharacteristics values of a strong neutrosophic bilinear operator on a strong neutrosophic bivector space $V = V_1 \cup V_2$ over a neutrosophic bifield $F = F_1 \cup F_2$.

**DEFINITION 2.3.40:** *Let $V = V_1 \cup V_2$ be a strong neutrosophic bivector space over the neutrosophic bifield $F = F_1 \cup F_2$ and let $T = T_1 \cup T_2$ be a bilinear operator on V; i.e., $T = T_1 \cup T_2$; $V = V_1 \cup V_2 \to V = V_1 \cup V_2$ and $T_i : V_i \to V_i$, $i = 1, 2$. This is the only way bilinear operator can be defined on V. A bicharacteristic value of T is a biscalar $c = c_1 \cup c_2$ ($c_i \in F_i$, $i = 1, 2$) in $F = F_1 \cup F_2$ such that there is a non zero bivector $\alpha = \alpha_1 \cup \alpha_2$ in $V = V_1 \cup V_2$ with $T\alpha = c\alpha$; i.e., $T\alpha = T_1\alpha_1 \cup T_2\alpha_2 = c_1\alpha_1 \cup c_2\alpha_2$; i.e., $T_i\alpha_i = c_i\alpha_i$ $i = 1, 2$. If $c = c_1 \cup c_2$ is a bicharacteristic value of $T = T_1 \cup T_2$ then*

i. *any $\alpha = \alpha_1 \cup \alpha_2$ such that $T\alpha = c\alpha$ is called the bicharacteristic bivector of $T = T_1 \cup T_2$ associated with the bicharacteristic value $c = c_1 \cup c_2$.*
ii. *The collection of all $\alpha = \alpha_1 \cup \alpha_2$ such that $T\alpha = c\alpha$ is called the bicharacteristic space associated with c.*

*If $T = T_1 \cup T_2$ is any bilinear operator on the bivector space $V = V_1 \cup V_2$. We call the bicharacteristic values associated with T to be bicharacteristic roots, bilatent roots bieigen values, biproper values or bispectral values.*

*These can be neutrosophic or real; will always be neutrosophic if $F = F_1 \cup F_2 = K_1 I \cup K_2 I$ where $K_1$ and $K_2$ are real fields.*

*These can be real or neutrosophic if $F = F_1 \cup F_2 = N(K_1) \cup N(K_2)$, $K_1$ and $K_2$ are real fields.*

If T is any bilinear operator and $c = c_1 \cup c_2$ in any biscalar the set of bivector $\alpha = \alpha_1 \cup \alpha_2$ such that $T\alpha = c\alpha$ is a strong neutrosophic bivector subspace of V. It is infact the binull space of the bilinear transformation $(T - cI_d) = (T_1 - c_1 I_{d_1}) \cup (T_2 - c_2 I_{d_2})$ where $I_{d_j}$ denotes the unit neutrosophic matrix for $j = 1$,

2. We call $c = c_1 \cup c_2$ the bicharacteristic value of $T = T_1 \cup T_2$ if this bispace is different from the bizero space $0 = 0 \cup 0$; that is $(T - c\, I_d) = (T_1 - c_1 I_{d_1}) \cup (T_2 - c_2 I_{d_2})$ fails to be one to one bilinear transformation that is when the bideterminant of $T - cI_d$ $= \det(T_1 - c_1 I_{d_1}) \cup \det(T_2 - c_2 I_{d_2}) = 0 \cup 0$.

We have the following theorem.

**THEOREM 2.3.29:** *Let $T = T_1 \cup T_2$ be a bilinear operator on a finite $(n_1, n_2)$ bidimensional strong neutrosophic bivector space $V = V_1 \cup V_2$ defined over the neutrosophic bifield $F = F_1 \cup F_2$ and let $c = c_1 \cup c_2$ be a biscalar in F. The following are equivalent.*
i. *$c = c_1 \cup c_2$ is a bicharacteristic value of $T = T_1 \cup T_2$.*
ii. *The bioperator $(T_1 - c_1 I_{d_1}) \cup (T_2 - c_2 I_{d_2}) = (T - c\, I_d)$ is bisingular or (nor biinvertible.*
iii. *$Det(T - cI_d) = 0 \cup 0$; i.e., $\det(T_1 - c_1 I_{d_1}) \cup \det(T_2 - c_2 I_{d_2})$*
   *$= 0 \cup 0$.*

This theorem is direct and the interested reader is expected to prove it.

Now we define the bicharacteristic value of a neutrosophic bimatrix $A = A_1 \cup A_2$ where each $A_i$ is a $n_i \times n_i$ neutrosophic matrix with entries from the neutrosophic field $F_i$, $i = 1, 2$, so that A is a neutrosophic bimatrix defined over the bifield $F = F_1 \cup F_2$. A bicharacteristic value of A in the bifield $F = F_1 \cup F_2$ is a biscalar $c = c_1 \cup c_2$ in $F = F_1 \cup F_2$ such that the bimatrix $A - cI_d$ $= (A_1 - c_1 I_{d_1}) \cup (A_2 - c_2 I_{d_2})$ is bisingular or not biinvertible.

$c = c_1 \cup c_2$ is a bicharacteristic value of $A = A_1 \cup A_2$ a $(n_1 \times n_1, n_2 \times n_2)$ neutrosophic bimatrix over the neutrosophic bifield $F = F_1 \cup F_2$ if and only if bidet $(A - cI_d) = 0 \cup 0$; i.e., $\det(A_1 - c_1 I_{d_1}) \cup \det(A_2 - c_2 I_{d_2}) = 0 \cup 0$; we form the bimatrix $(xI_d - A)$ $= (x I_{d_1} - A_1) \cup (x I_{d_2} - A_2)$. Clearly the bicharacteristic values of A in $F = F_1 \cup F_2$ are just biscalars $c = c_1 \cup c_2$ in $F_1 \cup F_2$ such that $f(c) = f_1(c_1) \cup f_2(c_2) = 0 \cup 0$. For this reason $f = f_1 \cup f_2$ is called the bicharacteristic polynomial of A. Clearly f is a neutrosophic bipolynomial of differ degrees in x over different neutrosophic fields. It is important to note that $f = f_1 \cup f_2$ is a bimonic bipolynomial which has bidegree exactly $(n_1, n_2)$. The bimonic neutrosophic bipolynomial is also a neutrosophic bipolynomial over $F = F_1 \cup F_2$.

We will illustrate this situation by some examples.

*Example 2.3.77:* Let

$$A = A_1 \cup A_2 = \begin{bmatrix} I & 0 & 1 \\ 0 & 1 & 0 \\ I & 0 & 0 \end{bmatrix} \cup \begin{bmatrix} 2 & I & 0 & I \\ I & I & 0 & 0 \\ 0 & 2 & 2I & 1 \\ 0 & 0 & 0 & 1 \end{bmatrix}$$

be a neutrosophic bimatrix of order $(3 \times 3, 4 \times 4)$ over the neutrosophic bifield $F = F_1 \cup F_2 = N(Z_2) \cup N(Z_3)$. The bicharacteristic neutrosophic bipolynomial associated with the neutrosophic bimatrix A is given by

$$(xI_d - A) = (xI_{3\times 3} - A_1) \cup (xI_{4\times 4} - A_2)$$

$$= \begin{bmatrix} x+I & 0 & 1 \\ 0 & x+1 & 0 \\ I & 0 & x \end{bmatrix} \cup \begin{bmatrix} x+1 & 2I & 0 & 2I \\ 2I & x+2I & 0 & 0 \\ 0 & 1 & x+I & 2 \\ 0 & 0 & 0 & x+2 \end{bmatrix}$$

is a neutrosophic bimatrix with neutrosophic polynomial entries.

$$\begin{aligned} f &= f_1 \cup f_2 \\ &= \det(xI_d - A) \\ &= \det(x\, I_{3\times 3} - A_1) \cup \det(x\, I_{4\times 4} - A_2) \\ &= \{(x + I)\, (x + 1)\, x + I\, (x + 1)\} \cup (x + 1)\, (x + 2I)\, (x + I) \\ &\quad (x + 2) + 2I\, (x + I)\, (x + 2)\} \\ &= \{x^3 + Ix^2 + Ix + x^2 + Ix + I\} \cup \{(x^2 + 2Ix + x + 2I) \\ &\quad (x^2 + 2I + Ix + 2I) + 2Ix^2 + I \\ &= \{x^3 + (I + 1)\, x^2 + I\} \cup \{x^4 + 2I + Ix^2 + 2x^2\}. \end{aligned}$$

Thus the bipolynomial is a monic neutrosophic polynomial of degree (3, 4) over the bifield $F = N(Z_2) \cup N(Z_3)$.

*Example 2.3.78:* Let

$$A = A_1 \cup A_2 = \begin{bmatrix} I & 0 \\ 2 & 2 \end{bmatrix} \cup \begin{bmatrix} I & 0 & 1 \\ 0 & 1 & I \\ 0 & 0 & 1 \end{bmatrix}$$

be a neutrosophic bimatrix with entries from the neutrosophic bifield $F = F_1 \cup F_2 = N(Z_3) \cup N(Z_2)$. The bicharacteristic neutrosophic bipolynomial associated with the neutrosophic bimatrix A is given by

$$(xI_d - A) = (x\,I_{2\times 2} - A_1) \cup (x\,I_{3\times 3} - A_2)$$

$$= \left\{ \begin{bmatrix} x & 0 \\ 0 & x \end{bmatrix} - \begin{bmatrix} I & 0 \\ 2 & 2 \end{bmatrix} \right\} \cup \left\{ \begin{bmatrix} x & 0 & 0 \\ 0 & x & 0 \\ 0 & 0 & x \end{bmatrix} - \begin{bmatrix} I & 0 & 1 \\ 0 & 1 & I \\ 0 & 0 & I \end{bmatrix} \right\}$$

$$= \begin{bmatrix} x+2I & 0 \\ 1 & x+1 \end{bmatrix} \cup \begin{bmatrix} x+I & 0 & 1 \\ 0 & x+1 & I \\ 0 & 0 & x+I \end{bmatrix}.$$

Let $f = f_1 \cup f_2 = \det(xI_d - A) = \det(xI_{2\times 2} - A_1) \cup \det(xI_{3\times 3} - A_2)$

$$= \begin{vmatrix} x+2I & 0 \\ 1 & x+1 \end{vmatrix} \cup \begin{vmatrix} x+I & 0 & 1 \\ 0 & x+1 & I \\ 0 & 0 & x+I \end{vmatrix}$$

$$\begin{aligned}
&= \{(x+2I)(x+1)\} \cup \{(x+I)^2 (x+1)\} \\
&= \{x^2 + 2Ix + x + 2I\} \cup \{(x^2 + 2I + I)(x+1) \\
&= (x^2 + 2Ix + x + 2I) \cup (x^2 + I)(x+1)\} \\
&= \{x^2 + (2I+1)x + 2I\} \cup \{x^3 + Ix + x^2 + I\}.
\end{aligned}$$

We see $\det(xI_d - A)$ is a neutrosophic bipolynomial which is monic neutrosophic bipolynomial of bidegree (2, 3) over the bifield $F = N(Z_3) \cup N(Z_2)$.

We now proceed onto define similar neutrosophic bimatries when the entries of these neutrosophic bimatrices are from the neutrosophic bifield $F = F_1 \cup F_2$.

**DEFINITION 2.3.41:** *Let $A = A_1 \cup A_2$ be a $(n_1 \times n_1, n_2 \times n_2)$ neutrosophic bimatrix over the neutrosophic bifield $F = F_1 \cup F_2$, that is each $A_i$ takes its entries from the neutrosophic field $F_i$, $i = 1, 2$. We say two neutrosophic bimatrices A and B of same order are similar if there exists a neutrosophic non invertible bimatrix $P = P_1 \cup P_2$ of $(n_1 \times n_1, n_2 \times n_2)$ order such that $B = P^{-1} A P$ where $P^{-1} = P_1^{-1} \cup P_2^{-1}$, $B = B_1 \cup B_2$ and*

$$B = P_1^{-1} A_1 P_1 \cup P_2^{-1} A_2 P_2.$$

Clearly
$$\begin{aligned} \det(xI_d - B) &= \det(xI_d - P^{-1} AP) \\ &= \det(P^{-1}(xI_d - A)P) \\ &= \det P^{-1} \cdot \det(xI_d - A) \det P \\ &= \det(xI_d - A) \\ &= \det(xI_{d_1} - A_1) \cup \det(xI_{d_2} - A_2). \end{aligned}$$

Thus
$$\det(xI_{d_1} - B_1) \cup \det(xI_{d_2} - B_2) = \det(xI_{d_1} - A_1) \cup \det(xI_{d_2} - A_2).$$

**DEFINITION 2.3.42:** *Let $T = T_1 \cup T_2$ be a linear bioperator on a strong neutrosophic bivector space $V = V_1 \cup V_2$ over the neutrosophic bifield $F = F_1 \cup F_2$. We say $T = T_1 \cup T_2$ is bidiagonalizable if there is a bibasis for $V = V_1 \cup V_2$ and for each bivector of which is a bicharacteristic bivecor of $T = T_1 \cup T_2$. Suppose $T = T_1 \cup T_2$ is a bidiagonalizable bilinear operator. Let*

$$\{C_1^1, \ldots, C_{k_1}^1\} \cup \{C_1^2, \ldots, C_{k_2}^2\}$$

*be the bidistinct bicharacteristic values of $T = T_1 \cup T_2$. Then there is a bibasis $B = B_1 \cup B_2$ in which T is represented by a bidiagonal matrix which has for its bidiagonal entries the* *scalars $C_i^t$ each repeated a certain number of times $t = 1, 2$. If $C_i^t$ is repeated $d_i^t$ times then the neutrosophic bimatrix has the biblock form*

$$[T]_B = [T_1]_{B_1} \cup [T_2]_{B_2}$$

$$= \begin{bmatrix} C_1^1 I_1^1 & 0 & \cdots & 0 \\ 0 & C_2^1 I_2^1 & \cdots & 0 \\ \vdots & \vdots & & \vdots \\ 0 & 0 & \cdots & C_{k_1}^1 I_{k_1}^1 \end{bmatrix} \cup \begin{bmatrix} C_1^2 I_1^2 & 0 & \cdots & 0 \\ 0 & C_2^2 I_2^2 & \cdots & 0 \\ \vdots & \vdots & & \vdots \\ 0 & 0 & \cdots & C_{k_2}^2 I_{k_2}^2 \end{bmatrix}.$$

*$I_j^t$ is the $d_j^t \times d_j^t$ identity matrix $t = 1, 2$.*

From this neutrosophic bimatrix we make the following observations.

First the bicharacteristic neutrosophic bipolynomial for $T = T_1 \cup T_2$ is the biproduct of bilinear factors $f = f_1 \cup f_2 = (x - C_1^1)^{d_1^1} \ldots (x - C_{k_1}^1)^{d_{k_1}^1} \cup (x - C_1^2)^{d_1^2} \ldots (x - C_{k_2}^2)^{d_{k_2}^2}$.

If the biscalar neutrosophic bifield $F = F_1 \cup F_2$ is bialgebraically closed, if each $F_i$ is algebraically closed for $i = 1, 2$; then every bipolynomial over $F = F_1 \cup F_2$ can be bifactored; however if $F = F_1 \cup F_2$ is not algebraically biclosed (bialgebraically closed) we are citing a special property of $T = T_1 \cup T_2$, when we say that its bicharacteristic polynomial does not have such a factorization.

The second thing to be noted is that $d_i^t$ is the number of times $C_i^t$ is repeated as a root of $f_t$ which is equal to the dimension of the space in $V_t$ of characteristic vectors associated with the characteristic value $C_i^t$; $i = 1, 2, \ldots, k_t$; $t = 1, 2$. This is because the binullity of a bidiagonal bimatrix is equal to the number of bizeros which has on its main bidiagonal and the neutrosophic bimatrix

$$[T - CI_d]_B = [T_1 - C_{1_1}^1 I_{d_1}]_{B_1} \cup [T_2 - C_{1_2}^2 I_{d_2}]_{B_2}$$

has $\left(d_{1_1}^1, d_{1_2}^2\right)$ bizeros on its main bidiagonal.

We give some results, the proof of which is direct and the interested reader can analyse them.

**THEOREM 2.3.30:** *Suppose that $T\alpha = C\alpha$ that is $(T_1 \cup T_2)(\alpha_1 \cup \alpha_2) = C_1\alpha_1 \cup C_2\alpha_2$ i.e., $T_1\alpha_1 \cup T_2\alpha_2 = C_1\alpha_1 \cup C_2\alpha_2$; $T = T_1 \cup T_2$ be a bilinear operator when the biscalar $C = C_1 \cup C_2 \in F = F_1 \cup F_2$ ($F = F_1 \cup F_2$ a neutrosophic bifield) and $\alpha = \alpha_1 \cup \alpha_2$ is a bivector from a strong neutrosophic bivector space $V = V_1 \cup V_2$ over $F = F_1 \cup F_2$. If $f = f_1 \cup f_2$ is any bipolynomial then $f(T)\alpha = f(C)\alpha$; i.e.,*

$$f_1(T_1)\alpha_1 \cup f_2(T_2)\alpha_2 = f_1(C_1)\alpha_1 \cup f_2(C_2)\alpha_2.$$

**THEOREM 2.3.31:** *Let $T = T_1 \cup T_2$ be a linear bioperator on the finite $(n_1, n_2)$ bidimensional strong neutrosophic bivector space $V = V_1 \cup V_2$ over the bifield $F = F_1 \cup F_2$. If $\{C_1^1, C_2^1, \ldots, C_{k_1}^1\} \cup \{C_1^2, C_2^2, \ldots, C_{k_2}^2\}$ be distinct bicharacteristic values of $T = T_1 \cup T_2$. Let $W_i = W_{i_1}^1 \cup W_{i_2}^2$ be the strong neutrosophic bisubspace of bicharacteristic bivectors associated with the bicharacteristic values $C_i = C_{i_1}^1 \cup C_{i_2}^2$. If $W = \{W_1^1 + \ldots + W_{k_1}^1\} \cup \{W_1^2 + \ldots + W_{k_2}^2\}$ the bidimension*

$$W = \{(\dim W_1^1 + \ldots + \dim W_{k_1}^1)\} \cup \{(\dim W_1^2 + \ldots + \dim W_{k_2}^2)\}$$
$$= \dim W^1 \cup \dim W^2.$$

*Infact if $B_{i_t}^t$ is the basis of $W_{i_t}^t$; $1 \le i_t \le k_t$, $t = 1, 2$, then $B = \{B_1^1, \ldots, B_{k_1}^1\} \cup \{B_1^2, \ldots, B_{k_2}^2\}$ is a bibasis of $W$.*

**THEOREM 2.3.32:** *Let $T = T_1 \cup T_2$ be a bilinear operator (linear bioperator) of a finite $(n_1, n_2)$ bidimension strong neutrosophic bivector space $V = V_1 \cup V_2$ over the bifield $F = F_1 \cup F_2$.*

*Let $\{C_1^1, \ldots, C_{k_1}^1\} \cup \{C_1^2, \ldots, C_{k_2}^2\}$ be the distinct bicharacteristic values of $T = T_1 \cup T_2$ and let $W_i = W_{i_1}^1 \cup W_{i_2}^2$ be* *the binull space of $T - C_iI_d = [T_1 - C^1_{l_1}I_{d_1}] \cup [T_2 - C^2_{l_2}I_{d_2}]$.*
*The following are equivalent*
*i. T is bidiagonalizable*
*ii. The bicharacteristic bipolynomial for $T = T_1 \cup T_2$ is*

$$f = f_1 \cup f_2$$
$$= (x - C^1_1)^{d^1_1} \ldots (x - C^1_{k_1})^{d^1_{k_1}} \cup (x - C^2_1)^{d^2_1} \ldots (x - C^2_{k_2})^{d^2_{k_2}}.$$

We will discuss more elaborately by giving proofs when V is a strong neutrosophic n-vector space over a neutrosophic n-field; n > 2. We define the notion of bipolynomial for the bioperator T : V → V.

**DEFINITION 2.3.43:** *Let $T = T_1 \cup T_2$ be a bilinear operator on a finite $(n_1, n_2)$ bidimensional strong neutrosophic bivector space $V = V_1 \cup V_2$ over the neutrosophic bifield $F = F_1 \cup F_2$. The biminimal neutrosophic bipolynomial for T is the unique monic bigenerator of the biideal of bipolynomials over the bifield $F = F_1 \cup F_2$ which biannihilate $T = T_1 \cup T_2$.*

The biminimal neutrosophic bipolyomial starts from the fact that the bigenerator of a neutrosophic bipolynomial biideal is characterized by being the bimonic bipolynomial of biminimum bidegree in the biideal that implies that the biminimal bipolynomial $p = p_1 \cup p_2$ for the bilinear operator $T = T_1 \cup T_2$ is uniquely determined by the following properties.

i. p is a bimonic neutrosophic bipolynomial over the biscalar neutrosophic bifield $F = F_1 \cup F_2$.
ii. $p(T) = p_1(T_1) \cup p_2(T_2) = 0 \cup 0$.
iii. No neutrosophic bipolynomial over the bifield $F = F_1 \cup F_2$ which biannihilates $T = T_1 \cup T_2$ has smaller bidegree than p = $p_1 \cup p_2$ has. $(n_1 \times n_1, n_2 \times n_2)$ to be the order of the neutrosophic bimatrix $A = A_1 \cup A_2$ over the neutrosophic bifield $F = F_1 \cup F_2$ where each $A_i$ has a $n_i \times n_i$ neutrosophic matrix with entries from the neutrosophic field $F_i$, which associated matrix of $T_i$; i = 1, 2.

The biminimal neutrosophic bipolynomial for $A = A_1 \cup A_2$ is defined in an analogous way as the unique bimonic generator of the biideal of all neutrosophic bipolynomial over the bifield $F = F_1 \cup F_2$ which biannihilate A.

If the neutrosophic linear bioperator $T = T_1 \cup T_2$ is represented by some bibasis by the neutrosophic bimatrix $A = A_1 \cup A_2$ then T and A have the same neutrosophic biminimal bipolynomial because $f(T) = f_1(T_1) \cup f_2(T_2)$ is represented in the bibasis by the neutrosophic bimatrix $f(A) = f_1(A_1) \cup f_2(A_2)$, so $f(T) = 0 \cup 0$ if and only if $f(A) = f_1(A_1) \cup f_2(A_2) = 0 \cup 0$, that is if and only if $f(T) = f_1(T_1) \cup f_2(T_2) = 0 \cup 0$.

So

$$\begin{aligned} f(P^{-1}AP) &= f_1(P_1^{-1}A_1P_1) \cup f_2(P_2^{-1}A_2P_2) \\ &= P_1^{-1}f_1(A_1)P_1 \cup P_2^{-1}f_2(A_2)P_2 \\ &= P^{-1}f(A)P \end{aligned}$$

for every neutrosophic bipolynomial $f = f_1 \cup f_2$.

Another important feature about the neutrosophic biminimal polynomials of neutrosophic bimatrices is that suppose $A = A_1 \cup A_2$ is a $(n_1 \times n_1, n_2 \times n_2)$ neutrosophic bimatrix with entries from the bifield $F = F_1 \cup F_2$. Suppose $K = K_1 \cup K_2$ is a neutrosophic bifield which contains the neutrosophic bifield $F = F_1 \cup F_2$; that is $K \supseteq F$ and $K_i \supseteq F_i$ for every i, i = 1, 2. $A = A_1 \cup A_2$ is a $(n_1 \times n_1, n_2 \times n_2)$ neutrosophic bimatrix over $F = F_1 \cup F_2$ or over $K = K_1 \cup K_2$ but we do not obtain two neutrosophic biminimal polynomial but only one neutrosophic minimal bipolynomial.

We now proceed on to prove one interesting theorem about the neutrosophic biminimal polynomials for T (or A).

**THEOREM 2.3.33:** *Let $T = T_1 \cup T_2$ be a neutrosophic linear bioperator on a $(n_1, n_2)$ bidimensional strong neutrosophic bivector space $V = V_1 \cup V_2$ (or let A be a $(n_1 \times n_1, n_2 \times n_2)$ neutrosophic bimatrix that is $A = A_1 \cup A_2$ where each $A_i$ is a $n_i \times n_i$ neutrosophic matrix with its entries from the neutrosophic field $F_i$ of $F = F_1 \cup F_2$ true for i = 1, 2). The bicharacteristic* *and biminimal neutrosophic bipolynomial for $T = T_1 \cup T_2$ (for $A = A_1 \cup A_2$) have the same biroots except for bimultiplicities.*

*Proof:* Let $p = p_1 \cup p_2$ be a neutrosophic biminimal bipolynomial for $T = T_1 \cup T_2$. Let $c = c_1 \cup c_2$ be a biscalar from the neutrosophic bifield $F = F_1 \cup F_2$. To prove $p(c) = p_1(c_1) \cup p_2(c_2) = 0 \cup 0$ if and only if $c = c_1 \cup c_2$ is the bicharacteristic value of T.

Suppose $p(c) = p_1(c_1) \cup p_2(c_2) = 0 \cup 0$ then $p = (x - c_1)q_1 \cup (x - c_2)q_2$ where $q = q_1 \cup q_2$ is the neutrosophic bipolynomial since bideg q < bideg p, the neutrosophic biminimal bipolynomial $p = p_1 \cup p_2$ tells us $q(T) = q_1(T_1) \cup q_2(T_2) \neq 0 \cup 0$. Choose the bivector $\beta = \beta_1 \cup \beta_2$ such that $q(T)\beta = q_1(T_1)\beta_1 \cup q_2(T_2)\beta_2 \neq 0 \cup 0$. Let $\alpha = q(T) p$ that is $\alpha = \alpha_1 \cup \alpha_2 = q_1(T_1)\beta_1 \cup q_2(T_2)\beta_2$.

Then
$$\begin{aligned} 0 \cup 0 &= p(T)\beta \\ &= p_1(T_1)\beta_1 \cup p_2(T_2)\beta_2 \\ &= (T - cI)\, q(T)\, \beta \\ &= (T_1 - c_1 I_1)q_1(T_1)\beta_1 \cup (T_2 - c_2 I_2)q_2(T_2)\beta_2 \\ &= (T_1 - c_1 I_1)\alpha_1 \cup (T_2 - c_2 I_2)\alpha_2 \end{aligned}$$
and thus $c = c_1 \cup c_2$ is a bicharacteristic value of $T = T_1 \cup T_2$.

Suppose $c = c_1 \cup c_2$ is the bicharacteristic value of the bilinear operator $T = T_1 \cup T_2$ say $T\alpha = c\alpha$; i.e., $T_1 \alpha_1 \cup T_2 \alpha_2 = c_1 \alpha_1 \cup c_2 \alpha_2$ with $\alpha \neq 0 \cup 0$.
From the earlier results we have $p(T)\alpha = p(c)\alpha$ that is
$$p_1(T_1)\alpha_1 \cup p_2(T_2)\alpha_2 = p_1(c_1)\alpha_1 \cup p_2(c_2)\alpha_2.$$
Since $p(T) = p_1(T_1) \cup p_2(T_2) = 0 \cup 0$ and $\alpha = \alpha_1 \cup \alpha_2 \neq 0 \cup 0$ we have $p_1(c_1) \cup p_2(c_2) = p(c) \neq 0 \cup 0$.

Let $T = T_1 \cup T_2$ be a bidiagonalizable bilinear operator and let $\{c_1^1,...,c_{k_1}^1\} \cup \{c_1^2,...,c_{k_2}^2\}$ be the bidistinct bicharacteristic values of T. Then the biminimal neutrosophic bipolynomial for $T = T_1 \cup T_2$ is the neutrosophic bipolynomial $p = p_1 \cup p_2 = (x - c_1^1) \ldots (x - c_{k_1}^1) \cup (x - c_1^2) \ldots (x - c_{k_2}^2)$.

If $\alpha = \alpha_1 \cup \alpha_2$ is a bicharacteristic bivector then one of the bioperators $\{(T_1 - c_1^1 I_1), ..., (T_1 - c_{k_1}^1 I_1)\} \cup \{(T_2 - c_1^2 I_2), ...,$ $(T_2 - c_{k_2}^2 I_2)\}$ send $\alpha = \alpha_1 \cup \alpha_2$ into $0 \cup 0$, thus resulting in $\{(T_1 - c_1^1 I_1), \ldots, (T_1 - c_{k_1}^1 I_1)\} \cup \{(T_2 - c_1^2 I_2), \ldots, (T_2 - c_{k_2}^2 I_2)\} = 0 \cup 0$ for every bicharacteristic bivector $\alpha = \alpha_1 \cup \alpha_2$.

Hence there exists a bibasis for the underlying bispace which consists of bicharacteristic vectors of $T = T_1 \cup T_2$. Hence $p(T) = p_1(T_1) \cup p_2(T_2) = \{(T_1 - c_1^1 I_1), \ldots, (T_1 - c_{k_1}^1 I_1)\} \cup \{(T_2 - c_1^2 I_2), \ldots, (T_2 - c_{k_2}^2 I_2)\} = 0 \cup 0$.

Thus we conclude if T is bidiagonlizable bilinear operator then the neutrosophic biminimal bipolynomial for $T = T_1 \cup T_2$ is a product of bidistinct bilinear factors.

Now we proceed onto prove the Cayley Hamilton theorem for strong neutrosophic bivector spaces of finite bidimension defined over the neutrosophic bifield of Type II.

**THEOREM 2.3.34: (Cayley Hamilton):** *Let $T = T_1 \cup T_2$ be a bilinear operator on a finite $(n_1, n_2)$ bidimensional strong neutrosophic bivector space defined over a neutrosophic bifield $F = F_1 \cup F_2$. If $f = f_1 \cup f_2$ is the bicharacteristic neutrosophic bipolynomial for T then $f(T) = f_1(T_1) \cup f_2(T_2) = 0 \cup 0$, in other words the biminimal neutrosophic bipolynomial bidivides the bicharacteristic neutrosophic bipolynomial for T.*

*Proof:* Let $K = K_1 \cup K_2$ be a bicommuting neutrosophic ring with biidentity $I_2 = (1, 1)$ consisting of all bipolynomial in $T = T_1 \cup T_2$. $K = K_1 \cup K_2$ is actually a bicommuting bialgebra with biidentity over the neutrosophic bifield $F = F_1 \cup F_2$ (that is both $F_1$ and $F_2$ are not pure).

Let $\{\alpha_1^1, \ldots, \alpha_{n_1}^1\} \cup \{\alpha_1^2, \ldots, \alpha_{n_2}^2\}$ be a bibasis for $V = V_1 \cup V_2$ and let $A = A^1 \cup A^2$ be a bimatrix which represents $T = T_1 \cup T_2$ in the given bibasis.

Then
$$T\alpha_i = T_1 \alpha_{i_1}^1 \cup T_2 \alpha_{i_2}^2$$
$$= \sum_{j_1=1}^{n_1} A_{j_1 i_1}^1 \alpha_{j_1}^1 \cup \sum_{j_2=1}^{n_2} A_{j_2 i_2}^2 \alpha_{j_2}^2 \ ;$$

$1 \leq j_i \leq n_{j_i}$ ; $i = 1, 2$. These biequations may be equivalently written in the form

$$\sum_{j_1=1}^{n_1}\left(\delta_{j_1 i_1} T_1 - A^1_{j_1 i_1} I_{i_1}\right)\alpha^1_{j_1} \cup \sum_{j_2=1}^{n_2}\left(\delta_{j_2 i_2} T_2 - A^2_{j_2 i_2} I_{i_2}\right)\alpha^2_{j_2} = 0 \cup 0.$$

Let $B = B^1 \cup B^2$ denote the element of $K_1^{n_1 \times n_1} \cup K_2^{n_2 \times n_2}$; i.e., $B^i$ is an element of $K_i^{n_i \times n_i}$ with entries $B^t_{i_t j_t} = \delta_{i_t j_t} T_t - A_{i_t j_t} I_t$ ;. $t = 1, 2$. When $n_t = 2$, $1 \leq j_t, i_t \leq n_t$;

$$B^t = \begin{bmatrix} T_t - A^t_{11} I_t & A^t_{21} I_t \\ -A^t_{12} I_t & T_t - A^t_{22} I_t \end{bmatrix}$$

and $\det B^t = \left(T_t - A^t_{11} I_t\right)\left(T_t - A^t_{22} I_t\right) - \left(A^t_{12} A^t_{21}\right) I_t = f_t(I_t)$ where $f_t$ is the neutrosophic characteristic polynomial associated with $T_t$, $t = 1, 2$. $f_t = x^2 -$ trace $A^t_x + \det A^t$. For case $n_t > 2$ it is clear that $\det B^t = f_t(T_t)$ since $f_t$ is the determinant of the neutrosophic matrix $xI_t - A_t$ whose entries are neutrosophic polynomial;
$(xI_t - A^t)_{i_t j_t} = \delta^x_{i_t j_t} - A^t_{i_t j_t}$.

We will shown $f_t(T_t) = 0$. In order that $f_t(T_t)$ is a zero operator it is necessary and sufficient that

$(\det B^t)_{\alpha^t_{k_t}} = 0$ for $k_t = 0, 1, \ldots, n_t$.

By definition of $B^t$ the vectors $\alpha^t_1 \cup \ldots \cup \alpha^t_{n_t}$ satisfy the equations;

$$\sum_{j_t=1}^{n_t} B^t_{i_t j_t} \alpha^t_{j_t} = 0; \ 1 \leq i_t \leq n_t.$$

When $n_t = 2$ we can write the above equation in the form

$$\begin{bmatrix} T_t - A^t_{11} I_t & -A^t_{21} I_t \\ -A^t_{12} I_t & T_t - A^t_{22} I_t \end{bmatrix} \begin{bmatrix} \alpha^t_1 \\ \alpha^t_2 \end{bmatrix} = \begin{bmatrix} 0 \\ 0 \end{bmatrix}.$$

In this case the usual adjoint $B^t$ is the neutrosophic matrix

$$\tilde{B}^t = \begin{bmatrix} T_t - A_{22}^t I_t & A_{21}^t I_t \\ A_{12}^t I & T_t - A_{22}^t I_t \end{bmatrix}$$

and

$$\tilde{B}^t B^t = \begin{bmatrix} \det B^t & 0 \\ 0 & \det B^t \end{bmatrix}.$$

Hence

$$\det B^t \begin{bmatrix} \alpha_1^t \\ \alpha_2^t \end{bmatrix} = \tilde{B}^t B^t \begin{bmatrix} \alpha_1^t \\ \alpha_2^t \end{bmatrix} =$$

$$\tilde{B}^t B^t \begin{bmatrix} \alpha_1^t \\ \alpha_2^t \end{bmatrix} = \tilde{B}^t B^t \begin{bmatrix} \alpha_1^t \\ \alpha_2^t \end{bmatrix} = \begin{bmatrix} 0 \\ 0 \end{bmatrix}.$$

In the general case $\tilde{B}^t = \text{adj } B^t$. Then $\sum_{j_t=1}^{n_t} \tilde{B}_{k_t i_t}^t B_{i_t j_t}^t \alpha_{j_t}^t = 0$

for each pair $k_t$, $i_t$ and summing on $i_t$ we have

$$0 = \sum_{i_t=1}^{n_t} \sum_{j_t=1}^{n_t} \tilde{B}_{k_t i_t}^t B_{i_t j_t}^t \alpha_{j_t}^t$$

$$= \sum_{i_t=1}^{n_t} \left( \sum_{j_t=1}^{n_t} \tilde{B}_{k_t i_t}^t B_{i_t j_t}^t \alpha_{j_t}^t \right).$$

Now $\tilde{B}^t B^t = (\det B_t) I_t$ so that

$$\sum_{i_t=1}^{n_t} \tilde{B}_{k_t i_t}^t B_{i_t j_t}^t = \delta_{k_t j_t} \det B^t$$

Therefore

$$0 = \sum_{j_t=1}^{n_t} \delta_{k_t j_t} \left( \det B^t \right) \alpha_{j_t}^t = (\det B^t) \alpha_{k_t}^t \ ; \ 1 \leq k_t \leq n_t.$$

Since this is true for each t; t = 1, 2; we have $0 \cup 0 = (\det B^1)$ $\alpha_{k_1}^1 \cup (\det B^2) \alpha_{k_2}^2$ ; $1 \leq k_i \leq n_i$; i = 1, 2.

The Cayley-Hamilton theorem is very important for it is useful in narrowing down the search for the biminimal neutrosophic bipolynomials of various bioperators.

If we know the neutrosophic bimatrix $A = A^1 \cup A^2$ which represents $T = T_1 \cup T_2$ in some ordered bibasis then we can compute the bicharacteristic neutrosophic bipolynomial $f = f_1 \cup f_2$. We know the biminimal neutrosophic polynomial $p = p_1 \cup p_2$ bidivides f that is each $p_i / f_i$; for $i = 1, 2$ (which we call as bidivides f) and that the two neutrosophic bipolynomials have the same biroots.

However we do not have a method of computing the roots even in case of polynomials so it is more difficult in case of finding the biroots of the neutrosophic bipolynomials. However if $f = f_1 \cup f_2$ factors as $f = (x - c_1^1)^{d_1^1} \ldots (x - c_{k_1}^1)^{d_{k_1}^1} \cup (x - c_1^2)^{d_1^2} \ldots (x - c_{k_2}^2)^{d_{k_2}^2}$ the distinct bisets $d_{i_t}^t \geq t$; $t = 1, 2, \ldots, k_t$ then $p = p_1 \cup p_2 = (x - c_1^1)^{r_1^1} \ldots (x - c_{k_1}^1)^{r_{k_1}^1} \cup (x - c_1^2)^{r_1^2} \ldots (x - c_{k_2}^2)^{r_{k_2}^2}$; $1 \leq r_j^t \leq d_j^t$.

We will illustrate this by a simple example.

*Example 2.3.79:* Let

$$A = A_1 \cup A_2 = \begin{bmatrix} 3 & 1 & I \\ 2 & 2I & 1 \\ 2I & 2 & 0 \end{bmatrix} \cup \begin{bmatrix} 0 & I \\ I & 1 \end{bmatrix}$$

be a neutrosophic bimatrix with entries from the neutrosophic bifield $F = F_1 \cup F_2 = N(Z_5) \cup N(Z_2)$. Clearly the bicharacteristic neutrosophic bipolynomial associated with the neutrosophic bimatrix A is given by

$$\begin{aligned} f &= f_1 \cup f_2. \\ &= \begin{vmatrix} x+2 & 4 & 4I \\ 3 & x+3I & 4 \\ 3I & 3 & x \end{vmatrix} \cup \begin{vmatrix} x & I \\ I & x+1 \end{vmatrix} \\ &= (x^3 + (3I + 1)x^2 + (4I + 1)x + 3I + 1) \cup (x(x+1) + I) \end{aligned}$$

is the biminimal neutrosophic bipolynomial of the neutrosophic bimatrix A.

Now we proceed onto define the notion of biinvariant subspaces or equivalently we may call them as invariant bisubspaces.

**DEFINITION 2.3.44:** *Let $V = V_1 \cup V_2$ be the strong neutrosophic bivector space over the bifield $F = F_1 \cup F_2$ of type II. Let $T = T_1 \cup T_2$ be a bilinear operator on V. If $W = W_1 \cup W_2$ is a strong neutrosophic bivector subspace of V we say W is biinvariant under T if each of the bivectors in W, i.e., for the bivector $\alpha = \alpha_1 \cup \alpha_2$ in W the bivector $T\alpha = T_1\alpha_1 \cup T_2\alpha_2$ is in W; i.e., each $T_i \alpha_i \in W_i$ for every $\alpha_i \in W_i$ under the operator $T_i$ for $i = 1, 2$, i.e., if T(W) is contained in W; that is $T_i(W_i) \subseteq W_i$ for $i = 1, 2$. i.e., $T(W) = T_1(W_1) \cup T_2(W_2) \subseteq W_1 \cup W_2$.*

The simple examples are we can say $V = V_1 \cup V_2$, the strong neutrosophic bivector space is invariant under a bilinear operator T on V. Similarly the zero subspace of a strong neutrosophic bivector space is invariant under T.

Now we proceed onto give the biblock neutrosophic matrix associated with a bioperator T of V.

Let $W = W_1 \cup W_2$ be a strong neutrosophic bivector subspace of the strong neutrosophic bivector space $V = V_1 \cup V_2$. Let $T = T^1 \cup T^2$ be a bioperator on V such that $W = W_1 \cup W_2$ is biinvariant under the bioperator T then $T = T^1 \cup T^2$ induces a bilinear operator; $T_w = T^1_{W_1} \cup T^2_{W_2}$ on the bisubspace W. This bilinear operator $T_w$ defined by $T_w(\alpha) = T(\alpha)$ for all $\alpha \in W$; i.e., if $\alpha = \alpha_1 \cup \alpha_2$ then
$$T_w(\alpha) = T_w(\alpha_1 \cup \alpha_2)$$
$$= T^1_{W_1}(\alpha_1) \cup T^2_{W_2}(\alpha_2).$$

Clearly $T_w$ is different from T as bidomain is W and not V. When $V = V_1 \cup V_2$ is a $(n_1, n_2)$ finite bidimensional the biinvariance of $W = W_1 \cup W_2$ under $T = T_1 \cup T_2$ has a simple neutrosophic bimatrix interpretation.

Let $B = B_1 \cup B_2 = \{\alpha_1^1,...,\alpha_{r_1}^1\} \cup \{\alpha_1^2,...,\alpha_{r_2}^2\}$ be a bibasis for W. The bidimension of W is $(r_1, r_2)$.

Let $A = [T]_B$, that is if $A = A_1 \cup A_2$ is the neutrosophic bimatrix such that $A = A_1 \cup A_2 = \left[T^1\right]_{B_1} \cup \left[T^2\right]_{B_2}$ so that

$$T_{\alpha_{j_t}^t}^t \sum_{i_t=1}^{n_t} A_{i_t j_t}^t \alpha_{i_t}^t$$

for $i = 1, 2$. Thus

$$T\alpha_j = T^1 \alpha_{j_1}^1 \cup T^2 \alpha_{j_2}^2$$

$$= \sum_{i_1=1}^{n_1} A_{i_1 j_1}^1 \alpha_{i_1}^1 \cup \sum_{i_2=1}^{n_2} A_{i_2 j_2}^2 \alpha_{i_2}^2 .$$

Since $W = W_1 \cup W_2$ is biinvariant under T the bivector $T\alpha_j$ belong to W for $j_1 < r_1$ and $j_2 < r_2$.

$$T\alpha_j = \sum_{i_1=1}^{r_1} A_{i_1 j_1}^1 \alpha_{i_1}^1 \cup \sum_{i_2=1}^{r_2} A_{i_2 j_2}^2 \alpha_{i_2}^2$$

that is $A_{i_k j_k}^k = 0$ if $j_k < r_k$ and $i_k > r_k$ for every $k = 1, 2$. Schematically A has the biblock

$$A = \begin{bmatrix} B & C \\ 0 & D \end{bmatrix} = \begin{bmatrix} B^1 & C^1 \\ 0 & D^1 \end{bmatrix} \cup \begin{bmatrix} B^2 & C^2 \\ 0 & D^2 \end{bmatrix}$$

where $B^t$ is a $r_t \times r_t$ neutrosophic matrix, $C^t$ is a $r_t \times (n_t - r_t)$ neutrosophic matrix and $D^t$ is a $(n_t - r_t) \times (n_t - r_t)$ neutrosophic matrix for $t = 1, 2$. $B = B_1 \cup B_2$ is the neutrosophic bimatrix induced by the bioperator $T_w$ on the bibasis $B' = B_1' \cup B_2'$. In view of the above properties we have the following Lemma.

**LEMMA 2.3.1:** *Let $W = W_1 \cup W_2$ be a biinvariant strong neutrosophic bisubspace of the bioperator $T = T_1 \cup T_2$ on the strong neutrosophic bivector space $V = V_1 \cup V_2$ over the neutrosophic bifield $F = F_1 \cup F_2$ which is not pure. The bicharacteristic neutrosophic bipolynomial for the birestriction* *operator $T_w = T_{1W_1} \cup T_{2W_2}$ bidivides the neutrosophic bicharacteristic polynomial for T. The biminimal neutrosophic bipolynomial for $T_w$ bidivides the biminimal neutrosophic polynomial for T.*

*Proof:* We know from the above results

$$A = \begin{bmatrix} B & C \\ 0 & D \end{bmatrix} = \begin{bmatrix} B^1 & C^1 \\ 0 & D^1 \end{bmatrix} \cup \begin{bmatrix} B^2 & C^2 \\ 0 & D^2 \end{bmatrix}$$

where
$$A = [T]_B = T_{1B_1} \cup T_{2B_2}$$

and
$$B = [T_w]_{B'} = B^1 \cup B^2 = [T_{1w_1}]_{B'_1} \cup [T_{2w_2}]_{B'_2}.$$

Because of the biblock form of the neutrosophic bimatrix

det (xI – A)
$\quad = \quad$ det $(xI_1 - A_1) \cup$ det $(xI_2 - A_2)$ (where $A = A_1 \cup A_2$)
$\quad = \quad$ det $(xI - B)$ det $(xI - D)$
$\quad = \quad$ {det $(xI_1 - B^1)$ det $(xI_1 - D^1)$ $\cup$
$\qquad$ det $(xI_2 - B^2)$ det $(xI_2 - D^2)$}.

That proves the statement about bicharacteristic neutrosophic polynomials. Notice that we used $I = I_1 \cup I_2$ to represent the biidentity matrix of the bituple of different sizes. The $k^{th}$ power of the neutrosophic bimatrix has the biblock form,

$$A^k = (A^1)^k \cup (A^2)^k$$

$$A^k = \begin{bmatrix} (B^1)^k & (C^1)^k \\ 0 & (D^1)^k \end{bmatrix} \cup \begin{bmatrix} (B^2)^k & (C^2)^k \\ 0 & (D^2)^k \end{bmatrix}$$

where $C^k = (C^1)^k \cup (C^2)^k$ is $\{(r_1 \times n_1 - r_1), (r_2 \times n_2 - r_2)\}$ bimatrix. Therefore any neutrosophic bipolynomial which biannihilates A also biannihilates B (and D too). So the biminimal neutrosophic bipolynomial for B bidivides the biminimal polynomial for A.

Let $T = T_1 \cup T_2$ be any linear bioperator on a $(n_1, n_2)$ finite dimensional space $V = V_1 \cup V_2$ over the neutrosophic bifield $F = F_1 \cup F_2$ (where both $F_1$ and $F_2$ are not pure neutrosophic fields). Let $W = W_1 \cup W_2$ be a strong neutrosophic bivector subspace of V spanned by all bicharacteristic bivectors of $T = T_1 \cup T_2$. Let $\{C_1^1, \ldots, C_{k_1}^1\} \cup \{C_1^2, \ldots, C_{k_2}^2\}$ be the bidistinct characteristic values of T. For each i let $W_i = W_{i_1}^1 \cup W_{i_2}^2$ be the strong neutrosophic bivector space associated with the bicharacteristic value $C_i = C_{i_1}^1 \cup C_{i_2}^2$ and let $B_i = \{B_i^1 \cup B_i^2\}$ be the ordered basis of $W_i$, i.e., $B_i^t$ is a basis of $W_i^t$.

$$B' = \{B_1^1, \ldots, B_{k_1}^1\} \cup \{B_1^2, \ldots, B_{k_2}^2\}$$

is a biordered bibasis for

$$W = \{W_1^1 + \ldots + W_{k_1}^1\} \cup \{W_1^2 + \ldots + W_{k_2}^2\} = W_1 \cup W_2.$$

In particular bidimension

$$= \{\dim W_1^1 + \ldots + \dim W_{k_1}^1\} \cup \{\dim W_1^2 + \ldots + \dim W_{k_2}^2\}.$$

We prove the result for one particular $W_i = \{W_1^i + \ldots + W_{k_i}^i\}$ and since $W_i$ is arbitrarily chosen the result is true for every i, i = 1, 2. $B_i' = \{\alpha_1^i, \ldots, \alpha_{r_i}^i\}$ so that the first few $\alpha^i$'s form a basis $B_i'$, the next $B_2'$. Then $T_i \alpha_j^t = t_j^i \alpha_j^t$; j = 1, 2, …, $r_i$ where $\{t_1^i, \ldots, t_{r_i}^i\} = \{C_1^i, \ldots, C_{k_i}^i\}$ where $C_j^i$ is repeated dim $W_j^i$ times j = 1, …, $r_i$. Now $W_i$ is invariant under $T_i$ since for each $\alpha^i$ in $W_i$, we have

$$\alpha^i = x_1^i \alpha_1^i + \ldots + x_{r_i}^i \alpha_{r_i}^i$$

$$T_i \alpha^i = t_1^i x_1^i \alpha_1^i + \ldots + t_{r_i}^i x_{r_i}^i \alpha_{r_i}^i.$$

Choose any other vector $\alpha_{r_i+1}^i, \ldots, \alpha_{n_i}^i$ in $V_i$ such that $B_i = \{\alpha_1^i, \ldots, \alpha_{n_i}^i\}$ is a basis for $V_i$. The matrix of $T_i$ relative to $B_i$ has the block form mentioned earlier and the neutrosophic matrix of the restriction operator relative to the basis $B_i'$ is

$$B^i = \begin{bmatrix} t_1^i & 0 & \cdots & 0 \\ 0 & t_2^i & \cdots & 0 \\ \vdots & \vdots & & \vdots \\ 0 & 0 & \cdots & t_{r_i}^i \end{bmatrix}.$$

The characteristic neutrosophic polynomial of $B^i$; i.e., of $T_{iw_i}$ is $g_i = g_i (x - C_1^i)^{e_1^i} \ldots (x - C_{k_i}^i)^{e_{k_i}^i}$ where $e_j^i = \dim W_j^i$; $j = 1, 2, \ldots, k_i$. Further more $g_i$ divides $f_i$, the characteristic neutrosophic polynomial for $T_i$. Therefore the multiplicity of $C_j^i$ as a root of $f_i$ is atleast $\dim W_j^i$. Thus $T_i$ is diagonalizable if and only if $r_i = n_i$, i.e., if and only $\{e_1^i + \ldots + e_{k_i}^i\} = n_i$. Since what we proved for $T_i$ is true for $T = T_1 \cup T_2$. Hence true for every $B^1 \cup B^2$.

We now proceed onto define T biconductor of $\alpha$ into $W = W_1 \cup W_2 \subseteq V_1 \cup V_2$.

**DEFINITION 2.3.45:** *Let $W = W_1 \cup W_2$ be a biinvariant strong neutrosophic bivector subspace for $T = T_1 \cup T_2$ and let $\alpha = \alpha_1 \cup \alpha_2$ be a bivector in the strong neutrosophic bivector space $V = V_1 \cup V_2$. The T-biconductor of $\alpha = \alpha_1 \cup \alpha_2$ into $W = W_1 \cup W_2$ is the biset $S_\tau(\alpha; W) = S_{\tau_1}(\alpha_1; W_1) \cup S_{\tau_2}(\alpha_2; W_2)$ which consists of all neutrosophic bipolynomials $g = g_1 \cup g_2$ over the neutrosophic bifield $F = F_1 \cup F_2$ such that $g(T)\alpha$ is in $W$; that is $g_1(T_1)\alpha_1 \cup g_2(T_2)\alpha_2 \in W_1 \cup W_2$.*

Since the bioperator T will be fixed throughout the discussions we shall usually drop the subscript T and write $S(\alpha; W) = S(\alpha_1; W_1) \cup S(\alpha_2; W_2)$. The authors usually call the collection of neutrosophic bipolynomials the bistuffer. We as in case of vector spaces prefer to call as biconductor that is the bioperator $g(T) = g_1(T_1) \cup g_2(T_2)$; slowly leads to the bivector $\alpha_1 \cup \alpha_2$ into $W = W_1 \cup W_2$. In the special case when $W = \{0\} \cup \{0\}$ the biconductor is called the T annihilator of $\alpha_1 \cup \alpha_2$.

We prove the following interesting lemma.

**LEMMA 2.3.2:** *If $W = W_1 \cup W_2$ is a strong neutrosophic biinvariant subspace for $T = T_1 \cup T_2$ then W is biinvariant under every neutrosophic bipolynomial in $T = T_1 \cup T_2$. Thus for each $\alpha = \alpha_1 \cup \alpha_2$ in $V = V_1 \cup V_2$ the biconductor $S(\alpha; W) = S(\alpha_1; W_1) \cup S(\alpha_2; W_2)$ is a biideal in the neutrosophic bipolynomial algebra $F[x] = F_1[x] \cup F_2[x]$ where $F_1$ and $F_2$ are neutrosophic fields and $F_1$ and $F_2$ are not pure neutrosophic fields.*

*Proof:* Given $W = W_1 \cup W_2 \subseteq V = V_1 \cup V_2$ is a strong neutrosophic bivector subspace over the neutrosophic bifield $F = F_1 \cup F_2$ (Both $F_1$ and $F_2$ are not pure neutrosophic), of the strong neutrosophic bivector space $V = V_1 \cup V_2$. If $\beta = \beta_1 \cup \beta_2$ is in $W = W_1 \cup W_2$, then $T\beta = T_1 \beta_1 \cup T_2 \beta_2$ is in $W = W_1 \cup W_2$. Thus $T(T\beta) = T^2 \beta = T_1^2 \beta_1 \cup T_2^2 \beta_2$ is in $W = W_1 \cup W_2$. By induction $T_\beta^k = T_1^{k_1} \beta_1 \cup T_2^{k_2} \beta_2$ is in $W = W_1 \cup W_2$ for every neutrosophic bipolynomial $f = f_1 \cup f_2$.

The definition $S(\alpha; W) = S(\alpha_1; W_1) \cup S(\alpha_2; W_2)$ is meaningful if $W = W_1 \cup W_2$ is any bisubset of W. If W is a strong neutrosophic bivector subspace then $S(\alpha; W)$ is a strong neutrosophic bisubspace of $F[x] = F_1[x] \cup F_2[x]$ because (cf + g)T = cf(T) + g(T); i.e.,$(c_1f_1 + g_1)T_1 \cup (c_2f_2 + g_2)T_2 = c_1f_1(T_1) + g_1(T_1) \cup c_2f_2(T_2) + g_2(T_2)$. If $W = W_1 \cup W_2$ is also biinvariant under $T = T_1 \cup T_2$ and let $g = g_1 \cup g_2$ be a neutrosophic bipolynomial in $S(\alpha; W) = S(\alpha_1; W_1) \cup S(\alpha_2; W_2)$; i.e., let $g(T)\alpha = g_1(T_1)\alpha_1 \cup g_2(T_2)\alpha_2$ be in $W = W_1 \cup W_2$ is any neutrosophic bipolynomial then f(T) g(T) $\alpha$ is in $W = W_1 \cup W_2$ that is $f(T)[g(T)\alpha] = f_1(T_1)[g_1(T_1)\, \alpha_1] \cup f_2(T_2)[g_2(T_2)\, \alpha_2]$ will be in $W = W_1 \cup W_2$. Since (fg)T = f(T)g(T) we have
$$(f_1g_1)T_1 \cup (f_2g_2)T_2 = f_1(T_1)g_1(T_1) \cup f_2(T_2)\, g_2(T_2);$$

(fg) $\in S(\alpha; W)$; that is $(f_ig_i) \in S(\alpha_i; W_i)$; i = 1, 2. Hence the claim.

The unique bimonic generator of the neutrosophic biideal $S(\alpha; W)$ is also called the T biconductor of $\alpha = \alpha_1 \cup \alpha_2$ in W (the T biannihilator in case $W = \{0\} \cup \{0\}$). The T biconductor of $\alpha$ into W is the bimonic neutrosophic bipolynomial g of least bidegree such that $g(T)\alpha = g_1(T_1)\alpha_1 \cup g_2(T_2)\alpha_2$ is in $W = W_1 \cup W_2$.

A neutrosophic bipolynomial $f = f_1 \cup f_2$ is in $S(\alpha; W) = S(\alpha_1; W_1) \cup S(\alpha_2; W_2)$ if and only if g bidivides f.

Note the biconductor $S(\alpha; W)$ always contains the bipolynomial for T, hence every T biconductor bidivides the biminimal polynomial for $T = T_1 \cup T_2$.

We prove the following lemma.

**LEMMA 2.3.3:** *Let $V = V_1 \cup V_2$ be a strong neutrosophic $(n_1, n_2)$ bidimensional bivector space over the neutrosophic bifield $F = F_1 \cup F_2$ (both $F_1$ and $F_2$ are not pure neutrosophic fields). Let $T = T_1 \cup T_2$ be a bilinear operator on V such that the neutrosophic biminimal polynomial for T is a product of bilinear factors $p = p_1 \cup p_2 = (x - c_1^1)^{r_1^1} \ldots (x - c_{k_1}^1)^{r_{k_1}^1} \cup (x - c_1^2)^{r_1^2} \ldots (x - c_{k_2}^2)^{r_{k_2}^2}$ ; $c_{t_i}^i \in F_i$; $1 \le t_i \le k_i$. $i = 1, 2$.*

*Let $W = W_1 \cup W_2$ be a strong neutrosophic proper bivector subspace of V ($V \ne W$) which is biinvariant under T. There exists a bivector $\alpha = \alpha_1 \cup \alpha_2$ in $V = V_1 \cup V_2$ such that*
i. *$\alpha$ is not in $W = W_1 \cup W_2$*
ii. *$(T - cI)\alpha = (T_1 - c_1 I_1)\alpha_1 \cup (T_2 - c_2 I_2)\alpha_2$*

*is in $W = W_1 \cup W_2$ for some bicharacteristic value of the bioperator T.*

*Proof:* (1) and (2) express that T biconductor of $\alpha = \alpha_1 \cup \alpha_2$ into $W = W_1 \cup W_2$ is a neutrosophic bilinear bipolynomial. Suppose $\beta = \beta_1 \cup \beta_2$ is any bivector in $V = V_1 \cup V_2$ which is not in $W = W_1 \cup W_2$. Let $g = g_1 \cup g_2$ be the T biconductor of $\beta$ in $W = W_1 \cup W_2$. Then g bidivides $p = p_1 \cup p_2$ the neutrosophic biminimal bipolynomial for $T = T_1 \cup T_2$. Since $\beta = \beta_1 \cup \beta_2$ is not in $W = W_1 \cup W_2$, the neutrosophic bipolynomial g is not constant. Therefore $g = g_1 \cup g_2 = (x - c_1^1)^{e_1^1} \ldots (x - c_{k_1}^1)^{e_{k_1}^1} \cup (x$ $- c_1^2)^{e_1^2} \ldots (x - c_{k_2}^2)^{e_{k_2}^2}$; where atleast one of the bipair of integers $e_i^1 \cup e_i^2$ is positive. Choose $j_t$ so that $e_{j_t}^t > 0$, t = 1, 2, then $(x - c_j) = (x - c_{j_1}^1) \cup (x - c_{j_2}^2)$ bidivides g. $g = (x - c_j)h$; i.e., $g = g_1 \cup g_2 = (x - c_{j_1}^1)h_1 \cup (x - c_{j_2}^2)h_2$. By the definition of g the bivector $\alpha = \alpha_1 \cup \alpha_2 = h_1(T_1)\beta_1 \cup h_2(T_2)\beta_2 = h(T)\beta$ cannot be in W. But $(T - c_jI) \alpha = (T - c_jI)h(T)\beta = g(T)\beta$ is in W

$(T_1 - c_{j_1}^1 I_1) \alpha_1 \cup (T_2 - c_{j_2}^2 I_2)$
$= (T_1 - c_{j_1}^1)h_1(T_1)\beta_1 \cup (T_2 - c_{j_2}^2)h_2(T_2)\beta_2$
$= g_1(T_1) \beta_1 \cup g_2(T_2) \beta_2$

with $g_i(T_i)\beta_i \in W_i$ for i = 1, 2.

Next we obtain the condition for T to be bitriangulable.

**THEOREM 2.3.35:** *Let $V = V_1 \cup V_2$ be a $(n_1, n_2)$ finite bidimensional strong neutrosophic bivector space over the bifield $F = F_1 \cup F_2$ ($F_1$ and $F_2$ are neutrosophic fields and they are not pure neutrosophic fields) and let $T = T_1 \cup T_2$ be a bilinear operator on $V = V_1 \cup V_2$. Then T is bitriangulable if and only if the biminimal neutrosophic bipolynomial for T is a biproduct of bilinear neutrosophic bipolynomials over the neutrosophic bifield $F = F_1 \cup F_2$.*

*Proof:* Suppose the biminimal neutrosophic bipolynomial $p = p_1 \cup p_2$, bifactors as $p = (x - c_1^1)^{r_1^1} \ldots (x - c_{k_1}^1)^{r_{k_1}^1} \cup (x - c_1^2)^{r_1^2} \ldots (x - c_{k_2}^2)^{r_{k_2}^2}$. By the repeated application of the lemma 2.3.3 we arrive at a bibasis $B = \{\alpha_1^1, \ldots, \alpha_{n_1}^1\} \cup \{\alpha_1^2, \ldots, \alpha_{n_2}^2\} = B_1 \cup B_2$ in which the neutrosophic bimatrix representing $T = T_1 \cup T_2$ is upper bitriangular

$$[T]_B = [T_1]_{B_1} \cup [T_2]_{B_2}$$

$$= \begin{bmatrix} a_{11}^1 & a_{12}^1 & \cdots & a_{1n_1}^1 \\ 0 & a_{22}^1 & \cdots & a_{2n_1}^1 \\ \vdots & \vdots & & \vdots \\ 0 & 0 & \cdots & a_{1n_1n_1}^1 \end{bmatrix} \cup \begin{bmatrix} a_{11}^2 & a_{12}^2 & \cdots & a_{1n_2}^2 \\ 0 & a_{22}^2 & \cdots & a_{2n_2}^2 \\ \vdots & \vdots & & \vdots \\ 0 & 0 & \cdots & a_{1n_2n_2}^2 \end{bmatrix}.$$

Merely $[T]_B$ = neutrosophic bitriangular bimatrix of $(n_1 \times n_1, n_2 \times n_2)$ order shows

$$\begin{aligned} T\alpha_j &= T_1 \alpha_{j_1}^1 \cup T_2 \alpha_{j_2}^2 \\ &= a_{1j_1}^1 \alpha_1^1 + \ldots + a_{j_1 j_1}^1 \alpha_{j_1}^1 \cup a_{1j_2}^2 \alpha_1^2 + \ldots + a_{j_2 j_2}^2 \alpha_{j_2}^2 \ldots \quad \text{(a)} \end{aligned}$$

$1 \leq j_i \leq n_i$, $i = 1, 2$ that is $T\alpha_j$ is in the strong neutrosophic bisubspace spanned by $\{\alpha_1^1, \ldots, \alpha_{j_1}^1\} \cup \{\alpha_1^2, \ldots, \alpha_{j_2}^2\}$. To find $\{\alpha_1^1, \ldots, \alpha_{j_1}^1\} \cup \{\alpha_1^2, \ldots, \alpha_{j_2}^2\}$ we start by applying the lemma to the bisubspace $W = W_1 \cup W_2 = \{0\} \cup \{0\}$ to obtain the bivector $\alpha_1^1 \cup \alpha_1^2$. Then apply lemma to $W_1^1 \cup W_1^2$ the bistrong neutrosophic bispace spanned by $\alpha^1 = \alpha_1^1 \cup \alpha_1^2$ and obtain $\alpha^2 = \alpha_2^1 \cup \alpha_2^2$. Next apply lemma to $W_2 = W_2^1 \cup W_2^2$. We have now obtained using the relation (a) the strong neutrosophic bivector space spanned by $\alpha^1$ and $\alpha^2$ and is biinvariant under T.

If T is bitriangulable then it is evident that the bicharacteristic neutrosophic bipolynomial for T has the form

$$\begin{aligned} f &= f_1 \cup f_2 \\ &= (x - c_1^1)^{d_1^1} \ldots (x - c_{k_1}^1)^{d_{k_1}^1} \cup (x - c_1^2)^{d_1^2} \ldots (x - c_{k_2}^2)^{d_{k_2}^2}. \end{aligned}$$

The bidiagonal entries $(a_{11}^1, \ldots, a_{1n_1}^1) \cup (a_{11}^2, \ldots, a_{1n_2}^2)$ are the bicharacteristic values with $c_j^t$ repeated $d_{jt}^t$ times. But if f can be bifactored so also is the biminimal bipolynomial p because p bidivides f.

The reader is expected prove the following corollary.

**COROLLARY 2.3.2:** *If $F = F_1 \cup F_2$ is a bialgebraically closed bifield. Every $(n_1 \times n_1, n_2 \times n_2)$ neutrosophic bimatrix over $F$ is similar over the bifield $F$ to a neutrosophic bitriangular bimatrix.*

**THEOREM 2.3.36:** *Let $V = V_1 \cup V_2$ be a $(n_1, n_2)$ bidimensional strong neutrosophic bivector space over the neutrosophic bifield $F = F_1 \cup F_2$ ($F_1$ and $F_2$ are not pure neutrosophic fields) and let $T = T_1 \cup T_2$ be a bilinear operator on $V = V_1 \cup V_2$. Then $T$ is bidiagonalizable if and only if the neutrosophic biminimal bipolynomial for T has the form*
$$p = p_1 \cup p_2$$
$$= (x - c_1^1) \ldots (x - c_{k_1}^1) \cup (x - c_1^2) \ldots (x - c_{k_2}^2).$$
*where $\{c_1^1, \ldots, c_{k_1}^1\} \cup \{c_1^2, \ldots, c_{k_2}^2\}$ are bidistinct element of $F = F_1 \cup F_2$.*

*Proof:* We know if $T = T_1 \cup T_2$ is bidiagonalizable its biminimal neutrosophic bipolynomial is a byproduct of bidistinct linear factors. Hence one way of the proof is clear.

To prove the converse let $W = W_1 \cup W_2$ be a strong neutrosophic bisubspace spanned by all the bicharacteritic bivectors of T and suppose $W \neq V$. Then we know by the properties of bilinear operator that there exists a bivector $\alpha = \alpha_1 \cup \alpha_2$ in $V = V_1 \cup V_2$ and not in $W = W_1 \cup W_2$ and the bicharacteristic value $c_j = c_{j_1}^1 \cup c_{j_2}^2$ of $T = T_1 \cup T_2$ such that the bivector

$\beta = (T - c_j I)\alpha$
$\quad = (T_1 - c_{j_1}^1 I_1)\alpha_1 \cup (T_2 - c_{j_2}^2 I_2)\alpha_2$
$\quad = \beta_1 \cup \beta_2$

lies in $W = W_1 \cup W_2$ where each $\beta_i \in W_i$, $i = 1, 2$. Since $\beta = \beta_1 \cup \beta_2$ is in W; $\beta_i = \beta_i^1 + \ldots + \beta_i^{k_i}$; $i = 1, 2$ with $T_i \beta_i^t = c_i^t \beta_i^t$, $t = 1, 2, \ldots, k_i$ and this is true for every $i = 1, 2$ and hence the bivector $h(T)\beta = \{h^1(c_1^1)\beta_1^1 + \ldots + h^1(c_{k_1}^1)\beta_{k_1}^1\} \cup \{h^2(c_1^2)\beta_1^2 + \ldots + h^2(c_{k_2}^2)\beta_{k_2}^2\}$ for every neutrosophic bipolyomial h. Now

$$\begin{aligned} p &= (x - c_j)q \\ &= p_1 \cup p_2 \\ &= (x - c_{j_1}^1)q_1 \cup (x - c_{j_2}^2)q_2 \end{aligned}$$

for some neutrosophic bipolynomial $q = q_1 \cup q_2$, Also $q - q(c_j) = (x - c_j)h$ that is

$$q_1 - q_1(c_{j_1}^1) \cup q_2 - q_2(c_{j_2}^2) = (x - c_{j_1}^1)h_{j_1}^1 \cup (x - c_{j_2}^2)h_{j_2}^2.$$

We have
$$\begin{aligned} q(t)\alpha - q(c_j)\alpha &= h(T)(T - c_j I)\alpha \\ &= h(T)\beta. \end{aligned}$$
But $h(T)\beta$ is in $W = W_1 \cup W_2$ and since
$$\begin{aligned} 0 &= p(T)\alpha \\ &= (T - c_j I)q(T)\alpha \\ &= p_1(T_1)\alpha_1 \cup p_2(T_2)\alpha_2 \\ &= (T_1 - c_{j_1}^1 I_1)q_1(T_1)\alpha_1 \cup (T_2 - c_{j_2}^2 I_2)q_2(T_2)\alpha_2 \end{aligned}$$

and the bivector $q(T)\alpha$ is in $W$, that is $q_1(T_1)\alpha_1 \cup q_2(T_2)\alpha_2$ is in $W = W_1 \cup W_2$. Therefore $q(c_j)\alpha = q_1(c_{j_1}^1)\alpha_1 \cup q_2(c_{j_2}^2)\alpha_2$ is in $W = W_1 \cup W_2$.

Since $\alpha = \alpha_1 \cup \alpha_2$ is not in $W = W_1 \cup W_2$, we have $q(c_j) = q_1(c_{j_1}^1) \cup q_2(c_{j_2}^2) = 0 \cup 0$. This contradicts the fact that $p$ has distinct roots.

Hence the claim.

We can now describe this more in terms of how the values are determined and its relation to Cayley Hamilton Theorem for strong neutrosophic bivector spaces of type II.

Suppose $T = T_1 \cup T_2$ is a bilinear operator of a strong neutrosophic bivector space of type II which is represented by the neutrosophic bimatrix $A = A_1 \cup A_2$ in some bibasis for which we wish to find out whether $T = T_1 \cup T_2$ is bidiagonalizable. We compute the bicharacteristic neutrosophic bipolynomial $f = f_1 \cup f_2$. If we can bifactor $f = f_1 \cup f_2$ as $(x - c_1^1)^{d_1^1} \ldots (x - c_{k_1}^1)^{d_{k_1}^1} \cup (x - c_1^2)^{d_1^2} \ldots (x - c_{k_2}^2)^{d_{k_2}^2}$, we have two different methods for finding whether or not $T$ is bidiagonlaizable. One method is to see whether for each i(t) (i(t) means i is independent on t) we can find a $d_i^t$ (t = 1, 2); $1 \le i \le k_t$ independent characteristic vectors associated with the characteristic value $c_i^t$. The other method is to check whether or not

$$(T - c_1 I) \cup (T - c_2 I) =$$
$$(T_1 - c_1^1 I_1) \ldots (T_1 - c_{k_1}^1 I_1) \cup (x - c_1^2 I_2) \ldots (x - c_{k_2}^2 I_2)$$

is the bizero operator.

**LEMMA 2.3.4:** *Let $V = V_1 \cup V_2$ be a strong neutrosophic bivector space of $(n_1, n_2)$ bidimensional over the neutrosophic bifield $F = F_1 \cup F_2$ ($F_1$ and $F_2$ are not neutrosophic pure) of type II.*

*Let $\{W_1^1, \ldots, W_{k_1}^1\} \cup \{W_1^2, \ldots, W_{k_2}^2\}$ be strong neutrosophic bivector subspace of V and let $W = W_1 \cup W_2 = \{W_1^1 + \ldots + W_{k_1}^1\} \cup \{W_1^2 + \ldots + W_{k_2}^2\}$. Then the following are equivalent.*

i. $\{W_1^1, \ldots, W_{k_1}^1\} \cup \{W_1^2, \ldots, W_{k_2}^2\}$ *are biindependent, that is* $\{W_1^t, \ldots, W_{k_t}^t\}$ *are independent for t = 1, 2.*

ii. *For each $j_t$, $2 \le j_t \le k_t$; t = 1, 2, we have* $W_{j_t}^t \{W_1^t + \ldots + W_{j_t-1}^t\} = \{0\}$ *for t = 1, 2.*

iii. *If $B_i^t$ is a bibasis for $W_i^t$, $1 \le i \le k_t$, t = 1, 2; then the bisequences $\{B_1^1, \ldots, B_{k_1}^1\} \cup \{B_1^2, \ldots, B_{k_2}^2\}$ is a bibasis for the strong neutrosophic bisubspace $W = W_1 \cup W_2 = \{W_1^1 + \ldots + W_{k_1}^1\} \cup \{W_1^2 + \ldots + W_{k_2}^2\}$.*

*Proof:* Assume (i) let $\alpha^t \in W_{j_t}^t \cap \{W_1^t + \ldots + W_{j_t-1}^t\}$ then there are vectors $(\alpha_1^t, \ldots, \alpha_{j_t-1}^t)$ with $\alpha_i^t \in W_i^t$ such that $(\alpha_1^t + \ldots + \alpha_{j_t-1}^t + \alpha^t) + \alpha^t = 0 + \ldots + 0 = 0$ and since $\{W_1^t, \ldots, W_{k_t}^t\}$ are independent it must be that $\alpha_1^t = \alpha_2^t = \ldots = \alpha_{j_t-1}^t = 0$. This is true for each t; t = 1, 2. Now let us observe that (ii) implies (i). Suppose $0 = \left(\alpha_1^t + \ldots + \alpha_{k_t}^t\right)$; $\alpha_i^t \in W_i^t$, i = 1, 2, ..., $k_t$ (we denote both the zero vector and zero scalar by 0). Let $j_t$ be the largest integer $i_t$ such that $\alpha_i^t \neq 0$. Then $0 = \alpha_1^t + \ldots + \alpha_{j_t}^t$; $\alpha_{j_t}^t \neq 0$ thus $\alpha_{j_t}^t = -\alpha_1^t - \ldots - \alpha_{j_{t-1}}^t$ is a non zero vector in $W_{j_t}^t \cap \{W_1^t + \ldots + W_{j_{t-1}}^t\}$.

Now that we know (i) and (ii) are the same let us see why (i) is equivalent to (iii). Assume (i). Let $B_i^t$ be a basis for $W_i^t$; $1 \leq i \leq k_t$ and let $B^t = \{B_1^t, \ldots, B_{k_t}^t\}$ true for each t, t = 1, 2.

Any linear relation between the vector in $B^t$ will have the form $\left(\beta_1^t + \ldots + \beta_{k_t}^t\right) = 0$ where $\beta_i^t$ is some linear combination of vectors in $B_i^t$. Since $\{W_1^t, W_2^t, \ldots, W_{k_t}^t\}$ are independent each of $\beta_i^t$ is 0. Since each $B_i^t$ is an independent relation.

The relation between vectors in $B^t$ is trivial. This is true for every t; t = 1, 2; so in $B = B^1 \cup B^2 = \{B_1^1, \ldots, B_{k_1}^1\} \cup \{B_1^2, \ldots, B_{k_2}^2\}$ every birelation in bivector in B is a trivial birelation. It is left for the reader to prove(c) implies (a).

If any of the conditions of the above lemma hold we say the bisum $W = \{W_1^1 + \ldots + W_{k_1}^1\} \cup \{W_1^2 + \ldots + W_{k_2}^2\}$; bidirect or that W is the bidirect sum of $\{W_1^1, \ldots, W_{k_1}^1\} \cup \{W_1^2 \ldots W_{k_2}^2\}$ and we write $W = \{W_1^1 \oplus \ldots \oplus W_{k_1}^1\} \cup \{W_1^2 \oplus \ldots \oplus W_{k_2}^2\}$. This bidirect sum will also be known as the biindependent sum or the biinterior direct sum of $\{W_1^1, \ldots, W_{k_1}^1\} \cup \{W_1^2, \ldots, W_{k_2}^2\}$. Let $V = V_1 \cup V_2$ be a strong neutrosophic bivector space over the neutrosophic bifield $F = F_1 \cup F_2$. A biprojection of V is a bilinear operator $E = E_1 \cup E_2$ on V such that $E^2 = E_1^2 \cup E_2^2 = E_1 \cup E_2 = E$.

Since E is a biprojetion. Let $R = R_1 \cup R_2$ be the birange of E and let $N = N_1 \cup N_2$ be the null bispace $E = E_1 \cup E_2$.

1. The bivector $\beta = \beta_1 \cup \beta_2$ is the birange R if and only if $E\beta = \beta$ that is $E_1 \beta_1 \cup E_2 \beta_2 = \beta_1 \cup \beta_2$. If $\beta = E\alpha$ then $E\beta = E^2\alpha = E\alpha = \beta$. Conversely if $\beta = E\beta$ then of course $\beta$ is in the birange of $E = E_1 \cup E_2$.
2. $V = R \oplus N$; that is $V = V_1 \cup V_2 = R_1 \oplus N_1 \cup R_2 \oplus N_2$; that is each $V_i = R_i \oplus N_i$; $i = 1, 2$.
3. The unique expression for $\alpha$ as a sum of bivector in R and N is $\alpha = E\alpha + (\alpha - E\alpha)$ that is $\alpha = \alpha_1 \cup \alpha_2 = E_1\alpha_1 + (\alpha_1 - E_1\alpha_1) \cup E_2\alpha_2 + (\alpha_2 - E_2\alpha_2)$.

From (1), (2) and (3) it is easy to verify, if $R = R_1 \cup R_2$ and $N = N_1 \cup N_2$ are strong neutrosophic bivector subspace of V such that $V = R \oplus N = R_1 \oplus N_1 \cup R_2 \oplus N_2$, there is one and only one biprojection operator $E = E_1 \cup E_2$ which has birange R and binull space N. That operator is called the biprojection on R along N.

Any biprojeciton $E = E_1 \cup E_2$ is trivially bidiagonalizable. If $\{\alpha_1^1,...,\alpha_{r_1}^1\} \cup \{\alpha_1^2,...,\alpha_{r_2}^2\}$ a bibasis of R and $\{\alpha_{r_1+1}^1,...,\alpha_{n_1}^1\} \cup \{\alpha_{r_2+1}^2,...,\alpha_{n_2}^2\}$ a bibasis for N then the bibasis $B = \{\alpha_1^1,...,\alpha_{n_1}^1\} \cup \{\alpha_1^2,...,\alpha_{n_2}^2\} = B_1 \cup B_2$, bidiagonalizes $E = E_1 \cup E_2$.

$$[E]_B = [E_1]_{B_1} \cup [E_2]_{B_2} = \begin{bmatrix} I_1 & 0 \\ 0 & 0 \end{bmatrix} \cup \begin{bmatrix} I_2 & 0 \\ 0 & 0 \end{bmatrix}$$

where $I_1$ is a $r_1 \times r_1$ identity matrix and $I_2$ is a $r_2 \times r_2$ identity matrix.

Projections can be used to describe the bidirect sum decomposition of the strong neutrosophic bivector space $V = V_1 \cup V_2$. For suppose $V = \{W_1^1 \oplus ... \oplus W_{k_1}^1\} \cup \{W_1^2 \oplus ... \oplus W_{k_2}^2\}$ for each j(t) we define $E_j^t$ on $V_t$. (t = 1, 2). Let $\alpha = \alpha_1 \cup \alpha_2$ be in $V = V_1 \cup V_2$ say $\alpha = \{\alpha_1^1 + ... + \alpha_{k_1}^1\} \cup \{\alpha_1^2 + ... + \alpha_{k_2}^2\}$ with $\alpha_i^t$ in $W_i^t$, $1 \leq i \leq k_t$ for t = 1, 2. Define $E_j^i \alpha^t = \alpha_j^t$ then $E_j^t$ is a well defined rule. It is easy to see that $E_i^t$ is linear and that range of $E_j^t$ is $W_j^t$ and $(E_j^t)^2 = E_j^t$. The null space of $E_j^t$ is the strong neutrosophic subspace $W_1^t + ... + W_{j-1}^t + W_{j+1}^t + ... + W_{k_t}^t$ for the statement $E_j^t \alpha^t = 0$ simply means $\alpha_j^t = 0$ that is $\alpha$ is actually a sum of vectors from the spaces $W_i^t$ with $i \neq j$. Interms of the projections $E_j^t$ we have $\alpha^t = E_1^t \alpha^t + ... + E_k^t \alpha^t$ for each $\alpha$ in V. The above equation implies $I_t = \{E_1^t + ... + E_{k_t}^t\}$. Note also that if $i \neq j$ then $E_i^t E_j^t = 0$ because the range of $E_j^t$ is the strong neutrosophic subspace $W_j^t$ which is contained in the null space of $E_i^t$. This is true for each t, t = 1, 2. Hence true on the strong neutrosophic bivector space
$$V = \{W_1^1 \oplus ... \oplus W_{k_1}^1\} \cup \{W_1^2 \oplus ... \oplus W_{k_2}^2\}$$

Now we prove an interesting result.

**THEOREM 2.3.37:** *Let $V = V_1 \cup V_2$ be a strong neutrosophic bivector space over the neutrosophic bifield $F = F_1 \cup F_2$ (both $F_1$ and $F_2$ are not pure neutrosophic) of type II and suppose $V = \{W_1^1 \oplus ... \oplus W_{k_1}^1\} \cup \{W_1^2 \oplus ... \oplus W_{k_2}^2\}$ then there exists $(k_1, k_2)$ bilinear operators $\{E_1^1, ..., E_{k_1}^1\} \cup \{E_1^2, ..., E_{k_2}^2\}$ on V such that*

  i. *Each $E_i^t$ is a projection, that is $(E_i^t)^2 = E_i^t$ for $t = 1, 2$; $1 \leq i \leq k_t$.*
  ii. *$E_i^t E_j^t = 0$ if $i \neq j$; $1 \leq i, j \leq k_t$ and $t = 1, 2$.*
  iii. *$I = I_1 \cup I_2 = \{E_1^1 + ... + E_{k_1}^1\} \cup \{E_1^2 + ... + E_{k_2}^2\}$.*
  iv. *The range of $E_i^t$ is $W_i^t$, for $i = 1, 2, ..., k_t$ and $t = 1, 2$.*

*Proof:* We are primarily interested in the bidirect sum bidecomposition $V = \{W_1^1 \oplus ... \oplus W_{k_1}^1\} \cup \{W_1^2 \oplus ... \oplus W_{k_2}^2\} = W_1 \cup W_2$ where each of the strong neutrosophic bivector bisubspaces $W_t$ is invariant under some given bilinear operator $T = T_1 \cup T_2$.

Given such a decomposition of V, T induces bilinear operators $T = T_1 \cup T_2$.

Given such a bidecomposition of $V = V_1 \cup V_2$, $T = T_1 \cup T_2$ induces bilinear operators $\left(T_i^1 \cup T_i^2\right)$ on $\left(W_i^1 \cup W_i^2\right)$ by restriction, the action of T is $\alpha$, is a bivector in V we have unique bivectors $\{\alpha_1^1,...,\alpha_{k_1}^1\} \cup \{\alpha_1^2,...,\alpha_{k_2}^2\}$ with $\alpha_i^t$ in $W_i^t$ such that $\alpha = \{\alpha_1^1+...+\alpha_{k_1}^1\} \cup \{\alpha_1^2+...+\alpha_{k_2}^2\}$ and then

$$T\alpha = \{T_1^1\alpha_1^1+...+T_{k_1}^1\alpha_{k_1}^1\} \cup \{T_1^2\alpha_1^2+...+T_{k_2}^2\alpha_{k_2}^2\}.$$

We shall describe this situation by saying that $T = T_1 \cup T_2$ is the bidirect sum of the operators $\{T_1^1,...,T_{k_1}^1\} \cup \{T_1^2,...,T_{k_2}^2\}$.

It must be remembered in using this terminology that the $T_i^t$ are not bilinear operators on the strong neutrosophic bivector space $V = V_1 \cup V_2$ but on the various strong neutrosophic bivector subspaces

$$W = W_1 \cup W_2.$$
$$= \{W_1^1 \oplus ... \oplus W_{k_1}^1\} \cup \{W_1^2 \oplus ... \oplus W_{k_2}^2\}$$

which enables us to associate with each $\alpha = \alpha_1 \cup \alpha_2$ in V a unique pair of k-tuple $\{\alpha_1^1,...,\alpha_{k_1}^1\} \cup \{\alpha_1^2,...,\alpha_{k_2}^2\}$ of vectors $\alpha_i^t \in W_i^t$, $i = 1, 2, ..., k_t$. $t = 1, 2$.

$\alpha = \{\alpha_1^1+...+\alpha_{k_1}^1\} \cup \{\alpha_1^2+...+\alpha_{k_2}^2\}$ is in such a way that we can carry out the bilinear operators on V by working in the individual strong neutrosophic bivector subspaces $W_i = W_i^1 + W_i^2$. The fact that each $W_i$ is biinvariant under T enable us to view the action of T as the independent action of the operators $T_i^t$ on the bisubspaces $W_i^t$; $i = 1, 2, ..., k_t$, $t = 1, 2$. Our purpose is to study T by finding biinvariant bidirect sum decompositions in which the $T_i^t$ operators of an elementary nature.

The following theorem is left as an exercise for the reader to prove.

**THEOREM 2.3.38:** *Let $T = T_1 \cup T_2$ be a bilinear operator on a strong neutrosophic bivector space $V = V_1 \cup V_2$ of type II over the neutronsophic bifield $F = F_1 \cup F_2$ ($F_1$ and $F_2$ are not pure neutrosophic fields). Let $(W_1^1, \ldots, W_{k_1}^1) \cup (W_1^2, \ldots, W_{k_2}^2)$ and $(E_1^1, \ldots, E_{k_1}^1) \cup (E_1^2, \ldots, E_{k_2}^2)$ be as before. Then a necessary and sufficient condition that each strong neutrosophic bivector subspace $W_i^t$ to be biinvariant under $T_i$ for $1 \le i \le k_t$; $t = 1, 2$ is that $E_i^t T_t = T_t E_i^t$ or $ET = TE$ for every $1 \le i \le k_t$ and $t = 1, 2$.*

Now we proceed on to define the notion of biprimary decomposition of strong neutrosophic bivector space $V = V_1 \cup V_2$ of $(n_1, n_2)$ dimension over the neutrosophic bifield $F = F_1 \cup F_2$ where $F_1$ and $F_2$ are not pure neutrosophic fields.

**THEOREM 2.3.39:** (Primary bidecomposition theorem): *Let $T = T_1 \cup T_2$ be a bilinear operator on a finite $(n_1, n_2)$ dimension strong neutrosophic bivector space $V = V_1 \cup V_2$ over the neutroscophic bifield $F = F_1 \cup F_2$ ($F_1$ and $F_2$ are not pure neutrosophic fields). Let $p = p_1 \cup p_2$ be the biminimal neutrosophic polynomial for $T = T_1 \cup T_2$.*

$p = p_{11}^{r_1^1} \ldots p_{1k_1}^{r_{k_1}^1} \cup p = p_{21}^{r_1^1} \ldots p_{2k_1}^{r_{k_1}^1}$ *where $p_{ti}^t$ are distinct irreducible monic neutrosophic polynomials over $F_t$; $i = 1, 2, \ldots, k_t$; $t = 1, 2$ and $r_i^t$ are positive integers. Let $W_i = W_i^1 \cup W_i^2$ be the null bispace of $p(T) = p_{1i}(T_i^1)^{r_i^1} \cup p_{2i}(T_i^2)^{r_i^2}$; $i = 1, 2$; then,*

  i. $W = W_1 \cup W_2 = (W_1^1 \oplus \ldots \oplus W_{k_1}^1) \cup (W_1^2 \oplus \ldots \oplus W_{k_2}^2)$
  ii. *each $W_i = W_i^1 + W_i^2$ is biinvariant under $T_i$; $i = 1, 2$.*
  iii. *If $T_i^r$ is the operator induced on $W_i^r$ by $T_i$ then the minimal neutrosophic polynomial for $T_i^r$ is $p_i^r$; $r = 1, 2, \ldots, k_i$, $i = 1, 2$.*

We prove the corollary to this theorem.

**COROLLARY 2.3.3:** *If $\{E_1^1,\ldots,E_{k_1}^1\} \cup \{E_1^2,\ldots,E_{k_2}^2\}$ are biprojections associated with the biprimary decomposition of T = $T_1 \cup T_2$ then each $E_i^t$ is a neutrosophic polynomial in T; $1 \le i \le k_t$; t = 1, 2 and accordingly if a linear operator S commutes with T then S commutes with each of the $E_i$; that is each strong neutrosophic subspace $W_i$ is invariant under S.*

*Proof:* For any bilinear operator defined on a strong neutrosophic bivector space defined over the neutrosophic bifield $F = F_1 \cup F_2$ ($F_1$ and $F_2$ pure neutrosophic fields) of type II, we can associate the notion of bidiagonal part of T and binilpotent part of T.

Consider the neutrosophic bimininal polynomial for $T = T_1 \cup T_2$ which is decomposed into first degree polynomials that is the case in which each $p_i$ is of the form $p_i^t = x - c_i^t$. Now the range of $E_i^t$ is the null space of $W_i^t$ of $(T_t - c_i^t I_t)^{r_i^t}$; we know by earlier results D is a bidiagonalizable part of T.

Let us look at the bioperator
$$N = T - D$$
$$N_1 \cup N_2 = (T_1 - D_1) \cup (T_2 - D_2)$$
$$T = (T_1 E_1^1 + \ldots + T_1 E_{k_1}^1) \cup (T_2 E_1^2 + \ldots + T_2 E_{k_2}^2)$$

and

$$D = D_1 \cup D_2 = (c_1 E_1^1 + \ldots + c_{k_1}^1 E_{k_1}^1) \cup (c_2 E_1^2 + \ldots + c_{k_2}^2 E_{k_2}^2)$$

so

$$N = N_1 \cup N_2$$
$$= \{(T_1 - c_1^1 I_1)E_1^1 + \ldots + (T_1 - c_{k_1}^1 I_1)E_{k_1}^1\} \cup$$
$$\{(T_2 - c_1^2 I_2)E_1^2 + \ldots + (T_2 - c_{k_2}^2 I_2)E_{k_2}^2\}.$$

Now

$$N^2 = N_1^2 \cup N_2^2$$
$$= (T_1 - c_1^1 I_1)^2 E_1^1 + \ldots + (T_1 - c_{k_1}^1 I_1)^2 E_{k_1}^1 \cup$$
$$(T_2 - c_1^2 I_2)^2 E_1^2 + \ldots + (T_2 - c_{k_2}^2 I_2)^2 E_{k_2}^2.$$

and in general

$$N^r = \{(T_1 - c_1^1 I_1)^{r_1} E_1^1 + \ldots + (T_1 - c_{k_1}^1 I_1)^{r_{k_1}} E_{k_1}^1\} \cup$$

$$\{(T_2 - c_1^2 I_2)^{r_2} E_1^2 + \ldots + (T_2 - c_{k_2}^2 I_2)^{r_{k_2}} E_{k_2}^2\}.$$

When $r \geq r_i$ for each i we have $N^r = 0$ i.e., $N_1^r \cup N_2^r = 0 \cup 0$; that is each of the bioperator $(T_t - c_1^t I_t)^{r_t} = 0$ on the range $E_i^t$; $1 \leq t \leq k_i$ and $i = 1, 2$. Thus $(T - cI)^r = 0$ for a suitable r.

Let $N = N_1 \cup N_2$ be a bilinear operator on a bivector space $V = V_1 \cup V_2$. We say N is binilpotent if there is some pair of integers $(r_1, r_2)$ such that $N_i^{r_i} = 0$ for $i = 1, 2$. We choose $r > r_i$; $i = 1, 2$ then $N^r = 0$, where $N = N_1 \cup N_2$.

In view of this we have the following theorem for strong neutrosophic bivector spaces of type II defined over the neutrosophic bifield $F = F_1 \cup F_2$ ($F_1$ and $F_2$ not pure neutrosophic).

**THEOREM 2.3.40:** *Let $T = T_1 \cup T_2$ be a bilinear operator on the $(n_1, n_2)$ finite bidimensional strong neutrosophic bivector space $V = V_1 \cup V_2$ over the neutrosophic bifield $F = F_1 \cup F_2$ (both $F_1$ and $F_2$ are not pure). Suppose that the biminimal neutrosophic polynomial for $T = T_1 \cup T_2$ decomposes over $F = F_1 \cup F_2$ into a biproduct of bilinear neutrosophic polynomials. Then there is a bi diagonalizable operator $N = N_1 \cup N_2$ on $V = V_1 \cup V_2$ such that,*
  *i. $T = D + N$ ; i.e.;*
    $T_1 \cup T_2 = D_1 \cup D_2 + N_1 \cup N_2$
    $= D_1 + N_1 \cup D_2 + N_2.$
  *ii. $DN = ND$ that is*
    $(D_1 \cup D_2)(N_1 \cup N_2) = D_1 N_1 \cup D_2 N_2$
    $= N_1 D_1 \cup N_2 D_2.$
*The bidiagonalizable operator $D = D_1 \cup D_2$ and the binilpotent operator $N = N_1 \cup N_2$ are uniquely determined by (i) and (ii) and each of them is a bipolynomial in $T_1$ and $T_2$.*

Consequent of the above theorem the following corollary is direct.

**COROLLARY 2.3.4:** *Let V be a finite bidimension strong neutrosophic bivector space over the special algebraically closed neutrosophic bifield $F = F_1 \cup F_2$. Then every bilinear operator $T = T_1 \cup T_2$ on $V = V_1 \cup V_2$ can be written as the sum of a bidiagonalizable operator $D = D_1 \cup D_2$ and a binilpotent operator $N = N_1 \cup N_2$ which commute. These bioperators D and N are unique and each is a bipolynomial in $(T_1, T_2)$.*

*Let $V = V_1 \cup V_2$ be a finite bidimensional strong neutrosophic bivector space over the neutrosophic bifield $F = F_1 \cup F_2$ and $T = T_1 \cup T_2$ be an aribitary and fixed bilinear operator on $V = V_1 \cup V_2$. If $\alpha = \alpha_1 \cup \alpha_2$ is a bivector in V then there is a smallest bisubspace of $V = V_1 \cup V_2$ which is biinvarient under $T = T_1 \cup T_2$ and contains $\alpha$. This strong neutrosophic bispace can be defined as the biintersection of all T-invariant strong neutrosophic bisubspaces which contain $\alpha$.*

If $W = W_1 \cup W_2$ be any strong neutrosophic bivector supspace of $W = W_1 \cup W_2$ of a strong neutrosophic bivector space $V = V_1 \cup V_2$ which is biinvariant under $\alpha = \alpha_1 \cup \alpha_2$; that is each $T_i$ in T is such that the strong neutrosophic subspace $W_i$ on $V_i$ is invariant under $T_i$ and contains $\alpha_i$; true for i = 1, 2.

Then $W = W_1 \cup W_2$ must also contain $T\alpha$; that is $T_i\alpha_i$ is in $W_i$ for each i = 1, 2; hence $T(T\alpha)$ is in W; that is $T_i(T_i\alpha_i) = T_i^2\alpha_i$ is in W and so on; that is $T_i^{m_i}(\alpha_i)$ is in $W_i$, for each i so that $T^m(\alpha) \in W$; i = 1, 2. W must contain $g(T)\alpha$ for every neutrosophic bipolynomial $g = g_1 \cup g_2$ over the neutrosophic bifield $F = F_1 \cup F_2$. The set of all bivectors of the form $g(T)\alpha = g_1(T_1)\alpha_1 \cup g_2(T_2)\alpha_2$ with $g = g_1 \cup g_2 \in F[x] = F_1[x] \cup F_2[x]$, is clearly biinvariant and is thus the smallest bi T-invariant (T-biinvariant) strong neutrosophic bisubspace which contains $\alpha = \alpha_1 \cup \alpha_2$.

In view of this we have the following definition.

**DEFINITION 2.3.46:** *Let $\alpha = \alpha_1 \cup \alpha_2$ be any bivector in a strong neutrosophic bivector space $V = V_1 \cup V_2$ over the neutrosophic bifield $F = F_1 \cup F_2$ ($F_1$ and $F_2$ are not pure neutrosophic). The*

*T-bicyclic strong neutrosophic bisubspace generated by $\alpha = \alpha_1 \cup \alpha_2$ is a strong neutrosophic bisubspace $Z(\alpha; T) = Z(\alpha_1; T_1) \cup Z(\alpha_2; T_2)$ of all bivectors $g(T)\alpha = g_1(T_1)\alpha_1 \cup g_2(T_2)\alpha_2$; $g = g_1 \cup g_2$ in $F[x] = F_1[x] \cup F_2[x]$ is a neutrosophic bipolynomial. If $Z(\alpha; T) = V$ then $\alpha$ is a bicyclic vector for T.*

Another way of describing this strong neutrosophic bisubspace $Z(\alpha; T)$ is that $Z(\alpha; T)$ is the strong neutrosophic bisubspace spanned by the bivectors $T_\alpha^k$; $k \geq 0$ and $\alpha$ is a bicyclic bivector for $T = T_1 \cup T_2$ if and only if these bivectors span V; that is each $T_{i\alpha_i}^{k_i}$ span $V_i$, $k_i \geq 0$ and thus $\alpha_i$ is a cyclic vector for $T_i$ if and only if these vectors span $V_i$, true for i = 1, 2.

We just caution the reader that the general bioperator $T = T_1 \cup T_2$ has no bicyclic bivector.

For any T the T bicyclic strong neutrosophic bisubspace generated by the bizero vector is the bizero strong neutrosophic bisubspace of V. The bispace $Z(\alpha;T) = Z(\alpha_1;T_1) \cup Z(\alpha_2;T_2)$ is (1, 1) dimensional if and only if $\alpha$ is a bicharacteristic vector for T. For the biidentity operator, every nonzero bivector generates a (1, 1) dimensional bicyclic strong neutrosophic bisubspace thus if bidimV > (1, 1) the biidentity operator has non cyclic vector.

For any T and $\alpha$ we shall be interested in the bilinear relation $c_0\alpha + c_1T\alpha + \ldots + c_kT\alpha^k = 0$ where $\alpha = \alpha_1 \cup \alpha_2$ so that
$$c_0^1\alpha_1 + c_1^1T_1\alpha_1 + \ldots + c_{k_1}^1 T_1\alpha_1^{k_1} = 0$$
and
$$c_0^2\alpha_2 + c_1^2T_2\alpha_2 + \ldots + c_{k_2}^2 T_2\alpha_2^{k_2} = 0$$
between the bivectors $T\alpha^j$, we shall be interested in the neutrosophic bipolynomial $g = g_1 \cup g_2$ where
$$g_i = c_0^i + c_1^i x + \ldots + c_{k_i}^i x^{k_i}$$
true for i = 1, 2, which has the property that $g(T)\alpha = 0$.

The set of all g in $F[x] = F_1[x] \cup F_2[x]$ such that $g(T)\alpha = 0$ is clearly a neutrosophic biideal in $F[x]$. It is also a non zero neutrosophic biideal in $F[x]$ because it contains biminimal bipolynomial $p = p_1 \cup p_2$ of the bioperator T. $p(T)\alpha = p_1(T_1)\alpha_1$ $\cup\ p_2(T_2)\alpha_2 = 0 \cup 0$; that is $p(T_1)\alpha_1 \cup p_2(T_2)\alpha_2 = 0 \cup 0$ for every $\alpha = \alpha_1 \cup \alpha_2$ in $V = V_1 \cup V_2$.

**DEFINITION 2.3.47:** *If $\alpha = \alpha_1 \cup \alpha_2$ is any bivector in strong neutrosophic bivector space $V = V_1 \cup V_2$ of type II defined over the neutrosophic bifield $F = F_1 \cup F_2$ ($F_1$ and $F_2$ are not pure neutrosophic). The T-annihilator ($T = T_1 \cup T_2$) of $\alpha = \alpha_1 \cup \alpha_2$ is the neutrosophic biideal $M(\alpha; T) = M(\alpha_1; T_1) \cup M(\alpha_2; T_2)$ in $F[x] = F_1[x] \cup F_2[x]$ consisting of all neutrosophic bipolynomials $g = g_1 \cup g_2$ over $F = F_1 \cup F_2$ such that $g(T) = g_1(T_1) \cup g_2(T_2) = 0 \cup 0$.*

*The unique neutrosophic monic bipolynomial $p\alpha = p_1\alpha_1 \cup p_2\alpha_2$ which bigenerates this biideal will also be called the bi T-annihilator of $\alpha$ or T bi annihilator of $\alpha$. The bi T-annihilator $p\alpha$ bidivides the neutrosophic biminimal bipolynomial of the bioperator $T = T_1 \cup T_2$. Clearly bidegree $(p\alpha) > (0, 0)$ unless $\alpha = \alpha_1 \cup \alpha_2$ is the zero bivector.*

**THEOREM 2.3.41:** *Let $\alpha = \alpha_1 \cup \alpha_2$ be any non zero bivector in $V = V_1 \cup V_2$; V a strong neutrosophic bivector space over the neutrosophic bifield $F = F_1 \cup F_2$ (both $F_2$ and $F_2$ are not pure neutrosophic).*

*Let $p_\alpha = p_{1\alpha_1} \cup p_{2\alpha_2}$ be the bi T annihilator of $\alpha = \alpha_1 \cup \alpha_2$.*

  i. *The bidegree of $p_\alpha$ is equal to the bidimension of the bicyclic strong neutrosophic bisubspace $Z(\alpha;T) = Z(\alpha_1;T_1) \cup Z(\alpha_2;T_2)$.*
  
  ii. *If the bidegree of $p_\alpha = p_{1\alpha_1} \cup p_{2\alpha_2}$ is $(k_1, k_2)$ then the bivectors $\alpha = \alpha_1 \cup \alpha_2$, $T\alpha = T_1\alpha_1 \cup T_2\alpha_2$, ..., $T_{\alpha_1}^{k_1-1} = T_1^{k_1-1}\alpha_1 \cup T_2^{k_2-1}\alpha_2$ form a bibasis for $Z(\alpha;T)$. That is $\{\alpha_1, T_1\alpha_1, T_1^2\alpha_1, ..., T_1^{k_1-1}\alpha_1\} \cup \{\alpha_2, T_2\alpha_2, T_2^2\alpha_2, ..., T_2^{k_2-1}\alpha_2\}$ form a bibasis for $Z(\alpha;T) = Z(\alpha_1;T_1) \cup Z(\alpha_2;T_2)$; that is $Z(\alpha_i, T_i)$ has $\{\alpha_i, T_i\alpha_i, ..., T_i^{k_i-1}\alpha_i\}$ as its basis; true for every $i = 1, 2$.*

iii.  If $S = S_1 \cup S_2$ is a bilinear operator on $Z(\alpha;T)$ induced by $T$, then by the biminimal neutrosophic polynomial for $S$ is $p_\alpha$.

*Proof:* Let $g = g_1 \cup g_2$ be a neutrosophic bipolynomial over the neutrosophic bifield $F = F_1 \cup F_2$, write $g = p_\alpha q + r$, that is $g_1 \cup g_2 = p_{1\alpha_1} q_1 + r_1 \cup p_{2\alpha_2} q_2 + r_2$; where $p_\alpha = p_{1\alpha_1} \cup p_{2\alpha_2}$ for $\alpha = \alpha_1 \cup \alpha_2$, $q = q_1 \cup q_2$ and $r = r_1 \cup r_2$ so $g_i = p_{i\alpha_i} q_i + r_i$ truce for $i = 1, 2$. Here either $r = 0 \cup 0$ or bidegree $r <$ bidegree $p_\alpha = (k_1, k_2)$. The neutrosophic bipolynomial $p_\alpha q = p_{1\alpha_1} q_1 \cup p_{2\alpha_2} q_2$ is in the T biannihilator of $\alpha = \alpha_1 \cup \alpha_2$ and so $g(T)\alpha = r(T)\alpha$, that is $g(T_1)\alpha_1 \cup g_2(T_2)\alpha_2 = r_1(t_1)\alpha_1 \cup r_2(T_2)\alpha_2$.

Since $r = r_1 \cup r_2 = 0 \cup 0$ or bidegree $r < (k_1, k_2)$ the bivector $r(T)\alpha = r_1 T_1(\alpha_1) \cup r_2 T_2(\alpha_2)$ is a bilinear combination of the bivectors $\alpha, T\alpha,\ldots,T^{k-1}\alpha$; that is a bilinear combination of bivectors $\alpha = \alpha_1 \cup \alpha_2$

$$T\alpha = T_1\alpha_1 \cup T_2\alpha_2,$$
$$T^2\alpha = T_1^2\alpha_1 \cup T_2^2\alpha_2,$$
$$T^3\alpha = T_1^3\alpha_1 \cup T_2^3\alpha_2, \ldots,$$
$$T^{k-1}\alpha = T_1^{k_1-1}\alpha_1 \cup T_2^{k_2-1}\alpha_2$$

and since $g(T)\alpha = g_1(T_1)\alpha_1 \cup g_2(T_2)\alpha_2$ is a typical bivector in $Z(\alpha; T)$; i.e., each $g_i(T_i)\alpha_i$ is a typical vector in $Z(\alpha_i; T_i)$ for each $g_i(T_i)\alpha_i$ is a typical vector in $Z(\alpha_I; T_i)$; $i = 1, 2$. This shows that these $(k_1, k_2)$ bivectors span $Z(\alpha; T)$.

These bivectors are certainly bilinearly independent because any non trivial bilinear relation between them would give us a non zero neutrosophic bipolynomial $g = g_1 \cup g_2$ such that $g(T)(\alpha) = g_1(T_1)(\alpha_1) \cup g_2(T_2)\alpha_2 = 0 \cup 0$ and bidegree $g <$ bidegree $p_\alpha$, which is absurd.

This proves (i) and (ii).

Let $S = S_1 \cup S_2$ be a bilinear operator on $(Z_\alpha; T)$ obtained by restricting $T$ to that strong neutrosophic bivector subspace. If $g = g_1 \cup g_2$ is any neutrosophic bifield $F = F_1 \cup F_2$ then $p_\alpha(S)g(T)\alpha = p_\alpha(T)g(T)\alpha$

that is

$$p_{1\alpha_1}(S_1)g_1(T_1)\alpha_1 \cup p_{2\alpha_2}(S_2)g_2(T_2)\alpha_2$$
$$= p_{1\alpha_1}(T_1)g_1(T_1)\alpha_1 \cup p_{2\alpha_2}(T_2)g_2(T_2)\alpha_2$$
$$= g(T)p_\alpha(T)\alpha$$
$$= g_1(T_1)p_{1\alpha_1}(T_1)\alpha_1 \cup g_2(T_2)p_{2\alpha_2}(T_2)\alpha_2$$
$$= g_1(T_1)(0) \cup g_2(T_2)(0)$$
$$= 0 \cup 0.$$

Thus the bioperator $p_\alpha S = p_1\alpha_1(S_1) \cup p_2\alpha_2(S_2)$ sends every bivector in $Z(\alpha;T) = Z(\alpha_1;T_1) \cup Z(\alpha_2; T_2)$ into $0 \cup 0$ and is the bizero operator on $Z(\alpha;T)$. Further more if $h = h_1 \cup h_2$ is a neutrosophic of bidegree less than $(k_1, k_2)$ we cannot have $h(S) = h_1(S_1) \cup h_2(S_2) = 0 \cup 0$ for then $h(S)\alpha = h_1(S_1)\alpha_1 \cup h_2(S_2)\alpha_2 = 0 \cup 0$; contradicting the definition of $p_\alpha$. This shows that $p_\alpha$ is the neutrosophic biminimal polynomial for S.

A particular consequences of this interesting theorem is that if $\alpha = \alpha_1 \cup \alpha_2$ happens to be a bicyclic vector for $T = T_1 \cup T_2$ then the neutrosophic biminimal bipolynomial for T have bidegree equal to the bidimension of the strong neutrosophic bivector space $V = V_1 \cup V_2$, hence by Cayley Hamilton theorem for the bivector spaces we have the neutrosophic biminimal polynomial for T is the bicharacteristic neutrosophic polynomial for T. We shall prove later that for any T there is a bivector $\alpha = \alpha_1 \cup \alpha_2$ in $V = V_1 \cup V_2$ which has the neutrosophic biminimal polynomial for $T = T_1 \cup T_2$ for its biannihilator.

It will then follow that $T = T_1 \cup T_2$ has a bicyclic vector if and only if the biminimal and the bicharacteristic neutrosophic polynomial for T are identical. We now study the general bioperator $T = T_1 \cup T_2$ by using the bioperator vector. Let us consider a bilinear operator $S = S_1 \cup S_2$ on the strong neutrosophic bivector space $W = W_1 \cup W_2$ of bidimension $(k_1, k_2)$ which is a cyclic bivector $\alpha = \alpha_1 \cup \alpha_2$.

By the above theorem just proved the bivectors $\alpha$, $S\alpha$, $S^2\alpha$, …, $S^{k-1}\alpha$; that is

$$\{\alpha_1, S_1\alpha_1, S_1^2\alpha_1, ..., S_1^{k_1-1}\alpha_1\}, \{\alpha_2, S_2\alpha_2, S_2^2\alpha_2, ..., S_2^{k_2-1}\alpha_2\}$$

form a bibasis for the bispace $W = W_1 \cup W_2$ and the annihilator $p_\alpha = p_{1\alpha_1} \cup p_{2\alpha_2}$ of $\alpha = \alpha_1 \cup \alpha_2$ is the biminimal neutrosophic bipolynomial for $S = S_1 \cup S_2$ (hence also the bicharacterstic neutrosophic bipolynomial for S).

If we let $\alpha^i = S^{i-1}\alpha$; that is $\alpha^i = \alpha_1^i \cup \alpha_2^i$ and $\alpha^i = S^{i-1}\alpha$ implies $\alpha_1^i = S_1^{i_1-1}\alpha_1$, $\alpha_2^i = S_2^{i_2-1}\alpha_2$; $1 \le i \le k_i - 1$ then the action of S on the bibasis $\{\alpha_1^1,...,\alpha_{k_1}^1\} \cup \{\alpha_1^2, \alpha_2^2,...,\alpha_{k_2}^2\}$ is $S\alpha^i = \alpha^{i+1}$ for i = 1, 2, …, k–1 that is $S_t\alpha_t^i = \alpha_t^{i+1}$ for i = 1, 2, …, $k_i$–1 and t = 1, 2. $S\alpha^k = -c_0\alpha^1 - ... - c_{k-1}\alpha^k$ that is $S_t\alpha_t^k = -c_0^t\alpha_t^1 - ... - c_{k_t-1}^t\alpha_t^k$ for t = 1, 2; where

$$p\alpha = \{c_0^1 + c_1^1 x + ... + c_{k_1-1}^1 x^{k_1-1} + x^{k_1}\} \cup$$
$$\{c_0^2 + c_2^2 x + ... + c_{k_2-1}^2 x^{k_2-1} + x^{k_2}\}.$$

The biexpression for $S\alpha_k$ follows from the fact $p_\alpha(S)\alpha = 0 \cup 0$; that is $S^k\alpha + c_{k-1}S^{k-1}\alpha + ... + c_1 S\alpha + c_0\alpha = 0 \cup 0$ that is

$$\{S_1^{k_1}\alpha_1 + c_{k_1-1}^1 S_1^{k_1-1}\alpha_1 + ... + c_1^1 S_1\alpha_1 + c_0^1\alpha_1\} \cup$$
$$\{S_2^{k_2}\alpha_2 + c_{k_2-1}^2 S_2^{k_2-1}\alpha_2 + ... + c_1^2 S_2\alpha_2 + c_0^2\alpha_2\} = 0 \cup 0.$$

This is given by the neutrosophic bimatrix $S = S_1 \cup S_2$ in the bibasis

$$B = B_1 \cup B_2$$
$$= \{\alpha_1^1,...,\alpha_{k_1}^1\} \cup \{\alpha_1^2,...,\alpha_{k_2}^2\}$$

$$= \begin{bmatrix} 0 & 0 & 0 & ... & 0 & -c_0^1 \\ 1 & 0 & 0 & ... & 0 & -c_1^1 \\ 0 & 1 & 0 & ... & 0 & -c_2^1 \\ \vdots & \vdots & \vdots & & \vdots & \vdots \\ 0 & 0 & 0 & ... & 1 & -c_{k_1-1}^1 \end{bmatrix} \cup \begin{bmatrix} 0 & 0 & 0 & ... & 0 & -c_0^2 \\ 1 & 0 & 0 & ... & 0 & -c_1^2 \\ 0 & 1 & 0 & ... & 0 & -c_2^2 \\ \vdots & \vdots & \vdots & & \vdots & \vdots \\ 0 & 0 & 0 & ... & 0 & -c_{k_2-1}^2 \end{bmatrix}.$$

The neutrosophic bimatrix is defined as the bicompanion bimatrix of the monic neutrosophic polynomial $p_\alpha = p_{1\alpha_1} \cup p_{2\alpha_2}$ or can also be represented with some flaw in notation as $p^1_{\alpha_1} \cup p^2_{\alpha_2}$ where $p = p^1 \cup p^2$.

Now we prove yet another interesting theorem.

**THEOREM 2.3.42:** *If $S = S_1 \cup S_2$ is a bilinear operator on a finite $(n_1, n_2)$ dimensional strong neutrosophic bivector space $W = W_1 \cup W_2$ then $S$ has a bicyclic bivector if and only if there is some bibasis for $W$ in which $S$ is represented by the bicompanion neutrosophic bimatrix of the neutrosophic biminimal polynomial for $S$.*

*Proof:* We just have noted that if $S = S_1 \cup S_2$ has a bycyclic bivector then there is such an ordered bibasis for $W = W_1 \cup W_2$. Conversely if there is some bibasis $\{\alpha^1_1,...,\alpha^1_{k_1}\} \cup \{\alpha^2_1,...,\alpha^2_{k_2}\}$ for $W$ in which $S$ is represented by the bicompanion neutrosophic biminimal polynomial, it is obvious that $\alpha^1_1 \cup \alpha^2_1$ is a bycyclic vector for $S$.

We give yet another interesting corollary.

**COROLLARY 2.3.5:** *If $A = A_1 \cup A_2$ be a bicompanion neutrosophic bimatrix of a bimonic neutrosophic bipolynomial $p = p_1 \cup p_2$ (each $p_i$ is monic) then $p$ is both the biminimal neutrosophic polynomial and the bicharacteristic neutrosophic bipolynomial of $A$.*

*Proof:* One way to see this is to let $S = S_1 \cup S_2$ a linear bioperator on $F_1^{k_1} \cup F_2^{k_2}$, which is represented by $A = A_1 \cup A_2$ in the bibasis. By applying the earlier theorem and the Cayley Hamilton theorem for bivector spaces. We give another method which is by direct calculation.

Now we proceed on to define the notion of bicyclic decomoposition or we can call it as cyclic bidecomposition and its birational form or equivalently rational biform.

Our aim is to show that any bilinear operator $T = T_1 \cup T_2$ of finite $(n_1, n_2)$ dimensional strong neutrosophic bivector space $V = V_1 \cup V_2$, there exists a biset of bivectors $\{\alpha_1^1,...,\alpha_{r_1}^1\},\{\alpha_1^2,...,\alpha_{r_2}^2\}$ in V such that

$$V = V_1 \cup V_2 =$$
$$Z(\alpha_1^1;T_1) \oplus ... \oplus Z(\alpha_{r_1}^1;T_1) \cup Z(\alpha_1^2;T_2) \oplus ... \oplus Z(\alpha_{r_2}^2;T_2).$$

In other words, we want to prove that V is a bidirect sum of bi T cyclic strong neutrosophic bivector subspaces. This will show that T is a bidirect sum of a bifinite number of bilinear operators each of which has a bicyclic bivector. The effect of this will be to reduce many problems about the general bilinear operator to similar problems about a linear bioperator which has a bicyclic bivector.

The bicyclic bidecomposition theorem is closely related to the problem in which T biinvariant bisubspaces $W = W_1 \cup W_2$ have the property that there exists a T biinvariant bisubspaces $W^1$ such that $V = W \oplus W^1$; that is

$$V = V_1 \cup V_2 = W_1 \oplus W_1^1 \cup W_2 \oplus W_2^1.$$

If $W = W_1 \cup W_2$ is any strong neutrosophic bisubspace of finite $(n_1, n_2)$ dimensional strong neutrosophic bivector space then there exists a strong neutrosophic bisubspace $W^1 = W_1^1 + W_2^1$ such that $V = W \oplus W^1$ that is,

$$V = V_1 \cup V_2 = W_1 \oplus W_1^1 \cup W_2 \oplus W_2^1$$

for each $V_i$ is a direct sum of $W_i$ and $W_i^1$, for i = 1, 2; that is $V_i = W_i \oplus W_i^1$. Usually, there are many such strong neutrosophic bivector spaces $W^1$ and each of this is called the bicomplementary to W.

We study the problem when a T biinvariant strong neutrosophic bisubspace has a complementary strong neutrosophic bisubspace which is also biinvariant under the same T.

Let us suppose that $V = W \oplus W^1$ that is $V = V_1 \cup V_2 = W_1 \oplus W_1^1 \cup W_2 \oplus W_2^1$ where both W and $W^1$ are strong neutrosophic biinvariant under T, then we study what special property is enjoyed by the strong neutrosophic bisubspace W. Each bivector $\beta = \beta_1 \cup \beta_2$ in $V = V_1 \cup V_2$ is of the form $\beta = \gamma +$ $\gamma^1$ where $\gamma$ is in W and $\gamma^1$ is in $W^1$ where $\gamma = \gamma_1 \cup \gamma_2$ and $\gamma^1 = \gamma_1^1 \cup \gamma_2^1$.

If $f = f_1 \cup f_2$ any neutrosophic bipolynomial over the scalar neutroscophic bipolynomial over the scalar neutrosophic bifield $F = F_1 \cup F_2$ then

$$f(T)\beta = f(T)\gamma + f(T)\gamma^1$$
$$= f_1(T_1)\beta_1 \cup f_2(T_2)\beta_2 = f(T)\gamma + f(T)\gamma^1$$
$$= f_1(T_1)\gamma_1 + f_1(T_1)\gamma_1^1 \cup f_2(T_2)\gamma_2 + f_2(T_2)\gamma_2^1.$$

Since W and $W^1$ are biinvariant under $T = T_1 \cup T_2$ the bivector $f(T)\gamma = f_1(T_1)\gamma_1 \cup f_2(T_2)\gamma_2$ is in $W = W_1 \cup W_2$ and $f(T)\gamma^1 = f_1(T_1)\gamma_1^1 \cup f_2(T_2)\gamma_2^1$ is in $W^1 = W_1^1 \cup W_2^1$. Therefore $f(T)\beta = f_1(T_1)\beta_1 \cup f_2(T_2)\beta_2$ is in W if and only if $f(T)\gamma^1 = 0 \cup 0$; that is $f_1(T_1)\gamma_1^1 \cup f_2(T_2)\gamma_2^1 = 0 \cup 0$. So if $f(T)\beta$ is in W then $f(T)\beta = f(T)\gamma$.

Now we define yet another new notion for bilinear operators on strong neutrosophic bivector spaces.

**DEFINITION 2.3.48:** *Let $T = T_1 \cup T_2$ be a bilinear operator on a strong neutrosophic bivector space $V = V_1 \cup V_2$ and $W = W_1 \cup W_2$ be a strong neutrosophic bivector subspace of V. We say W is bi T-admissable, if*
   *i. W is biinvariant under T,*
   *ii. $f(T)\beta$ is in W*
*for each $\beta \in V$ i.e., $f_1(T_1)\beta_1 \cup f_2(T_2)\beta_2$ is in $W = W_1 \cup W_2$ for every $\beta = \beta_1 \cup \beta_2$ in $V = V_1 \cup V_2$, there exists a bivector $\gamma = \gamma_1 \cup \gamma_2$ in $W_1 \cup W_2 = W$ such that $f(T)\beta = f(T)\gamma$ that is if W is biinvariant and has a bicomplementary biinvarient bisubspace then W is biadmissible.*

The biadmissibility characterizes those biinvariant bisubspaces which have bicomplementary biinvariant bisubspaces.

We see the biadmissibility property is involved in the bidecomposition of the bivector space $V = Z(\alpha_1;T) \oplus ... \oplus Z(\alpha_r,T) = Z(\alpha_1^1;T_1) \oplus ... \oplus Z(\alpha_{r_1}^1;T_1) \cup Z(\alpha_1^2;T_2) \oplus ... \oplus Z(\alpha_{r_2}^2;T_2)$

We arrive by some method or another we have selected bivectors $\{\alpha_1^1, \alpha_2^1, \ldots, \alpha_{r_1}^1\} \cup \{\alpha_1^2, \alpha_2^2, \ldots, \alpha_{r_2}^2\}$ and strong neutrosophic bisubspaces which is proper say

$$W_j = W_j^1 \cup W_j^2$$
$$= \{Z(\alpha_1^1; T_1) + \ldots + Z(\alpha_{j_1}^1; T_1)\} \cup \{Z(\alpha_1^2; T_2) + \ldots + Z(\alpha_{j_2}^2; T_2)\}.$$

We find the nonzero bivector $(\alpha_{j_1+1}^1 \cup \alpha_{j_2+1}^2)$ such that $W_j \cap (Z_{j+1}; T) = 0 \cup 0$ that is $W_j^1 \cap Z(\alpha_{j_1+1}^1; T_1) \cup (W_j^2) \cap Z(\alpha_{j_2+1}^2; T_2)$ $= 0 \cup 0$ because the strong neutrosophic bivector subspace $W_{j+1} = W_j \oplus Z(\alpha_{j+1}, T)$ that is
$$W_{j+1} = W_{j_1+1}^1 \cup W_{j_2+1}^2 = W_{j_1}^1 \oplus Z(\alpha_{j_1+1}; T_1) \cup W_{j_2}^2 \oplus Z(\alpha_{j_2+1}; T_2)$$
would be atleast one bidimensional nearer to exhausting V. But are we guaranteed of the existence of such $\alpha_{j+1} = \alpha_{j_1+1}^1 \cup \alpha_{j_2+1}^2$.

If $\{(\alpha_1^1, \ldots, \alpha_{j_1}^1) \cup (\alpha_2^2, \ldots, \alpha_{j_2}^2)\}$ have been choosen so that $W_j$ is T biadmissible strong neutrosophic bisubspace then it is rather easy to find a suitable $\alpha_{j_1+1}^1 \cup \alpha_{j_2+1}^2$.

Let $W = W_1 \cup W_2$ be a proper T biinvariant strong neutrosophic bisubspace. Let us find a non zero bivector $\alpha = \alpha_1 \cup \alpha_2$ such that $W \cap Z(\alpha;T) = \{0\} \cup \{0\}$; that is $W_1 \cap Z(\alpha_1;T_1) \cup W_2 \cap Z(\alpha_2;T_2) = \{0\} \cup \{0\}$. We can choose some bivector $\beta = \beta_1 \cup \beta_2$ which is not in $W = W_1 \cup W_2$; that is each $\beta_i$ is not in $W_i$, i=1, 2. Consider the T biconductor $S(\beta; W) = S(\beta_1; W_1) \cup S(\beta_2; W_2)$ which consists of all neutrosophic bipolynomials $g = g_1 \cup g_2$ such that $g(T)\beta = g_1(T_1)\beta_1 \cup g_2(T_2)\beta_2$ is in $W = W_1 \cup W_2$.

Recall that the neutrosophic bimonic polynomial $f = f_1 \cup f_2 = S(\beta; W)$; i.e., $f = f_1 \cup f_2 = S(\beta_1; W_1) \cup S(\beta_2; W_2)$ which bigenerate the neutrosophic biideals $S(\beta; W) = S(\beta_1; W_1) \cup S(\beta_2; W_2)$; that is each $f_i = S(\beta_i; W_i)$ generate the ideal $S(\beta_i; W_i)$ for $i = 1, 2$; that is $S(\beta; W)$ is also the T biconductor of $\beta$ into W. The bivector $f(T)\beta = f_1(T_1)\beta_1 \cup f_2(T_2)\beta_2$ is in $W = W_1 \cup W_2$. Now if W is T biadmissible there is a $\gamma = \gamma_1 \cup \gamma_2$ in W with $f(T)\beta = f(T)\gamma$. Let $\alpha = \beta - \gamma$ and let g be any neutrosophic bipolynomial since $\beta - \gamma$ is in W, $g(T)\beta$ will be in W if and only if $g(T)\alpha$ is in W; in other words $S(\alpha; W) = S(\beta; W)$. Thus the neutrosophic bipolynomial f is also the T biconductor of $\alpha$ into W.

But $f(T)\alpha = 0 \cup 0$. That tells us $f_1(T_1)\alpha_1 \cup f_2(T_2)\alpha_2 = 0 \cup 0$; that is $g(T)\alpha$ is in W if and only if $g(T)\alpha = g_1(T_1)\alpha_1 \cup g_2(T_2)\alpha_2 = 0 \cup 0$. The strong neutrosophic bisubspaces $Z(\alpha; T) = Z(\alpha_1; T_1) \cup Z(\alpha_2; T_2)$ and $W = W_1 \cup W_2$ are biindependent and f is the T biannihilator of $\alpha$.

Now we prove the cyclic decomposition theorem for $f_i$ linear operators on strong neutrosophic bivector spaces defined over the neutrosophic bifield $F = F_1 \cup F_2$ ($F_1$, $F_2$ are not pure neutrosophic fields) of type II.

**THEOREM 2.3.42:** *(Bicyclic decomposition theorem): Let $T = T_1 \cup T_2$ be a bilinear operator on a finite bidimensional $(n_1, n_2)$ strong neutrosophic bivector space $V = V_1 \cup V_2$ and let $W_0 = W_0^1 \cup W_0^2$ be a proper T biadmissible strong neutrosophic bivector subspace of V. There exists non zero bivectors $\{\alpha_1^1,...,\alpha_{r_1}^1\} \cup \{\alpha_1^2,...,\alpha_{r_2}^2\}$ in V with respective T biannihilators $\{p_1^1,...,p_{r_1}^1\} \cup \{p_1^2,...,p_{r_2}^2\}$ such that,*

    i.   $V = W_0 \oplus Z(\alpha_1;T) \oplus ... \oplus Z(\alpha_r;T)$
          $= W_0^1 \oplus Z(\alpha_1^1;T_1) \oplus ... \oplus Z(\alpha_{r_1}^1;T_1) \cup$
          $W_0^2 \oplus Z(\alpha_1^2;T_2) \oplus ... \oplus Z(\alpha_{r_2}^2;T_2)$

    ii.  $p_{k_r}^t$ *divides* $p_{k_r-1}^t$; *k=1, 2, ..., r and t = 1, 2.*

*Further more the integer r and the biannihilators $\{p_1^1,...,p_{r_1}^1\} \cup \{p_1^2,...,p_{r_2}^2\}$ are uniquely determined by (i) and (ii) and infact that no $\alpha_{k_r}^t$ is zero for t = 1, 2.*

*Proof:* The proof is given under four steps.

Take $W_0 = \{0\} \cup \{0\} = W_0^1 \cup W_0^2$; that is each $W_0^i = 0$ for i = 1, 2, although W does not produce any substantial simplification. Throughout the proof we shall abbreviate $f(T)\beta$ to $f\beta$ that is $f_1(T_1)\beta_1 \cup f_2(T_2)\beta_2$ to $f_1\beta_1 \cup f_2\beta_2$.

*Step 1:* There exists nonzero bivectors $\{\beta_1^1,...,\beta_{r_1}^1\} \cup \{\beta_1^2,...,\beta_{r_2}^2\}$ in the strong neutrosophic bivector space $V = V_1 \cup V_2$ such that

(i) $V = W_0 + Z(\beta_1;T) + ... + Z(\beta_r;T)$
$= W_0^1 + Z(\beta_1^1;T_1) + ... + Z(\beta_{r_1}^1;T_1) \cup$
$W_0^2 + Z(\beta_1^2;T_2) + ... + Z(\beta_{r_2}^2;T_2)$.

(ii) If $1 \leq k_i \leq r_i$; $i = 1, 2$ and $W_k = W_{k_1}^1 + W_{k_2}^2$
$= \{W_0^1 + Z(\beta_1^1;T_1) + ... + Z(\beta_{k_1}^1;T_1)\}$
$\cup \{W_0^2 + Z(\beta_1^2;T_2) + ... + Z(\beta_{k_2}^2;T_2)\}$

then the biconductor
$$p_k = p_{k_1}^1 \cup p_{k_2}^2 = S(\beta_{k_1}^1;W_{k_1-1}^1) \cup S(\beta_{k_2}^2;W_{k_2-1}^2)$$
has the maximum bidegree among all T biconductors into the strong neutrosophic bivector subspace
$$W_{k-1} = (W_{k_1-1}^1 \cup W_{k_2-1}^2)$$
that is for every $(k_1, k_2)$;
bidegree $p_k = \max_{\alpha^1 \text{ in } V_1} \deg\{S(\alpha^1;W_{k_1-1}^1)\} \cup \max_{\alpha^2 \text{ in } V_2} \deg\{S(\alpha^2;W_{k_2-1}^2)\}$.

This step depends upon only the fact that $W_0 = W_0^1 \cup W_0^2$ is a biinvariant strong neutrosophic bivector subspace. If $W = W_1 \cup W_2$ is a proper neutrosophic bi T-invariant bivector subspace
$$0 < \max_{\alpha} \text{bidegree } (S(\alpha; W)) \leq \text{bidim } V$$
that is
$$0 \cup 0 \cup < \max_{\alpha_1} \text{bidegree}(S(\alpha_1; W_1)) \cup$$
$$\max_{\alpha_2} \text{bidegree } (S(\alpha_2; W_2)) \leq (n_1, n_2)$$
and we can choose a bivector $\beta = \beta_1 \cup \beta_2$ so that bidegree $S(\beta; W) = \deg(S(\beta_1; W_1)) \cup \deg(S(\beta_2; W_2))$ attains the maximum. The strong neutrosophic bivector subspace $W + Z(\beta;T) = (W_1 + Z(\beta_1;T_1)) \cup (W_2 + Z(\beta_2;T_2))$ is then T biinvariant and has bidimension larger than bidimension W.

Apply this process to $W = W_0$ to obtain $\beta_1 = \beta_1^1 \cup \beta_2^1$. If $W_1 = W_0 + Z(\beta_1; T)$ that is
$$W_1^1 \cup W_2^1 = W_0^1 + Z(\beta_1^1; T_1) \cup W_0^2 + Z(\beta_2^1; T_2)$$
is still proper then apply the process to $W_1$ to obtain $\beta_2 = \beta_1^2 \cup \beta_2^2$.

Continue in that manner. Since bidim $W_k >$ bidim $W_k - 1$ that is
$$\mathrm{bidim}\, W_{k_1}^1 \cup \mathrm{bidim}\, W_{k_2}^2 > \mathrm{bidim}\, W_{k_1-1}^1 \cup \mathrm{bidim}\, W_{k_2-1}^2$$
we must reach $W_r = V$ that is $W_{r_1}^1 \cup W_{r_2}^2 = V_1 \cup V_2$ in not more than bidim V steps.

*Step 2:* Let $\{\beta_1^1,...,\beta_{r_1}^1\} \cup \{\beta_1^2,...,\beta_{r_2}^2\}$ be a biset of nonzero bivectors which satisfy the conditions (i) and (ii) of step 1. F$_i$x $(k_1, k_2)$; $1 \le k_i \le r_i$; $i = 1, 2$.

Let $\beta = \beta_1 \cup \beta_2$ be any bivector in the strong neutrosophic bivector space $V = V_1 \cup V_2$ and let $f = S(\beta; W_{k-1})$ that is $f_1 \cup f_2 = S(\beta_1; W_{k_1-1}^1) \cup S(\beta_2; W_{k_2-1}^2)$. If
$$f\beta = \beta_0 + \sum_{1 \le i \le k} g_i \beta_i$$
that is
$$f_1\beta_1 \cup f_2\beta_2 = (\beta_0^1 + \sum_{1 \le i_1 \le k_1} g_{i_1}^1 \beta_{i_1}^1) \cup (\beta_0^2 + \sum_{1 \le i_2 \le k_2} g_{i_2}^2 \beta_{i_2}^2)$$
$\beta_{i_t} \in W_{i_t}^t$; $t = 1, 2$, then $f = f_1 \cup f_2$ bidivides each neutrosophic bipolynomial $g_i = g_i^1 \cup g_i^2$ and $\beta_0 = f\gamma_0$ that is $\beta_0^1 \cup \beta_0^2 = f_1\gamma_0^1 \cup f_2\gamma_0^2$ where $\gamma_0 = \gamma_0^1 \cup \gamma_0^2 \in W_0 = W_0^1 \cup W_0^2$. If each $k_i = 1$ for $i = 1, 2$; this is just the statement that $W_0$ is T biadmissible. In order to prove this assertion for $(k_1, k_2) > (1, 1)$ apply the bidivision algorithms for the neutrosophic bipolynomials that is $g_i = fh_i + r_i$; $r_i = 0$ if bideg $r_i <$ bideg $f$ that is $g_{i_1}^1 \cup g_{i_2}^2 = (f_1 h_{i_1}^1 + r_{i_1}^1) \cup (f_2 h_{i_2}^2 + r_{i_2}^2)$. If bideg $r_i <$ bideg $f$. We wish to show that $r_i = 0 \cup 0$ for each $i = (i_1, i_2)$.

Let,
$$\gamma = \beta - \sum_{i=1}^{k-1} h_i \beta_i$$

that is
$$\gamma_1 \cup \gamma_2 = (\beta_1 - \sum_{i_1=1}^{k_1-1} h_{i_1}^1 \beta_{i_1}^1) \cup (\beta_2 - \sum_{i_2=1}^{k_2-1} h_{i_2}^2 \beta_{i_2}^2).$$

Since $\gamma - \beta$ is in $W_{k-1}$ that is $(\gamma_1 - \beta_1) \cup (\gamma_2 - \beta_2)$ is in $W_{k_1-1}^1 \cup W_{k_2-1}^2$. $S(\gamma_i; W_{k_i-1}^i) = S(\beta_i; W_{k_i-1}^i) = f_i$ ; $1 \leq i \leq 2$
that is
$$S(\gamma_1; W_{k_1-1}^1) \cup S(\gamma_2; W_{k_2-1}^2) = S(\beta_1; W_{k_1-1}^1) \cup S(\beta_2; W_{k_2-1}^2)$$
$$= f_1 \cup f_2.$$
$$S(\gamma; W_{k-1}) = S(\beta; W_{k-1}) = f.$$

Also
$$f\gamma = \beta_0 + \sum_{i-1}^{k-1} r_i \beta_i$$

that is
$$f_1\gamma_1 \cup f_2\gamma_2 = (\beta_0^1 + \sum_{i_1=1}^{k_1-1} r_{i_1}^1 \beta_{i_1}^1) \cup (\beta_0^2 + \sum_{i_2=1}^{k_2-1} r_{i_2}^2 \beta_{i_2}^2).$$

If $r_j = (r_{j_1}^1, r_{j_2}^2) \neq (0,0)$ we arrive at a contradiction. Let $j = (j_1, j_2)$ be the largest index $i = (i_1, i_2)$ for which $r_i = (r_i^1, r_i^2) \neq (0,0)$ then $f\gamma = \beta_0 + \sum_1^j r_i \beta_i$ ; $r_j \neq (0, 0)$ and bideg $r_j$ < bideg $f$. Let $p = S(\gamma; W_{j-1})$;
$$p_1 \cup p_2 = S(\gamma_1; W_{j_1-1}^1) \cup S(\gamma_2; W_{j_2-1}^2).$$
Since $W_{k-1} = W_{k_1-1}^1 \cup W_{k_2-1}^2$ contains $W_{j-1} = W_{j_1-1}^1 \cup W_{j_2-1}^2$ the biconductor $S(\gamma; W_{k-1})$ that is
$$f = f_1 \cup f_2 = S(\gamma_1; W_{k_1-1}^1) \cup S(\gamma_2; W_{k_2-1}^2)$$
must bidivide $p$. $p = fg$ that is $p_1 \cup p_2 = f_1g_1 \cup f_2g_2$. Apply $g(T) = g_1(T_1) \cup g_2(T_2)$ to both sides; that is
$$p\gamma = gf\gamma = gr_j\beta_j + g\beta_0 + \sum_{1 \leq i \leq j} gr_i\beta_i$$
that is,
$$p_1\gamma_1 \cup p_2\gamma_2 = (g_1f_1\gamma_1 \cup g_2f_2\gamma_2)$$

$$= g_1 r_{j_1}^1 \beta_{j_1}^1 + g_1 \beta_0^1 + \sum_{1 \leq i_1 \leq j_1} g_1 r_{i_1}^1 \beta_{i_1}^1 \cup g_2 r_{j_2}^2 \beta_{j_2}^2 + g_2 \beta_0^2 + \sum_{1 \leq i_2 \leq j_2} g_2 r_{i_2}^2 \beta_{i_2}^2 .$$

By definition, $p\gamma$ is in $W_{j-1}$ and the last two terms on the right side of the above equation are in $W_{j-1} = W_{j_1-1}^1 \cup W_{j_2-1}^2$.

Therefore
$$gr_j \beta_j = g_1 r_{j_1}^1 \beta_{j_1}^1 \cup g_2 r_{j_2}^2 \beta_{j_2}^2$$
is in $W_{j-1} = W_{j_1-1}^1 \cup W_{j_2-1}^2$. Now using condition (ii) of step 1 bideg $(gr_j) \geq \text{bideg}(S(\beta_j; W_{j-1}))$; that is

$\deg(g_1, r_{j_1}) \cup \deg(g_2, r_{j_2})$

$\quad \geq \quad \deg S(\beta_{j_1}; W_{j_1-1}^1) \cup \deg S(\beta_{j_2}; W_{j_2-1}^2)$

$\quad = \quad \text{bidegp}_j$

$\quad = \quad \deg p_{j_1}^1 \cup \deg p p_{j_2}^2 \geq \text{bidegree}(S(\gamma; W_{j-1}))$

$\quad = \quad \deg S(\gamma_1; W_{j_1-1}) \cup \deg S(\gamma_2; W_{j_2-1})$

$\quad = \quad \text{bidegree } p$

$\quad = \quad \deg p_1 \cup \deg p_2$

$\quad = \quad \text{bideg } fg$

$\quad = \quad \deg f_1 g_1 \cup \deg f_2 g_2.$

Thus bideg $r_j >$ bideg $f$; i.e., $\deg r_{j_1} \cup \deg r_{j_2} > \deg f_1 \cup \deg f_2$ and that contradicts the choice of $j = (j_1, j_2)$. We now know that $f = f_1 \cup f_2$ bidivides each $g_i = g_{i_1}^1 \cup g_{i_2}^2$ that is $f_{i_t}$ divides $g_{i_t}^t$; t =1, 2 and hence $\beta_0 = f\gamma$ that is $\beta_0^1 \cup \beta_0^2 = f_1 \gamma_1 \cup f_2 \gamma_2$. Since $W_0 = W_0^1 \cup W_0^2$ is T biadmissible (i.e., each $W_0^k$ is $T_k$ admissible for k = 1, 2); we have $\beta_0 = f\gamma_0$ where
$$\gamma_0 = \gamma_0^1 \cup \gamma_0^2 \in W_0 = W_0^1 \cup W_0^2$$
that is $\beta_0^1 \cup \beta_0^2 = f_1 \gamma_0^1 \cup f_2 \gamma_0^2$ where $\gamma_0 = W_0$. We make a mention that step 2 is a stronger form of the assertion that each of the strong neutrosophic vector bisubspaces $W_1 = W_1^1 \cup W_1^2$, $W_2 = W_2^1 \cup W_2^2$, ..., $W_r = W_r^1 \cup W_r^2$ is T biadmissible.

*Step 3:* There exists non zero bivectors
$$(\alpha_1^1,...,\alpha_{r_1}^1) \cup (\alpha_1^2,...,\alpha_{r_2}^2)$$
in $V = V_1 \cup V_2$ which satisfy condition (i) and (ii) of the theorem. Start with bivectors
$$\{\beta_1^1,\beta_2^1,...,\beta_{r_1}^1\} \cup \{\beta_1^2,\beta_2^2,...,\beta_{r_2}^2\}$$
as in step 1. Fix $k = (k_1, k_2)$ as $1 \le k_i \le r_i$; $i = 1, 2$. We apply step 2 to the bivector $\beta = \beta_1 \cup \beta_2 = \beta_{k_1}^1 \cup \beta_{k_2}^2 = \beta_k$ and T biconductor $f = f_1 \cup f_2 = p_{k_1}^1 \cup p_{k_2}^2 = p_k$. We obtain

$$p_k \beta_k = p_k \gamma_0 + \sum_{1 \le i \le k} p_k h_i \beta_i \ ;$$

that is,

$$p_{k_1}^1 \beta_{k_1}^1 \cup p_{k_2}^2 \beta_{k_2}^2$$
$$= (p_{k_1}^1 \gamma_0^1 + \sum_{1 \le i_1 \le k_1} p_{k_1}^1 h_{i_1}^1 \beta_{i_1}^1) \cup (p_{k_2}^2 \gamma_0^2 + \sum_{1 \le i_2 \le k_2} p_{k_2}^2 h_{i_2}^2 \beta_{i_2}^2).$$

where $\gamma_0 = \gamma_0^1 \cup \gamma_0^2$ is in $W_0 = W_0^1 \cup W_0^2$ and $\{h_1^1,...,h_{k_1-1}^1\} \cup \{h_1^2,...,h_{k_2-1}^2\}$ are neutrosophic bipolynomials. Let

$$\alpha_k = \beta_k - \gamma_0 - \sum_{1 \le i \le k} h_i \beta_i \ ;$$

i.e.,

$$\{\alpha_{k_1}^1 \cup \alpha_{k_2}^2\} = \left(\beta_{k_1}^1 - \gamma_0^1 - \sum_{1 \le i_1 < k_1} h_{i_1}^1 \beta_{i_1}^1\right) \cup \left(\beta_{k_2}^2 - \gamma_0^2 - \sum_{1 \le i_2 < k_2} h_{i_2}^2 \beta_{i_2}^2\right).$$

Since
$$\beta_k - \alpha_k = (\beta_{k_1}^1 - \alpha_{k_1}^1) \cup (\beta_{k_2}^2 - \alpha_{k_2}^2)$$
is in
$$W_{k-1} = W_{k_1-1}^1 \cup W_{k_2-1}^2 \ ;$$
is in
$$S(\alpha_k; W_{k-1}) = S(\beta_k; W_{k-1}) = p_k$$
$$= S(\alpha_{k_1}^1; W_{k_1-1}^1) \cup (\alpha_{k_2}^2; W_{k_2-1}^2)$$
$$= p_{k_1}^1 \cup p_{k_2}^2$$

and since $p_k\alpha_k = 0 \cup 0$ that is $p^1_{k_1}\alpha^1_{k_1} \cup p^2_{k_2}\alpha^2_{k_2} = 0 \cup 0$ we have,

$$W^1_{k_1-1} \cap Z(\alpha^1_{k_1};T_1) \cup W^2_{k_2-1} \cap Z(\alpha^2_{k_2};T_2) = \{0\} \cup \{0\}.$$

Because each $\alpha_k = \alpha^1_{k_1} \cup \alpha^2_{k_2}$ satisfies the above two equations just mentioned, it follows that,

$$W_k = W_0 \oplus Z(\alpha_1;T) \oplus ... \oplus Z(\alpha_k;T)$$

that is

$$W^1_{k_1} \cup W^2_{k_2} = W^1_0 \oplus Z(\alpha^1_1;T_1) \oplus ... \oplus Z(\alpha^1_{k_1};T_1)\} \cup$$
$$\{W^2_0 \oplus Z(\alpha^2_1;T_2) \oplus ... \oplus Z(\alpha^2_{k_2};T_2)\}$$

and that $p_k = p^1_{k_1} \cup p^2_{k_2}$ is the T biannihilator of $\alpha_k = \alpha^1_{k_1} \cup \alpha^2_{k_2}$ other words, bivectors $\{\alpha^1_1,...,\alpha^1_{r_1}\} \cup \{\alpha^2_1,...,\alpha^2_{r_2}\}$ define the same bisequence of strong neutrosophic bisubspaces $W_1 = W^1_1 \cup W^2_1$, $W_2 = W^1_2 \cup W^2_2$, …, as do the bivector $\{\beta^1_1 \cup ... \cup \beta^1_{r_1}\}$, $\{\beta^2_1 \cup ... \cup \beta^2_{r_2}\}$ and the T biconductors $p_k = S(\alpha_k;W_{k-1})$ that is $(p^1_{k_1} \cup p^2_{k_2}) = S(\alpha^1_{k_1};W^1_{k_1-1}) \cup S(\alpha^2_{k_2};W^2_{k_2-1})$ have the same maximality properties. The bivectors $\{\alpha^1_1,...,\alpha^1_{r_1}\} \cup \{\alpha^2_1,...,\alpha^2_{r_2}\}$ have the additional property that the strong neutrosophic bivector spaces $W_0 = \{W^1_0 \cup W^2_0\}$,

$$Z(\alpha_1;T) = Z(\alpha^1_1;T_1) \cup Z(\alpha^2_1;T_2)$$
$$Z(\alpha_2;T) = Z(\alpha^1_2;T_1) \cup Z(\alpha^2_2;T_2)$$

are biindependent. It is therefore easy to verify condition (ii) of the theorem. Since $(p^1_{i_1}\alpha^1_{i_1}) \cup (p^2_{i_2}\alpha^2_{i_2}) = 0 \cup 0$, we have the trivial relation

$$p_k\alpha_k = p^1_{k_1}\alpha^1_{k_1} \cup p^2_{k_2}\alpha^2_{k_2}$$
$$= (0 + p^1_1\alpha^1_1 + ... + p^1_{k_1-1}\alpha^1_{k_1-1}) \cup (0 + p^2_1\alpha^2_1 + ... + p^2_{k_2-1}\alpha^2_{k_2-1}).$$

Apply step 2 with $\{\beta^1_1,...,\beta^1_{k_1}\} \cup \{\beta^2_1,...,\beta^2_{k_2}\}$ replaced by $\{\alpha^1_1,...,\alpha^1_{k_1}\} \cup \{\alpha^2_1,...,\alpha^2_{k_2}\}$ and with $\beta = \beta_1 \cup \beta_2 = \alpha^1_{k_1} \cup \alpha^1_{k_2}$, $p_k$ bidivides each $p_i$; $i < k$ that is $(i_1, i_2) < (k_1, k_2)$; i.e., $p^1_{k_1} \cup p^2_{k_2}$ bidivides each $p^1_{i_1} \cup p^2_{i_2}$, i.e., each $p^t_{k_t}$ divides $p^t_{i_t}$ for $t = 1, 2$.

*Step 4:* The number $r = (r_1, r_2)$ and the neutrosophic bipolynomials $(p_1, ..., p_{r_1})$, $(p_2, ..., p_{r_2})$ are uniquely determined by the condition of the theorem. Suppose that in addition to the bivectors $\{\alpha_1^1, ..., \alpha_{r_1}^1\} \cup \{\alpha_1^2, ..., \alpha_{r_2}^2\}$ we have non zero bivectors

$$\{\gamma_1^1, ..., \gamma_{s_1}^1\} \cup \{\gamma_1^2, ..., \gamma_{s_2}^2\}$$

with respective T biannihilators

$$\{g_1^1, ..., g_{s_1}^1\} \cup \{g_1^2, ..., g_{s_2}^2\}$$

such that

$$V = W_0 \oplus Z(\gamma_1; T) \oplus ... \oplus Z(\gamma_s; T)$$

that is

$$V = V_1 \cup V_2 = W_0^1 \oplus Z(\gamma_1^1; T_1) \oplus ... \oplus Z(\gamma_{s_1}^1; T_1)\} \cup$$
$$\{W_0^2 \oplus Z(\gamma_1^2; T_2) \oplus ... \oplus Z(\gamma_{s_2}^2; T_2)\}$$

$g_{k_t}^t$ divides $g_{k_t - 1}^t$ for $t = 1, 2$ and $k_t = 1, 2, ..., s_t$. We shall show that $r = s$ that is $(r_1, r_2) = (s_1, s_2)$ and $p_i^t = g_i^t$; $1 \leq t \leq 2$; that is $p_i^1 \cup p_i^2 = g_i^1 \cup g_i^2$ for each i. We see that $p_1 = g_1 = p_1^1 \cup p_1^2 = g_1^1 \cup g_1^2$. The neutrosophic bipolynomial $g_1 = g_1^1 \cup g_1^2$ is determined by the above equation as the T biconductor of V into $W_0$; that is $V = V_1 \cup V_2$ into $W_0^1 \cup W_0^2$. Let $S(V; W_0) = S(V_1; W_0^1) \cup S(V_2; W_0^2)$ be the collection of all neutrosophic bipolynomials $f = f_1 \cup f_2$ such that $f\beta = f_1\beta_1 \cup f_2\beta_2$ is in $W_0$ for every $\beta = \beta_1 \cup \beta_2$ in V that is neutrosophic polynomials f such that the birange of $f(T) = $ range of $f_1(T_1) \cup$ range of $f_2(T_2)$ is contained in $W_1 = W_0^1 \cup W_0^2$; i.e., range of $f_i(T_i)$ is in $W_0^i$ for i=1, 2. Thus $S(V_i, W_i)$ is a non zero neutrosophic zero ideal in the neutrosophic polynomial algebra so that we see $S(V; W_0) = S(V_1; W_0^1) \cup S(V_2; W_0^2)$ is a non zero neutrosophic biideal in the neutrosophic bipolynomial algebra.

The neutrosophic polynomial $g_1^t$ is the monic generator of that neutrosophic ideal i.e., the bimonic neutrosophic polynomial $g_1 = g_1^1 \cup g_1^2$ is the neutrosophic monic bigenerator of that biideal. Each $\beta = \beta_1 \cup \beta_2$ in $V = V_1 \cup V_2$ has the form

$$\beta = (\beta_0^1 + f_1^1\gamma_1^1 + \ldots + f_{s_1}^1\gamma_{s_1}^1) \cup (\beta_0^2 + f_1^2\gamma_1^2 + \ldots + f_{s_2}^2\gamma_{s_2}^2)$$

and so

$$g_1\beta = g_1\beta_0 + \sum_1^s g_1 f_i \gamma_i$$

that is

$$g_1^1\beta_1 \cup g_1^2\beta_2 = \left[g_1^1\beta_0^1 + \sum_1^s g_1^1 f_{i_1}^1 \gamma_{i_1}^1\right] \cup \left[g_1^2\beta_0^2 + \sum_1^s g_1^2 f_{i_2}^2 \gamma_{i_2}^2\right].$$

Since each $g_i^t$ divides $g_1^t$ for t = 1, 2 we have $g_1\gamma_i = 0 \cup 0$ that is $g_1^1\gamma_{i_1}^1 \cup g_1^2\gamma_{i_2}^2 = 0 \cup 0$ for all i = ($i_1$, $i_2$) and $g_0\beta = g_1\beta_0$ is in $W_0 = W_0^1 \cup W_0^2$. Thus $g_i^t$ is in $S(V_t; W_0^t)$ for t = 1, 2; so $g_1 = g_1^1 \cup g_1^2$ is in $S(V; W_0) = S(V_1; W_0^1) \cup S(V_2; W_0^2)$.

Since each $g_i^t$ is monic, $g_1$ is a monic neutrosophic bipolynomial of least bidegree which sends $\gamma_1^t$ into $W_0^t$ so that $\gamma_1 = \gamma_1^1 \cup \gamma_1^2$ into $W_0 = W_0^1 \cup W_0^2$; we see that $g_1 = g_1^1 \cup g_1^2$ is the neutrosophic monic bipolynomial of least bidegree in the neutrosophic biideal $S(V; W_0)$. By the same argument $p_i$ is the bigenerator of the neutrosophic ideal so $p_1 = g_1$; that is $p_1^1 \cup p_1^2 = g_1^1 \cup g_1^2$.

If $f = f_1 \cup f_2$ is a neutrosophic bipolynomial and $W = W_1 \cup W_2$ is a strong neutrosophic bisubspace of $V = V_1 \cup V_2$ we shall employ the short hand fW for the set of all bivectors $f\alpha = f_1\alpha_1 \cup f_2\alpha_2$ with $\alpha = \alpha_1 \cup \alpha_2$ in $W = W_1 \cup W_2$.

The three facts can be proved by the reader.
(1) $fZ(\alpha;T) = Z(f_\alpha;T)$ that is
    $f_1(Z(\alpha_1;T_1)) \cup f_2(Z(\alpha_2;T_2)) = Z(f_1\alpha_1; T_1) \cup Z(f_2\alpha_2; T_2)$.
(2) $V = V_1 \oplus \ldots \oplus V_k = V_1^1 \oplus \ldots \oplus V_{k_1}^1 \cup V_1^2 \oplus \ldots \oplus V_{k_2}^2$ where each $V_t$ is a biinvariant under $T_i$; $1 \leq i \leq t$; t = 1, 2, then $fV = fV_1 \oplus fV_2$ that is
    $f_1V_1 \cup f_2V_2 = f_1V_1^1 \oplus \ldots \oplus f_1V_{k_1}^1 \cup f_2V_1^2 \oplus \ldots \oplus f_2V_{k_2}^2$.
(3) If $\alpha = \alpha_1 \cup \alpha_2$ and $\gamma = \gamma_1 \cup \gamma_2$ have the same T biannihilator then f$\alpha$ and f$\gamma$ have the same T biannihilator and hence bidim $Z(f\alpha;T)$ = bidim $Z(f\gamma;T)$ that is f$\alpha$ = $f_1\alpha_1 \cup f_2\alpha_2$ and $f\gamma = f_1\gamma_1 \cup f_2\gamma_2$ with dim $Z(f_1\alpha_1; T_1) \cup$ dim $Z(f_2\alpha_2; T_2)$ = dim $Z(f_1\gamma_1; T_1) \cup$ dim $Z(f_2\gamma_2; T_2)$.

Since we know $p_1 = g_1$ we know that $Z(\alpha_1\ T)$ and $Z(\gamma_1; T)$ have the same bidimension. Therefore bidim $W_0$ + bidim $Z(Y_1, T) <$ bidim V as before

$$\dim W_0^1 + \dim Z(\gamma_1^1; T_1) \cup \dim W_0^2 + \dim Z(\gamma_1^2; T_2)$$
$$\leq \dim V_1 \cup \dim V_2.$$

Now to check whether or not $p_2 = g_2$; $p_2^1 \cup p_2^2 = g_2^1 \cup g_2^2$. From the decomposition of $V = V_1 \cup V_2$ we obtain the two decomposition of the strong neutrosophic bivector subspace $p_2 V = p_2^1 V_1 \cup p_2^2 V_2$.

$$p_2 V = p_2^1 W_0 \oplus Z(p_2\alpha_1; T)$$

that is

$$p_2^1 V_1 \cup p_2^2 V_2 = p_2^1 W_0^1 \oplus Z(p_2^1\alpha_1^1; T_1) \cup p_2^2 W_0^2 \oplus Z(p_2^2\alpha_1^2; T_2).$$
$$p_2 V = p_2 W_0 \oplus Z(p_2\gamma_1; T) \oplus ... \oplus Z(p_2\gamma_s; T)$$

that is

$$p_2^1 V_1 \cup p_2^2 V_2 = p_2^1 W_0^1 \oplus Z(p_2^1\gamma_1^1; T_1) \oplus ... \oplus Z(p_2^1\gamma_{s_1}^1; T_1) \cup$$
$$p_2^2 W_0^2 \oplus Z(p_2^2\gamma_1^2; T_2) \oplus ... \oplus Z(p_2^2\gamma_{s_2}^2; T_2).$$

The proof follows from the fact if $r = (r_1, r_2) \geq (2, 2)$ then $p_2 = p_2^1 \cup p_2^2 = g_2 = g_2^1 \cup g_2^2$. We have made use of the facts (1) and (2) above and we have used the fact $p_2\alpha_i = p_2^1\alpha_{i_1}^1 \cup p_2^2\alpha_{i_2}^2 = 0 \cup 0$; $i = (i_1, i_2) > (2, 2)$. Since we know $p_1 = g_1$ fact (3) above tell us that, $Z(p_2\alpha_i; T) = Z(p_2^1\alpha_1^1; T_1) \cup Z(p_2^2\alpha_1^2; T_2)$
and
$$Z(p_2\gamma_i; T) = Z(p_2^1\gamma_1^1; T_1) \cup Z(p_2^2\gamma_1^2; T_2)$$

have the same bidimension. Hence it is apparent from above equalities that

$$\text{bidim } Z(p_2\gamma_i; T) = 0 \cup 0.$$

$\dim Z(p_2^1\gamma_{i_1}^1; T_1) \cup \dim Z(p_2^2\gamma_{i_2}^2; T_2) = 0 \cup 0$; $i = (i_1, i_2) \geq (2, 2)$.

We conclude $p_2\gamma_2 = (p_2^1\gamma_2^1) \cup (p_2^2\gamma_2^2) = 0 \cup 0$ and $g_2$ bidivides $p_2$; that is $g_2^t$ divides $p_2^t$ for t = 1, 2. The argument can be reserved to show that $p_2$ bidivides $g_2$; i.e., $p_2^t$ divides $g_2^t$ for each t; t = 1, 2. Hence $p_2 = g_2$.

We leave the following corollaries for the reader to prove.

**COROLLARY 2.3.6:** *If $T = T_1 \cup T_2$ is a bilinear operator on a finite $(n_1, n_2)$ bidimensional strong neutrosophic bivector space $V = V_1 \cup V_2$ then T-biadmissible strong neutrosophic vector bisubspace has a complementary strong neutrosophic bivector subspace which is also invariant under T.*

**COROLLARY 2.3.7:** *Let $T = T_1 \cup T_2$ be a bilinear operator on a finite $(n_1, n_2)$ strong neutrosophic bivector space $V = V_1 \cup V_2$.*
  i. *There exists bivectors $\alpha = \alpha_1 \cup \alpha_2$ in $V = V_1 \cup V_2$ such that the T biannihilator of $\alpha$ is the neutrosophic biminimal polynomial for T.*
  ii. *T has a bicyclic bivector if and only if the bicharacterstic and neutrosophic biminimal polynomials for T are identical.*

Now we proceed on to prove the generalized Cayley Hamilton theorem.

**THEOREM 2.3.44:** *(Generalized Cayley Hamilton theorem). Let $T = T_1 \cup T_2$ be a bilinear operator on a finite $(n_1, n_2)$ finite bidimensional strong neutrosophic bivector space $V = V_1 \cup V_2$ over a neutrosophic bifield $F = F_1 \cup F_2$ of type II (Both $F_1$ and $F_2$ are not pure neutrosophic). Let p and f be the biminimal bicharacterstic neutrosophic bipolynomials for T, respectively.*
  i. *p bidivides f i.e., $p = p_1 \cup p_2$ and $f = f^1 \cup f^2$ then $p_i$ divides $f^i$; i = 1, 2.*
  ii. *p and f have the same prime factors except for multiplicities.*
  iii. *If $p = f_1^{r_1}...f_k^{r_k}$ is the prime factorization of p then $f = f_1^{d_1}...f_k^{d_k}$ where $d_i$ is the bimultiplicity of $f_i(T)^{r_i}$ bidivided by the bidegree $f_i$. That is if*
  $$p = p_1 \cup p_2 = (f_1^1)^{r_1^1}...(f_{k_1}^1)^{r_{k_1}^1} \cup (f_1^2)^{r_1^2}...(f_{k_2}^2)^{r_{k_2}^2}$$

*then*

$$f = (f_1^1)^{d_1^1} \ldots (f_{k_1}^1)^{d_{k_1}^1} \cup (f_1^2)^{d_1^2} \ldots (f_{k_2}^2)^{d_{k_2}^2}; \; d_i^t = (d_1^t, \ldots, d_{k_t}^t)$$

*is the nullity of $f_i^t(T_t)^{r_i^t}$ which is bidivided by the bidegree $f_i^t$; i.e., ; $1 \leq i \leq k_t$; this is true for each t, t = 1, 2.*

*Proof:* The trivial case $V = \{0\} \cup \{0\}$ is obvious. To prove (i) and (ii) consider a bicyclic decomposition

$$V = Z(\alpha_1; T) \oplus \ldots \oplus Z(\alpha_r; T)$$
$$= Z(\alpha_1^1; T_1) \oplus \ldots \oplus Z(\alpha_{r_1}^1; T_1) \cup Z(\alpha_1^2; T_2) \oplus \ldots \oplus Z(\alpha_{r_2}^2; T_2).$$

By the second corollary $p_1 = p$. Let $S_i = S_i^1 \cup S_i^2$ be the birestriction of $T = T_1 \cup T_2$; i.e., each $S_i^s$ is the restriction of $T_s$ (for s = 1, 2, ..., $r_s$) to $Z(\alpha_i^s; T_s)$. Then $S_i$ has a bicyclic bivector so that $p_i = p_{i_1}^1 \cup p_{i_2}^1$ is both biminimal neutrosophic polynomial and the bicharacteristic neutrosophic polynomial for $S_i$. Therefore the neutrosophic bicharacteristic polynomial $f = f^1 \cup f^2$ is the byproduct $f = p_1^1 \ldots p_{r_1}^1 \cup p_1^2 \ldots p_{r_2}^2$. That is evident from earlier results that the neutrosophic bimatrix of T assumes a suitable bibasis.

Clearly $p_1 = p$ bidivides f; hence the claim (1). Obviously any prime bidivisor of p is a prime bidivisor of f. Conversely a prime bidivisor of $f = p_1^1 \ldots p_{r_1}^1 \cup p_1^2 \ldots p_{r_2}^2$ must bidivide one of the factor $p_i^t$ which in turn bidivides $p_1$.

Let $p = (f_1^1)^{r_1^1} \ldots (f_{k_1}^1)^{r_{k_1}^1} \cup (f_1^2)^{r_1^2} \ldots (f_{k_2}^2)^{r_{k_2}^2}$ be the prime bifactorization of p. We employ the biprimary decomposition theorem which tells $V_t^i = V_1^i \cup V_1^i$ is the binull space for $f_i^t(T_t)^{r_i^t}$ then

$$V = V_1 \oplus \ldots \oplus V_k = (V_1^1 \oplus \ldots \oplus V_{k_1}^1) \cup (V_1^2 \oplus \ldots \oplus V_{k_2}^2)$$

and $(f_i^t)^{r_i^t}$ is the neutrosophic minimal polynomial of the operator $T_t^i$ restricting $T_t$ to the strong neutrosophic invariant subspace $V_t^i$. This is true for each t; t = 1, 2. Apply part (ii) of the present theorem to the bioperator $T_t^i$. Since its neutrosophic minimal polynomial is a power of the prime $f_i^t$ the neutrosophic characteristic polynomial for $T_t^i$ has the form $(f_i^t)^{r_i^t}$ where $d_i^t > r_i^t$; $t = 1, 2$.

We have

$$d_i^t = \frac{\dim V_t^i}{\deg f_i^t}$$

for every $t = 1, 2$ and $\dim V_t^i$ = nullity $f_i^t(T_t)^{r_i^t}$ for every t; $t = 1, 2$. Since $T_t$ is the direct sum of operator $T_t^1,...,T_t^{k_1}$ the neutrosophic characteristic polynomial $f^t$ is the product $f^t = (f_1^t)^{d_1^t}...(f_{k_t}^t)^{d_{k_t}^t}$. Hence the claim.

The following corollary is left as an exercise for the reader.

**COROLLARY 2.3.8:** *Let $T = T_1 \cup T_2$ be a binilpotent operator of the strong neutrosophic bivector space of $(n_1, n_2)$ bidimension over the neutrosophic bifield $F = F_1 \cup F_2$ (both $F_1$ and $F_2$ are not pure neutrosophic fields) of type II then the bicharacteristic bipolynomial for T is $x^{n_1} \cup x^{n_2}$.*

Let us observe that the neutrosophic bimatrix analogue of the bicyclic decomposition theorem. If we have the bioperator $T = T_1 \cup T_2$ and the bidirect sum decomposition, let $B^i$ be the bicyclic ordered bibasis

$$\{\alpha_{i_1}^1, T_1\alpha_{i_1}^1,...,T_1^{k_{r_1}^1-1}\alpha_{i_1}^1\} \cup \{\alpha_{i_2}^2, T_2\alpha_{i_2}^2,...,T_2^{k_{r_2}^2-1}\alpha_{i_2}^2\}$$

for $Z(\alpha_i;T) = Z(\alpha_{i_1}^1;T_1) \cup Z(\alpha_{i_2}^2;T_2)$. Here $(k_{i_1}^1, k_{i_2}^2)$ denotes the bidimension of $Z(\alpha_i;T)$ that is the bidegree of the biannihilator $p_i = p_{i_1}^1 \cup p_{i_2}^2$.

The neutrosophic bimatrix of the induced operator $T_i$ in the bibasis $B_i$ is the bicompanion neutrosophic bimatrix of the neutrosophic bipolynomial $p_i$. Thus if we let B to be the bibasis for V.

which is the biunion of $B^i$ arranged in order $\{B_1^1 \ldots B_{r_1}^1\} \cup \{B_1^2 \ldots B_{r_2}^2\}$; then the neutrosophic bimatrix of T in the bibasis B will be $A = A_1 \cup A_2$.

$$= \begin{bmatrix} A_1^1 & 0 & \cdots & 0 \\ 0 & A_2^1 & \cdots & 0 \\ \vdots & & & \vdots \\ 0 & 0 & \cdots & A_{r_1}^1 \end{bmatrix} \cup \begin{bmatrix} A_1^2 & 0 & \cdots & 0 \\ 0 & A_2^2 & \cdots & 0 \\ \vdots & \vdots & & \vdots \\ 0 & 0 & \cdots & A_{r_2}^2 \end{bmatrix}$$

where $A_i^t$ is the $k_i^t \times k_i^t$ companion neutrosophic matrix of $p_i^t$ for $t = 1, 2$. A $(n_1 \times n_1, n_2 \times n_2)$ neutrosophic bimatrix A which is the bidirect sum of the neutrosophic bicompanion matrices of the non scalar monic neutrosophic bipolynomial $\{p_1^1 \ldots p_{r_1}^1\} \cup \{p_1^2 \ldots p_{r_2}^2\}$ such that $p_{i_t+1}^t$ divides $p_{i_t}^t$ for $i_t = 1, 2, \ldots, r_t - 1$ and $t = 1, 2$ will be defined as the rational biform or equivalently birational form.

**THEOREM 2.3.45:** *Let $F = F_1 \cup F_2$ be a neutrosophic bifield (Both $F_1$ and $F_2$ are not pure neutrosophic). Let $B = B_1 \cup B_2$ be a $(n_1 \times n_1, n_2 \times n_2)$ neutrosophic bimatrix over F. Then B is bisimilar over the bifield F to one and only one neutrosophic matrix in the rational form.*

*Proof:* We know from the usual neutrosophic square matrix every square matrix over a fixed neutrosophic field is similar to one and only one neutrosophic matrix which is in the rational form.

So the neutrosophic bimatrix $B = B_1 \cup B_2$ over the neutrosophic bifield $F = F_1 \cup F_2$ is such that each $B_i$ is a $n_i \times n_i$ neutrosophic square matrix over $F_i$; is similar to one and only one neutrosophic matrix which is in the rational form say $C_i$.

This is true for every i; i = 1, 2 so $B = B_1 \cup B_2$ is bisimilar over the field to one and only one bimatrix C which is in the rational biform.

The neutrosophic bipolynomials $\{p_1^1 \ldots p_{r_1}^1\} \cup \{p_1^2 \ldots p_{r_2}^2\}$ are called invariant bifactors or biinvariant factors for the neutrosophic bimatrix $B = B_1 \cup B_2$.

We shall just introduce the notion of biJordan form or Jordan biform for a strong neutrosophic bivector space of type II.

Suppose that $N = N_1 \cup N_2$ be a nilpotent bilinear operator on finite ($n_1$ $n_2$) bidimension strong neutrosophic bivector space $V = V_1 \cup V_2$ over a neutrosophic bifield $F = F_1 \cup F_2$ ($F_1$ and $F_2$ are not pure neutrosophic) of type II. Consider the bicyclic decomposition for N which we have described in the theorem. We have a pair of positive integers ($r_1$, $r_2$) and non zero bivector $\{\alpha_1^i, \alpha_2^i\}$ in V with biannihilators $\{p_1^1 \ldots p_{r_1}^1\} \cup \{p_1^2 \ldots p_{r_2}^2\}$ such that

$$V = Z(\alpha_1; N) \oplus \ldots \oplus Z(\alpha_r; N)$$
$$= Z(\alpha_1^1, N_1) \oplus \ldots \oplus (\alpha_{r_1}^1, N_1) \cup Z(\alpha_1^2, N_2) \oplus \ldots \oplus (\alpha_{r_2}^2, N_2)$$

and $p_{i_t+1}^t$ divides $p_{i_t}^t$ for $i_t = 1, 2, \ldots, r_t - 1$ and $t = 1, 2$. Since N is binilpotent and the biminimal neutrosophic polynomial is $x^{k_1} \cup x^{k_2}$ with $k_t \leq n_t$; $t = 1, 2$. Thus each $p_{i_t}^t$ is of the form $p_{i_t}^t = x^{k_i^t}$ and the bidivisibility condition says $k_1^t \geq k_2^t \geq \ldots \geq k_{r_t}^t$; $t = 1, 2$. Of course $k_1^t = k^t$ and $k_r^t \geq 1$.

The bicompanion neutrosophic bimatrix of $x^{k_{i_1}^1} \cup x^{k_{i_2}^2}$ is the $k_{i_r} \times k_{i_r}$ neutrosophic bimatrix. $A = A_{i_1}^1 \cup A_{i_2}^2$ with

$$A_{i_t}^t = \begin{bmatrix} 0 & 0 & \cdots & 0 & 0 \\ 1 & 0 & \cdots & 0 & 0 \\ 0 & 1 & \cdots & 0 & 0 \\ \vdots & & & & \vdots \\ 0 & 0 & \cdots & 1 & 0 \end{bmatrix};$$

$t = 1, 2$. Thus we from earlier results have a bibasis for $V = V_1 \cup V_2$ in which the neutrosophic bimatrix of N is the bidirect

sum of the elementary nilpotent neutrosophic bimatrices of sizes of $i_t$ which decreases as $i_t$ increases. One sees from this that associated with a binilpotent $(n_1 \times n_1, n_2 \times n_2)$ neutrosophic bimatrix is a positive pair of integers $(r_1, r_2)$ that is $\{k_1^1, ..., k_{r_1}^1\} \cup \{k_1^2, ..., k_{r_2}^2\}$ such that

$$\{k_1^1 + ... + k_{r_1}^1\} = n_1 \text{ and } \{k_1^2 + ... + k_{r_2}^2\} = n_2$$

and $k_{i_t}^t \geq k_{i_{t+1}}^t$ ; $t = 1, 2$ and $1 \leq i, i + 1 \leq r_t$ and these bisets of positive integers determine the birational form of the neutrosophic bimatrix that is they determine the neutrosophic bimatrix up to similarity.

Here is one thing, we like to mention about the binilpotent bioperator N.

The positive biinteger $(r_1, r_2)$ is precisely the binullity of N infact the strong neutrosophic binull space has a bibasis with $(r_1, r_2)$ bivectors $N_1^{k_{i_1}-1}\alpha_{i_1}^1 \cup N_2^{k_{i_2}-1}\alpha_{i_2}^2$. For let $\alpha = \alpha_1 \cup \alpha_2$ be in the strong neutrosophic binull space of N we write

$$\alpha = \left(f_1^1\alpha_1^1 + ... + f_{r_1}^1\alpha_{r_1}^1\right) \cup \left(f_1^2\alpha_1^2 + ... + f_{r_2}^2\alpha_{r_2}^2\right)$$

where $\left(f_{i_1}^1, f_{i_2}^2\right)$ is a neutrosophic bipolynomial the bidegree of which we may assume is less than $k_{i_1}, k_{i_2}$. Since $N\alpha = 0 \cup 0$ for each $i_r$ we have

$$\begin{aligned}
0 \cup 0 &= N(f_i, \alpha_i) \\
&= N_1\left(f_{i_t}, \alpha_{i_t}\right) \cup N_2\left(f_{i_t}, \alpha_{i_t}\right) \\
&= N_1 f_{i_1}(N_1)\alpha_{i_1} \cup N_2 f_{i_2}(N_2)\alpha_{i_2} \\
&= \left(xf_{i_1}\right)\alpha_{i_1} \cup \left(xf_{i_2}\right)\alpha_{i_2}.
\end{aligned}$$

Thus $\left(xf_{i_1}\right) \cup \left(xf_{i_2}\right)$ is bidivisible by $x^{k_{i_1}} \cup x^{k_{i_2}}$ and some bi deg $(f_{i_1}, f_{i_2}) > (k_{i_1}, k_{i_2})$ this imply that

$$f_{i_1} \cup f_{i_2} = C_{i_1}^1 x^{k_{i_1}-1} \cup C_{i_2}^2 x^{k_{i_2}-1}$$

where $C_{i_1}^1 \cup C_{i_2}^2$ is some biscalar, but then $\alpha = \alpha_1 \cup \alpha_2$.

$$C_1^1\left(x^{k_{i_1-1}}\alpha_1^1\right)+\ldots+C_{r_1}^1\left(x^{k_{i_1-1}}\alpha_{r_1}^1\right) \cup C_1^2\left(x^{k_{i_2-1}}\alpha_1^2\right)+\ldots+C_{r_2}^2\left(x^{k_{i_2-1}}\alpha_{r_2}^2\right)$$

which shows that all the bivectors form a bibasis for the strong neutrosophic binull space of $N = N_1 \cup N_2$.

Suppose T is a bilinear operator on $V = V_1 \cup V_2$ and that T factors over the neutrosophic bifield $F = F_1 \cup F_2$ as $f = f_1 \cup f_2$

$$= \left(x - C_1^1\right)^{d_1^1}\ldots\left(x - C_{k_1}^1\right)^{d_{k_1}^1} \cup \left(x - C_1^2\right)^{d_1^2}\ldots\left(x - C_{k_2}^2\right)^{d_{k_2}^2}$$

where $\{C_1^1\ldots C_{k_1}^1\} \cup \{C_1^2\ldots C_{k_2}^2\}$ are bidistinct bielement of $F = F_1 \cup F_2$ and $d_{i_t}^t \geq 1$; $t = 1,2$.

Then the neutrosophic biminimal polynomial for T will be

$$p = \left(x - C_1^1\right)^{r_1^1}\ldots\left(x - C_{k_1}^1\right)^{r_{k_1}^1} \cup \left(x - C_1^2\right)^{r_1^2}\ldots\left(x - C_{k_2}^2\right)^{r_{k_2}^2}$$

where $1 \leq r_{i_t}^t \leq d_{i_t}^t$; $t = 1, 2$.

If $W_{i_1}^1 \cup W_{i_2}^2$ is the strong neutrosophic binull space of

$$(T - C_i I)^{r_i} = \left(T_1 - C_1^1 I_1\right)^{r_{i_1}^1} \cup \left(T_2 - C_1^2 I_2\right)^{r_{i_2}^2}$$

then the biprimary decomposition theorem tells us that

$$V = V_1 \cup V_2 = \{W_1^1 \oplus \ldots \oplus W_{k_1}^1\} \cup \{W_1^2 \oplus \ldots \oplus W_{k_2}^2\}$$

and that the operator $T_{i_t}^t$ induced on $W_{i_t}^t$ defined by $T_t^{i_t}$ has neutrosophic biminimal polynomial $\left(x - C_{i_t}^t\right)^{r_i^t}$ for $t = 1, 2$; $1 \leq i_t \leq k_t$. Let $N_{i_t}^t$ be the bilinear operator on $W_{i_t}^t$ defined by $N_{i_t}^t = T_{i_t}^t - C_{i_t}^t I_t$; $1 \leq i_t \leq k_t$ then $N_{i_t}^t$ is binilpotent and has neutrosophic biminimal polynomial $x_{i_t}^{r_{i_t}^t}$. On $W_{i_t}^t$, $T_t$ acts like $N_{i_t}^t$ plus the scalar $C_{i_t}^t$ times the identity operator. Suppose we choose a bibasis for the strong neutrosophic bisubspace $W_{i_1}^1 \cup W_{i_2}^2$ corresponding to the bicyclic decomposition for the binilpotent $N_{i_t}^t$. Then the neutrosophic k matrix $T_{i_t}^t$ in this bibasis will be the bidirect sum of neutrosophic bimatrices;

$$\begin{bmatrix} C_1 & 0 & \cdots & 0 & 0 \\ 1 & C_1 & \cdots & 0 & 0 \\ \vdots & \vdots & & \vdots & \vdots \\ 0 & 0 & \cdots & C_1 & 0 \\ 0 & 0 & \cdots & 1 & C_1 \end{bmatrix} \cup \begin{bmatrix} C_2 & 0 & \cdots & 0 & 0 \\ 1 & C_2 & \cdots & 0 & 0 \\ \vdots & \vdots & & \vdots & \vdots \\ 0 & 0 & \cdots & C_2 & 0 \\ 0 & 0 & \cdots & 1 & C_2 \end{bmatrix}$$

each with $C = C_{i_t}^t$ for $t = 1, 2$. Further more the sizes of these neutrosophic bimatrices will decrease as one reads from left to right. A neutrosophic bimatirx of the form described above is called a bielementary Jordan bimatrix with bicharacteristic values $C_1 \cup C_2$.

Suppose we pull the bibasis for $W_{i_1}^1 \cup W_{i_2}^2$ together and obtain an biordered bibasis for $V = V_1 \cup V_2$. Let us describe the neutrosophic bimatrix A of T in the bibasis. A neutrosophic bimatrix A is the bidirect sum

$$A = \begin{bmatrix} A_1^1 & 0 & \cdots & 0 \\ 1 & A_2^1 & \cdots & 0 \\ \vdots & \vdots & & \vdots \\ 0 & 0 & \cdots & A_{k_1}^1 \end{bmatrix} \cup \begin{bmatrix} A_1^2 & 0 & \cdots & 0 \\ 1 & A_2^2 & \cdots & 0 \\ \vdots & \vdots & & \vdots \\ 0 & 0 & \cdots & A_{k_2}^2 \end{bmatrix}$$

of the $k_i$ sets of neutrosophic bimatrices

$$\left\{ A_1^1 ... A_{k_1}^1 \right\} \cup \left\{ A_1^2 ... A_{k_2}^2 \right\}.$$

Each

$$A_{i_1}^t = \begin{bmatrix} J_{t_1}^{i_t} & 0 & \cdots & 0 \\ 0 & J_{t_2}^{i_t} & \cdots & 0 \\ \vdots & \vdots & & \vdots \\ 0 & 0 & \cdots & J_{t_n}^{i_t} \end{bmatrix}$$

where each $J_{j_t}^{i_t}$ is an elementary Jordon neutrosophic matrix with characteristic value $C_{i_t}^t$; $1 \leq i_t \leq k_t$; $t = 1, 2$. Also within each $A_{i_t}^t$; the sizes of the neutrosophic matrices $J_{j_t}^{i_t}$ decrease as $j_t$ increases $1 \leq j_t \leq n_t$; $t = 1, 2$.

A $(n_1 \times n_1, n_2 \times n_2)$ neutrosophic bimatrix A which satisfies all the conditions described so far for some bisets of distinct $k_i$ scalars $\{C_1^1...C_{k_1}^1\} \cup \{C_1^2...C_{k_2}^2\}$ will be said to be in Jordan biform or biJordan form.

## 2.4 Neutrosophic Biinner Product Bivector Space

Now we proceed onto define the new notion of biinner product strong neutrosophic bivector space of type II and derive a few interesting properties about them.

**DEFINITION 2.4.1:** *Let $F = F_1 \cup F_2$ be a real neutrosophic bifield and $V = V_1 \cup V_2$ be a strong neutrosophic bivector space over the bineutrosophic bifield. An biinner product on V is a bifunction which assigns to each biordered pair of bivectors $\alpha = \alpha_1 \cup \alpha_2$ and $\beta = \beta_1 \cup \beta_2$ in V a biscalar $(\alpha/\beta) = (\alpha_1/\beta_1) \cup (\alpha_2/\beta_2)$ in $F = F_1 \cup F_2$ that is $(\alpha_i/\beta_i) \in F_i$, $i=1, 2$ in such a way that for all $\alpha = \alpha_1 \cup \alpha_2$, $\beta = \beta_1 \cup \beta_2$, and $\gamma = \gamma_1 \cup \gamma_2$ in $V = V_1 \cup V_2$ and for all biscalar $c = c_1 \cup c_2$ in $F_1 \cup F_2 = F$.*

i. *$(\alpha+\beta/\gamma) = (\alpha/\gamma) + (\beta/\gamma)$*
 *$(\alpha_1+\beta_1/\gamma_1) \cup (\alpha_2+\beta_2/\gamma_2)$*
 *$= (\alpha_1/\gamma_1) + (\beta_1/\gamma_1) \cup (\beta_2/\gamma_2) + (\beta_2/\gamma_2)$*
ii. *$(c\alpha/\beta) = c(\alpha/\beta)$*
 *that is $(c_1 \alpha_1/\beta_1) \cup (c_2 \alpha_2/\beta_2)$*
 *$= c_1 (\alpha_1/\beta_1) \cup c_2 (\alpha_2/\beta_2)$*
iii. *$(\alpha/\beta) = (\beta/\alpha)$*
iv. *$(\alpha/\alpha) = (\alpha_1/\alpha_1) \cup (\alpha_2/\alpha_2) > 0 \cup 0$ if $\alpha_i \neq 0$ for $i=1, 2$.*

*A strong neutrosophic bivector space $V = V_1 \cup V_2$ endowed with a biinner product is called the strong neutrosophic biinner product space over the real neutrosophic bifield $F = F_1 \cup F_2$.*

*Let $F = F_1 \cup F_2$ and for $V = F_1^{n_1} \cup F_2^{n_2}$ a strong neutrosophic bivector space over the real neutrosophic bifield F*

$= F_1 \cup F_2$, there is a biinner product called the bistandard inner product. It is defined for

$$\alpha = \alpha_1 \cup \alpha_2 = \left\{ x_1^1 \ldots x_{n_1}^1 \right\} \cup \left\{ x_1^2 \ldots x_{n_2}^2 \right\}$$

and

$$\beta = \beta_1 \cup \beta_2 = \left\{ y_1^1 \ldots y_{n_1}^1 \right\} \cup \left\{ y_1^2 \ldots y_{n_2}^2 \right\} \text{ by}$$

$$(\alpha/\beta) = \sum_{j_1} x_{j_1}^1 y_{j_1}^1 \cup \sum_{j_2} x_{j_2}^2 y_{j_2}^2 .$$

If $A = A_1 \cup A_2$ is a neutrosophic bimatrix over the bifield $F = F_1 \cup F_2$ where $A_i \in F_i^{n_i \times n_i}$ for $i = 1, 2$. $F_i^{n_i \times n_i}$ is a strong neutrosophic vector space over $F_i$; $i = 1, 2$. $V = F_1^{n_1 \times n_1} \cup F_2^{n_2 \times n_2}$ is a strong neutrosophic bivector space over the neutrosophic bifield $F = F_1 \cup F_2$ and $V$ is isomorphic to the strong neutrosophic bivector space $F_1^{n_1^2} \cup F_2^{n_2^2}$ in a natural way. It therefore follows;

$$(A/B) = \sum_{j_1 k_1} A_{j_1 k_1}^1 B_{j_1 k_1}^1 \cup \sum_{j_2 k_2} A_{j_2 k_2}^2 B_{j_2 k_2}^2$$

defines a biinner product on $V$. A strong neutrosophic bivector space over the neutrosophic bifield $F = F_1 \cup F_2$ (both $F_1$ and $F_2$ not pure neutrosophic) is known as the biinner product neutrosophic space or neutrosophic biinner product space.

We have the following interesting theorem.

**THEOREM 2.4.1:** *If $V = V_1 \cup V_2$ be a real biinner product neutrosophic space, then for any bivectors $\alpha = \alpha_1 \cup \alpha_2$ and $\beta = \beta_1 \cup \beta_2$ in $V$ and any scalar $c = c_1 \cup c_2$.*

  i.  $\| c\alpha \| = |c| \| \alpha \|$
      that is $\| c\alpha \| = \| c_1 \alpha_1 \| \cup \| c_2 \alpha_2 \|$
      $= |c_1| \| \alpha_1 \| \cup |c_2| \| \alpha_2 \|$;
  ii. $\|\alpha\| > (0 \cup 0)$ for $\alpha \neq 0$
      that is $\| \alpha_1 \| \cup \| \alpha_2 \| > (0, 0) = 0 \cup 0$;
  iii. $\|(\alpha/\beta)\| < \| \alpha \| \| \beta \|$
       that is $\|(\alpha_1/\beta_1)\| \cup \|(\alpha_2/\beta_2)\|$
       $= \| \alpha_1 \| \| \beta_1 \| \cup \| \alpha_2 \| \| \beta_2 \|$.

*However as in case of usual bivector spaces we in case of strong neutrosophic bivector spaces define the concept of biorthogonal bivectors.*

*If $\alpha, \beta \in V = V_1 \cup V_2$ be bivectors of a biinner product space, we can define*

$$\gamma = \beta - \frac{(\beta/\alpha)}{||\alpha||^2}\alpha \, ;$$

$$\gamma_1 \cup \gamma_2 = \left(\beta_1 - \frac{(\beta_1/\alpha_1)}{||\alpha_1||^2}\alpha_1\right) \cup \left(\beta_2 - \frac{(\beta_2/\alpha_2)}{||\alpha_2||^2}\alpha_2\right).$$

As in case of usual vector spaces we can in case of strong neutrosophic biinner product spaces define biorthogonality or biorthogonal bivectors.

Let $\alpha = \alpha_1 \cup \alpha_2$ and $\beta = \beta_1 \cup \beta_2$ be neutrosophic bivectors in a neutrosophic biinner product space $V = V_1 \cup V_2$.

Then $\alpha = \alpha_1 \cup \alpha_2$ is biorthogonal to $\beta = \beta_1 \cup \beta_2$ if $(\alpha/\beta) = \alpha_1 \beta_1 \cup \alpha_2 \beta_2 = 0 \cup 0$ that is $(\alpha_1/\beta_1) \cup (\alpha_2/\beta_2) = 0 \cup 0$. Since this implies $\beta = \beta_1 \cup \beta_2$ is biorthogonal to $\alpha = \alpha_1 \cup \alpha_2$.

It is left as an exercise for the reader to prove, the following theorem.

**THEOREM 2.4.2:** *A biorthogonal biset of non zero bivectors is bilinearly independent.*

**THEOREM 2.4.3:** *Let $V = V_1 \cup V_2$ be a strong neutrosophic biinner product space and let $\{\beta_1^1 \ldots \beta_{n_1}^1\} \cup \{\beta_1^2 \ldots \beta_{n_2}^2\}$ be any biindependent vector in V. Then one way to construct biorthogonal vectors;*

*$\{\alpha_1^1 \ldots \alpha_{n_1}^1\} \cup \{\alpha_1^2 \ldots \alpha_{n_2}^2\}$ in $V = V_1 \cup V_2$ is such that for each $k_i$, i=1, 2 the biset $\{\alpha_1^1 \ldots \alpha_{k_1}^1\} \cup \{\alpha_1^2 \ldots \alpha_{k_2}^2\}$ is the bibasis for the strong neutrosophic bisubspace spanned by $\{\beta_1^1 \ldots \beta_{k_1}^1\} \cup \{\beta_1^2 \ldots \beta_{k_2}^2\}$.*

*Proof:* The bivectors $\{\alpha_1^1 \ldots \alpha_{n_1}^1\} \cup \{\alpha_1^2 \ldots \alpha_{n_2}^2\}$ can be obtained by means of a construction analogous to Gram-Schmidt orthogonalization process called or defined as Gram-Schmidt biorthogonalization process.

First let $\alpha = \alpha_1^1 \cup \alpha_1^2$ and $\beta_1 = \beta_1^1 \cup \beta_1^2$. The other bivector is calculated using the rule

$$\gamma = \beta - \frac{(\beta/\alpha)}{\|\alpha\|^2} \alpha$$

$$\gamma = \gamma_1 \cup \gamma_2 = \left( \beta_1 - \frac{(\beta_1/\alpha_1)}{\|\alpha_1\|^2} \alpha_1 \right) \cup \left( \beta_2 - \frac{(\beta_2/\alpha_2)}{\|\alpha_2\|^2} \alpha_2 \right)$$

However we will indicate the proof of the result for any general n; n ≥ 3 in chapter 3 of this book.

We cannot define orthonormality as $i \in F_i$ as well as $i \in V_i$ $i = 1, 2$.

However we define just biapproximation in strong neutrosophic bivector spaces over neutrosophic bifields of type II.

**DEFINITION 2.4.2:** *Let $V = V_1 \cup V_2$ be a strong neutrosophic bivector space over the neutrosophic bifield $F = F_1 \cup F_2$ (Both $F_1$ and $F_2$ are not pure neutrosophic) of type II. Let $W = W_1 \cup W_2$ be a strong neutrosophic bivector subspace of V over the neutrosophic bifield $F = F_1 \cup F_2$.*

*Let $\beta = \beta_1 \cup \beta_2$ be a bivector in $V = V_1 \cup V_2$. To find the bibest approximation to $\beta = \beta_1 \cup \beta_2$ (or the best biapproximation to $\beta = \beta_1 \cup \beta_2$) in $W = W_1 \cup W_2$. This means to find a bivector $\alpha = \alpha_1 \cup \alpha_2$ for which $\|\beta - \alpha\| = \|\beta_1 - \alpha_1\| \cup \|\beta_2 - \alpha_2\|$ is as small as possible subject to the restriction that $\alpha = \alpha_1 \cup \alpha_2$ should belong to $W = W_1 \cup W_2$; that is to be more precise.*

*A best biapproximation to $\beta = \beta_1 \cup \beta_2$ in $W = W_1 \cup W_2$ is a bivector $\alpha = \alpha_1 \cup \alpha_2$ in W such that $\|\beta - \alpha\| \leq \|\beta - \gamma\|$ that is*

$\|\beta_1 - \alpha_1\| \cup \|\beta_2 - \alpha_2\| \le \|\beta_1 - \gamma_1\| \cup \|\beta_2 - \gamma_2\|$ for every bivector $\gamma = \gamma_1 \cup \gamma_2$ in W.

**THEOREM 2.4.4:** *Let $W = W_1 \cup W_2$ be a strong neutrosophic subbispace of a strong neutrosophic biinner product space $V = V_1 \cup V_2$ and $\beta = \beta_1 \cup \beta_2$ be in $V = V_1 \cup V_2$,*

i. *The bivector $\alpha = \alpha_1 \cup \alpha_2$ in W is a best biapproximation to $\beta = \beta_1 \cup \beta_2$ by bivectors in $W = W_1 \cup W_2$ if and only if $\beta - \alpha = (\beta_1 - \alpha_1) \cup (\beta_2 - \alpha_2)$ is biorthogonal to every vector in W. That is each $\beta_i - \alpha_i$ is orthogonal to every vector in $W_i$, true for $i = 1, 2$.*
ii. *If a best biapproximation to $\beta = \beta_1 \cup \beta_2$ by bivectors in $W = W_1 \cup W_2$ exists, it is unique.*

However we cannot define the notions and properties related to biorthonormality.

We now proceed onto define biorthogonal complement of a biset of bivectors in V.

**DEFINITION 2.4.3:** *Let $V = V_1 \cup V_2$ be a strong neutrosophic bivector space over the neutrosophic bifield $F = F_1 \cup F_2$ (Both $F_1$ and $F_2$ are not pure neutrosophic) of type II be a strong neutrosophic biinner product space.*

*Let $S = S_1 \cup S_2$ be any set of bivectors in V. The biorthogonal complement of S denoted by $S^{\perp(\perp)} = S_1^\perp \cup S_2^\perp$ is the set of all bivectors in V which are biorthogonal to every bivector in S.*

Properties related with the biorthogonal set is left as an exercise for the reader to derive.

The following results are simple and hence are left for the reader to prove.

**THEOREM 2.4.5:** *Let $V = V_1 \cup V_2$ be a strong neutrosophic biinner product space, $W = W_1 \cup W_2$ a finite dimensional strong neutrosophic bisubspace and $E = E_1 \cup E_2$ be the* *biorthogonal projection of V on W. Then the bimapping $\beta \to \beta - E\beta$; that is*

$$\beta_1 \cup \beta_2 \to (\beta_1 - E_1\beta_1) \cup (\beta_2 - E_2\beta_2)$$

*is the biorthogonal projection of V on W.*

**THEOREM 2.4.6:** *Let $W = W_1 \cup W_2$ be a finite $(n_1, n_2)$ bidimensional strong neutrosophic bisubspace of the strong neutrosophic biinner product space $V = V_1 \cup V_2$ of type II and let $E = E_1 \cup E_2$ be the biorthogonal projection of V on W.*

*Then $E = E_1 \cup E_2$ is an idempotent bilinear transformation of V onto W, $W^\perp$ is the null bispace of E and $V = W \oplus W^\perp$ that is*

$$\begin{aligned} V &= V_1 \cup V_2 \\ &= W_1 \oplus W_1^\perp \cup W_2 \oplus W_2^\perp. \end{aligned}$$

**THEOREM 2.4.7:** *Under the conditions of the above theorem $I - E = I_2 - E_1 \cup I_2 - E_2$ is the biorthogonal biprojection of V on $W^\perp$. It is a biidempotent bilinear transformation of V onto $W^\perp = W_1^\perp \cup W_2^\perp$ with binull space $W = W_1 \cup W_2$.*

**THEOREM 2.4.8:** *Let $\{\alpha_1^1 \ldots \alpha_{n_1}^1\} \cup \{\alpha_1^2 \ldots \alpha_{n_2}^2\}$ be a biorthogonal set of non zero bivectors in a strong neutrosophic biinner product space $V = V_1 \cup V_2$ over $F_1 \cup F_2$ of type II.*

*If $\beta = \beta_1 \cup \beta_2$ is any bivector in $V = V_1 \cup V_2$ then*

$$\sum_{k_1} \left( \frac{|(\beta_1 / \alpha_{k_1}^1)|^2}{\|\alpha_{k_1}^1\|^2} \right) \cup \sum_{k_2} \left( \frac{|(\beta_2 / \alpha_{k_2}^2)|^2}{\|\alpha_{k_2}^2\|^2} \right) \leq \|\beta_1\|_2 \cup \|\beta_2\|_2$$

*and equality holds if and only if*

$$\beta = \sum_{k_1} \left( \frac{|\beta_1 / \alpha_{k_1}^1|}{\|\alpha_{k_1}^1\|^2} \alpha_{k_1}^1 \right) \cup \sum_{k_2} \left( \frac{|\beta_2 / \alpha_{k_2}^2|}{\|\alpha_{k_2}^2\|^2} \alpha_{k_2}^2 \right) = \beta_1 \cup \beta_2.$$

Now we proceed onto define the notion of strong neutrosophic n-vector spaces of type II, $n \geq 3$ and neutrosophic n-vector space $n \geq 3$ in chapter three.

Several problems are proposed in chapter four of this book for the interested reader.

Chapter Three

# NEUTROSOPHIC n-VECTOR SPACES

In this chapter we for the first time introduce the notion of neutrosophic n-vector spaces of both type I and type II (n ≥ 3) and discuss some of the important properties about them. This chapter has three sections. Section one introduces the notion of neutrosophic n-vector spaces. Neutrosophic strong n-vector spaces are introduced in section two and neutrosophic n-vector spaces of type II is studied in section three.

## 3.1 Neutrosophic n-Vector Spaces

In this section we introduce strong neutrosophic n-vector spaces n ≥ 3 and illustrate it by some examples and discuss some of their properties.

**DEFINITION 3.1.1:** *Let $V = V_1 \cup V_2 \cup \ldots \cup V_n$ ($n \geq 3$) be such that each $V_i$ is a neutrosophic set and is a vector space over the* *real field F; $1 \leq i \leq n$. We call $V = V_1 \cup V_2 \cup ... \cup V_n$ to be a neutrosophic n-vector space over the field F.*

We illustrate this by some examples.

***Example 3.1.1:*** Let $V = V_1 \cup V_2 \cup V_3 \cup V_4 \cup V_5$

$$= \left\{ \begin{pmatrix} a_1 & a_2 & a_3 \\ a_4 & a_5 & a_6 \end{pmatrix} \middle| a_i \in QI; 1 \leq i \leq 6 \right\} \cup$$

$$\left\{ \begin{pmatrix} a & b \\ 0 & d \end{pmatrix} \middle| a,b,c,d \in QI \right\} \cup$$

$\{(a_1, a_2, a_3, a_4, a_5) \mid a_i \in N(Q); 1 \leq i \leq 5\} \cup$

$$\left\{ \begin{pmatrix} a_1 & a_2 \\ a_3 & a_4 \\ a_5 & a_6 \\ a_7 & a_8 \end{pmatrix} \middle| a_i \in N(Q); 1 \leq i \leq 8 \right\} \cup$$

{QI[x]; all polynomial in the variable x with coefficients from the neutrosophic field QI}; V is a neutrosophic 5-vector space over the real field Q.

***Example 3.1.2:*** Let $V = V_1 \cup V_2 \cup V_3 \cup V_4 =$

$$\left\{ \begin{pmatrix} a \\ b \\ c \end{pmatrix} \middle| a,b,c \in Z_7 I \right\} \cup \{[a, b, c, d, e] \mid a, b, c, d, e \in N(Z_7)\}$$

$\cup$ {$N(Z_7)[x]$; all polynomials in the variable x with coefficients from the neutrosophic field $N(Z_7)$} $\cup$

$$\left\{ \begin{pmatrix} a_1 & a_2 & a_3 \\ a_4 & a_5 & a_6 \\ a_7 & a_8 & a_9 \end{pmatrix} \middle| a_i \in Z_7 I; 1 \le i \le 9 \right\},$$

V is a neutrosophic 4-vector space over the real field $Z_7$.

*Example 3.1.3:* Let $V = V_1 \cup V_2 \cup V_3 =$

$$\left\{ \begin{pmatrix} a & b \\ c & d \end{pmatrix} \middle| a,b,c,d \in Z_{13} I \right\} \cup$$

$\{(a_1, a_2, a_3, a_4, a_5, a_6) \mid a_i \in N(Z_{13}); 1 \le i \le 6\} \cup$

$$\left\{ \begin{pmatrix} a_1 & a_2 \\ a_3 & a_4 \\ a_5 & a_6 \\ a_7 & a_8 \end{pmatrix} \middle| a_i \in Z_{13} I; 1 \le i \le 8 \right\},$$

V is a neutrosophic 3-vector space over the real field $Z_{13}$.

*Note:* Let $V = V_1 \cup V_2 \cup \ldots \cup V_n$ be a neutrosophic n-vector space over the real field F ($n \ge 3$). If n = 2 we get the neutrosophic bivector space. Having seen examples of neutrosophic n-vector spaces ($n \ge 3$) we now proceed onto define some substructures related with them.

**DEFINITION 3.1.2:** *Let $V = V_1 \cup V_2 \cup \ldots \cup V_n$ ($n \ge 3$) be a neutrosophic n-vector space over the real field F. Let $W = W_1 \cup W_2 \cup \ldots \cup W_n \subseteq V_1 \cup V_2 \cup \ldots \cup V_n$; if W is itself a neutrosophic n-vector space over the same real field F then we call W to be a neutrosophic n-vector subspace of V over the field F (each $W_i \subseteq V_i$ is a proper subspace different {0} space and $V_i$), i = 1, 2, …, n.*

We will illustrate this definition by some examples.

***Example 3.1.4:*** Let $V = V_1 \cup V_2 \cup V_3 \cup V_4 =$

$$\left\{ \begin{pmatrix} a_1 & a_2 & a_3 & a_4 \\ a_5 & a_6 & a_7 & a_8 \\ a_9 & a_{10} & a_{11} & a_{12} \end{pmatrix} \middle| a_i \in Z_{17}I; 1 \le i \le 12 \right\} \cup$$

$$\left\{ \begin{pmatrix} a_1 & a_2 \\ a_3 & a_4 \\ a_5 & a_6 \\ a_7 & a_8 \\ a_9 & a_{10} \end{pmatrix} \middle| a_i \in Z_{17}I; 1 \le i \le 10 \right\} \cup$$

$$\left\{ \begin{pmatrix} a_1 & a_2 & a_3 \\ a_4 & a_5 & a_6 \\ a_7 & a_8 & a_9 \end{pmatrix} \middle| a_i \in Z_{17}I; 1 \le i \le 9 \right\} \cup$$

$\{Z_{17}I[x]$; all polynomials the variable x with coefficients from $Z_{17}I\}$ be a neutrosophic 4-vector space over the real field $Z_{17}$. Take $W = W_1 \cup W_2 \cup W_3 \cup W_4 =$

$$\left\{ \begin{pmatrix} a & a & a & a \\ a & a & a & a \\ a & a & a & a \end{pmatrix} \middle| a \in Z_{17}I \right\} \cup \left\{ \begin{pmatrix} a & a \\ a & a \\ b & b \\ b & b \\ c & c \end{pmatrix} \middle| a,b,c \in Z_{17}I \right\} \cup$$

$$\left\{ \begin{pmatrix} a & b & c \\ 0 & d & e \\ 0 & 0 & f \end{pmatrix} \middle| a,b,c,d,e,f \in Z_{17}I \right\} \cup$$

$$\left\{ \sum_{i=0}^{\infty} a_i x^{2i} \middle| a_i \in Z_{17}I; 0 \le i \le \infty \right\}$$

$\subseteq V_1 \cup V_2 \cup V_3 \cup V_4$. Clearly W is a neutrosophic 4-vector subspace of V over the real field $Z_{17}$.

**Example 3.1.5:** Let $V = V_1 \cup V_2 \cup V_3 \cup V_4 \cup V_5 \cup V_6 =$

$$\{N(Q)\} \cup \left\{ \begin{pmatrix} a_1 & a_2 & a_3 & a_4 & a_5 & a_6 \\ a_7 & a_8 & a_9 & a_{10} & a_{11} & a_{12} \end{pmatrix} \middle| a_i \in QI; 1 \le i \le 12 \right\}$$

$$\cup \{(x_1, x_2, x_3, x_4, x_5, x_6, x_7) \mid x_i \in QI; 1 \le i \le 7\} \cup$$

$$\left\{ \begin{pmatrix} a_1 & a_2 & a_3 \\ a_4 & a_5 & a_6 \\ a_7 & a_8 & a_9 \\ a_{10} & a_{11} & a_{12} \\ a_{13} & a_{14} & a_{15} \\ a_{16} & a_{17} & a_{18} \\ a_{19} & a_{20} & a_{21} \end{pmatrix} \middle| a_i \in N(Q); 1 \le i \le 21 \right\} \cup$$

$$\left\{ \begin{pmatrix} a_1 & a_2 & a_3 & a_4 \\ a_5 & a_6 & a_7 & a_8 \\ a_9 & a_{10} & a_{11} & a_{12} \\ a_{13} & a_{14} & a_{15} & a_{16} \end{pmatrix} \middle| a_i \in N(Q); 1 \le i \le 21 \right\} \cup$$

$\{N(Q)[x]$; all polynomials in the variable x with coefficient from $N(Q)\}$ be a neutrosophic 6-vector space over the real field Q. Consider $W = W_1 \cup W_2 \cup W_3 \cup W_4 \cup W_5 \cup W_6 = \{\{QI\}\}$

$$\cup \left\{ \begin{pmatrix} a_1 & a_2 & a_3 & a_4 & a_5 & a_6 \\ a_1 & a_2 & a_3 & a_4 & a_5 & a_6 \end{pmatrix} \middle| a_i \in QI; 1 \le i \le 6 \right\} \cup$$

$$\{x, x, x, x, x, x, x) \mid x \in QI\} \cup$$

$$\left\{ \begin{pmatrix} a & a & a \\ b & b & b \\ a & a & a \\ b & b & b \\ c & c & c \\ b & b & b \\ c & c & c \end{pmatrix} \middle| a,b,c \in QI \right\} \cup$$

$$\left\{ \begin{pmatrix} a & a & a & a \\ b & b & b & b \\ c & c & c & c \\ d & d & d & d \end{pmatrix} \middle| a,b,c,d \in N(Q) \right\} \cup$$

{QI[x] | all polynomials in the variable x with coefficients from QI} $\subseteq V_1 \cup V_2 \cup V_3 \cup V_4 \cup V_5 \cup V_6$; W is a neutrosophic 6-vector subspace of V over the real field Q.

Now having seen the examples of neutrosophic n-vector subspace of a neutrosophic n-vector spaces we proceed onto define the new notion of sub neutrosophic n-vector subspace.

**DEFINITION 3.1.3:** *Let $V = V_1 \cup V_2 \cup ... \cup V_n$ ($n \geq 3$) be a neutrosophic n-vector space over the real field F. Let $W = W_1 \cup W_2 \cup ... \cup W_n \subseteq V_1 \cup V_2 \cup ... \cup V_n$ be such that W is a neutrosophic n-vector space over a proper subfield K of F then we define W to be a sub neutrosophic n-vector subspace of V over the subfield K of F.*

We will illustrate this situation by some examples.

*Example 3.1.6:* Let $V = V_1 \cup V_2 \cup V_3 \cup V_4 =$

$$\left\{ \begin{pmatrix} a_1 & a_2 & a_3 & a_4 \\ a_5 & a_6 & a_7 & a_8 \end{pmatrix} \middle| a_i \in RI; 1 \leq i \leq 8 \right\} \cup$$

$$\left\{\begin{pmatrix} a_1 & a_2 & a_3 \\ a_4 & a_5 & a_6 \\ a_7 & a_8 & a_9 \\ a_{10} & a_{11} & a_{12} \\ a_{13} & a_{14} & a_{15} \\ a_{16} & a_{17} & a_{18} \end{pmatrix} \middle| a_i \in Q\left(\sqrt{2}, \sqrt{3}, \sqrt{5}, \sqrt{7}, \sqrt{11}, \sqrt{13}, \sqrt{19}\right)I; \; 1 \le i \le 18 \right\}$$

$\cup$ {RI[x] ; all polynomials in the variable x with coefficients from RI} $\cup$ {5 × 5 neutrosophic matrices with entries from RI} be a neutrosophic 4-vector space over the real field F = $Q(\sqrt{2}, \sqrt{3}, \sqrt{5}, \sqrt{7}, \sqrt{11}, \sqrt{13}, \sqrt{19})$.

Take W = $W_1 \cup W_2 \cup W_3 \cup W_4$ =

$$\left\{\begin{pmatrix} a_1 & a_2 & a_3 & a_4 \\ a_5 & a_6 & a_7 & a_8 \end{pmatrix} \middle| \begin{array}{l} a_i \in Q\,(\sqrt{2}, \sqrt{3}, \sqrt{5}, \sqrt{7}, \sqrt{11}, \sqrt{13}, \\ \sqrt{19}, \sqrt{23}, \sqrt{41}, \sqrt{43}, \sqrt{53})I; \; 1 \le i \le 8 \end{array}\right\} \cup$$

$$\left\{\begin{pmatrix} a_1 & a_2 & a_3 \\ a_1 & a_2 & a_3 \\ a_1 & a_2 & a_3 \\ a_1 & a_2 & a_3 \\ a_1 & a_2 & a_3 \\ a_1 & a_2 & a_3 \end{pmatrix} \middle| a_1, a_2, a_3 \in Q\left(\sqrt{2}, \sqrt{3}, \sqrt{5}, \sqrt{7}, \sqrt{11}, \sqrt{13}, \sqrt{19}\right)I \right\}$$

$$\cup \left\{ \sum_{i=0}^{\infty} a_i x^{2i} \middle| a_i \in RI[x]; 0 \le i \le \infty \right\} \cup$$

{5 × 5 neutrosophic diagonal matrices with entries from RI} $\subseteq$ $V_1 \cup V_2 \cup V_3 \cup V_4$, W is a sub neutrosophic 4-vector subspace of V over the subfield K = $Q(\sqrt{2}, \sqrt{3}, \sqrt{5})$ $\subseteq$ F = $Q(\sqrt{2}, \sqrt{3}, \sqrt{5}, \sqrt{7}, \sqrt{11}, \sqrt{13}, \sqrt{19})$.

**Example 3.1.7:** Let $V = V_1 \cup V_2 \cup V_3 \cup V_4 \cup V_5 \cup V_6 \cup V_7$ be a neutrosophic 7-vector space over the real field R where

$$V_1 = \left\{ \begin{pmatrix} a & b & c \\ d & e & f \\ g & h & i \end{pmatrix} \middle| a,b,c,d,e,f,g,h,i \in N(R) \right\},$$

$V_2 = \{N(R)[x]$; all polynomials in the variable x with coefficients from $N(R)\}$,

$$V_3 = \left\{ \begin{pmatrix} a_1 & a_2 & a_3 & a_4 & a_5 & a_6 \\ a_7 & a_8 & a_9 & a_{10} & a_{11} & a_{12} \\ a_{13} & a_{14} & a_{15} & a_{16} & a_{17} & a_{18} \end{pmatrix} \middle| a_i \in N(R); 1 \leq i \leq 18 \right\},$$

$$V_4 = \left\{ \begin{pmatrix} a & b \\ c & d \end{pmatrix} \middle| a,b,c,d \in N(R) \right\}, \quad V_5 = \{N(R)\},$$

$$V_6 = \left\{ \begin{pmatrix} a_1 & a_2 \\ a_3 & a_4 \\ a_5 & a_6 \\ a_7 & a_8 \\ a_9 & a_{10} \\ a_{11} & a_{12} \\ a_{13} & a_{14} \end{pmatrix} \middle| a_i \in RI; 1 \leq i \leq 14 \right\}$$

and $V_7 = \{$all $9 \times 9$ neutrosophic matrices with entries from RI$\}$.
Take $W = W_1 \cup W_2 \cup W_3 \cup W_4 \cup W_5 \cup W_6 \cup W_7 =$

$$\left\{ \begin{pmatrix} a & b & c \\ d & e & f \\ g & h & i \end{pmatrix} \middle| a,b,c,d,e,f,g,h,i \in RI \right\} \cup$$

{RI[x] all polynomials in the variable x with coefficients from RI} ∪

$$\left\{ \begin{pmatrix} a & a & a & a & a & a \\ b & b & b & b & b & b \\ c & c & c & c & c & c \end{pmatrix} \middle| a,b,c \in RI \right\} \cup$$

$$\left\{ \begin{pmatrix} a & b \\ c & d \end{pmatrix} \middle| a,b,c,d \in QI \right\} \cup \{RI\} \cup$$

$$\left\{ \begin{pmatrix} a_1 & a_2 \\ a_3 & a_4 \\ a_5 & a_6 \\ a_7 & a_8 \\ a_9 & a_{10} \\ a_{11} & a_{12} \\ a_{13} & a_{14} \end{pmatrix} \middle| a_i \in QI; 1 \le i \le 14 \right\} \cup$$

{all 9 × 9 neutrosophic matrices with entries from N(Q)} ⊆ $V_1$ ∪ $V_2$ ∪ $V_3$ ∪ $V_4$ ∪ $V_5$ ∪ $V_6$ ∪ $V_7$ is a sub neutrosophic 7-vector subspace of V over the subfield Q of the field R.

We define a neutrosophic n-vector space which has no proper sub neutrosophic n-vector subspace to be a subsimple neutrosophic n-vector space.

We will illustrate this by examples.

*Example 3.1.8:* Let V = $V_1$ ∪ $V_2$ ∪ $V_3$ ∪ $V_4$ ∪ $V_5$ =

$$\left\{ \begin{pmatrix} a & b \\ c & d \end{pmatrix} \middle| a,b,c,d \in QI \right\} \cup \{QI\,[x]\} \cup$$

$$\left\{ \begin{pmatrix} a_1 & a_2 \\ a_3 & a_4 \\ a_5 & a_6 \\ a_7 & a_8 \\ a_9 & a_{10} \end{pmatrix} \middle| a_i \in QI; 1 \le i \le 10 \right\} \cup$$

{10 × 10 neutrosophic matrices with entries from QI} ∪ {($a_1$, $a_2$, $a_3$) | $a_i$ ∈ RI} be a neutrosophic 5-vector space over the real field Q. Clearly V has no sub neutrosophic 5-vector subspace as Q is a prime field that Q has no proper subfields. Hence V is a subsimple neutrosophic 5-vector space over Q.

*Example 3.1.9:* Let $V = V_1 \cup V_2 \cup V_3 \cup V_4 \cup V_5 \cup V_6 = \{Z_2I \times Z_2I \times Z_2I \times Z_2I\} \cup \{Z_2I[x]$ all polynomials in the variable x with coefficients from $Z_2I\} \cup$

$$\left\{ \begin{pmatrix} 0 & 0 \\ 0 & 0 \end{pmatrix}, \begin{pmatrix} I & I \\ I & I \end{pmatrix} I + I = 2I \equiv 0 (\mathrm{mod}\, 2) \right\} \cup$$

$$\left\{ \begin{pmatrix} a & a & a \\ a & a & a \\ a & a & a \end{pmatrix} \middle| a \in N(Z_2) \right\} \cup$$

{all 3 × 7 neutrosophic matrices with entries from $N(Z_2)$} ∪ {all 9 × 4 neutrosophic matrices with entries from $N(Z_2)$} be a neutrosophic 6-vector space over the real field $Z_2$. Further it can be easily verified V has no proper neutrosophic 6-vector subspace over $Z_2$. Since $Z_2$ is a prime field of characteristic two it has no proper subfields. This V is a subsimple neutrosophic 6-vector space over $Z_2$.

**THEOREM 3.1.1:** *Let $V = V_1 \cup V_2 \cup \ldots \cup V_n$ be a neutrosophic n-vector space over a real field F. If F is a prime field that is F has no proper subfields then V is a subsimple neutrosophic n-vector space over F.*

*Proof:* Given $V = V_1 \cup V_2 \cup \ldots \cup V_n$ is a neutrosophic n-vector space over the real field F, such that F is a prime field; that is F has no proper subfields. By the definition of such neutrosophic n-vector subspaces we see V does not have a proper sub neutrosophic n-vector subspace hence V is a subsimple neutrosophic n-vector space over F.

We define the notion of doubly simple neutrosophic n-vector space over a real field F.

Let $V = V_1 \cup V_2 \cup \ldots \cup V_n$ be a neutrosophic n-vector space over a field F. If V has no proper neutrosophic n-vector subspace over the field F then we all V to be a simple neutrosophic n-vector space over the field F.

We will first illustrate this situation by some examples.

*Example 3.1.10:* Let $V = V_1 \cup V_2 \cup V_3 \cup V_4 =$

$\left\{ \begin{pmatrix} 0 & 0 \\ 0 & 0 \end{pmatrix}, \begin{pmatrix} I & I \\ I & I \end{pmatrix}, \begin{pmatrix} 2I & 2I \\ 2I & 2I \end{pmatrix} \right\|$ elements of these matrices are from

$Z_3I \} \cup \left\{ \begin{pmatrix} 0 & 0 & 0 & 0 \\ 0 & 0 & 0 & 0 \end{pmatrix}, \begin{pmatrix} I & I & I & I \\ I & I & I & I \end{pmatrix}, \begin{pmatrix} 2I & 2I & 2I & 2I \\ 2I & 2I & 2I & 2I \end{pmatrix} \right\|$

elements of these $2 \times 4$ matrices are from neutrosophic field $Z_3I\} \cup \{5 \times 5$ neutrosophic matrices with entries from $Z_3I\} \times \{(a_1, a_2, a_3, a_4, a_5, a_6, a_7) \mid a_i \in N(Z_3); 1 \leq i \leq 7\}$ be a neutrosophic 4-vector space over the field $Z_3$. Clearly V has no neutrosophic 4-vector subspaces so V is a simple neutrosophic n-vector space.

*Example 3.1.11:* Let $V = \{(a\ a\ a\ a) \mid a \in Z_5I\} \cup$

$\left\{ \begin{pmatrix} a & a \\ 0 & a \end{pmatrix} \middle| a \in Z_5I \right\} \cup$

$$\left\{ \begin{pmatrix} a & a & a \\ a & a & a \\ a & a & a \\ a & a & a \\ a & a & a \\ a & a & a \\ a & a & a \end{pmatrix} \middle| a \in Z_5 I \right\} \cup$$

$$\left\{ \begin{pmatrix} a & 0 & 0 & 0 & 0 \\ a & a & 0 & 0 & 0 \\ a & a & a & 0 & 0 \\ a & a & a & a & 0 \\ a & a & a & a & a \end{pmatrix} \middle| a \in Z_5 I \right\}$$

$= V_1 \cup V_2 \cup V_3 \cup V_4$ be a neutrosophic 4-vector space over the real field $Z_5$. V is a simple neutrosophic 4-vector space.

We define doubly simple neutrosophic vector space.

**DEFINITION 3.1.4:** *Let $V = V_1 \cup V_2 \cup \ldots \cup V_n$ be a neutrosophic n-vector space over the real field F. Suppose V is a simple neutrosophic n-vector space as well as simple subneutrosophic bivector space then we call V to be a doubly simple neutrosophic n-vector space.*

We will illustrate this situation by some simple examples.

*Example 3.1.12:* Let $V = V_1 \cup V_2 \cup V_3 \cup V_4 \cup V_5 \cup V_6 \cup V_7$ be neutrosophic 7-vector space over the real field $Z_7$ where

$$V_1 = \left\{ \begin{pmatrix} a & a & a \\ a & a & a \end{pmatrix} \middle| a \in Z_7 I \right\},$$

$$V_2 = \left\{ \begin{pmatrix} a & a & a & a \\ a & a & a & a \\ a & a & a & a \\ a & a & a & a \\ a & a & a & a \\ a & a & a & a \\ a & a & a & a \\ a & a & a & a \end{pmatrix} \middle| a \in Z_7 I \right\}, V_3 = \left\{ \begin{pmatrix} a & 0 & 0 & 0 \\ a & a & 0 & 0 \\ a & a & a & 0 \\ a & a & a & a \end{pmatrix} \middle| a \in Z_7 I \right\}$$

$$V_4 = \{(a\ a\ a\ a\ a\ a\ a\ a) \mid a \in Z_7 I\},\ V_5 = \left\{ \begin{pmatrix} a \\ a \\ a \\ a \\ a \\ a \\ a \\ a \\ a \end{pmatrix} \middle| a \in Z_7 I \right\},$$

$$V_6 = \left\{ \begin{pmatrix} a & 0 \\ a & a \end{pmatrix} \middle| a \in Z_7 I \right\}$$

and

$$V_7 = \left\{ \begin{pmatrix} a & a & a & a & a & a \\ a & a & a & a & a & a \\ a & a & a & a & a & a \\ a & a & a & a & a & a \end{pmatrix} \middle| a \in Z_7 I \right\}.$$

It is easily verified V is a simple neutrosophic 7-vector space as each $V_i$ is a simple neutrosophic vector space for i = 1, 2, …, 7. Further $Z_7$ is a prime field so V has no subneutrosophic 7 vector subspaces. Thus V is a doubly simple neutrosophic 7-vector space over $Z_7$.

***Example 3.1.13:*** Let $V = V_1 \cup V_2 \cup V_3 =$

$$\left\{ \begin{pmatrix} a & b & c \\ d & e & f \\ g & h & i \end{pmatrix} \middle| a,b,c,d,e,f,g,h,i \in Z_{17}I \right\} \cup$$

$$\left\{ \begin{pmatrix} a & a & a & a & a \\ b & b & b & b & b \end{pmatrix} \middle| a,b \in Z_{17}I \right\} \cup \left\{ \begin{pmatrix} a & a \\ a & a \end{pmatrix} \middle| a \in Z_{17}I \right\}$$

be a neutrosophic 3-vector space over the real field $Z_{17}$. Clearly V is also a doubly simple neutrosophic bivector space over the field $Z_{17}$.

A neutrosophic n-vector space can have neutrosophic n-vector subspace still it can be a simple sub neutrosophic n-vector space. This is shown by some simple examples.

***Example 3.1.14:*** Let $V = V_1 \cup V_2 \cup V_3 \cup V_4 \cup V_5$ be a neutrosophic 5-vector space over the real field $F = Z_{11}$. Here

$$V_1 = \left\{ \begin{pmatrix} a & b \\ c & d \end{pmatrix} \middle| a,b,c,d \in Z_{11}I \right\}$$

$$V_2 = \left\{ \begin{pmatrix} a & a & a & a & a \\ a & a & a & a & a \end{pmatrix} \middle| a \in N(Z_{11}) \right\},$$

$$V_3 = \left\{ \begin{pmatrix} a & a \\ b & b \\ c & c \\ d & d \\ e & e \\ f & f \end{pmatrix} \middle| a,b,c,d,e,f \in N(Z_{11}) \right\}$$

$$V_4 = \left\{ \begin{pmatrix} a_1 & 0 & 0 & 0 \\ a_2 & a_3 & 0 & 0 \\ a_4 & a_5 & a_6 & 0 \\ a_7 & a_8 & a_9 & a_{10} \end{pmatrix} \middle| a_i \in N(Z_{11}); 1 \le i \le 10 \right\}$$

and $V_5 = \{Z_{17}I[x]$; all polynomials in the variable x with coefficients from $Z_{17}I$.

Take $W = W_1 \cup W_2 \cup W_3 \cup W_4 \cup W_5 =$

$$\left\{ \begin{pmatrix} a & a \\ a & a \end{pmatrix} \middle| a \in Z_{11}I \right\} \cup$$

$$\left\{ \begin{pmatrix} a & a & a & a & a \\ a & a & a & a & a \end{pmatrix} \middle| a \in Z_{11}I \right\} \cup$$

$$\left\{ \begin{pmatrix} a & a \\ a & a \\ a & a \\ a & a \\ a & a \\ a & a \end{pmatrix} \middle| a \in Z_{11}I \right\} \cup$$

$$\left\{ \begin{pmatrix} a_1 & 0 & 0 & 0 \\ a_2 & a_3 & 0 & 0 \\ a_4 & a_5 & a_6 & 0 \\ a_7 & a_8 & a_9 & a_{10} \end{pmatrix} \middle| a_i \in Z_{11}I; 1 \le i \le 10 \right\} \cup$$

$$\left\{ \sum_{i=0}^{\infty} a_i x^{2i} \middle| a_i \in Z_{11}I; 0 \le i \le \infty \right\}$$

⊆ $V_1 \cup V_2 \cup V_3 \cup V_4 \cup V_5$; W is a neutrosophic 5-vector subspace of V. So V is not a simple neutrosophic 5-vector space, however V has no subneutrosophic subvector space as $Z_{11}$ is a prime field so V is a subsimple neutrosophic 5-vector space over $Z_{11}$.

Thus V is not a doubly simple neutrosophic 5-vector space over the field $Z_{11}$.

Now we proceed onto define the notion of neutrosophic n-linear algebra $n \geq 3$.

**DEFINITION 3.1.5:** *Let $V = V_1 \cup V_2 \cup ... \cup V_n$ be a neutrosophic n-vector space over the real field F. If each $V_i$ is a neutrosophic linear algebra over the field F then we define V to be a neutrosophic n-linear algebra over the field F.*

We illustrate this situation by some simple examples.

*Example 3.1.15:* Let $V = V_1 \cup V_2 \cup V_3 \cup V_4 =$

$$\left\{ \begin{pmatrix} a & b & c \\ d & e & f \\ g & h & i \end{pmatrix} \middle| a,b,...,h,i \in Z_2 I \right\} \cup \left\{ \begin{pmatrix} a & b \\ c & d \end{pmatrix} \middle| a,b,c,d \in Z_2 I \right\}$$

$\cup \{(a_1, a_2, a_3, a_4) \mid a_i \in Z_2 I, 1 \leq i \leq 4\} \cup \{Z_2 I[x];$ all polynomials in the variable x with coefficients from $Z_2 I\}$ be a neutrosophic 4 linear algebra over the real field $Z_2 = \{0, 1\}$.

*Example 3.1.16:* Let $V = V_1 \cup V_2 \cup V_3 \cup V_4 \cup V_5 \cup V_6 =$

$$\left\{ \begin{pmatrix} a & b \\ c & d \end{pmatrix} \middle| a,b,c,d \in QI \right\} \cup$$

$\{N(Q)[x];$ all polynomials in the variable x with coefficients from the neutrosophic field QI$\} \cup \{5 \times 5$ neutrosophic matrices

with entries from QI} ∪ {7 × 7 neutrosophic upper triangular matrices with entries from QI} ∪

$$\left\{ \begin{pmatrix} a & a & a \\ a & a & a \\ a & a & a \end{pmatrix} \middle| a \in QI \right\} \cup$$

{N(Q)} is a neutrosophic 6-linear algebra over the real field Q.

Now we will state the following theorem. The reader is expected to prove it.

**THEOREM 3.1.2:** *Let $V = V_1 \cup V_2 \cup ... \cup V_n$ be a neutrosophic n-linear algebra defined over the real field F. Every neutrosophic n-linear algebra is a neutrosophic n-vector space. But in general a neutrosophic n-vector space need not be a neutrosophic n-linear algebra.*

We give an example of a neutrosophic n-vector space which is not a neutrosophic n-linear algebra.

*Example 3.1.17:* Let $V = V_1 \cup V_2 \cup V_3 \cup V_4 \cup V_5 =$

$$\left\{ \begin{pmatrix} a & a & a & a & a \\ a & a & a & a & a \\ a & a & a & a & a \end{pmatrix} \middle| a \in Z_7 I \right\} \cup$$

$$\left\{ \begin{pmatrix} a & b \\ c & d \\ e & f \\ g & h \\ a & b \end{pmatrix} \middle| a,b,c,d,e,f,g,h \in Z_7 I \right\} \cup$$

$$\left\{ \begin{pmatrix} a \\ b \\ c \\ d \\ e \end{pmatrix} \middle| a,b,c,d,e \in N(Z_7) \right\} \cup$$

$$\left\{ \begin{pmatrix} a & 0 & 0 \\ 0 & b & d \\ e & f & g \end{pmatrix} \middle| a,b,c,d,e,f,g \in N(Z_7) \right\} \cup$$

$$\left\{ \begin{pmatrix} 0 & x \\ y & 0 \end{pmatrix} \middle| x,y \in N(Z_7) \right\}$$

be a neutrosophic 5-vector space over the real field $Z_7$. Clearly V is not a neutrosophic 5-linear algebra over the real field $Z_7$. For we see in $V_1$ we cannot define product so $V_1$ is not a neutrosophic linear algebra over $Z_7$.

$$V_2 = \left\{ \begin{pmatrix} a & b \\ c & d \\ e & f \\ g & h \\ a & b \end{pmatrix} \middle| a,b,c,d,e,f,g,h \in Z_7 I \right\}$$

is not a neutrosophic linear algebra over $Z_7$ a product in $V_2$ cannot be defined only addition is valid. $V_4$ is a neutrosophic linear algebra over $Z_7$. However $V_5$ and $V_3$ are not neutrosophic linear algebras. Thus $V = V_1 \cup V_2 \cup V_3 \cup V_4 \cup V_5$ is not a neutrosophic 5-linear algebra over the field $Z_7$. Hence the claim.

We now proceed onto define the notion of neutrosophic n-linear subalgebra of a neutrosophic n-linear algebra.

**DEFINITION 3.1.6:** *Let $V = V_1 \cup V_2 \cup \ldots \cup V_n$ be a neutrosophic n-linear algebra over the real field F. Let $W = W_1 \cup W_2 \cup \ldots \cup W_n \subseteq V_1 \cup V_2 \cup \ldots \cup V_n$ be a neutrosophic n-linear algebra over the field F then we call W to be neutrosophic n-linear subalgebra of V over the field F.*

We will illustrate this situation by some examples.

*Example 3.1.18:* Let $V_1 \cup V_2 \cup V_3 \cup V_4 \cup V_5 =$

$$\left\{ \begin{pmatrix} a & b \\ c & d \end{pmatrix} \middle| a,b,c,d \in Z_{11}I \right\} \cup$$

$$\{(x_1, x_2, x_3, x_4, x_5) \mid x_i \in N(Z_{11}); 1 \leq i \leq 5\} \cup$$

$$\left\{ \begin{pmatrix} a & 0 & 0 \\ b & c & 0 \\ d & d & f \end{pmatrix} \middle| a,b,c,d,e,f,g \in Z_{11}I \right\} \cup$$

$\{Z_{11}I[x];$ all polynomials in the variable x with coefficients from $Z_{11}I\} \cup \{10 \times 10$ neutrosophic matrices with entries from $Z_{11}I\}$ be a neutrosophic 5-linear algebra over the real field $Z_{11}$. Take $W = W_1 \cup W_2 \cup W_3 \cup W_4 \cup W_5 =$

$$\left\{ \begin{pmatrix} a & a \\ a & a \end{pmatrix} \middle| a \in Z_{11}I \right\} \cup$$

$$\{(x_1, x_2, x_3, x_4, x_5) \mid x_i \in Z_{11}I; 1 \leq i \leq 5\} \cup$$

$$\left\{ \begin{pmatrix} a & 0 & 0 \\ a & a & 0 \\ a & a & a \end{pmatrix} \middle| a \in Z_{11}I \right\} \cup \left\{ \sum_{i=0}^{\infty} a_i x^{2i} \middle| a_i \in Z_{11}I; 0 \leq i \leq \infty \right\} \cup$$

{all 10 × 10 upper triangular matrices with entries from $Z_{11}I$} $\subseteq$ $V_1 \cup V_2 \cup V_3 \cup V_4 \cup V_5$ W is a neutrosophic 5-linear subalgebra of V over the real field $Z_{11}$.

***Example 3.1.19:*** Let $V = V_1 \cup V_2 \cup V_3 \cup V_4 =$

$$\left\{ \begin{pmatrix} a & 0 & 0 \\ b & b & 0 \\ c & c & c \end{pmatrix} \middle| a,b,c \in Z_{17}I \right\} \cup$$

$$\left\{ \sum_{i=0}^{\infty} a_i x^i \middle| a_i \in Z_{17}I; 0 \le i \le \infty \right\} \cup$$

$$\left\{ \begin{pmatrix} a & b & c \\ 0 & d & e \\ 0 & 0 & f \end{pmatrix} \middle| a,b,c,d,e,f \in Z_{17}I \right\} \cup$$

$$\left\{ \begin{pmatrix} a_1 & a_2 & a_3 & a_4 \\ a_5 & a_6 & a_7 & a_8 \\ a_9 & a_{10} & a_{11} & a_{12} \\ a_{13} & a_{14} & a_{15} & a_{16} \end{pmatrix} \middle| a_i \in N(Z_{17}); 1 \le i \le 16 \right\}$$

be a neutrosophic 4-linear algebra over the real field $Z_17$. Take $W = W_1 \cup W_2 \cup W_3 \cup W_4 =$

$$\left\{ \begin{pmatrix} a & 0 & 0 \\ a & a & 0 \\ a & a & a \end{pmatrix} \middle| a \in Z_{17}I \right\} \cup$$

$$\left\{ \sum_{i=0}^{\infty} a_i x^{2i} \middle| a_i \in Z_{17}I; 0 \le i \le \infty \right\} \cup$$

$$\left\{ \begin{pmatrix} a & a & a \\ 0 & a & a \\ 0 & 0 & a \end{pmatrix} \middle| a \in Z_{17}I \right\} \cup$$

$$\left\{ \begin{pmatrix} a & a & a & a \\ b & b & b & b \\ c & c & c & c \\ d & d & d & d \end{pmatrix} \middle| a,b,c,d \in Z_{17}I \right\}$$

$\subseteq V_1 \cup V_2 \cup V_3 \cup V_4$; W is a neutrosophic 4-linear subalgerba of V over the real field $Z_{17}$.

We see in general all neutrosophic n-linear algebras need not have neutrosophic n-linear subalgebras.

Suppose we have a neutrosophic n-linear algebra V which no proper neutrosophic n-linear subalgebra then we call V to be a simple neutrosophic n-linear algebra.

We will illustrate this situation by some examples.

*Example 3.1.20:* Let $V = V_1 \cup V_2 \cup V_3 \cup V_4 \cup V_5 =$

$$\left\{ \begin{pmatrix} a & a \\ a & a \end{pmatrix} \middle| a \in Z_{19}I \right\} \cup \left\{ \begin{pmatrix} a & a & a \\ 0 & a & a \\ 0 & 0 & a \end{pmatrix} \middle| a \in Z_{19}I \right\} \cup$$

$$\left\{ \begin{pmatrix} a & a & a & a & a \\ a & a & a & a & a \\ a & a & a & a & a \\ a & a & a & a & a \\ a & a & a & a & a \end{pmatrix} \middle| a \in Z_{19}I \right\} \cup$$

$\{(a\ a\ a\ a\ a\ a\ a\ a) \mid a \in Z_{19}\} \cup$

$$\left\{ \begin{pmatrix} a & 0 & 0 & 0 \\ a & a & 0 & 0 \\ a & a & a & 0 \\ a & a & a & a \end{pmatrix} \middle| a \in Z_{19}I \right\}$$

be a neutrosophic 5-linear algebra over the real field $Z_{19}$. We see V has no neutrosophic 5-linear subalgebra over the field $Z_{19}$. Thus V is a simple neutrosophic 5-linear algebra over $Z_{19}$.

*Example 3.1.21:* Let $V = V_1 \cup V_2 \cup V_3 \cup V_4 \cup V_5 \cup V_6 \cup V_7 \cup V_8 =$

$$\left\{ \begin{pmatrix} a & 0 & 0 & 0 \\ a & a & 0 & 0 \\ a & a & a & 0 \\ a & a & a & a \end{pmatrix} \middle| a \in Z_7 I \right\} \cup$$

$$\left\{ \sum_{i=0}^{\infty} a_i x^i \middle| a_i \in Z_7 I; 0 \leq i \leq \infty \right\} \cup \left\{ \begin{pmatrix} a & a & a \\ a & a & a \\ a & a & a \end{pmatrix} \middle| a \in N(Z_7 I) \right\} \cup$$

{$N(Z_7I)$} $\cup$ {all 9 × 9 upper triangular matrices with entries from $Z_7I$} $\cup$

$$\left\{ \begin{pmatrix} a & a & a & a & a & a \\ a & a & a & a & a & a \\ a & a & a & a & a & a \\ a & a & a & a & a & a \\ a & a & a & a & a & a \\ a & a & a & a & a & a \end{pmatrix} \middle| a \in Z_7 I \right\} \cup$$

{all 10 × 10 lower triangular neutrosophic matrices with entries with entries from $Z_7I$} $\cup$

$$\left\{ \begin{pmatrix} a & a \\ a & a \end{pmatrix} \middle| a \in Z_7 I \right\}$$

be a neutrosophic 8-linear algebra over the real field $Z_7$. Clearly V is a simple neutrosophic 8-linear algebra as the neutrosophic linear algebras $V_1$, $V_3$, $V_6$ and $V_8$ are simple neutrosophic linear algebras over the real field $Z_7$.

Now we proceed onto define yet another new substructures in neutrosophic n-linear algebras.

**DEFINITION 3.1.7:** *Let $V = V_1 \cup V_2 \cup ... \cup V_n$ be a neutrosophic n-linear algebra a real field F. Suppose $W = W_1 \cup W_2 \cup ... \cup W_n \subseteq V_1 \cup V_2 \cup ... \cup V_n$ be a proper n-subset of V such that W is a neutrosophic n-linear algebra over a proper subfield K of F then we define W to be subneutrosophic n-linear subalgebra of V over the subfield K of the field F.*

We will illustrate this by some examples.

*Example 3.1.22:* Let $V = V_1 \cup V_2 \cup V_3 \cup V_4 \cup V_5 \cup V_6 =$

$$\left\{ \begin{pmatrix} a & b \\ c & d \end{pmatrix} \middle| a,b,c,d \in RI \right\} \cup \{RI\} \cup \left\{ \begin{pmatrix} a & 0 & 0 \\ a & a & 0 \\ a & a & a \end{pmatrix} \middle| a \in N(R) \right\} \cup$$

$$\left\{ \sum_{i=0}^{\infty} a_i x^i \middle| a_i \in RI; 0 \leq i \leq \infty \right\} \cup$$

$$\left\{ \begin{pmatrix} a & a & a & a \\ b & b & b & b \\ c & c & c & c \\ d & d & d & d \end{pmatrix} \middle| a,b,c,d \in RI \right\} \cup$$

{11 × 11 neutrosophic matrices with entries from the neutrosophic field RI} be a neutrosophic 6-linear algebra over the real field R, the field of reals. Take $W = W_1 \cup W_2 \cup W_3 \cup W_4 \cup W_5 \cup W_6 =$

$$\left\{ \begin{pmatrix} a & a \\ a & a \end{pmatrix} \middle| a \in RI \right\} \cup \{QI\} \cup$$

$$\left\{ \begin{pmatrix} a & 0 & 0 \\ a & a & 0 \\ a & a & a \end{pmatrix} \middle| a \in RI \right\} \cup \left\{ \begin{pmatrix} a & a & a & a \\ a & a & a & a \\ a & a & a & a \\ a & a & a & a \end{pmatrix} \middle| a \in RI \right\} \cup$$

$$\left\{ \sum_{i=0}^{\infty} a_i x^{2i} \middle| a_i \in QI; 0 \leq i \leq \infty \right\} \cup$$

{all 11 × 11 neutrosophic matrices with entries from QI} $\subseteq V_1 \cup V_2 \cup V_3 \cup V_4 \cup V_5 \cup V_6$, W is a subneutrosophic 6 linear algebra over the real field $Q \subseteq R$.

***Example 3.1.23:*** Let $V = V_1 \cup V_2 \cup V_3 \cup V_4 \cup V_5 =$

$$\left\{ \begin{pmatrix} a & b \\ c & d \end{pmatrix} \middle| a,b,c,d \in RI \right\} \cup$$

$$\left\{ \begin{pmatrix} a & 0 & 0 \\ b & c & 0 \\ d & e & 0 \end{pmatrix} \middle| a,b,c,d,e \text{ are in } N(RI) \right\} \cup$$

$$\left\{ \sum_{i=0}^{\infty} a_i x^i \middle| a_i \in RI; 0 \leq i \leq \infty \right\}$$

∪ {N(R)} ∪ {7 × 7 neutrosophic matrices with entries from RI}
is a neutrosophic 5-linear algebra over the real field R.

Take $W = W_1 \cup W_2 \cup W_3 \cup W_4 \cup W_5 =$

$$\left\{ \begin{pmatrix} a & b \\ c & d \end{pmatrix} \middle| a,b,c,d \in QI \right\} \cup \left\{ \begin{pmatrix} a & 0 & 0 \\ b & c & 0 \\ d & e & 0 \end{pmatrix} \middle| a,b,c,d,e \in N(Q) \right\} \cup$$

$$\left\{ \sum_{i=0}^{\infty} a_i x^i \middle| a_i \in QI; 0 \le i \le \infty \right\} \cup$$

{N(Q)} ∪ {7 × 7 neutrosophic matrices with entries from QI} ⊆ $V_1 \cup V_2 \cup V_3 \cup V_4 \cup V_5$ is a subneutrosophic 5-linear subalgebra of V over the subfield Q of R,

Now if a neutrosophic n-linear algebra V has no proper subneutrosophic linear subalgebra over a subfield K of F (V is defined over F), then we call V to be subsimple. neutrosophic n-linear algebra.

We will illustrate this situation by some simple examples.

*Example 3.1.24:* Let $V = V_1 \cup V_2 \cup V_3 \cup V_4 =$

$$\left\{ \begin{pmatrix} a & b \\ c & d \end{pmatrix} \middle| a,b,c,d \in Z_7 I \right\} \cup$$

$$\left\{ \sum_{i=0}^{\infty} a_i x^i \middle| a_i \in Z_7 I; 0 \le i \le \infty \right\} \cup$$

$$\left\{ \begin{pmatrix} a_1 & 0 & 0 & 0 \\ a_2 & a_3 & 0 & 0 \\ a_4 & a_5 & a_6 & 0 \\ a_7 & a_8 & a_9 & a_{10} \end{pmatrix} \middle| a_i \in Z_7 I; 1 \le i \le 10 \right\} \cup$$

$\{N(Z_7)\}$ be a neutrosophic 4-linear algebra over the real field $Z_7$. Since $Z_7$ has no proper subfields that is as $Z_7$ is a prime field we see V has no subneutrosophic 4-linear subalgebras. Hence V is a subsimple 4-linear algebra.

***Example 3.1.25:*** Let $V = V_1 \cup V_2 \cup V_3 \cup V_4 \cup V_5 =$

$$\left\{ \begin{pmatrix} a_1 & a_2 & a_3 \\ a_4 & a_5 & a_6 \\ a_7 & a_8 & a_9 \end{pmatrix} \middle| a_i \in QI; 1 \leq i \leq 9 \right\} \cup$$

$$\left\{ \sum_{i=0}^{\infty} a_i x^i \middle| a_i \in QI; 0 \leq i \leq \infty \right\} \cup \{N(R)\} \cup$$

$$\left\{ \begin{pmatrix} a_1 & 0 & 0 & 0 \\ a_2 & a_3 & 0 & 0 \\ a_4 & a_5 & a_6 & 0 \\ a_7 & a_8 & a_9 & a_{10} \end{pmatrix} \middle| a_i \in RI \quad 0 \leq i \leq 10 \right\} \cup$$

{All $10 \times 10$ neutrosophic matrices with entries from RI} is a neutrosophic 5-linear algebra over the field Q. Clearly Q is a prime field so V has no subneutrosophic 5-linear algebra, hence V is a subsimple neutrosophic 5-linear algebra.

In view of this example we have nice theorem which gurantees the existence of subsimple neutrosophic n-linear algebras.

**THEOREM 3.1.3:** *Let $V = V_1 \cup V_2 \cup \ldots \cup V_n$ be a neutrosophic n-linear algebra over a real field F, where F is a prime field i.e., has no subfields then V is a subsimple neutrosophic n-linear algebra.*

*Proof:* Follows from the fact that $V = V_1 \cup \ldots \cup V_n$ is defined over the prime field F for a subneutrosophic n-linear algebra to exist we need the existence of a subfield in F. Hence V is a subneutrosophic simple n-linear algebra.

A simple neutrosophic n-linear algebra need not in general be a simple subneutrosophic n-linear algebra.

*Example 3.1.26:* Let $V = V_1 \cup V_2 \cup V_3 \cup V_4 \cup V_5 =$

$$\left\{ \begin{pmatrix} a & b \\ c & d \end{pmatrix} \middle| a,b,c,d \in N(Z_{11}) \right\} \cup$$

$$\left\{ \begin{pmatrix} a_1 & 0 & 0 \\ a_2 & a_3 & 0 \\ a_4 & a_5 & a_6 \end{pmatrix} \middle| a_i \in Z_{11}I; 1 \leq i \leq 6 \right\} \cup$$

$$\left\{ \sum_{i=0}^{\infty} a_i x^i \middle| a_i \in Z_{11}I; 0 \leq i \leq \infty \right\} \cup$$

$$\left\{ \begin{pmatrix} a_1 & a_2 & a_3 & a_4 \\ 0 & a_5 & a_6 & a_7 \\ 0 & 0 & a_8 & a_9 \\ 0 & 0 & 0 & a_{10} \end{pmatrix} \middle| a_i \in Z_{11}I; 1 \leq i \leq 10 \right\} \cup$$

{all $7 \times 7$ matrices with entries from $Z_{11}I$} be a neutrosophic 5-linear algebra over the real field $Z_{11}$.

Take $W = W_1 \cup W_2 \cup W_3 \cup W_4 \cup W_5 =$

$$\left\{ \begin{pmatrix} a & a \\ a & a \end{pmatrix} \middle| a \in Z_{11}I \right\} \cup$$

$$\left\{ \begin{pmatrix} a & 0 & 0 \\ a & a & 0 \\ a & a & a \end{pmatrix} \middle| a \in Z_{11}I \right\} \cup$$

$$\left\{\sum_{i=0}^{\infty} a_i x^{2i} \,\bigg|\, a_i \in Z_{11}I; 0 \leq i \leq \infty\right\} \cup$$

$$\left\{\begin{pmatrix} a & a & a & a \\ 0 & a & a & a \\ 0 & 0 & a & a \\ 0 & 0 & 0 & a \end{pmatrix} \,\bigg|\, a_i \in Z_{11}I\right\} \cup$$

$$\left\{\begin{pmatrix} a & 0 & 0 & 0 & 0 & 0 & 0 \\ a & a & 0 & 0 & 0 & 0 & 0 \\ a & a & a & 0 & 0 & 0 & 0 \\ a & a & a & a & 0 & 0 & 0 \\ a & a & a & a & a & 0 & 0 \\ a & a & a & a & a & a & 0 \\ a & a & a & a & a & a & a \end{pmatrix} \,\bigg|\, a \in Z_{11}I\right\}$$

$\subseteq V_1 \cup V_2 \cup V_3 \cup V_4 \cup V_5$.
W is a neutrosophic 5-linear subalgebra of V over the field $Z_{11}$. But V has no subneutrosophic 5-sublinear algebra.
    Hence the claim.

Now we define yet another new substructures.

**DEFINITION 3.1.8:** *Let $V = V_1 \cup V_2 \cup \ldots \cup V_n$ be a neutrosophic n-linear algebra over the real field F. Let $W = W_1 \cup W_2 \cup \ldots \cup W_n \subseteq V_1 \cup V_2 \cup \ldots \cup V_n$; be such that W is only a neutrosophic n-vector space over F and not a neutrosophic n-linear subalgebra of V; then we call W to be a neutrosophic pseudo n-vector subspace of V or W is a pseudo neutrosophic n-vector subspace of V.*

We will illustrate this situation by some simple examples.

**Example 3.1.27:** Let $V = V_1 \cup V_2 \cup V_3 \cup V_4 \cup V_5 =$

$$\left\{ \begin{pmatrix} a & b \\ c & d \end{pmatrix} \middle| a,b,c,d \in N(Q) \right\} \cup$$

$$\left\{ \begin{pmatrix} a & b & c \\ d & e & f \\ g & h & i \end{pmatrix} \middle| a,b,c,d,e,f,g,h,i \in QI \right\} \cup$$

$$\left\{ \sum_{i=0}^{\infty} a_i x^i \middle| a_i \in N(Q); 0 \le i \le \infty \right\} \cup$$

$\{N(Q)\} \cup \{5 \times 5$ neutrosophic upper triangular matrices with entries from $N(Q)\}$ be a neutrosophic 5-linear algebra over the real field Q.

Take $W = W_1 \cup W_2 \cup W_3 \cup W_4 \cup W_5 =$

$$\left\{ \begin{pmatrix} 0 & b \\ c & 0 \end{pmatrix} \middle| b,c \in N(Q) \right\} \cup \left\{ \begin{pmatrix} 0 & a & b \\ 0 & 0 & c \\ d & 0 & 0 \end{pmatrix} \middle| a,b,c,d \in QI \right\} \cup$$

$$\left\{ \sum_{i=0}^{8} a_i x^i \middle| a_i \in N(Q); 0 \le i \le 8 \right\} \cup \{QI\} \cup$$

$$\left\{ \begin{pmatrix} 0 & 0 & 0 & 0 & a \\ 0 & 0 & 0 & b & 0 \\ 0 & 0 & c & 0 & 0 \\ 0 & d & 0 & 0 & 0 \\ e & 0 & 0 & 0 & 0 \end{pmatrix} \middle| a,b,c,d,e \in N(Q) \right\}$$

$\subseteq V_1 \cup V_2 \cup V_3 \cup V_4 \cup V_5$; W is only a neutrosophic 5-vector space over the field Q; thus W is only a pseudo neutrosophic 5-vector subspace of V over Q.

***Example 3.1.28:*** Let $V_1 \cup V_2 \cup V_3 \cup V_4 \cup V_5 \cup V_6 =$

$$\left\{ \begin{pmatrix} a & b & c \\ d & e & f \\ g & h & i \end{pmatrix} \middle| a,b,c,d,e,f,g,h,i \in N(Z_2) \right\} \cup$$

$$\left\{ \sum_{i=0}^{\infty} a_i x^i \middle| a_i \in N(Z_2) \right\} \cup \left\{ \begin{pmatrix} a & b \\ a & b \end{pmatrix} \middle| a,b \in N(Z_2) \right\} \cup$$

{5 × 5 neutrosophic matrices with entries from $Z_2I$} ∪ {7 × 7 neutrosophic matrices with entries from $Z_2I$} ∪ {4 × 4 neutrosophic matrices with entries from $Z_2I$} be a neutrosophic 6-linear algebra over the real field $Z_2$.

Take $W = W_1 \cup W_2 \cup W_3 \cup W_4 \cup W_5 \cup W_6 =$

$$\left\{ \begin{pmatrix} 0 & 0 & a \\ 0 & b & 0 \\ c & 0 & 0 \end{pmatrix} \middle| a,b,c \in N(Z_2) \right\} \cup$$

$$\left\{ \sum_{i=0}^{29} a_i x^i \middle| a_i \in N(Z_2); 0 \le i \le 29 \right\} \cup \left\{ \begin{pmatrix} 0 & b \\ a & 0 \end{pmatrix} \middle| a,b \in N(Z_2) \right\} \cup$$

$$\left\{ \begin{pmatrix} 0 & 0 & 0 & 0 & a \\ 0 & 0 & 0 & b & 0 \\ 0 & 0 & c & 0 & 0 \\ 0 & d & 0 & 0 & 0 \\ e & 0 & 0 & 0 & 0 \end{pmatrix} \middle| a,b,c,d,e \in Z_2I \right\} \cup$$

$$\left\{ \begin{pmatrix} 0 & 0 & 0 & 0 & 0 & 0 & a \\ 0 & 0 & 0 & 0 & 0 & a & a \\ 0 & 0 & 0 & 0 & a & a & a \\ 0 & 0 & 0 & a & a & a & a \\ 0 & 0 & a & a & a & a & a \\ 0 & a & a & a & a & a & a \\ a & a & a & a & a & a & a \end{pmatrix} \middle| a \in Z_2 I \right\} \cup$$

$$\left\{ \begin{pmatrix} 0 & 0 & 0 & a_1 \\ 0 & 0 & a_2 & a_3 \\ 0 & a_4 & a_5 & a_6 \\ a_7 & a_8 & a_9 & a_{10} \end{pmatrix} \middle| a_i \in Z_2 I; 1 \leq i \leq 10 \right\}$$

$\subseteq V_1 \cup V_2 \cup V_3 \cup V_4 \cup V_5 \cup V_6$.

It is easily verified that W is only a neutrosophic 6-vector space over $Z_2$, so W is a pseudo neutrosophic 6-vector subspace of V over $Z_2$.

Now we proceed onto define pseudo subneutrosophic n-vector subspace of V.

**DEFINITION 3.1.9:** *Let $V = V_1 \cup V_2 \cup ... \cup V_n$ be a neutrosophic n-linear algebra over a real field F. Suppose $W = W_1 \cup W_2 \cup ... \cup W_n \subseteq V_1 \cup V_2 \cup ... \cup V_n$ be a neutrosophic n-vector space over a subfield K of F then we call W to be a pseudo subneutrosophic n-vector subspace of V over the subfield K of F.*

We will illustrate this by some simple examples.

*Example 3.1.29:* Let $V = V_1 \cup V_2 \cup V_3 \cup V_4 \cup V_5 =$

$$\left\{ \begin{pmatrix} a & b \\ c & d \end{pmatrix} \middle| a, b, c, d \in N(R) \right\} \cup \{N(R)\} \cup$$

$$\left\{ \begin{pmatrix} a_1 & a_7 & a_8 \\ a_2 & a_3 & a_9 \\ a_4 & a_5 & a_6 \end{pmatrix} \middle| a_i \in RI; 1 \le i \le 9 \right\} \cup$$

$$\left\{ \begin{pmatrix} a_1 & a_2 & a_3 & a_4 \\ a_{11} & a_5 & a_6 & a_7 \\ a_{12} & a_{13} & a_8 & a_9 \\ a_{14} & a_{15} & a_{16} & a_{10} \end{pmatrix} \middle| a_i \in N(R); 1 \le i \le 16 \right\} \cup$$

{N(R)[x]; all polynomials in the variable x with coefficients from N(R)} be a neutrosophic 5-linear algebra over the real field R. Take $W = W_1 \cup W_2 \cup W_3 \cup W_4 \cup W_5$

$$= \left\{ \begin{pmatrix} a & b \\ c & d \end{pmatrix} \middle| a,b,c,d \in QI \right\} \cup \{Q(R)\} \cup$$

$$\left\{ \begin{pmatrix} a_1 & 0 & 0 \\ a_2 & a_3 & 0 \\ a_4 & a_5 & a_6 \end{pmatrix} \middle| a_i \in QI; 1 \le i \le 6 \right\} \cup$$

$$\left\{ \begin{pmatrix} a_1 & a_2 & a_3 & a_4 \\ 0 & a_5 & a_6 & a_7 \\ 0 & 0 & a_8 & a_9 \\ 0 & 0 & 0 & a_{10} \end{pmatrix} \middle| a_i \in N(Q); 1 \le i \le 10 \right\} \cup$$

{N(Q)[x]; all polynomials in the variable x with coefficients from N(Q)} $\subseteq V_1 \cup V_2 \cup V_3 \cup V_4 \cup V_5$ is a subneutrosophic 5-linear subalgebra of V.

This will be different from pseudo subneutrosophic 5-vector subspace of V.

Take $W = W_1 \cup W_2 \cup W_3 \cup W_4 \cup W_5 =$

$$\left\{\begin{pmatrix} 0 & b \\ c & 0 \end{pmatrix} \middle| b,c \in QI\right\} \cup \{QI\} \cup \left\{\begin{pmatrix} 0 & 0 & a \\ 0 & b & 0 \\ c & 0 & 0 \end{pmatrix} \middle| a,b,c \in QI\right\} \cup$$

$$\left\{\begin{pmatrix} 0 & 0 & 0 & a \\ 0 & 0 & b & 0 \\ 0 & c & 0 & 0 \\ d & 0 & 0 & 0 \end{pmatrix} \middle| a,b,c,d \in QI\right\} \cup \left\{\sum_{i=0}^{41} a_i x^i \middle| a_i \in QI; 0 \le i \le 41\right\}$$

$\subseteq V_1 \cup V_2 \cup V_3 \cup V_4 \cup V_5$; W is only a neutrosophic 5-vector space over the field Q (Q a subfield R).

W is a pseudo subneutrosophic 5-vector subspace of V over the subfield Q of R.

*Example 3.1.30:* Let $V = V_1 \cup V_2 \cup V_3 \cup V_4 \cup V_5 \cup V_6 =$ {N(R)} $\cup$ {N(R)[x]; all neutrosophic polynomials in the variable x with coefficients from N (R)} $\cup$

$$\left\{\begin{pmatrix} a_1 & a_2 & a_3 \\ a_4 & a_5 & a_6 \\ a_7 & a_8 & a_9 \end{pmatrix} \middle| a_i \in N(R); 1 \le i \le 9\right\} \cup$$

$$\left\{\begin{pmatrix} a_1 & a_2 & a_3 & a_4 \\ a_5 & a_6 & a_7 & a_8 \\ a_9 & a_{10} & a_{11} & a_{12} \\ a_{13} & a_{14} & a_{15} & a_{16} \end{pmatrix} \middle| a_i \in RI; 1 \le i \le 16\right\} \cup$$

{All 8 × 8 matrices with entries from the neutrosophic field RI} $\cup$ {6 × 6 matrices with entries from the neutrosophic 6-linear algebra over the real field R} be a neutrosophic 6-linear algebra over the real field R. Take $W = W_1 \cup W_2 \cup W_3 \cup W_4 \cup W_5 \cup W_6 = \{QI\} \cup$

$$\left\{ \sum_{i=0}^{50} a_i x^i \,\middle|\, a_i \in QI; 0 \leq i \leq 50 \right\} \cup$$

$$\left\{ \begin{pmatrix} 0 & 0 & a \\ 0 & b & 0 \\ c & 0 & 0 \end{pmatrix} \,\middle|\, a,b,c \in QI \right\} \cup$$

$$\left\{ \begin{pmatrix} 0 & 0 & a & b \\ 0 & 0 & c & 0 \\ 0 & d & a & 0 \\ e & 0 & 0 & 0 \end{pmatrix} \,\middle|\, a,b,c,d,a,e \in QI \right\} \cup$$

$$\left\{ \begin{pmatrix} 0 & 0 & 0 & 0 & 0 & 0 & 0 & a \\ 0 & 0 & 0 & 0 & 0 & b & c & 0 \\ 0 & 0 & 0 & 0 & 0 & d & 0 & 0 \\ 0 & 0 & 0 & 0 & e & 0 & 0 & 0 \\ 0 & 0 & 0 & g & 0 & 0 & 0 & 0 \\ 0 & 0 & h & 0 & 0 & 0 & 0 & 0 \\ 0 & p & 0 & 0 & 0 & 0 & 0 & 0 \\ r & 0 & 0 & 0 & 0 & 0 & 0 & 0 \end{pmatrix} \,\middle|\, a,b,c,d,e,g,h,p,r \in QI \right\} \cup$$

$$\left\{ \begin{pmatrix} 0 & 0 & 0 & a & 0 & b \\ 0 & 0 & 0 & 0 & d & 0 \\ 0 & 0 & e & f & 0 & 0 \\ 0 & 0 & g & h & 0 & 0 \\ 0 & p & 0 & 0 & 0 & 0 \\ a & 0 & 0 & 0 & 0 & 0 \end{pmatrix} \,\middle|\, a,b,d,e,f,g,h,p,a \in QI \right\}$$

⊆ $V_1 \cup V_2 \cup V_3 \cup V_4 \cup V_5 \cup V_6$; W is a pseudo subneutrosophic 6-vector subspace of V over the real field Q.

If a neutrosophic n-linear algebra $V = V_1 \cup V_2 \cup \ldots \cup V_n$ does not contain any pseudo subneutrosophic n-vector subspace over a subfield K of F where V is defined over F; then we call V to be a pseudo simple subneutrosophic n-vector space over the field F.

We will illustrate this by some simple examples.

*Example 3.1.31:* Let $V = V_1 \cup V_2 \cup V_3 \cup V_4 \cup V_5 =$

$$\left\{ \begin{pmatrix} a & b \\ c & d \end{pmatrix} \middle| a,b,c,d \in Z_7I \right\} \cup \{N(Z_7)\} \cup$$

$$\left\{ \begin{pmatrix} a_1 & a_2 & a_3 & a_4 \\ a_5 & a_6 & a_7 & a_8 \\ a_9 & a_{10} & a_{11} & a_{12} \\ a_{13} & a_{14} & a_{15} & a_{16} \end{pmatrix} \middle| a_i \in Z_7I; 1 \leq i \leq 16 \right\} \cup$$

$$\left\{ \sum_{i=0}^{\infty} a_i x^i \middle| a_i \in Z_7I; 0 \leq i \leq \infty \right\} \cup$$

$$\left\{ \begin{pmatrix} a & b & c \\ d & e & f \\ g & h & i \end{pmatrix} \middle| a,b,c,d,e,f,g,h,i \in Z_7I \right\}$$

be a neutrosophic 5-linear algebra over the real field $Z_7$. Since $Z_7$ is a prime field it has no proper subfields. Hence V does not contain any pseudo subneutrosophic 5-vector subspace over $Z_7$. Hence V is a pseudo simple subneutrosophic 5-vector space over the field $Z_7$.

Now we proceed onto define linear transformation of neutrosophic n-vector space over the real field and discuss a few of its properties.

**DEFINITION 3.1.10:** *Let $V = V_1 \cup V_2 \cup ... \cup V_n$ be a neutrosophic n-vector space over a real field F and $W = W_1 \cup W_2 \cup ... \cup W_n$ be a neutrosophic n-vector space over the same real field F. Define $T : V \to W$. $T = T_1 \cup T_2 \cup ... \cup T_n : V = V_1 \cup V_2 \cup ... \cup V_n \to W = W_1 \cup W_2 \cup ... \cup W_n$ by $T(V_i) = W_j$ such that no two distinct $V_i$'s are mapped on to the same $W_j$; $1 \leq i, j \leq n$, where $T_i$ is a neutrosophic linear transformation from $V_i$ into $W_j$; $1 \leq i, j \leq n$, for i=1, 2, 3, ..., n. We call $T = T_1 \cup T_2 \cup ... \cup T_n$ to be a neutrosophic n-linear transformation of V into W.*

*If $W = V$ then we call T to be a neutrosophic n-linear operator on V. The set of all neutrosophic n-linear transformations of V into W, V and W defined over a real field F is denoted by*

*$N\,Hom_F(V, W) = \{$all neutrosophic n-linear transformations of V into W$\}$. $NHom_F(V, V) = \{$Collection of all neutrosophic n-linear operators of V into V$\}$.*

It is interesting and important to note that $V = V_1 \cup ... \cup V_n$ and $W = W_1 \cup ... \cup W_n$ are both defined over the same field F and both of them are only neutrosophic n-linear vector spaces.

We will illustrate by an example the neutrosophic n-linear transformation of V into W.

**Example 3.1.32:** Let $V = V_1 \cup V_2 \cup V_3 \cup V_4 \cup V_5 =$

$$\left\{ \begin{pmatrix} a & b \\ c & d \end{pmatrix} \middle| a,b,c,d \in Z_5 I \right\} \cup$$

$$\left\{ \begin{pmatrix} a & 0 & 0 \\ b & c & 0 \\ d & e & f \end{pmatrix} \middle| a,b,c,d,e,f \in (Z_5) \right\} \cup$$

$$\left\{ \begin{pmatrix} a & 0 & 0 & 0 \\ b & c & 0 & 0 \\ d & e & f & 0 \\ g & h & i & j \end{pmatrix} \middle| a,b,c,d,e,f,g,h,i,j \in Z_5I \right\} \cup$$

$$\left\{ \sum_{i=0}^{6} a_i x^i \middle| a_i \in Z_5I; 0 \le i \le 6 \right\} \cup$$

$$\left\{ \begin{pmatrix} a & 0 & 0 & 0 & 0 \\ 0 & b & 0 & 0 & 0 \\ 0 & 0 & c & 0 & 0 \\ 0 & 0 & 0 & d & 0 \\ 0 & 0 & 0 & 0 & e \end{pmatrix} \middle| a,b,c,d,e \in Z_5I \right\}$$

be a neutrosophic 5-vector space over the real field $Z_5$. $W = W_1 \cup W_2 \cup W_3 \cup W_4 \cup W_5 =$

$$\left\{ \begin{pmatrix} a & b & c \\ 0 & d & e \\ 0 & 0 & f \end{pmatrix} \middle| a,b,c,d,e,f \in Z_5I \right\} \cup$$

$$\{(a, b, c, d) \mid a, b, c, e \in Z_5I\} \cup$$

$$\left\{ \begin{pmatrix} a & b & c & d \\ 0 & 0 & 0 & 0 \\ e & f & g & h \\ 0 & 0 & 0 & 0 \end{pmatrix} \middle| a,b,c,d,e,f,g,h \in Z_5I \right\} \cup$$

$$\left\{ \sum_{i=0}^{9} a_i x^i \middle| a_i \in Z_5I; 0 \le i \le 9 \right\} \cup$$

$$\left\{ \begin{pmatrix} a & 0 & 0 & 0 & 0 & 0 & 0 \\ 0 & b & 0 & 0 & 0 & 0 & 0 \\ 0 & 0 & c & 0 & 0 & 0 & 0 \\ 0 & 0 & 0 & d & 0 & 0 & 0 \\ 0 & 0 & 0 & 0 & e & 0 & 0 \\ 0 & 0 & 0 & 0 & 0 & f & 0 \\ 0 & 0 & 0 & 0 & 0 & 0 & h \end{pmatrix} \middle| a,b,c,d,e,f,h \in Z_5 I \right\}$$

be a neutrosophic 5-vector space over the field $Z_5$.

Define $T : V \to W$ i.e., $T = T_1 \cup T_2 \cup T_3 \cup T_4 \cup T_5 : V = V_1 \cup V_2 \cup V_3 \cup V_4 \cup V_5 \to W = W_1 \cup W_2 \cup W_3 \cup W_4 \cup W_5$ where

$$T_1 : V_1 \to W_2,$$
$$T_2 : V_2 \to W_1,$$
$$T_3 : V_3 \to W_4,$$
$$T_4 : V_4 \to W_5$$

and

$$T_5 : V_5 \to W_3.$$

$$T_1 \begin{bmatrix} a & b \\ c & d \end{bmatrix} = (a, b, c, d);$$

$$T_2 \begin{bmatrix} a & 0 & 0 \\ b & c & 0 \\ d & e & f \end{bmatrix} = \begin{bmatrix} a & b & c \\ 0 & d & e \\ 0 & 0 & f \end{bmatrix},$$

$$T_3 \begin{bmatrix} a & 0 & 0 & 0 \\ b & c & 0 & 0 \\ d & e & f & 0 \\ g & h & i & j \end{bmatrix} = (a + bx + cx^2 + dx^3 + ex^4 + fx^5 + gx^6 + hx^7 + ix^8 + jx^9);$$

$$T_4\left[\sum_{i=0}^{6} a_i x^i\right] = \begin{bmatrix} a_0 & 0 & 0 & 0 & 0 & 0 & 0 \\ 0 & a_1 & 0 & 0 & 0 & 0 & 0 \\ 0 & 0 & a_2 & 0 & 0 & 0 & 0 \\ 0 & 0 & 0 & a_3 & 0 & 0 & 0 \\ 0 & 0 & 0 & 0 & a_4 & 0 & 0 \\ 0 & 0 & 0 & 0 & 0 & a_5 & 0 \\ 0 & 0 & 0 & 0 & 0 & 0 & a_6 \end{bmatrix}$$

and

$$T_5 \begin{bmatrix} a & 0 & 0 & 0 & 0 \\ 0 & b & 0 & 0 & 0 \\ 0 & 0 & c & 0 & 0 \\ 0 & 0 & 0 & d & 0 \\ 0 & 0 & 0 & 0 & e \end{bmatrix} = \begin{bmatrix} a & b & c & d \\ 0 & 0 & 0 & 0 \\ e & a & b & c \\ 0 & 0 & 0 & 0 \end{bmatrix}$$

It is easily verified that T is a neutrosophic 6-linear transformation of V into W.

**Example 3.1.33:** Let $V = V_1 \cup V_2 \cup V_3 \cup V_4 \cup V_5 \cup V_6 =$

$$\left\{ \begin{pmatrix} a & b \\ c & d \end{pmatrix} \middle| a, b, c, d \in N(Q) \right\} \cup$$

$$\{(a, b, c, d) \mid a, b, c, d \in N(Q)\} \cup$$

$$\left\{ \begin{pmatrix} a & 0 & 0 \\ b & c & 0 \\ d & e & f \end{pmatrix} \middle| a, b, c, d, e, f \in QI \right\} \cup$$

$$\left\{ \sum_{i=0}^{5} a_i x^i \middle| a_i \in QI; 0 \le i \le 5 \right\} \cup$$

$$\left\{ \begin{pmatrix} a & b & c & d \\ 0 & e & f & g \\ 0 & 0 & h & i \\ 0 & 0 & 0 & j \end{pmatrix} \middle| a,b,c,d,e,f,g,h,i,j \in QI \right\} \cup$$

$$\left\{ \begin{pmatrix} a & 0 & 0 & 0 \\ b & c & 0 & 0 \\ d & e & f & 0 \\ g & h & i & j \end{pmatrix} \middle| a,b,c,d,e,f,g,h,i,j \in QI \right\}$$

be a neutrosophic 6-vector space over the field Q. Define $T = T_1 \cup T_2 \cup T_3 \cup T_4 \cup T_5 \cup T_6 : V = V_1 \cup V_2 \cup V_3 \cup V_4 \cup V_5 \cup V_6 \to V = V_1 \cup V_2 \cup V_3 \cup V_4 \cup V_5 \cup V_6$ by

$$T_1 : V_1 \to V_2,$$
$$T_2 : V_2 \to V_1,$$
$$T_3 : V_3 \to V_4,$$
$$T_4 : V_4 \to V_3,$$
$$T_5 : V_5 \to V_6$$

and

$$T_6 : V_6 \to V_5$$

defined as follows:

$$T_1 \begin{pmatrix} a & b \\ c & d \end{pmatrix} = (a, b, c, d),$$

$$T_2 (a, b, c, d) = \begin{pmatrix} a & b \\ c & d \end{pmatrix},$$

$$T_3 \begin{pmatrix} a & 0 & 0 \\ b & c & 0 \\ d & e & f \end{pmatrix} = (a + bx + cx^2 + dx^3 + ex^4 + fx^5),$$

$$T_4\left(\sum_{i=0}^{5} a_i x^i\right) = \begin{pmatrix} a_0 & 0 & 0 \\ a_1 & a_2 & 0 \\ a_3 & a_4 & a_5 \end{pmatrix},$$

$$T_5 \begin{pmatrix} a & b & c & d \\ 0 & e & f & g \\ 0 & 0 & h & i \\ 0 & 0 & 0 & j \end{pmatrix} = \begin{bmatrix} a & 0 & 0 & 0 \\ b & c & 0 & 0 \\ d & e & f & 0 \\ g & h & i & j \end{bmatrix}$$

and

$$T_6 \begin{bmatrix} a & 0 & 0 & 0 \\ b & c & 0 & 0 \\ d & e & f & 0 \\ g & h & i & j \end{bmatrix} = \begin{pmatrix} a & b & c & d \\ 0 & e & f & g \\ 0 & 0 & h & i \\ 0 & 0 & 0 & j \end{pmatrix}.$$

It is easily verified that T is a neutrosophic 6-linear operator on V.

Now we proceed onto define other types of neutrosophic n-linear operators which will be know as the usual or common neutrosophic n-linear operators.

**DEFINITION 3.1.11:** *Let $V = V_1 \cup V_2 \cup ... \cup V_n$ be a neutrosophic n-vector space over a real field F.*

*Let $T: V \to V$ be a n-map such that $T = T_1 \cup T_2 \cup ... \cup T_n$ : $V = V_1 \cup V_2 \cup ... \cup V_n \to V_1 \cup V_2 \cup ... \cup V_n$ where $T_i : V_i \to V_i$; i=1, 2, ..., n if each $T_i$ is a linear operator then we define $T : V \to V$ to be a neutrosophic common n-linear operator on V or common neutrosophic n-linear operator on V.*

*We will denote the collection of all common neutrosophic n-linear operators on V by CN $Hom_F$ (V, V); clearly CN $Hom_F$ (V, V) is a neutrosophic $n^2$-subvector space of $NHom_F$ (V, V).*

*Further CN $Hom_F$ (V, V) = $Hom_F$ ($V_1$, $V_1$) $\cup$ $Hom_F$ ($V_2$, $V_2$) $\cup ... \cup Hom_F$ ($V_n$, $V_n$).*

We will illustrate this situation by an example.

**Example 3.1.34:** Let $V = V_1 \cup V_2 \cup V_3 \cup V_4 \cup V_5 \cup V_6 \cup V_7$

$$= \left\{ \begin{pmatrix} a & b \\ c & d \end{pmatrix} \middle| a,b,c,d \in Z_{11}I \right\} \cup$$

$$\{(a, b, c, d) \mid a, b, c, d \in Z_{11}I\} \cup$$

$$\left\{ \sum_{i=0}^{12} a_i x^i \middle| a_i \in Z_{11}I; 0 \leq i \leq 12 \right\} \cup$$

$$\left\{ \begin{pmatrix} a & b & c \\ d & e & f \\ g & h & i \end{pmatrix} \middle| a,b,c,d,e,f,g,h,i \in Z_{11}I \right\} \cup$$

$$\left\{ \begin{pmatrix} a & 0 & 0 & 0 & 0 \\ 0 & b & 0 & 0 & 0 \\ 0 & 0 & c & 0 & 0 \\ 0 & 0 & 0 & d & 0 \\ 0 & 0 & 0 & 0 & e \end{pmatrix} \middle| a,b,c,d,e \in Z_{11}I \right\} \cup$$

$$\left\{ \begin{pmatrix} 0 & 0 & a & b & d \\ 0 & 0 & 0 & e & f \\ 0 & 0 & 0 & 0 & g \\ a & 0 & 0 & 0 & 0 \\ b & e & 0 & 0 & 0 \\ d & f & g & 0 & 0 \end{pmatrix} \middle| a,b,e,f,g,d \in Z_{11}I \right\} \cup$$

$$\left\{ \begin{pmatrix} a_1 & a_2 & a_3 & a_4 \\ a_5 & a_6 & a_7 & a_8 \\ a_9 & a_{10} & a_{11} & a_{12} \end{pmatrix} \middle| a_i \in Z_{11}I; 1 \leq i \leq 12 \right\}$$

be a neutrosophic 7-vector space over the field $Z_{11}$.

Define $T : V_1 \cup V_2 \cup V_3 \cup V_4 \cup V_5 \cup V_6 \cup V_7 \to V_1 \cup V_2 \cup V_3 \cup V_4 \cup V_5 \cup V_6 \cup V_7$ where $T : T_1 \cup T_2 \cup \ldots \cup T_7$ such that $T_i : V_i \to V_i$; $i = 1, 2, \ldots, 7$.

$T_1 : V_1 \to V_1$ is a neutrosophic linear operator on $V_i$ defined by

$$T_1 \begin{pmatrix} a & b \\ c & d \end{pmatrix} = \left\{ \begin{pmatrix} a & b \\ 0 & d \end{pmatrix} \right\},$$

$T_2$ is a neutrosophic linear operator on $V_2$ defined by $T_2: V_2 \to V_2$ and $T_2(a, b, c, d) = (a, b, a, b)$.

$T_3$ is a neutrosophic linear operator on $V_3$, $T_3: V_3 \to V_3$ is given by

$$T_3 \left( \sum_{i=0}^{12} a_i x^i \right) = (a_0 + a_2 x^2 + a_4 x^4 + a_6 x^6 + a_8 x^8 + a_{10} x^{10} + a_{12} x^{12}).$$

$T_4 : V_4 \to V_4$ is a neutrosophic linear operator given by

$$T_4 \begin{pmatrix} a & b & c \\ d & e & f \\ g & h & i \end{pmatrix} = \left\{ \begin{pmatrix} a & b & c \\ 0 & e & f \\ 0 & 0 & i \end{pmatrix} \right\}$$

and $T_5 : V_5 \to V_5$ is a neutrosophic linear operator given by

$$T_5 \begin{pmatrix} a & 0 & 0 & 0 & 0 \\ 0 & b & 0 & 0 & 0 \\ 0 & 0 & c & 0 & 0 \\ 0 & 0 & 0 & d & 0 \\ 0 & 0 & 0 & 0 & e \end{pmatrix} = \begin{bmatrix} a & 0 & 0 & 0 & 0 \\ 0 & a+b & 0 & 0 & 0 \\ 0 & 0 & b+c & 0 & 0 \\ 0 & 0 & 0 & c+d & 0 \\ 0 & 0 & 0 & 0 & d+a \end{bmatrix}.$$

$T_6 : V_6 \to V_6$ is a neutrosophic linear operator on $V_6$ given by

$$T_6 \begin{pmatrix} 0 & 0 & a & b & d \\ 0 & 0 & 0 & e & f \\ 0 & 0 & 0 & 0 & g \\ a & 0 & 0 & 0 & 0 \\ b & e & 0 & 0 & 0 \\ d & f & g & 0 & 0 \end{pmatrix} = \begin{bmatrix} 0 & 0 & a+b & b+d & d+e \\ 0 & 0 & 0 & f+e & f+g \\ 0 & 0 & 0 & 0 & g \\ a+b & 0 & 0 & 0 & 0 \\ b+d & f+e & 0 & 0 & 0 \\ d+e & f+g & g & 0 & 0 \end{bmatrix}$$

and

$T_7 : V_7 \to V_7$ is a neutrosophic linear operator on $V_7$ given by

$$T_7 \begin{pmatrix} a_1 & a_2 & a_3 & a_4 \\ a_5 & a_6 & a_7 & a_8 \\ a_9 & a_{10} & a_{11} & a_{12} \end{pmatrix} = \begin{pmatrix} a_4 & a_3 & a_2 & a_1 \\ a_6 & a_5 & a_8 & a_7 \\ a_{12} & a_{11} & a_{10} & a_9 \end{pmatrix}.$$

It is easily verified that $T = T_1 \cup T_2 \cup T_3 \cup T_4 \cup T_5 \cup T_6 \cup T_7$ is a neutrosophic 7-linear operator on V.

Now we define two types of neutrosophic (m, n) linear transformation of a neutrosophic m-vector space into a neutrosophic n-vector space $m \neq n$ and $m > n$.

**DEFINITION 3.1.12:** *Let $V = V_1 \cup V_2 \cup ... \cup V_m$ be a neutrosophic m-vector space over the real field F and $W = W_1 \cup W_2 \cup ... \cup W_n$ be a neutrosophic n-vector space over the same field F; ($m \neq n$) and $m > n$).*

*Let $T = T_1 \cup T_2 \cup ... \cup T_m$ be a m-map from V into W such that $T_1 \cup T_2 \cup ... \cup T_m : V_1 \cup V_2 \cup ... \cup V_m \to W_1 \cup W_2 \cup ... \cup W_n$ given by $T_i : V_i \to W_j$, it is sure to happen that more than one $V_i$ is mapped onto a $W_j$, such that each $T_i$ is a neutrosophic linear transformation from $V_i$ to $W_j$; $1 \leq i \leq m$ and $1 \leq j \leq n$.*

*Then $T = T_1 \cup T_2 \cup ... \cup T_m$ is defined as a special (m, n) neutrosophic linear transformation of V to W.*

We will first illustrate this situation by an example.

*Example 3.1.35:* Let $V = V_1 \cup V_2 \cup V_3 \cup V_4 \cup V_5 \cup V_6 =$

$$\left\{ \begin{pmatrix} a & b \\ c & d \end{pmatrix} \middle| a, b, c, d \in QI \right\} \cup$$

$$\{(a, b, c, d, e, f) \mid a, b, c, d, e, f \in QI\} \cup$$

$$\left\{ \begin{pmatrix} a_1 & a_2 \\ a_3 & a_4 \\ a_5 & a_6 \end{pmatrix} \middle| a_i \in QI; 1 \le i \le 6 \right\} \cup \left\{ \begin{pmatrix} a & b \\ c & d \end{pmatrix} \middle| a, b, c, d \in QI \right\} \cup$$

$$\left\{ \begin{pmatrix} a_1 & a_2 & a_3 & a_4 \\ a_5 & a_6 & a_7 & a_8 \end{pmatrix} \middle| a_i \in QI; 1 \le i \le 8 \right\} \cup$$

$$\left\{ \begin{pmatrix} a_1 & a_2 \\ a_3 & a_4 \\ a_5 & a_6 \\ a_7 & a_8 \end{pmatrix} \middle| a_i \in QI; 1 \le i \le 8 \right\}$$

be a neutrosophic 6-vector space over the field Q.

Let $W = W_1 \cup W_2 \cup W_3 \cup W_4 =$

$$\left\{ \begin{pmatrix} a & b \\ c & d \\ e & f \end{pmatrix} \middle| a, b, c, d, e, f \in QI \right\} \cup$$

$$\left\{ \sum_{i=0}^{7} a_i x^i \middle| a_i \in QI; 0 \le i \le 7 \right\} \cup \{(a, b, c, d) \mid a, b, c, d \in QI\} \cup$$

$$\left\{ \begin{pmatrix} a_1 & a_2 & a_3 \\ a_4 & a_5 & a_6 \end{pmatrix} \middle| a_i \in QI; 1 \le i \le 6 \right\}$$

be a neutrosophic 4-vector space over the real field Q.

Define $T = T_1 \cup T_2 \cup T_3 \cup T_4 \cup T_5 \cup T_6 : V = V_1 \cup V_2 \cup V_3 \cup V_4 \cup V_5 \cup V_6 \to W = W_1 \cup W_2 \cup W_3 \cup W_4$ by

$$T_1 : V_1 \to W_3$$
$$T_2 : V_2 \to W_4$$
$$T_3 : V_3 \to W_1$$
$$T_4 : V_4 \to W_1,$$
$$T_5 : V_5 \to W_2$$

and $\qquad T_6 : V_6 \to W_3$

where each $T_i$ is a neutrosophic linear transformation from $V_i$ to $W_j$; $i = 1, 2, \ldots, 6$ and $j = 1, 2, 3, 4$.

$T_1 : V_1 \to W_3$ is defined by

$$T_1 \begin{pmatrix} a & b \\ c & d \end{pmatrix} = (a, b, c, d),$$

$T_1$ is a neutrosophic linear transformation from $V_1$ to $W_3$.

$T_2 : V_2 \to W_4$ is such that

$$T_2 (a, b, c, d, e, f) = \begin{pmatrix} a & b & c \\ d & e & f \end{pmatrix}.$$

Clearly $T_2$ is neutrosophic linear transformation from $V_2$ to $W_4$.

$T_3 : V_3 \to W_1$ is given by

$$T_3 \begin{pmatrix} a_1 & a_2 \\ a_3 & a_4 \\ a_5 & a_6 \end{pmatrix} = \begin{pmatrix} a_1 & a_3 \\ a_2 & a_5 \\ a_6 & a_4 \end{pmatrix}$$

$T_3$ is a neutrosophic linear transformation from $V_3$ to $W_1$.

$T_4 : V_4 \to W_1$ is defined by

$$T_4 \begin{pmatrix} a & 0 & 0 \\ b & c & 0 \\ d & e & f \end{pmatrix} = \begin{bmatrix} a+b & d \\ c & e+f \end{bmatrix}$$

$T_4$ is a neutrosophic linear transformation of $V_4$ to $W_1$.
$T_5 : V_5 \to W_2$ is such that

$$T_5 \begin{pmatrix} a_1 & a_2 & a_3 & a_4 \\ a_5 & a_6 & a_7 & a_8 \end{pmatrix} = a_1 + a_2x + a_3x^2 + a_4x^3 + a_5x^4 + a_6x^5 + a_7x^6 + a_8x^7$$

$T_5$ is a neutrosophic linear transformation of $V_5$ to $W_2$.

$T_6 : V_6 \to W_3$ is defined by

$$T_6 \begin{pmatrix} a_1 & a_2 \\ a_3 & a_4 \\ a_5 & a_6 \\ a_7 & a_8 \end{pmatrix} = (a_1 + a_2, a_3 + a_4, a_5 + a_6, a_7 + a_8).$$

Clearly $T_6$ is a neutrosophic linear transformation of $V_6$ to $W_3$. Thus $T = (T_1 \cup T_2 \cup T_3 \cup T_4 \cup T_5 \cup T_6)$ is a (6, 4) neutrosophic linear transformation of V to W.

**DEFINITION 3.1.13:** *Let $V = V_1 \cup V_2 \cup V_3 \cup ... \cup V_m$ be a neutrosophic m-vector space over a field F and $W = W_1 \cup W_2 \cup ... \cup W_n$ be a neutrosophic n-vector space defined over the same field F ($m \neq n$, $m < n$). Define a m-map, $T = T_1 \cup T_2 \cup ... \cup T_m$: $V = V_1 \cup V_2 \cup ... \cup V_m$ into $W_1 \cup W_2 \cup ... \cup W_n$ such that $T_i : V_i \to W_j$ where each $V_i$ is mapped into a distinct $W_j$; $1 \leq i \leq m$ and $1 \leq j \leq n$; where each $T_i$ is a neutrosophic linear transformation of $V_i$ to $W_j$.*

*We define $T = T_1 \cup T_2 \cup ... \cup T_m$ as a (m, n) neutrosophic linear transformation of V into W.*

We will illustrate this situation by an example.

***Example 3.1.36:*** Let $V = V_1 \cup V_2 \cup V_3 \cup V_4 \cup V_5 =$

$$\left\{ \begin{pmatrix} a_1 & a_2 & a_3 \\ a_4 & a_5 & a_6 \end{pmatrix} \middle| a_i \in Z_{17}I; 1 \leq i \leq 6 \right\} \cup$$

$$\left\{ \begin{pmatrix} a_1 & a_2 \\ a_3 & a_4 \\ a_5 & a_6 \\ a_7 & a_8 \\ a_9 & a_{10} \end{pmatrix} \middle| a_i \in Z_{17}I; 1 \le i \le 10 \right\} \cup$$

$$\left\{ \sum_{i=0}^{11} a_i x^i \middle| a_i \in Z_{17}I; 0 \le i \le 11 \right\} \cup \left\{ \begin{pmatrix} a & b \\ c & d \end{pmatrix} \middle| a,b,c,d \in Z_{17}I \right\} \cup$$

$$\left\{ \begin{pmatrix} a_1 & 0 & 0 & 0 \\ a_2 & a_3 & 0 & 0 \\ a_4 & a_5 & a_6 & 0 \\ a_7 & a_8 & a_9 & a_{10} \end{pmatrix} \middle| a_i \in Z_{17}I; 1 \le i \le 10 \right\}$$

be a neutrosophic 5-vector space over the field $Z_{17}$. Let $W = W_1 \cup W_2 \cup W_3 \cup W_4 \cup W_5 \cup W_6 \cup W_7 =$

$$\left\{ \begin{pmatrix} a & 0 & 0 \\ b & c & 0 \\ d & e & f \end{pmatrix} \middle| a,b,c,d,e,f \in Z_{17}I \right\} \cup$$

$$\left\{ \begin{pmatrix} a_1 & a_2 & a_3 & a_4 & a_5 & a_6 \\ a_7 & a_8 & a_9 & a_{10} & a_{11} & a_{12} \end{pmatrix} \middle| a_i \in Z_{17}I; 1 \le i \le 12 \right\} \cup$$

$$\left\{ \begin{pmatrix} a & d & e & b \\ g & h & i & j \\ k & l & m & n \\ s & p & q & r \end{pmatrix} \middle| a,b,d,e,g,h,i,j,k,l,m,n,s,p,q,r \in Z_{17}I \right\} \cup$$

$$\{(a, b, c, d, e) \mid a, b, c, d, e \in Z_{17}I\} \cup$$

$$\left\{ \begin{pmatrix} a_1 & 0 & 0 & 0 & 0 \\ 0 & a_2 & 0 & 0 & 0 \\ 0 & 0 & a_3 & 0 & 0 \\ 0 & 0 & 0 & a_4 & 0 \\ 0 & 0 & 0 & 0 & a_5 \end{pmatrix} \middle| a_i \in Z_{17}I; 1 \leq i \leq 5 \right\} \cup$$

$$\left\{ \sum_{i=0}^{9} a_i x^i \middle| a_i \in Z_{17}I; 0 \leq i \leq 9 \right\} \cup$$

$$\left\{ \begin{pmatrix} a_1 & a_2 & a_3 \\ a_4 & a_5 & a_6 \\ a_7 & a_8 & a_9 \\ a_{10} & a_{11} & a_{12} \end{pmatrix} \middle| a_i \in Z_{17}I; 1 \leq i \leq 12 \right\}$$

be a neutrosophic 7-vector space over the field $Z_{17}$. Define $T = T_1 \cup T_2 \cup T_3 \cup T_4 \cup T_5 : V = V_1 \cup V_2 \cup V_3 \cup V_4 \cup V_5 \to W = W_1 \cup W_2 \cup W_3 \cup W_4 \cup W_5 \cup W_6 \cup W_7$ such that

$$T_1 : V_1 \to W_1,$$
$$T_2 : V_2 \to W_6,$$
$$T_3 : V_3 \to W_2,$$
$$T_4 : V_4 \to W_3,$$

and
$$T_5 : V_5 \to W_5.$$

where

$T_1 : V_1 \to W_1$ is a neutrosophic linear transformation given by

$$T_1 \begin{pmatrix} a_1 & a_2 & a_3 \\ a_4 & a_5 & a_6 \end{pmatrix} = \left\{ \begin{pmatrix} a_1 & 0 & 0 \\ a_2 & a_3 & 0 \\ a_4 & a_5 & a_6 \end{pmatrix} \right\}$$

and

$T_2 : V_2 \to W_6$ is defined by

$$T_2 \begin{pmatrix} a_1 & a_2 \\ a_3 & a_4 \\ a_5 & a_6 \\ a_7 & a_8 \\ a_9 & a_{10} \end{pmatrix} = (a_1 + a_2 x + a_3 x^2 + a_4 x^3 + a_5 x^4 + a_6 x^5 + a_7 x^6 +$$

$$a_8 x^7 + a_9 x^8 + a_{10} x^9)$$

is a neutrosophic linear transformation of $V_2$ to $W_6$.

$T_3 : V_3 \to W_2$ is given by

$$T_3 \left( \sum_{i=0}^{11} a_i x^i \right) = \begin{pmatrix} a_1 & a_2 & a_3 & a_4 & a_5 & a_6 \\ a_7 & a_8 & a_9 & a_{10} & a_{11} & a_0 \end{pmatrix}.$$

$T_3$ is again a neutrosophic linear transformation from $V_3$ to $W_2$.

Consider $T_4 : V_4 \to W_3$ given by

$$T_4 \begin{pmatrix} a & b \\ c & d \end{pmatrix} = \begin{pmatrix} a & 0 & 0 & b \\ 0 & 0 & 0 & 0 \\ 0 & 0 & 0 & 0 \\ c & 0 & 0 & d \end{pmatrix}$$

$T_4$ is also a neutrosophic linear transformation from $V_4$ to $W_3$.

$T_5 : V_5 \to W_5$ defined by

$$T_5 \begin{pmatrix} a_1 & 0 & 0 & 0 \\ a_2 & a_3 & 0 & 0 \\ a_4 & a_5 & a_6 & 0 \\ a_7 & a_8 & a_9 & a_{10} \end{pmatrix} = \begin{pmatrix} a_1 + a_2 & 0 & 0 & 0 & 0 \\ 0 & a_3 + a_4 & 0 & 0 & 0 \\ 0 & 0 & a_5 + a_6 & 0 & 0 \\ 0 & 0 & 0 & a_7 + a_8 & 0 \\ 0 & 0 & 0 & 0 & a_9 + a_{10} \end{pmatrix}$$

is a neutrosophic linear transformation. Thus $T = T_1 \cup T_2 \cup T_3 \cup T_4 \cup T_5$ is a (5, 7) neutrosophic linear transformation of V to W.

Now having defined several types of neutrosophic n-linear transformations of neutrosophic n-vector spaces V and W we can define in a similar way all types of neutrosophic n-linear transformation for neutrosophic n-linear algebras with appropriate changes.
    We will only illustrate them by examples as modified definitions can be easily obtained by any reader.

***Example 3.1.37:*** Let $V = V_1 \cup V_2 \cup V_3 \cup V_4 =$

$$\left\{ \begin{pmatrix} a & b \\ c & d \end{pmatrix} \middle| a,b,c,d \in Z_{13}I \right\} \cup$$

$$\{(a_1, a_2, a_3, a_4, a_5, a_6) \mid a_i \in Z_{13}I; 1 \le i \le 6\} \cup$$

$$\left\{ \begin{pmatrix} a_1 & 0 & 0 \\ a_2 & a_3 & 0 \\ a_4 & a_5 & a_6 \end{pmatrix} \middle| a_i \in Z_{13}I; 1 \le i \le 6 \right\} \cup$$

$$\left\{ \begin{pmatrix} a_1 & a_2 & a_3 & a_4 \\ a_5 & a_6 & a_7 & a_8 \\ a_9 & a_{10} & a_{11} & a_{12} \\ a_{13} & a_{14} & a_{15} & a_{16} \end{pmatrix} \middle| a_i \in Z_{13}I; 1 \le i \le 16 \right\}$$

be a neutrosophic 4-linear algebra over the field $Z_{13}$. Let $W = W_1 \cup W_2 \cup W_3 \cup W_4 =$

$$\left\{ \begin{pmatrix} a_1 & a_2 & a_3 \\ 0 & a_4 & a_5 \\ 0 & 0 & a_6 \end{pmatrix} \middle| a_i \in Z_{13}I; 1 \le i \le 6 \right\} \cup$$

$$\left\{ \begin{pmatrix} a & 0 & 0 & 0 \\ 0 & b & 0 & 0 \\ 0 & 0 & c & 0 \\ 0 & 0 & 0 & d \end{pmatrix} \middle| a,b,c,d \in Z_{13}I \right\} \cup$$

$$\left\{ \begin{pmatrix} a & b \\ c & d \end{pmatrix} \middle| a,b,c,d \in Z_{13}I \right\} \cup$$

$$\left\{ \begin{pmatrix} a_1 & a_2 & a_3 & a_4 & a_5 \\ a_6 & a_7 & a_8 & a_9 & a_{10} \\ a_{11} & a_{12} & a_{13} & a_{14} & a_{15} \\ a_{16} & a_{17} & a_{18} & a_{19} & a_{20} \\ a_{21} & a_{22} & a_{23} & a_{24} & a_{25} \end{pmatrix} \middle| a_i \in Z_{13}I; 1 \le i \le 25 \right\}$$

be a neutrosophic 4-linear algebra over the field $Z_{13}$. Define a 4-map $T = T_1 \cup T_2 \cup T_3 \cup T_4 : V = V_1 \cup V_2 \cup V_3 \cup V_4 \to W_1 \cup W_2 \cup W_3 \cup W_4$ as follows.

$$T_1 : V_1 \to W_2,$$
$$T_2 : V_2 \to W_1,$$
$$T_3 : V_3 \to W_4$$

and
$$T_4 : V_4 \to W_3$$

so that each $T_i$ is a neutrosophic linear transformation.

$T_1 : V_1 \to W_2$ is such that

$$T_1 \begin{pmatrix} a & b \\ c & d \end{pmatrix} = \begin{pmatrix} a & 0 & 0 & 0 \\ 0 & b & 0 & 0 \\ 0 & 0 & c & 0 \\ 0 & 0 & 0 & d \end{pmatrix}$$

is a neutrosophic linear transformation of neutrosophic linear algebras $V_1$ into $W_2$.

$T_2 : V_2 \to W_1$ is such that

$$T_2(a_1, a_2, a_3, a_4, a_5, a_6) = \begin{pmatrix} a_1 & a_2 & a_3 \\ 0 & a_4 & a_5 \\ 0 & 0 & a_6 \end{pmatrix}$$

is a neutrosophic linear transformation of $V_2$ to $W_1$.

$T_3 : V_3 \to W_4$ is defined by

$$T_3 \begin{pmatrix} a_1 & 0 & 0 \\ a_2 & a_3 & 0 \\ a_4 & a_5 & a_6 \end{pmatrix} = \begin{pmatrix} a_1 & 0 & 0 & 0 & 0 \\ 0 & a_2 & 0 & 0 & 0 \\ 0 & 0 & a_3 & 0 & 0 \\ 0 & 0 & 0 & a_4 & 0 \\ 0 & 0 & 0 & 0 & a_5 \end{pmatrix}$$

is a neutrosophic linear transformation of $V_3$ into $W_4$.

$T_4 : V_4 \to W_3$  given by

$$T_4 \begin{bmatrix} a_1 & a_2 & a_3 & a_4 \\ a_5 & a_6 & a_7 & a_8 \\ a_9 & a_{10} & a_{11} & a_{12} \\ a_{13} & a_{14} & a_{15} & a_{16} \end{bmatrix} = \begin{bmatrix} a_1 + a_2 + a_3 + a_4 & a_5 + a_6 + a_7 + a_8 \\ a_9 + a_{10} + a_{11} + a_{12} & a_{13} + a_{14} + a_{15} + a_{16} \end{bmatrix}$$

is a neutrosophic linear transformation of $V_4$ into $W_3$.

Thus $T = T_1 \cup T_2 \cup T_3 \cup T_4$ is a neutrosophic 4-linear transformation the neutrosophic 4-linear algebra $V = V_1 \cup V_2 \cup V_3 \cup V_4$ into the neutrosophic 4-linear algebra $W = W_1 \cup W_2 \cup W_3 \cup W_4$.

We will now give an example of a neutrosophic n-linear operator.

*Example 3.1.38:* Let $V = V_1 \cup V_2 \cup V_3 \cup V_4 \cup V_5$ be a neutrosophic 5-linear algebra over the field $Z_{17}$ where

$$V_1 = \left\{ \begin{pmatrix} a & b \\ c & d \end{pmatrix} \middle| a,b,c,d \in Z_{17}I \right\},$$

$$V_2 = \{(a_1, a_2, a_3, a_4, a_5, a_6) \mid a_i \in Z_{17}I, 1 \le i \le 6\},$$

$$V_3 = \left\{ \begin{pmatrix} a_1 & 0 & 0 \\ a_2 & a_3 & 0 \\ a_4 & a_5 & a_6 \end{pmatrix} \middle| a_i \in Z_{17}I; 1 \le i \le 6 \right\},$$

$$V_4 = \left\{ \begin{pmatrix} a_1 & 0 & 0 & 0 \\ 0 & a_2 & 0 & 0 \\ 0 & 0 & a_3 & 0 \\ 0 & 0 & 0 & a_4 \end{pmatrix} \middle| a_i \in Z_{17}I; 1 \le i \le 4 \right\}$$

and
$V_5 = \{$all $6 \times 6$ neutrosophic matrices with entries from $Z_{17}I\}$.

Define $T = T_1 \cup T_2 \cup T_3 \cup T_4 \cup T_5 : V = V_1 \cup V_2 \cup V_3 \cup V_4 \cup V_5 \to V = V_1 \cup V_2 \cup V_3 \cup V_4 \cup V_5$; where

$$T_1 : V_1 \to V_4;$$
$$T_2 : V_2 \to V_3;$$
$$T_3 : V_3 \to V_2;$$
$$T_4 : V_4 \to V_1$$

and
$$T_5 : V_5 \to V_2$$

such that

$$T_1 \begin{pmatrix} a & b \\ c & d \end{pmatrix} = \begin{pmatrix} a_1 & 0 & 0 & 0 \\ 0 & a_2 & 0 & 0 \\ 0 & 0 & a_3 & 0 \\ 0 & 0 & 0 & a_4 \end{pmatrix};$$

$$T_2 (a_1, a_2, a_3, a_4, a_5, a_6) = \begin{pmatrix} a_1 & 0 & 0 \\ a_2 & a_3 & 0 \\ a_4 & a_5 & a_6 \end{pmatrix};$$

$$T_3 \begin{pmatrix} a_1 & 0 & 0 \\ a_2 & a_3 & 0 \\ a_4 & a_5 & a_6 \end{pmatrix} = (a_1, a_2, a_3, a_4, a_5, a_6);$$

$$T_4 \begin{pmatrix} a_1 & 0 & 0 & 0 \\ 0 & a_2 & 0 & 0 \\ 0 & 0 & a_3 & 0 \\ 0 & 0 & 0 & a_4 \end{pmatrix} = \begin{pmatrix} a_1 & a_2 \\ a_3 & a_4 \end{pmatrix}$$

and

$$T_5 \begin{pmatrix} a_1 & a_2 & a_3 & a_4 & a_5 & a_6 \\ a_7 & a_8 & a_9 & a_{10} & a_{11} & a_{12} \\ a_{13} & a_{14} & a_{15} & a_{16} & a_{17} & a_{18} \\ a_{19} & a_{20} & a_{21} & a_{22} & a_{23} & a_{24} \\ a_{25} & a_{26} & a_{27} & a_{28} & a_{29} & a_{30} \\ a_{31} & a_{32} & a_{33} & a_{34} & a_{35} & a_{36} \end{pmatrix} = (a_1, a_2, a_3, a_4, a_5, a_6).$$

It is easily verified $T = T_1 \cup T_2 \cup T_3 \cup T_4 \cup T_5$ is a neutrosophic 5-linear operator on V. Clearly T is not a usual neutrosophic 5-linear operator on V.

We will now illustrate by the example the usual neutrosophic n-linear operator on V.

*Example 3.1.39:* Let $V = V_1 \cup V_2 \cup V_3 \cup V_4 \cup V_5 \cup V_6 =$

$$\left\{ \begin{pmatrix} a & b \\ c & d \end{pmatrix} \middle| a,b,c,d \in Z_{19}I \right\} \cup$$

$$\{(a_1, a_2, a_3, a_4, a_5, a_6, a_7) \mid a_i \in Z_{19}I, 1 \leq i \leq 7\} \cup$$

$$\left\{ \begin{pmatrix} a_1 & 0 & 0 & 0 \\ a_2 & a_3 & 0 & 0 \\ a_4 & a_5 & a_6 & 0 \\ a_7 & a_8 & a_9 & a_{10} \end{pmatrix} \middle| a \in Z_{19}I; 1 \leq i \leq 10 \right\} \cup$$

{All 5×5 neutrosophic matrix with entries from $Z_{19}I$} ∪ {8 × 8 neutrosophic diagonal matrices with entries from $Z_{19}I$} ∪

$$\left\{ \begin{pmatrix} a_1 & a_2 & a_3 & a_4 \\ 0 & a_5 & a_6 & a_7 \\ 0 & 0 & a_8 & a_9 \\ 0 & 0 & 0 & a_{10} \end{pmatrix} \middle| a_i \in Z_{19}I; 1 \leq i \leq 10 \right\}$$

be a neutrosophic 6-linear algebra over the real field $Z_{19}$.

Define $T = T_1 \cup T_2 \cup T_3 \cup T_4 \cup T_5 \cup T_6 : V = V_1 \cup V_2 \cup V_3 \cup V_4 \cup V_5 \cup V_6 \to V_1 \cup V_2 \cup V_3 \cup V_4 \cup V_5 \cup V_6 = V$ such that $T_i : V_i \to V_i$ for $i = 1, 2, \ldots, 6$.

$T_1 : V_1 \to V_1$ such that

$$T_1 \begin{pmatrix} a & b \\ c & d \end{pmatrix} = \begin{pmatrix} b & a \\ d & c \end{pmatrix};$$

$T_2 : V_2 \to V_2$ is defined by

$T_2 (a_1, a_2, a_3, a_4, a_5, a_6) = (a_1 + a_2, a_2 + a_3, a_3 + a_4, a_4 + a_5, a_5 + a_6, a_6 + a_1)$,

$$T_3 \begin{pmatrix} a_1 & a_2 & a_3 & a_4 \\ 0 & a_5 & a_6 & a_7 \\ 0 & 0 & a_8 & a_9 \\ 0 & 0 & 0 & a_{10} \end{pmatrix} \to \begin{pmatrix} a_1 & 0 & 0 & 0 \\ 0 & a_3 & 0 & 0 \\ 0 & 0 & a_6 & 0 \\ 0 & 0 & 0 & a_{10} \end{pmatrix}$$

where $T_3 : V_3 \to V_3$ ;

$T_4 : V_4 \to V_4$ is such that

$$T_4 \begin{pmatrix} a_1 & a_2 & a_3 & a_4 & a_5 \\ a_6 & a_7 & a_8 & a_9 & a_{10} \\ a_{11} & a_{12} & a_{13} & a_{14} & a_{15} \\ a_{16} & a_{17} & a_{18} & a_{19} & a_{20} \\ a_{21} & a_{22} & a_{23} & a_{24} & a_{25} \end{pmatrix} \rightarrow \begin{pmatrix} a_1 & 0 & 0 & 0 & 0 \\ a_2 & a_3 & 0 & 0 & 0 \\ a_4 & a_5 & a_6 & 0 & 0 \\ a_7 & a_8 & a_9 & a_{10} & 0 \\ a_{11} & a_{12} & a_{13} & a_{14} & a_{15} \end{pmatrix}$$

$T_4$ is a neutrosophic linear operator on $V_4$.

$T_5 : V_5 \rightarrow V_5$ is such that $T_5$ maps any $8 \times 8$ matrix into the $8 \times 8$ diagonal matrix

$T_6 : V_6 \rightarrow V_6$ is such that

$$T_6 \begin{pmatrix} a_1 & a_2 & a_3 & a_4 \\ 0 & a_5 & a_6 & a_7 \\ 0 & 0 & a_8 & a_9 \\ 0 & 0 & 0 & a_{10} \end{pmatrix} \rightarrow \begin{pmatrix} a_1 & 0 & 0 & 0 \\ 0 & a_5 & 0 & 0 \\ 0 & 0 & a_8 & 0 \\ 0 & 0 & 0 & a_{10} \end{pmatrix}$$

$T_6$ is a neutrosophic linear operator on $V_6$. Thus $T = T_1 \cup T_2 \cup T_3 \cup T_4 \cup T_5 \cup T_6$ is a usual neutrosophic 6-linear operator on $V = V_1 \cup V_2 \cup V_3 \cup V_4 \cup V_5 \cup V_6$.

*Example 3.1.40:* Let $V = V_1 \cup V_2 \cup V_3 \cup V_4 \cup V_5 \cup V_6 \cup V_7$

$$= \left\{ \begin{pmatrix} a & b \\ c & d \end{pmatrix} \middle| a,b,c,d \in QI \right\} \cup$$

$$\left\{ \begin{pmatrix} a_1 & 0 & 0 \\ a_2 & a_3 & 0 \\ a_4 & a_5 & a_6 \end{pmatrix} \middle| a_i \in QI; 1 \leq i \leq 6 \right\} \cup$$

$$\left\{ \begin{pmatrix} a_1 & a_2 & a_3 \\ 0 & a_4 & a_5 \\ 0 & 0 & a_6 \end{pmatrix} \middle| a_i \in QI; 1 \le i \le 6 \right\} \cup$$

{All 5 × 5 neutrosophic matrices with entries from QI} ∪ {all 7 × 7 neutrosophic diagonal matrices with entries from QI} ∪

$$\left\{ \begin{pmatrix} a_1 & 0 & 0 & 0 \\ a_2 & a_3 & 0 & 0 \\ a_4 & a_5 & a_6 & 0 \\ a_7 & a_8 & a_9 & a_{10} \end{pmatrix} \middle| a_i \in QI; 1 \le i \le 10 \right\} \cup$$

$$\left\{ \begin{pmatrix} a_1 & a_2 & a_3 & a_4 \\ 0 & a_5 & a_6 & a_7 \\ 0 & 0 & a_8 & a_9 \\ 0 & 0 & 0 & a_{10} \end{pmatrix} \middle| a_i \in QI; 1 \le i \le 10 \right\}$$

be a neutrosophic 7-linear algebra over the field Q.

Let $W = W_1 \cup W_2 \cup W_3 \cup W_4 = \{(a, b, c, d, e, f) \mid a, b, c, d, e, f \in QI\} \cup$

$$\left\{ \begin{pmatrix} a & b \\ c & d \end{pmatrix} \middle| a, b, c, d \in QI \right\} \cup$$

{7 × 7 neutrosophic upper triangular matrices with entries from QI} ∪

$$\left\{ \begin{pmatrix} a_1 & 0 & 0 & 0 \\ a_2 & a_3 & 0 & 0 \\ a_4 & a_5 & a_6 & 0 \\ a_7 & a_8 & a_9 & a_{10} \end{pmatrix} \middle| a_i \in QI; 1 \le i \le 10 \right\}$$

be a neutrosophic 4-linear algebra over the real field Q.

Define $T = T_1 \cup T_2 \cup T_3 \cup T_4 \cup T_5 \cup T_6 \cup T_7 : V = V_1 \cup V_2 \cup V_3 \cup V_4 \cup V_5 \cup V_6 \cup V_7 \to W_1 \cup W_2 \cup W_3 \cup W_4$ as follows.

$T_1 : V_1 \to W_2$ where

$$T_1 \begin{pmatrix} a & b \\ c & d \end{pmatrix} = \begin{pmatrix} a+b & b+c \\ c+d & d+a \end{pmatrix},$$

$T_2 : V_2 \to W_1$ defined by

$$T_2 \begin{pmatrix} a_1 & 0 & 0 \\ a_2 & a_3 & 0 \\ a_4 & a_5 & a_6 \end{pmatrix} = (a_1, a_2, a_3, a_4, a_5, a_6),$$

$T_3 : V_3 \to W_1$ is defined by

$$T_3 \begin{pmatrix} a_1 & a_2 & a_3 \\ 0 & a_4 & a_5 \\ 0 & 0 & a_6 \end{pmatrix} = (a_1, a_2, a_3, a_4, a_5, a_6),$$

$T_4 : V_4 \to W_4$ is such that

$$T_4 \begin{pmatrix} a_1 & a_2 & a_3 & a_4 & a_5 \\ a_6 & a_7 & a_8 & a_9 & a_{10} \\ a_{11} & a_{12} & a_{13} & a_{14} & a_{15} \\ a_{16} & a_{17} & a_{18} & a_{19} & a_{20} \\ a_{21} & a_{22} & a_{23} & a_{24} & a_{25} \end{pmatrix} = \begin{pmatrix} a_1 & 0 & 0 & 0 \\ a_6 & a_7 & 0 & 0 \\ a_{11} & a_{12} & a_{13} & 0 \\ a_{21} & a_{22} & a_{24} & a_{25} \end{pmatrix},$$

$T_5 : V_5 \to W_3$ is defined by

$$T_5 \begin{pmatrix} a_1 & 0 & 0 & 0 & 0 & 0 & 0 \\ 0 & a_2 & 0 & 0 & 0 & 0 & 0 \\ 0 & 0 & a_3 & 0 & 0 & 0 & 0 \\ 0 & 0 & 0 & a_4 & 0 & 0 & 0 \\ 0 & 0 & 0 & 0 & a_5 & 0 & 0 \\ 0 & 0 & 0 & 0 & 0 & a_6 & 0 \\ 0 & 0 & 0 & 0 & 0 & 0 & a_7 \end{pmatrix} = \begin{pmatrix} a_1 & a_1 & a_1 & a_1 & a_1 & a_1 & a_1 \\ 0 & a_2 & a_2 & a_2 & a_2 & a_2 & a_2 \\ 0 & 0 & a_3 & a_3 & a_3 & a_3 & a_3 \\ 0 & 0 & 0 & a_4 & a_4 & a_4 & a_4 \\ 0 & 0 & 0 & 0 & a_5 & a_5 & a_5 \\ 0 & 0 & 0 & 0 & 0 & a_6 & a_6 \\ 0 & 0 & 0 & 0 & 0 & 0 & a_7 \end{pmatrix}$$

$T_6 : V_6 \to W_4$ defined by

$$T_6 \begin{pmatrix} a_1 & 0 & 0 & 0 \\ a_2 & a_3 & 0 & 0 \\ a_4 & a_5 & a_6 & 0 \\ a_7 & a_8 & a_9 & a_{10} \end{pmatrix} = \begin{pmatrix} a_{10} & 0 & 0 & 0 \\ a_9 & a_8 & 0 & 0 \\ a_7 & a_6 & a_5 & 0 \\ a_4 & a_3 & a_2 & a_1 \end{pmatrix}$$

and $T_7 : V_7 \to W_4$ is defined by

$$T_7 \begin{pmatrix} a_1 & a_2 & a_3 & a_4 \\ 0 & a_5 & a_6 & a_7 \\ 0 & 0 & a_8 & a_9 \\ 0 & 0 & 0 & a_{10} \end{pmatrix} = \begin{pmatrix} a_1 & 0 & 0 & 0 \\ a_2 & a_3 & 0 & 0 \\ a_4 & a_5 & a_6 & 0 \\ a_7 & a_8 & a_9 & a_{10} \end{pmatrix}$$

$T = T_1 \cup T_2 \cup T_3 \cup T_4 \cup T_5 \cup T_6 \cup T_7 : V = V_1 \cup V_2 \cup V_3 \cup V_4 \cup V_5 \cup V_6 \cup V_7 \to W_1 \cup W_2 \cup W_3 \cup W_4$ is a (7, 4) neutrosophic linear algebra transformation of V into W.

***Example 3.1.41:*** Let $V = V_1 \cup V_2 \cup V_3 \cup V_4 \cup V_5 =$

$$\left\{ \begin{pmatrix} a & b \\ c & d \end{pmatrix} \middle| a,b,c,d \in Z_7 I \right\} \cup$$

$\{(a_1, a_2, a_3, a_4, a_5, a_6, a_7, a_8) \mid a_i \in Z_7 I; 1 \le i \le 8\} \cup$

$$\left\{ \begin{pmatrix} a & b & c \\ d & e & f \\ g & h & i \end{pmatrix} \middle| a,b,c,d,e,f,g,h,i \in Z_7 I \right\} \cup$$

$$\left\{ \begin{pmatrix} a_1 & 0 & 0 & 0 \\ a_2 & a_3 & 0 & 0 \\ a_4 & a_5 & a_6 & 0 \\ a_7 & a_8 & a_9 & a_{10} \end{pmatrix} \middle| a_i \in Z_7I; 1 \leq i \leq 10 \right\} \cup$$

$\{5 \times 5$ diagonal neutrosophic matrices with entries from $Z_7I\}$ be a neutrosophic 5-linear algebra over the real field $Z_7$. Let $W = W_1 \cup W_2 \cup W_3 \cup W_4 \cup W_5 \cup W_6 \cup W_7 = \{8 \times 8$ neutrosophic diagonal matrices with entries from $Z_7I\} \cup$

$$\left\{ \begin{pmatrix} a & b \\ c & d \end{pmatrix} \middle| a,b,c,d \in Z_7I \right\} \cup$$

$$\left\{ \begin{pmatrix} a_1 & a_2 & a_3 & a_4 \\ 0 & a_5 & a_6 & a_7 \\ 0 & 0 & a_8 & a_9 \\ 0 & 0 & 0 & a_{10} \end{pmatrix} \middle| a_i \in Z_7I; 1 \leq i \leq 10 \right\} \cup$$

$\{(a_1, a_2, a_3, a_4, a_5) \mid a_i \in Z_7I; 1 \leq i \leq 5\} \cup \{4 \times 4$ neutrosophic diagonal matrices with entries from $Z_7I\} \cup$

$$\left\{ \begin{pmatrix} a & b & c \\ 0 & d & e \\ 0 & 0 & f \end{pmatrix} \middle| a,b,c,d,e,f, \in Z_7I \right\} \cup$$

$$\left\{ \begin{pmatrix} a_1 & 0 & 0 \\ a_2 & a_3 & 0 \\ a_4 & a_5 & a_6 \end{pmatrix} \middle| a_i \in Z_7I; 1 \leq i \leq 6 \right\}$$

be a neutrosophic 7-linear algebra over $Z_7$.

Define $T = T_1 \cup T_2 \cup T_3 \cup T_4 \cup T_5 : V = V_1 \cup V_2 \cup V_3 \cup V_4 \cup V_5 \to W_1 \cup W_2 \cup W_3 \cup W_4 \cup W_5 \cup W_6 \cup W_7$ as follows;

$$T_1 : V_1 \to W_5,$$
$$T_2 : V_2 \to W_1,$$
$$T_3 : V_3 \to W_7,$$
$$T_4 : V_4 \to W_3$$

and

$$T_5 : V_5 \to W_5$$

defined in the following way.

$T_1 : V_1 \to W_5$ is such that

$$T_1 \begin{pmatrix} a & b \\ c & d \end{pmatrix} = \begin{pmatrix} a & 0 & 0 & 0 \\ 0 & b & 0 & 0 \\ 0 & 0 & c & 0 \\ 0 & 0 & 0 & d \end{pmatrix},$$

$T_2 : V_2 \to W_1$ is defined as

$$T_2(a_1, a_2, a_3, a_4, a_5, a_6, a_7, a_8) = \begin{pmatrix} a_1 & 0 & 0 & 0 & 0 & 0 & 0 & 0 \\ 0 & a_2 & 0 & 0 & 0 & 0 & 0 & 0 \\ 0 & 0 & a_3 & 0 & 0 & 0 & 0 & 0 \\ 0 & 0 & 0 & a_4 & 0 & 0 & 0 & 0 \\ 0 & 0 & 0 & 0 & a_5 & 0 & 0 & 0 \\ 0 & 0 & 0 & 0 & 0 & a_6 & 0 & 0 \\ 0 & 0 & 0 & 0 & 0 & 0 & a_7 & 0 \\ 0 & 0 & 0 & 0 & 0 & 0 & 0 & a_8 \end{pmatrix},$$

$T_3 : V_3 \to W_7$ is given by

$$T_3 \begin{pmatrix} a & b & c \\ d & e & f \\ g & h & i \end{pmatrix} = \begin{pmatrix} a & 0 & 0 \\ b & c & 0 \\ d & e & f \end{pmatrix},$$

and $T_4 : V_4 \to W_3$ is defined as

$$T_4 \begin{pmatrix} a_1 & 0 & 0 & 0 \\ a_2 & a_3 & 0 & 0 \\ a_4 & a_5 & a_6 & 0 \\ a_7 & a_8 & a_9 & a_{10} \end{pmatrix} = \begin{pmatrix} a_1 & a_2 & a_3 & a_4 \\ 0 & a_5 & a_6 & a_7 \\ 0 & 0 & a_8 & a_9 \\ 0 & 0 & 0 & a_{10} \end{pmatrix},$$

$T_5 : V_5 \to W_4$ is expressed as

$$T_5 \begin{pmatrix} a_1 & 0 & 0 & 0 & 0 \\ 0 & a_2 & 0 & 0 & 0 \\ 0 & 0 & a_3 & 0 & 0 \\ 0 & 0 & 0 & a_4 & 0 \\ 0 & 0 & 0 & 0 & a_5 \end{pmatrix} = (a_1, a_2, a_3, a_4, a_5).$$

It is easily verified that $T = T_1 \cup T_2 \cup T_3 \cup T_4 \cup T_5$ is a (5, 7) neutrosophic linear transformation of V to W or neutrosophic (5, 7) linear transformation of V to W.

Now having seen several types of neutrosophic (m, n) linear transformation of neutrosophic m-linear algebra and neutrosophic n-linear algebra now we proceed onto define more properties about these neutrosophic n-linear transformation of neutrosophic n-linear algebra.

We have defined several subalgebraic structures of neutrosophic n-linear algebras (n-linear vector spaces) we now define subspace preserving n-linear operators in case of neutrosophic n-linear algebras (n-vector spaces).

**DEFINITION 3.1.14:** *Let $V = V_1 \cup V_2 \cup ... \cup V_n$ be a neutrosophic n-vector space over the real field F. Suppose $W = W_1 \cup W_2 \cup ... \cup W_n$ be a neutrosophic n-vector subspace of V over F. Let $T = T_1 \cup T_2 \cup ... \cup T_n$ be a neutrosophic n-linear operator on V. Suppose $T_i (W_i) \subseteq W_i$ for every $W_i$, $i = 1, 2, ..., n$ then we call $T = T_1 \cup T_2 \cup ... \cup T_n$ to be a vector subspace preserving neutrosophic n-linear operator on V.*

It is important to note that in general every neutrosophic n-linear operator on V need not preserve every neutrosophic n-vector subspace of V or even a single neutrosophic n-subspace of V.

We will say however the neutrosophic n-linear operator T which is the identity operator on V however preserves every neutrosophic n-vector subspace of V.

We will illustrate subspace preserving neutrosophic n-linear operator by an example.

*Example 3.1.42:* Let $V = V_1 \cup V_2 \cup V_3 \cup V_4 \cup V_5 =$

$$\left\{ \begin{pmatrix} a & b & c \\ d & e & f \\ g & h & i \end{pmatrix} \middle| a,b,c,d,e,f,g,h,i \in QI \right\} \cup$$

$$\left\{ \begin{pmatrix} a & b \\ c & d \end{pmatrix} \middle| a,b,c,d \in QI \right\} \cup$$

$$\left\{ \sum_{i=0}^{\infty} a_i x^i \middle| a_i \in QI; 0 \leq i \leq \infty \right\} \cup$$

$$\left\{ \begin{pmatrix} a & b & c & d \\ 0 & e & f & g \\ 0 & h & i & j \\ 0 & 0 & 0 & k \end{pmatrix} \middle| a,b,c,d,e,f,g,h,i,j,k \in QI \right\} \cup$$

$$\left\{ \begin{pmatrix} a_1 & a_2 & a_3 & a_4 \\ a_5 & a_6 & a_7 & a_8 \\ a_9 & a_{10} & a_{11} & a_{12} \end{pmatrix} \middle| a_i \in QI; 1 \leq i \leq 12 \right\}$$

be a neutrosophic 5-vector space over the field Q.
Take $W = W_1 \cup W_2 \cup W_3 \cup W_4 \cup W_5 =$

$$\left\{ \begin{pmatrix} a & a & a \\ a & a & a \\ a & a & a \end{pmatrix} \middle| a \in QI \right\} \cup$$

$$\left\{ \begin{pmatrix} a & b \\ 0 & c \end{pmatrix} \middle| a,b,c \in QI \right\} \cup \left\{ \sum_{i=0}^{\infty} a_i x^{2i} \middle| a_i \in QI; 0 \leq i \leq \infty \right\} \cup$$

$$\left\{ \begin{pmatrix} a & a & a & a \\ 0 & a & a & a \\ 0 & a & a & a \\ 0 & 0 & 0 & a \end{pmatrix} \middle| a \in QI \right\} \cup$$

$$\left\{ \begin{pmatrix} a_1 & a_2 & a_3 & a_4 \\ a_1 & a_2 & a_3 & a_4 \\ a_1 & a_2 & a_3 & a_4 \end{pmatrix} \middle| a_i \in QI; 1 \leq i \leq 4 \right\}$$

$\subseteq V_1 \cup V_2 \cup V_3 \cup V_4 \cup V_5$, be a neutrosophic 5-vector subspace of V over the field Q.

Define $T = T_1 \cup T_2 \cup T_3 \cup T_4 \cup T_5 : V = V_1 \cup V_2 \cup V_3 \cup V_4 \cup V_5 \to V = V_1 \cup V_2 \cup V_3 \cup V_4 \cup V_5$ as follows. $T_i : V_i \to V_i$, $i = 1, 2, \ldots, 5$ is a neutrosophic linear operator for each i and $T_i$ is defined as follows;

$T_1 : V_1 \to V_1$ is such that

$$T_1 \begin{pmatrix} a & b & c \\ d & e & f \\ g & h & i \end{pmatrix} = \begin{pmatrix} a & a & a \\ a & a & a \\ a & a & a \end{pmatrix}$$

$T_2 : V_2 \to V_2$ is given by
$$T_2 \begin{pmatrix} a & b \\ c & d \end{pmatrix} = \begin{pmatrix} a & b \\ 0 & d \end{pmatrix}.$$

$T_3 : V_3 \to V_3$ is defined by

$$T_3 \left( \sum_{i=0}^{\infty} a_i x^i \right) = \sum_{i=0}^{\infty} a_i x^{2i}$$

that is $Ix^i \to Ix^{2i}$ for every $i = 0, 1, \ldots, \infty$.

$T_4 : V_4 \to V_4$ is such that

$$T_4 \begin{pmatrix} a_1 & a_2 & a_3 & a_4 \\ 0 & a_5 & a_6 & a_7 \\ 0 & 0 & a_8 & a_9 \\ 0 & 0 & 0 & a_{10} \end{pmatrix} = \begin{pmatrix} a & a & a & a \\ 0 & a & a & a \\ 0 & 0 & a & a \\ 0 & 0 & 0 & a \end{pmatrix}$$

and $T_5 : V_5 \to V_5$ is defined by

$$T_5 \begin{pmatrix} a_1 & a_2 & a_3 & a_4 \\ a_5 & a_6 & a_7 & a_8 \\ a_9 & a_{10} & a_{11} & a_{12} \end{pmatrix} = \begin{pmatrix} a_1 & a_2 & a_3 & a_4 \\ a_1 & a_2 & a_3 & a_4 \\ a_1 & a_2 & a_3 & a_4 \end{pmatrix}.$$

It is easily verified that $T = T_1 \cup T_2 \cup T_3 \cup T_4 \cup T_5$ is a neutrosophic 5-linear operator on V. Further this T preserves the neutrosophic 5-vector subspace $W = W_1 \cup W_2 \cup W_3 \cup W_4 \cup W_5$. It is easily verified that $T_i (W_i) \subseteq W_i$ for $i = 1, 2, \ldots, 5$. Hence the claim.

We see in general all neutrosophic n-linear operators T on V need not preserve a neutrosophic n-vector subspace of V.

We will illustrate this situation by an example.

*Example 3.1.43:* Let $V = V_1 \cup V_2 \cup V_3 \cup V_4 =$

$$\left\{ \begin{pmatrix} a_1 & a_2 & a_3 \\ a_4 & a_5 & a_6 \\ a_7 & a_8 & a_9 \end{pmatrix} \middle| a_i \in Z_{11}I; 1 \leq i \leq 9 \right\} \cup$$

$$\left\{ \begin{pmatrix} a_1 & a_2 & a_3 & a_4 & a_5 \\ a_6 & a_7 & a_8 & a_9 & a_{10} \end{pmatrix} \middle| a_i \in Z_{11}I; 1 \leq i \leq 10 \right\} \cup$$

$$\left\{ \sum_{i=0}^{\infty} a_i x^i \middle| a_i \in Z_{11}I; 0 \leq i \leq \infty \right\} \cup$$

$$\left\{ \begin{pmatrix} a_1 & a_2 & a_3 \\ a_4 & a_5 & a_6 \\ a_7 & a_8 & a_9 \\ a_{10} & a_{11} & a_{12} \\ a_{13} & a_{14} & a_{15} \end{pmatrix} \middle| a_i \in Z_{11}I; 1 \leq i \leq 15 \right\}$$

be a neutrosophic 4-vector space over the field $Z_{11}$. Take $W = W_1 \cup W_2 \cup W_3 \cup W_4 =$

$$\left\{ \begin{pmatrix} 0 & 0 & 0 \\ a_1 & a_2 & a_3 \\ 0 & 0 & 0 \end{pmatrix} \middle| a_i \in Z_{11}I; 1 \leq i \leq 3 \right\} \cup$$

$$\left\{ \begin{pmatrix} 0 & 0 & 0 & 0 & 0 \\ a_1 & a_2 & a_3 & a_4 & a_5 \end{pmatrix} \middle| a_i \in Z_{11}I; 1 \leq i \leq 5 \right\} \cup$$

$$\left\{ \sum_{i=0}^{\infty} a_i x^{2i} \middle| a_i \in Z_{11}I; 0 \leq i \leq \infty \right\} \cup$$

$$\left\{ \begin{pmatrix} a_1 & 0 & a_6 \\ a_2 & 0 & a_7 \\ a_3 & 0 & a_8 \\ a_4 & 0 & a_9 \\ a_5 & 0 & a_{10} \end{pmatrix} \middle| a_i \in Z_{11}I; 1 \leq i \leq 10 \right\}$$

$\subseteq V_1 \cup V_2 \cup V_3 \cup V_4$ to be a neutrosophic 4-subspace of V over $Z_{11}$. Define $T = T_1 \cup T_2 \cup T_3 \cup T_4 : V = V_1 \cup V_2 \cup V_3 \cup V_4 \to V = V_1 \cup V_2 \cup V_3 \cup V_4$ as $T_i : V_i \to V_i$; i = 1 to 4 as follows.
$T_1 : V_1 \to V_1$; where

$$T_1 \begin{pmatrix} a_1 & a_2 & a_3 \\ a_4 & a_5 & a_6 \\ a_7 & a_8 & a_9 \end{pmatrix} = \begin{pmatrix} 0 & 0 & 0 \\ a_4 & a_5 & a_6 \\ 0 & 0 & 0 \end{pmatrix},$$

$T_2 : V_2 \to V_2$ as

$$T_2 \begin{pmatrix} a_1 & a_2 & a_3 & a_4 & a_5 \\ a_6 & a_7 & a_8 & a_9 & a_{10} \end{pmatrix} = \begin{pmatrix} 0 & 0 & 0 & 0 & 0 \\ a_6 & a_7 & a_8 & a_9 & a_{10} \end{pmatrix},$$

$T_3 : V_3 \to V_3$ is defined by

$$T_3 \left( \sum_{i=0}^{\infty} a_i x^i \right) = \sum_{i=0}^{\infty} a_i x^{2i},$$

that is $a_i x^i \mapsto a_i x^{2i}$ for $i = 0, 1, \ldots, \infty$ and $T_4 : V_4 \to V_4$ is such that

$$T_4 \begin{pmatrix} a_1 & a_{15} & a_6 \\ a_2 & a_{14} & a_7 \\ a_3 & a_{13} & a_8 \\ a_4 & a_{11} & a_9 \\ a_5 & a_{12} & a_{10} \end{pmatrix} = \begin{pmatrix} a_1 & 0 & a_6 \\ a_2 & 0 & a_7 \\ a_3 & 0 & a_8 \\ a_4 & 0 & a_9 \\ a_5 & 0 & a_{10} \end{pmatrix}.$$

It is easily verified $T = T_1 \cup T_2 \cup T_3 \cup T_4 : V \to V$ is a neutrosophic 4-linear operator on V and it preserve the neutrosophic 4-subspace $W = W_1 \cup W_2 \cup W_3 \cup W_4$ that is $T(W) = T_1(W_1) \cup T_2(W_2) \cup T_3(W_3) \cup T_4(W_4) \subseteq W_1 \cup W_2 \cup W_3 \cup W_4$.

Consider $P = P_1 \cup P_2 \cup P_3 \cup P_4 : V = V_1 \cup V_2 \cup V_3 \cup V_4 \to V_1 \cup V_2 \cup V_3 \cup V_4$ where $P_i : V_i \to V_i$; $i = 1, 2, 3, 4$ are neutrosophic linear operators defined as follows.

$P_1 : V_1 \to V_1$ such that

$$P_1 \begin{pmatrix} a_1 & a_2 & a_3 \\ a_4 & a_5 & a_6 \\ a_7 & a_8 & a_9 \end{pmatrix} = \begin{pmatrix} a_1 & 0 & 0 \\ 0 & a_5 & 0 \\ 0 & 0 & a_9 \end{pmatrix},$$

$P_2 : V_2 \to V_2$ is defined by

$$P_2 \begin{pmatrix} a_1 & a_2 & a_3 & a_4 & a_5 \\ a_6 & a_7 & a_8 & a_9 & a_{10} \end{pmatrix} = \begin{pmatrix} a_1 & 0 & a_3 & 0 & a_5 \\ a_6 & 0 & a_8 & 0 & a_{10} \end{pmatrix},$$

$P_3 : V_3 \to V_3$ is given by

$$P_3 \left( \sum_{i=0}^{\infty} a_i x^i \right) = a_1 x + a_3 x^3 + a_5 x^5 + \ldots + a_{2n+1} x^{2n+1} + \ldots$$

and $P_4 : V_4 \to V_4$ is defined by

$$P_4 \begin{pmatrix} a_1 & a_2 & a_3 \\ a_4 & a_5 & a_6 \\ a_7 & a_8 & a_9 \\ a_{10} & a_{11} & a_{12} \\ a_{13} & a_{14} & a_{15} \end{pmatrix} = \begin{pmatrix} 0 & a_2 & 0 \\ 0 & a_5 & 0 \\ 0 & a_8 & 0 \\ 0 & a_{11} & 0 \\ 0 & a_{14} & 0 \end{pmatrix}.$$

The neutrosophic 4-linear operator P on V does not preserve the neutrosophic 4-subspace $W = W_1 \cup W_2 \cup W_3 \cup W_4 \subseteq V_1 \cup V_2 \cup V_3 \cup V_4$. Thus $P = P_1 \cup P_2 \cup P_3 \cup P_4$ is neutrosophic 4-linear operator which does not preserve the subspace W.

It can be easily proved that for V and W any two neutrosophic n-vector spaces over a real field F if T and S are neutrosophic n-linear transformations of V to W then (T+S) is a neutrosophic n-linear transformation of V to W. Further if T is a neutrosophic n-linear transformation of V to W for any c an element of F, the function cT defined by $(cT)\alpha = cT\alpha$ is again a neutrosophic n-linear transformation of V to W. It is interesting to note the set of all neutrosophic n-linear transformations from V into W with addition and scalar multiplication defined above is a neutrosophic n-vector space over F.

Further as in case of usual n-vector spaces of type I we see in case of two neutrosophic n-vector spaces V and W over the field F both V and W are of n-finite dimension say $(n_1, n_2, \ldots, n_n)$ and $(m_1, m_2, \ldots, m_n)$ over the field F, and if $T = T_1 \cup T_2 \cup \ldots \cup T_n$ is a neutrosophic n-linear transformation of V into W, where $T_i : V_i \to W_j$ (That is no two $V_i$'s are mapped on to the same $W_j$) true for $1 \le i, j \le n$; then n rank T + n nullity T = n dim V; that is rank $T_1 \cup \ldots \cup$ rank $T_n$ + (nullity $T_1 \cup \ldots \cup$ nullity $T_n$) = dim $V_1 \cup \ldots \cup$ dim $V_n$; that is (rank $T_1$ + nullity $T_1$) $\cup$ (rank $T_2$ + nullity $T_2$) $\cup \ldots \cup$ (rank $T_n$ + nullity $T_n$) = dim $V_1 \cup$ dim $V_2 \cup \ldots \cup$ dim $V_n$.

We can also prove the result which is as follows:

Let $V = V_1 \cup \ldots \cup V_n$ and $W = W_1 \cup W_2 \cup \ldots \cup W_n$ be two neutrosophic n-vector spaces over the field F. Let
$$B = \left\{ \left( \alpha_1^1, \alpha_2^1, \ldots, \alpha_{n_1}^1 \right) \cup \left( \alpha_1^2, \alpha_2^2, \ldots, \alpha_{n_2}^2 \right) \cup \ldots \cup \left( \alpha_1^n, \alpha_2^n, \ldots, \alpha_{n_n}^n \right) \right\}$$
be a n-basis of $V = V_1 \cup V_2 \cup \ldots \cup V_n$; i.e., $\left( \alpha_1^i, \alpha_2^i, \ldots, \alpha_{n_i}^i \right)$ is a basis of $V_i$; i = 1, 2, …, n. Let
$$C = \left\{ \left( \beta_1^1, \beta_2^1, \ldots, \beta_{n_1}^1 \right) \cup \left( \beta_1^2, \beta_2^2, \ldots, \beta_{n_2}^2 \right) \cup \ldots \cup \left( \beta_1^n, \beta_2^n, \ldots, \beta_{n_n}^n \right) \right\}$$
be any n-vector in $W = W_1 \cup W_2 \cup \ldots \cup W_n$ then there is precisely only one neutrosophic linear n-transformation $T = T_1 \cup T_2 \cup \ldots \cup T_n$ from V onto W such that $T\alpha_j^i = \beta_j^i$; j = 1, 2, …, $n_i$; $1 \le i \le n$.

The following result is also true and it can be proved as in the case of n-vector spaces.

Let $V = V_1 \cup V_2 \cup \ldots \cup V_n$ be a neutrosophic $(n_1, n_2, \ldots, n_n)$ n-dimensional vector space over a field F. If $W = W_1 \cup W_2 \cup \ldots \cup W_n$ is a neutrosophic $(m_1, m_2, \ldots, m_n)$ n-dimensional vector space over the field F. Then $NL^n$ (V, W) = {collection of all neutrosophic n-linear transformations from V to W} is finite $(m_1n_1, m_2n_2, \ldots, m_nn_n)$ dimensional neutrosophic n-vector space over F. It is to be noted that as in case of n-vector spaces we in case of neutrosophic n-vector spaces also define neutrosophic n-linear transformations V to W where V is a neutrosophic n-vector space where as W is a neutrosophic m-vector space $m > n$.

Let V and W be two neutrosophic n-vector space and neutrosophic m-vector space respectively ($m > n$) over a real field F. Let $T = T_1 \cup T_2 \cup \ldots \cup T_n$ be a neutrosophic n-linear transformation of V into W such that $T_i : V_i \to W_j$, where each $V_i$ is mapped onto a distinct $W_j$ $1 \le i \le n$ and $1 \le j \le m$ that is no two $V_i$'s are mapped onto the same $W_j$ true for each i, $i = 1, 2, \ldots, n$. Then $NL^n(V, W)$ can be defined and in this case $NL^n(V, W)$ is finite dimensional n-space over F of n-dimension $\left(m_{i_1}n_1, m_{i_2}n_2, \ldots, m_{i_n}n_n\right)$ where $1 \le i_1, i_2, \ldots, i_n \le m$.

Now as in case of n-vector spaces we can define for neutrosophic n-vector spaces composition of neutrosophic n-linear transformations.

Let V, W and Z be three neutrosophic n-vector spaces over the field F, i.e., $V = V_1 \cup V_2 \cup \ldots \cup V_n$, $W = W_1 \cup W_2 \cup \ldots \cup W_n$ and $Z = Z_1 \cup Z_2 \cup \ldots \cup Z_n$.

Define $T = T_1 \cup T_2 \cup \ldots \cup T_n : V = V_1 \cup V_2 \cup \ldots \cup V_n \to W_1 \cup W_2 \cup \ldots \cup W_n$; $T_i : V_i \to W_j$; $i = 1, 2, \ldots, n$; $1 \le j \le n$. Let $P = P_1 \cup P_2 \cup \ldots \cup P_n : W = W_1 \cup W_2 \cup \ldots \cup W_n \to Z = Z_1 \cup Z_2 \cup \ldots \cup Z_n$; $P_j : W_j \to Z_k$; $j = 1, 2, \ldots, n$ and $1 \le k \le m$ so that no two subspaces $W_j$ are mapped on to same $Z_k$; $k = 1, 2, \ldots, n$. Now

$(P_j T_i) (c\alpha^i + \beta^i)$
$\quad = P_j [T_i (c\alpha^i + \beta^i)]$
$\quad = P_j T_i (c\alpha^i) + P_j (\beta^i)$
$\quad = P_j [c\omega^i + \delta^j]$ (as $T_i : V_i \to W_j$; $\delta^j, \omega^j \in W_j$)
$\quad = cP_j (\omega^j) + P_j (\delta^j)$
$\quad = ca^k + b^k$; $a^k, b^k \in Z_k$.

Thus $P_jT_i$ is a neutrosophic n-linear transformation from $W_j$ to $Z_k$. Hence the claim and the result is true for each i and j. Thus PT is a neutrosophic n-linear transformation from W to Z.

So PT = $(P_1 \cup P_2 \cup \ldots \cup P_n)(T_1 \cup T_2 \cup \ldots \cup T_n) = P_1 T_{i_1} \cup P_2 T_{i_2} \cup \ldots \cup P_n T_{i_n}$ where $(i_1, i_2, \ldots, i_n)$ is a permutation of $(1, 2, \ldots, n)$.

Now we for the notational convenience recall that if $V = V_1 \cup V_2 \cup \ldots \cup V_n$ is a neutrosophic n-vector space over the field F then $V_i$'s will be known as the component neutrosophic subvector space of V. $V_i$'s are also known as the component of V. Now we proceed onto give some properties of neutrosophic n-linear operators.

Let $V = V_1 \cup V_2 \cup \ldots \cup V_n$ be a neutrosophic n-vector space over the field F. Let $T = T_1 \cup T_2 \cup \ldots \cup T_n$ be a neutrosophic n-linear operator on V with $T_i : V_i \to V_i$, $i = 1, 2, \ldots, n$. Let $S = S_1 \cup S_2 \cup \ldots \cup S_n$ be another neutrosophic n-linear operator on V with $S_i : V_i \to V_i$, $i = 1, 2, \ldots, n$. Now ST and TS is again neutrosophic n-linear operators on V.

Thus the neutrosophic n-space of all neutrosophic n-linear operators has a product defined as composition.

In this case the neutrosophic n-linear operator TS is also defined. In general ST ≠ TS; that is ST – TS ≠ 0.

$NL^n(V, V)$ is a neutrosophic n-vector space of n-dimension $\left(n_1^2, n_2^2, \ldots, n_n^2\right)$ where the n-dimension of V is $(n_1, n_2, \ldots, n_n)$. All relations like n-nilpotent, n-diagonalizable can be defined in case of neutrosophic n-vector spaces of type I with appropriate modifications.

## 3.2 Neutrosophic Strong n-Vector Spaces

In this section we proceed onto define the new notion of neutrosophic strong n-vector spaces ($n \geq 3$) and discuss a few of their properties.

**DEFINITION 3.2.1:** *Let $V = V_1 \cup V_2 \cup \ldots \cup V_n$ be such that each $V_i$ is a neutrosophic vector space over the same neutrosophic field F then we call V to be a neutrosophic strong n-vector space or strong neutrosophic n-vector space.*

We will first illustrate this by some examples.

***Example 3.2.1:*** Let $V = V_1 \cup V_2 \cup V_3 \cup V_4 \cup V_5 =$

$$\left\{ \begin{pmatrix} a & b \\ c & d \end{pmatrix} \middle| a,b,c,d \in N(Q) \right\} \cup$$

$$\left\{ \begin{pmatrix} a_1 & a_2 & a_3 \\ a_4 & a_5 & a_6 \end{pmatrix} \middle| a_i \in QI; 1 \le i \le 6 \right\} \cup$$

$$\left\{ \begin{pmatrix} a_1 & a_2 \\ a_3 & a_4 \\ a_5 & a_6 \\ a_7 & a_8 \\ a_9 & a_{10} \end{pmatrix} \middle| a_i \in QI; 1 \le i \le 10 \right\} \cup$$

{QI[x]; all polynomials in the variable x with coefficients from QI} $\cup$ (QI $\times$ QI $\times$ QI $\times$ QI $\times$ QI) = {(a, b, c, d, e) | a, b, c, d, e $\in$ QI} be a neutrosophic strong 5-vector space over the neutrosophic field F = QI.

***Example 3.2.2:*** Let $V = V_1 \cup V_2 \cup V_3 \cup V_4 \cup V_5 \cup V_6 \cup V_7 =$

$$\left\{ \begin{pmatrix} a_1 & a_2 & a_3 & a_4 \\ a_5 & a_6 & a_7 & a_8 \end{pmatrix} \middle| a_i \in Z_2 I; 1 \le i \le 8 \right\} \cup$$

$$\left\{ \begin{pmatrix} a_1 & a_2 & a_3 \\ a_4 & a_5 & a_6 \\ a_7 & a_8 & a_9 \\ a_{10} & a_{11} & a_{12} \\ a_{13} & a_{14} & a_{15} \end{pmatrix} \middle| a_i \in N(Z_2); 1 \le i \le 15 \right\} \cup$$

$$\left\{ \begin{pmatrix} a_1 & a_2 & a_3 \\ a_4 & a_5 & a_6 \\ a_7 & a_8 & a_9 \end{pmatrix} \middle| a_i \in Z_2 I; 1 \le i \le 9 \right\} \cup$$

$\{(x_1, x_2, x_3, x_4, x_5, x_6, x_7) \mid x_i \in Z_2 I;\ 1 \le i \le 7\} \cup \{9 \times 9$ upper triangular matrices with entries from $N(Z_2)\} \cup \{N(Z_2)[x]$; all polynomials in the variable x with coefficients from $N(Z_2)\} \cup$

$$\left\{ \begin{pmatrix} a_1 & a_2 \\ a_3 & a_4 \\ a_5 & a_6 \\ a_7 & a_8 \\ a_9 & a_{10} \\ a_{11} & a_{12} \end{pmatrix} \middle| a_i \in Z_2 I; 1 \le i \le 12 \right\}$$

be a strong neutrosophic 7-vector space over the neutrosophic field $F = Z_2 I$.

Now having seen examples of strong neutrosophic n-vector spaces $n \ge 3$ we define substructures in them. It is both interesting and important to note that in a strong neutrosophic n-vector space V if $n = 2$ we get the strong neutrosophic bivector space defined and discussed in chapter two of this book. When $n = 3$ we call V to be a strong neutrosophic trivector space.

*Example 3.2.3:* Let $V = V_1 \cup V_2 \cup V_3 =$

$$\left\{ \begin{pmatrix} a & b \\ c & d \end{pmatrix} \middle| a, b, c, d \in RI \right\} \cup$$

$$\left\{ \begin{pmatrix} a_1 & a_2 & a_3 & a_4 & a_5 \\ a_6 & a_7 & a_8 & a_9 & a_{10} \end{pmatrix} \middle| a_i \in QI; 1 \le i \le 10 \right\} \cup$$

$$\left\{ \begin{pmatrix} a_1 & a_2 \\ a_3 & a_4 \\ a_5 & a_6 \\ a_7 & a_8 \end{pmatrix} \middle| a_i \in N(R); 1 \le i \le 8 \right\};$$

V is a strong neutrosophic trivector space (3-vector space) over the neutrosophic field F = QI.

**DEFINITION 3.2.2:** *Let $V = V_1 \cup V_2 \cup ... \cup V_n$ be a strong neutrosophic n-vector space over the neutrosophic field F. Let $W = W_1 \cup W_2 \cup ... \cup W_n \subseteq V_1 \cup V_2 \cup ... \cup V_n$ such that $\phi \ne W_i \not\subset V_i$; i = 1, 2, ..., n be a strong neutrosophic n-vector space over the same neutrosophic field F. Then we call W to be a strong neutrosophic n-vector subspace of V over F.*

We note that even if one $W_i = \phi$ or $W_i = V_i$ ($1 \le i \le n$) then we do not call $W = W_1 \cup W_2 \cup ... \cup W_n$ to be a strong neutrosophic n-vector subspace of V.

We will first illustrate this situation by some examples.

*Example 3.2.4:* Let $V = V_1 \cup V_2 \cup V_3 \cup V_4 \cup V_5 =$

$$\left\{ \begin{pmatrix} a & b \\ c & d \end{pmatrix} \middle| a,b,c,d \in Z_7 I \right\} \cup$$

$$\left\{ \begin{pmatrix} a_1 & a_2 & a_3 & a_4 \\ a_5 & a_6 & a_7 & a_8 \end{pmatrix} \middle| a_i \in Z_7 I; 1 \le i \le 8 \right\} \cup$$

$$\left\{ \sum_{i=0}^{\infty} a_i x^i \middle| a_i \in Z_7 I; i = 0,1,2,...,\infty \right\} \cup$$

$$\left\{ \begin{pmatrix} a & b & c \\ d & e & f \\ g & h & i \\ k & l & m \\ o & p & q \end{pmatrix} \middle| a, b, \ldots, p, q \in Z_7I \right\}$$

$\cup \ \{(x_1, x_2, x_3, x_4, x_5, x_6) \mid x_i \in N(Z_7); 1 \leq i \leq 6\}$ be a strong neutrosophic 5-vector space over the neutrosophic field $Z_7I$. Let $W = W_1 \cup W_2 \cup W_3 \cup W_4 \cup W_5 =$

$$\left\{ \begin{pmatrix} a & a \\ a & a \end{pmatrix} \middle| a \in Z_7I \right\} \cup \left\{ \begin{pmatrix} a & a & a & a \\ b & b & b & b \end{pmatrix} \middle| a, b \in Z_7I \right\} \cup$$

$$\left\{ \sum_{i=0}^{51} a_i x^i \middle| a_i \in Z_7I; 0 \leq i \leq 51 \right\} \cup$$

$$\left\{ \begin{pmatrix} a & a & a \\ a & a & a \\ a & a & a \\ a & a & a \\ a & a & a \end{pmatrix} \middle| a \in Z_7I \right\} \cup$$

$\{(x_1, x_2, x_3, x_4, x_5, x_6) \mid x_i \in Z_7I; 1 \leq i \leq 6\} \subseteq V_1 \cup V_2 \cup V_3 \cup V_4 \cup V_5$. It is easily verified that W is a strong neutrosophic 5-vector subspace of V over the neutrosophic field $Z_7I$.

*Example 3.2.5:* Let $V = V_1 \cup V_2 \cup V_3 \cup V_4 \cup V_5 \cup V_6 \cup V_7 =$

$$\left\{ \begin{pmatrix} a_1 & a_2 \\ a_3 & a_4 \end{pmatrix} \middle| a_i \in N(Q); 1 \leq i \leq 4 \right\} \cup$$

$$\left\{\sum_{i=0}^{\infty} a_i x^i \,\middle|\, a_i \in N(Q); i = 0, 1, 2, \ldots, \infty \right\};$$

all polynomials in the variable x with coefficients from the neutrosophic field N (Q)} $\cup$

$$\left\{ \begin{pmatrix} a_1 \\ a_2 \\ a_3 \\ a_4 \\ a_5 \end{pmatrix} \,\middle|\, a_i \in QI; 1 \leq i \leq 5 \right\} \cup$$

$$\left\{ \begin{pmatrix} a_1 & a_2 & a_3 & a_4 \\ a_5 & a_6 & a_7 & a_8 \end{pmatrix} \,\middle|\, a_i \in QI; 1 \leq i \leq 8 \right\} \cup$$

{all 8 × 8 neutrosophic matrices with entries from N(Q)} $\cup$ {All 4 × 4 lower triangular matrices with entries from N(Q)} $\cup$ {all 7 × 7 upper triangular matrices with entries from N(Q)} be a neutrosophic strong 7-vector space over the neutrosophic field QI. Let $W = W_1 \cup W_2 \cup W_3 \cup W_4 \cup W_5 \cup W_6 \cup W_7 =$

$$\left\{ \begin{pmatrix} a & a \\ a & a \end{pmatrix} \,\middle|\, a \in N(Q) \right\} \cup \left\{ \sum_{i=0}^{8} a_i x^i \,\middle|\, a_i \in N(Q); i = 0, 1, 2, \ldots, 8 \right\};$$

all polynomials in the variable x with coefficients from the neutrosophic field N(Q) of degree less than or equal to 8} $\cup$

$$\left\{ \begin{pmatrix} a \\ a \\ a \\ a \\ a \end{pmatrix} \,\middle|\, a \in QI \right\} \cup \left\{ \begin{pmatrix} a & a & a & a \\ b & b & b & b \end{pmatrix} \,\middle|\, a, b \in QI \right\} \cup$$

{All 8 × 8 neutrosophic matrices with entries from the neutrosophic field QI} ∪ {All 4 × 4 lower triangular matrices with entries from the neutrosophic field QI} ∪ {All 7 × 7 diagonal matrices with entries from the neutrosophic field N (Q)} ⊆ $V_1 \cup V_2 \cup V_3 \cup V_4 \cup V_5 \cup V_6 \cup V_7$; it is easily verified that W is a strong neutrosophic 7-vector subspace of V over the neutrosophic field QI.

We say a strong neutrosophic n-vector space is simple if V has no proper strong neutrosophic n-vector subspaces.

We will illustrate this by some examples.

***Example 3.2.6:*** Let $V = V_1 \cup V_2 \cup V_3 \cup V_4 =$

$$\left\{ \begin{pmatrix} a & a & a \\ a & a & a \end{pmatrix} \middle| a \in Z_3I \right\} \cup \left\{ \begin{pmatrix} a & a & a \\ a & a & a \\ a & a & a \\ a & a & a \\ a & a & a \end{pmatrix} \middle| a \in Z_3I \right\} \cup$$

{(a a a a a a) | $a \in Z_3I$} ∪ {All 8×8 neutrosophic diagonal matrices of the form $Z_3I$, that is

$$\left\{ \begin{pmatrix} a & 0 & 0 & 0 & 0 & 0 & 0 & 0 \\ 0 & a & 0 & 0 & 0 & 0 & 0 & 0 \\ 0 & 0 & a & 0 & 0 & 0 & 0 & 0 \\ 0 & 0 & 0 & a & 0 & 0 & 0 & 0 \\ 0 & 0 & 0 & 0 & a & 0 & 0 & 0 \\ 0 & 0 & 0 & 0 & 0 & a & 0 & 0 \\ 0 & 0 & 0 & 0 & 0 & 0 & a & 0 \\ 0 & 0 & 0 & 0 & 0 & 0 & 0 & a \end{pmatrix} \middle| a \in Z_3I \right\}$$

be a strong neutrosophic 4-vector space over the neutrosophic field $Z_3I$. Clearly V has no strong neutrosophic 4-vector subspace. Hence V is a simple neutrosophic strong 4-vector space or simple strong neutrosophic 4-vector space over $Z_3I$.

*Example 3.2.7:* Let $V = V_1 \cup V_2 \cup V_3 \cup V_4 \cup V_5 \cup V_6 \cup V_7 =$

$$\{(a\ a\ a\ a\ a\ a\ a\ a\ a) \mid a \in Z_7I\} \cup$$

$$\left\{ \begin{pmatrix} a & a \\ a & a \\ a & a \end{pmatrix} \middle| a \in Z_7I \right\} \cup$$

$$\left\{ \begin{pmatrix} a & 0 & 0 & 0 & 0 \\ a & a & 0 & 0 & 0 \\ a & a & a & 0 & 0 \\ a & a & a & a & 0 \\ a & a & a & a & a \end{pmatrix} \middle| a \in Z_7I \right\} \cup$$

$$\left\{ \begin{pmatrix} a & a & a & a & a & a & a \\ a & a & a & a & a & a & a \\ a & a & a & a & a & a & a \\ a & a & a & a & a & a & a \end{pmatrix} \middle| a \in Z_7I \right\} \cup$$

$$\left\{ \begin{pmatrix} a & a \\ a & a \\ a & a \\ a & a \\ a & a \\ a & a \end{pmatrix} \middle| a \in Z_7I \right\} \cup \left\{ \begin{pmatrix} a & 0 & 0 & 0 & 0 \\ 0 & a & 0 & 0 & 0 \\ 0 & 0 & a & 0 & 0 \\ 0 & 0 & 0 & a & 0 \\ 0 & 0 & 0 & 0 & a \end{pmatrix} \middle| a \in Z_7I \right\} \cup$$

$$\left\{ \begin{pmatrix} a & a & a \\ 0 & 0 & 0 \\ a & a & a \\ 0 & 0 & 0 \end{pmatrix} \middle| a \in Z_7 I \right\}$$

be a strong neutrosophic 7-vector space over the neutrosophic field $Z_7 I$. It is easily verified that V has no proper strong neutrosophic 7-vector subspace. Hence V is a simple strong neutrosophic 7-vector space over $Z_7 I$.

**DEFINITION 3.2.3:** *Let $V = V_1 \cup V_2 \cup V_3 \cup \ldots \cup V_n$ be a strong neutrosophic n-vector space over a neutrosophic field F. If each $V_i$ is a neutrosophic strong linear algebra over F then we call $V = V_1 \cup V_2 \cup \ldots \cup V_n$ to be a strong neutrosophic n-linear algebra over F.*

We will illustrate then by some examples.

***Example 3.2.8:*** Let $V = V_1 \cup V_2 \cup V_3 \cup V_4 \cup V_5 =$

$$\left\{ \begin{pmatrix} a & b \\ c & d \end{pmatrix} \middle| a,b,c,d \in N(Q) \right\} \cup \left\{ \sum_{i=0}^{\infty} a_i x^i \middle| a_i \in N(Q) \right\} \cup$$

$$\left\{ \begin{pmatrix} a_1 & a_2 & a_3 \\ a_4 & a_5 & a_6 \\ a_7 & a_8 & a_9 \end{pmatrix} \middle| a_i \in QI; 1 \leq i \leq 9 \right\} \cup$$

$$\left\{ \begin{pmatrix} a_1 & 0 & 0 & 0 \\ 0 & a_2 & 0 & 0 \\ 0 & 0 & a_3 & 0 \\ 0 & 0 & 0 & a_4 \end{pmatrix} \middle| a_i \in QI \right\} \cup$$

{All $8 \times 8$ upper triangular matrices with entries from $N(Q)$} be a strong neutrosophic 5-linear algebra over QI.

***Example 3.2.9:*** Let $V = V_1 \cup V_2 \cup V_3 \cup V_4 \cup V_5 \cup V_6 \cup V_7 =$

$$\left\{ \begin{pmatrix} a & a \\ a & a \end{pmatrix} \middle| a_i \in Z_{19}I \right\} \cup$$

$$\left\{ \begin{pmatrix} a & 0 & 0 \\ a & a & 0 \\ a & a & a \end{pmatrix} \middle| a \in N(Z_{19}) \right\} \cup$$

{All 9×9 upper triangular neutrosophic matrices with entries from $Z_{19}I$} $\cup$ {$(x_1, x_2, x_3, x_y) \mid x_i \in Z_{19}I; 1 \le i \le 4$} $\cup$ {all 5×5 lower triangular matrices with entries from $Z_{19}I$} $\cup$ {All 12 × 12 neutrosophic diagonal matrices with entries from $Z_{19}I$} $\cup$

$$\left\{ \sum_{i=0}^{\infty} a_i x^i \middle| a_i \in N(Z_{19}) \right\}.$$

V is a strong neutrosophic 7-linear algebra over the neutrosophic field $Z_{19}I$.

It is important and interesting to record that every strong neutrosophic n-linear algebra is a strong neutrosophic n-vector space but in general every strong neutrosophic n-vector space need not be a strong neutrosophic n-linear algebra. We leave the proof of this to the reader; however we give an example of a neutrosophic strong n-vector space which is not a neutrosophic strong n-linear algebra.

***Example 3.2.10:*** Let $V = V_1 \cup V_2 \cup V_3 \cup V_4 \cup V_5 =$

$$\left\{ \begin{pmatrix} a_1 & a_2 & a_3 & a_4 \\ a_5 & a_6 & a_7 & a_8 \end{pmatrix} \middle| a_i \in Z_{13}I; 1 \le i \le 8 \right\} \cup$$

$$\left\{ \begin{pmatrix} a_1 \\ a_2 \\ a_3 \\ a_4 \\ a_5 \end{pmatrix} \middle| a_i \in Z_{13}I; 1 \leq i \leq 5 \right\} \cup$$

$$\left\{ \begin{pmatrix} 0 & 0 & 0 & a_1 \\ 0 & 0 & a_2 & 0 \\ 0 & a_3 & 0 & 0 \\ a_4 & 0 & 0 & 0 \end{pmatrix} \middle| a_i \in Z_{13}I; 1 \leq i \leq 4 \right\} \cup$$

$$\left\{ \begin{pmatrix} a & 0 & b & 0 \\ a & c & 0 & d \\ e & 0 & f & 0 \\ a & h & a & i \end{pmatrix} \middle| a,b,e,d,f,h,i \in Z_{13}I \right\} \cup$$

$$\left\{ \sum_{i=0}^{19} a_i x^i \middle| a_i \in Z_{13}I \right\}$$

be a strong neutrosophic 5-vector space over the neutrosophic field $Z_{13}I$. Clearly V is not a strong neutrosophic 5-linear algebra over $Z_{13}I$ as none of $V_i$ is closed under multiplication; $1 \leq i \leq 5$.

Hence in general a strong neutrosophic n-vector space need not be a strong neutrosophic n-linear algebra.

Now we proceed onto define the notion of strong neutrosophic n-linear subalgebra of a strong neutrosophic n-linear algebra.

**DEFINITION 3.2.4:** *Let $V = V_1 \cup V_2 \cup ... \cup V_n$ be a strong neutrosophic n-linear algebra over the neutrosophic field F. Suppose $W = W_1 \cup W_2 \cup ... \cup W_n \subseteq V_1 \cup V_2 \cup ... \cup V_n$ ($W_i \not\subseteq V_i$; $W_i \neq \phi$ and $W_i \neq V_i$ for each i, $1 \leq i \leq n$) is such that W is* *strong neutrosophic n-linear algebra over the field F, then we call W to be a strong neutrosophic n-linear subalgebra of V.*

It is interesting and important to note that even if one of the $W_i$ is $\{0\}$ or $V_i$ then W is not a strong neutrosophic n-linear subalgebra of V over F.

We will illustrate this by some simple examples.

*Example 3.2.11:* Let $V = V_1 \cup V_2 \cup V_3 \cup V_4 \cup V_5 =$

$$\left\{ \begin{pmatrix} a & b \\ c & d \end{pmatrix} \middle| a,b,c,d \in N(Q) \right\} \cup$$

$$\{(a, b, c, d, e, f) \mid a, b, c, d, e, f \in N(Q)\} \cup$$

$$\left\{ \sum_{i=0}^{\infty} a_i x^i \middle| a_i \in N(Q) \right\} \cup \left\{ \begin{pmatrix} a & 0 & 0 \\ 0 & b & 0 \\ 0 & 0 & d \end{pmatrix} \middle| a,b,d \in N(Q) \right\} \cup$$

$$\left\{ \begin{pmatrix} a & 0 & 0 & 0 \\ b & c & 0 & 0 \\ d & e & f & 0 \\ g & h & j & k \end{pmatrix} \middle| a,b,c,d,e,f,g,h,j,k \in N(Q) \right\}$$

be a strong neutrosophic n-linear algebra over QI. Let $W = W_1 \cup W_2 \cup W_3 \cup W_4 \cup W_5 =$

$$\left\{ \begin{pmatrix} a & a \\ a & a \end{pmatrix} \middle| a \in QI \right\} \cup \{(a\ a\ a\ a\ a\ a) \mid a \in N(Q)\} \cup$$

$$\left\{ \sum_{i=0}^{\infty} a_i x^{2i} \middle| a_i \in N(Q) \right\} \cup$$

$$\left\{ \begin{pmatrix} a & 0 & 0 \\ 0 & a & 0 \\ 0 & 0 & a \end{pmatrix} \middle| a \in N(Q) \right\} \cup$$

$$\left\{ \begin{pmatrix} a & 0 & 0 & 0 \\ b & c & 0 & 0 \\ d & e & g & 0 \\ p & q & r & s \end{pmatrix} \middle| a,b,c,d,e,g,p,q,r,s \in QI \right\}$$

$$\subseteq V_1 \cup V_2 \cup V_3 \cup V_4 \cup V_5.$$

It is easily verified that W is a strong neutrosophic n-linear subalgebra of V over the field QI.

**Example 3.2.12:** Let $V = V_1 \cup V_2 \cup V_3 =$

$$\left\{ \begin{pmatrix} a_1 & a_2 & a_3 \\ a_4 & a_5 & a_6 \\ a_7 & a_8 & a_9 \end{pmatrix} \middle| a_i \in Z_3I; 1 \leq i \leq 9 \right\} \cup$$

{All $4 \times 4$ neutrosophic matrices with entries from $N(Z_3)$} $\cup$ {All $9 \times 9$ upper triangular matrices with entries from $N(Z_3)$}; V is a neutrosophic strong trilinear algebra over the neutrosophic field $Z_3I$. Take $W = W_1 \cup W_2 \cup W_3 =$

$$\left\{ \begin{pmatrix} a & a & a \\ a & a & a \\ a & a & a \end{pmatrix} \middle| a \in Z_3I \right\} \cup$$

{All $4 \times 4$ neutrosophic matrices with entries from $Z_3I$} $\cup$ {All $9 \times 9$ upper triangular matrices with entries from $Z_3I$} $\subseteq V_1 \cup V_2 \cup V_3$. W is a strong neutrosophic 3-linear subalgebra of V of V over the neutrosophic field $Z_3I$.

We call a strong neutrosophic n-linear algebra to be simple if $V = V_1 \cup V_2 \cup \ldots \cup V_n$ has no proper strong neutrosophic n-linear subalgebra.

We will illustrate this by some examples.

*Example 3.2.13:* Let $V = V_1 \cup V_2 \cup V_3 \cup V_4 \cup V_5 =$

$$\left\{ \begin{pmatrix} a & a \\ a & a \end{pmatrix} \middle| a \in Z_5 I \right\} \cup \{(a, a, a, a, a, a) \mid a \in Z_5 I\} \cup$$

$$\left\{ \begin{pmatrix} a & 0 & 0 & 0 \\ 0 & a & 0 & 0 \\ 0 & 0 & a & 0 \\ 0 & 0 & 0 & a \end{pmatrix} \middle| a \in Z_5 I \right\} \cup \left\{ \begin{pmatrix} a & 0 & 0 \\ a & a & 0 \\ a & a & a \end{pmatrix} \middle| a \in Z_5 I \right\} \cup$$

$$\left\{ \begin{pmatrix} a & 0 & 0 & 0 & 0 \\ 0 & a & 0 & 0 & 0 \\ 0 & 0 & a & 0 & 0 \\ 0 & 0 & 0 & a & 0 \\ 0 & 0 & 0 & 0 & a \end{pmatrix} \middle| a \in Z_5 I \right\}$$

be a strong neutrosophic 5-linear algebra over the neutrosophic field $Z_5 I$. Clearly V has no proper strong neutrosophic 5-linear subalgebra.

Thus V is a simple strong neutrosophic 5-linear algebra over $Z_5 I$.

*Example 3.2.14:* Let $V = V_1 \cup V_2 \cup V_3 \cup V_4 \cup V_5 \cup V_6 =$

$$\left\{ \begin{pmatrix} a & a & a \\ a & a & a \\ a & a & a \end{pmatrix} \middle| a \in Z_{11} I \right\} \cup \left\{ \begin{pmatrix} a & 0 \\ a & 0 \end{pmatrix} \middle| a \in Z_{11} I \right\} \cup$$

{(a 0 a 0 a 0 a) | a ∈ $Z_{11}I$} ∪ {9 × 9 diagonal matrices where all the diagonal entries are equal to (say) a; a ∈ $Z_{11}I$} ∪

$$\left\{ \begin{pmatrix} a & 0 & 0 & 0 \\ a & a & 0 & 0 \\ a & a & a & 0 \\ a & a & a & a \end{pmatrix} \middle| a \in Z_{11}I \right\} \cup$$

$$\left\{ \begin{pmatrix} a & a & a & a & a & a \\ 0 & a & a & a & a & a \\ 0 & 0 & a & a & a & a \\ 0 & 0 & 0 & a & a & a \\ 0 & 0 & 0 & 0 & a & a \\ 0 & 0 & 0 & 0 & 0 & a \end{pmatrix} \middle| a \in Z_{11}I \right\}$$

be a strong neutrosophic 6-linear algebra over the neutrosophic field $Z_{11}I$. It is easy to vertify that V has no proper strong neutrosophic 6-linear subalgebras; hence V is a simple strong neutrosophic linear algebra.

Now we proceed onto define the notion of pseudo strong neutrosophic n-linear subalgebra.

**DEFINITION 3.2.5:** *Let $V = V_1 \cup V_2 \cup V_3 \cup ... \cup V_n$ be a strong neutrosophic n-linear algebra over the neutrosophic field F.*

*Let $W = W_1 \cup W_2 \cup ... \cup W_n \subseteq V_1 \cup V_2 \cup ... \cup V_n$, where some $W_i$ ($W_i \neq (0)$ and $W_i \neq V_i$) contained in $V_i$ are just strong neutrosophic vector space over the neutrosophic field F; i = 1, 2, ..., n.*

*We define $W = W_1 \cup W_2 \cup ... \cup W_n \subseteq V_1 \cup V_2 \cup ... \cup V_n$ to be a pseudo strong neutrosophic linear subalgebra of V over F if some $W_i$'s are strong neutrosophic vector spaces and some $W_j$'s are strong neutrosophic linear algebras over the neutrosophic field F. $1 \leq i, j \leq n$ ($i \neq j$).*

We will illustrate this situation by some examples.

*Example 3.2.15:* Let $V = V_1 \cup V_2 \cup V_3 \cup V_4 =$

$$\left\{ \begin{pmatrix} a & b \\ c & d \end{pmatrix} \middle| a,b,c,d \in Z_{23}I \right\} \cup$$

$$\{(a_1, a_2, a_3, a_4, a_5, a_6) \mid a_i \in Z_{23}I;\ 1 \leq i \leq 6\} \cup$$

$$\left\{ \begin{pmatrix} a & b & c \\ d & e & f \\ g & h & i \end{pmatrix} \middle| a,b,c,d,e,f,g,h,i \in Z_{23}I \right\} \cup$$

$$\left\{ \sum_{i=0}^{\infty} a_i x^i \middle| a_i \in Z_{23}I;\ i = 0,1,2,...,\infty \right\}$$

be a strong neutrosophic 4-linear algebra over the neutrosophic field $Z_{23}I$. Choose $W = W_1 \cup W_2 \cup W_3 \cup W_4 =$

$$\left\{ \begin{pmatrix} 0 & a \\ b & 0 \end{pmatrix} \middle| a,b \in Z_{23}I \right\} \cup \{(a\ a\ a\ a\ a\ a) \mid a \in Z_{23}I\} \cup$$

$$\left\{ \begin{pmatrix} 0 & 0 & a \\ 0 & b & 0 \\ c & 0 & 0 \end{pmatrix} \middle| a,b,c,d \in Z_{23}I \right\} \cup$$

$$\left\{ \sum_{i=0}^{\infty} a_i x^{2i} \middle| a_i \in Z_{23}I;\ i = 0,1,2,...,\infty \right\}$$

$\subseteq V_1 \cup V_2 \cup V_3 \cup V_4$. It is easily verified that $W = W_1 \cup W_2 \cup W_3 \cup W_4$ is a pseudo strong neutrosophic linear subalgebra of V over the neutrosophic field $F = Z_{23}I$.

**Example 3.2.16:** Let $V = V_1 \cup V_2 \cup V_3 \cup V_4 \cup V_5 \cup V_6 =$

$$\left\{ \begin{pmatrix} a & b & c \\ d & e & f \\ g & h & i \end{pmatrix} \middle| a,b,c,d,e,f,g,h,i \in Z_{13}I \right\} \cup$$

$$\{(a_1, a_2, a_3, a_4, a_5, a_6, a_7, a_8) \mid a_i \in Z_{13}I;\ 1 \leq i \leq 8\} \cup$$

$$\left\{ \begin{pmatrix} a_1 & 0 & 0 & 0 \\ a_2 & a_3 & 0 & 0 \\ a_4 & a_5 & a_6 & 0 \\ a_7 & a_8 & a_9 & a_{10} \end{pmatrix} \middle| a_i \in Z_{13}I; 1 \leq i \leq 10 \right\} \cup$$

$$\left\{ \begin{pmatrix} a_1 & a_2 \\ a_3 & a_4 \end{pmatrix} \middle| a_i \in Z_{13}I; 1 \leq i \leq 4 \right\} \cup$$

$$\left\{ \sum_{i=0}^{\infty} a_i x^i \middle| a_i \in Z_{13}I; 1 \leq i \leq \infty; \right.$$

that is all polynomials in the variable x with coefficients from $Z_{13}I\} \cup$

$$\left\{ \begin{pmatrix} a_1 & 0 & 0 & 0 & 0 \\ 0 & a_2 & 0 & 0 & 0 \\ 0 & 0 & a_3 & 0 & 0 \\ 0 & 0 & 0 & a_4 & 0 \\ 0 & 0 & 0 & 0 & a_5 \end{pmatrix} \middle| a_i \in Z_{13}I; 1 \leq i \leq 5 \right\}$$

be a strong neutrosophic 6-linear algebra over the neutrosophic field $Z_{13}I$.

Choose $W = W_1 \cup W_2 \cup W_3 \cup W_4 \cup W_5 \cup W_6 =$

$$\left\{ \begin{pmatrix} 0 & 0 & a \\ 0 & b & 0 \\ d & 0 & 0 \end{pmatrix} \middle| a, b, d \in Z_{13}I \right\} \cup$$

$\{(a_1, 0, a_2, 0, a_3, 0, a_4, 0) \mid a_i \in Z_{13}I; 1 \leq i \leq 4\} \cup$

$$\left\{ \begin{pmatrix} a & 0 & 0 & 0 \\ a & a & 0 & 0 \\ a & a & a & 0 \\ a & a & a & a \end{pmatrix} \middle| a \in Z_{13}I \right\} \cup$$

$$\left\{ \begin{pmatrix} 0 & a \\ b & 0 \end{pmatrix} \middle| a, b \in Z_{13}I \right\} \cup$$

$$\left\{ \sum_{i=0}^{27} a_i x^i \middle| a_i \in Z_{13}I; 0 \leq i \leq 27 \right. ;$$

all polynomials in the variable x with coefficients from $Z_{13}I\} \cup$

$$\left\{ \begin{pmatrix} a & 0 & 0 & 0 & 0 \\ 0 & a & 0 & 0 & 0 \\ 0 & 0 & a & 0 & 0 \\ 0 & 0 & 0 & a & 0 \\ 0 & 0 & 0 & 0 & a \end{pmatrix} \middle| a \in Z_{13}I \right\}$$

$\subseteq V_1 \cup V_2 \cup V_3 \cup V_4 \cup V_5 \cup V_6$.

It is easily verified that W is a pseudo strong neutrosophic 6-linear subalgebra of V over the field $Z_{13}I$.

Now we proceed onto give an example and then define the notion of strong pseudo neutrosophic n-vector space of a strong neutrosophic n-linear algebra over the neutrosophic field F.

**Example 3.2.17:** Let $V = V_1 \cup V_2 \cup V_3 \cup V_4 \cup V_5 =$

$$\left\{ \begin{pmatrix} a & b \\ c & d \end{pmatrix} \middle| a,b,c,d \in N(Q) \right\} \cup$$

$\left\{ \sum_{i=0}^{\infty} a_i x^i \right.$ ; all polynomials in the variable x with coefficients from the neutrosophic field QI; $0 \le a_i \le \infty \} \cup$

$$\left\{ \begin{pmatrix} a & a & a \\ a & a & a \\ a & a & a \end{pmatrix} \middle| a \in N(Q) \right\} \cup$$

$$\left\{ \begin{pmatrix} a_1 & a_2 & a_3 & a_4 \\ a_5 & a_6 & a_7 & a_8 \\ a_9 & a_{10} & a_{11} & a_{12} \\ a_{13} & a_{14} & a_{15} & a_{16} \end{pmatrix} \middle| a_i \in QI; 1 \le i \le 16 \right\} \cup$$

{all $8 \times 8$ neutrosophic matrices with entries from QI} be a strong neutrosophic 5-linear algebra over the neutrosophic field QI.

Let $W = W_1 \cup W_2 \cup W_3 \cup W_4 \cup W_5 =$

$$\left\{ \begin{pmatrix} 0 & a \\ b & 0 \end{pmatrix} \middle| a,b \in N(Q) \right\} \cup$$

$\left\{ \sum_{i=0}^{10} a_i x^i \right.$ ; all polynomials in the variable x with coefficients from the neutrosophic field QI of degree less than equal to 10; $a_i \in QI$; $i = 0, 1, 2, ..., 10\} \cup$

$$\left\{\begin{pmatrix} 0 & 0 & a \\ 0 & a & 0 \\ a & 0 & 0 \end{pmatrix} \middle| a \in QI\right\} \cup$$

$$\left\{\begin{pmatrix} 0 & 0 & 0 & a \\ 0 & 0 & a & 0 \\ 0 & a & 0 & 0 \\ a & 0 & 0 & 0 \end{pmatrix} \middle| a \in QI\right\} \cup$$

$$\left\{\begin{pmatrix} 0 & 0 & 0 & 0 & 0 & 0 & 0 & a_1 \\ 0 & 0 & 0 & 0 & 0 & 0 & a_2 & 0 \\ 0 & 0 & 0 & 0 & 0 & a_3 & 0 & 0 \\ 0 & 0 & 0 & 0 & a_4 & 0 & 0 & 0 \\ 0 & 0 & 0 & a_5 & 0 & 0 & 0 & 0 \\ 0 & 0 & a_6 & 0 & 0 & 0 & 0 & 0 \\ 0 & a_7 & 0 & 0 & 0 & 0 & 0 & 0 \\ a_8 & 0 & 0 & 0 & 0 & 0 & 0 & 0 \end{pmatrix} \middle| a_i \in QI; 1 \le i \le 8\right\}$$

$\subseteq V_1 \cup V_2 \cup V_3 \cup V_4 \cup V_5$. W is a strong neutrosophic 5-vector space over the neutrosophic field QI. We call this strong neutrosophic 5-vector space as strong pseudo neutrosophic 5-vector space of V over the neutrosophic field QI.

We now give the formal definition of this new notion.

**DEFINITION 3.2.6:** *Let $V = V_1 \cup V_2 \cup V_3 \cup \ldots \cup V_n$ be a strong neutrosophic n-linear algebra over the neutrosophic field F. Let $W = W_1 \cup W_2 \cup \ldots \cup W_n \subseteq V_1 \cup V_2 \cup \ldots \cup V_n$ be such that each $W_i \subseteq V_i$ is different from (0) and $V_i$ for i = 1, 2, ..., n and each $W_i$ is only a strong neutrosophic vector space over the field F and is not a strong neutrosophic linear algebra over F.*

*We define $W = W_1 \cup W_2 \cup \ldots \cup W_n$ to be a strong pseudo neutrosophic n-vector space of V over the field F.*

We will give an example of this concept.

***Example 3.2.18:*** Let $V = V_1 \cup V_2 \cup V_3 \cup V_4 =$

$$\left\{ \begin{pmatrix} a & b & c \\ d & e & f \\ g & h & i \end{pmatrix} \middle| a,b,c,d,e,f,g,h,i \in Z_{17}I \right\} \cup$$

$$\left\{ \begin{pmatrix} a & b \\ c & d \end{pmatrix} \middle| a,b,c,d \in Z_{17}I \right\} \cup$$

$\left\{ \sum_{i=0}^{\infty} a_i x^i \right.$; all polynomials in the variable x with coefficients from the neutrosophic field $Z_{17}I$; $a_i \in Z_{17}I$; $0 \le i \le \infty\} \cup$ {all 6 × 6 matrices with entries from $Z_{17}I$} be a strong neutrosophic 4-linear algebra over the neutrosophic field $Z_{17}I$. Take $W = W_1 \cup W_2 \cup W_3 \cup W_4 =$

$$\left\{ \begin{pmatrix} 0 & 0 & a \\ 0 & b & 0 \\ c & 0 & 0 \end{pmatrix} \middle| a,b,c \in Z_{17}I \right\} \cup \left\{ \begin{pmatrix} 0 & a \\ b & 0 \end{pmatrix} \middle| a,b \in Z_{17}I \right\} \cup$$

$\left\{ \sum_{i=0}^{5} a_i x^i \right.$; all polynomials in variable x with coefficients from $Z_{17}I$ of degree less than or equal to five $a_i \in Z_{17}I$; $0 \le i \le 5\} \cup$

$$\left\{ \begin{pmatrix} 0 & 0 & 0 & 0 & 0 & a_1 \\ 0 & 0 & 0 & 0 & a_2 & 0 \\ 0 & 0 & 0 & a_3 & 0 & 0 \\ 0 & 0 & a_4 & 0 & 0 & 0 \\ 0 & a_5 & 0 & 0 & 0 & 0 \\ a_6 & 0 & 0 & 0 & 0 & 0 \end{pmatrix} \middle| a_i \in Z_{17}I; 1 \le i \le 6 \right\}$$

$\subseteq V_1 \cup V_2 \cup V_3 \cup V_4$. It is easy to verify each $W_i \subseteq V_i$ is only a strong neutrosophic vector space over $Z_{17}I$. $i = 1, 2, 3, 4$. So $W = W_1 \cup W_2 \cup W_3 \cup W_4 \subseteq V_1 \cup V_2 \cup V_3 \cup V_4$ is only a pseudo strong neutrosophic 4-vector space of V over the field $Z_{17}I$.

**DEFINITION 3.2.7:** *Let $V = V_1 \cup V_2 \cup ... \cup V_n$ be a strong neutrosophic n-linear algebra over the neutrosophic field F. Let $W = W_1 \cup W_2 \cup ... \cup W_n \subseteq V_1 \cup V_2 \cup ... \cup V_n$ be a neutrosophic n-linear algebra over the neutrosophic field $K \subseteq F$ (K is only a proper subfield of F, $K \ne F$). We define W to be a pseudo strong neutrosophic n-linear subalgebra of V over the real field $K \subseteq F$.*

We will illustrate this situation by some examples.

*Example 3.2.19:* Let $V = V_1 \cup V_2 \cup V_3 \cup V_4 \cup V_5 =$

$$\left\{ \begin{pmatrix} a & b \\ c & d \end{pmatrix} \middle| a,b,c,d \in N(Z_{23}) \right\} \cup$$

$\{(a_1, a_2, a_3, a_4, a_5, a_6, a_7) \mid a_i \in N(Z_{23})\ 1 \le i \le 7\} \cup$

$\left\{ \sum_{i=0}^{\infty} a_i x^i \right.$ ; all polynomials in the variable x with coefficients from the neutrosophic field $N(Z_{23})\} \cup$

$$\left\{ \begin{pmatrix} a_1 & a_2 & a_3 \\ a_4 & a_5 & a_6 \\ a_7 & a_8 & a_9 \end{pmatrix} \middle| a_i \in N(Z_{23}); 1 \le i \le 9 \right\} \cup$$

{9 × 9 matrices with entries from $N(Z_{23})$} be a strong neutrosophic 5-linear algebra over the neutrosophic field $N(Z_{23})$.
Take $W = W_1 \cup W_2 \cup W_3 \cup W_4 \cup W_5 =$

$$\left\{ \begin{pmatrix} a & b \\ c & d \end{pmatrix} \middle| a,b,c,d \in Z_{23}I \right\} \cup \{(a\ a\ a\ a\ a\ a\ a) \mid a \in Z_{23}I\} \cup$$

$$\left\{ \sum_{i=0}^{\infty} a_i x^i \middle| a_i \in Z_{23}I; 0 \le i \le \infty \right\} \cup$$

$$\left\{ \begin{pmatrix} a_1 & 0 & 0 \\ a_2 & a_3 & 0 \\ a_4 & a_5 & 0 \end{pmatrix} \middle| a_i \in N(Z_{23}) \right\} \cup$$

{all $9 \times 9$ lower triangular matrices with entries from $N(Z_{23}) \subseteq V_1 \cup V_2 \cup V_3 \cup V_4 \cup V_5$ is a neutrosophic 5-linear algebra over the field $Z_{23} \subseteq N(Z_{23})$. Thus W is a pseudo strong neutrosophic 5-linear subalgebra of V over the subfield $Z_{23} \subseteq N(Z_{23})$.

*Example 3.2.20:* Let $V = V_1 \cup V_2 \cup V_3$ = {All $5 \times 5$ neutrosophic matrices with entries from $N(R)$} $\cup$ {All polynomial $\sum_{i=0}^{\infty} a_i x^i$ with coefficients from $N(R)$ in the variable x} $\cup$ {$(x_1, x_2, x_3, x_4, x_5, x_6, x_7, x_8, x_9, x_{10}) \mid x_i \in N(R); 1 \le i \le 10$} be a strong neutrosophic 3-linear algebra over the neutrosophic field $N(R)$. Take $W = W_1 \cup W_2 \cup W_3$ = {All $5 \times 5$ diagonal matrices with entries from $N(R)$} $\cup$

$$\left\{ \sum_{i=0}^{\infty} a_i x^{2i} \middle| a_i \in N(R) \right\} \cup$$

{$(a\ a\ a\ a\ a\ a\ a\ a\ a\ a) \mid a \in N(R)$} $\subseteq V_1 \cup V_2 \cup V_3$; W is a neutrosophic 3-linear algebra over the field $R \subseteq N(R)$. Thus W is a pseudo strong neutrosophic 3-linear subalgebra of V over the field $R \subseteq N(R)$.

**DEFINITION 3.2.8:** *Let $V = V_1 \cup V_2 \cup ... \cup V_n$ be a strong neutrosophic n-linear algebra over the neutrosophic field K. Let $W = W_1 \cup W_2 \cup ... \cup W_n \subseteq V_1 \cup V_2 \cup ... \cup V_n$ be a n-vector space over a real field F, $F \subseteq K$; we call $W = W_1 \cup W_2 \cup ... \cup W_n \subseteq V_1 \cup V_2 \cup ... \cup V_n$ to be a pseudo n-vector subspace of V over the real subfield F of K.*

We will illustrate this situation by some examples.

*Example 3.2.21:* Let $V = V_1 \cup V_2 \cup V_3 \cup V_4 \cup V_5 \cup V_6 \cup V_7 =$

$$\left\{ \begin{pmatrix} a & b \\ c & d \end{pmatrix} \middle| a,b,c,d \in N(Z_2) \right\} \cup$$

$\left\{ \sum_{i=0}^{\infty} a_i x^i \right.$; all polynomials in the variable x with coefficients from $N(Z_2)$; $a_i \in N(Z_2)$ and $i = 0, 1, 2, ..., \infty \}$ ∪ {All 5 × 5 matrices with entries from $N(Z_2)$} ∪ {All 7 × 7 matrices with entries from $N(Z_2)$} ∪ {All 4 × 4 matrices with entries from $N(Z_2)$} ∪ {All 6 × 6 matrices with entries from $N(Z_2)$} ∪

$$\left\{ \begin{pmatrix} a & b & c \\ d & e & f \\ g & h & i \end{pmatrix} \middle| a,b,c,d,e,f,g,h,i \in N(Z_2) \right\}$$

be a strong neutrosophic 7-linear algebra over the neutrosophic field $N(Z_2)$. Let $W = W_1 \cup W_2 \cup W_3 \cup W_4 \cup W_5 \cup W_6 \cup W_7$

$$= \left\{ \begin{pmatrix} 0 & a \\ b & 0 \end{pmatrix} \middle| a,b \in Z_2 \right\} \cup$$

$\left\{ \sum_{i=0}^{6} a_i x^i \right.$; all polynomials in the variable x with coefficients from the real field $Z_2$ of degree less than or equal to 6} ∪

$$\left\{ \begin{pmatrix} 0 & 0 & 0 & 0 & a \\ 0 & 0 & 0 & a & 0 \\ 0 & 0 & a & 0 & 0 \\ 0 & a & 0 & 0 & 0 \\ a & 0 & 0 & 0 & 0 \end{pmatrix} \middle| a \in Z_2 \right\} \cup$$

$$\left\{ \begin{pmatrix} 0 & 0 & 0 & 0 & 0 & 0 & a_1 \\ 0 & 0 & 0 & 0 & 0 & a_2 & 0 \\ 0 & 0 & 0 & 0 & a_3 & 0 & 0 \\ 0 & 0 & 0 & a_4 & 0 & 0 & 0 \\ 0 & 0 & a_5 & 0 & 0 & 0 & 0 \\ 0 & a_6 & 0 & 0 & 0 & 0 & 0 \\ a_7 & 0 & 0 & 0 & 0 & 0 & 0 \end{pmatrix} \middle| a_i \in Z_2; 1 \leq i \leq 7 \right\} \cup$$

$$\left\{ \begin{pmatrix} 0 & 0 & 0 & a \\ 0 & 0 & b & 0 \\ 0 & c & 0 & 0 \\ d & 0 & 0 & 0 \end{pmatrix} \middle| a,b,c,d \in Z_2 \right\} \cup$$

$$\left\{ \begin{pmatrix} 0 & 0 & 0 & 0 & 0 & 0 & a_1 \\ 0 & 0 & 0 & 0 & 0 & a_2 & 0 \\ 0 & 0 & 0 & 0 & a_3 & 0 & 0 \\ 0 & 0 & 0 & a_4 & 0 & 0 & 0 \\ 0 & 0 & a_5 & 0 & 0 & 0 & 0 \\ 0 & a_6 & 0 & 0 & 0 & 0 & 0 \\ a_7 & 0 & 0 & 0 & 0 & 0 & 0 \end{pmatrix} \middle| a_i \in Z_2; 1 \leq i \leq 7 \right\} \cup$$

$$\left\{ \begin{pmatrix} 0 & 0 & a \\ 0 & b & 0 \\ d & 0 & 0 \end{pmatrix} \middle| a,b,d \in Z_2 \right\}$$

$\subseteq V_1 \cup V_2 \cup V_3 \cup V_4 \cup V_5 \cup V_6 \cup V_7$. W is only a 7-vector space over the field $Z_2 \subseteq N(Z_2)$. Thus W is a pseudo 7-vector subspace of V over the real field $Z_2 \subseteq N(Z_2)$.

*Example 3.2.22:* Let $V = V_1 \cup V_2 \cup V_3 \cup V_4 \cup V_5 =$

$$\left\{ \begin{pmatrix} a & b \\ c & d \end{pmatrix} \middle| a,b,c,d \in N(Z_7) \right\} \cup$$

$\{(a_1, a_2, a_3, a_4, a_5) \mid a_i \in N(Z_7); 1 \leq i \leq 5\} \cup$ {All 4×4 matrices with entries from $N(Z_7)$} $\cup$ {All 6 × 6 matrices with entries from $N(Z_7)$} $\cup$ {All 5 × 5 matrices with entries from $N(Z_7)$} be a strong neutrosophic 5-linear algebra over the neutrosophic field $N(Z_7)$.

Take $W = W_1 \cup W_2 \cup W_3 \cup W_4 \cup W_5 =$

$$\left\{ \begin{pmatrix} 0 & a \\ b & 0 \end{pmatrix} \middle| a,b \in Z_7 \right\} \cup$$

$\{(a_1 \; 0 \; a_2 \; 0 \; a_3) \mid a_i \in Z_7; 1 \leq i \leq 3\} \cup$

$$\left\{ \begin{pmatrix} 0 & 0 & 0 & a \\ 0 & 0 & b & 0 \\ 0 & c & 0 & 0 \\ d & 0 & 0 & 0 \end{pmatrix} \middle| a,b,c,d \in Z_7 \right\} \cup$$

$$\left\{ \begin{pmatrix} 0 & 0 & 0 & 0 & 0 & a \\ 0 & 0 & 0 & 0 & b & 0 \\ 0 & 0 & 0 & c & 0 & 0 \\ 0 & 0 & d & 0 & 0 & 0 \\ 0 & e & 0 & 0 & 0 & 0 \\ p & 0 & 0 & 0 & 0 & 0 \end{pmatrix} \middle| a,b,c,d,e,p \in Z_7 \right\} \cup$$

$$\left\{ \begin{pmatrix} 0 & 0 & 0 & 0 & a \\ 0 & 0 & 0 & 0 & 0 \\ 0 & 0 & b & 0 & 0 \\ 0 & 0 & 0 & 0 & 0 \\ c & 0 & 0 & 0 & 0 \end{pmatrix} \middle| a, b, c \in Z_7 \right\}$$

$\subseteq V_1 \cup V_2 \cup V_3 \cup V_4 \cup V_5$. Clearly W is only a 7-vector space over the real field $Z_7$. Thus W is a pseudo 7-vector subspace of V over the real field $Z_7$ of $N(Z_7)$.

**DEFINITION 3.2.9:** *Let $V = V_1 \cup V_2 \cup ... \cup V_n$ be a strong neutrosophic n-linear algebra over the neutrosophic field F. Let $W = W_1 \cup W_2 \cup ... \cup W_n \subseteq V_1 \cup V_2 \cup ... \cup V_n$ be a n-linear algebra over the real field K; $K \not\subseteq F$. We define W to be a pseudo n-linear subalgebra of V over the real field $K \subseteq F$.*

We will illustrate this situation by some examples.

*Example 3.2.23:* Let $V = V_1 \cup V_2 \cup V_3 \cup V_4 =$

$$\left\{ \begin{pmatrix} a & b \\ c & d \end{pmatrix} \middle| a, b, c, d \in N(Z_{11}) \right\} \cup$$

$\{a_1, a_2, a_3, a_4, a_5, a_6) \mid a_i \in N(Z_{11}); 1 \leq i \leq 6\} \cup$

$$\left\{ \begin{pmatrix} a & b & c \\ d & e & f \\ g & h & u \end{pmatrix} \middle| a, b, c, d, e, f, g, h, i \in N(Z_{11}) \right\} \cup$$

$$\left\{ \sum_{i=0}^{\infty} a_i x^i \middle| a_i \in N(Z_{11}) \right\}$$

be a strong neutrosophic 4-linear algebra over the neutrosophic field $N(Z_{11})$.

Let $W = W_1 \cup W_2 \cup W_3 \cup W_4 =$

$$\left\{ \begin{pmatrix} a & b \\ c & d \end{pmatrix} \middle| a,b,c,d \in Z_{11} \right\} \cup \{(a, a, a, b, b, b) \mid a, b \in Z_{11}\} \cup$$

$$\left\{ \begin{pmatrix} a & a & a \\ b & b & b \\ c & c & c \end{pmatrix} \middle| a,b,c \in Z_{11} \right\} \cup$$

$$\left\{ \sum_{i=0}^{\infty} a_i x^{2i} \middle| a_i \in Z_{11} \right\}$$

$\subseteq V_1 \cup V_2 \cup V_3 \cup V_4$.

It can be easily verified that W is a pseudo 4-linear subalgebra of V over the field $Z_{11} \subseteq N(Z_{11})$.

**Example 3.2.24:** Let $V = V_1 \cup V_2 \cup V_3 \cup V_4 \cup V_5 \cup V_6 =$

$$\left\{ \begin{pmatrix} a & b \\ c & d \end{pmatrix} \middle| a,b,c,d \in N(Q) \right\} \cup$$

$\{(a, b, c) \mid c, a, b \in N(Q)\} \cup$

$\left\{ \sum_{i=0}^{\infty} a_i x^i \middle| a_i \in N(Q) \right.$ ; collection of all polynomials in the variable

x with coefficients from the neutrosophic field N(Q)$\} \cup$

$$\left\{ \begin{pmatrix} a & b & 0 & 0 \\ c & d & 0 & 0 \\ 0 & 0 & e & f \\ 0 & 0 & g & h \end{pmatrix} \middle| a,b,c,d,e,f,g,h \in N(Q) \right\} \cup$$

$$\left\{ \begin{pmatrix} a & 0 & 0 & 0 & 0 \\ b & d & 0 & 0 & 0 \\ e & f & g & 0 & 0 \\ h & i & j & k & 0 \\ l & m & n & p & q \end{pmatrix} \middle| a,b,d,e,f,g,h,i,j,k,l,m,n,p,q \in N(Q) \right\}$$

$$\cup \left\{ \begin{pmatrix} a & 0 & 0 & 0 & 0 & 0 \\ 0 & b & 0 & 0 & 0 & 0 \\ 0 & 0 & c & 0 & 0 & 0 \\ 0 & 0 & 0 & d & 0 & 0 \\ 0 & 0 & 0 & 0 & e & 0 \\ 0 & 0 & 0 & 0 & 0 & g \end{pmatrix} \middle| a,b,c,d,e,g \in N(Q) \right\}$$

be a strong neutrosophic 6-linear algebra over the neutrosophic field N (Q). Consider $W = W_1 \cup W_2 \cup W_3 \cup W_4 \cup W_5 \cup W_6 =$

$$\left\{ \begin{pmatrix} a & 0 \\ 0 & b \end{pmatrix} \middle| a,b \in Q \right\} \cup \{(a, b, 0) \mid a, b \in Q\} \cup$$

$$\left\{ \sum_{i=0}^{\infty} a_i x^i \middle| a_i \in Q; 0 \leq i \leq \infty \right\} \cup \left\{ \begin{pmatrix} a & b & 0 & 0 \\ c & d & 0 & 0 \\ 0 & 0 & 0 & 0 \\ 0 & 0 & 0 & 0 \end{pmatrix} \middle| a,b,c,d \in Q \right\} \cup$$

$$\left\{ \begin{pmatrix} a & 0 & 0 & 0 & 0 \\ 0 & b & 0 & 0 & 0 \\ 0 & 0 & c & 0 & 0 \\ 0 & 0 & 0 & d & 0 \\ 0 & 0 & 0 & 0 & e \end{pmatrix} \middle| a,b,c,d,e \in Q \right\} \cup$$

$$\left\{ \begin{pmatrix} a & 0 & 0 & 0 & 0 & 0 \\ 0 & b & 0 & 0 & 0 & 0 \\ 0 & 0 & c & 0 & 0 & 0 \\ 0 & 0 & 0 & d & 0 & 0 \\ 0 & 0 & 0 & 0 & e & 0 \\ 0 & 0 & 0 & 0 & 0 & f \end{pmatrix} \middle| a,b,c,d,e,f,g \text{ are in } Q \right\}$$

$\subseteq V_1 \cup V_2 \cup V_3 \cup V_4 \cup V_5$.

It is easy to prove, W is a 5-linear algebra over Q, the real field of rationals. Hence W is a pseudo 5-linear subalgebra of V over the real field Q.

Now as in case of strong neutrosophic bivector spaces we can define in case of strong neutrosophic n-vector spaces V and W defined over the same neutrosophic field F; where $V = V_1 \cup V_2 \cup \ldots \cup V_n$ and $W = W_1 \cup W_2 \cup \ldots \cup W_n$; a strong neutrosophic n-linear transformation T from V to W such that $T = T_1 \cup T_2 \cup \ldots \cup T_n : V = V_1 \cup V_2 \cup \ldots \cup V_n \to W = W_1 \cup W_2 \cup \ldots \cup W_n$ with $T_i : V_i \to W_j$; $1 \leq i, j \leq n$ so that no two $V_i$'s are mapped on to the same $W_j$. We denote the collection of all strong neutrosophic n-linear transformations of V to W by $SNHom_F (V, W)$. Like in case of strong neutrosophic bivector spaces we can define strong neutrosophic n-linear operator for strong neutrosophic n-vector space V defined over the field K.

That is if V = W then the strong neutrosophic n-linear transformation will be known as strong neutrosophic n-linear operator on V. $SNHom_K (V, V)$ denotes the set of all strong neutrosophic n-linear operators on V.

Now as in case of usual neutrosophic n-vector spaces over the real field F we can define special (m, n) linear transformation where $V = V_1 \cup V_2 \cup \ldots \cup V_m$ and $W = W_1 \cup W_2 \cup \ldots \cup W_n$, m < n and (m, n) linear transformations when m > n. All properties derived for neutrosophic n-vector spaces (n-linear algebras) defined over a real field can be derived with appropriate modifications in case of strong neutrosophic n-vector spaces (n-linear algebras) defined over the neutrosophic field.

Interested reader can construct examples.
Now we proceed onto define n-basis.

**DEFINITION 3.2.10:** *Let $V = V_1 \cup V_2 \cup ... \cup V_n$ be a strong neutrosophic n-vector space over the neutrosophic field K. A proper n-subset $S = S_1 \cup S_2 \cup ... \cup S_n \subseteq V_1 \cup V_2 \cup ... \cup V_n$ is said to be a n-basis of V if S a n-linearly independent n-set and each $S_i \subseteq V_i$ generates $V_i$ and $S_i$ is a basis of $V_i$ true for each $i = 1, 2, ..., n$.*

*If the n-set $X = X_1 \cup X_2 \cup ... \cup X_n \subseteq V_1 \cup V_2 \cup ... \cup V_n$ is such that each $X_i$ is a linearly independent subset of $V_i$; $i = 1, 2, ..., n$ then we say $X = X_1 \cup X_2 \cup ... \cup X_n$ is a n-linearly independent n-subset of V.*

*A n-basis $S = S_1 \cup S_2 \cup ... \cup S_n \subseteq V_1 \cup V_2 \cup ... \cup V_n$ is always a n-linearly independent n- subset of V over the field F.*

However every n-linearly independent n-subset of V need not be a n-basis of V. If a n-subset $Y = Y_1 \cup Y_2 \cup ... \cup Y_n \subseteq V_1 \cup V_2 \cup ... \cup V_n$ is not a n-linearly independent n-subset of V then we define Y to be a n-linearly dependent n-subset of V.

Interested reader can give examples of these concepts.
We can as in case of neutrosophic n-vector spaces define the notion of n-kernel of a n-linear transformation.

**DEFINITION 3.2.11:** *Let $V = V_1 \cup V_2 \cup ... \cup V_n$ and $W = W_1 \cup W_2 \cup ... \cup W_n$ be two strong neutrosophic n-vector spaces over the neutrosophic field F. Let $T = T_1 \cup T_2 \cup ... \cup T_n : V = V_1 \cup V_2 \cup ... \cup V_n \to W = W_1 \cup W_2 \cup ... \cup W_n$ ; $T_i : V_i \to W_j$; $i = 1, 2, 3, ..., n$ and $j = 1, 2, ..., n$ such that no two $V_i$'s are mapped onto the same $W_j$. The n-kernel of $T = T_1 \cup T_2 \cup ... \cup T_n$ denoted by*

$$kerT = kerT_1 \cup kerT_2 \cup ... \cup kerT_n$$

*where $ker\, T_i = \{v^i \in V_i \mid T_i(v^i) = 0\}$; $i = 1, 2, ..., n$.*

*Thus $kerT = \{(v^1, v^2, ..., v^n) / T(v^1, v^2, ..., v^n) = \{T_1(v^1) \cup T_2(v^2) \cup ... \cup T_n(v^n) = 0 \cup 0 \cup ... \cup 0\}$.*

It is left as a simple exercise for the reader to prove kerT is a proper neutrosophic n-subgroup of V. Further kerT is a strong neutrosophic n-vector subspace of V.

## 3.3 Neutrosophic n-Vector Spaces of Type II

In this section we proceed onto define the new notion of neutrosophic n-vector spaces of type II. We discuss several interesting results about them.

**DEFINITION 3.3.1:** *Let $V = V_1 \cup V_2 \cup ... \cup V_n$; where each $V_i$ is a neutrosophic vector space over $F_i$; $V_i \not\subseteq V_j$ and $V_j \not\subseteq V_i$ (if $i \neq j$, $i \leq i, j \leq n$) and $F_i \not\subseteq F_j$ as well as $F_j \not\subseteq F_i$ ($i \neq j$, $1 \leq i, j \leq n$). We define $V = V_1 \cup V_2 \cup ... \cup V_n$ to be a neutrosophic n-vector space over the real n field $F = F_1 \cup F_2 \cup ... \cup F_n$ of type II.*

We will illustrate this situation by some examples.

***Example 3.3.1:*** Let $V = V_1 \cup V_2 \cup V_3 \cup V_4 \cup V_5 =$

$$\left\{ \begin{pmatrix} a & b \\ c & d \end{pmatrix} \middle| a,b,c,d \in Z_7 I \right\} \cup$$

$$\left\{ \begin{pmatrix} a_1 & a_2 & a_3 \\ a_4 & a_5 & a_6 \end{pmatrix} \middle| a_i \in Z_{11}I; 1 \leq i \leq 6 \right\} \cup$$

$$\left\{ \begin{pmatrix} a_1 & a_2 \\ a_3 & a_4 \\ a_5 & a_6 \\ a_7 & a_8 \\ a_9 & a_{10} \\ a_{11} & a_{12} \end{pmatrix} \middle| a_i \in N(Q); 1 \leq i \leq 12 \right\} \cup$$

$$\left\{ \sum_{i=0}^{12} a_i x^i \text{ ; all polynomials in the variable x with coefficients} \right.$$

from the neutrosophic field $Z_2I$ of degree less than or equal to 12 $a_i \in Z_2I; 0 \le i \le 12 \} \cup$

$$\left\{ \begin{pmatrix} a & 0 & a \\ 0 & b & 0 \\ d & 0 & a \end{pmatrix} \middle| a,b,c,d \in Z_{29}I \right\}$$

be a neutrosophic 5-vector space over the 5-field $F = Z_7 \cup Z_{11} \cup Q \cup Z_2 \cup Z_{29}$ of type II.

***Example 3.3.2:*** Let $V = V_1 \cup V_2 \cup V_3 \cup V_4 \cup V_5 \cup V_6 \cup V_7 =$

$$\left\{ \begin{pmatrix} a \\ b \\ c \\ d \\ e \end{pmatrix} \middle| a,b,c,d,e \in N(Z_2) \right\} \cup$$

$$\left\{ \begin{pmatrix} a_1 & a_2 & a_3 & a_4 \\ a_5 & a_6 & a_7 & a_8 \end{pmatrix} \middle| a_i \in N(Z_3); 1 \le i \le 8 \right\} \cup$$

$$\left\{ \sum_{i=0}^{20} a_i x^i \text{ ; all polynomials in the variable x of degree less than} \right.$$

or equal to 20 with coefficients from the field $Z_{23}I \} \cup$

$$\left\{ \begin{pmatrix} a & b & c & d \\ 0 & g & e & f \\ 0 & 0 & h & i \\ 0 & 0 & 0 & k \end{pmatrix} \middle| a,b,c,d,g,e,f,h,i,k \in Z_{11}I \right\} \cup$$

$$\left\{ \begin{pmatrix} a & b & c \\ b & a & c \\ c & a & b \\ b & c & a \\ c & b & a \\ a & c & b \end{pmatrix} \middle| a,b,c \in N(Q) \right\} \cup$$

$$\left\{ \begin{pmatrix} 0 & 0 & 0 & a \\ 0 & 0 & b & 0 \\ 0 & c & 0 & 0 \\ d & 0 & 0 & 0 \end{pmatrix} \middle| a,b,c,d \in Z_{17}I \right\} \cup$$

$$\left\{ \begin{pmatrix} a_1 & a_2 \\ a_3 & a_4 \\ a_5 & a_6 \\ a_7 & a_8 \\ a_9 & a_{10} \end{pmatrix} \middle| a_i \in N(Z_{31}); 1 \leq i \leq 10 \right\}$$

be a neutrosophic 7-vector space over the 7-field $F = Z_2 \cup Z_3 \cup Z_{23} \cup Z_{11} \cup Q \cup Z_{17} \cup Z_{31}$ of type II.

Even if we do not mention the word type II by the context one of easily understand what type of n-spaces are under study.

**DEFINITION 3.3.2:** *Let $V = V_1 \cup V_2 \cup … \cup V_n$ be a neutrosophic n-vector space over the n-field $F = F_1 \cup F_2 \cup … \cup F_n$. Let $W = W_1 \cup W_2 \cup … \cup W_n \subseteq V_1 \cup V_2 \cup … \cup V_n$, to be a neutrosophic n-vector space over the n-field $F = F_1 \cup F_2 \cup … \cup F_n$, then we define $W_1 \cup W_2 \cup … \cup W_n$ to be a neutrosophic n-vector subspace of V of type II over the n-field F.*

We illustrate this situation by some simple examples.

***Example 3.3.3:*** Let $V = V_1 \cup V_2 \cup V_3 \cup V_4 =$

$$\left\{ \begin{pmatrix} a_1 & a_2 & a_3 & a_4 & a_5 \\ a_6 & a_7 & a_8 & a_9 & a_{10} \end{pmatrix} \middle| a_i \in N(Q); 1 \le i \le 10 \right\} \cup$$

$$\left\{ \begin{pmatrix} a_1 & a_2 \\ a_3 & a_4 \\ a_5 & a_6 \\ a_7 & a_8 \\ a_9 & a_{10} \\ a_{11} & a_{12} \end{pmatrix} \middle| a_i \in Z_{11}I; 1 \le i \le 12 \right\} \cup$$

$$\left\{ \sum_{i=0}^{12} a_i x^i \right.; \text{ all polynomials of degree less than or equal to 12}$$

with coefficients from the neutrosophic field $N(Z_5)$; $a_i \in N(Z_5)$; $0 \le i \le 12 \} \cup$

$$\left\{ \begin{pmatrix} 0 & 0 & 0 & a \\ 0 & 0 & b & 0 \\ 0 & c & 0 & 0 \\ d & 0 & 0 & 0 \end{pmatrix} \middle| a,b,c,d \in N(Z_3) \right\}$$

be a neutrosophic 4-vector space over the 4-field $F = Q \cup Z_{11} \cup Z_5 \cup Z_3$ of type II.

Choose $W = W_1 \cup W_2 \cup W_3 \cup W_4 =$

$$\left\{ \begin{pmatrix} a_1 & a_2 & a_3 & a_4 & a_5 \\ a_6 & a_7 & a_8 & a_9 & a_{10} \end{pmatrix} \middle| a_i \in QI; 1 \le i \le 10 \right\} \cup$$

$$\left\{ \begin{pmatrix} a & a \\ a & a \\ a & a \\ a & a \\ a & a \\ a & a \end{pmatrix} \middle| a \in Z_{11}I \right\} \cup$$

$\left\{ \sum_{i=0}^{12} a_i x^i \right.$ ; all polynomials in the variable x of degree less than or equal to 12 with coefficients from $Z_5 I\}\ \cup$

$$\left\{ \begin{pmatrix} 0 & 0 & 0 & a \\ 0 & 0 & a & 0 \\ 0 & a & 0 & 0 \\ a & 0 & 0 & 0 \end{pmatrix} \middle| a \in N(Z_3) \right\}$$

$\subseteq V_1 \cup V_2 \cup V_3 \cup V_4$; W is a neutrosophic 4-vector subspace of V over the 4-field F of type II.

*Example 3.3.4:* Let $V = V_1 \cup V_2 \cup V_3 \cup V_4 \cup V_5 \cup V_6 =$

$$\left\{ \begin{pmatrix} a_1 & a_2 \\ a_3 & a_4 \\ a_5 & a_6 \\ a_7 & a_8 \end{pmatrix} \middle| a_i \in Z_7 I; 1 \leq i \leq 8 \right\} \cup$$

$\left\{ \sum_{i=0}^{30} a_i x^i \right.$ ; all polynomials of degree less than or equal to thirty with coefficients from the field $N(Z_{13})$; $a_i \in N(Z_{13})$; $0 \leq i \leq 30\}$

$$\cup \left\{ \begin{pmatrix} a_1 & a_2 & a_3 & a_4 & a_5 \\ a_6 & a_7 & a_8 & a_9 & a_{10} \end{pmatrix} \middle| a_i \in N(Z_{17}); 1 \le i \le 10 \right\} \cup$$

$$\left\{ \begin{pmatrix} 0 & 0 & a \\ 0 & b & 0 \\ d & e & f \end{pmatrix} \middle| a,b,d,e,f \in Z_{19}I \right\} \cup$$

$$\left\{ \begin{pmatrix} a_1 \\ a_2 \\ a_3 \\ a_4 \\ a_5 \end{pmatrix} \middle| a_i \in N(Q); 1 \le i \le 5 \right\} \cup$$

{All $10 \times 19$ matrices with entries from $N(Z_2)$} be a neutrosophic 6-vector space over the 6-field $F = F_1 \cup F_2 \cup F_3 \cup F_4 \cup F_5 \cup F_6 = Z_7 \cup Z_{13} \cup Z_{17} \cup Z_{19} \cup Q \cup Z_2$ of type II. Take $W = W_1 \cup W_2 \cup W_3 \cup W_4 \cup W_5 \cup W_6 =$

$$\left\{ \begin{pmatrix} a & a \\ a & a \\ a & a \\ a & a \end{pmatrix} \middle| a \in Z_7 I \right\} \cup$$

$\left\{ \sum_{i=0}^{10} a_i x^i \right.$; all polynomials of degree less than or equal to 10 with coefficients from $Z_{13}I$ in the variable x; $a_i \in Z_{13}I$; $0 \le i \le 10\} \cup$

$$\left\{ \begin{pmatrix} a & a & a & a & a \\ a & a & a & a & a \end{pmatrix} \middle| a \in N(Z_{17}) \right\} \cup$$

$$\left\{ \begin{pmatrix} 0 & 0 & a \\ 0 & a & 0 \\ a & 0 & 0 \end{pmatrix} \middle| a \in Z_{19}I \right\} \cup$$

$$\left\{ \begin{pmatrix} a \\ a \\ a \\ a \\ a \\ a \end{pmatrix} \middle| a \in QI \right\} \cup$$

{All 10×19 matrices with entries from $Z_2I$} $\subseteq V_1 \cup V_2 \cup V_3 \cup V_4 \cup V_5 \cup V_6$ is a neutrosophic 6-vector subspace of V over the 6-field F of type II.

**DEFINITION 3.3.3:** *Let $V = V_1 \cup V_2 \cup V_3 \cup \ldots \cup V_n$ be a neutrosophic n-vector space over the n-field $F = F_1 \cup F_2 \cup \ldots \cup F_n$. If $W = W_1 \cup W_2 \cup \ldots \cup W_n \subseteq V_1 \cup V_2 \cup \ldots \cup V_n$ and $K_1 \cup K_2 \cup \ldots \cup K_n = K \subseteq F = F_1 \cup F_2 \cup \ldots \cup F_n$. If W is a neutrosophic n-vector space over the n-field K then we call W to be a special subneutrosophic n-vector subspace of V over the n-subfield K of F of type II.*

We will illustrate this by some examples and counter examples.

*Example 3.3.5:* Let $V = V_1 \cup V_2 \cup V_3 \cup V_4 \cup V_5 =$

$$\left\{ \begin{pmatrix} a & a & a & a & a \\ b & b & b & b & b \end{pmatrix} \middle| a,b \in N(Z_5) \right\} \cup$$

$$\left\{ \begin{pmatrix} a_1 & a_2 & a_3 \\ a_4 & a_5 & a_6 \\ a_7 & a_8 & a_9 \\ a_{10} & a_{11} & a_{12} \end{pmatrix} \middle| a_i \in N(Z_7); 1 \leq i \leq 12 \right\} \cup$$

$$\left\{ \begin{pmatrix} a_1 \\ a_2 \\ a_3 \\ a_4 \\ a_5 \\ a_6 \\ a_7 \end{pmatrix} \middle| a_i \in Z_{11}I; 1 \le i \le 7 \right\} \cup$$

$$\left\{ \begin{pmatrix} 0 & 0 & 0 & 0 & a \\ 0 & 0 & 0 & b & 0 \\ 0 & 0 & c & 0 & 0 \\ 0 & d & 0 & 0 & 0 \\ e & 0 & 0 & 0 & 0 \end{pmatrix} \middle| a,b,c,d,e \in Z_{19}I \right\} \cup$$

$$\left\{ \begin{pmatrix} a_1 & a_2 & a_3 & 0 & 0 \\ a_4 & a_5 & a_6 & 0 & 0 \\ a_7 & a_8 & a_9 & 0 & 0 \\ 0 & 0 & 0 & a_{10} & a_{11} \\ 0 & 0 & 0 & a_{12} & a_{13} \end{pmatrix} \middle| a_i \in N(Q); 1 \le i \le 13 \right\}$$

be a neutrosophic 5-vector space over the 5-field $F = Z_5 \cup Z_7 \cup Z_{11} \cup Z_{19} \cup Q$. We see each of the fields are prime so F has no 5-subfield. Thus V has no special subneutrosophic 5-vector subneutrosophic 5-vector subspace.

*Example 3.3.6:* Let $V = V_1 \cup V_2 \cup V_3 \cup V_4 =$

$$\left\{ \begin{pmatrix} a & b \\ c & d \end{pmatrix} \middle| a,b,c,d \in N(Q(\sqrt{2}, \sqrt{3}, \sqrt{19})) \right\} \cup$$

$$\left\{ \begin{pmatrix} a & 0 & 0 \\ 0 & b & 0 \\ 0 & 0 & d \end{pmatrix} \middle| a, b, d \in N(Q(\sqrt{17}, \sqrt{3}, \sqrt{5}, \sqrt{13}, \sqrt{11})) \right\} \cup$$

$$\left\{ \begin{pmatrix} a \\ b \\ c \\ d \\ e \end{pmatrix} \middle| a, b, c, d, e \in N(Q(\sqrt{23}, \sqrt{29}, \sqrt{11}, \sqrt{7}, \sqrt{2})) \right\} \cup$$

$$\left\{ \begin{pmatrix} a & b & c & d & e \\ f & g & h & i & j \end{pmatrix} \middle| \begin{array}{l} a, \ldots, j \in N(Q(\sqrt{2}, \sqrt{23}, \sqrt{19}, \\ \sqrt{17}, \sqrt{41}, \sqrt{43}, \sqrt{53})) \end{array} \right\}$$

be a neutrosophic 4-vector space over the 4-field

$$F = F_1 \cup F_2 \cup F_3 \cup F_4$$
$$= Q(\sqrt{2}, \sqrt{3}, \sqrt{19}) \cup Q(\sqrt{17}, \sqrt{3}, \sqrt{5}, \sqrt{13}, \sqrt{11}) \cup$$
$$Q(\sqrt{23}, \sqrt{29}, \sqrt{11}, \sqrt{7}, \sqrt{2}) \cup$$
$$Q(\sqrt{2}, \sqrt{23}, \sqrt{19}, \sqrt{17}, \sqrt{41}, \sqrt{43}, \sqrt{53})$$

of type II. Consider $W = W_1 \cup W_2 \cup W_3 \cup W_4 =$

$$\left\{ \begin{pmatrix} a & a \\ a & a \end{pmatrix} \middle| a \in N(Q(\sqrt{2}, \sqrt{3}, \sqrt{19})) \right\} \cup$$

$$\left\{ \begin{pmatrix} a & 0 & 0 \\ 0 & a & 0 \\ 0 & 0 & a \end{pmatrix} \middle| a \in N(Q(\sqrt{17}, \sqrt{3}, \sqrt{5}, \sqrt{13}, \sqrt{11})) \right\} \cup$$

$$\left\{\begin{pmatrix} a \\ b \\ a \\ b \\ a \end{pmatrix} \middle| a,b \in N(Q(\sqrt{23},\sqrt{29},\sqrt{11},\sqrt{7},\sqrt{2})) \right\} \cup$$

$$\left\{\begin{pmatrix} a & b & c & d & e \\ a & b & c & d & e \end{pmatrix} \middle| \begin{array}{l} a,b,c,d,e \in N(Q(\sqrt{2},\sqrt{23}, \\ \sqrt{19},\sqrt{17},\sqrt{41},\sqrt{43},\sqrt{53})) \end{array}\right\}$$

$\subseteq V_1 \cup V_2 \cup V_3 \cup V_4$. W is a neutrosophic 4-vector space over the 4-field

$$K = K_1 \cup K_2 \cup K_3 \cup K_4$$
$$= Q(\sqrt{3},\sqrt{19}) \cup Q(\sqrt{17},\sqrt{5},\sqrt{13}) \cup$$
$$Q(\sqrt{23},\sqrt{29},\sqrt{11},\sqrt{2}) \cup Q(\sqrt{41},\sqrt{43})$$
$$\subseteq F_1 \cup F_2 \cup F_3 \cup F_4.$$

K is clearly a 4-subfield of F. Thus W is a special sub neutrosophic 4-vector subspace of V over $K = K_1 \cup K_2 \cup K_3 \cup K_4$.

*Example 3.3.7:* Let $V = V_1 \cup V_2 \cup V_3 =$

$$\left\{\begin{pmatrix} a_1 & a_2 & a_3 \\ a_4 & a_5 & a_6 \end{pmatrix} \middle| a_i \in N\left(\frac{Z_2[x]}{\langle x^2+x+1 \rangle}\right)\right\}$$

($\langle x^2 + x + 1 \rangle$ denotes the ideal generated by $x^2 + x +1$, $1 \le i \le 6$)

$$\cup \left\{\begin{pmatrix} a \\ b \\ c \\ d \\ e \end{pmatrix} \middle| a,b,c,d,e \in N\left(\frac{Z_3[x]}{\langle x^3+x^2+x+2 \rangle}\right)\right\} \cup$$

$$\left\{ \begin{pmatrix} 0 & 0 & 0 & a \\ 0 & 0 & b & 0 \\ 0 & d & 0 & 0 \\ c & 0 & 0 & 0 \end{pmatrix} \middle| a,b,c,d \in N\left(\frac{Z_5[x]}{\langle x^2+2 \rangle}\right) \right\}$$

be a neutrosophic 3-vector space over the 3-field

$$F = F_1 \cup F_2 \cup F_3$$

$$= \frac{Z_2[x]}{\langle x^2+x+1 \rangle} \cup \frac{Z_3[x]}{\langle x^3+x^2+x+2 \rangle} \cup \frac{Z_5[x]}{\langle x^2+2 \rangle}.$$

Take $W = W_1 \cup W_2 \cup W_3 =$

$$\left\{ \begin{pmatrix} a & b & c \\ a & b & a \end{pmatrix} \middle| a,b,c \in N\left(\frac{Z_2[x]}{\langle x^2+x+1 \rangle}\right) \right\} \cup$$

$$\left\{ \begin{pmatrix} a \\ a \\ a \\ a \\ a \end{pmatrix} \middle| a \in N\left(\frac{Z_3[x]}{\langle x^3+x^2+x+2 \rangle}\right) \right\} \cup$$

$$\left\{ \begin{pmatrix} 0 & 0 & 0 & a \\ 0 & 0 & a & 0 \\ 0 & a & 0 & 0 \\ a & 0 & 0 & 0 \end{pmatrix} \middle| a \in N\left(\frac{Z_5[x]}{\langle x^2+2 \rangle}\right) \right\}$$

$\subseteq V_1 \cup V_2 \cup V_3$; be a neutrosophic 3-vector space over the 3-field $K = Z_2 \cup Z_3 \cup Z_5 \subseteq F_1 \cup F_2 \cup F_3$. Clearly W is a special subneutrosophic 3-vector subspace of V over the 3-subfield $K = Z_2 \cup Z_3 \cup Z_5 \subseteq F$.

Now a neutrosophic n-vector space $V = V_1 \cup V_2 \cup \ldots \cup V_n$ over a n-field $F = F_1 \cup F_2 \cup \ldots \cup F_n$ of type II is said to be n-simple if V has no proper special subneutrosophic n-vector subspace over the n-subfield $K = K_1 \cup K_2 \cup \ldots \cup K_n \subseteq F = F_1 \cup F_2 \cup \ldots \cup F_n$.

*Example 3.3.8:* Let $V = V_1 \cup V_2 \cup V_3 \cup V_4 \cup V_5 \cup V_6 \cup V_7 =$

$$\left\{ \begin{pmatrix} a & a & b & c \\ d & e & f & g \end{pmatrix} \middle| a,b,c,d,e,f,g \in N(Z_2) \right\} \cup$$

$$\left\{ \begin{pmatrix} a & b \\ c & d \\ e & f \\ g & h \\ i & j \end{pmatrix} \middle| a,b,c,d,e,f,g,h,i,j \in N(Z_3) \right\} \cup$$

$$\left\{ \begin{pmatrix} 0 & 0 & 0 & a \\ 0 & 0 & b & 0 \\ 0 & c & 0 & 0 \\ d & 0 & 0 & 0 \end{pmatrix} \middle| a,b,c,d \in N(Z_5) \right\} \cup$$

$\left\{ \sum_{i=0}^{6} a_i x^i \right.$; all polynomials in the variable x with coefficients from $N(Z_{17})$; $a_i \in N(Z_{17})$; $0 \leq i \leq 6\} \cup$

$$\left\{ \begin{pmatrix} a & b & c & d & e \\ f & g & h & i & j \\ a & b & c & d & e \end{pmatrix} \middle| a,b,c,d,e,f,g,h,i,j \in N(Z_7) \right\} \cup$$

$$\left\{ \begin{pmatrix} a_1 \\ a_2 \\ a_3 \\ a_4 \\ a_5 \\ a_6 \\ a_7 \end{pmatrix} \middle| a_i \in N(Z_{13}); 1 \leq i \leq 7 \right\} \cup$$

$$\left\{ \begin{pmatrix} a_1 & a_2 & a_3 & a_4 & a_{17} \\ a_5 & a_6 & a_7 & a_8 & a_{18} \\ a_9 & a_{10} & a_{11} & a_{12} & a_{19} \\ a_{13} & a_{14} & a_{15} & a_{16} & a_{20} \end{pmatrix} \middle| a_i \in N(Z_{19}); 1 \leq i \leq 20 \right\}$$

be a neutrosophic 7-vector space over the 7-field $F = F_1 \cup F_2 \cup \ldots \cup F_n = Z_2 \cup Z_3 \cup Z_5 \cup Z_{17} \cup Z_7 \cup Z_{13} \cup Z_{19}$. Clearly every $F_i$ in F is a prime field. So F has no proper 7-subfield. Hence even if V has a proper 7-subset $W = W_1 \cup W_2 \cup W_3 \cup W_4 \cup W_5 \cup W_6 \cup W_7 \subseteq V_1 \cup V_2 \cup V_3 \cup V_4 \cup V_5 \cup V_6 \cup V_7$ which is a neutrosophic 7-vector subspace yet it cannot become a special subneutrosophic 7-vector subspace of V as F has no proper 7-subfield.

Inview of this we has the following theorem the proof of which is straight forward.

**THEOREM 3.3.1:** *Let $V = V_1 \cup V_2 \cup \ldots \cup V_n$ be any neutrosophic n-vector space over the n-field $F = F_1 \cup F_2 \cup \ldots \cup F_n$ where each $F_i$ is a prime field; $1 \leq i \leq n$. V is a simple neutrosophic n-vector space over the n-field F.*

**DEFINITION 3.3.4:** *Let $V = V_1 \cup V_2 \cup \ldots \cup V_n$ be a neutrosophic n-vector space over a n-field $F = F_1 \cup F_2 \cup \ldots \cup F_n$ where some of the $F_i$'s are prime fields and some of the $F_j$'s are non prime fields; $1 \leq i, j \leq n$. Let $K = K_1 \cup K_2 \cup \ldots \cup K_n \subseteq F_1 \cup F_2 \cup \ldots \cup F_n$ where some of the $K_i$'s are equal to $F_i$ ($K_i = F_i$) and* *some $K_j$'s are proper subfields of $F_j$ for $1 \leq i, j \leq n$. We call $K = K_1 \cup K_2 \cup ... \cup K_n \subseteq F_1 \cup F_2 \cup ... \cup F_n = F$ to be a quasi n-field. Let $W = W_1 \cup W_2 \cup ... \cup W_n \subseteq V_1 \cup V_2 \cup ... \cup V_n$ be such that W is a neutrosophic n-vector space over $K = K_1 \cup K_2 \cup ... \cup K_n \subseteq F_1 \cup F_2 \cup ... \cup F_n$.*

*We call W as a quasi special neutrosophic n-vector subspace of V over the quasi n-field $K_1 \cup K_2 \cup ... \cup K_n$.*

We will illustrate this situation by an example.

***Example 3.3.9:*** Let $V = V_1 \cup V_2 \cup V_3 \cup V_4 =$

$$\left\{ \begin{pmatrix} a_1 & a_2 & a_3 & a_4 & a_5 \\ a_6 & a_7 & a_8 & a_9 & a_{10} \end{pmatrix} \middle| a_i \in N(R); 1 \leq i \leq 10 \right\} \cup$$

$$\left\{ \begin{pmatrix} a_1 & a_2 & a_3 \\ a_4 & a_5 & a_6 \\ a_7 & a_8 & a_9 \\ a_{10} & a_{11} & a_{12} \\ a_{13} & a_{14} & a_{15} \\ a_{16} & a_{17} & a_{18} \\ a_{19} & a_{20} & a_{21} \end{pmatrix} \middle| a_i \in N(Z_7); 1 \leq i \leq 21 \right\} \cup$$

$$\left\{ \begin{pmatrix} 0 & 0 & 0 & a \\ 0 & 0 & b & 0 \\ 0 & c & 0 & 0 \\ d & 0 & 0 & 0 \end{pmatrix} \middle| a, b, c, d \in N\left\langle \frac{Z_3[x]}{\langle x^3 + x^2 + x + 2 \rangle} \right\rangle \right\} \cup$$

$\left\{ \sum_{i=0}^{24} a_i x^i \right.$ ; all polynomials in the variable x of degree less than or equal to 24 with coefficients from $N(Z_{17})$ $a_i \in N(Z_{17})$; $0 \leq i \leq 24\}$ be a neutrosophic 4-vector space over the 4-field

$$F = F_1 \cup F_2 \cup F_3 \cup F_4$$
$$= R \cup Z_7 \cup \left\langle \frac{Z_3[x]}{\langle x^3 + x^2 + x + 2 \rangle} \right\rangle \cup Z_{17}.$$

Take $W = W_1 \cup W_2 \cup W_3 \cup W_4 =$

$$\left\{ \begin{pmatrix} a & a & a & a & a \\ b & b & b & b & b \end{pmatrix} \middle| a, b \in N(Z_7) \right\} \cup$$

$$\left\{ \begin{pmatrix} a & a & a \\ a & a & a \\ a & a & a \\ b & b & b \\ a & a & a \\ a & a & a \\ a & a & a \end{pmatrix} \middle| a, b \in N(Z_7) \right\} \cup$$

$$\left\{ \begin{pmatrix} 0 & 0 & 0 & a \\ 0 & 0 & b & 0 \\ 0 & c & 0 & 0 \\ d & 0 & 0 & 0 \end{pmatrix} \middle| a, b, c, d \in N(Z_3) \right\} \cup$$

$\left\{ \sum_{i=0}^{15} a_i x^i \right.$; all polynomials in the variable x with coefficients from the neutrosophic field $N(Z_{17})$ of degree less than or equal to 17$\} \subseteq V_1 \cup V_2 \cup V_3 \cup V_4$ ; Take $K = K_1 \cup K_2 \cup K_3 \cup K_4 = Q \cup Z_7 \cup Z_3 \cup Z_{17} \subseteq F_1 \cup F_2 \cup F_3 \cup F_4$. Clearly W is a neutrosophic 4-vector space over the 4-field $K = K_1 \cup K_2 \cup K_3 \cup K_4$. Thus W is a quasi special neutrosophic 4-vector subspace of V over the 4-quasi field K.

The notion of n-basis and n-linearly independent elements can be defined as in case of neutrosophic n-vector spaces of type I here each $S_i \subseteq V_i$ is a basis of $V_i$ over $F_i$; $i = 1, 2, \ldots, n$ where $V = V_1 \cup V_2 \cup \ldots \cup V_n$ is defined over the n-field $F = F_1 \cup F_2 \cup \ldots \cup F_n$.

The reader is expected to construct examples. Also the notion of n-linear transformation of type II neutrosophic n-vector spaces can be defined as in case of type I n-vector spaces with necessary changes. Further the notion of kernel T, T a n-linear transformation of the neutrosophic n-vector space V to a neutrosophic n-vector space W defined over the same n-field $F = F_1 \cup \ldots \cup F_n$ is defined as follows.

Let $T = T_1 \cup T_2 \cup \ldots \cup T_n : V = V_1 \cup V_2 \cup \ldots \cup V_n \to W = W_1 \cup W_2 \cup \ldots \cup W_n$ be a map such that $T_i : V_i \to W_j$ is a linear transformation and no two $V_i$ is mapped onto the same $W_j$; $1 \le i, j \le n$. $i = 1, 2, \ldots, n$. The n-kernel of T denoted by $\ker T = \ker T_1 \cup \ker T_2 \cup \ldots \cup \ker T_n$ where $\ker T_i = \{v^i \in V_i \mid T_i(v^i) = \bar{0}\}$, $i = 1, 2, \ldots, n$. Thus

$$\ker T = \ker T_1 \cup \ker T_2 \cup \ldots \cup \ker T_n$$
$$= \{(v^1, v^2, \ldots, v^n) \in V_1 \cup V_2 \cup \ldots \cup V_n \, / \, T(v^1, v^2, \ldots, v^n)$$
$$= T_1(v^1) \cup T_2(v^2) \cup \ldots \cup T_n(v^n)$$
$$= 0 \cup 0 \cup \ldots \cup 0\}.$$

It is easily verified that kerT is a neutrosophic n-vector subspace of V over the n-field $F = F_1 \cup F_2 \cup \ldots \cup F_n$ of type II.

We can prove if V and W are n-finite dimensional and T a n-linear transformation then n rank T + n nullity T = $(n_1, n_2, \ldots, n_n)$ dim V = n-dimension of V.

(rank $T_1 \cup$ rank $T_2 \cup \ldots \cup$ rank $T_n$) + nullity $T_1 \cup \ldots \cup$ nullity $T_n$ = dim $V_1 \cup$ dim $V_2 \cup \ldots \cup$ dim $V_n$ (dim $V_i$ over the field $F_i$, $i = 1, 2, \ldots, n$)

Thus (rank $T_1$ + nullity $T_1$) $\cup$ (rank $T_2$ + nullity $T_2$) $\cup \ldots \cup$ rank $T_n$ + nullity $T_n$ = $(n_1, n_2, \ldots, n_n)$.

Further if $V = V_1 \cup V_2 \cup \ldots \cup V_n$ and $W = W_1 \cup W_2 \cup \ldots \cup W_n$ be two neutrosophic n-vector spaces over the n-field $F = F_1 \cup F_2 \cup \ldots \cup F_n$ of type II and if $S = S_1 \cup S_2 \cup \ldots \cup S_n$ and $T = T_1 \cup T_2 \cup \ldots \cup T_n$ are two n-linear transformations of V to W then the n function (T+S) =

$(T_1 \cup T_2 \cup \ldots \cup T_n + S_1 \cup S_2 \cup \ldots \cup S_n) = T_1 + S_1 \cup T_2 + S_2 \cup \ldots \cup T_n + S_n$ is defined by

$(T + S)\alpha$
$= ((T_1 + S_1) \cup (T_2 + S_2) \cup \ldots \cup (T_n + S_n)] (\alpha_1 \cup \alpha_2 \cup \ldots \cup \alpha_n)$
$= (T_1 + S_1) \alpha_1 \cup (T_2 + S_2) \alpha_2 \cup \ldots \cup (T_n + S_n) \alpha_n$
$= (T_1\alpha_1 + S_1\alpha_1) \cup (T_2\alpha_2 + S_2\alpha_2) \cup \ldots \cup (T_n\alpha_n + S_n\alpha_n)$

is a neutrosophic n-linear transformation from $V = V_1 \cup V_2 \cup \ldots \cup V_n$ to $W_1 \cup W_2 \cup \ldots \cup W_n = W$ and $\alpha_1 \cup \alpha_2 \cup \ldots \cup \alpha_n \in V_1 \cup V_2 \cup \ldots \cup V_n$. Also if $c = c_1 \cup c_2 \cup c_3 \cup \ldots \cup c_n \in F_1 \cup F_2 \cup \ldots \cup F_n$ then $(c_1 \cup c_2 \cup \ldots \cup c_n)(T_1 \cup T_2 \cup \ldots \cup T_n)(\alpha_1 \cup \alpha_2 \cup \ldots \cup \alpha_n) = c_1 T_1 \alpha_1 \cup c_2 T_2 \alpha_2 \cup \ldots \cup c_n T_n \alpha_n$.

Thus the set of all n-linear transformations of V to V with n-addition and n-scalar multiplication defined above is again a neutrosophic n-vector space of type II over the n-field $F = F_1 \cup F_2 \cup \ldots \cup F_n$.

Thus $NL(V, W) = NL^1(V_1, W_1) \cup NL^2(V_2, W_2) \cup \ldots \cup NL^n(V_n, W_n)$ is a neutrosophic n-vector space over the n-field $F = F_1 \cup F_2 \cup \ldots \cup F_n$ where $V_i$ and $W_i$ are neutrosophic vector spaces defined over the field $F_i$, $i = 1, 2, \ldots, n$.

Now we proceed onto define the new notion of neutrosophic n-linear algebra of type II.

**DEFINITION 3.3.5:** *Let $V = V_1 \cup V_2 \cup \ldots \cup V_n$ be a neutrosophic n-vector space over the n-field $F = F_1 \cup F_2 \cup \ldots \cup F_n$ of type II. If each $V_i$ is a neutrosophic linear algebra over $F_i$ for $i = 1, 2, \ldots, n$; then we call V to be a neutrosophic n-linear algebra over the n-field F of type II.*

We will illustrate this situation by some examples.

*Example 3.3.10:* Let $V = V_1 \cup V_2 \cup V_3 \cup V_4 =$

$$\left\{ \begin{pmatrix} a & b \\ c & d \end{pmatrix} \middle| a,b,c,d \in Z_7 I \right\} \cup$$

$$\left\{ \sum_{i=0}^{\infty} a_i x^i \text{ ; all polynomials in the variable x with coefficients} \right.$$

from $N(Z_{11})$ that is $a_i \in (Z_{11})$; $0 \le i \le \infty\} \cup$

$$\left\{ \begin{pmatrix} a & 0 & 0 & 0 \\ b & c & 0 & 0 \\ d & e & g & 0 \\ p & q & r & s \end{pmatrix} \middle| a,b,c,d,e,g,p,q,r,s \in N(Q) \right\} \cup$$

$$\left\{ \begin{pmatrix} a & b & c \\ 0 & d & e \\ 0 & 0 & f \end{pmatrix} \middle| a,b,c,d,e,f \in N(Z_2) \right\}$$

be a neutrosophic 4-linear algebra over the 4-field $F = F_1 \cup F_2 \cup F_3 \cup F_4 = Z_7 \cup Z_{11} \cup Q \cup Z_2$.

*Example 3.3.11:* Let $V = V_1 \cup V_2 \cup V_3 \cup V_4 \cup V_5 \cup V_6 =$

$$\left\{ \sum_{i=0}^{\infty} a_i x^i \text{ ; all polynomials in the variable x with coefficients} \right.$$

from $N(Z_{19})$. $a_i \in N(Z_{19})$; $0 \le i \le \infty\} \cup$

$$\left\{ \begin{pmatrix} a & 0 \\ b & d \end{pmatrix} \middle| a,b,d \in Z_{13}I \right\} \cup$$

$$\left\{ \begin{pmatrix} a & b & c \\ d & e & f \\ g & h & i \end{pmatrix} \middle| a,b,c,d,e,f,g,h,i \in N(Q) \right\} \cup$$

{All $10 \times 10$ upper triangular matrices with entries from the neutrosophic field $N(Z_{23})\} \cup$ {All $12 \times 12$ diagonal matrices with entries from $Z_2I$} $\cup$ {All $5 \times 5$ lower triangular matrices with entries from the neutrosophic field. N $(Z_{41})$}; is a neutrosophic 6-linear algebra over the 6-field $F = F_1 \cup F_2 \cup F_3 \cup F_4 \cup F_5 \cup F_6 = Z_{19} \cup Z_{13} \cup Q \cup Z_{23} \cup Z_2 \cup Z_{41}$ of type II.

Now we proceed onto define the substructure in neutrosophic n-linear algebras of type II.

**DEFINITION 3.3.6:** *Let $V = V_1 \cup V_2 \cup ... \cup V_n$ be a neutrosophic n-linear algebra over the n-field $F = F_1 \cup F_2 \cup ... \cup F_n$. Let $W = W_1 \cup W_2 \cup ... \cup W_n \subseteq V_1 \cup V_2 \cup ... \cup V_n$ ; if W itself is a neutrosophic n-linear algebra of type II over the n-field $F = F_1 \cup F_2 \cup ... \cup F_n$ and $W_i \neq \{0\}$ or $W_i \neq V_i$ for every i, i = 1, 2, ..., n. We call W to be a neutrosophic n-linear subalgebra of V over the n-field F of type II. If V has no neutrosophic n-linear subalgebras then we define V to be a n-simple neutrosophic n-linear algebra over F of type II.*

We will illustrate both the situations by some simple examples.

***Example 3.3.12:*** Let $V = V_1 \cup V_2 \cup V_3 \cup V_4 \cup V_5 \cup V_6 \cup V_7$

$$= \left\{ \begin{pmatrix} a & b \\ c & d \end{pmatrix} \middle| a, b, c, d \in N(Q) \right\} \cup$$

$\{(a_1, a_2, a_3, a_4, a_5) \mid a_i \in Z_{11}I; 1 \leq i \leq 5\} \cup$

$\left\{ \sum_{i=0}^{\infty} a_i x^i \right.$; all polynomials in the variable x with coefficients from the neutrosophic field $N(Z_{19}) \cup$ {All $10 \times 10$ neutrosophic matrices with entries form $Z_{41}I$} $\cup$ {All $8 \times 8$ upper triangular matrices with entries from $N(Z_2)$} $\cup$ {All $5 \times 5$ low triangular matrices with entries from $N(Z_5)$} $\cup$ {$3 \times 3$ matrices with entries from $N(Z_{31})$} be a neutrosophic 7-linear algebra over the 7-field $F = F_1 \cup ... \cup F_7 = Q \cup Z_{11} \cup Z_{19} \cup Z_{41} \cup Z_2 \cup Z_5 \cup Z_{31}$. Consider $W = W_1 \cup W_2 \cup W_3 \cup W_4 \cup W_5 \cup W_6 \cup W_7 =$

$$\left\{\begin{pmatrix} a & b \\ c & d \end{pmatrix} \middle| a,b,c,d \in QI\right\} \cup \{(a\ a\ a\ a\ 0) \mid a \in Z_{11}I\} \cup$$

$\left\{\sum_{i=0}^{\infty} a_i x^{2i}\right.$ ; all polynomials of even degree in the variable x with coefficients from $Z_{19}I\}$ $\cup$ {All $10 \times 10$ neutrosophic upper triangular matrices with entries from $Z_{41}I\}$ $\cup$ {All $8 \times 8$ upper triangular matrices with entries from $Z_2I\}\cup$ {All $5 \times 5$ lower triangular matrices with entries from $Z_5I\}$ $\cup$ {$3 \times 3$ matrices with entries from $Z_{31}I\}$ $\subseteq V_1 \cup V_2 \cup V_3 \cup V_4 \cup V_5 \cup V_6 \cup V_7$; it is easily verified W is a neutrosophic 7-linear subalgebra of V over the 7-field F.

***Example 3.3.13:*** Let $V = V_1 \cup V_2 \cup V_3 \cup V_4 \cup V_5 =$

$$\left\{\begin{pmatrix} a & a \\ a & a \end{pmatrix} \middle| a \in Z_2I\right\} \cup \{(a\ a\ a\ a\ a\ a) \mid a \in Z_3I\} \cup$$

$$\left\{\begin{pmatrix} a & 0 & 0 \\ a & a & 0 \\ a & a & a \end{pmatrix} \middle| a \in Z_5I\right\} \cup \left\{\begin{pmatrix} a & 0 & 0 & 0 & 0 \\ 0 & a & 0 & 0 & 0 \\ 0 & 0 & a & 0 & 0 \\ 0 & 0 & 0 & a & 0 \\ 0 & 0 & 0 & 0 & a \end{pmatrix} \middle| a \in Z_7I\right\} \cup$$

$$\left\{\begin{pmatrix} a & a & a & a \\ 0 & a & a & a \\ 0 & 0 & a & a \\ 0 & 0 & 0 & a \end{pmatrix} \middle| a \in Z_{13}I\right\}$$

be a neutrosophic 5-linear algebra over the 5-field $F = Z_2 \cup Z_3 \cup Z_5 \cup Z_7 \cup Z_{13}$. It is easy to verify V has no proper neutrosophic 5-sublinear algebras of type II over F in V. Thus V is a simple neutrosophic 5-linear algebra over F.

**DEFINITION 3.3.7:** *Let $V = V_1 \cup V_2 \cup ... \cup V_n$ be a neutrosophic n-linear algebra over the n-field $F = F_1 \cup F_2 \cup ... \cup F_n$ of type II. Let $W = W_1 \cup W_2 \cup ... \cup W_n \subseteq V_1 \cup V_2 \cup ... \cup V_n$ be such that W is a neutrosophic n-vector space over the n-field F and no $W_i \subseteq V_i$ is a neutrosophic linear algebra over $F_i$; $W_i \neq V_i$, $i = 1, 2, ..., n$. We define W to be a pseudo neutrosophic n-vector subspace of V over the n-field F.*

We will illustrate this situation by some examples.

*Example 3.3.14:* Let $V = V_1 \cup V_2 \cup V_3 \cup V_4 =$

$$\left\{ \begin{pmatrix} a & b & c \\ d & e & f \\ g & h & i \end{pmatrix} \middle| a,b,c,d,e,f,g,h,i \in Z_7I \right\} \cup$$

$\left\{ \sum_{i=0}^{\infty} a_i x^i \right.$; all polynomials in the variable x with coefficients from $Z_2I\} \cup$

$$\left\{ \begin{pmatrix} a & b & c & d \\ a & e & f & g \\ a & a & p & g \\ a & a & a & r \end{pmatrix} \middle| a,b,c,d,e,f,g,p,q,r \in N(Q) \right\} \cup$$

{All $7 \times 7$ neutrosophic matrices with entries from $N(Z_{11})$} be a neutrosophic 4-linear algebra over 4-field $F = Z_7 \cup Z_2 \cup Q \cup Z_{11}$. Let $W = W_1 \cup W_2 \cup W_3 \cup W_4 =$

$$\left\{ \begin{pmatrix} 0 & 0 & b \\ 0 & c & 0 \\ d & 0 & 0 \end{pmatrix} \middle| b,c,d \in Z_7I \right\} \cup$$

$$\left\{ \sum_{i=0}^{12} a_i x^i \text{ ; all polynomials of degree less than or equal to 12} \right.$$

with coefficients from $Z_2 I$; $a_i \in Z_2 I$; $0 \le i \le 12\} \cup$

$$\left\{ \begin{pmatrix} 0 & 0 & 0 & d \\ 0 & 0 & f & 0 \\ 0 & a & 0 & 0 \\ a & 0 & 0 & 0 \end{pmatrix} \middle| f, d, a \in N(Q) \right\} \cup$$

$$\left\{ \begin{pmatrix} 0 & 0 & 0 & 0 & 0 & 0 & a \\ 0 & 0 & 0 & 0 & 0 & b & 0 \\ 0 & 0 & 0 & 0 & d & 0 & 0 \\ 0 & 0 & 0 & e & 0 & 0 & 0 \\ 0 & 0 & f & 0 & 0 & 0 & 0 \\ 0 & g & 0 & 0 & 0 & 0 & 0 \\ h & 0 & 0 & 0 & 0 & 0 & 0 \end{pmatrix} \middle| a,b,c,d,e,f,g,h \in N(Z_{11}) \right\}$$

$\subseteq V_1 \cup V_2 \cup V_3 \cup V_4$ be a neutrosophic 4-vector space over F $= Z_7 \cup Z_2 \cup Q \cup Z_{11}$. Clearly W is a pseudo neutrosophic 4-vector subspace of V over the 4-field $F = Z_7 \cup Z_2 \cup Q \cup Z_{11}$ of type II.

*Example 3.3.15:* Let $V = V_1 \cup V_2 \cup V_3 \cup V_4 \cup V_5 \cup V_6 =$

$$\left\{ \begin{pmatrix} a & b \\ c & d \end{pmatrix} \middle| a, b, c, d \in N(Q) \right\} \cup$$

$\left\{ \sum_{i=0}^{\infty} a_i x^i \text{ ; all polynomials in the variable x with coefficients} \right.$ from $N(Z_{11})\} \cup$ {all $8 \times 8$ neutrosophic matrices with entries from $N(Z_{29})\} \cup$ {all $7 \times 7$ upper triangular matrices with entries from $N(Z_3)\} \cup \{4 \times 4$ matrices with entries from $Z_{23}I\} \cup \{5 \times 5$ neutrosophic matrices with entries from $N(Z_{41})\}$ be a neutrosophic 6-linear algebra over the 6-field algebra over the 6-field $F = Q \cup Z_{11} \cup Z_{29} \cup Z_3 \cup Z_{23} \cup Z_{41}$ of type II. Let $W = W_1 \cup W_2 \cup W_3 \cup W_4 \cup W_5 \cup W_6 =$

$$\left\{ \begin{pmatrix} 0 & a \\ b & 0 \end{pmatrix} \middle| a, b \in QI \right\} \cup$$

$\left\{ \sum_{i=0}^{5} a_i x^i \right.$ ; all polynomials in the variable x of degree less than or equal to 5 with coefficients from $N(Z_{11})\} \cup$

$$\left\{ \begin{pmatrix} 0 & 0 & 0 & 0 & 0 & 0 & 0 & a \\ 0 & 0 & 0 & 0 & 0 & 0 & b & 0 \\ 0 & 0 & 0 & 0 & 0 & c & 0 & 0 \\ 0 & 0 & 0 & 0 & d & 0 & 0 & 0 \\ 0 & 0 & 0 & 0 & 0 & 0 & 0 & 0 \\ 0 & 0 & g & 0 & 0 & 0 & 0 & 0 \\ 0 & f & 0 & 0 & 0 & 0 & 0 & 0 \\ p & 0 & 0 & 0 & 0 & 0 & 0 & 0 \end{pmatrix} \middle| a,b,c,d,g,f,p \in Z_{29}I \right\} \cup$$

$$\left\{ \begin{pmatrix} 0 & 0 & 0 & 0 & 0 & 0 & a \\ 0 & 0 & 0 & 0 & 0 & b & 0 \\ 0 & 0 & 0 & 0 & c & 0 & 0 \\ 0 & 0 & 0 & g & 0 & 0 & 0 \\ 0 & 0 & 0 & 0 & 0 & 0 & 0 \\ 0 & 0 & 0 & 0 & 0 & 0 & 0 \\ 0 & 0 & 0 & 0 & 0 & 0 & 0 \end{pmatrix} \middle| a,b,c,g \in Z_3I \right\} \cup$$

$$\left\{ \begin{pmatrix} 0 & 0 & a & b \\ 0 & 0 & c & 0 \\ 0 & d & g & 0 \\ g & 0 & 0 & f \end{pmatrix} \middle| a,b,c,d,g,f \in Z_3I \right\} \cup$$

$$\left\{ \begin{pmatrix} 0 & 0 & a & b & c \\ 0 & 0 & 0 & d & 0 \\ 0 & 0 & e & f & 0 \\ 0 & g & 0 & 0 & h \\ g & 0 & 0 & 0 & p \end{pmatrix} \middle| a,b,c,d,e,f,g,h,p \in Z_{41}I \right\}$$

$\subseteq V_1 \cup V_2 \cup V_3 \cup V_4 \cup V_5 \cup V_6$ is a pseudo neutrosophic 6-vector subspace of V over $F = Q \cup Z_{11} \cup Z_{29} \cup Z_3 \cup Z_{43} \cup Z_{41}$ of type II.

**DEFINITION 3.3.8:** *Let $V = V_1 \cup V_2 \cup ... \cup V_n$ be a n-finite neutrosophic n-vector space (n-linear algebra) of type II over the n-field $F = F_1 \cup F_2 \cup ... \cup F_n$. Suppose $W = W_1 \cup W_2 \cup ... \cup W_n \subseteq V_1 \cup V_2 \cup ... \cup V_n$ be a neutrosophic n-vector subspace (n-linear subalgebra) of V of n-dimension ($n_1 - 1$, $n_2 - 1$, ..., $n_n - 1$) over the n-field F of type II, where n-dimension of V is ($n_1$, $n_2$, ..., $n_n$). Then we define W to be a neutrosophic hyper n-space (n-algebra) of V.*

The reader is requested to give examples of the above definition.

We define neutrosophic n-polynomial ring or neutrosophic polynomial n-ring over the n-field $F = F_1 \cup F_2 \cup ... \cup F_n$ to be $F_1[x] \cup F_2[x] \cup ... \cup F_n[x] = F[x]$ where $F_1, F_2, ..., F_n$ are n distinct neutrosophic fields.

*Example 3.3.16:* $F[x] = N(Z_7)[x] \cup Z_{11}I[x] \cup RI[x] \cup Z_2I[x] \cup N[Z_{13}][x] \cup N(Z_{47})[x]$ is a 6-polynomial neutrosophic ring.

**DEFINITION 3.3.9:** *Let $V = V_1 \cup V_2 \cup ... \cup V_n$ be a neutrosophic n-vector space (linear algebra) over the n-field $F = F_1 \cup F_2 \cup ... \cup F_n$. Let $W = W_1 \cup W_2 \cup ... \cup W_n \subseteq V_1 \cup V_2 \cup ... \cup V_n$ be such that W is only just a n-vector space (linear algebra) over the n-field F then we all W to be pseudo n-vector space (n-linear algebra) of V over the n-field $F = F_1 \cup F_2 \cup ... \cup F_n$.*

We will illustrate this situation by some simple examples.

***Example 3.3.17:*** Let V = $V_1 \cup V_2 \cup V_3 \cup V_4 \cup V_5$ =

$$\left\{ \begin{pmatrix} a_1 & a_2 & a_3 \\ a_4 & a_5 & a_6 \end{pmatrix} \middle| a_i \in N(Q); 1 \leq i \leq 6 \right\} \cup$$

$$\left\{ \begin{pmatrix} a_1 & a_2 \\ a_3 & a_4 \\ a_5 & a_6 \\ a_7 & a_8 \\ a_9 & a_{10} \end{pmatrix} \middle| a_i \in N(Z_{11}); 1 \leq i \leq 10 \right\} \cup$$

$$\left\{ \begin{pmatrix} 0 & 0 & 0 & a_1 \\ 0 & 0 & a_2 & a_3 \\ 0 & a_4 & a_5 & a_6 \\ a_7 & a_8 & a_9 & a_{10} \end{pmatrix} \middle| a_i \in N(Z_2); 1 \leq i \leq 10 \right\} \cup$$

$\left\{ \sum_{i=0}^{19} a_i x^i \right.$ ; all polynomials in the variable x with coefficients from $N(Z_{17})$ of degree less than or equal to 19; $a_i \in N(Z_{17})$; $0 \leq i \leq 19 \} \cup$

$$\left\{ \begin{pmatrix} a_1 & a_2 & a_5 & a_6 \\ a_3 & a_4 & a_7 & 0 \\ a_8 & a_9 & 0 & 0 \\ a_{10} & 0 & 0 & 0 \end{pmatrix} \middle| a_i \in N(Z_{43}); 1 \leq i \leq 10 \right\}$$

be a neutrosophic 5-vector space over the 5-field F = $F_1 \cup F_2 \cup F_3 \cup F_4 \cup F_5 = Q \cup Z_{11} \cup Z_2 \cup Z_{17} \cup Z_{43}$ of type II. Let W = $W_1 \cup W_2 \cup W_3 \cup W_4 \cup W_5$ =

$$\left\{ \begin{pmatrix} a_1 & a_2 & a_3 \\ a_1 & a_2 & a_3 \end{pmatrix} \middle| a_i \in Q; 1 \le i \le 3 \right\} \cup$$

$$\left\{ \begin{pmatrix} a & a \\ b & b \\ c & c \\ d & d \\ e & e \end{pmatrix} \middle| a,b,c,d,e \in Z_{11} \right\} \cup$$

$$\left\{ \begin{pmatrix} 0 & 0 & 0 & a \\ 0 & 0 & b & b \\ 0 & c & c & c \\ d & d & d & d \end{pmatrix} \middle| a,b,c,d \in Z_2 \right\} \cup$$

$\left\{ \sum_{i=0}^{10} a_i x^i \right.$ ; all polynomials in the variable x of degree less than or equal to 10 with coefficients from $Z_{17}\} \cup$

$$\left\{ \begin{pmatrix} a & a & a & a \\ b & b & b & 0 \\ c & c & 0 & 0 \\ d & 0 & 0 & 0 \end{pmatrix} \middle| a,b,c,d \in Z_{43} \right\}$$

$\subseteq V_1 \cup V_2 \cup V_3 \cup V_4 \cup V_5$ be a pseudo 5-vector space of V over the 5-field F.

*Example 3.3.18:* Let $V = V_1 \cup V_2 \cup V_3 \cup V_4 \cup V_5 \cup V_6 =$

$$\left\{ \begin{pmatrix} a & b \\ c & d \end{pmatrix} \middle| a,b,c,d \in N(Z_2) \right\} \cup$$

{All $9 \times 9$ neutrosophic matrices with entries from $N(Z_{11})$} $\cup$
{all $7 \times 7$ upper triangular neutrosophic matrices with entries from $N(Z_{13})$} $\cup$

$$\left\{ \begin{pmatrix} a & b & 0 & 0 \\ c & d & 0 & 0 \\ 0 & 0 & e & f \\ 0 & 0 & g & h \end{pmatrix} \middle| a,b,c,d,e,f,g,h \in N(Z_{47}) \right\} \cup$$

$$\left\{ \sum_{i=0}^{\infty} a_i x^i \; ; \text{ all polynomial in the variable x with coefficients from} \right.$$
$N(Z_{53})$} $\cup$

$$\left\{ \begin{pmatrix} a_1 & a_2 & a_3 & 0 & 0 & 0 \\ a_4 & a_5 & a_6 & 0 & 0 & 0 \\ a_7 & a_8 & a_9 & 0 & 0 & 0 \\ 0 & 0 & 0 & a_{10} & a_{11} & a_{12} \\ 0 & 0 & 0 & a_{13} & a_{14} & a_{15} \\ 0 & 0 & 0 & a_{16} & a_{17} & a_{18} \end{pmatrix} \middle| a_i \in N(Z_5); 1 \le i \le 18 \right\}$$

be a neutrosophic 6-linear algebra over the field $F = Z_2 \cup Z_{11} \cup Z_{13} \cup Z_{47} \cup Z_{53} \cup Z_5$ of type II. Consider $W = W_1 \cup W_2 \cup W_3 \cup W_4 \cup W_5 \cup W_6 =$

$$\left\{ \begin{pmatrix} a & b \\ c & d \end{pmatrix} \middle| a,b,c,d \in Z_2 \right\} \cup$$

{All $9 \times 9$ matrices with entries from $Z_{11}$} $\cup$ {all $7 \times 7$ upper triangular matrices with entries from $Z_{13}$} $\cup$

$$\left\{\begin{pmatrix} a & a & 0 & 0 \\ a & a & 0 & 0 \\ 0 & 0 & b & b \\ 0 & 0 & b & b \end{pmatrix} \middle| a,b \in Z_{47}\right\} \cup$$

$$\left\{\sum_{i=0}^{\infty} a_i x^i \text{ ; all polynomial in the variable x with coefficients from } Z_{53}\right\} \cup$$

$$\left\{\begin{pmatrix} a & a & a & 0 & 0 & 0 \\ a & a & a & 0 & 0 & 0 \\ a & a & a & 0 & 0 & 0 \\ 0 & 0 & 0 & b & b & b \\ 0 & 0 & 0 & b & b & b \\ 0 & 0 & 0 & b & b & b \end{pmatrix} \middle| a,b \in Z_5\right\}$$

$\subseteq V_1 \cup V_2 \cup V_3 \cup V_4 \cup V_5 \cup V_6$; is a 6-linear algebra over F. Thus W is a pseudo 6 linear algebra of V over the 6-field F of type II.

A neutrosophic n-vector space (n-linear algebra) $V = V_1 \cup V_2 \cup \ldots \cup V_n$ of type II over the n-field $F = F_1 \cup F_2 \cup \ldots \cup F_n$ is said to be pseudo simple n-vector space (n-linear algebra) if V does not contain any pseudo n-vector subspace (n-linear subalgebra).

In view of this we give the following theorem which guarantees the existence of pseudo simple n-vector spaces (n-linear algebras).

**THEOREM 3.3.2:** *Let $V = V_1 \cup V_2 \cup \ldots \cup V_n$ be a neutrosophic n-vector space (n-linear algebra) over the n-field $F = F_1 \cup F_2 \cup \ldots \cup F_n$ of type II. If each $V_i$ is defined over a field $K_i$ with $K_i = F_i$ I ($F_i$ real prime field) for $i = 1, 2, \ldots, n$. V is a pseudo simple n-vector space (n-linear algebra) over F of type II.*

*Proof:* Given $V = V_1 \cup V_2 \cup \ldots \cup V_n$ is a neutrosophic n-vector space over the n-field $F = F_1 \cup F_2 \cup \ldots \cup F_n$ where each $V_i$ is defined over $F_iI$; $i = 1, 2, \ldots, n$. Thus there does not exist a $W_i \subseteq V_i$ where $W_i$ is a real vector space over $F_i$ as $F_i \not\subseteq F_iI$ for $i = 1, 2, \ldots, n$. Thus $V = V_1 \cup V_2 \cup \ldots \cup V_n$ does not contain any n-vector subspace. So V is a pseudo simple n-vector space.

We will illustrate this by some simple examples.

*Example 3.3.19:* Let $V = V_1 \cup V_2 \cup V_3 \cup V_4 \cup V_5 \cup V_6 \cup V_7$

$$= \left\{ \begin{pmatrix} a & b & c \\ d & e & f \end{pmatrix} \middle| a,b,c,d,e,f \in Z_2I \right\} \cup$$

$$\left\{ \begin{pmatrix} a_1 & a_2 & a_3 \\ a_4 & a_5 & a_6 \\ a_7 & a_8 & a_9 \\ a_{10} & a_{11} & a_{12} \\ a_{13} & a_{14} & a_{15} \end{pmatrix} \middle| a_i \in Z_7I; 1 \le i \le 15 \right\} \cup$$

$$\left\{ \begin{pmatrix} a_1 & a_2 & a_3 & a_4 & a_5 & a_6 & a_7 \\ a_8 & a_9 & a_{10} & a_{11} & a_{12} & a_{13} & a_{14} \end{pmatrix} \middle| a_i \in Z_{11}I; 1 \le i \le 14 \right\} \cup$$

$$\left\{ \begin{pmatrix} a_1 \\ a_2 \\ a_3 \\ a_4 \\ a_5 \\ a_6 \\ a_7 \\ a_8 \end{pmatrix} \middle| a_i \in Z_{19}I; 1 \le i \le 8 \right\} \cup$$

$$\left\{ \sum_{i=0}^{8} a_i x^i \text{ ; all polynomials in the variable x of degree less than} \right.$$

or equal to 8 with coefficients from the field QI; $a_i \in QI$; $0 \le i \le 8$} $\cup$

$$\left\{ \begin{pmatrix} 0 & 0 & 0 & a \\ 0 & 0 & b & d \\ 0 & c & d & f \\ g & h & p & q \end{pmatrix} \middle| a,b,d,c,f,g,h,p,q \in Z_3 I \right\} \cup$$

$$\left\{ \begin{pmatrix} 0 & 0 & 0 & 0 & 0 & a_1 \\ 0 & 0 & 0 & 0 & a_2 & a_3 \\ 0 & 0 & 0 & a_4 & 0 & 0 \\ 0 & 0 & a_5 & 0 & a_6 & 0 \\ 0 & a_7 & 0 & 0 & 0 & 0 \\ a_8 & 0 & 0 & 0 & 0 & a_9 \end{pmatrix} \middle| a_i \in Z_{31} I; 1 \le i \le 9 \right\}$$

be a neutrosophic 7-vector space over the 7-field $F = Z_2 \cup Z_7 \cup Z_{11} \cup Z_{19} \cup Q \cup Z_3 \cup Z_{31}$ of type II. It is easily verified that V has no proper pseudo 7-vector subspace over the 7-field F. Hence V is a pseudo simple 7-vector space.

***Example 3.3.20:*** Let $V = V_1 \cup V_2 \cup V_3 \cup V_4 =$

$$\left\{ \begin{pmatrix} a & b \\ c & d \end{pmatrix} \middle| a,b,c,d \in QI \right\} \cup$$

$\{(a_1, a_2, a_3, a_4, a_5) \mid a_i \in Z_7 I; 1 \le i \le 5\} \cup$

$$\left\{ \begin{pmatrix} a_1 & 0 & 0 & 0 \\ a_2 & a_3 & 0 & 0 \\ a_4 & a_5 & a_6 & 0 \\ a_7 & a_8 & a_9 & a_{10} \end{pmatrix} \middle| a_i \in Z_{19} I; 1 \le i \le 10 \right\} \cup$$

$$\left\{ \begin{pmatrix} a & b & c \\ d & e & f \\ g & h & i \end{pmatrix} \middle| a,b,c,d,e,f,i \in Z_{47}I \right\}$$

be a neutrosophic 4-linear algebra over the 4-field $F = Q \cup Z_7 \cup Z_{19} \cup Z_{47}$. It is easily verified that V has no pseudo 4-linear subalgebra over F. Thus V is a pseudo simple 4-linear algebra.

**DEFINITION 3.3.10:** *Let $V = V_1 \cup V_2 \cup ... \cup V_n$ be a neutrosophic n-vector space (n-linear algebra) over the n-field $F = F_1 \cup F_2 \cup ... \cup F_n$. Suppose $W = W_1 \cup W_2 \cup ... \cup W_n \subseteq V_1 \cup V_2 \cup ... \cup V_n$ ($W_i \neq V_i$, $W_i \neq \{0\}$ for $i = 1, 2, ..., n$) is a neutrosophic n-vector space (n-linear algebra) over the n-subfield $K = K_1 \cup K_2 \cup ... \cup K_n \subseteq F_1 \cup F_2 \cup ... \cup F_n$ (Each $K_i$ is a proper subfield of $F_i$; $i = 1, 2, ..., n$) the we call W to be a subneutrosophic n-vector subspace (n-linear subalgebra) of V over the n-subfield $K \subset F$ of type II. If V has no proper subneutrosophic n-vector subspace (n-linear subalgebra) over a proper n-subfield of F then we call V to be a n-simple subneutrosophic n-vector space (n-linear algebra).*

Now we can define strong neutrosophic n-vector spaces of type II and derive for them also some interesting properties with appropriate modifications.

**DEFINITION 3.3.11:** *Let $V = V_1 \cup V_2 \cup ... \cup V_n$ be such that each $V_i$ is a strong neutrosophic vector space over the neutrosophic field $F_i$, $i = 1, 2, ..., n$ then we call $V = V_1 \cup V_2 \cup ... \cup V_n$ to be a strong neutrosophic n-vector space over the neutrosophic n-field $F = F_1 \cup F_2 \cup ... \cup F_n$.*

We will illustrate them by some examples.

***Example 3.3.21:*** Let $V = V_1 \cup V_2 \cup V_3 \cup V_4 =$
$$\left\{ \begin{pmatrix} a & a & a & a & a \\ b & b & b & b & b \end{pmatrix} \middle| a,b \in N(Z_{11}) \right\} \cup$$

$$\left\{ \begin{pmatrix} a & a \\ b & b \\ c & c \\ d & d \end{pmatrix} \middle| a,b,c,d \in N(Z_7) \right\} \cup$$

$$\left\{ \begin{pmatrix} 0 & 0 & 0 & a_1 \\ 0 & 0 & a_2 & 0 \\ 0 & a_3 & 0 & 0 \\ a_4 & 0 & 0 & 0 \end{pmatrix} \middle| a_i \in N(Z_{47}); 1 \leq i \leq 4 \right\} \cup$$

$$\left\{ \begin{pmatrix} a & b & c & d & e & g \\ a & b & c & d & e & g \\ a & b & c & d & e & g \end{pmatrix} \middle| a,b,c,d,e,g \in N(Z_{19}) \right\}$$

be a strong neutrosophic 4-vector space over the 4 neutrosophic field $F = Z_{11}I \cup N(Z_7) \cup Z_{47}I \cup N(Z_{19})$.

*Example 3.3.22:* Let $V = V_1 \cup V_2 \cup V_3 =$

$$\left\{ \begin{pmatrix} a \\ b \\ c \\ d \\ e \\ f \end{pmatrix} \middle| a,b,c,d,e,f \in N(Q) \right\} \cup$$

$$\left\{ \begin{pmatrix} a & b & a & b & a & b \\ b & c & b & c & b & c \\ c & d & c & d & c & d \end{pmatrix} \middle| a,b,c,d \in N(Z_{41}) \right\} \cup$$

$$\left\{ \begin{pmatrix} a & b & c \\ d & e & f \\ a & b & c \\ d & e & f \\ a & a & a \end{pmatrix} \middle| a,b,c,d,e,f \in Z_{11}I \right\}$$

be a strong neutrosophic trivector space over the neutrosophic trifield $F = N(Q) \cup Z_{41}I \cup Z_{11}I$.

We can define the notion of strong neutrosophic n-linear algebra.

**DEFINITION 3.3.12:** *Let $V = V_1 \cup V_2 \cup ... \cup V_n$ be a strong neutrosophic n-vector space over the neutrosophic n-field $F = F_1 \cup F_2 \cup ... \cup F_n$. If each $V_i$ is strong neutrosophic linear algebra over the neutrosophic field $F_i$, $i = 1, 2, ..., n$ then we call V to be a strong neutrosophic n-linear algebra over the neutrosophic n-field $F = F_1 \cup F_2 \cup ... \cup F_n$.*

We will illustrate this situation by some examples.

*Example 3.3.23:* Let $V = V_1 \cup V_2 \cup V_3 \cup V_4 \cup V_5 =$

$$\left\{ \begin{pmatrix} a & b \\ c & d \end{pmatrix} \middle| a,b,c,d \in N(Q) \right\} \cup$$

$\left\{ \sum_{i=0}^{\infty} a_i x^i \right.$ ; all neutrosophic polynomials in the variable x with coefficients from the neutrosophic field $Z_{11}I$; $a_i \in Z_{11}I$; $0 \leq i \leq \infty \} \cup \{(a_1, a_2, a_3, .... a_{20}) \mid a_i \in Z_3I; 1 \leq i \leq 20\} \cup$

$$\left\{ \begin{pmatrix} a_1 & 0 & 0 \\ a_2 & a_3 & 0 \\ a_4 & a_5 & a_6 \end{pmatrix} \middle| a_i \in Z_2I; 1 \leq i \leq 6 \right\} \cup$$

$$\left\{ \begin{pmatrix} a_1 & 0 & 0 & 0 & 0 \\ 0 & a_2 & 0 & 0 & 0 \\ 0 & 0 & a_3 & 0 & 0 \\ 0 & 0 & 0 & a_4 & 0 \\ 0 & 0 & 0 & 0 & a_5 \end{pmatrix} \middle| a_i \in Z_{17}I; 1 \le i \le 5 \right\}$$

be a strong neutrosophic 5-linear algebra over the neutrosophic 5-field $F = QI \cup Z_{11}I \cup Z_3I \cup Z_2I \cup Z_{17}I$ of type II.

***Example 3.3.24:*** Let $V = V_1 \cup V_2 \cup V_3 \cup V_4 = $ {All $10 \times 10$ upper triangular neutrosophic matrices with entries from the neutrosophic field $Z_2I$} $\cup$ {set of all $15 \times 15$ neutrosophic diagonal matrices with entries from the neutrosophic field $Z_3I$} $\cup$ {set of all $3 \times 3$ lower triangular matrices with entries from the neutrosophic field $N(Q)$} $\cup$ {$(a_1, a_2, a_3, a_4, a_5, a_6, a_7, a_8) \mid a_i \in Z_5I$} be a strong neutrosophic 4-linear algebra over the neutrosophic 4-field $F = Z_2I \cup Z_3I \cup N(Q) \cup Z_5I$.

We as in case of neutrosophic n-vector spaces (n-linear algebras) of type II define in case of strong neutrosophic n-vector spaces (n-linear algebras) the notion of strong neutrosophic n-vector subspaces (n-linear subalgebras) of type II.

Recall in $F = F_1 \cup F_2 \cup \ldots \cup F_n$ and if some of the $F_i$'s are neutrosophic fields and some of the $F_j$'s of real fields; $1 \le i, j \le n$ then we call F to be a quasi neutrosophic n-field. We shall just give some examples of them.

$$F = F_1 \cup F_2 \cup F_3 \cup F_4 \cup F_5 \cup F_6$$
$$= Z_7I \cup Z_2 \cup N(Q) \cup R \cup Z_{11}I \cup N(Z_{23})$$

is a quasi neutrosophic 6-field.

$$K = K_1 \cup K_2 \cup K_3 \cup K_4 \cup K_5 \cup K_6 \cup K_7 \cup K_8$$
$$= Q(\sqrt{2})I \cup Z_3I \cup N(Z_{29}) \cup Z_{17}I \cup Q(\sqrt{7}\sqrt{11}) \cup$$
$$N(Q(\sqrt{19}\sqrt{23}\sqrt{3})) \cup N(Z_{47}) \cup Z_{43}I$$

is a quasi neutrosophic 8-field.

Using the notion of quasi neutrosophic n-field we can define the new notion of quasi strong neutrosophic n-vector spaces (n-linear algebras) over quasi n-neutrosophic fields.

**DEFINITION 3.3.13:** *Let $V = V_1 \cup V_2 \cup ... \cup V_n$ be such that some $V_i$'s are vector spaces over the real field $F_i$ and some of the $V_j$'s are strong neutrosophic vector spaces over the neutrosophic field $F_j$ ($i \neq j$, $1 \leq i, j \leq n$). We define $V_1 \cup V_2 \cup ... \cup V_n$ to be a quasi strong neutrosophic n-vector space over the quasi neutrosophic n-field $F = F_1 \cup F_2 \cup ... \cup F_n$.*

We will illustrate this and the substructures mentioned earlier by some examples.

*Example 3.3.25:* Let $V = V_1 \cup V_2 \cup V_3 \cup V_4 \cup V_5 \cup V_6 =$

$$\left\{ \begin{pmatrix} a & b & a \\ c & d & b \end{pmatrix} \middle| a,b,c,d \in Z_5 I \right\} \cup$$

$$\left\{ \begin{pmatrix} a_1 \\ a_2 \\ a_3 \\ a_4 \\ a_5 \end{pmatrix} \middle| a_i \in Z_7 I \right\} \cup$$

$$\left\{ \begin{pmatrix} 0 & 0 & 0 & 0 & a_1 \\ 0 & 0 & 0 & a_2 & a_4 \\ 0 & 0 & a_3 & a_1 & a_6 \\ 0 & a_4 & a_6 & a_7 & a_8 \\ a_9 & a_1 & a_2 & a_3 & a_4 \end{pmatrix} \middle| a_i \in Z_{13} I; 1 \leq i \leq 6 \right\} \cup$$

$$\left\{ \begin{pmatrix} a_1 & a_2 & a_3 & a_4 \\ a_5 & a_6 & a_7 & a_8 \end{pmatrix} \middle| a_i \in N(Q); 1 \leq i \leq 8 \right\} \cup$$

$$\left\{\sum_{i=0}^{12} a_i x^i \text{ ; all neutrosophic polynomials in the variable x with}\right.$$

coefficients from the neutrosophic field $Z_{41}I$ of degree less than or equal to 12; $a_i \in Z_{41}I$; $0 \le i \le 12\} \cup$

$$\left\{\begin{pmatrix} a_1 & a_2 \\ a_3 & a_4 \\ a_5 & a_6 \\ a_7 & a_8 \\ a_9 & a_{10} \\ a_{11} & a_{12} \\ a_{13} & a_{14} \end{pmatrix} \middle| a_i \in Z_{11}I; 1 \le i \le 14 \right\}$$

be a quasi strong neutrosophic 6-vector space over the quasi neutrosophic 6-field $F = Z_5 \cup Z_7I \cup Z_{13} \cup N(Q) \cup Z_{41} \cup Z_{11}I$.

*Example 3.3.26:* Let $V = V_1 \cup V_2 \cup V_3 \cup V_4 =$

$$\left\{\begin{pmatrix} a_1 & a_2 & a_3 \\ a_4 & a_5 & a_6 \\ a_7 & a_8 & a_9 \\ a_{10} & a_{11} & a_{12} \\ a_{13} & a_{14} & a_{15} \\ a_{16} & a_{17} & a_{18} \\ a_{19} & a_{20} & a_{21} \end{pmatrix} \middle| a_i \in N(Q); 1 \le i \le 21 \right\} \cup$$

$$\left\{\begin{pmatrix} a & b & c & d & e & f & g \\ b & c & d & e & g & f & a \\ a & b & d & f & g & a & e \end{pmatrix} \middle| a,b,c,d,e,f,g \in N(Z_{41}) \right\} \cup$$

$$\left\{ \begin{pmatrix} 0 & 0 & 0 & a \\ 0 & 0 & b & 0 \\ 0 & c & 0 & 0 \\ d & 0 & 0 & 0 \end{pmatrix} \middle| a,b,c,d \in N(Z_{41}) \right\} \cup$$

$$\left\{ \begin{pmatrix} a & b & c & d & e \\ a & b & a & c & a \\ d & e & d & e & c \end{pmatrix} \middle| a,b,c,d,e \in N(Z_{17}I) \right\}$$

be a quasi strong neutrosophic 4-vector space over the quasi neutrosophic 4-field $F = Q \cup Z_5I \cup Z_41 \cup Z_{17}I$.

*Example 3.3.27:* Let $V = V_1 \cup V_2 \cup V_3 \cup V_4 \cup V_5 =$

$$\left\{ \begin{pmatrix} a & b \\ c & d \end{pmatrix} \middle| a,b,c,d \in N(Z_{11}) \right\} \cup$$

$\{(a_1, a_2, a_3, a_4, a_5, a_6, a_7) \mid a_i \in Z_{17}I; 1 \leq i \leq 7\} \cup \{$All $7 \times 7$ upper triangular matrices with entries from $N(Q)\} \cup \{$All $8 \times 8$ lower triangular neutrosophic matrices with entries from $N(Z_{23})\} \cup$ $\left\{ \sum_{i=0}^{\infty} a_i x^i \right.$ ; all polynomials in the variable x with coefficients from the neutrosophic field $N(Z_{41})$; $a_i \in N(Z_{41})$; $0 \leq i \leq \infty\}$ be a strong quasi neutrosophic 5-linear algebra over the quasi neutrosophic is field $F = Z_{11} \cup Z_{17}I \cup N(Q) \cup Z_{23} \cup N(Z_{41})$.

*Example 3.3.28:* Let $V = V_1 \cup V_2 \cup V_3 \cup V_4 = \{$All 10×10 neutrosophic matrices with entries from $N(Z_{53})\} \cup$

$$\left\{ \begin{pmatrix} a_1 & 0 & 0 \\ a_2 & a_3 & 0 \\ a_1 & a_2 & a_3 \end{pmatrix} \middle| a_i \in Z_{11}I; 1 \leq i \leq 3 \right\} \cup$$

{All 5 × 5 diagonal matrices with entries from $N(Z_{41})$} $\cup$

$$\left\{ \begin{pmatrix} a & b & 0 & 0 \\ c & d & 0 & 0 \\ 0 & 0 & b & a \\ 0 & 0 & d & c \end{pmatrix} \middle| a,b,c,d \in N(Z_2) \right\}$$

be a strong quasi neutrosophic 4-linear algebra over the quasi neutrosophic 4-field $F = N(Z_{53}) \cup Z_{11} \cup N(Z_{41}) \cup Z_2$.

***Example 3.3.29:*** Let $V = V_1 \cup V_2 \cup V_3 \cup V_4 \cup V_5 =$

$$\left\{ \begin{pmatrix} a_1 & a_2 & a_3 \\ a_4 & a_5 & a_6 \end{pmatrix} \middle| a_i \in N(Z_3); 1 \le i \le 6 \right\} \cup$$

$\left\{ \sum_{i=0}^{15} a_i x^i \right.$; all neutrosophic polynomials in the variable x of degree less than or equal to fifteen with coefficients from $N(Q)$; $a_i \in N(Q); 0 \le i \le 15 \}$ $\cup$ {All 5 × 5 neutrosophic matrices with entries from $Z_7 I$} $\cup$

$$\left\{ \begin{pmatrix} a & b \\ c & d \\ e & f \\ g & h \\ i & j \end{pmatrix} \middle| a,b,c,d,e,f,g,h,i,j \in Z_{11} I \right\} \cup$$

$$\left\{ \begin{pmatrix} 0 & 0 & 0 & a \\ 0 & 0 & b & c \\ 0 & c & d & a \\ a & b & c & d \end{pmatrix} \middle| a,b,c,d \in Z_{17} I \right\}$$

be a strong neutrosophic 5-vector space over the 5-neutrosophic field $F = Z_3I \cup N(Q) \cup Z_7I \cup Z_{11}I \cup Z_{17}I$. $W = W_1 \cup W_2 \cup W_3 \cup W_4 \cup W_5 =$

$$\left\{ \begin{pmatrix} a & a & a \\ b & b & b \end{pmatrix} \middle| a, b \in Z_3I \right\} \cup$$

$\left\{ \sum_{i=0}^{6} a_i x^i \right.$ ; all polynomial in the variable x with coefficients from the field $N(Q)$ of degree less than or equal to 6; $0 \leq i \leq 6$; $a_i \in N(Q)\} \cup \{$All $5 \times 5$ upper triangular matrices with entries from the field $Z_7I\} \cup$

$$\left\{ \begin{pmatrix} a & a \\ a & a \\ a & a \\ a & a \\ b & b \end{pmatrix} \middle| a, b \in Z_{11}I \right\} \cup$$

$$\left\{ \begin{pmatrix} 0 & 0 & 0 & a \\ 0 & 0 & b & 0 \\ 0 & d & 0 & 0 \\ e & 0 & 0 & 0 \end{pmatrix} \middle| a, b, d, e \in Z_{17}I \right\}$$

$\subseteq V_1 \cup V_2 \cup V_3 \cup V_4 \cup V_5$ is a strong neutrosophic 5-vector subspace over the 5-neutrosophic field F.

*Example 3.3.30:* Let $V = V_1 \cup V_2 \cup V_3 \cup V_4 \cup V_5 \cup V_6 = \{(a_1, a_2, a_3, a_4, a_5, a_6) \mid a_i \in Z_7I; 1 \leq i \leq 6\} \cup \{$All polynomials in the variable x with coefficients from $N(Q)$; $N(Q)[x]\} \cup$

$$\left\{ \begin{pmatrix} a & b & c \\ d & e & f \\ g & h & k \end{pmatrix} \middle| a, b, c, d, e, f, g, h, k \in Z_{11}I \right\} \cup$$

{All 7 × 7 neutrosophic matrices with entries from $Z_5I$} ∪ {All 9 × 9 upper triangular matrices with entries from $N(Z_{23})$} ∪

$$\left\{ \begin{pmatrix} a & b & 0 & 0 \\ d & e & 0 & 0 \\ 0 & 0 & a & a \\ 0 & 0 & a & a \end{pmatrix} \middle| a, b, e, b \in Z_{29}I \right\}$$

be a strong neutrosophic 6-linear algebra over the neutrosophic 6-field $F = Z_7I \cup QI \cup Z_{11}I \cup Z_5I \cup N(Z_{23}) \cup Z_{29}I$. Let $W = W_1 \cup W_2 \cup W_3 \cup W_4 \cup W_5 \cup W_6 = \{(a, a, a, b, b, a) \mid a, b \in Z_7I\}$ ∪ {All polynomials in the variable x with coefficients from QI; that is QI[x]} ∪

$$\left\{ \begin{pmatrix} a & a & a \\ a & a & a \\ b & a & b \end{pmatrix} \middle| a, b \in Z_{11}I \right\} \cup$$

{7 × 7 neutrosophic upper triangular matrices with entries from $Z_5I$} ∪ {All 9 × 9 diagonal matrices with entries form $N(Z_{23})$}

$$\cup \left\{ \begin{pmatrix} a & a & 0 & 0 \\ a & a & 0 & 0 \\ 0 & 0 & a & a \\ 0 & 0 & a & a \end{pmatrix} \middle| a \in Z_{29}I \right\}$$

⊆ $V_1 \cup V_2 \cup V_3 \cup V_4 \cup V_5 \cup V_6$ is a strong neutrosophic 6-linear subalgebra of V over the neutrosophic 6-field F.

Thus we have given examples of strong neutrosophic n-vector subspaces (n-linear subalgebras) over a neutrosophic n-field F.

As in case of neutrosophic n-vector spaces (n-linear algebras) one can in case of strong neutrosophic n-vector spaces

(n-linear algebras) define the notion of n-linearly independent n-subset, n-basis, strong neutrosophic n-linear transformations from a strong neutrosophic n-vector space (n-linear algebra) into a strong neutrosophic n-vector space (n-linear algebra) W both defined on the same n-neutrosophic field F. We can also define and study the properties enjoyed by strong neutrosophic n-linear operators on a strong neutrosophic n-vector space (n-linear algebra) V. Concepts like strong neutrosophic eigen values, eigen vectors, etc of strong neutrosophic n-linear operators can be obtained after suitable modifications.

Further almost all theorems derived in case of strong neutrosophic bivector spaces can be derived for strong neutrosophic n-vector spaces. Hence we leave this task for the interested reader. Only in case of strong neutrosophic n-vector spaces one can define n-linear functions, strong neutrosophic dual n-vector spaces and prove $(V^*)^* = V$. That is if

$$\begin{aligned}
V^* &= V_1^* \cup V_2^* \cup \ldots \cup V_n^* \\
&= (V_1 \cup V_2 \cup \ldots \cup V_n)^* \\
(V^*)^* &= (V_1^* \cup V_2^* \cup \ldots \cup V_n^*)^* \\
&= (V_1^*)^* \cup (V_2^*)^* \cup \ldots \cup (V_n^*)^* \\
&= V_1 \cup V_2 \cup \ldots \cup V_n.
\end{aligned}$$

As $(V_i^*)^* = V_i$ for each i, i = 1, 2, …, n.

The reader is expected to prove the following theorem.

**THEOREM 3.3.3:** *Let $T = T_1 \cup T_2 \cup \ldots \cup T_n$ be a n-linear operator on a finite $(n_1, n_2, \ldots, n_n)$ dimension strong neutrosophic n-vector space $V = V_1 \cup V_2 \cup \ldots \cup V_n$ over the n-field $F = F_1 \cup F_2 \cup \ldots \cup F_n$. Let*
$$C = \{C_1^1, C_2^1, \ldots, C_{k_1}^1\} \cup \{C_1^2, C_2^2, \ldots, C_{k_2}^2\} \cup \ldots \cup \{C_1^n, C_2^n, \ldots, C_{k_n}^n\}$$
*be distinct n-characteristic values of $T = T_1 \cup T_2 \cup \ldots \cup T_n$ and let $W_i = W_{i_1}^1 \cup W_{i_2}^2 \cup \ldots \cup W_{i_n}^n$ be the null n-space (n-null space) of*
$$T - CI_d = [T_1 - C_{i_1}^1 I_{d_1}] \cup [T_2 - C_{i_2}^2 I_{d_2}] \cup \ldots \cup [T_n - C_{i_n}^n I_{d_n}];$$

*The following are equaivalent*
*(i) T is n-diagonalizable*
*(ii) The n-characeristic n-polynomial for $T = T_1 \cup T_2 \cup ... \cup T_n$ is*

$$f = f_1 \cup f_2 \cup ... \cup f_n$$
$$= (x - C_1^1)^{d_1^1} ... (x - C_{k_1}^1)^{d_{k_1}^1} \cup (x - C_1^2)^{d_1^2} ...(x - C_{k_2}^2)^{d_{k_2}^2}$$
$$\cup ... \cup (x - C_1^n)^{d_1^n} ... (x - C_{k_n}^n)^{d_{k_n}^n}.$$

We define the notion of the n-ideal generated by n-polynomials which n-annihilate $T = T_1 \cup T_2 \cup ... \cup T_n$.

**DEFINITION 3.3.14:** *Let $T = T_1 \cup T_2 \cup ... \cup T_n$ be a n-linear operator of the finite $(n_1, n_2, ..., n_n)$ dimensional strong neutrosophic n-vector space $V = V_1 \cup V_2 \cup ... \cup V_n$ over the neutrosophic n-field $F = F_1 \cup ... \cup F_n$. The n-minimal neutrosophic n-polynomial for T is the unique monic n-generator of the n-ideal of n-polynomials over the n-field $F = F_1 \cup F_2 \cup ... \cup F_n$ which n-annihilate $T = T_1 \cup T_2 \cup ... \cup T_n$.*

*The n-minimal neutrosophic n-polynomial $p = p_1 \cup p_2 \cup ... \cup p_n$ for the n-linear operator $T = T_1 \cup T_2 \cup ... \cup T_n$ is uniquely determined by the following properties.*

*i. p is a n-monic neutrosophic n-polynomial over the n-field $F = F_1 \cup F_2 \cup ... \cup F_n$.*
*ii. $p(T) = p_1(T_1) \cup p_2(T_2) \cup ... \cup p_n(T_n) = 0 \cup 0 \cup ... \cup 0$.*
*iii. No neutrosophic n-polynomial over the n-field $F = F_1 \cup F_2 \cup ... \cup F_n$ which n-annihilates $T = T_1 \cup T_2 \cup ... \cup T_n$ has smaller n-degree than $p = p_1 \cup p_2 \cup ... \cup p_n$, has.*
*iv. $(n_1 \times n_1, n_2 \times n_2, ..., n_n \times n_n)$ to be the order of the neutrosophic n-matrix $A = A_1 \cup A_2 \cup ... \cup A_n$ over the n-field $F = F_1 \cup F_2 \cup ... \cup F_n$ (Each $A_i$ is the neutrosophic matrix of order $n_i \times n_i$ associated with $T_i$ with entries from the field $F_i$, $i = 1, 2, ..., n$).*

The reader to expected to derive these facts and also obtain all the related results like Cayley Hamilton Theorem, n-projection primary n-decomposition theorem, n-cyclic decomposition theorem, Generalized Cayley Hamilton theorem and so on. All these concepts can be extended appropriately from the results proved in case of strong neutrosophic bivector spaces over bifields.

The reader can define and derive n-Jordan form or Jordan n-form analogous to bi Jordan form or Jordan biform.

The notion of n-inner product on strong neutrosophic n-vector spaces of type II is an important and an interesting notion.

**DEFINITION 3.3.15:** *Let $F = F_1 \cup F_2 \cup F_3 \cup \ldots \cup F_n$ be a real neutrosophic n-field and $V = V_1 \cup V_2 \cup \ldots \cup V_n$ be a strong neutrosophic n-vector space over F. A n-inner product of V is a n-function which assigns to each n-ordered pair of n-vectors $\alpha = \alpha_1 \cup \alpha_2 \cup \ldots \cup \alpha_n$ and $\beta = \beta_1 \cup \beta_2 \cup \ldots \cup \beta_n$ in V a n-scalar $(\alpha/\beta) = (\alpha_1/\beta_1) \cup (\alpha_2/\beta_2) \cup \ldots \cup (\alpha_n/\beta_n)$ in $F = F_1 \cup F_2 \cup \ldots \cup F_n$; that is $(\alpha_i / \beta_i) \in F_i$; $i = 1, 2, \ldots, n$. $(\alpha_i, \beta_i \in V_i)$ in such a way that for all $\alpha = \alpha_1 \cup \alpha_2 \cup \ldots \cup \alpha_n$, $\beta = \beta_1 \cup \beta_2 \cup \ldots \cup \beta_n$ and $\gamma = \gamma_1 \cup \gamma_2 \cup \ldots \cup \gamma_n$ in $V = V_1 \cup V_2 \cup \ldots \cup V_n$ and for all n-scalars $c = c_1 \cup c_2 \cup \ldots \cup c_n$ in $F = F_1 \cup F_2 \cup \ldots \cup F_n$.*

(a) $(\alpha + \beta/\gamma) = (\alpha/\gamma) + (\beta/\gamma)$
  $(\alpha_1 + \beta_1/\gamma_1) \cup (\alpha_2 + \beta_2/\gamma_2) \cup \ldots \cup (\alpha_n + \beta_n/\gamma_n)$
  $= (\alpha_1/\gamma_1) + (\beta_1/\gamma_1) \cup (\alpha_2/\gamma_2) + (\beta_2/\gamma_2) \cup \ldots \cup (\alpha_n/\gamma_n) + (\beta_n/\gamma_n)$.

(b) $(c\alpha/\beta) = c(\alpha/\beta)$ that is
  $(c_1\alpha_1/\beta_1) \cup (c_2\alpha_2/\beta_2) \cup \ldots \cup (c_n\alpha_n/\beta_n)$
  $= c_1(\alpha_1/\beta_1) \cup c_2(\alpha_2/\beta_2) \cup \ldots \cup c_n(\alpha_n/\beta_n)$.

(c) $(\alpha/\beta) = (\beta/\alpha)$ that is
  $(\alpha_1 \cup \alpha_2 \cup \ldots \cup \alpha_n / \beta_1 \cup \beta_2 \cup \ldots \cup \beta_n) =$
  $(\beta_1 \cup \beta_2 \cup \ldots \cup \beta_n / \alpha_1 \cup \alpha_2 \cup \ldots \cup \alpha_n )$.

(d) $(\alpha/\alpha) = (\alpha_1 \cup \alpha_2 \cup \ldots \cup \alpha_n / \alpha_1 \cup \alpha_2 \cup \ldots \cup \alpha_n)$
  $= (\alpha_1/\alpha_1) \cup (\alpha_2/\alpha_2) \cup \ldots \cup (\alpha_n/\alpha_n) > (0 \cup 0 \cup \ldots \cup 0)$
  *if $\alpha_i \neq 0$ for $i = 1, 2, \ldots, n$.*

*A strong neutrosophic n-vector space endowed with a n-linear product is defined as a strong neutrosophic n-inner product space over the real neutrosophic n-field $F = F_1 \cup F_2 \cup ... \cup F_n$.*

*Let $V = F_1^{n_1} \cup F_2^{n_2} \cup ... \cup F_n^{n_n}$ be a strong neutrosophic n-vector space over the real neutrosophic n-field $F = F_1 \cup F_2 \cup ... \cup F_n$, there is a standard n-inner product called the n-standard inner product. It is defined for*

$$\alpha = \alpha_1 \cup \alpha_2 \cup ... \cup \alpha_n$$
$$= (x_1^1 ... x_{n_1}^1) \cup (x_1^2 ... x_{n_2}^2) \cup ... \cup (x_1^n ... x_{n_n}^n)$$

*and*

$$\beta = \beta_1 \cup \beta_2 \cup ... \cup \beta_n$$
$$= (y_1^1 ... y_{n_1}^1) \cup (y_1^2 ... y_{n_2}^2) \cup ... \cup (y_1^n ... y_{n_n}^n)$$

*by*

$$(\alpha/\beta) = \sum_{j_1} x_{j_1}^1 y_{j_1}^1 \cup \sum_{j_2} x_{j_2}^2 y_{j_2}^2 \cup ... \cup \sum_{j_n} x_{j_n}^n y_{j_n}^n.$$

*If $A = A_1 \cup A_2 \cup ... \cup A_n$ is a neutrosophic n-matrix over the n-field $F = F_1 \cup F_2 \cup ... \cup F_n$, where $A_i \in F_i^{n_i \times n_i}$, $i = 1, 2, ..., n$. $F_i^{n_i \times n_i}$ is a strong neutrosophic vector space over $F_i$; $i = 1, 2, ..., n$. $V = F_1^{n_1 \times n_1} \cup F_2^{n_2 \times n_2} \cup ... \cup F_n^{n_n \times n_n}$ is a strong neutrosophic n-vector space over the n-field $F = F_1 \cup F_2 \cup ... \cup F_n$ and $V = F_1^{n_1 \times n_1} \cup F_2^{n_2 \times n_2} \cup ... \cup F_n^{n_n \times n_n}$ is a strong neutrosophic n-vector space over the n-field $F = F_1 \cup F_2 \cup ... \cup F_n$ is isomorphic to the strong neutrosophic n-vector space $F_1^{n_1^2} \cup F_2^{n_2^2} \cup ... \cup F_n^{n_n^2}$ in a natural way.*

$$(A/B) = \sum_{j_1 k_1} A_{j_1 k_1}^1 B_{j_1 k_1}^1 \cup \sum_{j_2 k_2} A_{j_2 k_2}^2 B_{j_2 k_2}^2 \cup ... \cup \sum_{j_n k_n} A_{j_n k_n}^n B_{j_n k_n}^n$$

*defines a n-inner product on V. A strong neutrosophic n-vector space over the neutrosophic n-field $F = F_1 \cup F_2 \cup ... \cup F_n$ on which is defined a n-linear product is known as the n-inner product neutrosophic space or neutrosophic n-inner product space.*

The reader is expected to prove the following theorem.

**THEOREM 3.3.4:** *If $V = V_1 \cup V_2 \cup \ldots \cup V_n$ be a n-inner product neutrosophic space then for any n-vectors $\alpha = \alpha_1 \cup \alpha_2 \cup \ldots \cup \alpha_n$ and $\beta = \beta_1 \cup \beta_2 \cup \ldots \cup \beta_n$ in V and any scalar $c = c_1 \cup c_2 \cup \ldots \cup c_n$ in $F = F_1 \cup F_2 \cup \ldots \cup F_n$*

i. $||c\alpha|| = |c| \, ||\alpha||$ that is
   $||c\alpha|| = ||c_1 \alpha_1|| \cup \ldots \cup ||c_n \alpha_n||$
   $= |c_1| \, ||\alpha_1|| \cup \ldots \cup |c_n| \, ||\alpha_n||$;

ii. $||\alpha|| > 0 \cup 0 \cup \ldots \cup 0$ for $\alpha \neq 0$,
    that is $||\alpha_1|| \cup ||\alpha_2|| \cup \ldots \cup ||\alpha_n|| > (0, 0, \ldots, 0)$
    $= 0 \cup 0 \cup \ldots \cup 0$

iii. $||(\alpha/\beta)|| < ||\alpha|| \, ||\beta||$ that is
     $||(\alpha_1/\beta_1)|| \cup ||(\alpha_2/\beta_2)|| \cup \ldots \cup ||(\alpha_n/\beta_n)||$
     $\leq ||\alpha_1|| \, ||\beta_1|| \cup ||\alpha_2|| \, ||\beta_2|| \cup \ldots \cup ||\alpha_n|| \, ||\beta_n||$.

As in case of strong neutrosophic bivector spaces we can define in case of strong neutrosophic n-vector spaces the notion of n-orthogornal n-vectors.

If $\alpha, \beta \in V = V_1 \cup V_2 \cup \ldots \cup V_n$ be a n-vector in a n-inner product space we can define

$$\gamma = \beta - \frac{(\beta/\alpha)}{||\alpha||^2}\alpha \; ;$$

that is

$$\gamma_1 \cup \gamma_2 \cup \ldots \cup \gamma_n$$

$$= \left(\beta_1 - \frac{(\beta_1/\alpha_1)}{||\alpha_1||^2}\alpha_1\right) \cup \left(\beta_2 - \frac{(\beta_2/\alpha_2)}{||\alpha_2||^2}\alpha_2\right)$$

$$\cup \ldots \cup \left(\beta_n - \frac{(\beta_n/\alpha_n)}{||\alpha_n||^2}\alpha_n\right).$$

We say $\alpha = \alpha_1 \cup \alpha_2 \cup \ldots \cup \alpha_n$ is n-orthogonal to $\beta = \beta_1 \cup \beta_2 \cup \ldots \cup \beta_n$ if

$$(\alpha \mid \beta) = (\alpha_1 \mid \beta_1) \cup (\alpha_2 \mid \beta_2) \cup \ldots \cup (\alpha_n \mid \beta_n)$$
$$= 0 \cup 0 \cup \ldots \cup 0.$$

This clearly implies $\beta = \beta_1 \cup \beta_2 \cup \ldots \cup \beta_n$ is n-orthogonal to $\alpha = \alpha_1 \cup \alpha_2 \cup \ldots \cup \alpha_n$.

The reader can easily prove the following theorem.

**THEOREM 3.3.5:** *A n-orthogonal n-set of non zero n-vectors is n-linearly independent.*

**THEOREM 3.3.6:** *Let $V = V_1 \cup V_2 \cup \ldots \cup V_n$ be a strong neutrosophic n-inner product space and let $(\beta_1^1, \ldots, \beta_{n_1}^1) \cup (\beta_1^2, \ldots, \beta_{n_2}^2) \cup \ldots \cup (\beta_1^n, \ldots, \beta_{n_n}^n)$ be any n-independent vectors in V. Then one way to construct n-orthogonal vectors $(\alpha_1^1, \ldots, \alpha_{n_1}^1) \cup (\alpha_1^2, \ldots, \alpha_{n_2}^2) \cup \ldots \cup (\alpha_1^n, \ldots, \alpha_{n_n}^n)$ in $V = V_1 \cup V_2 \cup \ldots \cup V_n$ is such that for each $k_i = 1, 2, \ldots, n_i$ the set $(\alpha_1^1, \ldots, \alpha_{k_1}^1) \cup (\alpha_1^2, \ldots, \alpha_{k_2}^2) \cup \ldots \cup (\alpha_1^n, \ldots, \alpha_{k_n}^n)$ is a n-basis for the neutrosophic n-vector subspace spanned by $(\beta_1^1, \ldots, \beta_{k_1}^1) \cup (\beta_1^2, \ldots, \beta_{k_2}^2) \cup \ldots \cup (\beta_1^n, \ldots, \beta_{k_n}^n)$.*

*Proof:* The n-vectors $(\alpha_1^1, \ldots, \alpha_{n_1}^1) \cup (\alpha_1^2, \ldots, \alpha_{n_2}^2) \cup \ldots \cup (\alpha_1^n, \ldots, \alpha_{n_n}^n)$ will be obtained by means of a construction known as Gram-Schmidt n-orthogonalization process.

First let $\alpha_1 = \alpha_1^1 \cup \alpha_1^2 \cup \ldots \cup \alpha_1^n = \beta_1 = \beta_1^1 \cup \beta_1^2 \cup \ldots \cup \beta_1^n$. The other n-vector are given inductively as follows.

Suppose $\alpha_1, \alpha_2, \ldots, \alpha_{m_i}$ ($1 \le m_i < n_i$) have been choosen so that for every $k_i$, $\{\alpha_1, \alpha_2, \ldots, \alpha_{k_i}\}$; $1 \le k_i \le m_i$ is an orthogonal n-basis for the n-subspace of V that is spanned by $\beta_1 \ldots \beta_k$. To construct the next n-vector $\alpha_{m_i+1}$ let

$$\alpha_{m_i+1} = \alpha^1_{m_i+1} \cup \alpha^2_{m_i+1} \cup \ldots \cup \alpha^n_{m_i+1}$$

$$\beta_{m_i+1} = \sum_{k=1}^{m} \frac{(\beta_{m+1}/\alpha_k)\alpha_k}{\|\alpha_k\|^2}$$

$$= \beta^1_{m_{i_1}+1} \cup \beta^2_{m_{i_2}+1} \cup \ldots \cup \beta^n_{m_{i_n}+1}$$

$$-\sum_{k_{i_1}} \frac{(\beta^1_{m_{i_1}+1}/\alpha^1_{k_{i_1}})}{\|\alpha^1_{k_{i_1}}\|^2}\alpha^1_{k_{i_1}} \cup \sum_{k_{i_2}} \frac{(\beta^2_{m_{i_2}+1}/\alpha^2_{k_{i_2}})}{\|\alpha^2_{k_{i_2}}\|^2}\alpha^2_{k_{i_2}} \cup \ldots \cup$$

$$\sum_{k_{i_n}} \frac{(\beta^n_{m_{i_n}+1}/\alpha^n_{k_{i_n}})}{\|\alpha^n_{k_{i_n}}\|^2}\alpha^n_{k_{i_n}} = \left(\beta^1_{m_{i_1}+1} - \sum \frac{(\beta^1_{m_{i_1}+1}/\alpha^1_{k_{i_1}})}{\|\alpha^1_{k_{i_1}}\|^2}\alpha^1_{k_{i_1}}\right) \cup \ldots \cup$$

$$\left(\beta^n_{m_{i_n}+1} - \sum \frac{(\beta^n_{m_{i_n}+1}/\alpha^n_{k_{i_n}})}{\|\alpha^n_{k_{i_n}}\|^2}\alpha^n_{k_{i_n}}\right).$$

Then $\alpha_{m_i+1} \neq 0$ i.e., $\alpha^i_{m+1} \neq 0$; for other wise $\beta^i_{m+1}$ is a linear combination of $\alpha^i_1 \cup \alpha^i_2 \cup \ldots \cup \alpha^i_{m+1}$; $i = 1, 2, \ldots, n$. Further more it can be verified $(\alpha^i_{m_i+1}/\alpha^i_j) = 0$; $1 \leq j \leq m_i$ true for $i = 1, 2, \ldots, n$.

Hence true for $(\alpha_{m+1}/\alpha_j) = 0 \cup \ldots \cup 0$.

Therefore

$\{\alpha^1_1, \ldots, \alpha^1_{m_1+1}\} \cup \{\alpha^2_1, \ldots, \alpha^2_{m_2+1}\} \cup \ldots \cup \{\alpha^n_1, \ldots, \alpha^n_{m_n+1}\}$

is an n-orthogonal set consisting $(m_1 + 1, m_2 + 1, \ldots, m_n + 1)$ non zero n-vectors in the n-subspace spanned by

$\{\beta^1_1, \ldots, \beta^1_{m_1+1}\} \cup \{\beta^2_1, \ldots, \beta^2_{m_2+1}\} \cup \ldots \cup \{\beta^n_1, \ldots, \beta^n_{m_n+1}\}$.

In particular for $m = 4$ we have

$$\alpha^1_1 \cup \alpha^2_1 \cup \ldots \cup \alpha^n_1 = \beta^1_1 \cup \beta^2_1 \cup \ldots \cup \beta^n_1.$$

$$\alpha_2^1 \cup \alpha_2^2 \cup \ldots \cup \alpha_2^n = \beta_2^1 \cup \beta_2^2 \cup \ldots \cup \beta_2^n - \frac{\left(\beta_2^1 / \alpha_1^1\right)}{\|\alpha_1^1\|^2} \alpha_1^1 \cup$$

$$\frac{\left(\beta_2^2 / \alpha_1^2\right)}{\|\alpha_1^2\|^2} \alpha_1^2 \cup \ldots \cup \frac{\left(\beta_2^n / \alpha_1^n\right)}{\|\alpha_1^n\|^2} \alpha_1^n$$

$$= \beta_2^1 - \frac{\left(\beta_2^1 / \alpha_1^1\right)}{\|\alpha_1^1\|^2} \alpha_1^1 \cup \beta_2^2 - \frac{\left(\beta_2^2 / \alpha_1^2\right)}{\|\alpha_1^2\|^2} \alpha_1^2 \cup \ldots \cup \beta_2^n - \frac{\left(\beta_2^n / \alpha_1^n\right)}{\|\alpha_1^n\|^2} \alpha_1^n$$

and so on.

Interested reader on similar lines can construct $\alpha_3 = \alpha_3^1 \cup \alpha_3^2 \cup \ldots \cup \alpha_3^n$ interms of $\beta_3$, $\alpha_2$ and $\alpha_1$ and so on.

Now we define the notion of best n-approximation in case of strong neutrosophic n-vector spaces over the neutrosophic n-field F.

**DEFINITION 3.3.16:** *Let $V = V_1 \cup V_2 \cup \ldots \cup V_n$ be a strong neutrosophic n-vector space over the neutrosophic n-field $F = F_1 \cup F_2 \cup \ldots \cup F_n$ (All $F_i$'s are not pure neutrosophic for $i = 1, 2, \ldots, n$) of type II.*

*Let $W = W_1 \cup W_2 \cup \ldots \cup W_n$ be a strong neutrosophic n-vector subspace of V over F. Let $\beta = \beta_1 \cup \beta_2 \cup \ldots \cup \beta_n$ be a n-vector in $V = V_1 \cup V_2 \cup \ldots \cup V_n$, $\beta \notin W$ (i.e., $\beta_i \notin W_i$ for $i = 1, 2, \ldots, n$).*
*To find the n best approximation (best n-approximation) to $\beta = \beta_1 \cup \beta_2 \cup \ldots \cup \beta_n$ in $W = W_1 \cup W_2 \cup \ldots \cup W_n$.*
*That is to find a n-vector $\alpha = \alpha_1 \cup \alpha_2 \cup \ldots \cup \alpha_n$ for which*

$$\|\beta-\alpha\| = \|\beta_1 - \alpha_1\| \cup \|\beta_2 - \alpha_2\| \cup \ldots \cup \|\beta_n - \alpha_n\|$$

*is as small as possible subject to the restriction $\alpha = \alpha_1 \cup \alpha_2 \cup \ldots \cup \alpha_n$ is in $W = W_1 \cup W_2 \cup \ldots \cup W_n$ (that is each $\|\beta_i - \alpha_i\|$ is* *as small as possible subject to the restriction that $\alpha_i$ should belong to $W_i$, $i = 1, 2, ..., n$).*

*To be more precise a best n-approximation to $\beta = \beta_1 \cup \beta_2 \cup ... \cup \beta_n$ in $W = W_1 \cup W_2 \cup ... \cup W_n$ is a n-vector $\alpha = \alpha_1 \cup \alpha_2 \cup ... \cup \alpha_n$ in $W$ such that*

$$||\beta - \alpha|| < ||\beta - \gamma||$$

*that is*

$$||\beta_1 - \alpha_1|| \cup ||\beta_2 - \alpha_2|| \cup ... \cup ||\beta_n - \alpha_n||$$
$$\leq ||\beta_1 - \gamma_1|| \cup ||\beta_2 - \gamma_2|| \cup ... \cup ||\beta_n - \gamma_n||$$

*for every n-vector $\gamma = \gamma_1 \cup \gamma_2 \cup ... \cup \gamma_n$ in $W$.*

It is important to note that as in case of n-vector spaces of type II, we as in case of strong neutrosophic n-vector spaces define the notion of n-orthogonality.

However the interested reader can prove the following theorem.

**THEOREM 3.3.7:** *Let $W = W_1 \cup W_2 \cup ... \cup W_n$ be a strong neutrosophic n-vector subspace of a strong neutrosophic n-inner product space $V = V_1 \cup V_2 \cup ... \cup V_n$. Let $\beta = \beta_1 \cup \beta_2 \cup ... \cup \beta_n \in V = V_1 \cup V_2 \cup ... \cup V_n$; $\beta_i \in V_i$; $i = 1, 2, ..., n$.*

i. *The n-vector $\alpha = \alpha_1 \cup \alpha_2 \cup ... \cup \alpha_n$ in $W$ is a best n-approximation to $\beta = \beta_1 \cup \beta_2 \cup ... \cup \beta_n$ by n-vectors in $W = W_1 \cup W_2 \cup ... \cup W_n$ if and only if $\beta - \alpha = (\beta_1 - \alpha_1) \cup (\beta_2 - \alpha_2) \cup ... \cup (\beta_n - \alpha_n)$ is n-orthogonal to every n-vector in $W$.*

   *That is each $\beta_i - \alpha_i$ is orthogonal to every vector in $W_i$; true for $i = 1, 2, ..., n$.*

ii. *If the best n-approximation to $\beta = \beta_1 \cup \beta_2 \cup ... \cup \beta_n$ in $W = W_1 \cup W_2 \cup ... \cup W_n$ exists then it is unique.*

Now we proceed onto define the notion of n-orthogonal complement of a n-set of n-vectors in $V = V_1 \cup V_2 \cup ... \cup V_n$.

**DEFINITION 3.3.17:** *Let $V = V_1 \cup V_2 \cup ... \cup V_n$ be a strong neutrosophic inner produce space of type II defined over the neutrosophic n-field $F = F_1 \cup F_2 \cup ... \cup F_n$.*

*Let $S = S_1 \cup S_2 \cup ... \cup S_n$ be any n-set of n-vectors in $V = V_1 \cup V_2 \cup ... \cup V_n$. The n-orthogonal complement of S denoted by $S^\perp = S_1^\perp \cup S_2^\perp \cup ... \cup S_n^\perp$ is the set of all n-vectors in V which are n-orthogonal to every n-vector in S.*

The reader is expected to prove the following theorems.

**THEOREM 3.3.8:** *Let $V = V_1 \cup V_2 \cup ... \cup V_n$ be a strong neutrosophic n-inner product space. Let $W = W_1 \cup W_2 \cup ... \cup W_n$ be a finite dimensional strong neutrosophic n-vector subspace of V and $E = E_1 \cup E_2 \cup ... \cup E_n$ be the n-orthogonal projection of V on W.*

*Then the n-mapping $\beta \to (\beta - E\beta)$; that is*
*$\beta_1 \cup \beta_2 \cup ... \cup \beta_n \to (\beta_1 - E\beta_1) \cup (\beta_2 - E\beta_2) \cup ... \cup (\beta_n - E\beta_n)$*
*i.e., each $\beta_i \to (\beta_i - E\beta_i)$ for $i = 1, 2, ..., n$ is the n-orthogonal projection of V on W.*

**THEOREM 3.3.9:** *Let $W = W_1 \cup W_2 \cup ... \cup W_n$ be a finite $(n_1, n_2, ..., n_n)$ n-dimensional strong neutrosophic n-vector subspace of the strong neutrosophic n-inner product space $V = V_1 \cup V_2 \cup ... \cup V_n$ of type II.*

*Let $E = E_1 \cup E_2 \cup ... \cup E_n$ be an n-idempotent n-linear transformation of V onto W. $W^\perp = W_1^\perp \cup W_2^\perp \cup ... \cup W_n^\perp$ is the null n-subspace of $E = E_1 \cup E_2 \cup ... \cup E_n$ and $V = W \oplus W^\perp$ that is*

$$V = V_1 \cup V_2 \cup ... \cup V_n$$
$$= W_1 \oplus W_1^\perp \cup ... \cup W_n \oplus W_n^\perp.$$

**THEOREM 3.3.10:** *Under the conditions of the above theorem $I - E = I_1 \cup I_2 \cup ... \cup I_n - (E_1 \cup E_2 \cup ... \cup E_n) = I_1 - E_1 \cup I_2 - E_2 \cup ... \cup I_n - E_n$ is the n-orthogonal n-projection of V on $W^\perp$. It is a n-idempotent n-linear transformation of V on to $W^\perp$ with n-null space $W = W_1 \cup W_2 \cup ... \cup W_n$.*

**THEOREM 3.3.11:** *Let*

$$\{\alpha_1^1,...,\alpha_{n_1}^1\} \cup \{\alpha_1^2,...,\alpha_{n_2}^2\} \cup ... \cup \{\alpha_1^n,...,\alpha_{n_n}^n\}$$

*be a n-orthogonal set of non zero n-vectors in a strong n-inner product space* $V = V_1 \cup V_2 \cup ... \cup V_n$ *over* $F = F_1 \cup F_2 \cup ... \cup F_n$ *of type II.*

*If* $\beta = \beta_1 \cup \beta_2 \cup ... \cup \beta_n$ *is any n-vector in* $V = V_1 \cup V_2 \cup ... \cup V_n$ *then*

$$\sum_{k_1}\left(\frac{|(\beta_1/\alpha_{k_1}^1)|^2}{||\alpha_{k_1}^1||^2}\right) \cup$$

$$\sum_{k_2}\left(\frac{|(\beta_2/\alpha_{k_2}^2)|^2}{||\alpha_{k_2}^2||^2}\right) \cup ... \cup$$

$$\sum_{k_n}\left(\frac{|(\beta_n/\alpha_{k_n}^2)|^2}{||\alpha_{k_n}^2||^2}\right)$$

$$\leq ||\beta_1||^2 \cup ||\beta_2||^2 \cup ... \cup ||\beta_n||^2$$

*and equality holds if and only if*

$$\beta_1 \cup \beta_2 \cup ... \cup \beta_n$$

$$= \sum_{k_1}\left(\frac{(\beta_1/\alpha_{k_1}^1)}{||\alpha_{k_1}^1||^2}\alpha_{k_1}^1\right) \cup$$

$$\sum_{k_2}\left(\frac{(\beta_2/\alpha_{k_2}^2)}{||\alpha_{k_2}^2||^2}\alpha_{k_2}^2\right) \cup ... \cup$$

$$\sum_{k_n} \left( \frac{(\beta_n / \alpha_{k_n}^2)}{\| \alpha_{k_n}^2 \|^2} \alpha_{k_n}^2 \right).$$

Results on neutrosophic bivector spaces (bilinear algebras) discussed and derived in Chapter 2 can be derived for neutrosophic n-vector spaces (n-linear algebras).

Further all results true in case of n-linear algebras of type II can be derived in case of neutrosophic n-linear algebras of type II with appropriate modifications.

**Chapter Four**

# APPLICATIONS OF NEUTROSOPHIC n-LINEAR ALGEBRAS

In this chapter we just suggest the possible applications of the neutrosophic n-linear algebras of type I and II (n ≥ 2), strong neutrosophic n-linear algebras of type I and II and quasi strong neutrosophic linear algebras of type II.

These neutrosophic n-linear algebras over the neutrosophic n-fields or over the real n-field can be used in neutrosophic fuzzy models like Neutrosophic Cognitive Maps (NCMs) and when we have multiexperts we can use the neutrosophic n-matrices and model n- NCMs (n ≥ 2).

These neutrosophic n-matrices can also be used to model Neutrosophic Fuzzy Relational maps, when n-experts give their opinion on any real world problem. Use of these neutrosophic n-matrices will save time and economy.

These neutrosophic n-matrices can be used in n-models whenever the concept of indeterminacy is present.

The n-NCMs (i.e., NCMs constructed using neutrosophic n-matrices which gives the NCM model of n-experts $n \geq 2$) can be used in creating metabolic regulatory n-Network models. Also multi expert NCM models can be used to find the driving speed vehicles of any one in free way. These n-NCM models will be very useful in Medical diagnostics. n-NCMs using neutrosophic n-matrices with entries from [0, 1] can be used in diagnosis and study of specific language impairment.

These structures will be best suited for web mining n-inferences and in robotics.

The strong neutrosophic n-linear operators when analyzing the eigen values and eigen vectors in any real models where indeterminacy is dominant can be used.

We see pure complex numbers ni, (i is such that $i^2 = -1$ and $n \in R$) at an even stages (powers) merges with real but when we use the indeterminant 'I' they at no point of time merge with reals. Thus if the presence of indeterminacy prevails use of neutrosophic models is more appropriate.

The n-NCM models can be used in legal rules when several lawyers give their opinion about a case.

These models will also be better suited in analysis of Business Performance Assessment as in business always a factor of indeterminacy is present.

Chapter Five

# SUGGESTED PROBLEMS

In this chapter we have given over eighty problems for the reader to solve. This will help one to understand the concepts introduced in this book.

1. Let $V = V_1 \cup V_2 =$

$$\left\{ \begin{pmatrix} a & b \\ c & d \end{pmatrix} \middle| a,b,c,d \in Z_{13}I \right\} \cup \left\{ \begin{pmatrix} a_1 & a_2 & a_3 \\ a_4 & a_5 & a_6 \end{pmatrix} \right\}$$

where $a_i \in N(Z_{13})$; $1 \leq i \leq 6$} be a neutrosophic bivector space over the field $Z_{13}$.
   a. Find neutrosophic bivector subspaces of V.
   b. Find $NHom_{Z_{13}}$ (V, V).
   c. Can V have special neutrosophic bivector subspaces over a subfield of $Z_{13}$?

2. Let $V = V_1 \cup V_2$ be a neutrosophic bivector space over a real field F. Develop some interesting properties of V.

3. Let $V = V_1 \cup V_2$ and $W = W_1 \cup W_2$ be two neutrosophic bivector spaces over the field F. Find the algebraic structure of $NHom_F(V, W)$.

4. Let $V = V_1 \cup V_2 =$

$$\left\{ \begin{pmatrix} a & b & c \\ 0 & d & e \\ 0 & 0 & f \end{pmatrix} \middle| a,b,c,d,e,f \in Z_{17}I \right\} \cup$$

$$\left\{ \begin{pmatrix} a_1 & a_2 \\ a_3 & a_4 \\ a_5 & a_6 \\ a_7 & a_8 \end{pmatrix} \middle| a_i \in N(Z_{17}) \right\}$$

be a neutrosophic bivector space over the real field $Z_{17}$. Find atleast one neutrosophic linear bioperator on V which is inverible. Does there exist a neutrosophic linear bioperator T on V which is inverible but T has a non trivial kernel?

5. Let $V = V_1 \cup V_2 = \{Z_{11}I[x]\} \cup$

$$\left\{ \begin{pmatrix} a_1 & a_2 & a_3 & a_4 & a_5 \\ a_6 & a_7 & a_8 & a_9 & a_{10} \end{pmatrix} \middle| a_i \in Z_{17}I; 1 \le i \le 10 \right\}$$

be a neutrosophic bivector space over the field $Z_{11}$. Find a bibasis of V over $Z_{11}$. Define a bilinear bioperator T on V which is not inverible.

6. Let $V = V_1 \cup V_2$, $W = W_1 \cup W_2$ and $S = S_1 \cup S_2$ be three neutrosophic bivector spaces over the real field F. Suppose T: $V \to W$ is a neutrosophic bilinear bitransformation, $P : W \to$

S is a neutrosophic bilinear transformation. Will TP: V→ S be a neutrosophic bilinear transformation from V to S?(we define (TP) (v) = P(T(v)) | v ∈ V so T(v) ∈ W; P((T(v)) ∈ S).

7. Let $V = V_1 \cup V_2$ and $W = W_1 \cup W_2$ be two finite bidimensional neutrosophic bivector spaces over the real field F. Prove if $T = T_1 \cup T_2 : V \to W$ is a linear bitransformation then nullity T + rank T = bidim V, that is nullity $(T_1 \cup T_2)$ + rank $(T_1 \cup T_2)$ = dim $V_1 \cup$ dim $V_2$.

8. Define neutrosophic hyperbisubspace of V. Illustrate this concept by an example.

9. Does there exist a neutrosophic bivector space which has no neutrosophic hyper bisubspace? Justify your claim.

10. Can you characterize those neutrosophic bivector spaces which has no neutrosophic hyper bisubspace?

11. Let $V = V_1 \cup V_2 = \{Z_{17}I\} \cup \{Z_{17}I[x]\}$ be a neutrosophic bivector space over the field $Z_{17}$. Does V have neutrosophic hyper bisubspace?

12. Obtain some interesting properties about the special strong neutrosophic bivector spaces over the neutrosophic bifield $F = F_1 \cup F_2$.

13. Let $V = V_1 \cup V_2 =$

$$\left\{ \begin{pmatrix} a & b & c \\ d & e & f \end{pmatrix} \middle| a,b,c,d,e,f \in Z_{11}I \right\} \cup$$

$$\left\{ \begin{pmatrix} a_1 & a_2 \\ a_3 & a_4 \\ a_5 & a_6 \\ a_7 & a_8 \end{pmatrix} \middle| a_i \in N(Q); 1 \le i \le 8 \right\}$$

be a strong neutrosophic bivector space over the neutrosophic bifield $F = Z_{11}I \cup N(Q) = F_1 \cup F_2$.
   a. Find a bibasis of V over the bifield $F = F_1 \cup F_2$.
   b. Can V have proper strong neutrosophic bivector subspaces $W_1$, $W_2$ which are such that $W_1 \cong W_2$?
   c. Is V bisimple?
   d. Find a linear bioperator on V which is non biinvertible.
   e. Find a linear bioperator on V which is biinvertible.
   f. Can V have proper hyper bisubspace?

14. Let $V = V_1 \cup V_2 = $ {collection of all $7 \times 7$ neutrosophic matrices with entries from the neutrosophic field $Z_5I$} $\cup$ {Collection of all $9 \times 9$ neutrosophic matrices with entries from the neutrosophic field QI} be a strong neutrosophic bivector space over the neutrosophic bifield $F = F_1 \cup F_2 = Z_5I \cup QI$.
   a. Is V a strong neutrosophic bilinear algebra over F?
   b. Can V have proper strong neutrosophic bilinear subalgebras over $F = F_1 \cup F_2$?
   c. Find SNH (V, V).
   d. Define a linear bioperator T on V which is biinvertible.
   e. Find a bibasis for V.
   f. Can V have pseudo strong bivector subspaces?
   g. Find for a bilinear operator T on V the associated bicharacteristic values and bicharacteristic vectors.
   h. Is V bisimple?
   i. Can V have pseudo real bilinear subalgebras?

15. Obtain some interesting properties about SNH (V, W) where V and W are strong neutrosophic bivector spaces defined over the pure neutrosophic bifield $F = F_1 \cup F_2$.

16. Will every strong neutrosophic bivector space have strong neutrosophic hyper space?

17. Let $B = B_1 \cup B_2 = $

$$\begin{pmatrix} 3I & 0 & 1 \\ 0 & 5I & 4 \\ 0 & 2I & 2 \end{pmatrix} \cup \begin{pmatrix} I & 0 & 1 & 2 \\ 0 & 2I & I & 0 \\ 0 & 1+I & 0 & 2+I \\ 2I+1 & 0 & 0 & I \end{pmatrix}$$

be a neutrosophic bimatrix with entries from the neutrosophic bifield $F = F_1 \cup F_2 = N(Z_5) \cup N(Z_3)$. Find the bicharacteristic values associated with B. Find the bicharacteristic neutrosophic bipolynomial associated with it.

18. Let $V = V_1 \cup V_2 =$

$$\left\{ \begin{pmatrix} a_1 & a_2 \\ a_3 & a_4 \\ a_5 & a_6 \\ a_7 & a_8 \\ a_9 & a_{10} \end{pmatrix} \middle| a_i \in Z_{17}I; 1 \leq i \leq 10 \right\} \cup$$

$$\left\{ \begin{pmatrix} a_1 & a_2 & a_3 \\ a_4 & a_5 & a_6 \end{pmatrix} \middle| a_i \in Z_5I; 1 \leq i \leq 10 \right\}$$

be a strong neutrosophic bivector space over the neutrosophic bifield $F = F_1 \cup F_2 = Z_{17} I \cup Z_5 I$. $W = W_1 \cup W_2 = \{(a, b, c, d, e, f, g) \mid a, b, c, d, e, f, g \in Z_{17}I\} \cup$

$$\left\{ \begin{pmatrix} a_1 & a_2 & a_3 \\ a_4 & a_5 & a_6 \\ a_7 & a_8 & a_9 \\ a_{10} & a_{11} & a_{12} \end{pmatrix} \middle| a_i \in Z_5I; 1 \leq i \leq 12 \right\}$$

be a strong neutrosophic bivector space over the same neutrosophic bifield $F = F_1 \cup F_2 = Z_{17}I \cup Z_5I$.

a. Find SNH(V, W) = $L^2$(V, W). Is $L^2$(V, W) a strong neutrosophic bivector space over the bifield $Z_{17}I \cup Z_5I$?
b. Find a T : V → W such that T is not biinvertible.
c. Can T : V → W be such that, T is biinvertible?
d. Find a S : W → W such that S is biinvertible.
e. Find $L^2$(W, V).
f. Find a bibasis for V.
g. Find a bibasis for W.
h. What is the bidimension of V?
i. What is the bidimension of $L^2$(V, W) over $Z_{17}I \cup Z_5I$?
j. Is V and W bisimple strong neutrosophic bivector spaces?
k. Find SNH (V, V) = $L^2$(V, V).
l. Find SNH (W, W) = $L^2$(W, W).

19. Let A = $A_1 \cup A_2$ =

$$\begin{pmatrix} I & 0 & 3 \\ 7I & 6 & 2 \\ 0 & 0 & 4I \end{pmatrix} \cup \begin{pmatrix} 4I & 0 & 2 & I \\ 3I & I & 0 & 7 \\ 0 & I & 0 & 0 \\ 8 & 0 & 7I & 4I \end{pmatrix}$$

be a neutrosophic bimatrix with entries from the neutrosophic bifield $N(Z_{11}) \cup N(Q) = F_1 \cup F_2 = F$.
a. Find the bicharacteristic bipolynomial associated with A.
b. Find the bicharacteristic bieigen values of A.
c. Is A bidiagonalizable over the bifield $F = F_1 \cup F_2$?

20. Obtain some interesting properties about bipolynomial biideals.

21. Let V = $V_1 \cup V_2$ =

$$\left\{ \sum_{i=0}^{\infty} a_i x^{2i} \middle| a_i \in Z_7I \right\} \cup \left\{ \sum_{i=0}^{19} a_i x^{2i} \middle| a_i \in N(Q); 0 \le i \le 19 \right\}$$

be a strong neutrosophic bivector space over the neutrosophic bifield F = $Z_7I \cup N(Q)$. Define T : V → V and find the bicharacteristic values associated with T. Find the bidimension of V over $F = Z_7I \cup N(Q)$.

22. Let $V = V_1 \cup V_2 =$

$$\left\{ \begin{pmatrix} a & a & a & a & a \\ a & a & a & a & a \end{pmatrix} \middle| a \in N(R) \right\} \cup$$

$$\left\{ \begin{pmatrix} a & b \\ b & a \\ a & b \\ b & a \end{pmatrix} \middle| a, b \in Z_7I \right\}$$

be a strong neutrosophic bivector space over the neutrosophic bifield $F = F_1 \cup F_2 = N(R) \cup Z_7I$.
 a. Find a bibasis of V relative to the bifield $F = R(I) \cup Z_7I$.
 b. Find a bibasis of V relative to the bifield $K = K_1 \cup K_2 = RI \cup Z_7I$.
 c. Find the bibasis of V relative to the bifield $S = S_1 \cup S_2 = N(Q) \cup Z_7I$ where V is a strong neutrosophic bivector space defined over $S = S_1 \cup S_2$.
 d. Find the bibasis of V relative to the bifield $P = P_1 \cup P_2 = QI \cup Z_7I$.

Compare the bibasis of V when defined over 4 different fields and establish that the bidimension of a strong neutrosophic bivector space is dependent on the neutrosophic bifield over which the bispace is defined.

23. Does there exist a strong neutrosophic bivector space whose bidimension is independent of the neutrosophic bifield over which it is defined?

24. Let $V =$

$$\left\{ \begin{pmatrix} a & b \\ c & d \end{pmatrix} \middle| a, b, c, d \in Z_{11}I \right\} \cup$$

$$\left\{ \begin{pmatrix} a & b & c \\ 0 & d & e \\ 0 & 0 & f \end{pmatrix} \middle| a,b,c,d,e,f \in Z_7I \right\}$$

be a strong neutrosophic bilinear algebra defined over the neutrosophic bifield $F = Z_{11}I \cup Z_7I$.

a. Find a bibasis of V and the bidimension of V over F.
b. Find the bidimension of NL(V,V) = {all linear bioperators on V} over $F = Z_{11}I \cup Z_7I$.
c. Is V pseudo bisimple? Justify your claim.

25. Let $V = V_1 \cup V_2 = \{(a_1, a_2, a_3, a_4, a_5) \mid a_i \in N(Q)\} \cup$

$$\left\{ \begin{pmatrix} a_1 \\ a_2 \\ a_3 \\ a_4 \\ a_5 \\ a_6 \end{pmatrix} \middle| a_i \in Z_{11}I \right\}$$

be a strong neutrosophic bivector space over the neutrosophic bifield $F = F_1 \cup F_2 = QI \cup Z_{11}I$ and $W = W_1 \cup W_2 =$

$$\left\{ \begin{pmatrix} a & b & c \\ 0 & d & e \\ 0 & 0 & 0 \end{pmatrix} \middle| a,b,c,d,e \in QI \right\} \cup$$

$$\left\{ \begin{pmatrix} a_1 & a_2 & a_3 \\ a_4 & a_5 & a_6 \end{pmatrix} \middle| a_i \in Z_{11}I \right\}$$

be a strong neutrosophic bivector space over the same bifield $F = F_1 \cup F_2 = QI \cup Z_{11}I$.
Find $L^2(V, W)$ and $L^2(W, V)$.

Determine the bidimension of $L^2(V, W)$ and $L^2(W, V)$ over $F = QI \cup Z_{11}I$.

26. Let $V = V_1 \cup V_2 =$

$$\left\{ \begin{pmatrix} a & b \\ c & d \end{pmatrix} \middle| a, b, c, d \in QI \right\} \cup$$

$$\left\{ \sum_{i=0}^{9} a_i x^i \middle| a_i \in Z_{13}I; 0 \leq i \leq 9 \right\}$$

be a strong neutrosophic bivector space over the neutrosophic bifield $F = QI \cup Z_{13}I$. Find the bidimension of $L^2(V, V)$ over $F = QI \cup Z_{13}I$.

27. Let $V = V_1 \cup V_2 =$

$$\left\{ \begin{pmatrix} a & b \\ c & d \end{pmatrix} \middle| a, b, c, d \in QI \right\} \cup$$

$$\left\{ \sum_{i=0}^{9} a_i x^i \middle| a_i \in QI; 0 \leq i \leq 9 \right\}$$

be a neutrosophic space over the neutrosophic field $F = QI$.
a. Find a bibasis for V.
b. Find the bidimension of $L^2(V, V)$.

28. Let $V = V_1 \cup V_2 =$

$$\left\{ \begin{pmatrix} a & b \\ c & d \\ e & f \end{pmatrix} \middle| a, b, c, d, e, f \in Z_3I \right\} \cup$$

{(a, a, a, a, b, b, b, b) | a, b ∈ $Z_2I$} be a strong neutrosophic bivector space over the neutrosophic bifield $F = Z_3I \cup Z_2I$.
   a. Find a bibasis for V.
   b. What is the bidimension of V over F?
   c. Find a $T = T_1 \cup T_2 : V = V_1 \cup V_2 \to V_1 \cup V_2$ such that all the bicharacteristic values associated with T are distinct.
   d. Find $L^2(V, V)$.
   e. What is the bidimension of $L^2(V, V)$ over F?

29. Let $V = V_1 \cup V_2 = \{QI\} \cup \{Z_{11}I\}$ be a strong neutrosophic bilinear algebra defined over the neutrosophic bifield $F = QI \cup Z_{11}I$.
   a. What is the bidimension of V?
   b. Suppose $V = V_1 \cup V_2$ is only a strong neutrosophic bivector space over $F = QI \cup Z_{11}I$, what is the bidimension of V?
   c. Does the bidimension in general dependent on its algebraic structure?

30. Let $V = V_1 \cup V_2 = N(Q) \cup Z_{19}I$ be a strong neutrosophic bivector space over the neutrosophic bifield $F = QI \cup Z_{19}I$.
   a. What is bidimension of V?
   b. If $V = V_1 \cup V_2$ is realized as a strong neutrosophic bilinear algebra over the bifield $F = QI \cup Z_{19}I$. What is its bidimension?
   c. Compare and find if any difference exists.

31. Let $V = V_1 \cup V_2 = \{QI\} \cup \{(a, b, c) | a, b, c, \in Z_{17}I\}$ be neutrosophic bivector space over the real bifield $F = Q \cup Z_{17}$.
   a. Find a bibasis for V over F.
   b. Find the bidimension of V over F.
   c. Find $L^2(V, V)$.
   d. What is the bidimension of $L^2(V, V)$?

32. Let $V = V_1 \cup V_2 =$
$$\left\{ \sum_{i=0}^{\infty} a_i x^i \middle| a_i \in Z_{11}I; i = 0, 1, 2, \ldots, \infty \right\}$$

$$\cup \left\{ \sum_{i=0}^{\infty} a_i x^i \,\middle|\, a_i \in Z_{17} I; i = 0, 1, 2, \ldots, \infty \right\}$$

be a neutrosophic bivector space over the bifield $F = Z_{11} \cup Z_{17}$.
a. What is the bidimension of V over $F = Z_{11} \cup Z_{17}$?
b. Suppose the same V is defined to be a strong neutrosophic bivector space over the field $K = Z_{11}I \cup Z_{17}I$, what is the bidimension of V over the bifield $K = Z_{11}I \cup Z_{17}I$?
c. Find $L_F^2(V, V)$ and $L_k^2(V, V)$.
d. What is the bidimensions of $L_F^2(V, V)$?

33. Let $V = V_1 \cup V_2$ and $W = W_1 \cup W_2$ be two $(n_1, n_2)$ and $(m_1, m_2)$ bidimensional strong neutrosophic bivector spaces over the bifield $F = F_1 \cup F_2$. Let $C^*$ and $B^*$ be the dual bibasis of V and W of C and B respectively. If A is a neutrosophic bimatrix of $T = T_1 \cup T_2$, a bilinear transformation of V to W relative to the bibasis C and B and $T^t$ relative to $C^*$ and $B^*$ respectively. Obtain some interesting relations between T and $T^t$.

34. Obtain some interesting properties / results about bidiagonalizable bilinear operator.

35. Find the bipolynomial for the neutrosophic bimatrix $A = A_1 \cup A_2 =$

$$\begin{pmatrix} 3I & 0 & 4I & 1 \\ 0 & 7I & 3 & 0 \\ 2 & 1 & 4 & 3 \\ 10I & 0 & 9 & I \end{pmatrix} \cup \begin{pmatrix} I & 0 & 1 \\ 0 & I & 1 \\ 0 & 1 & I \end{pmatrix}$$

where A is defined over the neutrosophic bifield $F = F_1 \cup F_2 = N(Z_{11}) \cup N(Z_2)$.

36. Illustrate by an example that birank $T^t$ = birank T.

37. Let $V = V_1 \cup V_2 =$

$$\left\{ \begin{pmatrix} a & b \\ c & d \end{pmatrix} \middle| a,b,c,d \in Z_{11}I \right\} \cup$$

$$\left\{ \begin{pmatrix} a_1 & a_2 & a_3 \\ a_4 & a_5 & a_6 \end{pmatrix} \middle| a_i \in N(Q), 1 \leq i \leq 6 \right\}$$

be a strong neutrosophic bivector space over the bifield $F = F_1 \cup F_2$ and let $T = T_1 \cup T_2$ be a bilinear operator on V defined by $T = T_1 \cup T_2 : V = V_1 \cup V_2 \to V = V_1 \cup V_2$ where $T_1 : V_1 \to V_1$ and $T_2 : V_2 \to V_2$ such that

$$T_1 \begin{pmatrix} a & b \\ c & d \end{pmatrix} = \begin{pmatrix} a & b \\ 0 & d \end{pmatrix}$$

and

$$T_2 \begin{pmatrix} a_1 & a_2 & a_3 \\ a_4 & a_5 & a_6 \end{pmatrix} = \begin{pmatrix} a_1 & 0 & a_3 \\ 0 & a_5 & 0 \end{pmatrix}.$$

Find $T^t$. Is birank $T^t$ = birank T?

38. Let $V = V_1 \cup V_2 =$

$$\left\{ \begin{pmatrix} a & b & c \\ d & e & f \\ g & h & i \end{pmatrix} \middle| a,b,c,d,e,f,g,h,i \in Z_{29}I \right\} \cup$$

$$\left\{ \begin{pmatrix} a_1 & a_2 & a_3 & a_4 \\ a_5 & a_6 & a_7 & a_8 \end{pmatrix} \middle| a_i \in N(Z_2); 1 \leq i \leq 8 \right\}$$

and

$$W = \left\{ \begin{pmatrix} a_1 & a_2 \\ a_3 & a_4 \\ a_5 & a_6 \\ a_7 & a_8 \end{pmatrix} \middle| a_i \in Z_{29}I; 1 \leq i \leq 8 \right\} \cup$$

$$\left\{ \begin{pmatrix} a_1 & a_2 & a_3 \\ a_4 & a_5 & a_6 \\ a_7 & a_8 & a_9 \\ a_{10} & a_{11} & a_{12} \end{pmatrix} \middle| a_i \in N(Z_2); 1 \leq i \leq 12 \right\}$$

be two strong neutrosophic bivector spaces over the neutrosophic bifield $F = F_1 \cup F_2 = Z_{29}I \cup N(Z_2)$. Define a $T = T_1 \cup T_2 : V = V_1 \cup V_2 \to W = W_1 \cup W_2$ such that $T^t$, its bitranspose of T.

a. Prove birange of $T^t$ is the biannihilator of the binull space of T.
b. Prove birank $T^t$ = birank T.

39. For the V and W given in problem (38) Find $L^2$ (V, W) and $L^2$ (W, V).

40. For V and W given in problem (38) find a T such that (a) T is biinverible (b) T is binon singular.

41. For a give strong neutrosophic bivector space $V = V_1 \cup V_2 =$

$$\left\{ \begin{pmatrix} a & b & c \\ d & e & f \\ g & h & i \end{pmatrix} \middle| a,b,c,d,e,f,g,h,i \in N(Q) \right\} \cup$$

$$\left\{ \begin{pmatrix} a_1 & a_2 & a_3 \\ a_4 & a_5 & a_6 \end{pmatrix} \middle| a_i \in N(Z_{29}); 1 \leq i \leq 6 \right\}$$

defined over the neutrosophic bifield $F = N(Q) \cup N(Z_{29})$. Define a bilinear functional $f = f_1 \cup f_2$ from V into F.

a. Find for a bibasis B of $V = V_1 \cup V_2$, the bibasis $B^*$ of $V^* = V_1^* \cup V_2^*$.
b. Prove $V^{**} = V$ and bidimension V = bidimension $V^*$.

42. Obtain some interesting properties about bilinear functionals defined on a strong neutrosophic bivector space over a neutrosophic bifield.

43. Find for $V = V_1 \cup V_2 =$

$$\left\{ \begin{pmatrix} a & b \\ c & d \end{pmatrix} \middle| a,b,c,d \in N(Q) \right\} \cup$$

$$\left\{ \begin{pmatrix} a_1 & a_2 & a_3 & a_4 \\ a_5 & a_6 & a_7 & a_8 \end{pmatrix} \middle| a_i \in N(Z_{11}); 1 \le i \le 8 \right\}$$

a strong neutrosophic bivector space defined over the neutrosophic bifield $F = N(Q) \cup N(Z_{11})$, a T bidiagonalizable linear bioperator on V. If $B = B_1 \cup B_2$ is a bibasis prove $[T]B = [T_1]_{B_1} \cup [T_2]_{B_2}$.

44. Let $V = V_1 \cup V_2 =$

$$\left\{ \begin{pmatrix} a_1 & a_2 \\ a_3 & a_4 \\ a_5 & a_6 \\ a_7 & a_8 \end{pmatrix} \middle| a_i \in N(Q); 1 \le i \le 8 \right\} \cup$$

$$\left\{ \begin{pmatrix} a_1 & a_2 & a_3 & a_4 & a_5 \\ a_6 & a_7 & a_8 & a_9 & a_{10} \end{pmatrix} \middle| a_i \in N(Z_{17}); 1 \le i \le 8 \right\}$$

be a strong neutrosophic bivector space over the neutrosophic bifield $F = F_1 \cup F_2 = N(Q) \cup N(Z_{17})$. Find $V^*$. Prove $V^{**} = V$. Find two distinct strong neutrosophic bivector subspaces $W_1$ and $W_2$ in V and find $W_1^o$ and $W_2^o$.

45. Let $V = V_1 \cup V_2 =$

$$\left\{ \begin{pmatrix} a & b \\ c & d \end{pmatrix} \middle| a,b,c,d \in N(Q) \right\} \cup \{Z_{11}I\},$$

$$W = \left\{ \begin{pmatrix} a & b & e \\ c & d & f \end{pmatrix} \middle| a,b,c,d,e,f \in N(Q) \right\} \cup$$

$$\{(a, a, a, b, b) \mid a, b \in Z_{11}I\} = W_1 \cup W_2$$

and $P = P_1 \cup P_2 =$

$$\left\{ \begin{pmatrix} a & b \\ c & d \\ e & f \end{pmatrix} \middle| a,b,c,d,e,f \in N(Q) \right\} \cup$$

$$\left\{ \begin{pmatrix} a & b & c & d & e \\ f & g & h & i & j \end{pmatrix} \middle| a,b,c,d,e,f,g,h,i,j \in Z_{11}I \right\}$$

be three strong neutrosophic bivector spaces over the bifield F $= QI \cup Z_{11}I$.

a. Find linear bioperator $T = T_1 \cup T_2 : V_1 \cup V_2 = V \to W = W_1 \cup W_2$ and $S = S_1 \cup S_2 : W = W_1 \cup W_2 \to P = P_1 \cup P_2$.
b. Find $L^2(V, W)$, $L^2(W, P)$ and $L^2(V, P)$ and their bidimensions.

46. Let $V = V_1 \cup V_2 =$

$$\left\{ \begin{pmatrix} a & b \\ c & d \end{pmatrix} \middle| a,b,c,d \in Z_7I \right\} \cup$$

$$\left\{ \begin{pmatrix} a & b \\ c & d \end{pmatrix} \middle| a,b,c,d \in N(Q) \right\}$$

be a strong neutrosophic bilinear algebra over the neutrosophic bifield $F = F_1 \cup F_2 = Z_7I \cup N(Q)$.

a. What is the bidimension of V?
   b. Find a bibasis of V.
   c. What is the bidimension of $L^2(V, W) = L^2(V_1, W_1) \cup L^2(V_2, W_2)$?
   d. Find a linear bifunctional $f = f_1 \cup f_2 : V = V_1 \cup V_2 \to F_1 \cup F_2$ and find bikernel $f = \text{kernel } f_1 \cup \text{kernel } f_2$.

47. Let $F = F_1 \cup F_2$ be a neutrosophic bifield. Show that the neutrosophic biideal generated by finite number of neutrosophic bipolynomial $f^1$, $f^2$ where $f^1 = f_1^1 \cup f_2^1$ and $f^2 = f_1^2 \cup f_2^2$ in $F[x] = F_1[x] \cup F_2[x]$ is the intersection of all neutrosophic biideals in $F[x]$.

48. Let $(n_1, n_2)$ be a biset of positive integers and $F = F_1 \cup F_2$ be a neutrosophic bifield, let $W = W_1 \cup W_2$ be the set of all bivectors $\left(x_1^1,...,x_{n_1}^1\right) \cup \left(x_1^2,...,x_{n_3}^2\right)$ in $F^{n_1} \cup F^{n_2}$ such that $x_1^1 + x_2^1 +...+ x_{n_1}^1 = 0$ and $x_1^2 + x_2^2 +...+ x_{n_2}^2 = 0$.
   a. Prove $W^o = W_1^o \cup W_2^o$ consists of all bilinear functionals $f = f^1 \cup f^2$ of the form
   $$f_1\left(x_1^1, x_2^1,..., x_{n_1}^1\right) \cup f_2\left(x_1^2, x_2^2,..., x_{n_2}^2\right) = c_1 \sum_{j=1}^{n_1} x_j^1 \cup c_2 \sum_{j=1}^{n_2} x_j^2.$$
   b. Show that the bidual space $W^*$ of $W$ can be naturally identified with the bilinear functionals
   $$f_1\left(x_1^1, x_2^1,..., x_{n_1}^1\right) \cup f_2\left(x_1^2, x_2^2,..., x_{n_2}^2\right)$$
   $$= \left(c_1^1 x_1^1 +...+ c_{n_1}^1 x_{n_1}^1\right) \cup \left(c_1^2 x_1^2 +...+ c_{n_2}^2 x_{n_2}^2\right)$$
   on $F_1^{n_1} \cup F_2^{n_2}$ which satisfy $c_1^i +...+ c_{n_1}^i = 0$ for $i = 1, 2$.

49. Let $W = W_1 \cup W_2$ be a strong neutrosophic bisubspace of a finite $(n_1, n_2)$ bidimensional bivector space over $V = V_1 \cup V_2$ and if $g = g^1 \cup g^2 = \left\{\left(g_1^1,...,g_{r_1}^1\right)\right\} \cup \left\{\left(g_1^2,...,g_{r_2}^2\right)\right\}$ is a bibasis for $W^o = W_1^o \cup W_2^o$ then prove

$$W = \bigcap_i N_{g_i} = \bigcap_{i_1=1}^{r_1} N^1_{g_{i_1}} \cup \bigcap_{i_2=1}^{r_2} N^2_{g_{i_2}}$$

where $\{(N^1_1,...,N^1_{r_1})\} \cup \{(N^2_1,...,N^2_{r_2})\}$ is the biset of binull space of bilinear functionals

$$f = f^1 \cup f^2 = \{(f^1_1, f^1_2,...,f^1_{r_1})\} \cup \{(f^2_1, f^2_2,...,f^2_{r_2})\} \text{ and}$$

$$g = g^1 \cup g^2 = \{(g^1_1,...,g^1_{r_1})\} \cup \{(g^2_1,...,g^2_{r_2})\}$$

and is the bilinear combination of the bilinear functionals $f = f^1 \cup f^2$.

50. Let $V = V_1 \cup V_2 \cup V_3 =$

$$\left\{ \begin{pmatrix} a & b & e \\ c & d & f \end{pmatrix} \middle| a,b,c,d,e,f \in Z_7I \right\} \cup$$

$$\left\{ \begin{pmatrix} a_1 & a_2 \\ a_3 & a_4 \\ a_5 & a_6 \\ a_7 & a_8 \end{pmatrix} \middle| a_i \in Z_{11}I; 1 \le i \le 8 \right\} \cup$$

$$\left\{ \begin{pmatrix} a_1 & a_2 & a_3 & a_4 & a_5 & a_6 \\ a_7 & a_8 & a_9 & a_{10} & a_{11} & a_{12} \end{pmatrix} \middle| a_i \in Z_{13}I; 1 \le i \le 12 \right\}$$

be a neutrosophic trivector space over the 3-field $F = Z_7 \cup Z_{11} \cup Z_{13}$.

a. Find a tribasis of V.
b. Find neutrosophic trivector subspaces of V.
c. What is the 3-dimension of V?
d. Define a neutrosophic trilinear operator T on V which is non invertible (if $T = T_1 \cup T_2 \cup T_3$ then $T^{-1} = T_1^{-1} \cup T_2^{-1} \cup T_3^{-1}$) Show $T^{-1}$ does not exists for the T defined.

51. Let $V = V_1 \cup V_2 \cup V_3 \cup V_4 \cup V_5 =$

$$\left\{ \begin{pmatrix} a & b & c \\ b & a & c \end{pmatrix} \middle| a, b, c \in Z_7I \right\} \cup$$

$$\left\{ \begin{pmatrix} a_1 & a_2 \\ a_3 & a_4 \\ a_5 & a_6 \\ a_7 & a_8 \end{pmatrix} \middle| a_i \in Z_{11}I; 1 \leq i \leq 8 \right\} \cup$$

$$\left\{ \begin{pmatrix} a_1 \\ a_2 \\ a_3 \\ a_4 \\ a_5 \\ a_6 \end{pmatrix} \middle| a_i \in Z_{17}I; 1 \leq i \leq 6 \right\} \cup$$

$$\left\{ \begin{pmatrix} a_1 & a_2 & a_3 & a_4 & a_5 \\ a_6 & a_7 & a_8 & a_9 & a_{10} \\ a_{11} & a_{12} & a_{13} & a_{14} & a_{15} \end{pmatrix} \middle| a_i \in N(Q) \right\} \cup$$

$\left\{ \sum_{i=0}^{25} a_i x^i \right.$; all polynomials in the variable x with coefficients from the neutrosophic field $N(Z_{19})$; $a_i \in N(Z_{19})$; $0 \leq i \leq 25 \}$ be a strong neutrosophic 5-vector space over the neutrosophic 5-field $F = Z_7I \cup Z_{11}I \cup Z_{17}I \cup QI \cup N(Z_{19})$.

a. Find a strong neutrosophic 5-vector subspace of V.
b. Is V pseudo simple?
c. Can on V be defined a strong neutrosophic 5-linear operator T so that T is invertible?
d. What is the 5-dimension of V?
e. Find a 5-basis of V.

f. Find a 5-linearly independent 5-subset of V which is not a 5-basis.

52. Obtain some important properties about $SNHom_F$ (V, V); V is a strong neutrosophic n-vector space over a neutrosophic n-field $F = F_1 \cup F_2 \cup \ldots \cup F_n$. What is the algebraic structure of $SNHom_F$ (V, V)?

53. Characterize those strong neutrosophic n-vector spaces which are simple.

54. Give an example of a strong neutrosophic 5-vector space which is pseudo simple.

55. Prove in case of a finite n-vector space V of type II, where dim $V = (n_1, n_2, \ldots, n_n)$. Rank T + nullity T = dim V.

56. Derive primary decomposition theorem for strong neutrosophic n-vector space over the neutrosophic n-field.

57. Let $V = V_1 \cup V_2 \cup V_3 \cup V_4 \cup V_5 =$

$$\left\{ \begin{pmatrix} a & b & 0 & 0 \\ c & d & 0 & 0 \\ 0 & 0 & a & d \\ 0 & 0 & b & c \end{pmatrix} \middle| a,b,c,d \in N(Q) \right\} \cup$$

$$\left\{ \begin{pmatrix} a_1 & a_2 & a_3 \\ a_4 & a_5 & a_6 \\ a_7 & a_8 & a_9 \end{pmatrix} \middle| a_i \in N(Z_{29}); 1 \le i \le 9 \right\} \cup$$

$$\left\{ \begin{pmatrix} a & b \\ c & d \end{pmatrix} \middle| a,b,c,d \in N(Z_{11}) \right\} \cup$$

$\{(a_1, a_2, \ldots, a_9) \mid a_i \in N(Z_{23}); 1 \leq i \leq 29)\} \cup \{$All $5 \times 5$ upper triangular matrices with entries from $N(Z_{13})\}$ be a strong neutrosophic 5-linear algebra over the 5-field, $F = QI \cup Z_{29}I \cup Z_{11}I \cup N(Z_{23}) \cup Z_{13}I$.

Let $W = W_1 \cup W_2 \cup W_3 \cup W_4 \cup W_5 =$

$$\left\{ \begin{pmatrix} a & a & 0 & 0 \\ b & b & 0 & 0 \\ 0 & 0 & a & a \\ 0 & 0 & b & b \end{pmatrix} \middle| a, b \in QI \right\} \cup$$

$$\left\{ \begin{pmatrix} a & a & a \\ a & a & a \\ a & a & a \end{pmatrix} \middle| a \in N(Z_{29}) \right\} \cup \left\{ \begin{pmatrix} a & a \\ b & b \end{pmatrix} \middle| a, b \in N(Z_{11}) \right\} \cup$$

$$\{(a, a, a, a, a, a, a, a, a) \mid a \in N(Z_{23})\} \cup$$

$$\left\{ \begin{pmatrix} a & a & a & a & a \\ 0 & a & a & a & a \\ 0 & 0 & a & a & a \\ 0 & 0 & 0 & a & a \\ 0 & 0 & 0 & 0 & a \end{pmatrix} \middle| a \in N(Z_{13}) \right\}$$

$\subseteq V_1 \cup V_2 \cup V_3 \cup V_4 \cup V_5$ be a strong neutrosophic 5-linear subalgebra over V the 5-field F.

For $\beta = \beta_1 \cup \beta_2 \cup \beta_3 \cup \beta_4 \cup \beta_5 =$

$$\begin{pmatrix} 7 & 2I & 0 & 0 \\ 5+I & 0 & 0 & 0 \\ 0 & 0 & 7 & 0 \\ 0 & 0 & 2I & 5+I \end{pmatrix} \cup \left\{ \begin{pmatrix} 0 & 3I & 0 \\ 7+I & 0 & I \\ 1 & 0 & 2+I \end{pmatrix} \right\} \cup$$

$$\left\{ \begin{pmatrix} 3 & 3+8I \\ 10I & 0 \end{pmatrix} \right\} \cup$$

$$\{(0, I, 0, 3I, 7+I, 0, 1, 2I+1, 0)\} \cup$$

$$\left\{ \begin{pmatrix} 9 & 1 & 2I & 9 & I \\ 0 & I & 0 & 1 & 2I \\ 0 & 0 & 7I & 0 & 0 \\ 0 & 0 & 0 & 2I & -4 \\ 0 & 0 & 0 & 0 & I \end{pmatrix} \right\}$$

$\in V_1 \cup V_2 \cup V_3 \cup V_4 \cup V_5$, find $\alpha \in W$ such that $= \alpha_1 \cup \alpha_2 \cup \alpha_3 \cup \alpha_4 \cup \alpha_5 \in W_1 \cup W_2 \cup W_3 \cup W_4 \cup W_5 = W$ is the best 5-approximation of $\beta$. Prove $\beta - \alpha$ is 5-orthogonal to every 5-vector in W, that is $\beta_1 - \alpha_1 \cup \beta_2 - \alpha_2 \cup \beta_3 - \alpha_3 \cup \beta_4 - \alpha_4 \cup \beta_5 - \alpha_5$ is 5-orthogonal to every 5-vector in $W = W_1 \cup W_2 \cup W_3 \cup W_4 \cup W_5$. Prove $\beta_i - \alpha_i$ is orthogonal to every vector in $W_i$; $i = 1, 2, 3, 4, 5$. Find $W^\perp$. Prove $V = W \oplus W^\perp$ where $W^\perp = W_1^\perp \cup W_2^\perp \cup W_3^\perp \cup W_4^\perp \cup W_5^\perp$, that is $W_1 \oplus W_1^\perp \cup W_2 \oplus W_2^\perp \cup W_3 \oplus W_3^\perp \cup W_4 \oplus W_4^\perp \cup W_5 \oplus W_5^\perp = V_1 \cup V_2 \cup V_3 \cup V_4 \cup V_5 = V$.

58. Prove if $V = V_1 \cup V_2 \cup \ldots \cup V_n$ is a strong neutrosophic n-vector space over the n-field $F = F_1 \cup F_2 \cup \ldots \cup F_n$ of finite $(n_1, n_2, \ldots, n_n)$ dimension over F. $T = T_1 \cup T_2 \cup \ldots \cup T_n$ be a n-linear operator on V.
Prove there exists a n-set $\{\alpha_1^1, \ldots, \alpha_{r_1}^1\} \cup \{\alpha_1^2, \ldots, \alpha_{r_2}^2\} \cup \ldots \cup \{\alpha_1^n, \ldots, \alpha_{r_n}^n\}$ in V such that $V = V_1 \cup V_2 \cup \ldots \cup V_n = Z(\alpha_1^1; T_1) \oplus \ldots \oplus Z(\alpha_{r_1}^1; T_1) \cup Z(\alpha_1^2; T_2) \oplus \ldots \oplus Z(\alpha_{r_2}^2; T_2) \cup \ldots \cup Z(\alpha_1^n; T_n) \oplus \ldots \oplus Z(\alpha_{r_n}^n; T_n)$; i.e., V is the n-direct sum of n-cyclic strong neutrosophic n-vector subspaces.

59. State and prove Generalized Cayley Hamilton Theorem for a finite $(n_1, \ldots, n_n)$ dimensional strong neutrosophic n-vector space $V = V_1 \cup V_2 \cup \ldots \cup V_n$ over the neutrosophic n-field $F = F_1 \cup F_2 \cup \ldots \cup F_n$ after appropriate changes.

60. Define n-projections associated with the n-primary decomposition of $T = T_1 \cup T_2 \cup \ldots \cup T_n$.

61. Let $T = T_1 \cup T_2 \cup \ldots \cup T_n$ be a n-linear operator on the strong neutrosophic n-vector space $V = V_1 \cup V_2 \cup \ldots \cup V_n$ over the neutrosophic n-field $F = F_1 \cup F_2 \cup \ldots \cup F_n$ ($F_i$'s are not pure neutrosophic; $i = 1, 2, \ldots, n$)
    Let
    $$\{W_1^1, \ldots, W_{k_1}^1\} \cup \{W_1^2, \ldots, W_{k_2}^2\} \cup \ldots \cup \{W_1^n, \ldots, W_{k_n}^n\}$$
    and
    $$\{E_1^1, \ldots, E_{k_1}^1\} \cup \{E_1^2, \ldots, E_{k_2}^2\} \cup \ldots \cup \{E_1^n, \ldots, E_{k_n}^n\};$$
    where $\{W_1^i, \ldots, W_{k_i}^i\}$ are independent for $i = 1, 2, \ldots, n$. $E = E_1 \cup E_2 \cup \ldots \cup E_n$ is a n-projection operator on V such that $E^2 = E$ that is $E^2 = (E_1 \cup \ldots \cup E_n)^2 = E_1^2 \cup \ldots \cup E_n^2 = E_1 \cup \ldots \cup E_n$ (That is each $E_i$ is a projection of $V_i$ such that $E_i^2 = E_i$, $i = 1, 2, \ldots, n$).
    Then a necessary and sufficient condition that each strong neutrosophic n-vector subspace $W_i^t$ to be invariant under $T_i$ for $1 \le i \le k_t$; $t = 1, 2, \ldots, n$ is that $E_i^t T_t = T_t E_i^t$ or $ET = TE$ for every $1 \le i \le k_t$; $t = 1, 2, \ldots, n$.

62. Let $T = T_1 \cup T_2 \cup \ldots \cup T_n$ be a n-linear operator on a $(n_1, n_2, \ldots, n_n)$ finite n-dimensional strong neutrosophic n-vector space $V = V_1 \cup V_2 \cup \ldots \cup V_n$ over the neutrosophic n-field $F = F_1 \cup F_2 \cup \ldots \cup F_n$ ($F_i$'s are not pure neutrosophic; $i = 1, 2, \ldots, n$). Suppose that the n-minimal neutrosophic polynomial for $T = T_1 \cup T_2 \cup \ldots \cup T_n$ decomposes over $F = F_1 \cup F_2 \cup \ldots \cup F_n$ into a n-product of n-linear neutrosophic polynomials. Then there is a n-diagonalizable operator $N = N_1 \cup N_2 \cup \ldots \cup N_n$ on $V = V_1 \cup V_2 \cup \ldots \cup V_n$ such that

a.  T = D + N that is
    $T_1 \cup T_2 \cup \ldots \cup T_n$
    $= D_1 \cup D_2 \cup \ldots \cup D_n + (N_1 \cup N_2 \cup \ldots \cup N_n)$.
    $= D_1 + N_1 \cup D_2 + N_2 \cup \ldots \cup D_n + N_n$.
b.  DN = ND that is
    $(D_1 \cup D_2 \cup \ldots \cup D_n)(N_1 \cup N_2 \cup \ldots \cup N_n)$
    $= D_1 N_1 \cup D_2 N_2 \cup \ldots \cup D_n N_n$
    $= N_1 D_1 \cup N_2 D_2 \cup \ldots \cup N_n D_n$
    $= ND$.

The n-diagonalizable operator $D = D_1 \cup D_2 \cup \ldots \cup D_n$ and the n-nilpotent operator $N = N_1 \cup N_2 \cup \ldots \cup N_n$ are uniquely determined by (a) and (b) and each of them is a n-polynomial in $T_1, T_2, \ldots, T_n$. Prove.

63. Prove $S(\beta;W) = S(\beta_1;W_1) \cup S(\beta_2; W_2) \cup \ldots \cup S(\beta_n; W_n)$ is the n-conductor of T where T is a n-linear operator on the strong neutrosophic n-vector space $V = V_1 \cup V_2 \cup \ldots \cup V_n$ and $W = W_1 \cup W_2 \cup \ldots \cup W_n \subseteq V_1 \cup V_2 \cup \ldots \cup V_n$ is a proper T-n-invariant neutrosophic n-vector subspace of V.
    Prove some interesting results about these structures like relating it with neutrosophic n-ideals. Hence or other wise prove the n-cyclic decomposition theorem.

64. If $T = T_1 \cup T_2 \cup \ldots \cup T_n$ is a n-linear operator of a finite $(n_1, n_2, \ldots, n_n)$ dimension strong neutrosophic n-vector space $V = V_1 \cup V_2 \cup \ldots \cup V_n$ over the n-field $F = F_1 \cup F_2 \cup \ldots \cup F_n$. Prove T is n-diagonalizable if and only if the n-characteristic n-polynomial $T = T_1 \cup T_2 \cup \ldots \cup T_n$ is $f = f_1 \cup f_2 \cup \ldots \cup f_n$
    $= (x - c_1^1)^{d_1^1} \ldots (x - c_{k_1}^1)^{d_{k_1}^1} \cup \ldots \cup (x - c_1^n)^{d_1^n} \ldots (x - c_{k_1}^n)^{d_{k_1}^n}$
    under the usual notations.

65. Obtain some interesting properties about quasi strong neutrosophic n-vector spaces over quasi neutrosophic n-field $F = F_1 \cup F_2 \cup \ldots \cup F_n$.

66. Let $V = V_1 \cup V_2 \cup \ldots \cup V_6 =$

$$\left\{ \begin{pmatrix} a_1 & a_2 & a_3 \\ a_4 & a_5 & a_6 \end{pmatrix} \middle| a_i \in Z_7 I \right\} \cup$$

$$\left\{ \begin{pmatrix} a & a \\ b & b \\ c & c \\ d & d \\ e & e \end{pmatrix} \middle| a,b,c,d,e \in Z_{11}I \right\} \cup$$

{All 5 × 5 lower triangular matrices with entries from the field N($Z_2$)} $\cup$

$$\left\{ \begin{pmatrix} 0 & 0 & 0 & a \\ 0 & 0 & b & 0 \\ 0 & c & 0 & 0 \\ d & 0 & 0 & 0 \end{pmatrix} \middle| a,b,c,d \in Z_3 I \right\} \cup$$

$\left\{ \sum_{i=0}^{4} a_i x^i \right.$ ; all polynomials in the variable x of degree less than or equal to 4 with coefficients from $Z_5 I\} \cup$

$$\left\{ \begin{pmatrix} a_1 & a_2 & a_3 & a_4 & a_5 \\ a_6 & a_7 & a_8 & a_9 & a_{10} \end{pmatrix} \middle| a_i \in Z_{23}I \right\}$$

be a strong neutrosophic 6-vector space over the neutrosophic 6-field $F = Z_7 I \cup Z_{11}I \cup N(Z_2) \cup Z_3 I \cup Z_5 I \cup Z_{23}I$.

a. Find a 6-basis of V.
b. Is V n-finite?
c. Find at least two strong neutrosophic 6-subspaces of V.
d. Write V as a direct sum of neutrosophic strong 6-vector subspaces.

e. If $F = Z_7I \cup Z_{11}I \cup N(Z_2) \cup Z_3I \cup Z_5I \cup Z_{23}I$ is changed to a 6-field $K = Z_7 \cup Z_{11} \cup Z_2 \cup Z_3 \cup Z_5 \cup Z_{23}$. Find a 6-basis.
f. Does the change of 6-field affect the structure of V? Justify your claim.

67. Find some interesting properties about neutrosophic n-linear algebras.

68. Can Cayley Hamilton theorem hold good for neutrosophic n-vector spaces defined over a real n-field? Justify your answer!

69. For the strong neutrosophic 4-linear algebra given by $V = V_1 \cup V_2 \cup V_3 \cup V_4 =$

$$\left\{ \begin{pmatrix} a & b \\ c & d \end{pmatrix} \middle| a,b,c,d \in Z_2I \right\} \cup$$

$$\{(a_1, a_2, a_3) \mid a_i \in Z_3I; 1 \leq i \leq 3\} \cup$$

$$\left\{ \begin{pmatrix} a_1 & 0 & 0 \\ 0 & a_2 & 0 \\ 0 & 0 & a_3 \end{pmatrix} \middle| a_i \in Z_5I; 1 \leq i \leq 3 \right\} \cup$$

{All $4 \times 4$ upper triangular neutrosophic matrices with entries from $Z_7I$} defined over the neutrosophic 4-field $F = Z_2I \cup Z_3I \cup Z_5I \cup Z_7I$.

a. Find a 4-basis for V.
b. Define a strong neutrosophic linear operator T on V and for that T find the neutrosophic 4-characteristic polynomial, neutrosophic 4-eigen values and 4-eigen vectors.
c. If F is replaced by $K = Z_2 \cup Z_3 \cup Z_5 \cup Z_7$ will the 4-basis be different?
d. Find $SNHom_K(V, V)$ and $SNHom_F(V, V)$. What is the difference between them as algebraic structures?

e. Is 4-rank T + 4-nullity T = 4-dim V? Justify your claim (T:V →V is a neutrosophic strong linear operator on V).
f. If V is assumed only as a neutrosophic strong 4-vector space over the neutrosophic 4-field, what will be 4-basis of V? Will the 4-basis of V differ? Justify / substantiate your claim.

70. Let $V = V_1 \cup V_2 \cup V_3 \cup V_4 =$

$$\begin{pmatrix} I & 0 \\ 1 & 3I \end{pmatrix} \cup \begin{pmatrix} I & 2I & 1 & 0 \\ 0 & 3I & I & 6I \\ 6 & 0 & 3I & 0 \\ 2I & 1 & 0 & 1 \end{pmatrix} \cup$$

$$\begin{pmatrix} I & 0 & 1 \\ 0 & 1 & I \\ 1 & 0 & I \end{pmatrix} \cup \begin{pmatrix} I & 0 & 0 & 0 & I \\ 2 & 0 & I & 1 & 0 \\ 0 & 1 & 1 & 0 & I \\ 0 & 0 & I & 2 & 0 \\ 0 & I & 0 & 0 & 2I \end{pmatrix}$$

be a neutrosophic 4-matrix with entries from the 4-field $F = F_1 \cup F_2 \cup F_3 \cup F_4 = Z_5I \cup Z_7I \cup Z_2I \cup Z_3I$ respectively. Find the 4-characteristic neutrosophic 4-polynomial associated with the neutrosophic 4-matrix V. Can this have neutrosophic 4-eigen values? Justify your claim.

71. For the example 2.3.72 given chapter two find $SNHom_F$ (V, W). Find a $T : V \to W$ so that $kerT = (0) \cup (0)$.

72. Obtain some interesting and special features enjoyed by quasi neutrosophic n-vector spaces.

73. If $L = L^1_{\alpha_1} \cup L^2_{\alpha_2} \cup ... \cup L^n_{\alpha_n}$ is a n-linear function induced by $\alpha = \alpha^1 \cup \alpha^2 \cup ... \cup \alpha^n$ in $V = V_1 \cup V_2 \cup ... \cup V_n$, a strong neutrosophic n-vector space over the n-field $F = F_1 \cup F_2 \cup ...$

$\cup\, F_n$. Is $\alpha = \alpha^1 \cup \alpha^2 \cup \ldots \cup \alpha^n \mapsto L_\alpha = L^1_{\alpha_1} \cup L^2_{\alpha_2} \cup \ldots \cup L^n_{\alpha_n}$ a n-isomorphism of $V = V_1 \cup V_2 \cup \ldots \cup V_n$ onto $V^{**} = V^{**}_1 \cup V^{**}_2 \cup \ldots \cup V^{**}_n$ ? Justify your claim.

74. Study the properties enjoyed by SNL $(V_1, F_1) \cup$ SNL $(V_2, F_2) \cup \ldots \cup$ SNL $(V_n, F_n)$ where $V = V_1 \cup V_2 \cup \ldots \cup V_n$ is a strong neutrosophic n-vector space defined over the neutrosophic n-field $F = F_1 \cup F_2 \cup \ldots \cup F_n$.

75. Find a set $X = X_1 \cup X_2 \cup X_3 \cup X_4 \subseteq V = V_1 \cup V_2 \cup V_3 \cup V_4$

$$= \left\{ \begin{pmatrix} a_1 & a_2 & a_3 & a_4 \\ a_5 & a_6 & a_7 & a_8 \end{pmatrix} \middle| a_i \in Z_7 I; 1 \le i \le 8 \right\} \cup$$

{All $7 \times 3$ neutrosophic matrices with entries from the neutrosophic field $Z_2 I$} $\cup$

$\left\{ \sum_{i=0}^{5} a_i x^i \right.$ ; all neutrosophic polynomials of degree less than or equal to 5 with coefficients from $Z_5 I$ in the variable x$\} \cup$

$$\left\{ \begin{pmatrix} 0 & 0 & 0 & a \\ 0 & 0 & b & 0 \\ 0 & c & 0 & 0 \\ d & 0 & 0 & 0 \end{pmatrix} \middle| a,b,c,d \in N(Q) \right\},$$

a strong neutrosophic 4-vector space over the 4-field $F = F_1 \cup F_2 \cup F_3 \cup F_4 = Z_7 I \cup Z_2 I \cup Z_5 I \cup QI$; a 4-linearly independent 4-set of V which is not a 4-basis of V.
 a. Find a 4-basis of V.
 b. What is the 4-dimension of V?
 c. Define a invertible 4-linear operator on V.
 d. Find $\text{SNHom}_F(V, V)$. What is the dimension of $\text{SNHom}_F(V, V)$?

e. Does SNL(V, F) = SNL($V_1$, $F_1$) $\cup$ SNL($V_2$, $F_2$) $\cup$ SNL($V_3$, $F_3$) $\cup$ SNL($V_4$, $F_4$) exist? Justify your claim.

76. Let $V = V_1 \cup V_2 \cup \ldots \cup V_n$ be neutrosophic strong n-linear algebra over the n-field $F = F_1 \cup F_2 \cup \ldots \cup F_n$. Consider a n-basis $\{\alpha_1^1,\ldots,\alpha_{n_1}^1\} \cup \{\alpha_1^2,\ldots,\alpha_{n_2}^2\} \cup \ldots \cup \{\alpha_1^n,\ldots,\alpha_{n_n}^n\}$ of V over F. If $W = W_1 \cup W_2 \cup \ldots \cup W_n$ is a strong neutrosophic n-vector space over the same F and if
$$\beta = \{\beta_1^1,\ldots,\beta_{n_1}^1\} \cup \{\beta_1^2,\ldots,\beta_{n_2}^2\} \cup \ldots \cup \{\beta_1^n,\ldots,\beta_{n_n}^n\}$$
be any n-vector in W. Prove there exists precisely a n-linear transformation $T = T_1 \cup T_2 \cup \ldots \cup T_n$ from V into W such that $T_i(\alpha_j^i) = \beta_j^i$ for $j = 1, 2, \ldots, n_i$ and $i = 1, 2$.

77. Prove if $V = V_1 \cup V_2 \cup \ldots \cup V_n$ and $W = W_1 \cup W_2 \cup \ldots \cup W_n$ are two strong neutrosophic n-vector spaces over the same neutrosophic n-field $F = F_1 \cup F_2 \cup \ldots \cup F_n$ of type II. If $T = T_1 \cup T_2 \cup \ldots \cup T_n$ is a n-linear transformation of V into W then prove the following are equivalent.
    a. $T = T_1 \cup T_2 \cup \ldots \cup T_n$ is n-invertible.
    b. $T = T_1 \cup T_2 \cup \ldots \cup T_n$ is n-non singular.
    c. $T = T_1 \cup T_2 \cup \ldots \cup T_n$ is onto that is the n-range of $T = T_1 \cup T_2 \cup \ldots \cup T_n$ is $W = W_1 \cup W_2 \cup \ldots \cup W_n$.

78. Prove every ($n_1, n_2, \ldots, n_n$) dimensional strong neutrosophic n-vector space $V = V_1 \cup V_2 \cup \ldots \cup V_n$ over the neutrosophic n-field $F = F_1 \cup F_2 \cup \ldots \cup F_n$ is n-isomorphism to $F_1^{n_1} \cup F_2^{n_2} \cup \ldots \cup F_n^{n_n}$.

79. Let $V = V_1 \cup V_2 \cup \ldots \cup V_n$ be a finite ($n_1, n_2, \ldots, n_n$) n-dimensional strong neutrosophic n-vector space over the neutrosophic n-field $F = F_1 \cup F_2 \cup \ldots \cup F_n$.
Let
$$B = B_1 \cup B_2 \cup \ldots \cup B_n$$
$$= \{\alpha_1^1,\ldots,\alpha_{n_1}^1\} \cup \{\alpha_1^2,\ldots,\alpha_{n_2}^2\} \cup \ldots \cup \{\alpha_1^n,\ldots,\alpha_{n_n}^n\}$$

be a n-basis of $V = V_1 \cup V_2 \cup \ldots \cup V_n$. There is a unique n-dual basis (dual n-basis)
$$B = B_1^* \cup B_2^* \cup \ldots \cup B_n^*$$
$$= \{f_1^1, f_2^1, \ldots, f_{n_1}^1\} \cup \{f_1^2, f_2^2, \ldots, f_{n_2}^2\} \cup \ldots \cup \{f_1^n, f_2^n, \ldots, f_{n_n}^n\}$$
for $V^* = V_1^* \cup V_2^* \cup \ldots \cup V_n^*$ such that $f_i^k(\alpha_j) = \delta_{ij}^k$.

Prove for each n-linear functional $f = f_1 \cup f_2 \cup \ldots \cup f_n$ we have
$$f = \sum_{k=1}^{n_i} f^i(\alpha_k^i) f_k^i$$
that is
$$f = \left( \sum_{k=1}^{n_1} f^1(\alpha_k^1) f_k^1 \right) \cup \left( \sum_{k=1}^{n_2} f^2(\alpha_k^2) f_k^2 \right) \cup \ldots \cup \left( \sum_{k=1}^{n_n} f^n(\alpha_k^n) f_k^n \right)$$

and for each n – vector $\alpha = \alpha_1 \cup \alpha_2 \cup \ldots \cup \alpha_n$ in $V = V_1 \cup V_2 \cup \ldots \cup V_n$ we have
$$\alpha = \left( \sum_{k=1}^{n_1} f_k^1(\alpha^1)(\alpha_k^1) \right) \cup \left( \sum_{k=1}^{n_2} f_k^2(\alpha^2)(\alpha_k^2) \right)$$
$$\cup \ldots \cup \left( \sum_{k=1}^{n_n} f_k^n(\alpha^n)(\alpha_k^n) \right).$$

80. Obtain some important properties about n-best approximations on strong neutrosophic n-vector space over the n-field $F = F_1 \cup F_2 \cup \ldots \cup F_n$.

81. Let $V = V_1 \cup V_2 \cup V_3 \cup V_4 = \{(a_1, a_2, \ldots, a_9) \mid a_i \in Z_7I, 1 \leq i \leq 9\} \cup \{(a_1, a_2, \ldots, a_{20}) \mid a_i \in N(Z_2), 1 \leq i \leq 20\} \cup \{(a_1, a_2, a_3, a_4, a_5) \mid a_i \in N(Z_{11}), 1 \leq i \leq 5\} \cup \{(a_1, a_2, \ldots, a_8) \mid a_i \in N(Z_5), 1 \leq i \leq 8\}$ be a strong neutrosophic inner product 4-space over the 4-field $Z_7I \cup Z_2I \cup Z_{11}I \cup Z_5I$.

Let $S = S_1 \cup S_2 \cup S_3 \cup S_4 = \{(a_1, a_2, 0, 0, 0, 0, a_7, a_8, 0), (0, 0, a_3, a_4, a_5, 0, 0, 0, a_9) (0, 0, 0, 0, 0, a_6, a_7, a_8, 0) \mid a_i \in Z_7I; 1 \leq i \leq 9\} \cup \{(a_1, a_2, a_3, a_4, a_5, 0, 0, \ldots, 0), (0, 0, 0, 0, 0, 0, 0, a_8,$

$a_9$, $a_{10}$, $a_{11}$, $a_{12}$, 0, 0, ..., 0) | $a_i \in Z_2I$; $1 \le i \le 5$; $i = 8, 9, 10, 11, 12$) $\cup$ {$(a_1, a_2, 0, 0, a_5)$ $(0, 0, 0, a_4, a_5)$} $\cup$ {$(a_1, a_2, 0, 0, a_5, a_6, 0, 0) \cup (a_1, a_2, 0, 0, 0, 0, a_7, a_8)$ | $a_i \in N(Z_5)$} $\subseteq V_1 \cup V_2 \cup V_3 \cup V_4$ be any 4-set of 4-vectors in V. Find the 4-orthogonal complement of S denoted by $S^\perp = S_1^\perp \cup S_2^\perp \cup S_3^\perp \cup S_4^\perp$.

82. Derive Cayley Hamilton theorem and Primary n-decomposition theorem for strong neutrosophic n-vector space V defined over the n-field F.

# FURTHER READING

1. Abraham, R., *Linear and Multilinear Algebra*, W. A. Benjamin Inc., 1966.

2. Albert, A., *Structure of Algebras*, Colloq. Pub., 24, Amer. Math. Soc., 1939.

3. Berlekamp, E.R., *Algebraic Coding Theory,* Mc Graw Hill Inc, 1968.

4. Birkhoff, G., and MacLane, S., *A Survey of Modern Algebra*, Macmillan Publ. Company, 1977.

5. Birkhoff, G., *On the structure of abstract algebras,* Proc. Cambridge Philos. Soc., 31  433-435, 1995.

6. Bruce, Schneier., *Applied Cryptography,* Second Edition, John Wiley, 1996.

7. Burrow, M., *Representation Theory of Finite Groups*, Dover Publications, 1993.

8. Charles W. Curtis, *Linear Algebra – An introductory Approach*, Springer, 1984.

9. Dubreil, P., and Dubreil-Jacotin, M.L., *Lectures on Modern Algebra*, Oliver and Boyd., Edinburgh, 1967.

10. Gel'fand, I.M., *Lectures on linear algebra*, Interscience, New York, 1961.

11. Greub, W.H., *Linear Algebra,* Fourth Edition, Springer-Verlag, 1974.

12. Halmos, P.R., *Finite dimensional vector spaces*, D Van Nostrand Co, Princeton, 1958.

13. Hamming, R.W., *Error Detecting and error correcting codes*, Bell Systems Technical Journal, 29, 147-160, 1950.

14. Harvey E. Rose, Linear Algebra, Bir Khauser Verlag, 2002.

15. Herstein I.N., *Abstract Algebra,* John Wiley,1990.

16. Herstein, I.N., and David J. Winter, *Matrix Theory and Linear Algebra*, Maxwell Pub., 1989.

17. Herstein, I.N., *Topics in Algebra,* John Wiley, 1975.

18. Hoffman, K. and Kunze, R., *Linear algebra*, Prentice Hall of India, 1991.

19. Hummel, J.A., *Introduction to vector functions,* Addison-Wesley, 1967.

20. Jacob Bill, *Linear Functions and Matrix Theory ,* Springer-Verlag, 1995.

21. Jacobson, N., *Lectures in Abstract Algebra*, D Van Nostrand Co, Princeton, 1953.

22. Jacobson, N., *Structure of Rings*, Colloquium Publications, 37, American Mathematical Society, 1956.

23. Johnson, T., *New spectral theorem for vector spaces over finite fields $Z_p$ ,* M.Sc. Dissertation, March 2003 (Guided by Dr. W.B. Vasantha Kandasamy).

24. Katsumi, N., *Fundamentals of Linear Algebra*, McGraw Hill, New York, 1966.

25. Kostrikin, A.I, and Manin, Y. I., *Linear Algebra and Geometry*, Gordon and Breach Science Publishers, 1989.

26. Lang, S., *Algebra*, Addison Wesley, 1967.

27. Lay, D. C., *Linear Algebra and its Applications*, Addison Wesley, 2003.

28. Mac William, F.J., and Sloane N.J.A., *The Theory of Error Correcting Codes*, North Holland Pub., 1977.

29. Pettofrezzo, A. J., *Elements of Linear Algebra,* Prentice-Hall, Englewood Cliffs, NJ, 1970.

30. Pless, V.S., and Huffman, W. C., *Handbook of Coding Theory,* Elsevier Science B.V, 1998.

31. Roman, S., *Advanced Linear Algebra*, Springer-Verlag, New York, 1992.

32. Rorres, C., and Anton H., *Applications of Linear Algebra*, John Wiley & Sons, 1977.

33. Semmes, Stephen, *Some topics pertaining to algebras of linear operators*, November 2002. http://arxiv.org/pdf/math.CA/0211171

34. Shannon, C.E., *A Mathematical Theory of Communication,* Bell Systems Technical Journal, 27, 379-423 and 623-656, 1948.

35. Shilov, G.E., *An Introduction to the Theory of Linear Spaces,* Prentice-Hall, Englewood Cliffs, NJ, 1961.

36. Smarandache, Florentin (editor), *Proceedings of the First International Conference on Neutrosophy, Neutrosophic Logic, Neutrosophic set, Neutrosophic probability and Statistics,* December 1-3, 2001 held at the University of New Mexico, published by Xiquan, Phoenix, 2002.

37. Smarandache, Florentin, *A Unifying field in Logics: Neutrosophic Logic, Neutrosophy, Neutrosophic set, Neutrosophic probability*, second edition, American Research Press, Rehoboth, 1999.

38. Smarandache, Florentin, *An Introduction to Neutrosophy*, http://gallup.unm.edu/~smarandache/Introduction.pdf

39. Smarandache, Florentin, *Collected Papers II*, University of Kishinev Press, Kishinev, 1997.

40. Smarandache, Florentin, *Neutrosophic Logic, A Generalization of the Fuzzy Logic*, http://gallup.unm.edu/~smarandache/NeutLog.txt

41. Smarandache, Florentin, *Neutrosophic Set, A Generalization of the Fuzzy Set*, http://gallup.unm.edu/~smarandache/NeutSet.txt

42. Smarandache, Florentin, *Neutrosophy : A New Branch of Philosophy,* http://gallup.unm.edu/~smarandache/Neutroso.txt

43. Smarandache, Florentin, *Special Algebraic Structures, in Collected Papers III*, Abaddaba, Oradea, 78-81, 2000.

44. Thrall, R.M., and Tornkheim, L., *Vector spaces and matrices*, Wiley, New York, 1957.

45. VaN Lint, J.H., *Introduction to Coding Theory,* Springer, 1999.

46. Vasantha Kandasamy and Rajkumar, R. Use of best approximations in algebraic bicoding theory, Varahmihir Journal of Mathematical Sciences, 6, 509-516, 2006.

47. Vasantha Kandasamy and Thiruvegadam, N., Application of pseudo best approximation to coding theory, *Ultra Sci.,* 17, 139-144, 2005.

48. Vasantha Kandasamy, W.B., *Bialgebraic structures and Smarandache bialgebraic structures,* American Research Press, Rehoboth, 2003.

49. Vasantha Kandasamy, W.B., Bivector spaces, *U. Sci. Phy. Sci.*, 11, 186-190 1999.

50. Vasantha Kandasamy, W.B., *Linear Algebra and Smarandache Linear Algebra*, Bookman Publishing, 2003.

51. Vasantha Kandasamy, W.B., Smarandache, Florentin and K. Ilanthenral, *Introduction to bimatrices,* Hexis, Phoenix, 2005.

52. Vasantha Kandasamy, W.B., Smarandache, Florentin and K. Ilanthenral, *Introduction to Linear Bialgebra,* Hexis, Phoenix, 2005.

53. Vasantha Kandasamy, W.B., and Florentin Smarandache, Basic Neutrosophic Algebraic Structures and their Applications to Fuzzy and Neutrosophic Models, Hexis, Church Rock, 2005.

54. Vasantha Kandasamy, W.B., and Florentin Smarandache, *n-linear algebra of type I and its applications*, InfoLearnQuest, Ann Arbor, 2008.

55. Vasantha Kandasamy, W.B., and Florentin Smarandache, *n-linear algebra of type II,* InfoLearnQuest, Ann Arbor, 2008.

56. Vasantha Kandasamy, W.B., Smarandache, Florentin and K. Ilanthenral, *Special Set Linear Algebra and Special Fuzzy Linear Algebra*, Editura CuArt, 2010.

57. Voyevodin, V.V., *Linear Algebra*, Mir Publishers, 1983.

58. Zelinksy, D., *A first course in Linear Algebra*, Academic Press, 1973.

# INDEX

## B

C

# ABOUT THE AUTHORS

**Dr.W.B.Vasantha Kandasamy** is an Associate Professor in the Department of Mathematics, Indian Institute of Technology Madras, Chennai. In the past decade she has guided 13 Ph.D. scholars in the different fields of non-associative algebras, algebraic coding theory, transportation theory, fuzzy groups, and applications of fuzzy theory of the problems faced in chemical industries and cement industries.

She has to her credit 646 research papers. She has guided over 68 M.Sc. and M.Tech. projects. She has worked in collaboration projects with the Indian Space Research Organization and with the Tamil Nadu State AIDS Control Society. She is presently working on a research project funded by the Board of Research in Nuclear Sciences, Government of India. This is her $48^{th}$ book.

On India's 60th Independence Day, Dr.Vasantha was conferred the Kalpana Chawla Award for Courage and Daring Enterprise by the State Government of Tamil Nadu in recognition of her sustained fight for social justice in the Indian Institute of Technology (IIT) Madras and for her contribution to mathematics. The award, instituted in the memory of Indian-American astronaut Kalpana Chawla who died aboard Space Shuttle Columbia, carried a cash prize of five lakh rupees (the highest prize-money for any Indian award) and a gold medal.
She can be contacted at vasanthakandasamy@gmail.com
Web Site: http://mat.iitm.ac.in/home/wbv/public_html/

**Dr. Florentin Smarandache** is a Professor of Mathematics at the University of New Mexico in USA. He published over 75 books and 150 articles and notes in mathematics, physics, philosophy, psychology, rebus, literature.

In mathematics his research is in number theory, non-Euclidean geometry, synthetic geometry, algebraic structures, statistics, neutrosophic logic and set (generalizations of fuzzy logic and set respectively), neutrosophic probability (generalization of classical and imprecise probability). Also, small contributions to nuclear and particle physics, information fusion, neutrosophy (a generalization of dialectics), law of sensations and stimuli, etc. He can be contacted at smarand@unm.edu